\documentclass[11pt,a4paper,reqno]{amsart}%
\usepackage{amsthm,amsmath,amsfonts,amssymb,xcolor,amsxtra,appendix,bookmark,dsfont,bm,mathrsfs}
\usepackage[latin1]{inputenc}
\usepackage{xcolor} 
\usepackage{amssymb}\usepackage{graphicx}
\usepackage{tikz}
\usepackage{mathtools}
\usepackage{enumitem}
\usepackage{hyperref}
\usepackage[latin1]{inputenc}
\usepackage{color} 
\usepackage{amssymb}\usepackage{graphicx}
\usepackage{tikz}
\usepackage{hyperref}
\theoremstyle{plain}
\newtheorem{theorem}{Theorem}
\newtheorem{lemma}[theorem]{Lemma}

\newtheorem{proposition}[theorem]{Proposition}

\newtheorem{definition}{Definition}

\theoremstyle{remark}
\newtheorem{remark}[theorem]{Remark}

\allowdisplaybreaks

\title[Wave turbulence theory for    ZK equations]{On the wave turbulence theory for a stochastic   KdV type equation}

\author[G. Staffilani]{Gigliola Staffilani
}
\address{Department of Mathematics, Massachusetts Institute of Technology, Cambridge, MA 02139, USA}
\email{gigliola@math.mit.edu} 
\thanks{ G.S. is  funded in part by  the NSF grant DMS-1764403, DMS-2052651 and the Simons Foundation}

\author[M.-B. Tran]{Minh-Binh Tran}
\address{Department of Mathematics, Texas A\&M University, College Station, TX 77843, USA}
\email{minhbinh@tamu.edu} 
\thanks{M.-B. T is  funded in part by  the  NSF Grants DMS-1854453, DMS-2204795, DMS-2305523,    Humboldt Fellowship,   NSF CAREER  DMS-2044626/DMS-2303146.}

\begin{document}
	\date{\today}

	\begin{abstract} Starting from the   Zakharov-Kuznetsov (ZK) equation, a multidimensional KdV type equation  on a hypercubic lattice, we provide a derivation of the 3-wave kinetic equation. We  show that the two point correlation function  can be asymptotically expressed as  the solution of the 3-wave  kinetic equation at the kinetic limit under very general assumptions: the initial condition is out of equilibrium, the dimension  is  $d\ge 2$, the smallness of the nonlinearity $\lambda$ is allowed to be independent of the size $D$ of the lattice.   To the best of  our knowledge, the work provides the first rigorous derivation of  nonlinear 3-wave kinetic equations. Also this is the first derivation   for wave kinetic equations  in the lattice setting and  out-of-equilibrium.  
	\end{abstract}
	
	\maketitle

	\tableofcontents
	
	\section{Introduction}\label{intro}

	Given a periodic global solution $u(t,x)$ to a certain dispersive equation one is  interested in analyzing the migration of its energy from low to high frequencies (forward cascade), more precisely one wants to  study the behavior  {of the energy spectrum given by the magnitude of the Fourier coefficients $\hat u(t,k)$} as $t$ and $k$ become large. 
	Here we take the view point that results on forward cascade are based on two different approaches. The first approach, arguably introduced by Bourgain in \cite{bourgain1996growth}, is centered on the asymptotic analysis of the norm  $H^s, s\gg1,$ of the solution $u$ itself, since $\|u(t)\|_{H^s}^2=\sum_{k}|\hat u(t,k)|^2\langle k\rangle^{2s}$ and hence its growth would show that migration to higher frequencies occurred.  In within this approach great progress has been made 
	using techniques that are more commonly employed in a PDE  context, see for example \cite{bourgain1999growth,bourgain1999growtha}, \cite{carles2012energy}, \cite{colliander2010transfer}, \cite{deng2017strichartz}, \cite{hani2015modified}, \cite{kuksin1997oscillations},  \cite{majda1997one},  \cite{sohinger2011bounds}, \cite{staffilani1997growth},  \cite{staffilani2020stability},  or techniques that are framed in within a dynamical systems approach, such as in \cite{berti2019long}, \cite{guardia2015growth}, \cite{haus2015growth}.  The second approach is based on finding an effective equation, referred to as the wave kinetic equation,  for the expectation (with respect to a certain probability  density) of $|\hat u(t,k)|^2$ by approximating the original dispersive equation in various ways, and then implementing a limit process  to finally arrive at a new equation that in an appropriate sense is satisfied by the expectation of $|\hat u(t,k)|^2$. 
	
	The second approach is directly connected to the so-called wave turbulence theory in continuum mechanics 
	(see \cite{taoweakturbulence}). Wave turbulence  describes the dynamics of both classical and non-classical nonlinear waves  out of  thermal equilibrium. In spite of the enormous diversity of wave fields describing the processes of random wave interactions  in nature, there is a common mathematical framework that models and describes the dynamics of spectral energy transfer through  probability densities associated with weakly non-linear interactions in quantum or classical wave systems. Those probability densities are solutions of wave kinetic equations, whose nonlocal interaction operators are of kinetic-type. Among the several types of wave kinetic equations, there are two common ones: 3-wave and 4-wave equations  associated with  weakly nonlinear wave systems with quadratic and cubic nonlinearities respectively. Denoting $\lambda>0$ the  parameter that describes the weak interactions of the wave system under consideration, as the collision rate  is $\lambda^{2}$, it is expected that the associated wave kinetic equation can be derived at the {\it van Hove limit} or {\it the kinetic 
		limit  }
	\begin{equation}
		\label{VanHove}
		t\backsim \lambda^{-2}.
	\end{equation}
	The derivation of  kinetic equations has been one of the  central questions in the wave turbulence theory, a theory that has its origin in the works of  Peierls \cite{Peierls:1993:BRK,Peierls:1960:QTS}, Brout-Prigogine \cite{brout1956statistical},  Zaslavskii-Sagdeev \cite{zaslavskii1967limits}, Hasselmann \cite{hasselmann1962non,hasselmann1974spectral},  Benney-Saffman-Newell \cite{benney1969random,benney1966nonlinear},  Zakharov \cite{zakharov2012kolmogorov}, and has been playing  important roles in a vast range of physical applications: water surface gravity and capillary waves; inertial waves due to rotation and internal waves on density stratifications, which are important in planetary atmospheres and oceans; Alfv\'en wave turbulence in solar wind; planetary Rossby waves, which are important for the weather and climate evolutions; waves in Bose-Einstein condensates (BECs) and in nonlinear optics; waves in plasmas of fusion devices; and many others, as discussed in the books of  Zakharov, Lvov, Falkovich \cite{zakharov2012kolmogorov} and Nazarenko
	\cite{Nazarenko:2011:WT}, and the review paper of Newell and Rumpf \cite{newell2011wave}. 
	
	Despite the growing number of physical applications of wave turbulence theory, the mathematically rigorous derivation of wave kinetic equations for different types of weakly nonlinear dispersive wave systems has  always been a challenging problem.

	In rigorously deriving wave kinetic equations,  the work  of Lukkarinen and Spohn \cite{LukkarinenSpohn:WNS:2011} for the lattice cubic nonlinear Schr\"odinger equation (NLS) is pioneering. In \cite{LukkarinenSpohn:WNS:2011},   the random initial condition is distributed according to the corresponding Gibbs measure at statistical equilibrium. The equation is put on a hypercubic lattice and kept deterministic. The linearization of  the 4-wave kinetic operator can then  be observed via the use of a special time correlation function at the kinetic  limit.

	To derive kinetic equations from wave equations, the wave equations can be  put either on  hypercubic lattices or on continuum domains. The lattice and continuum settings have different challenges and difficulties (see the discussions in  \cite{chen2005localization,chen2005localization2,lukkarinen2007asymptotics} and Subsection \ref{Sec:Novelty}). 
	
	While the work of Lukkarinen and Spohn is in the lattice setting,	the work of Buckmaster-Germain-Hani-Shatah \cite{buckmaster2019onthe,buckmaster2019onset}, Deng-Hani \cite{deng2019derivation}, and  Collot-Germain \cite{collot2019derivation,collot2020derivation} give  rigorous derivations of the 4-wave kinetic equations  from the cubic NLS equation in the continuum setting,  with random initial data, at limits closed to \eqref{VanHove} and the equation is put out of statistical equilibrium.  Attempts to derive the 4-wave kinetic equation from the stochastic NLS in the continuum setting have been carried out  by Dymov, Kuksin and collaborators in \cite{dymov2019formal,dymov2019formal2,dymov2020zakharov,dymov2021large}.
	
	Finally, let us also mention the (CR) equation,  which is derived from the deterministic cubic nonlinear Schr\"odinger equation in~\cite{buckmaster2016effective,faou2016weakly} and studied in~\cite{buckmaster2016analysis,germain2015continuous,germain2016continuous},  and which is a Hamiltonian equation whose nonlinearity is given by one component of the 4-wave kinetic operator. 
	
	Let us also comment that the Boltzmann equation is probably the most famous collisional kinetic model. The derivation of the Boltzmann equation is also a very important and active research line  \cite{bodineau2016brownian,gallagher2013newton,lanford1975time}.
	
	In our present work, we propose a derivation of the full 3-wave kinetic equation in the same lattice setting considered previously by Lukkarinen and Spohn \cite{LukkarinenSpohn:WNS:2011}.  We will start  with the ZK  equation in $d$-dimension ($d\ge2$), which has many physical applications concerning  drift waves in fusion plasmas  and
	Rossby waves in geophysical fluids  \,\cite{Nazarenko:2011:WT}[Section 6.2, Equation (6.1)] and ionic-sonic waves in a magnetized plasma \cite{lannes2013cauchy,zakharov1974three}, 
	\begin{equation}\label{KleinGordonNoise}
		\begin{aligned}
			\mathrm{d}\psi(x,t) \ & = \ -\Delta\partial_{x_1} \psi(x,t) \mathrm{d}t  \ + \ \lambda\partial_{x_1}\Big(\psi^2(x,t)\Big)\mathrm{d}t,\\
			\psi(x,0) \ & = \ \psi(x),
		\end{aligned}
	\end{equation}
with a random noise that will be specified later. The parameter $1>\lambda>0$ is  a real   constant. The constant $\lambda$ is the parameter describing the weak interactions of the nonlinear wave system and as mentioned above 
	we will  send $\lambda$ to $0$. {\it Notice that the ZK  equation  has normally served as the first example for which a 3-wave kinetic equation is derived \cite{Nazarenko:2011:WT}.}
	Our equation is defined on a hypercubic lattice $\Lambda$, of size $D$ that will be defined more precisely in the next section. 
	
	\medskip
	Let us now give a very rough statement of the main theorem that will be properly stated in Theorem \ref{TheoremMain} in Section \ref{Sec:Main}  after all the necessary notations are introduced. 
	\begin{theorem}[Rough Statement of the Main Theorem \ref{TheoremMain}]\label{TheoremMainRough}
		Suppose that $d\ge 2$. The two point correlation function of the solution   for   \eqref{KleinGordonNoise} with a stochasticity on a hypercubic lattice, after we take the limit $D\to\infty$,  can be asymptotically expressed via the solution of a 3-wave  kinetic equation at the kinetic time \eqref{VanHove} in the resonance broadening sense, under general assumptions on the initial data.  
	\end{theorem}

	We  now present a list of remarks concerning the main theorem and its proof. 

	The  proof of our main result is built by using as a starting point the   Erdos-Yau techniques \cite{erdHos2000linear} (see also \cite{erdHos2002linear,yau1998scaling}), introduced    in their study of the linear Schr\"odinger equation with a weak random potential in both lattice and continuum settings, and later developments by Erdos-Salmhofer-Yau \cite{erdHos2008quantum} and Lukkarinen-Spohn \cite{LukkarinenSpohn:WNS:2011}. One main component of the proof is the expansion of the solution as a  power series in terms of the initial data,  and the organization of this expansion in Feynman diagrams. Since the Duhamel formula has the advantage that it can be stopped at a different number of terms
	depending on the collision history of each elementary wave function, one imposes a {\it stopping rule}, first introduced in \cite{erdHos2000linear} and widely used later \cite{butz2015kinetic,chen2005localization,chen2006convergence,erdHos2008quantum, LukkarinenSpohn:KLW:2007,LukkarinenSpohn:WNS:2011}, in which the expansion is only done up to a given number of layers $\mathfrak{N}$, which is a function of the smallness parameter $\lambda$.  Most of the Feynman diagrams, after being integrated out, produce positive powers $\lambda^\theta, \theta>0$ of the small parameter $\lambda$ and hence become very small as $\lambda$ approaches $0$. The remaining  diagrams have very  special  structures: they are self-repeated. This phenomenon was first discovered for both lattice and continuum settings by Erdos-Yau  \cite{erdHos2000linear} in the context of the linear Schr\"odinger equation with random potential, and later  developed by Lukkarinen-Spohn \cite{LukkarinenSpohn:WNS:2011} to the more sophisticated case of the lattice nonlinear Schr\"odinger equation with random initial condition. The fact that most of the Feynman diagrams, after being integrated out, are negligible in the limit $\lambda\to 0$, is  often referred to as the {\it suppression of crossings} (see (ii) below for more details). 
	Following the discussion by  Erdos for the linear random Schr\"odinger equation \cite{erdos2010lecture,erdos2020}, we call the self-repeating structures/self-repeated diagrams   {\it ladder (leading) diagrams}. The combination of the self-repeating structures gives us a power series form of  solutions to the 3-wave kinetic equation. However, to achieve the full power of the techniques we use we need to overcome  several obstacles, the   resolution of which  require interdisciplinary tools. {\it To the best of our knowledge  the present work is  the first rigorous derivation of  nonlinear 3-wave kinetic equations, and also the first out-of-equilibrium derivation for wave kinetic equations in the lattice setting.  } 
	In the next subsection we  highlight the novelties of our approach.

	\subsection{Novelties of the approach}\label{Sec:Novelty}
	~~~\\~~~~

	{\bf (i) Clustering/moment estimates and the advantage of the Liouville equation. } One of the main  obstacles in controlling the Feynman diagrams  concerns the control of the error terms, which involve the full original time evolution. Erdos and Yau, who work in the linear setting, control the error terms using the unitarity of the evolution.  In the work of Lukkarinen and Spohn, since the equation is nonlinear, the control comes from imposing the assumption in  \cite{LukkarinenSpohn:WNS:2011}[Assumption 2.1] on the $l^1$-clustering estimate at equilibrium, (see also \cite{lukkarinen2009not}[Assumption A2] in the 
	context of  weakly interacting lattice quantum fluids).  The $l^1$-clustering estimate is a well-known and difficult problem in statistical physics that concerns bounds for the cumulants (see \cite{ginibre1971some,ruelle1964cluster,ruelle1999statistical,salmhofer2009clustering} and the references therein). Indeed, one of the main technical obstacles that forces Lukkarinen and Spohn to put the system at equilibrium is the difficulty in having the $l^1$-clustering estimate out of equilibrium.  While \cite{LukkarinenSpohn:WNS:2011}[Assumption 2.1] has been proved only for the zero boundary conditions and at equilibrium case 
	(see \,\cite{abdesselam2009clustering}), establishing{ such a clustering estimate} out of equilibrium  for the   ZK equation under consideration is still an open problem. To overcome this technical difficulty, we focus on the analysis of the Liouville equation of the density function under the presence of the noise.

	\medskip
	
	{\bf (ii) New types of crossing estimates - Challenges of the lattice setting. } As discussed above, in rigorously deriving kinetic equations from wave systems \cite{bal2003self,basile2016thermal,komorowski2020high,ryzhik1996transport}, one important key step is to show that most of the Feynman diagrams, after being integrated out, are negligible in the limit $\lambda\to 0$, leading to the dominance of the ladder diagrams \cite{chen2005localization,chen2006convergence,erdHos2008quantum,erdHos2000linear, LukkarinenSpohn:KLW:2007,LukkarinenSpohn:WNS:2011}.  
	An important element in all of these results is an estimate proving that all so-called crossing graphs are suppressed.
	For the   linear Schr\"odinger equation with a random potential, the estimate takes the form
	\begin{equation}
		\begin{aligned}\label{CrossingIntro1}
			& \ \sup_{(\alpha_1,\alpha_2,\alpha_3)\in\mathbb{R}^3} 
			\int_{(\mathbb{T}^d)^2}\frac{\mathrm{d}k_1\mathrm{d}k_2}{|\alpha_1-\omega(k_1)+{\bf i}\lambda^2||\alpha_2-\omega(k_2)+{\bf i}\lambda^2||\alpha_3-\omega(k_1-k_2+k_0)+{\bf i}\lambda^2|} \lesssim \langle \ln \lambda\rangle^{\gamma_a}\lambda^{\gamma_b},
		\end{aligned}
	\end{equation}
	for some constants $\gamma_a,\gamma_b\in\mathbb{R}$. Estimates of the type \eqref{CrossingIntro1} are often referred to as the ``crossing estimates''. As discussed in  \cite{lukkarinen2007asymptotics}, the validity of the corresponding estimate in the  continuum Schr\"odinger  setting, with { the dispersion relation}
	\begin{equation}\label{NLSDispersionCont}\omega(k)=|k|^2, \ \ \ \  k\in\mathbb{R}^d\end{equation}
	is fairly straightforward to prove, but the  lattice   case  turns out to be much more involved since
	\begin{equation}\label{NLSDispersionDiscrete}\omega(k)= \sin^2(2\pi   k^1) + \cdots + \sin^2(2\pi   k^d), \quad \mbox{ for } k=(k^1,\cdots ,k^d).\end{equation}
	The key problem of the lattice dispersion relation, which does not appear in the continuum case, is that there exist critical energy values where the energy surface has segments of zero Gauss curvature. This technical issue of the lattice Schr\"odinger  dispersion relation  is unavoidable for topological reasons, and has been first observed by Bourgain \cite{bourgain2002random}. Since, in $3$ dimension, \eqref{NLSDispersionDiscrete} is a Morse function on the torus, the first Betti number of the $3$-torus is $3$, which implies that the level surfaces of \eqref{NLSDispersionDiscrete} make a transition from the topological sphere to a genus $3$ surface, and back to the topological sphere. The bound in \eqref{CrossingIntro1} has been proved 
	to hold for the lattice Schr\"odinger dispersion relation    with  $\gamma_b = -4/5$ and $\gamma_a = 2$ by Chen  \cite{chen2005localization} (see also \cite{chen2005localization2,chen2006convergence,chen2011boltzmann}) and with $\gamma_b = -3/4$ and $\gamma_a = 6$ by Erdos-Salmhofer-Yau  \cite{erdHos2008quantum}. In \cite{erdHos2007decay}, the ``four denominator estimate'', which involves four resolvent terms
	instead of three as in \eqref{CrossingIntro1} has been studied. 
	In a later important work  \cite{lukkarinen2007asymptotics}, Lukkarinen has considered the case of more general dispersion relations $\omega$. He discovered that for the general set up, for  small $\omega$, each of the factors in \eqref{CrossingIntro1} is sharply concentrated around some level set of $\omega$, while  the arguments of $\omega$ are not allowed to vary independently of each other, and the magnitude of the integral for small $\lambda$ is thus determined by the overlap of the different level sets depending on the constants $\alpha_i$. Therefore, to establish \eqref{CrossingIntro1} it is important to consider the worst case scenario for the level sets, and  estimate the overlap between the { three level sets} as $k_1$ and $k_2$ are varied. {\it It has been proved by Lukkarinen \cite{lukkarinen2007asymptotics} that an analytic dispersion relation suppresses crossings if and only if it is not a constant on any affine hyperplane.  A counterexample, in which the crossing estimate \eqref{CrossingIntro1} fails to hold true, has also been introduced, which unfortunately covers the lattice ZK dispersion relation 	\begin{equation}\label{KDVDispersion}
			\omega_k \ = \ \omega(k) \ = \ \sin(2\pi   k^1)\Big[\sin^2(2\pi   k^1) + \cdots + \sin^2(2\pi   k^d)\Big], \quad\mbox{ for } k=(k^1,\cdots,k^d).
		\end{equation} To be more precise, the suppression of crossings does not hold for the dispersion relation \eqref{KDVDispersion},  due to the fact that the lattice ZK dispersion relation 
		vanishes, i.e.  $\omega(k)=0$, on the affine hyperplane $k=(0,k^2,\cdots,k^d)$.} Let us also comment that the setting in which the topological issue discussed above is avoidable and  only the low energy component is picked up is that of a lattice discretized nonlinear dispersive equation for which the lattice spacing converges to zero in an appropriate limit, such as in the work by Kirkpatrick-Lenzmann-Staffilani \cite{kirkpatrick2013continuum}.
	
	In the context of the lattice nonlinear Schr\"odinger equation considered in \cite{LukkarinenSpohn:WNS:2011}, the crossing estimate takes the form
	\begin{equation}
		\begin{aligned}\label{CrossingIntro2}
			& \ \sup_{(\alpha_1,\alpha_2)\in\mathbb{R}^2} \int_{(\mathbb{T}^d)^2}\mathrm{d}k_1\mathrm{d}k_2
			\frac{1}{\Big|\alpha_1\pm  \omega(k_1)\pm  \omega(k_2)\pm  \omega(k_3)+{\bf i}\lambda^2\Big|}\\
			&\ \ \ \times \frac{1}{\Big|\alpha_2\pm  \omega(k_1)\pm  \omega(k_2)\pm  \omega(k_3)\pm  \omega(k_1+k_4)  \pm \omega(k_1+k_5)+{\bf i}\lambda^2\Big|} \lesssim \langle \ln \lambda\rangle^{\gamma_a}\lambda^{\gamma_b},
		\end{aligned}
	\end{equation}
	where $\omega$ is the nearest neighboring dispersion relation of the Schr\"odinger operator defined in  \eqref{NLSDispersionDiscrete}. 
	Crossing estimates of type \eqref{CrossingIntro2}  have only two     denominators instead of three as in \eqref{CrossingIntro1} and they are  clearly more complicated than those of  type \eqref{CrossingIntro1}. Therefore, even though the strategy is still to show the dominance of the leading (ladder) diagrams, the classification of other types of graphs (crossing and nested diagrams) is much more complicated and involved. Due to  the above topological issue of the lattice Schr\"odinger dispersion relation \eqref{NLSDispersionDiscrete} observed by Bourgain, the result of \cite{LukkarinenSpohn:WNS:2011} is restricted to dimension $d\ge 4$.

	At this point it should be clear to the reader that in the context of the lattice ZK equation, the crossing estimates are more delicate due to the singularities coming from the dispersion relation. However, the quadratic nonlinearity in our ZK equation also creates an extra problem in obtaining the suppression of crossings, as explained below. One would like to prove  crossing estimates for the lattice ZK equation of  the following form
	\begin{equation}
		\begin{aligned}\label{CrossingIntro3}
			& \ \sup_{(\alpha_1,\alpha_2)\in\mathbb{R}^2}\int_{\mathbb{T}^d}\mathrm{d}k_1
			\frac{1}{\Big|\alpha_1\pm  \omega(k_1)\pm  \omega(k_1+k_2) +{\bf i}\lambda^2\Big|}\\
			&\ \ \ \times \frac{1}{\Big|\alpha_2\pm  \omega(k_1)\pm  \omega(k_1+k_2)\pm  \omega(k_1+k_3) +{\bf i}\lambda^2\Big|} \lesssim \langle \ln \lambda\rangle^{\gamma_a}\lambda^{\gamma_b}.
		\end{aligned}
	\end{equation}
	{\it To the best of  our knowledge, this is the first time that a crossing estimate of  type \eqref{CrossingIntro3} is considered.}
	{\it In comparison to the linear Schr\"odinger \eqref{CrossingIntro1} and the nonlinear Schr\"odinger crossing estimates \eqref{CrossingIntro2}, one can see immediately that the  lattice ZK crossing estimate  \eqref{CrossingIntro3} would be  harder to obtain as it only contains two denominators instead of three denominators and one integration in $k_1$ instead of two integrations in $k_1,k_2$.} The loss of one integration is due to the nature of the quadratic nonlinearity. We will explain more below that the crossing  estimate \eqref{CrossingIntro3} is harder to obtain in general and fails to hold true for the specific case of the lattice ZK dispersion relation.

An ``easy-to-see'' technical difficulty in  obtaining the crossing estimate \eqref{CrossingIntro3} can be explained as follows. To obtain the crossing estimate \eqref{CrossingIntro1} for the NLS, one could use an $L^3$ estimate
\begin{equation}
	\begin{aligned}\label{CrossingIntro2:1}
		& \ \Big|\int_{(\mathbb{T}^d)^2}\mathrm{d}k_1\mathrm{d}k_2
		e^{{\bf i} t(\pm  \omega(k_1)\pm  \omega(k_2)\pm  \omega(k_3))+{\bf i} s( \omega(k_1+k_4)  \pm \omega(k_1+k_5)) }\Big|
		\lesssim	& \|p_t\|_3^2\|\mathscr{K}(\pm t,\pm s,\pm s, k_4,k_5)\|_3,
	\end{aligned}
\end{equation}
thanks to the presence of the double integral $\int_{(\mathbb{T}^d)^2}\mathrm{d}k_1\mathrm{d}k_2$, where 
$$\mathscr{K}(x,t_0,t_1,t_2, k,k_*):=e^{-{\bf i}d(t_0+t_1+t_2)}\prod_{i=1}^d\int_{0}^{2\pi}\frac{\mathrm{d} p}{2\pi}e^{{\bf i}px^i+{\bf i}(t_0\cos(2\pi p)+t_1\cos(2\pi p+k^i)+t_2\cos(2\pi p + k_*^i ))},$$
with $k=(k^1,\cdots,k^d)$ and $k_*=(k^1_*,\cdots,k^d_*),$ and $p_t=\mathscr{K}(x,t,0,0, 0,0)$.  As a result, in order to obtain the crossing estimate \eqref{CrossingIntro3}, which only involves one integration,  a straightforward bound would be an $L^2$ estimate, which is definitely not satisfactory. Indeed, the $L^2$ norm is conserved and has no decay in time, {\it which leads to the divergence of the sum of all the   diagrams. We will discuss later in (iv) that the noise has no influence on pairing graphs, and therefore the restoration of  the convergence of those diverging
	     diagrams is indeed a challenging issue for any dispersion relations. } 

However, in the special case of  the ZK dispersion relation \eqref{KDVDispersion} under consideration, the problem is even more serious.
As discussed above, Lukkarinen \cite{lukkarinen2007asymptotics} has pointed out a counterexample, in which the  crossing estimate \eqref{CrossingIntro1} of the linear case fails to hold true for the lattice ZK dispersion relation \eqref{KDVDispersion}. {\it   As thus, the harder crossing estimate \eqref{CrossingIntro3} also fails to hold true for the lattice ZK dispersion relation. A main portion of our work is dedicated to  establishing several different and new types  of crossing estimates, which are more flexible than those used in the previous works,  for the singular  lattice ZK dispersion relation  under the low dimensional assumption $d\ge2$. These  novel types of crossing estimates allow us to go around the situation encountered previously in Lukkarinen's counterexample \cite{lukkarinen2007asymptotics}  and are embedded into new (and sophisticated) types of graph estimates. This part of the paper  is very much involved and way more complicated than the original situation of the lattice Schr\"odinger operator encountered by Bourgain \cite{bourgain2002random}.} Let us also point out that, similar in spirit but  simpler dispersive estimates have also been studied in a recent work of Grande-Kurianski-Staffilani \cite{grande2021nonlinear}.

\medskip

{\bf (iii) Resonance broadening/Quasi-resonance. }
{\it Unlike the Schr\"odinger dispersion relation, the  lattice ZK dispersion relation  not only creates major obstacles in obtaining the  crossing estimates but also has another serious problem: it prevents the convergence of the leading diagrams}.  In most formal derivations of wave kinetic equations written in physics books \cite{Nazarenko:2011:WT}, it is common practice to assume that

\begin{equation}
	\begin{aligned}\label{PhysicalDelta}
		&	\int_{\mathbb{T}^{2d}}\mathrm{d}k_2\mathrm{d}k_3\delta (\omega(k_3)+\omega(k_2)-\omega(k_2+k_3))F(k_2+k_3,k_2,k_3)\\
		\ = \ &  \int_{\mathbb{R}}\mathrm{d}s\int_{\mathbb{T}^{2d}} \mathrm{d}k_2\mathrm{d}k_3 e^{-{\bf i}s(\omega(k_3)+\omega(k_2)-\omega(k_2+k_3))}F(k_2+k_3,k_2,k_3),\end{aligned}
\end{equation}
for any test function $F(k_2+k_3,k_2,k_3)\in C^\infty(\mathbb{T}^{2d})$ and for any dispersion relation $\omega$. Unfortunately, it turns out that \eqref{PhysicalDelta} only holds true under some assumptions on $\omega$, while for most dispersion relations the quantity $ \delta (\omega(k_3)+\omega(k_2)-\omega(k_2+k_3))$ cannot be defined as a positive  measure. 
Let us define  \begin{equation}\label{Eq:CollisionWT6}
	\begin{aligned}
		&{\delta}_\ell (\omega(k_3)+\omega(k_2)-\omega(k_1)) \ := \ \frac{	\ell}{  \ell^2+(\omega(k_3)+\omega(k_2)-\omega(k_1))^2 },\ \ \ \ \ell >0.
\end{aligned}\end{equation} 
One writes
\begin{equation}
	\label{remark:TheoremMain:1}\begin{aligned}
		&	\int_{\mathbb{T}^{2d}}\mathrm{d}k_2\mathrm{d}k_3{\delta}_\ell (\omega(k_3)+\omega(k_2)-\omega(k_2+k_3))F(k_2+k_3,k_2,k_3)\\
		\ = \ & \int_{\mathbb{R}}\mathrm{d}se^{-|s|\ell} \int_{\mathbb{T}^{2d}} \mathrm{d}k_2\mathrm{d}k_3 e^{-{\bf i}s(\omega(k_3)+\omega(k_2)-\omega(k_2+k_3))}F(k_2+k_3,k_2,k_3).\end{aligned}
\end{equation}

When   $\omega$ is sufficiently good,  the oscillatory integral $\int_{\mathbb{T}^{2d}}\mathrm{d}k_2\mathrm{d}k_3$ $ e^{-{\bf i}s(\omega(k_3)+\omega(k_2)-\omega(k_2+k_3))}F(k_2+k_3,k_2,k_3)$ produces sufficient decay in $s$, yielding the convergence of $\delta_\ell (\omega(k_3)+\omega(k_2)-\omega(k_2+k_3))$ to the positive  measure $\delta (\omega(k_3)+\omega(k_2)-\omega(k_2+k_3))$ in the limit $\ell\to 0$ (see  \cite{LukkarinenSpohn:KLW:2007}[Proposition A.1] and \cite{LukkarinenSpohn:WNS:2011}[Proposition 7.4] for the proofs concerning the Schr\"odinger case). 

When $\omega$ is the lattice ZK dispersion relation, indeed  the delta function $\delta (\omega(k_3)+\omega(k_2)-\omega(k_2+k_3))$ cannot be defined as a positive  measure, yielding the divergence of the leading graphs, that contain oscillatory integrals of the form (with $\ell=0$) $\int_{\mathbb{T}^{2d}}\mathrm{d}k_2\mathrm{d}k_3$ $ e^{-{\bf i}s(\omega(k_3)+\omega(k_2)-\omega(k_2+k_3))}$ $F(k_2+k_3,k_2,k_3).$  

In our derivation of the 3-wave kinetic equation from the lattice ZK equation, for some diagrams, the quantity $F(k_2+k_3,k_2,k_3)$ is also allowed to vanish on the hyperplane $k=(0,k^2,\cdots,k^d)$, yielding the hope that one can use  $F(k_2+k_3,k_2,k_3)$ to control some singularities coming from the oscillatory integrals. However, even with the help of $F(k_2+k_3,k_2,k_3)$, the oscillatory integrals are so singular that the convergence of the leading diagrams still cannot be restored. In general, the lattice ZK dispersion relation creates serious divergent issues for both leading and non-leading diagrams. {For this reason in our prof we replace $\delta (\omega(k_3)+\omega(k_2)-\omega(k_2+k_3))$  with $\delta_\ell (\omega(k_3)+\omega(k_2)-\omega(k_2+k_3))$ for $\ell>0$ arbitrarely small.)}

The replacement of $\delta (\omega(k_3)+\omega(k_2)-\omega(k_2+k_3))$ by $\delta_\ell (\omega(k_3)+\omega(k_2)-\omega(k_2+k_3))$   is often referred to as {\it resonance broadening/quasi-resonance} and is justified when $\omega$ is singular or degenerate (which is the case of the lattice ZK dispersion relation under consideration) in the physical context \cite{chekhlov96,huang2000,lee2007formation,l1997statistical,lvov2004hamiltonian,lvov2010oceanic,lvov2012resonant,remmel2014nonlinear,smith2001,smith2005near,smith1999transfer,smith2002generation}. Recently, in the context of oceanography the idea of  resonance broadening, has  also been used to obtain  a well-posedness result  for a  3-wave kinetic equation  by   \cite{GambaSmithBinh}. {\it To the best of our knowledge, our work provides the first rigorous justification of the resonance broadening, discussed previously in  physical works.}

\medskip
{{{\bf (iv) Justifying the addition of the noise - Can a standard additive noise be used?} 
	Firstly, we need to justify why the noise is needed. From a technical point of view, in  Lukkarinen's counterexample,  the crossing estimate \eqref{CrossingIntro1} fails to hold true for the lattice ZK dispersion relation \eqref{KDVDispersion} while  \eqref{CrossingIntro1}  holds true for both the continuum and lattice Schr\"odinger dispersion  relations. { It is therefore straightforward that  the harder crossing estimate \eqref{CrossingIntro3} fails to hold true for the lattice ZK dispersion relation \eqref{KDVDispersion}, unfortunately.} From a deeper point of view, we know that if $k=(k^1,k^2,..,k^d)$ and  $k^1=0$, then $\omega(k)=0$ when $\omega$ is of the type \eqref{KDVDispersion}. We call the degenerate surface for which $k^1=0$ the {\it ghost manifold}. If we take $k_1, \cdots, k_m$ in the ghost manifold, it follows that $\omega(k_1)=\cdots=\omega(k_m)=0.$ Moreover, the sum vector $k_1 + \cdots + k_m$ is also in the ghost manifold, leading to $$\omega(k_1 + \cdots + k_m)=\omega(k_1)+\cdots+\omega(k_m).$$  This shows that on the ghost manifold, not just $3$-wave interactions, but any $m$-wave interactions are also allowed, with $m\ge3$. {\it Due to the resonance broadening effect    discussed above, all quasi-resonance $m$-wave interactions can also happen in a small neighborhood near the ghost manifold.   Therefore,  the appearance of the ghost manifold, which does not exist in the Schr\"odinger case  \cite{lukkarinen2007asymptotics,LukkarinenSpohn:WNS:2011},  destroys the structure of 3-wave interactions and as thus there is no 3-wave kinetic equation for the ZK equation without noise. The introduction of the noise into the equation, which vanishes in the limit $\lambda\to0$,  has the  role of only dealing with this ghost manifold: it  identifies those waves that  accidentally fall into the ghost manifold, and it  kicks them out of the set by only randomizing their phases and not their amplitudes. Except  for the special role of controlling the singular dispersion relation  on the ghost manifold, the noise   has no influence on most of the Feynman diagrams, including the most important ones: the leading diagrams.  Thus, unlike the Schr\"odinger case, a similar result on the derivation of the 3-wave kinetic equation in its current strong form should not be expected from the lattice ZK equation under very general assumptions on the initial data  without the noise.}}
	
	Secondly, one may wonder if a standard additive noise can be used in the ZK equation. Indeed, all types of noise will have  diffusive effects and therefore, they will all introduce diffusive terms into the Liouville equation. However, as the 3-wave kinetic equation conserves energy and momentum, it would be important that the noise  does not inject energy into the wave system and the conservation laws are preserved under the effect of the noise. A standard additive noise is therefore not preferable as it introduces extra energy into the equation. Another question is if one could replace the convolutive noise $\partial_{x_1}\phi\odot\mathrm{d}W(t)$ { that we use} by the multiplicative noise $\partial_{x_1}\phi\circ\mathrm{d}W(t)$, similarly to what has been done in \cite{hannani2021stochastic} by Hannani and Olla, as the multiplicative noise also preserves all the conservation laws. Indeed, as the multiplicative noise affects not only the phases, a new additional collision operator  { would be introduced into the kinetic equation in this case,} and it { would change the structure} of the kinetic wave equation. Thus, we are not going to pursue the study in this direction here.  Among all the noises,  our  convolutive noise has the special effect that it randomizes only the phases without injecting any energy into the wave system, and therefore, it is quite weak but sufficiently good to handle the ghost manifold.

    	When adding noises in any weak wave systems, another question one might ask is whether or not the effect of the noise dominates the effect of the weak nonlinearity.   {\it We will prove later in the manuscript that the influence of the noise \eqref{KleinGordonNoise} almost disappears in all  pairing graphs as well as in all  leading (ladder) diagrams. As discussed above, those  diagrams are where the nonlinearity affects the most, thus there is no competition between the weak nonlinearity and the weak noise. Therefore, it is not a concern whether or not the noise is strong (or weak) in comparison with the  nonlinearity.} Indeed, {except for the special } effect on the singularities of the dispersion relation, the noise does not affect the structure of the Feynman diagrams, and  the rest of the proof holds independently of the noise.

\medskip

To conclude this part of the introduction, we would like to remark that a first version of  our result was announced via various talks of the authors and a first draft of the manuscript was circulating among some selected researchers for feedback since December 2020. 
As we were finalizing our manuscript based on the feedbacks received, being unaware of the full content of our manuscript,  Deng and Hani  posted a deep work \cite{deng2021full}, in which the kinetic time limit was also  reached for the cubic NLS  equation in the continuum setting. While the approach of \cite{deng2021full} is based on the expansion of the solution, our approach  involves the Liouville   equation of the density function\footnote{This is similar to what has been done often in the physical literature (see \cite{akhiezer2013methods,reichl2016modern}).}. As discussed above, the Schr\"odinger dispersion relation in the continuum setting has very different challenges and difficulties in comparison with the ZK dispersion relation in the lattice setting. Moreover,  the kinetic limit of Deng and Hani is obtained via the assumption that the smallness of the nonlinearity $\lambda$ depends on the size of the domain  of the dispersive equation, while our kinetic limit holds when  $\lambda$ is kept independent of the size of the domain. Hence, the  two definitions of the kinetic limits of the two results, obtained at the same time, give  quite broad  ranges of parameters at which one could obtain  wave kinetic equations. Moreover, in \cite{ampatzoglou2021derivation}, a derivation of the inhomogeneous 4-wave kinetic equations in  limits close to the kinetic time from the quadratic NLS equation has been provided. Finally, we want to mention the breakthrough new work of Lukkarinen and Vuoksenmaa  \cite{LukkarinenVuoksenmaa}     for the NLS equation in the lattice setting based on the powerful tool of the cumulant expansion.  In this highly important work, the desirable derivation of the homogeneous 4-wave kinetic equation has been provided under the  condition that $\lambda$ is independent of the size of the domain and the noise is not needed.  
\medskip

In the next section, we are going to introduce the framework for our problem, the complete statement of the main theorem and remarks on the main result.

\bigskip

\bigskip

{\bf Acknowledgments:}   The authors would like to express their gratitude to  Tristan Buckmaster, Thomas Chen, J. Robert Dorfman, Laszlo Erdos, Jani Lukkarinen, Yuri Lvov,  Sergey Nazarenko, Alan Newell, Stefano Olla,  Yves Pomeau, Linda Reichl, Lenya Ryzhik, Avy Soffer, Eitan Tadmor  and Hong-Tzer Yau  for  the inspiring  discussions on  wave turbulence and kinetic theories. They would also like to thank   Erwan Faou for the improvement  of Assumption (B),  Alessio Figalli for the  fruitful email exchanges on optimal transport theory and equation \eqref{FokkerPlanck}, Herbert Spohn for pointing out  useful references, Zaher Hani and Yu Deng for the comparisons between theirs and our   results.  
\medskip

\section{Settings and main result}\label{Sec:Main} 
\subsection{The Set Up }\label{Sec:Setting}

Our strategy to   derive  the 3-wave kinetic equation from the stochastic equation of \eqref{KleinGordonNoise} is based on the following lattice setting, similar to \cite{LukkarinenSpohn:WNS:2011}:
\begin{itemize}
	\item[(i)] Firstly, we write \eqref{KleinGordonNoise} using a finite difference scheme and derive a kinetic equation for the  equation on the lattice.
	\item[(ii)] Secondly, we send the size $D$ of the domain to infinity,   and obtain the desired continuous kinetic equation at the kinetic limit $t=\mathcal{O}(\lambda^{-2})$. 
\end{itemize}

We now introduce our finite volume mesh, namely
\begin{equation}
	\label{Lattice}
	\Lambda \ = \ \Lambda(D) \ = \ \left\{0,1 \dots,2D
	\right\}^d, 
\end{equation}
for some constant $D\in\mathbb{N}$.   

The dynamics for the discretized  equation   reads
\begin{equation}
	\label{LatticeDynamics}\begin{aligned}
		\mathrm{d}\psi(x,t) \ = &  \ \ \sum_{y\in  \Lambda} O_1(x-y)\psi(y,t)\mathrm{d}t \ + \ \lambda\sum_{y\in  \Lambda} O_2(x-y)  \psi^2(y,t)\mathrm{d}t, \\
		\psi(x,0) \ & = \ \psi_0(x), \ \forall (x,t)\in\Lambda\times \mathbb{R}_+,\end{aligned}
\end{equation}
where $O_1(x-y)$ and $O_2(x-y)$ are finite difference operators that we will express below in the Fourier space. The boundary condition is assumed to be periodic. To obtain the lattice dynamics,  we  introduce  the Fourier transform 
\begin{equation}
	\label{Def:Fourier}\hat  \psi(k)=\sum_{x\in\Lambda} \psi(x) e^{-2\pi {\bf i} k\cdot x}, \quad k\in \Lambda_* = \Lambda_*(D) =  \left\{-\frac{D}{2D+1},\cdots,0,\cdots,\frac{D}{2D+1}\right\}^d,
\end{equation}
at the end of this standard procedure, \eqref{LatticeDynamics}  can  be rewritten in the Fourier space as a system of ODEs
\begin{equation}
	\label{LatticeDynamicsFourier}\begin{aligned}
		\frac{\mathrm{d}\hat\psi(k,t)}{\mathrm{d} t} \ & = \ {\bf i} \omega(k)\hat\psi(k,t) \ + \     {\bf i}\lambda \bar\omega(k)\frac{1}{|\Lambda_*|^2}\sum_{k=k_1+k_2;k_1,k_2\in\Lambda_*}\hat\psi(k_1,t)\hat\psi(k_2,t),\\
		\hat\psi(k,0) \ & = \ \hat\psi_0(k). \end{aligned}
\end{equation}

In the above expression, the conservation of momentum is understood in terms of modulo the lattice $\Lambda_*$, which can be written as 
\begin{equation}
	\label{ModeZ}
	k=k_1+k_2 \mod{\Lambda_*} \Longleftrightarrow \vec{{V}}\in \mathbb{Z}^d \mbox{ such that } k=k_1+k_2 +  \vec{{V}}.
\end{equation}
We also set 
\begin{equation} \label{positiveDomain}
	\Lambda_*^{\pm}:=\{k =(k_1,\dots,k^d)\in \Lambda_* |  \pm k^1>0 \}.
\end{equation}

The  dispersion relation  takes the  discretized form 
\begin{equation}
	\label{NearestNeighbord}
	\omega_k \ = \ \omega(k) \ = \ \sin(2\pi   k^1)\Big[\sin^2(2\pi   k^1) + \cdots + \sin^2(2\pi   k^d)\Big],
\end{equation}
with $k=(k^1,\cdots,k^d)$. Moreover, we  simply set
\begin{equation}
	\label{omegao}
	\bar\omega(k) \ =    \sin(2\pi   k^1).
\end{equation}

We randomize the above dynamics by introducing a random noise
\begin{equation}
	\label{LatticeDynamicsFourierZ}\begin{aligned}
		\mathrm{d}\hat\psi(k,t) \ & = \ {\bf i} \omega(k)\hat\psi(k,t)\mathrm{d}t \ + \ {\bf i}
		\sum_{l \in \Lambda^*}
		\mathfrak{g}(k,l)\hat\psi(k,t)\circ\mathrm{d}W_l \\
		&\ \ \  \ +  {\bf i}\frac{1}{|\Lambda_*|^2}\lambda \bar\omega(k)\sum_{k=k_1+k_2;k_1,k_2\in\Lambda^*}\hat\psi(k_1,t)\hat\psi(k_2,t),\\
		\hat\psi(k,0) \ & = \ \hat\psi_0(k) , \end{aligned}
\end{equation}
where $W_k(t)$ is a sequence of independent real Wiener processes on a filter probability space $(\Omega,{\bf{F}},({\bf{F}})_{t\ge 0},{\bf{P}})$, the values of $\mathfrak{g}(k,l)$ will be specified later.

We also define the mesh size to be \begin{equation}\label{Mesh} h^d=\left(\frac{1}{2D+1}\right)^d.\end{equation} 

In this convention, the inverse Fourier transform is 
\begin{equation}\label{Def:FourierInverse}
	f(x)= h^d\sum_{k\in\Lambda_*} \hat  f(k) e^{2\pi {\bf i} k\cdot x}.
\end{equation}
We also use the following notations
\begin{equation}
	\label{Shorthand1}
	\int_{\Lambda}\mathrm{d}x \ = \  \sum_{x\in\Lambda},\ \ \ \ 	\int_{\Lambda_*}\mathrm{d}k \ = \  \frac{1}{|\Lambda_*|}\sum_{k\in\Lambda_*},\ \ \ \ 
	\langle f, g\rangle \ = \  \sum_{x\in\Lambda}f(x)^* g(x),
\end{equation}
where   if $z\in \mathbb{C}$, then $z^*$ is the complex conjugate, as well as the  Japanese bracket
$
\langle x\rangle \ = \ \sqrt{1+|x|^2}, \ \ \forall x\in\mathbb{R}^d.
$	Moreover, for any $N\in\mathbb{N}\backslash\{0\}$, we define the delta function $\delta_N$ on $(\mathbb{Z}/N)^d$ as
\begin{equation}
	\label{Def:Delta}\delta_N(k) = |N|^d\mathbf{1}(k \mbox{ mod } 1 \ = \ 0), \ \ \ \forall k\in (\mathbb{Z}/N)^d.
\end{equation}
In our computations, we omit the sub-index $N$ and simply write
\begin{equation}
	\label{Def:Delta1}\delta(k) = |N|^d\mathbf{1}(k \mbox{ mod } 1 \ = \ 0), \ \ \ \forall k\in (\mathbb{Z}/N)^d.
\end{equation} 
We  restrict the frequency domain to
\begin{equation}\begin{aligned}
		\Lambda^* = \Lambda^*(D) = &\  \left\{-\frac{D}{2D+1},\cdots,-\frac{1}{2D+1},\frac{1}{2D+1},\cdots,\frac{D}{2D+1}\right\}\\
		&\times \left\{-\frac{D}{2D+1},\cdots,0,\cdots,\frac{D}{2D+1}\right\}^{d-1}. \end{aligned}
\end{equation}
We also set
\begin{equation}
	\label{Shorthand1a}
	\int_{\Lambda^*}\mathrm{d}k \ = \  \frac{1}{|\Lambda^*|}\sum_{k\in\Lambda^*}.
\end{equation}

Here we also assume  that the solution $\psi$ is  real-valued, since $\omega_{-k}=-\omega_k$, and it follows that 
(see \,\cite{Nazarenko:2011:WT})
\begin{equation}
	\label{psikpsiminusk}\psi(-k,t)=\psi^*(k,t),\ \ \ \ \ \forall k\in\Lambda^*,\ \ \  \forall t\in\mathbb{R}_+.
\end{equation}

Next we will define $\mathfrak{g}$. Let $\varphi^{\Lambda_*^+}$ be a one-to-one mapping from $\Lambda_*^+$ to $\{1,\cdots,|\Lambda_*^+|\}$ and $\phi_{\Lambda_*^+}=(\varphi^{\Lambda_*^+})^{-1}$ be its inverse. We define any vector $v$ in $\mathbb{Z}^{|\Lambda_*^+|}$ to be $v=(v_{\phi_{\Lambda_*^+}(1)},\cdots,v_{\phi_{\Lambda_*^+}(|\Lambda_*^+|)})$.  

We define 
\begin{equation}
	\begin{aligned}\label{VectorSpace1}
 \mathcal{V}_\mathscr{E}^1 = \Big\{v=(v_{\phi_{\Lambda_*^+}(1)},\cdots,v_{\phi_{\Lambda_*^+}(|\Lambda_*^+|)})\in \mathbb{R}^{|\Lambda_*^+|}~~~ & \ \Big| ~~~ \exists (-k_1,k_2,k_3)\in \Lambda^*_+, -k_1+k_2+k_3=0, \mbox{ such that }\\
& -v_{k_1}=v_{k_2}=v_{k_3}=1 \mbox{ and } v_{k}=0 \mbox{ if } k\ne k_1,k_2,k_3 \Big\}.
\end{aligned}
\end{equation}
Let $\mathscr{V}_\mathscr{E}^1$ be the vector space spanned by $\mathcal{V}_\mathscr{E}^1$ and let  $\mathscr{V}_\mathscr{E}^2$ be the vector space that satisfies
\begin{equation}
	\label{VectorSpace2}
\mathcal{V}_\mathscr{E}^1	\bigoplus \mathcal{V}_\mathscr{E}^2 \ = \  \mathbb{R}^{|\Lambda_*^+|},
\end{equation}
which means the intersection between $\mathcal{V}_\mathscr{E}^1, \mathcal{V}_\mathscr{E}^2$ is only the zero-vector.
We denote by $\{e_1,\cdots, e_{\mathscr{M}_\mathscr{E}}\}\subset \mathbb{Z}^{|\Lambda_*^+|}$ an orthonormal basis of $\mathcal{V}_\mathscr{E}^2$.
We  set 
$	\mho_i    \ge \mathscr{C}_\mho =\lambda^{c_\mho}
\mathscr{C}^o_\mho h^{-d}>0$ with $\mathscr{C}^o_\mho$ being a  positive   constant and $c_\mho$ being a sufficiently small positive constant. The matrix $\mathscr{E}=\Big(\mathscr{E}(k,k')\Big)=\Big(\mathscr{E}(\varphi_{\Lambda_*^+}(k),\varphi_{\Lambda_*^+}(k'))\Big)$ can now be specified via its spectral representation
\begin{equation}
	\label{VectorSpace3}
	\mathscr{E} \ := \  	  \sum_{i=1}^{\mathscr{M}_\mathscr{E}} \mho_ie_ie_i^T,
\end{equation}
where $e_i^T$ denotes the transpose of $e_i$. Then, it follows straightforwardly that the matrix $\mathscr{E}=\Big(\mathscr{E}(k,k')\Big)=\Big(\mathscr{E}(\varphi_{\Lambda_*^+}(k),\varphi_{\Lambda_*^+}(k'))\Big)$ is positive semi-definite. 

 Therefore, we can define
$\tilde{\mathfrak{g}}: \Lambda_*^+\times \Lambda_*^+ \to \mathbb{R}$, such that 
$
	\sum_{l \in \Lambda_*^+} \tilde{\mathfrak{g}}(k,l)
	\tilde{\mathfrak{g}}(k',l) = \mathscr{E}(k,k'). 
$
Next, we extend $\mathscr{E}(k,k')$ to $\Lambda^*\times \Lambda^*$, by setting $\mathscr{E}(k,k') =$ $ \mathscr{E}((k^1,\cdots,k^d),$ $(k'^1,k'^2,\cdots,k'^d))$ $ =\mathscr{E}((|k^1|,k^2,\cdots,k^d),$ $(|k'^1|,k'^2,\cdots,k'^d))$. Therefore,  we can define
$\mathfrak{g}: \Lambda^* \times \Lambda^* \to \mathbb{R}$, such that 
\begin{equation}
	\sum_{l \in \Lambda^*} \mathfrak{g}(k,l)
	\mathfrak{g}(k',l) = \mathscr{E}(k,k'), \forall k,k'\in\Lambda^*.
\end{equation}

We now set $$a_k \ = \ \frac{\hat{\psi}(k)}{ \sqrt{|\bar\omega(k)|}},$$


Let ${a},{a}^*$ denote the vectors $(a_k)_{k\in\Lambda^*}$, $(a_k^*)_{k\in\Lambda^*}$, and let us  set 
\begin{equation}
	\label{Kernel}
	\mathcal{M}(k,k_1,k_2) \ = \ 2\mathrm{sign}(k^1) \mathcal{W}(k,k_1,k_2),
\end{equation}
with 
\begin{equation}\label{Kernel2}  \mathcal{W} \ = \sqrt{{|\bar\omega(k)\bar\omega(k_1)\bar\omega(k_2)|}}.\end{equation}

By defining 
\begin{equation}
	\label{ShortenNotationofA}
	{a}(k,1,t) \ = \ a^*_k(t), \mbox{ and } {a}(k,-1,t) \ = \ a_k(t),
\end{equation}
we rewrite the system  as
\begin{equation}
	\begin{aligned}\label{StartPoint0}
		\mathrm{d}{a}(k,\sigma,t)\ & =  \ -{\bf i}\sigma\omega(k){a}(k,\sigma,t) \mathrm{d}t\ - \ {\bf i}\sigma\lambda\int_{\Lambda^*}\mathrm{d}k_1\int_{\Lambda^*}\mathrm{d}k_2\delta( k-k_1-k_2)\times\\
		& \ \  \ \ \  \ \times \mathcal{M}(k,k_1,k_2){a}(k_1,\sigma,t) {a}(k_2,\sigma,t)\mathrm{d}t-\  \ {\bf i}  \sigma\sum_{l \in \Lambda^*}
		\mathfrak{g}(k,l){a}(k,\sigma,t)\circ\mathrm{d}W_l,
		\\
		{a}(k,-1,0) \ & = \ a_0(k), \ \ \ \ \forall (k,t)\in\Lambda^*\times \mathbb{R}_+, \sigma\in\{\pm1\},
	\end{aligned}
\end{equation}
In order to absorb the quantity $-{\bf i}\sigma\omega(k){a}(k,\sigma,t)$ on the right hand side of \eqref{StartPoint0}, we set
\begin{equation}
	\label{HatA}
	\alpha(k,\sigma,t) \ = \ {a}(k,\sigma,t)e^{{\bf i}\sigma \omega(k)t}.
\end{equation}
For  sake of simplicity, we also denote ${a}(k,\sigma,t)$ and $\alpha(k,\sigma,t)$ as ${a}_t(k,\sigma)$ and $\alpha_t(k,\sigma)$. The following system can be now derived for the new quantities $\alpha_t(k,\sigma)$
\begin{equation}
	\begin{aligned}\label{StartPoint}
	&	\mathrm{d}\alpha_t(k,\sigma) \  =  - \ {\bf i}\sigma\lambda \int_{\Lambda^*}\mathrm{d}k_1\int_{\Lambda^*}\mathrm{d}k_2\delta(k-k_1-k_2)\times\\
		& \  \times \mathcal{M}(k,k_1,k_2)\alpha_t(k_1,\sigma) \alpha_t(k_2,\sigma)e^{{\bf i}t\sigma(-\omega(k_1)-\omega(k_2)+\omega(k))}\mathrm{d}t\ -\ {\bf i} \sum_{l \in \Lambda^*}
		\mathfrak{g}(k,l)\alpha(k,\sigma,t)\circ\mathrm{d}W_l,
		\\
		\alpha_0(k,-1) \ & = \ a_0(k), \ \ \ \ \forall (k,t)\in\Lambda^*\times \mathbb{R}_+.
	\end{aligned}
\end{equation}
From \eqref{psikpsiminusk}, we also deduce that $\forall k\in\Lambda^*, \  \forall t\in\mathbb{R}_+, \ \forall \sigma\in\{\pm1\}$
\begin{equation}
	\label{akaminusk}a_{-k}(t)=a_k^*(t), \ \ \ \alpha(k,\sigma,t) = \alpha^*(k,-\sigma,t) = \alpha^*(-k,\sigma,t).
\end{equation}
We also  define the standard normed spaces
\begin{equation}\label{Def:Norm1}L^{p}(\mathbb{T}^d):= \left\{F(k): \mathbb{T}^d\to \mathbb{R}\ \Big| \ \|F\|_{L^p}=\|F\|_{p}=\left[\int_{\mathbb{T}^d}\mathrm{d}k |F(k)|^p\right]^{\frac{1}{p}}<\infty\right\}, \ \ \ \mbox{ for } p\in[1,\infty),\end{equation}

\begin{equation}\label{Def:Norm2}L^{\infty}(\mathbb{T}^d):= \left\{F(k): \mathbb{T}^d\to \mathbb{R}\ \Big| \ \|F\|_{L^\infty}= \|F\|_{\infty}=\sup_{k\in\mathbb{T}^d} |F(k)|<\infty\right\},\end{equation}
\begin{equation}\label{Def:Norm3} l^{p}(\mathbb{Z}^d):= \left\{F(x): \mathbb{Z}^d\to \mathbb{R}\ \Big| \ \|F\|_{l^p}= \|F\|_{p}=\left[\sum_{x\in\mathbb{Z}^d} |F(x)|^p\right]^{\frac{1}{p}}<\infty\right\}, \ \ \ \mbox{ for } p\in[1,\infty),\end{equation}
and
\begin{equation}\label{Def:Norm4}l^{\infty}(\mathbb{Z}^d):= \left\{F(x): \mathbb{Z}^d\to \mathbb{R}\ \Big| \ \|F\|_{l^\infty}=\|F\|_{\infty}=\sup_{x\in\mathbb{Z}^d} |F(x)|<\infty\right\}.\end{equation}

\subsection{Liouville equation, two-point  correlation function and weak turbulence theory}  
Let us introduce the real processes $B_{1,k}$ and $B_{2,k}$ such that $a_k=B_{1,k} + {\bf i}B_{2,k}$. Since $a_k$, $B_{1,k}$ and $B_{2,k}$ are random variables, we  use the variables $b_{1,k}$,  $b_{2,k}$ and $\tilde{a}_k=b_{1,k} +{\bf i}b_{2,k}$,  to present their roles in the density function. 
We suppose that the initial data are random variables   $(B_{1,k}(0,\varpi),B_{2,k}(0,\varpi))$ 
(see \,\cite{Nazarenko:2011:WT}) on the probability space $(\Omega,{\bf{P}})$, with  probability density functions $\varrho(0,b_1,b_2)$ in which $(b_1,b_2)$ $=$ $(b_{1,k},b_{2,k})_{k\in\Lambda^*}$.

It is then straightforward from the standard It\^o calculus that the total density function $\varrho(t,b_1,b_2)$  satisfies the Liouville equation (the Heisenberg picture)
\begin{equation}
	\label{Liouville1}
	\frac{\partial}{\partial t}\varrho \ =  \mathbf{R}\varrho \ - \{\mathcal{H},\varrho\}
\end{equation}
where  the form of 
$ \mathbf{R}[\varrho]$ will be specified later and
 the Poisson bracket is  defined as
$$\{A,B\} \ = \ \sum_{k\in\Lambda^*}\partial_{b_{1,k}}A\partial_{b_{2,k}}B- \partial_{b_{2,k}}A\partial_{b_{1,k}}B.$$

By the change of variables, 
\begin{equation}\label{ChangofVariable}
	b_{1,k}+{\bf i}b_{2,k}\ =\ {a}_k \ = \ \sqrt{2c_{1,k}}e^{{\bf i}c_{2,k}},
\end{equation} with $c_{1,k}\in \mathbb{R}_+$ and $c_{2,k}\in [-\pi,\pi]$, we have an equation of $\varrho(t,c_1,c_2)$ with the new variables $c_{1,k},c_{2,k}$
\begin{equation}
	\begin{aligned}\label{FokkerPlanck}
		\partial_t\varrho
		=\  &\mathbf{R}[\varrho]\
		 \ -\sum_{k\in\Lambda^*}\omega_k\partial_{c_{2,k}}\varrho \
		+\ \sum_{k\in\Lambda^*}\lambda \mathfrak{H}^a(k)\partial_{c_{1,k}}\varrho+ \sum_{k\in\Lambda^*}\lambda\mathfrak{H}^b(k)\partial_{c_{2,k}}\varrho,
	\end{aligned}
\end{equation}
in which
in which, we specify 
$ \mathbf{R}[\varrho]$ under the new variables as follows
\begin{equation} \label{NOISEgenerator}
	\mathbf{R}[\varrho]= \sum_{k,k' \in \Lambda^* } 
	\mathscr{E}(k,k')\frac{\partial^2}{\partial_{c_{2,k}} \partial_{c_{2,k'}}}\varrho, 
\end{equation}
and
$$\mathfrak{H}^a({k})= \int_{\Lambda^*}\mathrm{d}k_1\int_{\Lambda^*}\mathrm{d}k_2\mathcal{M}(k,k_1,k_2)\sqrt{2c_{1,k_1}c_{1,k_2}c_{1,k}}\Big[\delta(k-k_1-k_2){\sin(c_{2,k_1}+c_{2,k_2}-c_{2,k})}$$
$$ +\ 2\delta(k+k_1-k_2){\sin(-c_{2,k_1}+c_{2,k_2}-c_{2,k})}\Big],$$
and
$$\mathfrak{H}^b({k})=-\int_{\Lambda^*}\mathrm{d}k_1\int_{\Lambda^*}\mathrm{d}k_2\mathcal{M}(k,k_1,k_2)\sqrt{c_{1,k_1}c_{1,k_2}}\Big[\delta(k-k_1-k_2){\sin(c_{2,k_1}+c_{2,k_2})}$$
$$ +\ 2\delta(k+k_1-k_2){\sin(-c_{2,k_1}+c_{2,k_2})}\Big]\frac{\sin(c_{2,k})}{\sqrt{2c_{1,k}}}$$
$$-\int_{\Lambda^*}\mathrm{d}k_1\int_{\Lambda^*}\mathcal{M}(k,k_1,k_2)\delta(k-k_1-k_2)\sqrt{c_{1,k_1}c_{1,k_2}}\Big[\delta(k-k_1-k_2){\cos(c_{2,k_1}+c_{2,k_2})}$$
$$ +\ 2\delta(k+k_1-k_2){\cos(-c_{2,k_1}+c_{2,k_2})}\Big]\frac{\cos(c_{2,k})}{\sqrt{2c_{1,k}}}.$$
By defining the new Hamiltonian
\begin{equation}
	\label{DensityHamiltonian}\begin{aligned}
		\hat{\mathfrak{H}}({k})\ = &\  \lambda\int_{\Lambda^*}\mathrm{d}k_1\int_{\Lambda^*}\mathrm{d}k_2\mathcal{M}(k,k_1,k_2)\sqrt{2c_{1,k_1}c_{1,k_2}c_{1,k}}\Big[\delta(k-k_1-k_2){\cos(c_{2,k}-c_{2,k_1}-c_{2,k_2})}\\
		& \ + \ 2\delta(k+k_1-k_2){\cos(c_{2,k}+c_{2,k_1}-c_{2,k_2})}\Big]\ + \ \omega_k c_{1,k},
	\end{aligned}
\end{equation}
we have  the Liouville equation  
\begin{equation}
	\begin{aligned}\label{FokkerPlanck2}
		\partial_t\varrho\
-\	\mathbf{R}[\varrho]\	+\ &\sum_{k\in\Lambda^*}\Big[\Big[\hat{\mathfrak{H}}(k),\varrho\Big]\Big]_k=0,
	\end{aligned}
\end{equation}
in which
\begin{equation}\label{DensityHamiltonian1}\Big[\Big[\hat{\mathfrak{H}}(k),\varrho\Big]\Big]_k \ = \ \partial_{c_{2,k}}\hat{\mathfrak{H}}(k)\partial_{c_{1,k}}\varrho \ - \ \partial_{c_{1,k}}\hat{\mathfrak{H}}(k) \partial_{c_{2,k}}\varrho.\end{equation}

	

	\begin{definition}\label{def:distinct}

		For any given observable $F_+:\mathbb{R}^{2|\Lambda_*^+|}\to \mathbb{C}$,  and for any  random 
		variables $B_1=(B_{1,k})_{k\in\Lambda_*^+}$, $B_2=(B_{2,k})_{k\in\Lambda_*^+}$  on   $(\Omega,{\bf{F}},({\bf{F}})_{t\ge 0},{\bf{P}})$, we  define the average to be
		\begin{equation}\label{Averaget+}
			\langle F_+(B_1,B_2)\rangle^+ \ = \  \langle F_+(B_1,B_2)\rangle_{t}^+ = \Big\langle{F_+(B_1,B_2)}\Big\rangle_{\bf{P}}^+ \coloneqq \int_{\mathbb{R}^{2|\Lambda^
					+_*|}}
			\mathrm{d}b_1\mathrm{d}b_2 F_+(b_1,b_2)\varrho(t,b_1,b_2),
		\end{equation}	
		in which $\langle\rangle_{{\bf{P}}}^+$ is the expectation with respect to $
		{\bf{P}}$.  The initial density measure at time $t=0$ of the system to be the density function of standard independent Gaussians  
		\begin{equation}
			\label{Propo:ExampleMeasure:1}
			{\varrho}(0,b_1,b_2) \ = \ \frac{1}{Z_0}\exp\left(-\sum_{k\in\Lambda^+_*}\frac{b_{1,k}^2+b_{2,k}^2}{\gimel(k)}\right) \ = \ \frac{1}{Z_0}\exp\left(-\sum_{k\in\Lambda^+_*}\frac{2c_{1,k}}{\gimel(k)}\right)
		\end{equation}
		in which 
		$$Z_0=\prod_{k\in\Lambda^+_*}\pi\gimel(k),$$
		 for some function $\gimel(k)>0$.
		Since $a_k=a_{-k}^*$,  the above definition can be extended to 
		$\Lambda^*$ as follows:  for $B_1=(B_{1,k})_{k \in \Lambda^*}$, 
		$B_2 = (B_{2,k})_{k \in \Lambda^*} $, we set
		\begin{equation*}
			B_1^{\pm}\coloneqq(B_{1,k})_{k \in 
				\Lambda^{\pm}_*} \quad \text{and} \quad B_2^{\pm} \coloneqq (B_{2,k})_{k \in 
				\Lambda^{\pm}_*}.
		\end{equation*}
		Then for any given observable $F: \mathbb{R}^{2 \Lambda^*} \to \mathbb{C}$,
		as $$F(B_1,B_2)=F(B_1^+,B_2^+,B_1^-,B_2^-),$$ and 
		we define
		\begin{equation}
			\label{Averaget}\begin{aligned}
				\langle F(B_1,B_2)\rangle_{t} \coloneqq
				\int_{\mathbb{R}^{2|\Lambda^
						+_*|}}
				\mathrm{d}b_1^+\mathrm{d}b_2^+ F(b_1^+,b_2^+,b_1^+,-b_2^+)\varrho(t,b_1,b_2).\end{aligned}
		\end{equation}
		In this case, we can consider the initial measure as
		\begin{equation}
			\label{Propo:ExampleMeasure:2}
			{\varrho}(0,b_1,b_2) \ = \ {\varrho}(0,c_1) \ = \  \frac{1}{Z_0}\exp\left(-\sum_{k\in\Lambda^*}\frac{b_{1,k}^2+b_{2,k}^2}{\gimel(k)}\right) \ = \ \frac{1}{Z_0}\exp\left(-\sum_{k\in\Lambda^*}\frac{2c_{1,k}}{\gimel(k)}\right)
		\end{equation}
			in which 
		$Z_0=\prod_{k\in\Lambda^*}\pi\gimel(k), $and under the additional condition that  ${\varrho}(0,c_1)$ is symmetric  in the sense ${\varrho}(0,\cdots,c_{1,(-k^1,k^2,\cdots,k^d)},$ $\cdots,c_{1,(k^1,k^2,\cdots,k^d)},\cdots)$ $={\varrho}(0,\cdots,c_{1,(k^1,k^2,\cdots,k^d)},$ $\cdots,c_{1,(-k^1,k^2,\cdots,k^d)},\cdots)$. 
			Recalling that $h$ is the mesh size, we assume  $\gimel(k)=h^{-d}\bar\gimel(k)$ with $0< \bar\gimel(k)<C_{\gimel}/|k^1|$ for some constant $C_{\gimel}>0$ and $k=(k^1,\cdots,k^d)$. We put   $\mathcal{A}_{x}^o=\int_{\Lambda^*}\mathrm{d}k{a_{k}}
		e^{{\bf i}2\pi k\cdot x}$. We then have
		\begin{equation}
			\begin{aligned}
				\label{Propo:ExampleMeasure:2}
				{\varrho}(0,c_1,c_2) \ = \ & \frac{1}{Z_0}\exp\left(-\sum_{k\in\Lambda^*}h^d\frac{2{c}_{1,k}}{\bar\gimel(k)}\right)
				\ = \ & \frac{1}{Z_0}\exp\left(-\sum_{x\in\Lambda }|\mathcal{A}_{x}^o*\Re_x|^2\right),
			\end{aligned}
		\end{equation}
		
		in which $\Re_x=\int_{\Lambda^*}\mathrm{d}k\sqrt{1/\bar\gimel(k)}
		e^{{\bf i}2\pi k\cdot x}.$
		This means even $\varrho(t)$ is a probability density over $\mathbb{R}^{|\Lambda_*^+|}$, thanks to the extension based on    symmetrization, it can be considered as a probability density over $\mathbb{R}^{|\Lambda_*|}$.

				The set $\{k_1,\cdots,k_n\}\subset \Lambda^*$ is called admissible for any $k\in \{k_1,\cdots,k_n\}$, if either $-k$ is also in $\{k_1,\cdots,k_n\}$ or there are $k',k''\in \{k_1,\cdots,k_n\}$ such that $k+k'+k''=0$. Otherwise, the set is not admissible. The set $\{(k_1,\sigma_1),\cdots,(k_n,\sigma_n)\}\subset \Lambda^*\times\{\pm1\}$ is called admissible if $\{k_1\sigma_1,\cdots k_n\sigma_n\}\subset \Lambda^*$ is admissible.
\end{definition}



{  
	For multi-indices $\mathcal{V}_1, \mathcal{V}_2 \in \mathbb{N}^{|\Lambda^*|}$, we define the multinomial expression
	\begin{equation} \label{defaV}
		({a})^{\mathcal{V}_1} \coloneqq \prod_{k \in \Lambda_*} ({a}_k)^{\mathcal{V}_{1,
				\varphi_{\Lambda^*}(k)}}, \quad ({a}^*)^{\mathcal{V}_2} \coloneqq \prod_{k \in \Lambda_*} ({a}^*_k)^{\mathcal{V}_{2,
				\varphi_{\Lambda^*}(k)}},
	\end{equation}
	in which $\mathcal{V} \in \mathbb{N}^{|\Lambda^*|}$ and $k \in \Lambda^*$, for simplicity, we denote $\mathcal{V}_k \coloneqq \mathcal{V}_{	\varphi_{\Lambda^*}{(k)}}$. Then, we obtain
	\begin{multline}\label{Noiseaction0}
		\mathbf{R}\left(({a})^{\mathcal{V}_1} ({a}^*)^{\mathcal{V}_2}\right) = - 2\sum_{k_1,k_2 \in \Lambda_*^+}\Big( \mathscr{E}(k_1,k_2)
		(\mathcal{V}_{1,k_1}-\mathcal{V}_{1,-k_1}+\mathcal{V}_{2,-k_1}-\mathcal{V}_{2,k_1}) \\
		\times 	(\mathcal{V}_{1,k_2}-\mathcal{V}_{1,-k_2}+\mathcal{V}_{2,-k_2}-\mathcal{V}_{2,k_2})
		({a})^{\mathcal{V}_1} ({a}^*)^{\mathcal{V}_2}\Big).
	\end{multline}
}


	Let us now  define 
	\begin{equation}
		\label{omegaoprime}
		\bar\omega'(k)=\frac{h^{-d}}{\tilde\omega(k)},\ \ \ \ \tilde\omega(k)=|k^1|, \ \ \ \mbox{ with } k=(k^1,\cdots,k^d),
	\end{equation}
	and the invariant measure
	\begin{equation}
		\label{Invariant}
		\tilde{\varrho}(b_1,b_2) \ = \ \frac{1}{Z}\exp\left(-\sum_{k\in\Lambda^*}\frac{b_{1,k}^2+b_{2,k}^2}{\bar\omega'(k)}\right),
	\end{equation}
	in which 
	$Z=\prod_{k\in\Lambda^*}\pi\bar\omega'(k).$  
	
	\begin{proposition}\label{Propo:ExpectationEqui}
		The following identity holds true 
		\begin{equation}
			\begin{aligned}
				& \mathbf{R}\tilde{\varrho}\ = \  \left\{\mathcal{H},\exp\left(-\sum_{k\in\Lambda^*}\frac{b_{1,k}^2+b_{2,k}^2}{\bar\omega'(k)}\right)\right\} \ = \  0.
			\end{aligned}
		\end{equation}

		Consider the   backward ZK equation in the Fourier space (see \eqref{LatticeDynamicsFourier})	\begin{equation}
			\label{EquationTrajectoryFourier}\begin{aligned}
				\frac{\mathrm{d}\bar{\wp}(k,s)}{\mathrm{d} s} \ & = -\ {\bf i} \omega(k)\bar{\wp}(k,s)  \ - \     {\bf i}\lambda \bar\omega(k)\frac{1}{|\Lambda_*|^2}\sum_{k=k_1+k_2;k_1,k_2\in\Lambda_*}\bar{\wp}(k_1,s)\bar{\wp}(k_2,s), \end{aligned}
		\end{equation}
		the density distribution can be computed as follows 
		\begin{equation}
			\begin{aligned}
				\label{Propo:ExpectationEqui:3}
				&\varrho(0,\mathcal{B}_{1}(0,t,b_1,b_2),\mathcal{B}_{2}(0,t,b_1,b_2)) \ =  \ \varrho(t,b_1,b_2),
			\end{aligned}
		\end{equation}
		in which $	\wp_k(s)=[\mathcal{B}_{1,k}(t-s,t) + \mathbf{i}\mathcal{B}_{1,k}(t-s,t)],$ and $	\bar{\wp}_k(s)=[\mathcal{B}_{1,k}(t-s,t) + \mathbf{i}\mathcal{B}_{1,k}(t-s,t)]\sqrt{|\bar\omega(k)|}.$ 
		
		Define 	\begin{equation}\label{Gaussian}
		\mathfrak{P} \ = \ \exp\left( \sum_{k\in\Lambda^*}c_{\mathfrak{P}}c_{1,k}\tilde\omega(k)h^{d}\right),	\end{equation}
	for any $n\in\mathbb{N},$ and $c_{\mathfrak{P}}$ is a positive constant.
	We have the  ``$L^1$-density identity''
	\begin{equation}\label{Propo:ExpectationEqui:1}
		\begin{aligned}
			\partial_t\int_{(\mathbb{R}_+\times[-\pi,\pi])^{|\Lambda^*|}}\mathrm{d}c_{1}\mathrm{d}c_{2}\, \varrho\, \mathfrak{P}  \
			= \ & 0.
		\end{aligned}
	\end{equation}

	Let $n$ be an arbitrary positive natural number and $\{k_{1},\cdots,k_{n}\}$ be a subset of $\Lambda^*$. 	We then have
	\begin{equation}\label{Propo:ExpectationEqui:2}
		\begin{aligned}
	\int_{\Lambda_*^n}{\mathrm{d}k_1\cdots\mathrm{d}k_n}	[\tilde\omega(k_1)\cdots\tilde\omega(k_n)]	\left \langle |a_{k_{1}}|^2\cdots|a_{k_{n}}|^2 \right\rangle_t \
			\lesssim  &\ \	h^{-dn}\mathcal{C}_o^n,
		\end{aligned}
	\end{equation}
	where the constant in the inequality \eqref{Propo:ExpectationEqui:2} is universal.
	\end{proposition}
	We prove this proposition in Subsection \ref{proofs} below.
	%
	

	

\medskip

{\bf Prediction from Weak Turbulence Theory.}
The standard expectation from 
physicists working on the derivation of the kinetic wave equation goes as follows. Consider the two-point  correlation function 
\begin{equation}
	\label{Expectation}
	f(k,t)\delta_{k=k'} \ = \ \langle {\alpha}_{t}(k,-1){\alpha}_{t}(k',1)\rangle.
\end{equation}
In the limit of  $D\to \infty$, $\lambda\to 0$ and $t=\lambda^{-2}\tau=\mathcal{O}(\lambda^{-2})$, the two-point correlation function 
$f(k,t)$ has the limit 
$$\lim_{\lambda\to 0, D\to\infty} f(k,\lambda^{-2}\tau) = f^\infty(k,\tau)$$
which solves the 3-wave  equation  (see \,\cite{Nazarenko:2011:WT})
\begin{equation}\label{Eq:WT0}
	\frac{\partial}{\partial \tau}f^\infty(k,\tau)=  \mathcal{C}\big(f^\infty  \big)(k,\tau)
\end{equation}
with the collision operator 
\begin{equation}
	\begin{aligned}\label{Eq:CollisionWT1}
		&\mathcal{C}(f^\infty)(k_1)=   \int_{(\mathbb{T}^d)^{2}}
		\mathrm{d}k_2\mathrm{d}k_3|\mathcal{M}(k_1,k_2,k_3)|^2  \frac{1}{\pi}\delta(\omega(k_3)+\omega(k_2)-\omega(k_1))\\
		&\times\delta(k_2+k_3-k_1) \Big( f_2^\infty f^\infty_3-f^\infty_1f^\infty_2\mathrm{sign}(k_1^1)\mathrm{sign}(k_3^1)-f^\infty_1f^\infty_3\mathrm{sign}(k_1^1)\mathrm{sign}(k_2^1)\Big),
\end{aligned}\end{equation}
in which $\mathbb{T}$ is the periodic torus $[-1/2,1/2]$. 
Here we have introduced the shorthand notation $f^\infty_j=f^\infty(k_j), j=1,2,3$. We also set $\mathbb{T}^d_+=\{k=(k^1,\cdots,k^d)\in\mathbb{T}^d\ | \  k^1\ge 0 \}$.

%
For the ZK equation, as discussed in (iii) of the introduction, the delta function $\delta(\omega(k_3)+\omega(k_2)-\omega(k_1))$ appearing in the collision operator \eqref{Eq:CollisionWT1} is not well-defined as a positive measure, and the resonance broadening is needed. We say that a function $f_\ell^\infty$ solves
the ``resonance broadening'' 3-wave  equation  if and only if
\begin{equation}\label{Eq:WT1}
	\frac{\partial}{\partial \tau}f^\infty_\ell(k,\tau)=  \mathcal{C}_\ell\big(f^\infty_\ell  \big)(k,\tau),
\end{equation}
with the collision operator 
\begin{equation}
	\begin{aligned}\label{Eq:CollisionWT5}
		&\mathcal{C}_\ell (f^\infty)(k_1)=   \int_{(\mathbb{T}^d)^{2}}
		\mathrm{d}k_2\mathrm{d}k_3|\mathcal{M}(k_1,k_2,k_3)|^2  \frac{1}{\pi}\delta_\ell (\omega(k_3)+\omega(k_2)-\omega(k_1))\\
		&\times\delta(k_2+k_3-k_1) \Big( f_2^\infty f^\infty_3-f^\infty_1f^\infty_2\mathrm{sign}(k_1^1)\mathrm{sign}(k_3^1)-f^\infty_1f^\infty_3\mathrm{sign}(k_1^1)\mathrm{sign}(k_2^1)\Big),
\end{aligned}\end{equation}
in which  $\delta_\ell$ is defined in \eqref{Eq:CollisionWT6}.
Note that when $\ell=0$, the resonance broadening $\delta_\ell (\omega(k_3)+\omega(k_2)-\omega(k_1))$ returns to $\delta (\omega(k_3)+\omega(k_2)-\omega(k_1))$ (see \cite{chekhlov96,huang2000,lee2007formation,l1997statistical,lvov2004hamiltonian,lvov2010oceanic,lvov2012resonant,remmel2014nonlinear,smith2001,smith2005near,smith1999transfer,smith2002generation}).

 The function $f_\ell^\infty(k,\tau)$ solves \eqref{Eq:CollisionWT5} with the initial condition  $f_\ell^\infty(k,0)=\lim_{D\to\infty} f(k,0)$  and  can be written in term of Taylor expansions as
\begin{equation}
 \label{Eq:CollisionWT5:a} f^\infty_\ell(k,\tau)  \ = \  \sum_{q=0}^\infty \frac{\tau^q}{q!}\mathcal{C}_\ell^q[\lim_{D\to\infty}f(k,0)],\end{equation}	where the operators $\mathcal{C}_\ell^q$ are explicit. 

{\bf Analytic expression of the leading diagrams}
Next, we present the expression of the leading quantities in $\lim_{D\to\infty}f(\tau)$. We split $\lim_{D\to\infty}f(\tau)$ as the sum of a leading part and non-leading part
\begin{equation}
	\label{LeadingDiagram:2} \lim_{D\to\infty}f(\tau)  : = \ f_{nonleading}(\tau)  \ + \ f_{leading}(\tau),\end{equation}  
and
\begin{equation}
	\label{LeadingDiagram:1}
	f_{leading}(\tau) \ :=\  	\sum_{q=0}^{\mathfrak{N}/4} \mathfrak{C}^\lambda_{q,\infty,0},
\end{equation}  
where $\mathfrak{N}$ is defined in \eqref{Def:Para0} below and  $\mathfrak{C}^\lambda_{q,\infty,\ell}$ can be defined for any $\ell\ge 0$ as

\begin{equation}
	\begin{aligned}\label{FinalProof:E10a:1}
		\mathfrak{C}^\lambda_{q,\infty,\ell}(\tau)\ = \ &(-1)^q \lambda^{-2q}
		\int_{(\mathbb{R}_+)^{\{1,\cdots,q\}}}\!\mathrm{d} s_{1} \cdots \mathrm{d} s_{2 q-1} 
		[\mathcal{C}^\lambda_{2,\infty,\ell}(\lambda^{-2 }s_1)\cdots\mathcal{C}^\lambda_{q+1,\infty,\ell}(\lambda^{-2 }s_{2q-1})]\\
		&\times
		\mathbf{1}\left(\sum_{i=1}^{q} s_{2i-1} \le \tau \right)
		\frac{1}{q!} \left(\tau-\sum_{i=1}^{q} s_{2i-1}\right)^q\,,
	\end{aligned}
\end{equation}
where
\begin{equation}
	\begin{aligned}\label{FinalProof:E10:a:1}
		&\mathcal{C}^\lambda_{m+1,\infty,\ell}(s,k_1',\cdots,k'_i,\cdots,k'_{m})=   \sum_{i=1}^{m}\mathcal{C}^\lambda_{i,m+1,\infty,\ell}(s,k_{i}')
\end{aligned}\end{equation}
with
\begin{equation}
	\begin{aligned}\label{FinalProof:E9:1}
		&\mathcal{C}^\lambda_{i,m+1,\infty,\ell}(s,k_1',\cdots,k'_i,\cdots,k'_{m})=   \iint_{(\mathbb{T}^d)^2}
		\mathrm{d}k'\mathrm{d}k_{m+1}'|\mathcal{M}(k_{i}',k',k_{m+1}')|^2 \\
		&\times \Big[e^{{\bf i}s(\omega(k')+\omega(k_{m+1}')-\omega(k_{i}'))-s\ell} +e^{-{\bf i}s(\omega(k')+\omega(k_{m+1}')-\omega(k_{i}'))-s\ell}\Big]\\
		&\times\delta(k'+k_{m+1}'-k_{i}')\Big( \mathfrak{L}_{m+1}^\infty(k_1',\cdots,k',\cdots,k'_{m+1})-\mathfrak{L}_{m+1}^\infty(k_1',\cdots,k'_i,\cdots,k'_{m+1})\mathrm{sign}(k_{i}')\mathrm{sign}(k_{m+1}')\\
		& -\mathfrak{L}_{m+1}^\infty(k_1',\cdots,k'_i,\cdots,k')\mathrm{sign}(k_{i}')\mathrm{sign}(k')\Big),
\end{aligned}\end{equation}
in which $k'$ takes the position of $k'_i$ in $\mathfrak{L}_{m+1}(k_1',\cdots,k',\cdots,k'_{m+1})$   and of $k'_{m+1}$ in $\mathfrak{L}_{m+1}(k_1',\cdots,$ $ k'_i,\cdots,k')$. We also have
\begin{equation}
	\mathfrak{L}_{m+1}^\infty(k_1',\cdots,k'_i,\cdots,k'_{m}) \ = \ \mathcal{C}^\lambda_{m+2,\infty,\ell}(s,k_1',\cdots,k'_i,\cdots,k'_{m+1}),
\end{equation}
and the $q+1$-correlation function of the first time slice is now defined to be
\begin{equation}
	\mathfrak{L}_{q+1}^\infty(k_{1}',\cdots,k_{q+1}'):=\prod_{i=1}^{q+1} \lim_{D\to\infty}{f}(k_{i}',0),\end{equation}

{\it In formal derivations, and in the rigorous derivation of wave kinetic equations coming from the nonlinear Schr\"odinger equation \cite{lukkarinen2009not,LukkarinenSpohn:2011:WNS} $	f_{leading}(\tau)$ is referred to as the leading diagrams, and when $\lambda\to 0$,  will give the Taylor series of the solution of \eqref{Eq:WT0}. We refer to \cite{lukkarinen2009not,LukkarinenSpohn:2011:WNS}  for the similar expressions of the leading diagrams in the nonlinear Schr\"odinger case. As  discussed in the introduction, the delta function $\delta(\omega(k_3)+\omega(k_2)-\omega(k_1))$ is not well-defined as a positive measure.   The limit 
\begin{equation}\lim_{\lambda\to 0}f_{leading}(\tau),\end{equation}
is therefore not well  defined.}  We modify $	f_{leading}(\tau)$, just to make sure that the limit $\lambda\to0$ can be taken by adding the resonance broadening $\ell>0$
\begin{equation}
	\label{LeadingDiagram:3}
	f_{\ell,leading}(\tau) \ :=\  	\sum_{q=0}^{\mathfrak{N}/4} \mathfrak{C}^\lambda_{q,\infty,\ell},
\end{equation}  
 and then add the modified leading part $f_{\ell,leading}$ with the standard non-leading part $f_{nonleading}$
\begin{equation}
	\label{LeadingDiagram:4} f_{\ell}(\tau)  : = \ 	f_{\ell,leading}(\tau) + 	f_{nonleading}(\tau).\end{equation}  
Our main result then shows that the non-leading term $ f_{nonleading}(\tau)$ vanishes in the kinetic limit and the  broadened leading part $f_{\ell,leading}$ converges to the Taylor series \eqref{Eq:CollisionWT5:a}, which is the solution of \eqref{Eq:WT1}.  {\it To the best of our knowledge, the theorem provides the first rigorous justification of the resonance broadening, discussed previously in the  physics literature  (see \cite{chekhlov96,huang2000,lee2007formation,l1997statistical,lvov2004hamiltonian,lvov2010oceanic,lvov2012resonant,remmel2014nonlinear,smith2001,smith2005near,smith1999transfer,smith2002generation}).}

We nam now state the main theorem.
\begin{theorem}\label{TheoremMain}
	Suppose that $d\ge 2$. Let us write  $t = \tau\lambda^{-2}$,  where  $\tau>0$.  Under the settings of Definition \ref{def:distinct},   and the noise constructed in Section \ref{Sec:Setting}, there exist $1>T_*>0$, such that for  $0<\tau<T_*$ and for any $\ell>0$
	
	\begin{itemize}
		
			\item[(i)] The non-leading part of $f(\tau$) vanishes in the kinetic limit
		\begin{equation}\label{TheoremMain:1}	\begin{aligned} &\lim_{\lambda\to 0}  \mathscr{W}\left(\Big| f(\tau)_{nonleading} \Big|>\theta \right)=  0,
		\end{aligned}	\end{equation}
		for any $\theta>0$. 		Here $ \mathscr{W}$ represents the Lebesgue measure.

		\item[(ii)] $f_{\ell}(\tau)$ converges to $f_{\ell}^\infty(\tau)$ in the kinetic limit, which is the solution of \eqref{Eq:WT1}, defined in \eqref{Eq:CollisionWT5:a} 
		\begin{equation}\label{TheoremMain:2}	\begin{aligned} &\lim_{\lambda\to 0}  \mathscr{W}\left(\Big| f_\ell(\tau)-f^\infty_\ell(\tau) \Big|>\theta \right)=  0,
		\end{aligned}	\end{equation}
		for any $\theta>0$.

	\end{itemize}  
	
\end{theorem}

\section{Duhamel expansions}

\subsection{Duhamel expansions of the Fourier modes of the solution} \label{Sec:Duha}

For two vectors $\mathcal{V}_1=(\mathcal{V}_{1,k})_{k\in\Lambda^*}\in\mathbb{N}^{|\Lambda^*|}, $ and $\mathcal{V}_2=(\mathcal{V}_{2,k})_{k\in\Lambda^*}\in\mathbb{N}^{|\Lambda^*|}$,  we denote
\begin{equation}\label{AverageVectorShorthand}
	  \quad \left\langle a^{\mathcal{V}_{1}} (a^*)^{\mathcal{V}_{2}}\right\rangle_t \ = \ \left\langle \prod_{k\in \Lambda^* }a_k^{\mathcal{V}_{1,k}} (a_k^*)^{\mathcal{V}_{2,k}}\right\rangle_t.
\end{equation}
First, we will study the evolution of the momenta $a^{\mathcal{V}_1}(a^*)^{\mathcal{V}_2}$ in which we adopt the notations of \eqref{AverageVectorShorthand}.

We then obtain an  equation for $\langle {a}^{\mathcal{V}_1}({a}^*)^{\mathcal{V}_2}\rangle$
\begin{equation}
	\label{KolmogorovEquation2}
	\begin{aligned}
		\frac{\partial}{\partial t}\left\langle {a}^{\mathcal{V}_1}({a}^*)^{\mathcal{V}_2}\right\rangle \ = \ & \left({\bf i}\sum_{k\in\Lambda^*} \omega_k(\mathcal{V}_{1}-\mathcal{V}_{2})- \mathscr{T}_{(\mathcal{V}_{1},-\mathcal{V}_{2})}\right)\left\langle {a}^{\mathcal{V}_1}({a}^*)^{\mathcal{V}_2}\right\rangle\\
		& + \lambda\,  \sum_{\mathcal{W}_1,\mathcal{W}_2}\mathcal{U}_{\mathcal{W}_1,\mathcal{W}_2}\left\langle {a}^{\mathcal{W}_1}({a}^*)^{\mathcal{W}_2}\right\rangle,
	\end{aligned}
\end{equation}
where $(\mathcal{V}_{1},-\mathcal{V}_{2})$ is a vector in $\mathbb{N}^{|\Lambda^*|}$ whose components are composed by those of $\mathcal{V}_{1},-\mathcal{V}_{2}$; $\mathscr{T}_{(\mathcal{V}_{1},-\mathcal{V}_{2})}$ is a non-negative constant depending on $\mathcal{V}_{1},-\mathcal{V}_{2}$; $\mathcal{U}_{\mathcal{W}_1,\mathcal{W}_2}$ is a coefficient depending on $\mathcal{M}$; and $\mathcal{W}_1,\mathcal{W}_2$ are new multiindices depending on $\mathcal{V}_1$ and $\mathcal{V}_2$. 
{We also have 
\begin{equation}
	\label{KolmogorovEquation2:explicit}
	\begin{aligned}
		\lambda\sum_{\mathcal{W}_1,\mathcal{W}_2}\mathcal{U}_{\mathcal{W}_1,\mathcal{W}_2}\left\langle {a}^{\mathcal{W}_1}({a}^*)^{\mathcal{W}_2}\right\rangle=&-\lambda\sum_{k\in\Lambda^*}\left\langle\mathcal{H}^1_2(k)\partial_{b_{1,k}}\Big( {a}^{\mathcal{V}_1}( {a}^*)^{\mathcal{V}_2}\Big)\right\rangle\ \\
		&+ \ \lambda\sum_{k\in\Lambda^*} \left\langle\mathcal{H}_2^2(k)\partial_{b_{2,k}}\Big( {a}^{\mathcal{V}_1}( {a}^*)^{\mathcal{V}_2}\Big)\right\rangle\\
		=\ &\sum_{k\in\Lambda^*}\left\langle-\lambda\mathcal{H}^1_2(k)\partial_{b_{1,k}}\Big( {a}^{\mathcal{V}_1}( {a}^*)^{\mathcal{V}_2}\Big)\ + \ \lambda\mathcal{H}_2^2(k)\partial_{b_{2,k}}\Big( {a}^{\mathcal{V}_1}( {a}^*)^{\mathcal{V}_2}\Big)\right\rangle,
	\end{aligned}
\end{equation}}
where
\begin{equation}\begin{aligned}
		\mathcal{H}^1_2(k)
		=\ &\frac12\left[\iint_{(\Lambda^*)^2}\mathrm{d}k_1\mathrm{d}k_2\mathcal{M}(k,k_1,k_2)\delta(k-k_1-k_2)(b_{1,k_1}b_{2,k_2}+b_{1,k_2}b_{2,k_1})\right.\\
		\ &\left.+2\iint_{(\Lambda^*)^2}\mathrm{d}k_1\mathrm{d}k_2\mathcal{M}(k,k_1,k_2)\delta(k+k_1-k_2)(b_{1,k_1}b_{2,k_2}-b_{1,k_2}b_{2,k_1})\right]\\
		\ =\ &\iint_{(\Lambda^*)^2}\mathrm{d}k_1\mathrm{d}k_2\mathcal{M}(k,k_1,k_2)\delta(k-k_1-k_2)(b_{1,k_1}b_{2,k_2}+b_{1,k_2}b_{2,k_1})\\
		\ =\ &\iint_{(\Lambda^*)^2}\mathrm{d}k_1\mathrm{d}k_2\mathcal{M}(k,k_1,k_2)\delta(k-k_1-k_2)\mathrm{Im}[{a}_{k_1}{a}_{k_2}],\end{aligned}
\end{equation}

and
\begin{equation}\begin{aligned} \mathcal{H}^2_2(k)
		=\ &\frac12\left[\iint_{(\Lambda^*)^2}\mathrm{d}k_1\mathrm{d}k_2\mathcal{M}(k,k_1,k_2)\delta(k-k_1-k_2)(b_{1,k_1}b_{1,k_2}-b_{2,k_1}b_{2,k_2})\right.\\
		& \left.+2\iint_{(\Lambda^*)^2}\mathrm{d}k_1\mathrm{d}k_2\mathcal{M}(k,k_1,k_2)\delta(k+k_1-k_2)(b_{1,k_1}b_{1,k_2}+b_{2,k_1}b_{2,k_2})\right]\\
		\ =\ &\iint_{(\Lambda^*)^2}\mathrm{d}k_1\mathrm{d}k_2\mathcal{M}(k,k_1,k_2)\delta(k-k_1-k_2)(b_{1,k_1}b_{1,k_2}-b_{2,k_1}b_{2,k_2})\\
		\ =\ &\iint_{(\Lambda^*)^2}\mathrm{d}k_1\mathrm{d}k_2\mathcal{M}(k,k_1,k_2)\delta(k-k_1-k_2)\mathrm{Re}[ {a}_{k_1} {a}_{k_2}].\end{aligned}
\end{equation}

In the above expression \eqref{KolmogorovEquation2:explicit}, we observe that $ \partial_{b_{1,k}}\Big( {a}^{\mathcal{V}_1}( {a}^*)^{\mathcal{V}_2}\Big)$ and $\partial_{b_{2,k}}\Big( {a}^{\mathcal{V}_1}( {a}^*)^{\mathcal{V}_2}\Big)$ will be zero if the components $\mathcal{V}_{1,k}$ and  $\mathcal{V}_{2,k}$ are $0$ in the vectors $\mathcal{V}_{1},\mathcal{V}_{2}$. If either $\mathcal{V}_{1,k}$ or  $\mathcal{V}_{2,k}$ is different from $0$, the quantities $\mathfrak{H}^1_1(k)$ and $\mathfrak{H}^1_2(k)$ {contain} the delta function $\delta(k-k_1-k_2)$, which will later defined to be the decomposition of the momentum $k$ into two momenta $k_1$ and $k_2$ in our graph representation. In this procedure, the component ${a}_k$ in $ {a}^{\mathcal{V}_1}( {a}^*)^{\mathcal{V}_2}$ is then replaced  by
\begin{equation}
	\label{KolmogorovEquation2:explicit:a}
	\begin{aligned}
		 {a}_k\ \longrightarrow\ & -\mathcal{H}^1_2(k)\partial_{b_{1,k}} {a}_k \ + \  \mathcal{H}_2^2(k)\partial_{b_{2,k}} {a}_k\\
		\ = \ & -\lambda\iint_{(\Lambda^*)^2}\mathrm{d}k_1\mathrm{d}k_2\mathcal{M}(k,k_1,k_2)\delta(k-k_1-k_2)\mathrm{Im}[ {a}_{k_1} {a}_{k_2}]\\
		&\ \ \ \ + \lambda\mathbf{i}\iint_{(\Lambda^*)^2}\mathrm{d}k_1\mathrm{d}k_2\mathcal{M}(k,k_1,k_2)\delta(k-k_1-k_2)\mathrm{Re}[ {a}_{k_1} {a}_{k_2}]
		\\
		\ = \ & \lambda\mathbf{i}\iint_{(\Lambda^*)^2}\mathrm{d}k_1\mathrm{d}k_2\mathcal{M}(k,k_1,k_2)\delta(k-k_1-k_2) {a}_{k_1} {a}_{k_2},
	\end{aligned}
\end{equation}
and  the component $ {a}_k^*$ in {${a}^{\mathcal{V}_1}( {a}^*)^{\mathcal{V}_2}$} is then replaced by
\begin{equation}
	\label{KolmogorovEquation2:explicit:b}
	\begin{aligned}
		 {a}_k^*& \longrightarrow   -\lambda\mathbf{i}\iint_{(\Lambda^*)^2}\mathrm{d}k_1\mathrm{d}k_2\mathcal{M}(k,k_1,k_2)\delta(k-k_1-k_2) {a}^*_{k_1} {a}^*_{k_2},
	\end{aligned}
\end{equation}
to form  the last term $\lambda\sum_{\mathcal{W}_1,\mathcal{W}_2}\mathcal{U}_{\mathcal{W}_1,\mathcal{W}_2}\left\langle {a}^{\mathcal{W}_1}({a}^*)^{\mathcal{W}_2}\right\rangle$, on the right hand side of \eqref{KolmogorovEquation2}.

Let us now define the new average 
\begin{equation}
	\label{Def:NewAverA}
	\langle\langle a^{\mathcal{V}_1}(a^*)^{\mathcal{V}_2}\rangle\rangle \ = \ e^{-{\bf i}t\sum_{k\in\Lambda^*} \omega_k(\mathcal{V}_{1,k}-\mathcal{V}_{2,k})+t   \mathscr{T}_{(\mathcal{V}_{1},-\mathcal{V}_{2})}}\langle a^{\mathcal{V}_1}(a^*)^{\mathcal{V}_2}\rangle,
\end{equation}
 
and equivalently, the new average for $\alpha$, where $\alpha$ was defined in \eqref{HatA}, takes the form 

\begin{equation}
	\label{Def:NewAverAlpha}
	\begin{aligned}
		& \langle \alpha(-1)^{\mathcal{V}_1} \alpha(1)^{\mathcal{V}_2}\rangle_*\ := \ \langle \alpha(-1)^{\mathcal{V}_1} \alpha(1)^{\mathcal{V}_2}\rangle_{*,t} \ := \ \langle\langle a^{\mathcal{V}_1}(a^*)^{\mathcal{V}_2}\rangle\rangle\\
		& \ = \ e^{t \sum_{k\in\Lambda^*}\mathscr{T}_{(\mathcal{V}_{1,k},-\mathcal{V}_{2,k})}}\langle \alpha(-1)^{\mathcal{V}_1}\alpha(1)^{\mathcal{V}_2}\rangle.
	\end{aligned}
\end{equation}
We obtain the equation 
\begin{equation}
	\label{KolmogorovEquation2:a}
	\begin{aligned}
		\frac{\partial}{\partial t}\left\langle {\alpha}^{\mathcal{V}_1}({\alpha}^*)^{\mathcal{V}_2}\right\rangle_* \ = \ & \sum_{\mathcal{W}_1,\mathcal{W}_2}\mathcal{U}_{\mathcal{W}_1,\mathcal{W}_2}e^{-{\bf i}t\sum_{k\in\Lambda^*} \omega_k(\mathcal{V}_{1,k}-\mathcal{V}_{2,k})+t  \mathscr{T}_{(\mathcal{V}_{1},-\mathcal{V}_{2})}}\\
		&\times e^{{\bf i}t\sum_{k\in\Lambda^*} \omega_k(\mathcal{W}_{1,k}-\mathcal{W}_{2,k})-t  \mathscr{T}_{(\mathcal{W}_{1},-\mathcal{W}_{2})}}\left\langle {\alpha}^{\mathcal{W}_1}({\alpha}^*)^{\mathcal{W}_2}\right\rangle_*.
	\end{aligned}
\end{equation}

 Thus, in terms of $\alpha$ and the average $\langle \cdot \cdot\rangle_*$, we obtain the equation

\begin{equation}
	\begin{aligned}\label{StartPointAverageAlphaDuhamel}
		& \frac{\partial}{\partial t}\langle \alpha_t(k,1),\alpha_t(k,-1)\rangle_* \  =  \ {\bf i}\lambda\int_{(\Lambda^*)^2}\mathrm{d}k_1\mathrm{d}k_2\delta(-k+k_1+k_2)\\
		& \ \times \mathcal{M}(k,k_1,k_2)e^{{\bf i}t(\omega(k_1)+\omega(k_2)-\omega(k))-t \mathscr{T}(\mathbf{1}_{-k},\mathbf{1}_{k_1},\mathbf{1}_{k_2})}\langle\alpha_t(k,-1) \alpha_t(k_1,1)\alpha_t(k_2,1)\rangle_*,
	\end{aligned}
\end{equation}
where $\mathbf{1}_{-k},\mathbf{1}_{k_1},\mathbf{1}_{k_2}$ are vectors in $\mathbb{N}^{|\Lambda^*|}$ that take values $1$ at the coordinates $-k,k_1,k_2$ respectively and $0$ elsewhere.

We will continue to expand \eqref{StartPointAverageAlphaDuhamel} using the Duhamel expansions. However, the expansion strategy leads to a technical difficulty in treating the wave numbers near   the singular manifold $\mathfrak{S}$, which appears when performing  the estimates on the oscillatory integrals {  and that is defined below}.  As a result, we introduce the following proposition about the existence of  appropriate  cut-off functions that will isolate the singular manifold $\mathfrak{S}$.

For $\aleph^*,\aleph_*=\pm 1,$ let $V=(V_1,\cdots,V_d),W=(W_1,\cdots,W_d)$ be two vectors in $[-\pi,\pi]^d$, and set


	\begin{equation}
	\label{Lemm:Angle:1}\begin{aligned}
		\mathfrak{C}_{\aleph_1^j,\aleph_2^j}^1 \ = \ & -\aleph'\sin(V_1)\sin(2V_j),
		\\
		\mathfrak{C}_{\aleph_1^j,\aleph_2^j}^3 \ = \ &  -\aleph''\sin(W_1)\cos(2W_j),\\
		\mathfrak{C}_{\aleph_1^j,\aleph_2^j}^2 \ = \  & -\sin(2V_j-2W_j)\sin(V_1-W_1)  \ - \aleph'\sin(W_1)\sin(2W_j) \\
		& - \ \aleph''\sin(V_1)\sin(2V_j).
	\end{aligned}
\end{equation}
We consider the equation
\begin{equation}
	\label{Lemm:Angle:2:a}\begin{aligned}
		\mathfrak{C}_{\aleph_1^j,\aleph_2^j}^1r_*^2+\mathfrak{C}_{\aleph_1^j,\aleph_2^j}^2r_*+\mathfrak{C}_{\aleph_1^j,\aleph_2^j}^3 \ = \ 0.
	\end{aligned}
\end{equation}
In the case that \eqref{Lemm:Angle:2:a} has two real roots (or a double real root), we denote these roots by $\tilde{r}_1,\tilde{r}_2$. In the case that \eqref{Lemm:Angle:2:a}   has two complex roots, we denote the real part of the complex root by $\tilde{r}_3$. We set

\begin{equation}
	\begin{aligned}\label{Lemm:Bessel2:2bb8:1}
		\cos(\Upsilon_*(\tilde{r}_i,\aleph^*,\aleph_*))
		\ = \ & {\Big|\aleph^*\tilde{r}_i+\aleph_*+\tilde{r}_i\cos(V_1)e^{{\bf i}2V_j}+\cos(W_1)e^{{\bf i}2W_j}\Big|}\\
		&\ \times \Big[\Big|\aleph^*\tilde{r}_i+\aleph_*+\tilde{r}_i\cos(V_1)e^{{\bf i}2V_j}+\tilde{r}_i\cos(W_1)e^{{\bf i}2W_j}\Big|^2\\
		& \  +\Big|\tilde{r}_i\sin(V_1)e^{{\bf i}2V_j}+\sin(W_1)e^{{\bf i}2W_j}\Big|^2\Big]^{-\frac12},\\
		\sin(\Upsilon_*(\tilde{r}_i,\aleph^*,\aleph_*))\
		= \ & {\Big|\tilde{r}_i\sin(V_1)e^{{\bf i}2V_j}+\sin(W_1)e^{{\bf i}2W_j}\Big|}\\
		&\ \times \Big[\Big|\aleph^*\tilde{r}_i+\aleph_*+\tilde{r}_i\cos(V_1)e^{{\bf i}2V_j}+\cos(W_1)e^{{\bf i}2W_j}\Big|^2\\
		& \ +\Big|\tilde{r}_i\sin(V_1)e^{{\bf i}2V_j}+\sin(W_1)e^{{\bf i}2W_j}\Big|^2\Big]^{-\frac12},
	\end{aligned}
\end{equation}
for $i=1,2,3$.
 We define \begin{equation}\label{Setstar}\begin{aligned}\mathfrak{S}(V,W)=\Big\{2\pi k=(2\pi k^1,\cdots,2\pi k^d)\in [-\pi,\pi]^d \Big|&\  2\pi k^1= \Upsilon_*(\tilde{r}_i,\aleph^*,\aleph_*)\pm \pi; i=1,2,3; \aleph^*,\aleph_*=\pm 1;\\
 		&\ k^j=0,\pm\frac14,\pm\frac12, j=1,\cdots,d \Big\}.\end{aligned}.\end{equation}


We have the following proposition, whose proof is quite standard and is therefore omitted. 
\begin{proposition}\label{Propo:Phi3A} Let $k_1,k_2$ be two momenta in the delta function $\delta(k_0-k_1-k_2)$ and $k_1',k_2'$ be two momenta in the delta function $\delta(k_0'-k_1'-k_2')$.  We assume that $\delta(k_0-k_1-k_2)$ appears after  $\delta(k_0'-k_1'-k_2')$ in the Duhamel expansions.  
	We set $V^1=k_1-k_1'$,  $V^2=k_1+k_1'$, $V^3=k_1+k_2'$,  $V^4=k_1-k_2'$, $W^1=k_2-k_1'$,  $W^2=k_2+k_1'$, $W^3=k_2+k_2'$,  $W^4=k_2-k_2'$ and denote $V^i=(V^1_1,\cdots,V^1_d)$ and $W^i=(W^1_1,\cdots,W^1_d)$. We put $\mathfrak{S}_*=\cup_{i,i'=1}^4\mathfrak{S}(V^i,W^{i'})$.

	For any $\eth>0$,   there exists a  cut-off function $$\Psi _1( k_1,k_2,k_1',k_2'): (\Lambda^*)^4\to [0,1],$$  
	and the smooth   version   
	$${\Psi}_1( k_1,k_2,k_1',k_2'): \mathbb{T}^{4d}\to [0,1],$$ such that $\Psi_1( k_1,k_2,k_1',k_2')= \Psi_1(k_1,k_0,k_3,k_0'),$ on $(\Lambda^*)^4$. And ${\Psi}_1( k_1,k_2,k_1',k_2')=1$ when $	|V^i_j|, \Big|V^i_j\pm\frac{1}{2}\Big|, |W^i_j|, \Big|W_j^i\pm \frac{1}{2}\Big|, |V_j^i\pm W_j^{i'}|   >|\ln\lambda|^{-\eth}$ and $\Psi_2(k_1,k_2,k_3, k_0')=0$ when either $  	|V^i_j|, \Big|V^i_j\pm\frac{1}{2}\Big|, |W^i_j|, \Big|W_j^i\pm \frac{1}{2}\Big|$, or $ |V_j^i\pm W_j^{i'}| <|\ln\lambda|^{-\eth}/2$, for $j=1,\cdots,d$, $i,i'=1,2$.

	Moreover, there also exists a  cut-off function
	$$\Psi _2 (k_1,k_2,k_1',k_2'): (\Lambda^*)^4\to [0,1],$$  
	and the smooth rescaled version     
	$${\Psi}_2( k_1,k_2,k_1',k_2'): \mathbb{T}^{4d}\to [0,1],$$ such that $\Psi_2(k_1,k_2,k_1',k_2')= \Psi_2(k_1,k_2,k_1',k_2')$ on $(\Lambda^*)^4$. Moreover, ${\Psi}_2(k_1,k_2,k_1',k_2')=1$ when $d( k_1,\mathfrak{S}_*)>|\ln\lambda|^{-\eth}$ and $d( k_2,\mathfrak{S}_*)>|\ln\lambda|^{-\eth}$, and ${\Psi}_2(k_1,k_2,k_1',k_2') =0$ when $d( k_1,\mathfrak{S}_*)<|\ln\lambda|^{-\eth}/2,$ or $d( k_2,\mathfrak{S}_*)<|\ln\lambda|^{-\eth}/2.$ 
	
	We finally set $\Psi_3=\Psi_1\Psi_2+\sum_{i=1}^4\mathbf{1}_{V^i=0}+\sum_{i=1}^4\mathbf{1}_{W^i=0}$ and we extend $\Psi_3$ to have a  continuum version  from $\Lambda^*$ to    $\mathbb{T}$.

\end{proposition}

\begin{definition}\label{def:Phi3}
	For any momenta $k_1,k_2\in\Lambda^*$ that are associated to a delta function $\delta(\sigma_0k_0+\sigma_1k_1+\sigma_2k_2)$ in the Duhamel expansions, we define the set 
	\begin{equation}\begin{aligned}
			\mathcal{S}_{(\sigma_0,k_0,\sigma_1, k_1,\sigma_2, k_2 )} : = & \Big\{(\sigma_0'k_0',\sigma_1'k_1',\sigma_2'k_2) ~~ | ~~ \delta(\sigma_0'k_0'+\sigma_1'k_1'+\sigma_2'k_2')\\
			& \mbox{ appears before } \delta(\sigma_0k_0+\sigma_1k_1+\sigma_2k_2) \mbox{  in the Duhamel expansions}\Big\},\end{aligned}
	\end{equation}
	and put
	
	\begin{equation}\begin{aligned}		\label{def:Phi3:1}\Phi_4(\sigma_0,k_0,\sigma_1, k_1,\sigma_2, k_2 ) \ = \ & \prod_{(\sigma_0'k_0',\sigma_1'k_1',\sigma_2'k_2) \in \mathcal{S}_{(\sigma_0,k_0,\sigma_1, k_1,\sigma_2, k_2 )} }\\
			& \Psi_3(-\sigma_0 \sigma_1 k_1,-\sigma_0 \sigma_2 k_2 ,-\sigma_0'\sigma_1' k_1',-\sigma_0'\sigma_2'k_2').\end{aligned}	
	\end{equation}
	
	Moreover,  we also have their extended  versions
	${\Phi}_4$ from $\Lambda^*$ to  $\mathbb{T}$ as in Proposition \ref{Propo:Phi3A} .
\end{definition}  
The following lemma is useful in contructing the next cut-off function.
\begin{lemma}
	\label{Lemma:Cutoffn}
	For any $m\in\mathbb{N}$, $m\ge 1$, we set 
	\begin{equation}
		\label{Lemma:Cutoffn:1}\begin{aligned}
			\mathfrak{S}_{2m,\Lambda_*} \ :=\ & \Big\{(k_1,\cdots,k_{2m})\in \times\Lambda_*^{2m}~~\Big|~~\\
			& \int_0^{1}\mathrm d \tau'[\tilde\omega(k_1)\cdots\tilde\omega(k_{2m})]^\frac12	 h^{dm}\langle |a_{k_{1}}|\cdots|a_{k_{{2m}}}|\rangle_{\tau'\lambda^{-2}} > \ln|\ln|\ln\lambda|||\Big\},
			\end{aligned}
	\end{equation}
then
	\begin{equation}\label{Lemma:Cutoffn:2}
	\mathscr{W}\Big(\lim_{D\to \infty }\mathfrak{S}_{2m,\Lambda_*}\Big) \lesssim
	\mathcal{C}_o^m\ln|\ln|\ln\lambda|||^{-1}.\end{equation}
\end{lemma}
\begin{proof}
	From \eqref{Propo:ExpectationEqui:2}, we find
		\begin{equation}\label{Lemma:Cutoffn:E1}
		\begin{aligned}
		&	\int_{\Lambda_*^{2m}}{\mathrm{d}k_1\cdots\mathrm{d}k_n}	[\tilde\omega(k_1)\cdots\tilde\omega(k_{2m})]^\frac12	\left \langle |a_{k_{1}}|\cdots|a_{k_{{2m}}}| \right\rangle_t \\
		\lesssim 	&\ \	\left[\int_{\Lambda_*^{2m}}{\mathrm{d}k_1\cdots\mathrm{d}k_n}	[\tilde\omega(k_1)\cdots\tilde\omega(k_{2m})]	\left \langle |a_{k_{1}}|^2\cdots|a_{k_{{2m}}}| ^2\right\rangle_t\right]^\frac12 \\
			\lesssim  &\ \	h^{-dm}\mathcal{C}_o^m,
		\end{aligned}
	\end{equation}
 which implies
	\begin{equation}\label{Lemma:Cutoffn:E3}
	\begin{aligned}
		&	\int_{\mathfrak{S}_{2m,\Lambda_*}}{\mathrm{d}k_1\cdots\mathrm{d}k_n}	[\tilde\omega(k_1)\cdots\tilde\omega(k_{2m})]^\frac12	 h^{dm}\left\langle |a_{k_{1}}|\cdots|a_{k_{{2m}}}| \right\rangle_t 
		\lesssim  &\ \	\mathcal{C}_o^m,
	\end{aligned}
\end{equation}
yielding 
	\begin{equation}\label{Lemma:Cutoffn:E4}
\mathscr{W}\Big(\lim_{D\to \infty }\mathfrak{S}_{2m,\Lambda_*}\Big) \lesssim
	\mathcal{C}_o^m\ln|\ln|\ln\lambda|||^{-1}.\end{equation}
\end{proof}
\begin{definition}\label{def:hbar}
	For $(k_1,\cdots,k_{2m})\in [0,1]\times\Lambda_*^{2m}$, we define the cut-off function
	\begin{equation}\begin{aligned}		\label{def:hbar:1}\hbar_{2m}(k_1,\cdots,k_{2m}) \ := \  1-\chi_{\mathfrak{S}_{2m,\Lambda_*}} \ = \ & 0 \mbox{ when } (k_1,\cdots,k_{2m})\in \mathfrak{S}_{2m,\Lambda_*}\\
	= \ & 1 \mbox{ otherwise}.	\end{aligned}	
	\end{equation}
Moreover,  we also have their extended  versions
$\hbar_{\tau,2m},$ on  $[0,1]\times\mathbb{T}^{2dm}$ as in Proposition \ref{Propo:Phi3A} .
\end{definition}

{ 


Now, we follow the soft partial time integration technique in the work of Lukkarinen-Spohn \cite{LukkarinenSpohn:WNS:2011}, that  takes its  origin from  the work of Erdos-Yau   \cite{erdHos2000linear}, to control the multi-layer Duhamel expansion. In shortening the notations, we set
\begin{equation}
	\label{Def:Omega}\begin{aligned}
		\mathfrak{X}(\sigma,k,\sigma',k_1,\sigma'',k_2) \ = & \ \sigma\omega(k) \ + \ \sigma' \omega(k_1) \ + \ \sigma'' \omega(k_2).
	\end{aligned}
\end{equation}
We observe that using \eqref{StartPoint} and  \eqref{Def:Omega} we have}
\begin{equation}
	\begin{aligned} \label{eq:duhamelproductlambda}
	\frac{	\mathrm{d}}{	\mathrm{d}t}\Big[e^{\varsigma	 t}\prod_{i=1}^2 & \alpha_t(k_i,\sigma_i) \Big]\ = \  \\
		=  & \varsigma	 e^{\varsigma t}\prod_{i=1}^n \alpha_t(k_i,\sigma_i) \ - \ {\bf i}\lambda \sum_{j=1}^2\sigma_j\prod_{i=1,i\ne j}^{n} \alpha_t(k_i,\sigma_i) \Big[ \iint_{(\Lambda^*)^2} \mathrm{d}k'_1 \mathrm{d}k'_2\\
		&\  \times  \delta(-k_j+k_1'+k_2') \mathcal{M}(k_j,k_1',k_2')\times \exp\Big[\varsigma t-{\bf i}t\mathfrak{X}(\sigma_j,k_j,\sigma_j,k_1',\sigma_j,k_2')\Big] \alpha_t(k_1',\sigma_j) \alpha_t(k_2',\sigma_j)\Big], 
\end{aligned}\end{equation}
where $\varsigma>0$ is a control parameter to  be specified later.

\medskip
{\it The full Duhamel expansion.} 
By repeating the expansion process $\mathfrak{N}$ times, using both the soft partial time and modified soft partial time integrations, we then obtain a multi-layer expression  in which  the time interval $[0,t]$ is divided into $\mathfrak{N}+1$ time slices $[0,s_0],$ $[s_0,s_0+s_1]$, $\dots$, $[s_0+\cdots+s_{\mathfrak{N}-1},t]$ and $t=s_0+\cdots+s_{\mathfrak{N}}.$ Let us write down the final result of this process in the following schematic manner:
\begin{equation}
	\begin{aligned} \label{eq:fullDuhamel}
		& \delta_{k=k'}\langle\alpha_t(k,-1)\alpha_t(k',1)\rangle_*\  
		=  \  \sum_{n=0}^{\mathfrak{N}-1}\langle\mathcal{F}_{n}^0(t,k,-1,\Gamma)[\alpha_0] \rangle_*\ + \ \sum_{n=0}^{\mathfrak{N}-1}\varsigma_n\int_0^t\mathrm{d}s \langle\mathcal{F}_{n}^1(s,t,k,-1,\Gamma)[\alpha_s]\rangle_*\\
		&\ +\ \sum_{n=1}^{\mathfrak{N}}\int_0^t\mathrm{d}s \langle\mathcal{F}_{n}^2(s,t,k,-1,\Gamma)[\alpha_s]\rangle_* \ + \ \int_0^t\mathrm{d}s \langle\mathcal{F}_{\mathfrak{N}}^3(s,t,k,-1,\Gamma)[\alpha_s]\rangle_*.
\end{aligned}\end{equation} 

{In the formulation \eqref{eq:fullDuhamel}, $\Gamma$ denotes the soft partial time integration vector.
\begin{equation}
	\label{SolfPartTimeVector}(\varsigma_0,\cdots,\varsigma_{\mathfrak N-1}) \mbox{ in } \mathbb{R}_+^{\mathfrak N}.
\end{equation}
}
We set
\begin{equation}
	\label{Def:Para0}
	\mathfrak{N}
	=\max\left(1,	\left[\frac{\mathfrak{N}_0|\ln\lambda|}{\ln\langle\ln\lambda\rangle}	\right]^{c_{\mathfrak{N}}}\right), 
\end{equation}
where $[x]$ is the integer part that satisfies $[x]\le x <[x]+1$ and $\mathfrak{N}_0$ is any number in $\mathbb{R}_+$, for $0<c_{\mathfrak{N}}\le\frac12$. And we 
 define the ``stopping rule''


We also set the parameter that controls the partial time integration of \eqref{eq:duhamelproductlambda}
\begin{equation}
	\label{Def:Para2}
	\varsigma' = \lambda^2\mathfrak{N}^{\wp	},
\end{equation}
in which $\wp$ is a positive constant. The components of the soft partial time integration vector are defined as 
\begin{equation}
	\label{Def:Para3}
	\varsigma_n = 0 \mbox{ when } 0\le n<[\mathfrak{N}/4], \mbox{ and }\varsigma_n = \varsigma'   \mbox{ when } [\mathfrak{N}/4]\le n\le \mathfrak{N}. 
\end{equation}

 {The first term in \eqref{eq:fullDuhamel}} has the following explicit form   (for $n\ge 1$)
\begin{equation}
	\begin{aligned} \label{eq:DefE0}
		&\langle\mathcal{F}_{n}^0(t,k_{n,1},\sigma_{n,1},\Gamma)[\alpha]\rangle_*\\
		\ = \  &  
		(-{\bf i}\lambda)^n\sum_{\substack{\rho_i\in\{1,\cdots,n-i+2\},\\ i\in\{0,\cdots,n-1\}}}\sum_{\substack{\bar\sigma\in \{\pm1\}^{\mathcal{I}_n},\\ \sigma_{i,\rho_i}+\sigma_{i-1,\rho_i}+\sigma_{i-1,\rho_i+1}\ne \pm3,
				\\ \sigma_{i-1,\rho_i}\sigma_{i-1,\rho_i+1}= 1}} \int_{(\Lambda^*)^{\mathcal{I}_n}}\mathrm{d}\bar{k}\Delta_{n,\rho}(\bar{k},\bar\sigma)h^{d}\\
		&\ \ \ \ \ \ \ \ \ \ \ \ \ \times \prod_{i=1}^n\Big[\sigma_{i,\rho_i}\mathcal{M}( k_{i,\rho_i}, k_{i-1,\rho_i}, k_{i-1,\rho_i+1})\left\langle\prod_{i=1}^{n+2}\alpha(k_{0,i},\sigma_{0,i})\right\rangle_*\\
		&\ \ \ \ \ \ \ \ \ \ \ \ \ \times \Phi_{1,i}( \sigma_{i-1,\rho_i},k_{i-1,\rho_i}, \sigma_{i-1,\rho_i+1},k_{i-1,\rho_i+1})\Big]\int_{(\mathbb{R}_+)^{\{0,\cdots,n\}}}\mathrm{d}\bar{s} \delta\left(t-\sum_{i=0}^ns_i\right)\prod_{i=0}^{n}e^{-s_i\varsigma_{n-i}}\\
		&\ \ \ \ \ \ \ \ \ \ \ \ \   \times\prod_{i=1}^{n}e^{-{\bf i}t_i(s)\mathfrak{X}(\sigma_{i,\rho_i}, k_{i,\rho_i},\sigma_{i-1,\rho_i},k_{i-1,\rho_i},\sigma_{i-1,\rho_i},k_{i-1,\rho_i+1})}\prod_{i=0}^{n}e^{-s_i \tau_{i}},
\end{aligned}\end{equation} 
in which for each delta function   
and \begin{equation}
	\label{CutoffPhii}\begin{aligned}
		&  \Phi_{1,i}( \sigma_{i-1,\rho_i},k_{i-1,\rho_i}, \sigma_{i-1,\rho_i+1},k_{i-1,\rho_i+1}) \ = \  \Phi_{1}( \sigma_{i-1,\rho_i},k_{i-1,\rho_i}, \sigma_{i-1,\rho_i+1},k_{i-1,\rho_i+1})\\ & \ =\ \Phi_{4}( \sigma_{i-1,\rho_i},k_{i-1,\rho_i}, \sigma_{i-1,\rho_i+1},k_{i-1,\rho_i+1}) \hbar_{s_{i-1}\lambda^2,i}(k_{n-i,1},\cdots,k_{n-i,i+2}) \mbox{ for } 1\le i\le n-1, \mbox{ $i$ is even}, \\
		&  \Phi_{1,i} \ = \  1 \mbox{ for } 1\le i\le n-1, \mbox{ $i$ is odd}, \\
		& \Phi_{1,i} + \Phi_{0,i} \ = \ 1, \ \ \  \Phi_{1} + \Phi_{0} \ = \ 1.  \end{aligned}
\end{equation}
We also have the following equivalent form of expressing the total phases 
\begin{equation}\label{GraphSec:E1}
	\begin{aligned}
		& \prod_{i=0}^{n}e^{-s_i\varsigma_{n-i}}
		\prod_{i=1}^{n}e^{-\mathbf{i}t_i(s)\mathfrak{X}(\sigma_{i,\rho_i}, k_{i,\rho_i},\sigma_{i-1,\rho_i},k_{i-1,\rho_i},\sigma_{i-1,\rho_i},k_{i-1,\rho_i+1})}\\
		\ =\ &\prod_{i=0}^{n}e^{-s_i\varsigma_{n-i}}
		\prod_{i=1}^{n}e^{-\mathbf{i}\left(\sum_{j=1}^{i-1}s_j\right)\mathfrak{X}(\sigma_{i,\rho_i}, k_{i,\rho_i},\sigma_{i-1,\rho_i},k_{i-1,\rho_i},\sigma_{i-1,\rho_i},k_{i-1,\rho_i+1})}
		\ =\   \prod_{i=0}^n e^{-\mathbf{i} s_i \vartheta_i},\end{aligned}
\end{equation}
where 
\begin{equation}\label{GraphSec:E2}
	\begin{aligned}
		\vartheta_i \ := \ & \sum_{l=i+1}^n \mathfrak{X}(\sigma_{l,\rho_l},k_{l,\rho_l},\sigma_{l-1,\rho_l},  k_{l-1,\rho_l},\sigma_{l-1,\rho_l+1}, k_{l-1,\rho_l+1}) - {\bf i}\varsigma_{n-i},
	\end{aligned}
\end{equation} 
in which the phase is defined in \eqref{Def:Omega}.  We set $\vartheta_n=0$.

Note that $\Phi_{1,i}(n)$ are  functions of $n$. However, in our computations, we simply write $\Phi_{1,i}$ since they are associated with Feynman diagrams that have $n$-layers.  

Moreover, we have used the notations for the set of indices
\begin{equation}
	\label{IndexSet2}\mathcal{I}_{n}  \ := \ \{(j,l) \  | \ 0\le j\le n-1, 1\le l \le n-j+2\}. \end{equation}

\begin{definition}[Phase Regulators]\label{Def:PhaseRegulator}
	We define the phase regulators to be the quantities
	\begin{equation}
		\label{Def:TauGen}
		\tau_i \ = \ \mathscr{T}_{\left(\mathbf{1}_{k_{i,j}\sigma_{i,j}}\right)_{j=1}^{n-i+2}} \quad \mbox{ for } i\in\{0,\cdots,n-1\},\ \ \ \tau_n=0,
	\end{equation}
\end{definition}
where $\left(\mathbf{1}_{k_{i,j}\sigma_{i,j}}\right)_{j=1}^{n-i+2}$ denotes a vector  in $\mathbb{N}^{|\Lambda^*|}$ that take values $1$ at the coordinates $ k_{i,j}$ respectively and $0$ elsewhere and we have used { the same definition of} $\mathscr{T}_{(\mathcal{V}_{1},-\mathcal{V}_{2})}$ in   \eqref{KolmogorovEquation2}. 
We also define
\begin{equation}
	\label{Def:Ti}t_i(s)=\sum_{j=1}^{i-1}s_j,\end{equation}
$\bar{s}$, $\bar{k}$ and $\bar\sigma$ are vectors representing all of the quantities $s_i$, $k_{i,j},\sigma_{i,j}$ appearing in the integration. The number $\rho_i$ encodes the position where the splitting happens hence the use of
\begin{equation}
	\label{Def:Rho}\delta\left(\sigma_{i,\rho_i}k_{i,\rho_i}+\sigma_{i-1,\rho_i}k_{i-1,\rho_i}+\sigma_{i-1,\rho_i+1}k_{i-1,\rho_i+1}\right).
\end{equation}
In the above summation we only allow $$(\sigma_{i,\rho_i},\sigma_{i-1,\rho_i},\sigma_{i-1,\rho_i+1})\ne (+1,+1,+1),(-1,-1,-1),\ \ \ \sigma_{i-1,\rho_i}\sigma_{i-1,\rho_i+1}= 1,$$ that means we have either $k_{i,\rho_i}-k_{i-1,\rho_i}-k_{i-1,\rho_i+1}=0$, or $-k_{i,\rho_i}+k_{i-1,\rho_i}+k_{i-1,\rho_i+1}=0$. The case $k_{i,\rho_i}+k_{i-1,\rho_i}+k_{i-1,\rho_i+1}=0$ is not allowed. Moreover, $$\mathcal{F}_{0}^0(t,k,\sigma,\Gamma)[\alpha]=e^{-\varsigma_0 t}\langle\alpha(k,\sigma)\alpha(k,-\sigma)\rangle_*$$
and finally
the function $\Delta_{n,\rho}$ contains all of the $\delta$-functions. More precisely,
\begin{equation}
	\begin{aligned} \label{eq:DefDelta}
		\Delta_{n,\rho}(\bar{k},\bar\sigma) 
		=  &\  
		\prod_{i=1}^{n}\Big\{\prod_{l=1}^{\rho_i-1}\Big[\delta(k_{i,l}-k_{i-1,l})\mathbf{1}(\sigma_{i,l}=\sigma_{i-1,l})\Big] \\
		&\times \mathbf{1}(\sigma_{i,\rho_i}=-\sigma_{i-1,\rho_i}) \delta\left(\sigma_{i,\rho_i}k_{i,\rho_i}+\sigma_{i-1,\rho_i}k_{i-1,\rho_i}+\sigma_{i-1,\rho_i+1}k_{i-1,\rho_i+1}\right)\\
		&\times\prod_{l=\rho_i+1}^{n-i+2}\Big[\delta(k_{i,l}-k_{i-1,l+1})\mathbf{1}(\sigma_{i,l}=\sigma_{i-1,l
			+1})\Big]\mathbf{1}(k_{n,1}=k_{n,2})\mathbf{1}(\sigma_{n,1}+\sigma_{n,2}=0)\Big\},
\end{aligned}\end{equation} 
where $\bar{k}$ and $\bar\sigma$ represent all of the quantities $k_{i,j},\sigma_{i,j}$ appearing in the formula.
This function makes sure that 
{ 
$$k_{i,j}=k_{i-1,j}, \quad \sigma_{i,j}=\sigma_{i-1,j}, \quad \mbox{ for } j\in\{1,\cdots,\rho_i-1\},$$  
$$k_{i,j}=k_{i-1,j+1}, \quad \sigma_{i,j}=\sigma_{i-1,j+1}, \quad \mbox{ for  }j\in\{\rho_i+1,\cdots,{n-i}+2\},$$  
$$  k_{n,1}=k_{n,2}, \quad \sigma_{n,1}+\sigma_{n,2}=0,\quad \sigma_{i,\rho_i}k_{i,\rho_i}+\sigma_{i-1,\rho_i}k_{i-1,\rho_i}+\sigma_{i-1,\rho_i+1}k_{i-1,\rho_i+1}=0,\quad  \sigma_{i,\rho_i}=-\sigma_{i-1,\rho_i},$$ }
as discussed in the previous subsection. 

%
%
%
%

We define the last  { term in \eqref{eq:fullDuhamel}} $\mathcal{F}_n^3$  as:
\begin{equation}
	\begin{aligned} \label{eq:DefE3}
		& \langle\mathcal{F}^3_{n}(s_0,t,k_{n,1},\sigma_{n,1},\Gamma)[\alpha]\rangle_*\\ 
		= & \  
		(-{\bf i}\lambda)^n\sum_{\substack{\rho_i\in\{1,\cdots,n-i+2\},\\ i\in\{0,\cdots,n-1\}}}\sum_{\substack{\bar\sigma\in \{\pm1\}^{\mathcal{I}_n},\\ \sigma_{i,\rho_i}+\sigma_{i-1,\rho_i}+\sigma_{i-1,\rho_i+1}\ne \pm3,
				\\ \sigma_{i-1,\rho_i}\sigma_{i-1,\rho_i+1}= 1}} \int_{(\Lambda^*)^{\mathcal{I}_n}}\mathrm{d}\bar{k}\Delta_{n,\rho}(\bar{k},\bar\sigma)h^{d}\\
		&\ \ \ \ \ \ \ \ \ \ \ \ \ \times \prod_{i=1}^n\Big[\sigma_{i,\rho_i}\mathcal{M}( k_{i,\rho_i}, k_{i-1,\rho_i}, k_{i-1,\rho_i+1})\Phi_{1,i}(\sigma_{i,\rho_i},k_{i,\rho_i}, \sigma_{i-1,\rho_i},k_{i-1,\rho_i}, \sigma_{i-1,\rho_i+1},k_{i-1,\rho_i+1})\Big]\\
		&\ \ \ \ \ \ \ \ \ \ \ \ \ \times \left\langle\prod_{i=1}^{n+2}\alpha(k_{0,i},\sigma_{0,i})\right\rangle_* \int_{(\mathbb{R}_+)^{\{1,\cdots,n\}}}\mathrm{d}\bar{s} \delta\left(t-\sum_{i=0}^ns_i\right)\prod_{i=1}^{n}e^{-s_i \varsigma_{n-i}}\prod_{i=0}^{n}e^{-s_i \tau_{i}}\\
		&\ \ \ \ \ \ \ \ \ \ \ \ \   \times\prod_{i=1}^{n}e^{-{\bf i}t_i(s)\mathfrak{X}(\sigma_{i,\rho_i}, k_{i,\rho_i},\sigma_{i-1,\rho_i},k_{i-1,\rho_i},\sigma_{i-1,\rho_i},k_{i-1,\rho_i+1})}.
\end{aligned}\end{equation} 
The only difference between the two formulas \eqref{eq:DefE0} and \eqref{eq:DefE3} is the integration with respect to $\mathrm{d}s$. In \eqref{eq:DefE0} the integration is taken over $(\mathbb{R}_+)^{\{0,\cdots,n\}}$ and in \eqref{eq:DefE3} it is over $(\mathbb{R}_+)^{\{1,\cdots,n\}}.$ 

\smallskip
{  The second term in  \eqref{eq:fullDuhamel}, namely  $\mathcal{F}_n^1$ is defined using the last term}, as follows. We set, 
$$\langle\mathcal{F}^1_{0}(s,t,k,\sigma,\Gamma)[\alpha]\rangle_*=e^{-\varsigma_0 (t-s)}\langle\alpha(k,\sigma) \alpha(k,-\sigma)\rangle_*,$$ 
and for $n>0$, we set
\begin{equation}
	\label{eq:DefE1}
	\langle\mathcal{F}^1_n(s,t,k_{n,1},\sigma_{n,1},\Gamma)[\alpha]\rangle_* \ = \ \int_0^{t-s}\mathrm{d}r e^{-r\varsigma_n}\langle\mathcal{F}^3_{n}(s+r,t,k_{n,1},\sigma_{n,1},\Gamma)[\alpha]\rangle_*.
\end{equation}

And,  
\begin{equation}
	\begin{aligned} \label{eq:DefE2}
		&\langle\mathcal{F}^2_{n}(s_0,t,k_{n,1},\sigma_{n,1},\Gamma)[\alpha]\rangle_*\\
		=  \ &
		(-{\bf i}\lambda)^n\sum_{\substack{\rho_i\in\{1,\cdots,n-i+2\},\\ i\in\{0,\cdots,n-1\}}}\sum_{\substack{\bar\sigma\in \{\pm1\}^{\mathcal{I}_n},\\ \sigma_{i,\rho_i}+\sigma_{i-1,\rho_i}+\sigma_{i-1,\rho_i+1}\ne \pm3,
				\\ \sigma_{i-1,\rho_i}\sigma_{i-1,\rho_i+1}= 1}} \int_{(\Lambda^*)^{\mathcal{I}_n}}\mathrm{d}\bar{k}\Delta_{n,\rho}(\bar{k},\bar\sigma)h^{d}\\
		&  \ \times \sigma_{1,\rho_1}\mathcal{M}( k_{1,\rho_1}, k_{0,\rho_1}, k_{0,\rho_1+1})\Phi_{0,1}(\sigma_{0,\rho_1}, k_{0,\rho_1},\sigma_{0,\rho_1+1},k_{0,\rho_1+1})\\
		& \ \times \prod_{i=2}^n\Big[\sigma_{i,\rho_i}\mathcal{M}( k_{i,\rho_i}, k_{i-1,\rho_i}, k_{i-1,\rho_i+1})\Phi_{1,i}(\sigma_{i-1,\rho_i}, k_{i-1,\rho_i},\sigma_{i-1,\rho_i+1},k_{i-1,\rho_i+1})\Big] \\
		&\ \times   \left\langle\prod_{i=1}^{n+2}\alpha(k_{0,i},\sigma_{0,i})\right\rangle_*\int_{(\mathbb{R}_+)^{\{1,\cdots,n\}}}\mathrm{d}\bar{s} \delta\left(t-\sum_{i=0}^ns_i\right)\prod_{i=1}^{n}e^{-s_i \varsigma_{n-i} }\prod_{i=0}^{n}e^{-s_i \tau_{i}}\prod_{i=0}^{n}e^{-s_i \tau_{i}}\\
		& \ \times\prod_{i=1}^{n}e^{-{\bf i}t_i(s)\mathfrak{X}(\sigma_{i,\rho_i}, k_{i,\rho_i},\sigma_{i-1,\rho_i}, k_{i-1,\rho_i},\sigma_{i-1,\rho_i+1},  k_{i-1,\rho_i+1})}.
\end{aligned}\end{equation}
 

The multi-layer Duhamel expansion can then be expressed via Feynman diagrams, as we will explain in Section \ref{Duhamel} below. 

A straightforward computation gives the following { propositions. }

\begin{proposition}[First term]\label{Proposition:Ampl1} For  $\mathfrak{N}\ge 1$, the terms associated to $\mathcal{F}_n^0$ in equation \eqref{eq:fullDuhamel} can be written as
	\begin{equation}\label{Proposition:Ampl:1}
		Q^{1} 
		\ = \  \sum_{n=1}^{\mathfrak{N}-1}\mathcal{G}^{0}_{n }( t,k,\sigma,\Gamma)
	\end{equation}
	where
	\begin{equation}
		\begin{aligned} \label{eq:DefFmain}
			&\mathcal{G}^{0}_{n }( t,k,\sigma,\Gamma)\  
			=  \  
			(-{\bf i}\lambda)^n\sum_{\substack{\rho_i\in\{1,\cdots,n-i+2\},\\ i\in\{0,\cdots,n-1\}}}\sum_{\substack{\bar\sigma\in \{\pm1\}^{\mathcal{I}_n},\\ \sigma_{i,\rho_i}+\sigma_{i-1,\rho_i}+\sigma_{i-1,\rho_i+1}\ne \pm3,
					\\ \sigma_{i-1,\rho_i}\sigma_{i-1,\rho_i+1}= 1}} h^{d} \\
			&\  \times \int_{(\Lambda^*)^{\mathcal{I}_n}}\mathrm{d}\bar{k}\Delta_{n,\rho}(\bar{k},\bar\sigma)\left\langle\prod_{i=1}^{n+2}\alpha(k_{0,i},\sigma_{0,i})\right\rangle_{0} \prod_{i=1}^n\Big[\sigma_{i,\rho_i}\mathcal{M}( \sigma_{i,\rho_i}, k_{i,\rho_i},\sigma_{i-1,\rho_i}, k_{i-1,\rho_i},\sigma_{i-1,\rho_i+1}, k_{i-1,\rho_i+1})\\
			& \ \times \Phi_{1,i}(\sigma_{i-1,\rho_i}, k_{i-1,\rho_i},\sigma_{i-1,\rho_i+1},k_{i-1,\rho_i+1})\Big]  \int_{(\mathbb{R}_+)^{\{0,\cdots,n\}}}\mathrm{d}\bar{s} \delta\left(t-\sum_{i=0}^ns_i\right)\prod_{i=0}^{n}e^{-s_i\varsigma_{n-i}}\prod_{i=0}^{n}e^{-s_i \tau_{i}}\\
			&\ \times \prod_{i=1}^{n}e^{-{\bf i}t_i(s)\mathfrak{X}( \sigma_{i,\rho_i}, k_{i,\rho_i},\sigma_{i-1,\rho_i}, k_{i-1,\rho_i},\sigma_{i-1,\rho_i+1}, k_{i-1,\rho_i+1})}\mho(\{(k_{0,1},\sigma_{0,1}),\cdots,(k_{0,n},\sigma_{0,n}),(k_{1,\rho_1},\cdots,k_{n,\rho_n})\}),
	\end{aligned}\end{equation}
in which $\mho(\{(k_{0,1},\sigma_{0,1}),\cdots,(k_{0,n},\sigma_{0,n}),(k_{1,\rho_1},\cdots,k_{n,\rho_n})\})=1$ if the set $\{(k_{0,1},\sigma_{0,1}),\cdots,(k_{0,n},\sigma_{0,n})\}$ is admissible and  all $k_{i,\rho_i}$ are different from $0$; otherwise, the function is $0$.
\end{proposition}

\begin{proposition}[Third term]\label{Proposition:Ampl2} The terms associated to $\mathcal{F}_n^2$   can be written as 
	\begin{equation}\label{Def:OperatorQErrors:Bis2}
		\begin{aligned}
			 Q^{3}  := & \ \sum_{n=1}^{\mathfrak N}\int_0^t\mathrm{d}s \mathcal{F}_{n}^2(s,t,k,\sigma,\Gamma)[\alpha_s] \ 
			=     \ \sum_{n=1}^{\mathfrak N}\int_0^t\mathrm{d}s\mathcal{G}^{2,pair}_{n}(s,t,k,\sigma,\Gamma)+ \ \sum_{n=1}^{\mathfrak N}\int_0^t\mathrm{d}s\mathcal{G}^{2,nonpair}_{n}(s,t,k,\sigma,\Gamma),
		\end{aligned}
	\end{equation}
	{ 
	$$=:Q^{3,pair} \ + \ Q^{3,nonpair}$$}
	where
	\begin{equation}
		\begin{aligned} \label{eq:DefF2}
			&\mathcal{G}^{2,pair}_{n}(s_0,t,k,\sigma,\Gamma)\  
			=  \ (-{\bf i}\lambda)^n\sum_{\substack{\rho_i\in\{1,\cdots,n-i+2\},\\ i\in\{0,\cdots,n-1\}}}\sum_{\substack{\bar\sigma\in \{\pm1\}^{\mathcal{I}_n},\\ \sigma_{i,\rho_i}+\sigma_{i-1,\rho_i}+\sigma_{i-1,\rho_i+1}\ne \pm3,
					\\ \sigma_{i-1,\rho_i}\sigma_{i-1,\rho_i+1}= 1}}\\
			&\ \times \int_{(\Lambda^*)^{\mathcal{I}_n}}\mathrm{d}\bar{k}\Delta_{n,\rho}(\bar k,\bar\sigma)\sigma_{1,\rho_1}\mathcal{M}( k_{1,\rho_1}, k_{0,\rho_1}, k_{0,\rho_2}) \Phi_{0,1}(\sigma_{0,\rho_1} ,k_{0,\rho_1},\sigma_{0,\rho_2}, k_{0,\rho_2})h^{d}\\
			&\ \times\prod_{i=2}^n\Big[\sigma_{i,\rho_i}\mathcal{M}( k_{i,\rho_i}, k_{i-1,\rho_i}, k_{i-1,\rho_i+1})\Phi_{1,i}(\sigma_{i-1,\rho_i} k_{i-1,\rho_i},\sigma_{i-1,\rho_i+1}k_{i-1,\rho_i+1})\Big]\\
			&\ \times e^{{\bf i}s_0\vartheta_0} \left\langle\prod_{i=1}^{n+2}\alpha(k_{0,i},\sigma_{0,i})\right\rangle_{s_0}  \int_{(\mathbb{R}_+)^{\{1,\cdots,n\}}}\mathrm{d}\bar{s} \delta\left(t-\sum_{i=0}^ns_i\right)\prod_{i=1}^{n}e^{-s_i\varsigma_{n-i}}\prod_{i=0}^{n}e^{-s_i \tau_{i}}\\
			&\  \times \prod_{i=1}^{n}e^{-{\bf i}t_i(s)\mathfrak{X}( \sigma_{i,\rho_i},k_{i,\rho_i},\sigma_{i-1,\rho_i}, k_{i-1,\rho_i},\sigma_{i-1,\rho_i+1},k_{i-1,\rho_i+1})}\mho(\{(k_{0,1},\sigma_{0,1}),\cdots,(k_{0,n},\sigma_{0,n}),(k_{1,\rho_1},\cdots,k_{n,\rho_n})\}),
	\end{aligned}\end{equation}
in which $\mho(\{(k_{0,1},\sigma_{0,1}),\cdots,(k_{0,n},\sigma_{0,n}),(k_{1,\rho_1},\cdots,k_{n,\rho_n})\})=1$ if the set $\{(k_{0,1},\sigma_{0,1}),\cdots,(k_{0,n},\sigma_{0,n})\}$ is admissible and  all $k_{i,\rho_i}$ are different from $0$; otherwise, the function is $0$.
We split
\begin{equation}
	\label{lemma:Q3FinalEstimate:2}
	Q^3 \ = \ Q^{3,pair} \ + \ Q^{3,nonpair} 
\end{equation}
where 	$Q^{3,pair}$ is the sum over expansions that are admissible  and all $k_{i,\rho_i}$ are different from $0$.
\end{proposition}

\begin{proposition}[Last term]\label{Proposition:Ampl3}
	The component associated to $\mathcal{F}_{\mathfrak{N}}^3$ can be expressed as
	\begin{equation}\label{Def:OperatorQErrors:Bis3}
		\begin{aligned}
			Q^4= & \ \int_0^t\mathrm{d}s \mathcal{F}_{\mathfrak{N}}^3(s,t,k,\sigma,\Gamma)[\alpha_s] \
			=   \  \int_0^t\mathrm{d}s\mathcal{G}^{3,pair}_{\mathfrak{N}}(s,t,k,\sigma,\Gamma)\
			+   \  \int_0^t\mathrm{d}s\mathcal{G}^{3,nonpair}_{\mathfrak{N}}(s,t,k,\sigma,\Gamma),
		\end{aligned}
	\end{equation}
	where
	\begin{equation}
		\begin{aligned} \label{eq:DefF3}
			&\mathcal{G}^{3}_{n}(s_0,t,k,\sigma,\Gamma)\  
			=  \  (-{\bf i}\lambda)^n\sum_{\substack{\rho_i\in\{1,\cdots,n-i+2\},\\ i\in\{0,\cdots,n-1\}}}\sum_{\substack{\bar\sigma\in \{\pm1\}^{\mathcal{I}_n},\\ \sigma_{i,\rho_i}+\sigma_{i-1,\rho_i}+\sigma_{i-1,\rho_i+1}\ne \pm3,
					\\ \sigma_{i-1,\rho_i}\sigma_{i-1,\rho_i+1}= 1}}h^{d}\\
			& \times \int_{(\Lambda^*)^{\mathcal{I}_n}}\mathrm{d}\bar{k}\Delta_{n,\rho}(\bar{k},\bar\sigma) \left\langle\prod_{i=1}^{n+2}\alpha(k_{0,i},\sigma_{0,i})\right\rangle_{s_0}\\
			& \times  \prod_{i=1}^n\Big[\sigma_{i,\rho_i}\mathcal{M}( k_{i,\rho_i}, k_{i-1,\rho_i}, k_{i-1,\rho_i+1}) \Phi_{1,i}(\sigma_{i-1,\rho_i}, k_{i-1,\rho_i}, \sigma_{i-1,\rho_i+1},k_{i-1,\rho_i+1}) \Big] \\
			&  \times\int_{(\mathbb{R}_+)^{\{1,\cdots,n\}}}\mathrm{d}\bar{s} \delta\left(t-\sum_{i=0}^ns_i\right)\prod_{i=1}^{n}e^{-s_i\varsigma_{n-i}}\prod_{i=0}^{n}e^{-s_i \tau_{i}}\\
			&  \times \prod_{i=1}^{n}e^{-{\bf i}t_i(s)\mathfrak{X}(\sigma_{i,\rho_i},k_{i,\rho_i},\sigma_{i-1,\rho_i}, k_{i-1,\rho_i}, \sigma_{i-1,\rho_i+1}, k_{i-1,\rho_i+1})  }\mho(\{(k_{0,1},\sigma_{0,1}),\cdots,(k_{0,n},\sigma_{0,n}),(k_{1,\rho_1},\cdots,k_{n,\rho_n})\}),
	\end{aligned}\end{equation}
in which $\mho(\{(k_{0,1},\sigma_{0,1}),\cdots,(k_{0,n},\sigma_{0,n}),(k_{1,\rho_1},\cdots,k_{n,\rho_n})\})=1$ if the set $\{(k_{0,1},\sigma_{0,1}),\cdots,(k_{0,n},\sigma_{0,n})\}$ is admissible and  all $k_{i,\rho_i}$ are different from $0$; otherwise, the function is $0$.
	We split
	\begin{equation}
		\label{Proposition:Ampl3:2}
		Q^4 \ = \ Q^{4,pair} \ + \ Q^{4,nonpair} 
	\end{equation}
	where 	$Q^{4,pair}$ is the sum over expansions that are admissible  and all $k_{i,\rho_i}$ are different from $0$.
\end{proposition}

\begin{proposition}[Second term]\label{Proposition:Ampl4} The terms associated to $\mathcal{F}_n^1$ can be written as \begin{equation}\label{Def:OperatorQErrors:Bis3}
		\begin{aligned}
			Q^2 = & \sum_{n=0}^{\mathfrak{N}-1}\int_0^t\mathrm{d}s\varsigma_n \mathcal{F}_{n}^1(s,t,k,\sigma,\Gamma)[\alpha_s] 
		\\	= \ & \sum_{n=0}^{\mathfrak{N}-1} \int_0^t\mathrm{d}s\varsigma_n\mathcal{G}^{1,pair}_{n}(s,t,k,\sigma,\Gamma)+\sum_{n=0}^{\mathfrak{N}-1} \int_0^t\mathrm{d}s\varsigma_n\mathcal{G}^{1,nonpair}_{n}(s,t,k,\sigma,\Gamma),
		\end{aligned}
	\end{equation}
	where
	\begin{equation}
		\begin{aligned} \label{eq:DefF3}
			&\mathcal{G}^{1,pair}_{n}(s,t,k,\sigma,\Gamma)\  
			= \varsigma_n\int_0^{t-s}\mathrm{d}r e^{-r\varsigma_n}\mathcal{G}^{3,pair}_{n}(s+r,t,k,\sigma,\Gamma)[\alpha_r].
		\end{aligned}
	\end{equation}
	We split
\begin{equation}
	\label{Proposition:Ampl4:2}
	Q^2 \ = \ Q^{2,pair} \ + \ Q^{2,nonpair} 
\end{equation}
where 	$Q^{2,pair}$ is the sum over expansions that are admissible  and all $k_{i,\rho_i}$ are different from $0$.
\end{proposition}

\subsection{Estimates of the density function}\label{proofs}	
Below, we give a proof of Proposition \ref{Propo:ExpectationEqui}.

\begin{proof}
	The proof of the first claim follows from a straightforward computation
	\begin{equation*}
		\begin{aligned} \left\{\mathcal{H},\exp\left(-\sum_{k\in\Lambda^*}\frac{b_{1,k}^2+b_{2,k}^2}{\bar\omega'(k)}\right)\right\}
			\ = & \ - \exp\left(-\sum_{k\in\Lambda^*}\frac{b_{1,k}^2+b_{2,k}^2}{\bar\omega'(k)}\right)\left\{\mathcal{H},\sum_{k\in\Lambda^*}\frac{b_{1,k}^2+b_{2,k}^2}{\bar\omega'(k)}\right\}
			\ =  \ 0.\end{aligned}
	\end{equation*}
	Moreover, we also have
	$$\left(b_{2,k}\frac{\partial }{\partial b_{1,k}}-b_{1,k}\frac{\partial }{\partial b_{2,k}}\right)\left(-\sum_{k\in\Lambda^*}\frac{b_{1,k}^2+b_{2,k}^2}{\bar\omega'(k)}\right)=0,$$
	leading to
	$\mathbf{R}\tilde{\rho} =  0.$

	 	Let us now consider the trajectories of $c_{1,k}$, $c_{2,k}$,
	 \begin{equation}
	 	\begin{aligned}
	 		\label{Trajectory:1}
	 		\partial_s \mathcal{X}_{1,k}(s,t) \ = & \  \lambda\sum_{k_1,k_2\in\Lambda^*}\mathcal{M}(k,k_1,k_2)\sqrt{2\mathcal{X}_{1,  k_1}2\mathcal{X}_{1, k_2}2\mathcal{X}_{1,k}}\\
	 		&\times\Big[\delta(k- k_1- k_2){\sin(\mathcal{X}_{2,k}- \mathcal{X}_{2,  k_1}- \mathcal{X}_{2, k_2})}\Big], \ \ \ \mathcal{X}_{1,k}(t,t) \ = \ c_{1,k},\\
	 		\partial_s \mathcal{X}_{2,k}(s,t)  \ = & \ \ - \ \lambda\sum_{k_1,k_2\in\Lambda^*}\mathcal{M}(k,k_1,k_2)\sqrt{2\mathcal{X}_{1,  k_1}2\mathcal{X}_{1, k_2}\mathcal{X}_{1,k}^{-1}}\\
	 		&\times\Big[\delta(k-k_1-k_2){\cos(\mathcal{X}_{2,k}- \mathcal{X}_{2,  k_1}- \mathcal{X}_{2, k_2})}\Big]\ - \ \omega_k, \ \ \ \mathcal{X}_{2,k}(t,t) \ = \ c_{2,k}.
	 	\end{aligned}
	 \end{equation}
	 We have
	 \begin{equation}
	 	\begin{aligned}
	 		\label{Trajectory:2}
	 		&\varrho(0,\mathcal{X}_{1}(0,t,c_1,c_2)) \ =  \ \varrho(t,c_1,c_2).
	 	\end{aligned}
	 \end{equation}

 { As in } \eqref{Trajectory:2}, we also have the trajectories of $b_{1,k},b_{2,k}$

 {Similar with \eqref{Trajectory:2}, we also have the trajectories of $b_{1,k},b_{2,k}$
 \begin{equation}
 	\begin{aligned}
 		\label{Trajectory:3}
 		\partial_s \mathcal{B}_{1,k}(s,t) \ = & \  \lambda\sum_{k_1,k_2\in\Lambda^*}\mathcal{M}(k,k_1,k_2)\delta(k-k_1-k_2)(\mathcal B_{1,k_1}\mathcal B_{2,k_2}+\mathcal B_{1,k_2}\mathcal B_{2,k_1})+\mathcal B_{2,k}\omega_k,\\
 		\ \ \ \mathcal{B}_{1,k}(t,t) \ = 	& \ b_{1,k},\\
 		\partial_s \mathcal{B}_{2,k}(s,t)  \ = & \ \ - \ \lambda\sum_{k_1,k_2\in\Lambda^*}\mathcal{M}(k,k_1,k_2)\delta(k-k_1-k_2)(\mathcal B_{1,k_1}\mathcal B_{1,k_2}-\mathcal B_{2,k_1}\mathcal B_{2,k_2})+\mathcal B_{1,k}\omega_k,\\
 		\ \ \ \mathcal{B}_{2,k}(t,t) \ = & \ b_{2,k}.
 	\end{aligned}
 \end{equation}
 We have
 \begin{equation}
 	\begin{aligned}
 		\label{Trajectory:4}
 		&\varrho(0,\mathcal{B}_{1}(0,t,b_1,b_2),\mathcal{B}_{2}(0,t,b_1,b_2)) \ =  \ \varrho(t,b_1,b_2).
 	\end{aligned}
 \end{equation}
 Setting $	\wp_k(s)=[\mathcal{B}_{1,k}(t-s,t) + \mathbf{i}\mathcal{B}_{1,k}(t-s,t)],$ and $	\bar{\wp}_k(s)=[\mathcal{B}_{1,k}(t-s,t) + \mathbf{i}\mathcal{B}_{1,k}(t-s,t)]\sqrt{|\bar\omega(k)|},$ we obtain

 \begin{equation}
 	\label{EquationTrajectoryFouriera}\begin{aligned}
 		\frac{\mathrm{d}\bar{\wp}(k,s)}{\mathrm{d} s} \ & = -\ {\bf i} \omega(k)\bar{\wp}(k,s)  \ - \     {\bf i}\lambda \bar\omega(k)\frac{1}{|\Lambda_*|^2}\sum_{k=k_1+k_2;k_1,k_2\in\Lambda_*}\bar{\wp}(k_1,s)\bar{\wp}(k_2,s),\end{aligned}
 \end{equation}
which finishes the proof of \eqref{Propo:ExpectationEqui:3}.}

Choosing $\mathfrak{P}$, defined in \eqref{Gaussian} as a test function in \eqref{FokkerPlanck2}, we find, by  integrating  by parts multiple times with the notice that $\partial_{c_{2,k}}\hat{\mathfrak{H}}(k)=0$ when $c_{1,k}=0$
\begin{equation} 
	\begin{aligned}
		0 \ = \ 	& \int_{(\mathbb{R}_+\times[-\pi,\pi])^{|\Lambda^*|}}\mathrm{d}c_{1}\mathrm{d}c_{2}\partial_t\varrho\mathfrak{P}\\
		&	+\ \sum_{k\in\Lambda^*}\int_{(\mathbb{R}_+\times[-\pi,\pi])^{|\Lambda^*|}}\mathrm{d}c_{1}\mathrm{d}c_{2}\Big[\Big[\hat{\mathfrak{H}}(k),\varrho\Big]\Big]_k\mathfrak{P}\ 
		- \sum_{k,k'\in\Lambda^*}\mathcal{E}(k,k')\int_{(\mathbb{R}_+\times[-\pi,\pi])^{|\Lambda^*|}}\mathrm{d}c_{1}\mathrm{d}c_{2}\partial_{c_{2,k}c_{2,k'}}\varrho\mathfrak{P}\\
		\ = \ 	& \partial_t\int_{(\mathbb{R}_+\times[-\pi,\pi])^{|\Lambda^*|}}\mathrm{d}c_{1}\mathrm{d}c_{2}\varrho\mathfrak{P}\ 
		+\ \sum_{k\in\Lambda^*}\int_{(\mathbb{R}_+\times[-\pi,\pi])^{|\Lambda^*|}}\mathrm{d}c_{1}\mathrm{d}c_{2}\partial_{c_{2,k}}\hat{\mathfrak{H}}(k)\partial_{c_{1,k}}\varrho\mathfrak{P}
		\\
		&	-\ \sum_{k\in\Lambda^*}\int_{(\mathbb{R}_+\times[-\pi,\pi])^{|\Lambda^*|}}\mathrm{d}c_{1}\mathrm{d}c_{2}\partial_{c_{1,k}}\hat{\mathfrak{H}}(k)\partial_{c_{2,k}}\varrho\mathfrak{P}
		\\
		\ = \ 	& \partial_t\int_{(\mathbb{R}_+\times[-\pi,\pi])^{|\Lambda^*|}}\mathrm{d}c_{1}\mathrm{d}c_{2}\varrho\mathfrak{P}\ 
		-\ \sum_{k\in\Lambda^*}\int_{(\mathbb{R}_+\times[-\pi,\pi])^{|\Lambda^*|}}\mathrm{d}c_{1}\mathrm{d}c_{2}\partial_{c_{2,k}}\hat{\mathfrak{H}}(k)\varrho\partial_{c_{1,k}}\mathfrak{P}
		\\
		& + \sum_{k\in\Lambda^*}\int_{(\mathbb{R}_+)^{|\Lambda^*|-1}\times[-\pi,\pi]^{|\Lambda^*|}}\prod_{k'\in\Lambda^*\backslash\{k\}}\mathrm{d}c_{1,k'}\mathrm{d}c_{2}\partial_{c_{2,k}}\hat{\mathfrak{H}}(k)\varrho\mathfrak{P}\Big|_{c_{1,k}=0}^{c_{1,k}=\infty}\\
		&	-\ \sum_{k\in\Lambda^*}\int_{(\mathbb{R}_+\times[-\pi,\pi])^{|\Lambda^*|}}\mathrm{d}c_{1}\mathrm{d}c_{2}\partial_{c_{1,k}c_{2,k}}\hat{\mathfrak{H}}(k)\varrho\mathfrak{P} 
		\ 	+\ \sum_{k\in\Lambda^*}\int_{(\mathbb{R}_+\times[-\pi,\pi])^{|\Lambda^*|}}\mathrm{d}c_{1}\mathrm{d}c_{2}\partial_{c_{1,k}c_{2,k}}\hat{\mathfrak{H}}(k)\varrho\mathfrak{P}\\
		&
		\ 	-\ \sum_{k\in\Lambda^*}\int_{(\mathbb{R}_+\times[-\pi,\pi])^{|\Lambda^*|}}\mathrm{d}c_{1}\prod_{k'\in\Lambda^*\backslash\{k\}}\mathrm{d}c_{2,k'}\partial_{c_{1,k}}\hat{\mathfrak{H}}(k)\varrho\mathfrak{P}\Big|_{c_{2,k}=-\pi}^{c_{2,k}=\pi}\\
		\ = \ 	& \partial_t\int_{(\mathbb{R}_+\times[-\pi,\pi])^{|\Lambda^*|}}\mathrm{d}c_{1}\mathrm{d}c_{2}\varrho\mathfrak{P}\ 
		-\ \sum_{k\in\Lambda^*}\int_{(\mathbb{R}_+\times[-\pi,\pi])^{|\Lambda^*|}}\mathrm{d}c_{1}\mathrm{d}c_{2}\partial_{c_{2,k}}\hat{\mathfrak{H}}(k)\varrho\partial_{c_{1,k}}\mathfrak{P}
		\\
		&	-\ \sum_{k\in\Lambda^*}\int_{(\mathbb{R}_+\times[-\pi,\pi])^{|\Lambda^*|}}\mathrm{d}c_{1}\mathrm{d}c_{2}\partial_{c_{1,k}c_{2,k}}\hat{\mathfrak{H}}(k)\varrho\mathfrak{P} 
		\ 	+\ \sum_{k\in\Lambda^*}\int_{(\mathbb{R}_+\times[-\pi,\pi])^{|\Lambda^*|}}\mathrm{d}c_{1}\mathrm{d}c_{2}\partial_{c_{1,k}c_{2,k}}\hat{\mathfrak{H}}(k)\varrho\mathfrak{P},\end{aligned}
\end{equation}
which, after the simplification of the last two similar terms, leads to
\begin{equation}
	\begin{aligned}
		0=\	& \partial_t\int_{(\mathbb{R}_+\times[-\pi,\pi])^{|\Lambda^*|}}\mathrm{d}c_{1}\mathrm{d}c_{2}\varrho\mathfrak{P}\
		-\ \sum_{k\in\Lambda^*}\int_{(\mathbb{R}_+\times[-\pi,\pi])^{|\Lambda^*|}}\mathrm{d}c_{1}\mathrm{d}c_{2}\partial_{c_{2,k}}\hat{\mathfrak{H}}(k)\varrho\partial_{c_{1,k}}\mathfrak{P}.
	\end{aligned}
\end{equation}
Let us study the second term on the right hand side of the above equation using \eqref{DensityHamiltonian}:  
\begin{equation}
	\begin{aligned}
		&	 \sum_{k\in\Lambda^*}\int_{(\mathbb{R}_+\times[-\pi,\pi])^{|\Lambda^*|}}\mathrm{d}c_{1}\mathrm{d}c_{2}\partial_{c_{2,k}}\hat{\mathfrak{H}}(k)\varrho\partial_{c_{1,k}}\mathfrak{P}\\
		= \ &\lambda \int_{(\mathbb{R}_+\times[-\pi,\pi])^{|\Lambda^*|}}\mathrm{d}c_{1}\mathrm{d}c_{2}\varrho\left\{\sum_{k\in\Lambda^*}\partial_{c_{1,k}}\mathfrak{P} \sum_{k_1',k_2'\in\Lambda^*}\mathcal{M}(k,k_1',k_2')\right.\\
		&\times 2 \sqrt{2c_{1,k_1'}c_{1,k_2'}c_{1,k}}\Big[\delta(k-k_1'-k_2'){\sin(c_{2,k_1'}+c_{2,k_2'}-c_{2,k})}\Big]\Big\}
		\\ 
		= \ & \sum_{k\in\Lambda^*} c_{\mathfrak{P}}\lambda  h^{2d}\tilde\omega(k)\int_{(\mathbb{R}_+\times[-\pi,\pi])^{|\Lambda^*|}}\mathrm{d}c_{1}\mathrm{d}c_{2}\varrho\left\{\mathfrak{P}  \sum_{k_1',k_2'\in\Lambda^*}\mathcal{M}(k,k_1',k_2')\right.\\
		&\times 2 \sqrt{2c_{1,k_1'}c_{1,k_2'}c_{1,k}}\Big[\delta(k-k_1'-k_2'){\sin(c_{2,k_1'}+c_{2,k_2'}-c_{2,k})}\Big]\Big\}\\
		= \ & \sum_{k\in\Lambda^*} h^{2d}c_{\mathfrak{P}}\lambda  \tilde\omega(k)\int_{(\mathbb{R}_+\times[-\pi,\pi])^{|\Lambda^*|}}\mathrm{d}c_{1}\mathrm{d}c_{2}\varrho\mathfrak{P}\Big\{ \sum_{k_1',k_2'\in\Lambda^*}|\mathcal{M}(k,k_1',k_2')|\mathrm{sign}k^1\\
		&\times 2 \sqrt{2c_{1,k_1'}c_{1,k_2'}c_{1,k}}\Big[\delta(k-k_1'-k_2'){\sin(c_{2,k_1'}+c_{2,k_2'}-c_{2,k})}\Big]\Big\}.
	\end{aligned}
\end{equation}
 By a rotation of $k,k_1',k_2'$, we have
\begin{equation}
	\begin{aligned}
		& \sum_{k,k_1',k_2'\in\Lambda^*}|\mathcal{M}(k,k_1',k_2')|\mathrm{sign}k^1 2\sqrt{2c_{1,k_1'}c_{1,k_2'}c_{1,k}}\Big[\delta(k-k_1'-k_2'){\sin(c_{2,k_1'}+c_{2,k_2'}-c_{2,k})}\Big]\tilde\omega(k)\\
		=\ &\mathfrak{M}_1\ + \ \mathfrak{M}_2 \ + \ \mathfrak{M}_3, 
	\end{aligned}
\end{equation}
in which
\begin{equation}
	\begin{aligned}
		\mathfrak{M}_1: =\ 	&	\sum_{k,k_1',k_2'\in\Lambda^*}|\mathcal{M}(k,k_1',k_2')|\mathrm{sign}k^1 \sqrt{2c_{1,k_1'}c_{1,k_2'}c_{1,k}}\Big[\delta(k-k_1'-k_2'){\sin(c_{2,k_1'}+c_{2,k_2'}-c_{2,k})}\Big]\tilde\omega(k),\\
		\mathfrak{M}_2 :=\ 
		& \sum_{k,k_1',k_2'\in\Lambda^*}|\mathcal{M}(k,k_1',k_2')|\mathrm{sign}(k_1')^1 \sqrt{2c_{1,k_1'}c_{1,k_2'}c_{1,k}}\Big[\delta(k_1'-k-k_2'){\sin(-c_{2,k_1'}+c_{2,k_2'}+c_{2,k})}\Big]{\tilde\omega(k'_1)},\\
		\mathfrak{M}_3 :=\  & \sum_{k,k_1',k_2'\in\Lambda^*}|\mathcal{M}(k,k_1',k_2')|\mathrm{sign}(k_2')^1 \sqrt{2c_{1,k_1'}c_{1,k_2'}c_{1,k}}\Big[\delta(k_2'-k_1'-k){\sin(c_{2,k_1'}-c_{2,k_2'}+c_{2,k})}\Big]{\tilde\omega(k'_2)}.
	\end{aligned}
\end{equation}
We observe that  $\mathfrak{M}_1+\mathfrak{M}_2+\mathfrak{M}_3=0$. Hence,
we finally obtain 
\begin{equation}
	\begin{aligned}
		\partial_t\int_{(\mathbb{R}_+\times[-\pi,\pi])^{|\Lambda^*|}}\mathrm{d}c_{1}\mathrm{d}c_{2}\varrho\mathfrak{P} \
		= \ & 0,
	\end{aligned}
\end{equation}
which finishes the proof of \eqref{Propo:ExpectationEqui:1}.

To prove \eqref{Propo:ExpectationEqui:2}, we compute
\begin{equation}\label{Propo:ExpectationEqui:EA}
	\begin{aligned}
		&\sum_{k_1,\cdots,k_n\in\Lambda_*}\left|\Big\langle |a_{k_{1}}|^2\cdots |a_{k_{n}}|^2\Big\rangle_t\right| \\
		= & \  \left|\int_{\mathbb{R}^{2n}}\prod_{i=1}^{n}\mathrm{d}b_{1,k_{i}}\mathrm{d}b_{2,k_{i}}\sum_{k_1,\cdots,k_n\in\Lambda_*}\big|(b_{1,k_{i}}+{\bf i}b_{2,k_{i}})\big|^2 \int_{\mathbb{R}^{2(|\Lambda^*|-n)}}\sum_{k_1,\cdots,k_n\in\Lambda_*}\prod_{k\in\Lambda^*\backslash\{k_{1},\cdots,k_{n}\}}\mathrm{d}b_{1,k}\mathrm{d}b_{2,k}\varrho(t)\right|\\	\lesssim  &\ \left|	\int_{(\mathbb{R}_+\times[-\pi,\pi])^{|\Lambda^*|}}\mathrm{d}c_{1}\mathrm{d}c_{2}\prod_{j=1}^{n}\sum_{k_1,\cdots,k_n\in\Lambda_*}\big|2c_{1,k_{i}}\big|\varrho(t)\right|\\\
		\
		\lesssim   &  \  h^{-2dn}	C^{{n}}[\tilde\omega(k_1)\cdots\tilde\omega(k_n)]^{-1} \left|\int_{(\mathbb{R}_+\times[-\pi,\pi])^{|\Lambda^*|}}\mathrm{d}c_{1}\mathrm{d}c_{2}\mathfrak{P}\varrho(t)\right|\
		\\
		\lesssim    &  \  h^{-2nd}	C^{{n}} [\tilde\omega(k_1)\cdots\tilde\omega(k_n)]^{-1} \left|\int_{(\mathbb{R}_+\times[-\pi,\pi])^{|\Lambda^*|}}\mathrm{d}c_{1}\mathrm{d}c_{2}\mathfrak{P}\varrho(0)\right|,
	\end{aligned}
\end{equation}
where we have used \eqref{Propo:ExpectationEqui:1}.

	Inserting the form of $\varrho(0)$ (see Definition \ref{def:distinct})  into \eqref{Propo:ExpectationEqui:EA}, we then find
\begin{equation}\label{Propo:ExampleMeasure:E3}
	\begin{aligned}
		&\sum_{k_1,\cdots,k_n\in\Lambda_*}\Big\langle |a_{k_{1}}|^2\cdots |a_{k_{n}}|^2\Big\rangle_t\\
		\lesssim &\ h^{-2dn}C^{{n}}[\tilde\omega(k_1)\cdots\tilde\omega(k_n)]^{-1}   \left[(2\pi)^{|\Lambda^*|}\int_{\mathbb{R}_+^{|\Lambda^*|}}\mathrm{d}c_{1}\left(\prod_{k\in\Lambda^*}\frac{e^{ \tilde\omega(k)h^{2d}c_{\mathfrak{P}}c_{1,k}-\frac{2c_{1,k}}{\gimel(k)}}}{\pi\gimel(k) }\right)\right].
	\end{aligned}
\end{equation}
Since 
\begin{equation}\label{IdentityExponential}
	\int_0^\infty\mathrm{d}x x^ne^{-ax} \ = \ \frac{n!}{a^{n+1}},
\end{equation} the last term \eqref{Propo:ExampleMeasure:E3} can be computed explicitly
as
\begin{equation*}\begin{aligned}
		& (2\pi)^{|\Lambda^*|}\int_{\mathbb{R}_+^{|\Lambda^*|}}\mathrm{d}c_{1}\left(\prod_{k\in\Lambda^*}\frac{e^{ h^{2d}c_{\mathfrak{P}}\tilde\omega(k)c_{1,k}-\frac{2c_{1,k}}{\gimel(k)}}}{\pi\gimel(k) }\right)\\
		= &\ \int_{\mathbb{R}_+^{|\Lambda^*|}}\prod_{k\in\Lambda^*}\mathrm{d}(2c_{1,k})\prod_{k\in\Lambda^*}\left(\frac{e^{2c_{1,k}\big(\frac{\tilde\omega(k)h^{2d}c_{\mathfrak{P}}}{2}-\frac{1}{\gimel(k)}\big)}}{\big(-\frac{h^{2d}c_{\mathfrak{P}\tilde\omega(k)}}{2}+\frac{1}{\gimel(k)}\big)^{-1}}\frac{1}{\gimel(k)\big(-\frac{h^{2d}c_{\mathfrak{P}}\tilde\omega(k)}{2}+\frac{1}{\gimel(k)}\big)}\right)\\
 		\lesssim & \ \prod_{k\in\Lambda^*}\frac{2}{2-c_{\mathfrak{P}}h^{2d}\tilde\omega(k)\gimel(k)}.
	\end{aligned}
\end{equation*}

This quantity can be bounded as (see Definition \ref{def:distinct})
\begin{equation}\label{Lemma:BoundL1Identity:E2}
	\begin{aligned}
\prod_{k\in\Lambda^*}\frac{2}{2-c_{\mathfrak{P}}h^{2d}\tilde\omega(k)\gimel(k)}\ \le	 & \ 	\left|\frac{2}{2-c_{\mathfrak{P},1}h^d}\right|^{|\Lambda^*|}\ 
		\le  	\  	\left|\frac{2}{2-c_{\mathfrak{P},1}|\Lambda^*|^{-1}}\right|^{|\Lambda^*|}\ \le \ c_{\mathfrak{P}}', 
	\end{aligned}
\end{equation}
for some constant $c_{\mathfrak{P}}'>0$.
Inequality  \eqref{Propo:ExpectationEqui:2} then follows.
\end{proof}

\section{Dispersive estimates}
\subsection{Dispersive estimates}\label{Subsec:DispersiveEstimates}

Let $V=(V_1,\cdots,V_d),W=(W_1,\cdots,W_d)$ be two vectors in $[-\pi,\pi]^d$ and $t_0,t_1,t_2$ be real numbers. We set
\begin{equation}\label{Lemm:Bessel}
	\begin{aligned}
		& {t_0}+{t_1}\cos(V_1)e^{{\bf i}2V_j}+{t_2}\cos(W_1)e^{{\bf i}2W_j} \
		= \ \Big|{t_0}+{t_1}\cos(V_1)e^{{\bf i}2V_j}+{t_2}\cos(W_1)e^{{\bf i}2W_j}\Big|e^{{\bf i}\aleph
			_1^j},\\
		& \mbox{ and }\\
		& {t_1}\sin(V_1)e^{{\bf i}2V_j}+{t_2}\sin(W_1)e^{{\bf i}2W_j} 
		= \  \ \Big|{t_1}\sin(V_1)e^{{\bf i}2V_j}+{t_2}\sin(W_1)e^{{\bf i}2W_j}\Big|e^{{\bf i}\aleph_2^j},
	\end{aligned}
\end{equation}
with $\aleph_1^j=\aleph_1^j(V,W),\aleph_2^j=\aleph_2^j(V,W)\in [-\pi,\pi]$, $j=2,\cdots,d$.

The following lemma gives an estimate on the angles $\aleph_1^j,\aleph_2^j$.
\begin{lemma}\label{Lemm:Angle}
	Suppose that \begin{equation}
		\label{Lemm:Angle:0} t_0=t_1\aleph' + t_2 \aleph''\end{equation} with $\aleph',\aleph''\in\{\pm 1\}$. 
	The following estimate then holds true
	\begin{equation}
		\label{Lemm:Angle:00}\begin{aligned}
			\frac{1}{|1-|\cos(\aleph_1^j-\aleph_2^j)||^\frac12} 
			\ \lesssim \ & \frac{t_1^2+t_2^2}{|\mathfrak{C}_{\aleph_1^j,\aleph_2^j}^1t_1^2+\mathfrak{C}_{\aleph_1^j,\aleph_2^j}^2t_1t_2+\mathfrak{C}_{\aleph_1^j,\aleph_2^j}^3t_2^2|}+1,\end{aligned}
	\end{equation}
	where the constants on the right hand side are universal and
	\begin{equation}
		\label{Lemm:Angle:1}\begin{aligned}
			\mathfrak{C}_{\aleph_1^j,\aleph_2^j}^1 \ = \ & -\aleph'\sin(V_1)\sin(2V_j),
			\\
			\mathfrak{C}_{\aleph_1^j,\aleph_2^j}^3 \ = \ &  -\aleph''\sin(W_1)\cos(2W_j),\\
			\mathfrak{C}_{\aleph_1^j,\aleph_2^j}^2 \ = \  & -\sin(2V_j-2W_j)\sin(V_1-W_1)  \ - \aleph'\sin(W_1)\sin(2W_j) \\
			& - \ \aleph''\sin(V_1)\sin(2V_j).
		\end{aligned}
	\end{equation}
	Setting $r_*=t_1/t_2$, we consider the equation
	\begin{equation}
		\label{Lemm:Angle:2}\begin{aligned}
			\mathfrak{C}_{\aleph_1^j,\aleph_2^j}^1r_*^2+\mathfrak{C}_{\aleph_1^j,\aleph_2^j}^2r_*+\mathfrak{C}_{\aleph_1^j,\aleph_2^j}^3 \ = \ 0.
		\end{aligned}
	\end{equation}
	In the case that \eqref{Lemm:Angle:2} has two real roots (or a double real root), we denote these roots by $\tilde{r}_1,\tilde{r}_2$. In the case that \eqref{Lemm:Angle:2} has two complex roots, we denote the real part of the complex root by $\tilde{r}_3$. Let $\epsilon_{r_*}\in(-1,1)$ be a small number. We consider the case when $r_*=(1+\epsilon_{r_*})\tilde r_i$, $i=1,2,3$ and define the function

	\begin{equation}
		\label{Lemm:Angle:3}\begin{aligned}
			& g_{t_1,t_2}(r_*) \ = \ \\
			&  \frac{\Big|{r_*}\sin(V_1)e^{{\bf i}2V_j}+ \sin(W_1)e^{{\bf i}2W_j}\Big|}{\Big[\Big|r_*\Big(\cos(V_1)e^{{\bf i}2V_j}+\aleph'\Big)+ \Big(\cos(W_1)e^{{\bf i}2W_j}+\aleph'\Big)\Big|^2+\Big|{r_*}\sin(V_1)e^{{\bf i}2V_j}+ \sin(W_1)e^{{\bf i}2W_j}\Big|^2\Big]^\frac12},
		\end{aligned}
	\end{equation}
	and
	\begin{equation}
		\label{Lemm:Angle:3a}\begin{aligned}
			& f_{t_1,t_2}(r_*) \ = \ \\
			&  \frac{\Big|r_*\Big(\cos(V_1)e^{{\bf i}2V_j}+\aleph'\Big)+ \Big(\cos(W_1)e^{{\bf i}2W_j}+\aleph''\Big)\Big|}{\Big[\Big|r_*\Big(\cos(V_1)e^{{\bf i}2V_j}+\aleph'\Big)+ \Big(\cos(W_1)e^{{\bf i}2W_j}+\aleph''\Big)\Big|^2+\Big|{r_*}\sin(V_1)e^{{\bf i}2V_j}+ \sin(W_1)e^{{\bf i}2W_j}\Big|^2\Big]^\frac12}.
		\end{aligned}
	\end{equation}
	We suppose that 
	\begin{equation}
		\label{Lemm:Angle:4}
		|V_j|, \Big|V_j-\frac{\pi}{2}\Big|, |V_j-{\pi}|,|W_j|, \Big|W_j-\frac{\pi}{2}\Big|, |W_j-{\pi}|\ge \langle\ln\lambda\rangle^{-c_{V,W}},
	\end{equation}
	and
	\begin{equation}
		\label{Lemm:Angle:5}
		|V_j\pm W_j| \ge \langle\ln\lambda\rangle^{-c_{V,W}},
	\end{equation}
	for some constant $c_{V,W}>0$ and for all $j=1,\cdots,d$. We thus have
	\begin{equation}
		\label{Lemm:Angle:6}\begin{aligned}
			&\Big|\frac{\mathrm{d}}{\mathrm{d} r_*}	|f_{t_1,t_2}(r_*)|^2\Big|	 \ \lesssim \langle\ln\lambda\rangle^{\mathfrak{C}_{\aleph_1^j,\aleph_2^j}^4},
		\end{aligned}
	\end{equation} 
	\begin{equation}
		\label{Lemm:Angle:7}\begin{aligned}
			&\Big|\frac{\mathrm{d}}{\mathrm{d} r_*}	|g_{t_1,t_2}(r_*)|^2\Big|	 \ \lesssim \langle\ln\lambda\rangle^{\mathfrak{C}_{\aleph_1^j,\aleph_2^j}^4},
		\end{aligned}
	\end{equation} 
	and
	\begin{equation}
		\label{Lemm:Angle:8}\begin{aligned}
			&\Big|\frac{\mathrm{d}}{\mathrm{d} r_*}	[f_{t_1,t_2}(r_*)g_{t_1,t_2}(r_*)]\Big|	 \ \lesssim \langle\ln\lambda\rangle^{\mathfrak{C}_{\aleph_1^j,\aleph_2^j}^4},
		\end{aligned}
	\end{equation} 
	for some explicit constant $\mathfrak{C}_{\aleph_1^j,\aleph_2^j}^4>0$.
	
	In the special case that $t_2=0,$ $t_1=-t_0$, we bound\begin{equation}
		\label{Lemm:Angle:9}
		[1-|\cos(\aleph_1^j-\aleph_2^j)|]^{-\frac12}\ 
		\lesssim  \ \frac{1}{|\sin(V_1)\sin(2V_j)|}+1.
	\end{equation}

\end{lemma}
\begin{proof}
	
	We compute
	\begin{equation}
		\label{Lemm:Angle:E1}
		\tan(\aleph_1^j) \ = \  \frac{t_1\cos(V_1)\sin(2V_j)+t_2\cos(W_1)\sin(2W_j)}{t_1[\cos(V_1)\cos(2V_j)+\aleph']+t_2[\cos(W_1)\cos(2W_j)+\aleph'']}\ =: \ \frac{At_1+Bt_2}{Ct_1+Dt_2},
	\end{equation}
	and
	\begin{equation}
		\label{Lemm:Angle:E2}
		\tan(\aleph_2^j) \ = \  \frac{t_1\sin(V_1)\sin(2V_j)+t_2\sin(W_1)\sin(2W_j)}{t_1\sin(V_1)\cos(2V_j)+t_2\sin(W_1)\cos(2W_j)}\ =: \ \frac{A't_1+B't_2}{C't_1+D't_2},
	\end{equation}
	then
	\begin{equation}
		\label{Lemm:Angle:E3}\begin{aligned}
			& \tan(\aleph_1^j) - \tan(\aleph_2^j) \ = \   \frac{At_1+Bt_2}{Ct_1+Dt_2} - \frac{A't_1+B't_2}{C't_1+D't_2}\\
			\ = \ & \frac{(AC'-CA')t_1^2+(BC'-CB'+AD'-A'D)t_1t_2+(BD'-B'D)t_2^2}{(Ct_1+Dt_2)(C't_1+D't_2)}.
		\end{aligned}
	\end{equation}

	We then bound
	\begin{equation}
		\label{Lemm:Angle:E4}\begin{aligned}
			& \frac{1}{|\tan(\aleph_1^j) - \tan(\aleph_2^j)|}\\ \ = \ &  
			\Big|\frac{(Ct_1+Dt_2)(C't_1+D't_2)}{(AC'-CA')t_1^2+(BC'-CB'+AD'-A'D)t_1t_2+(BD'-B'D)t_2^2}\Big|\\
			\ \le \ &  
			\frac{5(|t_1|+|t_2|)^2}{|(AC'-CA')t_1^2+(BC'-CB'+AD'-A'D)t_1t_2+(BD'-B'D)t_2^2|} \ =: \ E.
		\end{aligned}
	\end{equation}
	Similarly, we also compute
	\begin{equation}
		\label{Lemm:Angle:E5}\begin{aligned}
			& 
			\cot(\aleph_1^j) - \cot(\aleph_2^j) \ = \   \frac{Ct_1+Dt_2}{At_1+Bt_2} - \frac{C't_1+D't_2}{A't_1+B't_2}\\
			\ = \ & \frac{-(AC'-CA')t_1^2-(BC'-CB'+AD'-A'D)t_1t_2-(BD'-B'D)t_2^2}{(At_1+Bt_2)(A't_1+B't_2)},
		\end{aligned}
	\end{equation}
	and bound
	\begin{equation}
		\label{Lemm:Angle:E6}\begin{aligned}
			& \frac{1}{|\cot(\aleph_1^j) - \cot(\aleph_2^j)|}
			\ \le \ 
			\frac{5(|t_1|+|t_2|)^2}{|(AC'-CA')t_1^2+(BC'-CB'+AD'-A'D)t_1t_2+(BD'-B'D)t_2^2|} = {E}.
		\end{aligned}
	\end{equation}
	Combining \eqref{Lemm:Angle:E4} and \eqref{Lemm:Angle:E6}, we find
	\begin{equation}
		\label{Lemm:Angle:E7}
		|\sin(\aleph_1^j-\aleph_2^j)|\ \ge \ \frac{|\cos(\aleph_1^j)\cos(\aleph_2^j)|}{E}, \ \ \ \mbox{ and    } |\sin(\aleph_1^j-\aleph_2^j)|\ \ge \ \frac{|\sin(\aleph_1^j)\sin(\aleph_2^j)|}{E},
	\end{equation}
	leading to
	\begin{equation}
		\label{Lemm:Angle:E8}
		|\sin(\aleph_1^j-\aleph_2^j)|\ \ge \ \frac{|\cos(\aleph_1^j)\cos(\aleph_2^j)+\sin(\aleph_1^j)\sin(\aleph_2^j)|}{2E}\ = \ \frac{|\cos(\aleph_1^j-\aleph_2^j)|}{2E}.
	\end{equation}
	We then find $
	\label{Lemm:Angle:E9}
	|\cos(\aleph_1^j-\aleph_2^j)|^2\Big(1+\frac{1}{4E^2}\Big)\ \le  \  1,$
	yielding
	$
	\label{Lemm:Angle:E10}
	|\cos(\aleph_1^j-\aleph_2^j)|\ \le  \  \frac{2E}{\sqrt{4E^2+1}},$
	then
	\begin{equation}
		\label{Lemm:Angle:E10:a}
		1-|\cos(\aleph_1^j-\aleph_2^j)|\ \ge \ 1- \frac{2E}{\sqrt{4E^2+1}}=\frac{1}{\sqrt{4E^2+1}[\sqrt{4E^2+1}+2E]}\ \ge \ \frac{1}{{4E^2+1}}.
	\end{equation}
	Finally, we find
	\begin{equation}
		\label{Lemm:Angle:E11}\begin{aligned}
			& \frac{1}{|1-|\cos(\aleph_1^j-\aleph_2^j)||^\frac12} \\
			\ \lesssim \ & \sqrt{\frac{4(|t_1|+|t_2|)^4}{|(AC'-CA')t_1^2+(BC'-CB'+AD'-A'D)t_1t_2+(BD'-B'D)t_2^2|^2}+1}\\
			\ \lesssim \ & \frac{|t_1|^2+|t_2|^2}{|(AC'-CA')t_1^2+(BC'-CB'+AD'-A'D)t_1t_2+(BD'-B'D)t_2^2|}+1.\end{aligned}
	\end{equation}
	We also observe that
	\begin{equation}
		\label{Lemm:Angle:E12}\begin{aligned}
			& AC'-CA'\ = \  \cos(V_1)\sin(2V_j)\sin(V_1)\cos(2V_j)\\ 
			& \ - \ [\cos(V_1)\cos(2V_j)+\aleph']\sin(V_1)\sin(2V_j) \ = \ -\aleph'\sin(V_1)\sin(2V_j),
			\\
			& BD'-B'D \ = \  \cos(W_1)\sin(2W_j)\sin(W_1)\cos(2W_j) \\
			&\ - \ [\cos(W_1)\cos(2W_j)+\aleph'']\sin(W_1)\sin(2W_j)  \ = \  -\aleph''\sin(W_1)\sin(2W_j)\\
			&  BC'-CB'+AD'-A'D\ = \  \cos(W_1)\sin(2W_j)\sin(V_1)\cos(2V_j)\\
			& \ - \  [\cos(V_1)\cos(2V_j)+\aleph']\sin(W_1)\sin(2W_j) + \cos(V_1)\sin(2V_j)\sin(W_1)\cos(2W_j) \\
			& \ - \ \sin(V_1)\sin(2V_j)[\cos(W_1)\cos(2W_j)+\aleph'']\\
			& =\ \sin(2W_j)\cos(2V_j)\sin(V_1-W_1) \ - \ \cos(2W_j)\sin(2V_j)\sin(V_1-W_1)\\
			& \ \ - \aleph'\sin(W_1)\sin(2W_j) \ - \ \aleph''\sin(V_1)\sin(2V_j)\\
			& =\ -\sin(2V_j-2W_j)\sin(V_1-W_1)  \ - \aleph'\sin(W_1)\sin(2W_j) \ - \ \aleph''\sin(V_1)\sin(2V_j).
		\end{aligned}
	\end{equation}
	Combining 
	\eqref{Lemm:Angle:E11} and \eqref{Lemm:Angle:E12} yields the first conclusion of the lemma.
	
	Now, by definition, we compute
	\begin{equation}
		\label{Lemm:Angle:E13}\begin{aligned}
			&	|f_{t_1,t_2}(r_*)|^2	 \ = \  \left|\frac{|Ar_*+B+{\bf i } (Cr_*+D)|}{[|r_*A+B+ {\bf i}(Cr_*+D)|^2+|r_*A'+B'+ {\bf i}(C'r_*+D')|^2]^\frac12}\right|^2\\
			& \ = \ \frac{(A^2+C^2)r_*^2+2(AB+CD)r_*+B^2+D^2}{(A^2+A'^2+C^2+C'^2)r_*^2+2(AB+CD+A'B'+C'D')r_*+B^2+B'^2+D^2+D'^2}\\
			& \ = \ : \frac{M_1r_*^2+M_2r_*+M_3}{N_1r_*^2+N_2r_*+N_3},
		\end{aligned}
	\end{equation}
	which implies
	\begin{equation}
		\label{Lemm:Angle:E14}\begin{aligned}
			&\frac{\mathrm{d}}{\mathrm{d} r_*}	|f_{t_1,t_2}(r_*)|^2	 \ = \ \frac{(M_1N_2-M_2N_1)r_*^2+2(M_1N_3-N_1M_3)r_*+M_2N_3-N_2M_3}{|N_1r_*^2+N_2r_*+N_3|^2}.
		\end{aligned}
	\end{equation}
	
	Let us consider the case when $r_*=(1+\epsilon_{r_*})\tilde{r}_1$. We estimate using \eqref{Lemm:Angle:4}-\eqref{Lemm:Angle:5} 
	\begin{equation}
		\label{Lemm:Angle:E15}\begin{aligned}
			&|(M_1N_2-M_2N_1)r_*^2+2(M_1N_3-N_1M_3)r_*+M_2N_3-N_2M_3 | \lesssim |r_*^2+1|\\
			&	\lesssim\ \left|\frac{\mathfrak{C}_{\aleph_1^j,\aleph_2^j}^2+\sqrt{|\mathfrak{C}_{\aleph_1^j,\aleph_2^j}^2|^2-4\mathfrak{C}_{\aleph_1^j,\aleph_2^j}^1\mathfrak{C}_{\aleph_1^j,\aleph_2^j}^3}}{2\mathfrak{C}_{\aleph_1^j,\aleph_2^j}^1}\right|^2+1\ \lesssim\ \frac{1}{|\mathfrak{C}_{\aleph_1^j,\aleph_2^j}^1|^2} \lesssim \langle\ln\lambda\rangle^{\mathfrak{C}_{r_*}^1},
		\end{aligned}
	\end{equation}
	for some constant $\mathfrak{C}_{r_*}^1>0$. Similarly, when $r_*=(1+\epsilon_{r_*})\tilde{r}_2$ and $r_*=(1+\epsilon_{r_*})\tilde{r}_3$, the same estimate \eqref{Lemm:Angle:E15} also holds true. 
	
	Next, we estimate the denominator of \eqref{Lemm:Angle:E14}
	\begin{equation}
		\label{Lemm:Angle:E16}\begin{aligned}
			&N_1r_*^2+N_2r_*+N_3\  = \ \\
			= & \ (A^2+A'^2+C^2+C'^2)r_*^2+2(AB+CD+A'B'+C'D')r_*+B^2+B'^2+D^2+D'^2\\
			= & \ (A^2+A'^2+C^2+C'^2)\Big[r_*+\frac{AB+CD+A'B'+C'D'}{A^2+A'^2+C^2+C'^2}\Big]^2 \\
			& \  \  \  +\frac{(A^2+A'^2+C^2+C'^2)(B^2+B'^2+D^2+D'^2)-(AB+CD+A'B'+C'D')^2}{|A^2+A'^2+C^2+C'^2|}\\
			\ge & \ \frac{(A^2+A'^2+C^2+C'^2)(B^2+B'^2+D^2+D'^2)-(AB+CD+A'B'+C'D')^2}{|A^2+A'^2+C^2+C'^2|}.
		\end{aligned}
	\end{equation}
	We compute
	\begin{equation}
		\label{Lemm:Angle:E17}\begin{aligned}
			& A^2+A'^2+C^2+C'^2 \  = \ |\cos(V_1)\sin(2V_j)|^2 + |\sin(V_1)\cos(2V_j)|^2 \\
			&  +  |\cos(V_1)\cos(2V_j)+\aleph'|^2 +|\sin(V_1)\sin(2V_j)|^2 \lesssim 1,
		\end{aligned}
	\end{equation}
	which implies
	\begin{equation}
		\label{Lemm:Angle:E18}\begin{aligned}
			&N_1r_*^2+N_2r_*+N_3\  \gtrsim \ \\
			\gtrsim  & \ (A^2+A'^2+C^2+C'^2)(B^2+B'^2+D^2+D'^2)-(AB+CD+A'B'+C'D')^2.
		\end{aligned}
	\end{equation}
	Using the identity
	\begin{equation}\label{Lemm:Angle:E19}\begin{aligned}
			\left(\sum_{i=1}^4a_i^2\right)\left(\sum_{i=1}^4b_i^2\right)-\left(\sum_{i=1}^4a_ib_i\right)^2
			\ = \ & \frac12\left[\sum_{i,j=1}^4(a_ib_j-a_jb_i)^2\right],
		\end{aligned}
	\end{equation}
	for all $a_i,b_i\in\mathbb{R}$, $i=1,2,3,4$, we find by \eqref{Lemm:Angle:4}-\eqref{Lemm:Angle:5}
	\begin{equation}
		\label{Lemm:Angle:E20}\begin{aligned}
			N_1r_*^2+N_2r_*+N_3\  \gtrsim  &\ 
			[\cos(V_1)\sin(2V_j)\sin(W_1)\sin(2W_j) 	- \sin(V_1)\sin(2V_j)\cos(W_1)\sin(2W_j)]^2	\\
			\  \gtrsim  &\ 
			[\sin(V_1-W_1)\sin(2V_j)\sin(2W_j)]^2	\
			\  \gtrsim  \	\langle\ln\lambda\rangle^{-\mathfrak{C}_{r_*}^2},
		\end{aligned}
	\end{equation}
	for some constant $\mathfrak{C}_{r_*}^2>0$. Therefore
	\begin{equation}
		\label{Lemm:Angle:E21}\begin{aligned}
			&\Big|\frac{\mathrm{d}}{\mathrm{d} r_*}	|f_{t_1,t_2}(r_*)|^2\Big|	 \ \lesssim \langle\ln\lambda\rangle^{\mathfrak{C}_{\aleph_1^j,\aleph_2^j}^4},
		\end{aligned}
	\end{equation} for some constant $\mathfrak{C}_{\aleph_1^j,\aleph_2^j}^4>0$. Thus \eqref{Lemm:Angle:6} is proved. The  inequalities \eqref{Lemm:Angle:7} and  \eqref{Lemm:Angle:8} can be proved by similar arguments. 
	
	In the special case that $t_2=0,$ $t_1=-t_0$, we have 
	{$$
		\frac{t_1\cos(V_1)\sin(2V_j)}{t_1[\cos(V_1)\cos(2V_j)-1]} = \frac{At_1+Bt_2}{Ct_1+Dt_2},\quad 
		\frac{t_1\sin(V_1)\sin(2V_j)}{t_1\sin(V_1)\cos(2V_j)} =  \frac{A't_1+B't_2}{C't_1+D't_2},
		$$}
	then $A=\cos(V_1)\sin(2V_j),$ $C = \cos(V_1)\cos(2V_j)-1,$ {$A'= \sin(V_1)\sin(2V_j),$} $C'=\sin(V_1)\cos(2V_j). $  
	Thus $AC'-CA'=\sin(2V_j)$, and we bound, following \eqref{Lemm:Angle:E11}
	\begin{equation}
		\label{Lemm:Angle:E23}
		[1-|\cos(\aleph_1^j-\aleph_2^j)|]^{-\frac12}\ 
		\lesssim \ \frac{1}{|\sin(V_1)\sin(2V_j)|}+1.
	\end{equation}
	
\end{proof}
We set

\begin{equation*}
	\begin{aligned}\label{Lemm:Bessel2:2bb8:1}
		\cos(\Upsilon_*(\tilde{r}_i,\aleph^*,\aleph_*))
		\ = \ & {\Big|\aleph^*\tilde{r}_i+\aleph_*+\tilde{r}_i\cos(V_1)e^{{\bf i}2V_j}+\cos(W_1)e^{{\bf i}2W_j}\Big|}\\
		&\ \times \Big[\Big|\aleph^*\tilde{r}_i+\aleph_*+\tilde{r}_i\cos(V_1)e^{{\bf i}2V_j}+\cos(W_1)e^{{\bf i}2W_j}\Big|^2\\
		& \  +\Big|\tilde{r}_i\sin(V_1)e^{{\bf i}2V_j}+\sin(W_1)e^{{\bf i}2W_j}\Big|^2\Big]^{-\frac12},\\
		\sin(\Upsilon_*(\tilde{r}_i,\aleph^*,\aleph_*))\
		= \ & {\Big|\tilde{r}_i\sin(V_1)e^{{\bf i}2V_j}+\sin(W_1)e^{{\bf i}2W_j}\Big|}\\
		&\ \times \Big[\Big|\aleph^*\tilde{r}_i+\aleph_*+\tilde{r}_i\cos(V_1)e^{{\bf i}2V_j}+\cos(W_1)e^{{\bf i}2W_j}\Big|^2\\
		& \ +\Big|\tilde{r}_i\sin(V_1)e^{{\bf i}2V_j}+\sin(W_1)e^{{\bf i}2W_j}\Big|^2\Big]^{-\frac12},
	\end{aligned}
\end{equation*}
for $i=1,2,3$.
We define the singular manifold set, that handle the singular set that come from the oscillatory integrals in the rest of this section \begin{equation*}\label{Setstar}\begin{aligned}\mathfrak{S}(V,W)=\Big\{(\xi_1,\cdots,\xi_d)\in [-\pi,\pi]^d \Big|&\ \xi_1= \Upsilon_*(\tilde{r}_i,\aleph^*,\aleph_*)\pm \pi; i=1,2,3; \aleph^*,\aleph_*=\pm 1;\\
		&\ \xi_j/(2\pi)=0,\pm\frac14,\pm\frac12, j=1,\cdots,d \Big\}.\end{aligned}.\end{equation*}

We assume in some parts of the next lemmas that
\begin{equation}\label{Lemm:Bessel:aa}
	\begin{aligned}
		& |\cos(\aleph^j
		_1-\aleph^j_2)
		|<1, \ \ \ \ \   \Big|{t_0}+{t_1}\cos(V_1)e^{{\bf i}2V_j}+{t_2}\cos(W_1)e^{{\bf i}2W_j}\Big|>0,\\
		& \mbox{ and }\  \Big|{t_1}\sin(V_1)e^{{\bf i}2V_j}+{t_2}\sin(W_1)e^{{\bf i}2W_j}\Big|>0, \ \ \ \ j=2,\cdots,d.
	\end{aligned}
\end{equation}

For any $\eth>0,\lambda>0$ being small constants, there exists a smooth cut-off function $\check\Psi _1( \xi,V,W): [-\pi,\pi]^{3d}\to [0,1],$  
such that $\check\Psi_1( \xi,V,W) =1$ when $	|V^i_j|, \Big|V^i_j\pm\frac{1}{2}\Big|, |W^i_j|, \Big|W_j^i\pm \frac{1}{2}\Big|, |V_j^i\pm W_j^{i'}|   >|\ln\lambda|^{-\eth}$ and $\check\Psi _1( \xi,V,W) =0$ when either $  	|V^i_j|, \Big|V^i_j\pm\frac{1}{2}\Big|, |W^i_j|, \Big|W_j^i\pm \frac{1}{2}\Big|$, or $ |V_j^i\pm W_j^{i'}| <|\ln\lambda|^{-\eth}/2$, for $j=1,\cdots,d$, $i,i'=1,2$. 
Moreover, there also exists a  smooth cut-off function   
$\check{\Psi}_2( \xi,V,W) : [-\pi,\pi]^{3d}\to [0,1],$ such that   $\check{\Psi}_2=1$ when   $d( \xi,\mathfrak{S})>|\ln\lambda|^{-\eth}$, and $\check{\Psi}_2 =0$ when $d( \xi,\mathfrak{S})<|\ln\lambda|^{-\eth}/2$. 
We finally set $\check\Psi=\check\Psi_1\check\Psi_2$.

Let us consider the dispersion relation
\begin{equation}
	\label{Nearest}
	\omega(k) \ = \sin(2\pi k_1)\Big[\sin^2(2\pi k^1) + \cdots + \sin^2(2\pi k^d)\Big],
\end{equation}
where $k=(k^1,\cdots,k^d)\in\mathbb{T}^d$. By defining $\xi=(\xi_1,\cdots,\xi_d)=2\pi k\in [-\pi,\pi]^d$, we then simplify \eqref{Nearest} as
\begin{equation}
	\label{Nearest1}
	\omega(\xi) \ = \sin(\xi_1)\Big[\sin^2(\xi_1) + \cdots + \sin^2(\xi_d)\Big].
\end{equation}

We now consider the  functional
\begin{equation}
	\label{extendedBesselfunctions1a}
	\mathfrak{F}(m,t_0,t_1,t_2)  =  \int_{[-\pi,\pi]^d} \mathrm{d}\xi  e^{{\bf i}m\cdot \xi} e^{{\bf i}t_0\omega(\xi)+{\bf i}t_1\omega(\xi+V)+{\bf i}t_2\omega(\xi+W)},
\end{equation}
and
\begin{equation}
	\label{extendedBesselfunctions1a:1}
	\tilde{\mathfrak{F}}(m,t_0,t_1,t_2)  =  \int_{[-\pi,\pi]^d} \mathrm{d}\xi  \check\Psi(\xi,V,W)e^{{\bf i}m\cdot \xi} e^{{\bf i}t_0\omega(\xi)+{\bf i}t_1\omega(\xi+V)+{\bf i}t_2\omega(\xi+W)},
\end{equation}
for $m\in\mathbb{Z}^d$.

Let us rewrite \eqref{extendedBesselfunctions1a}-\eqref{extendedBesselfunctions1a:1} as follows

\begin{equation}
	\label{extendedBesselfunctions1}
	\begin{aligned}
		&\mathfrak{F}(m,t_0,t_1,t_2)  =   \int_{-\pi}^\pi\mathrm{d}\xi_1 \exp\Big({{\bf i}m_1\xi_1}\Big)\\
		&\ \ \times\exp\Big({\bf i}t_0\sin^3(\xi_1)+{\bf i}t_1\sin^3(\xi_1+V_1)+{\bf i}t_2\sin^3(\xi_1+W_1)\Big)\times\\
		&\ \  \times\Big[\prod_{j=2}^d\int_{-\pi}^\pi\mathrm{d}\xi_j \exp\Big({{\bf i}m_j\xi_j}\Big) \exp\Big({\bf i}t_0\sin(\xi_1)\sin^2(\xi_j)\\
		&\ \ +{\bf i}t_1\sin(\xi_1+V_1)\sin^2(\xi_j+V_j) +{\bf i}t_2\sin(\xi_1+W_1)\sin^2(\xi_j+W_j)\Big)\Big],
	\end{aligned}
\end{equation}
and
\begin{equation}
	\label{extendedBesselfunctions1}
	\begin{aligned}
		&\tilde{\mathfrak{F}}(m,t_0,t_1,t_2)  =   \int_{-\pi}^\pi\mathrm{d}\xi_1 \exp\Big({{\bf i}m_1\xi_1}\Big)\\
		&\ \ \times\exp\Big({\bf i}t_0\sin^3(\xi_1)+{\bf i}t_1\sin^3(\xi_1+V_1)+{\bf i}t_2\sin^3(\xi_1+W_1)\Big)\times\\
		&\ \  \times\Big[\prod_{j=2}^d\int_{-\pi}^\pi\mathrm{d}\xi_j \check\Psi(\xi,V,W)\exp\Big({{\bf i}m_j\xi_j}\Big) \exp\Big({\bf i}t_0\sin(\xi_1)\sin^2(\xi_j)\\
		&\ \ +{\bf i}t_1\sin(\xi_1+V_1)\sin^2(\xi_j+V_j) +{\bf i}t_2\sin(\xi_1+W_1)\sin^2(\xi_j+W_j)\Big)\Big].
	\end{aligned}
\end{equation}

\begin{lemma}\label{Lemm:Bessel2} There exists a universal constant $\mathfrak{C}_{\mathfrak{F},1}$ independent of $t_0,t_1,t_2,\aleph_1,\aleph_2$, such that 
	\begin{equation}
		\label{Lemm:Bessel2:1}
		\|\mathfrak{F}(\cdot,t_0,t_1,t_2) \|_{l^2} \ = \ \left( \sum_{m\in\mathbb{Z}^d}|\mathfrak{F}(m,t_0,t_1,t_2) |^2\right)^\frac12 \ \le \ \mathfrak{C}_{\mathfrak{F},1},
	\end{equation}
	and similarly
	\begin{equation}
		\label{Lemm:Bessel2:1:1}
		\|\tilde{\mathfrak{F}}(\cdot,t_0,t_1,t_2) \|_{l^2}  \ \le \ \mathfrak{C}_{\mathfrak{F},1},
	\end{equation}
\end{lemma}
\begin{proof}
	By the Plancherel theorem, we obtain
	\begin{equation}
		\label{Lemm:Bessel2:2}
		\begin{aligned}
			\sum_{m\in\mathbb{Z}^d}|\mathfrak{F}(m,t_0,t_1,t_2) |^2 
			\ = \ &  \int_{[-\pi,\pi]^d} \mathrm{d}\xi \Big|e^{{\bf i}t_0\omega(\xi)+{\bf i}t_1\omega(\xi+V)+{\bf i}t_2\omega(\xi+W)}\Big|^2,
		\end{aligned}
	\end{equation}
	which is a  bounded quantity. The conclusion \eqref{Lemm:Bessel2:1} of the lemma then follows. The second inequality \eqref{Lemm:Bessel2:1:1} can be proved by the same argument.
\end{proof} 
\begin{lemma}\label{Lemm:Bessel3}
	Under assumption \eqref{Lemm:Bessel:aa}, there exists a universal constant $\mathfrak{C}_{\mathfrak{F},4}>0$ independent of $t_0,t_1,t_2,\aleph_1,\aleph_2$, such that 
	\begin{equation}
		\label{Lemm:Bessel3:1}\begin{aligned}
			&\|\mathfrak{F}(\cdot,t_0,t_1,t_2) \|_{l^4}  \le  \mathfrak{C}_{\mathfrak{F},4}\prod_{j=2}^{d}\Big\langle\min\Big\{\Big|{t_0}+{t_1}\cos(V_1)e^{{\bf i}2V_j}+{t_2}\cos(W_1)e^{{\bf i}2W_j}\Big|,\\
			& \ \Big|{t_1}\sin(V_1)e^{{\bf i}2V_j}+{t_2}\sin(W_1)e^{{\bf i}2W_j}\Big|\Big\}[1-|\cos(\aleph^j
			_1-\aleph^j
			_2)|]^{\frac{1}{2}}\Big\rangle^{-(\frac{1}{8}-)}\\
			& \ \times \min\Big\{1,\sum_{j=2}^d\Big[\Big|{t_0}+{t_1}\cos(V_1)e^{{\bf i}2V_j}  +{t_2}\cos(W_1)e^{{\bf i}2W_j}\Big|\\
			&\ \ \ +\Big|{t_1}\sin(V_1)e^{{\bf i}2V_j}+{t_2}\sin(W_1)e^{{\bf i}2W_j}\Big|\Big]\\
			&\ \ \ \times\Big[\Big|{t_0}+{t_1}e^{{\bf i}3V_1}+{t_2}e^{{\bf i}3W_1}\Big|+(2d+1)\Big|{t_0}+{t_1}e^{{\bf i}V_1}+{t_2}e^{{\bf i}W_1}\Big|\Big]^{-1}\Big\}^\frac14.\end{aligned}
	\end{equation}
	In addition \eqref{Lemm:Bessel3:1} also holds true for $\tilde{\mathfrak{F}}$.
	Suppose further that \eqref{Lemm:Angle:0},\eqref{Lemm:Angle:4} and \eqref{Lemm:Angle:5} hold true, and 	\begin{equation}
		\label{Lemm:Bessel3:1:a}r_*=t_1/t_2=(1+\epsilon_{r_*})\tilde{r}_l\end{equation} for $l=1,2,3$ in which $\tilde{r}_l$ are defined in Lemma \ref{Lemm:Angle}. When \begin{equation}
		\label{Lemm:Bessel3:1:b} |\epsilon_{r_*}|=|\epsilon_{r_*}'|\langle
		\ln\lambda\rangle^{-c}\end{equation}  for an explicit constant $c>0$ depending only on the cut-off functions, and $\epsilon_{r_*}'$ is sufficiently small but independent of $\lambda$ and the cut-off functions,  then we have the estimate \begin{equation}
		\label{Lemm:Bessel3:2}\begin{aligned}
			&\|\tilde{\mathfrak{F}}(\cdot,t_0,t_1,t_2) \|_{l^4}  \le  \mathfrak{C}_{\mathfrak{F},4}\langle\ln\lambda\rangle^{\mathfrak{C}'_{\mathfrak{F},4}}\prod_{j=2}^{d}\Big\langle\Big\{\Big|{t_0}+{t_1}\cos(V_1)e^{{\bf i}2V_j}+{t_2}\cos(W_1)e^{{\bf i}2W_j}\Big|+\\
			& + \ \Big|{t_1}\sin(V_1)e^{{\bf i}2V_j}+{t_2}\sin(W_1)e^{{\bf i}2W_j}\Big|\Big\}|\cos(\aleph^j
			_1-\aleph^j
			_2)|^{\frac{1}{2}}\Big\rangle^{-(\frac18-)},\end{aligned}
	\end{equation}
	for universal constants $\mathfrak{C}_{\mathfrak{F},4},\mathfrak{C}_{\mathfrak{F},4}'>0.$
\end{lemma}
\begin{remark}\label{remark:Bessel2}
	The independence of $\mathfrak{C}_{\mathfrak{F},4}, \mathfrak{C}_{\mathfrak{F},4}'$ on $\aleph_1,\aleph_2$ is important, as we will later integrate $\aleph_1,\aleph_2$ out in our applications. The key difference in \eqref{Lemm:Bessel3:1} and \eqref{Lemm:Bessel3:2} is the alternative use of $1-|\cos(\aleph^j
	_1-\aleph^j
	_2)|$ and $|\cos(\aleph^j
	_1-\aleph^j
	_2)|$.
\end{remark}
\begin{proof}
	We observe that
	\begin{equation}
		\label{Lemm:Bessel2:1a}
		\begin{aligned}
			&|\mathfrak{F}(m,t_0,t_1,t_2) |^2  \ =  \ \int_{[-\pi,\pi]^d}\int_{[-\pi,\pi]^d}\mathrm{d}\xi\mathrm{d}\eta  e^{{\bf i}m\cdot\xi}e^{-{\bf i}m\cdot\eta}\\
			& \ \times e^{{\bf i}t_0\omega(\xi)+{\bf i}t_1\omega(\xi+V)+{\bf i}t_2\omega(\xi+W)}e^{-{\bf i}t_0\omega(\eta)-{\bf i}t_1\omega(\eta+V)-{\bf i}t_2\omega(\eta+W)},
		\end{aligned}
	\end{equation} 
	which, by the change of variable $\xi-\eta\to \xi$, gives
	\begin{equation}
		\label{Lemm:Bessel2:1a}
		\begin{aligned}
			& |\mathfrak{F}(m,t_0,t_1,t_2) |^2
			\ =  \ \int_{[-\pi,\pi]^d}\mathrm{d}\xi e^{{\bf i}m\cdot\xi} \int_{-\pi}^{\pi}\mathrm{d}\eta_1 \exp\Big({\bf i}t_0\sin^3(\xi_1+\eta_1)\\
			&   +{\bf i}t_1\sin^3(\xi_1+\eta_1+V_1)+{\bf i}t_2\sin^3(\xi_1+\eta_1+W_1)\\
			& -{\bf i}t_0\sin^3(\eta_1)-{\bf i}t_1\sin^3(\eta_1+V_1)-{\bf i}t_2\sin^3(\eta_1+W_1)\Big)\\
			& \  \times\Big[\prod_{j=2}^d\int_{-\pi}^\pi\mathrm{d}\eta_j \exp\Big({\bf i}t_0\sin(\xi_1+\eta_1)\sin^2(\xi_j+\eta_j) - {\bf i}t_0\sin(\eta_1)\sin^2(\eta_j) \\
			&\ \ +{\bf i}t_1\sin(\xi_1+\eta_1+V_1)\sin^2(\xi_j+\eta_j+V_j) +{\bf i}t_2\sin(\xi_1+\eta_1+W_1)\sin^2(\xi_j+\eta_j+W_j)\\
			&\ \ -{\bf i}t_1\sin(\eta_1+V_1)\sin^2(\eta_j+V_j)-{\bf i}t_2\sin(\eta_1+W_1)\sin^2(\eta_j+W_j)\Big)\Big]. 
		\end{aligned}
	\end{equation} 
	By the Plancherel theorem, we find
	\begin{equation}
		\label{Lemm:Bessel2:2}\begin{aligned}
			& \|\mathfrak{F}(\cdot,t_0,t_1,t_2) \|_{l^4}^4\ = \ \sum_{m\in\mathbb{Z}^d}|\mathfrak{F}(m,t_0,t_1,t_2) |^4\\  
			\ = \ & \int_{[-\pi,\pi]^d}\mathrm{d}\xi \Big| \int_{-\pi}^{\pi}\mathrm{d}\eta_1  \exp\Big({\bf i}t_0\sin^3(\xi_1+\eta_1)\\
			& +{\bf i}t_1\sin^3(\xi_1+\eta_1+V_1)  +{\bf i}t_2\sin^3(\xi_1+\eta_1+W_1)-{\bf i}t_0\sin^3(\eta_1)\\
			& -{\bf i}t_1\sin^3(\eta_1+V_1)-{\bf i}t_2\sin^3(\eta_1+W_1)\Big)\\
			& \  \times\Big[\prod_{j=2}^d\int_{-\pi}^\pi\mathrm{d}\eta_j \exp\Big({\bf i}t_0\sin(\xi_1+\eta_1)\sin^2(\xi_j+\eta_j)\\
			&\ \ +{\bf i}t_1\sin(\xi_1+\eta_1+V_1)\sin^2(\xi_j+\eta_j+V_j)\\
			&\ \  +{\bf i}t_2\sin(\xi_1+\eta_1+W_1)\sin^2(\xi_j+\eta_j+W_j)- {\bf i}t_0\sin(\eta_1)\sin^2(\eta_j)\\
			&\ \  -{\bf i}t_1\sin(\eta_1+V_1)\sin^2(\eta_j+V_j) -{\bf i}t_2\sin(\eta_1+W_1)\sin^2(\eta_j+W_j)\Big)\Big]\Big|^2.\end{aligned}
	\end{equation}
	Similarly, we also find
	\begin{equation}
		\label{Lemm:Bessel2:2:1}\begin{aligned}
			& \|\tilde{\mathfrak{F}}(\cdot,t_0,t_1,t_2) \|_{l^4}^4\ = \ \sum_{m\in\mathbb{Z}^d}|\tilde{\mathfrak{F}}(m,t_0,t_1,t_2) |^4\\  
			\ = \ & \int_{[-\pi,\pi]^d}\mathrm{d}\xi \Big| \int_{-\pi}^{\pi}\mathrm{d}\eta_1  \check\Psi(\xi+\eta,V,W)\check\Psi(\eta,V,W)\exp\Big({\bf i}t_0\sin^3(\xi_1+\eta_1)\\
			& +{\bf i}t_1\sin^3(\xi_1+\eta_1+V_1)  +{\bf i}t_2\sin^3(\xi_1+\eta_1+W_1)-{\bf i}t_0\sin^3(\eta_1)\\
			& -{\bf i}t_1\sin^3(\eta_1+V_1)-{\bf i}t_2\sin^3(\eta_1+W_1)\Big)\\
			& \  \times\Big[\prod_{j=2}^d\int_{-\pi}^\pi\mathrm{d}\eta_j \exp\Big({\bf i}t_0\sin(\xi_1+\eta_1)\sin^2(\xi_j+\eta_j)\\
			&\ \ +{\bf i}t_1\sin(\xi_1+\eta_1+V_1)\sin^2(\xi_j+\eta_j+V_j)\\
			&\ \  +{\bf i}t_2\sin(\xi_1+\eta_1+W_1)\sin^2(\xi_j+\eta_j+W_j)- {\bf i}t_0\sin(\eta_1)\sin^2(\eta_j)\\
			&\ \  -{\bf i}t_1\sin(\eta_1+V_1)\sin^2(\eta_j+V_j) -{\bf i}t_2\sin(\eta_1+W_1)\sin^2(\eta_j+W_j)\Big)\Big]\Big|^2.\end{aligned}
	\end{equation}
	
	{\bf (i) First, we prove \eqref{Lemm:Bessel3:1} for $\mathfrak{F}$.} The proof of \eqref{Lemm:Bessel3:1} for $\tilde{\mathfrak{F}}$ can be done by precisely the same argument.
	Let us define for $2\le j\le d$
	\begin{equation}
		\begin{aligned}
			\mathfrak{A}_j(\xi_1,\eta_1,\xi_j) \ = & \ \int_{-\pi}^\pi\mathrm{d}\eta_j e^{{\bf i}\mathfrak{B}_j(\xi_1,\eta_1,\xi_j,\eta_j)},
		\end{aligned}
	\end{equation}
	in which
	\begin{equation}
		\begin{aligned}
			& \mathfrak{B}_j(\xi_1,\eta_1,\xi_j,\eta_j) =  \ t_0\sin(\xi_1+\eta_1)\sin^2(\xi_j+\eta_j)\\
			&\ +t_1\sin(\xi_1+\eta_1+V_1)\sin^2(\xi_j+\eta_j+V_j)\\
			&\  +t_2\sin(\xi_1+\eta_1+W_1)\sin^2(\xi_j+\eta_j+W_j)- t_0\sin(\eta_1)\sin^2(\eta_j)\\
			&\  -t_1\sin(\eta_1+V_1)\sin^2(\eta_j+V_j)-t_2\sin(\eta_1+W_1)\sin^2(\eta_j+W_j),
		\end{aligned}
	\end{equation}
	and
	\begin{equation}
		\begin{aligned}
			& \mathfrak{B}_1(\xi_1,\eta_1) =  \ t_0\sin^3(\xi_1+\eta_1) +t_1\sin^3(\xi_1+\eta_1+V_1)\\
			&  +t_2\sin^3(\xi_1+\eta_1+W_1)-t_0\sin^3(\eta_1)-t_1\sin^3(\eta_1+V_1)-t_2\sin^3(\eta_1+W_1),
		\end{aligned}
	\end{equation}
	we find
	\begin{equation}
		\label{Lemm:Bessel2:2a}\begin{aligned}
			\sum_{m\in\mathbb{Z}^d}|\mathfrak{F}(m,t_0,t_1,t_2) |^2
			\ = \ & \int_{[-\pi,\pi]^d}\mathrm{d}\xi \Big|\int_{-\pi}^{\pi}\mathrm{d}\eta_1 \prod_{j=2}^d\mathfrak{A}_j(\xi_1,\eta_1,\xi_j) e^{{\bf i}\mathfrak{B}_1(\xi_1,\eta_1) }\Big|^2.\end{aligned}
	\end{equation}

	We now study the oscillatory integrals $\mathfrak{A}_j$ by writing the phase $\mathfrak{B}_j$ as 
	\begin{equation}
		\begin{aligned}
			&\mathfrak{B}_j(\xi_1,\eta_1,\xi_j,\eta_j)
			=   \frac{t_0}{2}\sin(\xi_1+\eta_1)[1-\cos(2\xi_j+2\eta_j)]\\
			&\  + \ \frac{t_1}{2}\sin(\xi_1+\eta_1+V_1)[1-\cos(2\xi_j+2\eta_j+2V_j)]\\
			&\  +\ \frac{t_2}{2}\sin(\xi_1+\eta_1+W_1)[1-\cos(2\xi_j+2\eta_j+2W_j)]\\
			&\ - \ \frac{t_0}{2}\sin(\eta_1)[1-\cos(2\eta_j)]\
			-\frac{t_1}{2}\sin(\eta_1+V_1)[1-\cos(2\eta_j+2V_j)]\\
			&\ -\ \frac{t_2}{2}\sin(\eta_1+W_1)[1-\cos(2\eta_j+2W_j)],
		\end{aligned}
	\end{equation}
	which could be split as the sum of
	\begin{equation}
		\begin{aligned}
			& \mathfrak{B}^a_j(\xi_1,\eta_1,\xi_j,\eta_j)
			\ =   \ -\frac{t_0}{2}\sin(\xi_1+\eta_1)\cos(2\xi_j+2\eta_j)\\
			&\  - \ \frac{t_1}{2}\sin(\xi_1+\eta_1+V_1)\cos(2\xi_j+2\eta_j+2V_j)\\
			&\  - \ \frac{t_2}{2}\sin(\xi_1+\eta_1+W_1)\cos(2\xi_j+2\eta_j+2W_j)\ + \ \frac{t_0}{2}\sin(\eta_1)\cos(2\eta_j)\\
			&\
			+\frac{t_1}{2}\sin(\eta_1+V_1)\cos(2\eta_j+2V_j)\ + \ \frac{t_2}{2}\sin(\eta_1+W_1)\cos(2\eta_j+2W_j),
		\end{aligned}
	\end{equation}
	and
	\begin{equation}
		\begin{aligned}
			&\mathfrak{B}^b_j(\xi_1,\eta_1)
			=   \frac{t_0}{2}\sin(\xi_1+\eta_1)  +  \frac{t_1}{2}\sin(\xi_1+\eta_1+V_1)\\
			&  + \frac{t_2}{2}\sin(\xi_1+\eta_1+W_1)
			-  \frac{t_0}{2}\sin(\eta_1)\
			-\frac{t_1}{2}\sin(\eta_1+V_1) - \frac{t_2}{2}\sin(\eta_1+W_1).
		\end{aligned}
	\end{equation}
	The oscillatory integral $\mathfrak{A}_j$ can now be written
	\begin{equation}
		\mathfrak{A}_j \ = \ e^{{\bf i}\mathfrak{B}^b_j }\int_{-\pi}^\pi\mathrm{d}\eta_j e^{{\bf i}\mathfrak{B}^a_j}.
	\end{equation}
	Let us express $\mathfrak{B}^a_j$ into the following form
	\begin{equation}
		\begin{aligned}
			& \mathfrak{B}^a_j(\xi_1,\eta_1,\xi_j,\eta_j)
			\ =   \ -\ \mathrm{Re}\Big[e^{{\bf i}2\xi_j+{\bf i}2\eta_j}\Big(\frac{t_0}{2}\sin(\xi_1+\eta_1)+ \frac{t_1}{2}\sin(\xi_1+\eta_1+V_1)e^{{\bf i}2V_j}\\
			&\  + \ \frac{t_2}{2}\sin(\xi_1+\eta_1+W_1)e^{{\bf i}2W_j}\Big)\Big]\\
			& + \ \mathrm{Re}\Big[e^{{\bf i}2\eta_j}\Big(
			\frac{t_0}{2}\sin(\eta_1)\
			+\frac{t_1}{2}\sin(\eta_1+V_1)e^{{\bf i}2V_j} +  \frac{t_2}{2}\sin(\eta_1+W_1)e^{{\bf i}2W_j}\Big)\Big].
		\end{aligned}
	\end{equation}
	Setting $$\mathfrak{C}_j^1=\frac{t_0}{2}\sin(\xi_1+\eta_1)+ \frac{t_1}{2}\sin(\xi_1+\eta_1+V_1)e^{{\bf i}2V_j}\  + \ \frac{t_2}{2}\sin(\xi_1+\eta_1+W_1)e^{{\bf i}2W_j}$$
	and
	$$\mathfrak{C}_j^2=\frac{t_0}{2}\sin(\eta_1)\
	+\frac{t_1}{2}\sin(\eta_1+V_1)e^{{\bf i}2V_j} +  \frac{t_2}{2}\sin(\eta_1+W_1)e^{{\bf i}2W_j},$$
	we find
	\begin{equation}
		\begin{aligned}\label{204new}
			& \mathfrak{B}^a_j(\xi_1,\eta_1,\xi_j,\eta_j)
			\ =   \ -\ \mathrm{Re}\Big[e^{{\bf i}2\xi_j+{\bf i}2\eta_j}\mathfrak{C}_j^1\Big]\ + \ \mathrm{Re}\Big[e^{{\bf i}2\eta_j}\mathfrak{C}_j^2\Big],
		\end{aligned}
	\end{equation}
	and define
	\begin{equation}
		\mathfrak{A}_j^a \ = \ \int_{-\pi}^\pi\mathrm{d}\eta_j e^{{\bf i}\mathfrak{B}^a_j}.
	\end{equation}
	We combine the phases ${\mathfrak{B}^b_j }$ and ${\mathfrak{B}_1}$  
	\begin{equation}
		\begin{aligned}\label{206new}
			& \mathfrak{B}^a_1=\mathfrak{B}_1 + \sum_{j=2}^{d}\mathfrak{B}^b_j \  =  \ t_0\Big[\sin^3(\xi_1+\eta_1)+\frac{d-1}{2}\sin(\xi_1+\eta_1)\Big]\\
			& +t_1\Big[\sin^3(\xi_1+\eta_1+V_1)+\frac{d-1}{2}\sin(\xi_1+\eta_1+V_1)\Big]\\
			&  +t_2\Big[\sin^3(\xi_1+\eta_1+W_1)+\frac{d-1}{2}\sin(\xi_1+\eta_1+W_1)\Big]\\
			& -t_0\Big[\sin^3(\eta_1)+\frac{d-1}{2}\sin(\eta_1)\Big]-t_1\Big[\sin^3(\eta_1+V_1)+\frac{d-1}{2}\sin(\eta_1+V_1)\Big]
			\\
			&-t_2\Big[\sin^3(\eta_1+W_1)+\frac{d-1}{2}\sin(\eta_1+W_1)\Big].
		\end{aligned}
	\end{equation}
	We then write
	\begin{equation}
		\label{Lemm:Bessel2:2a:b1}\begin{aligned}
			\sum_{m\in\mathbb{Z}^d}|\mathfrak{F}(m,t_0,t_1,t_2) |^2
			\ = \ & \int_{[-\pi,\pi]^d}\mathrm{d}\xi \Big|\int_{-\pi}^{\pi}\mathrm{d}\eta_1 \prod_{j=2}^d\mathfrak{A}_j^a(\xi_1,\eta_1,\xi_j) e^{{\bf i}\mathfrak{B}^a_1(\xi_1,\eta_1) }\Big|^2.\end{aligned}
	\end{equation}		
	We now set
	\begin{equation}
		\mathfrak{C}_j^1=\Xi_j^1 e^{{\bf i}2\Upsilon_j^1}, \ \ \ \ \mathfrak{C}_j^2=\Xi_j^2 e^{{\bf i}2\Upsilon_j^2},
	\end{equation}
	with $\Xi_j^1,\Xi_j^2\in\mathbb{R}_+$ and $\Upsilon_j^1,\Upsilon_j^2\in[-\pi,\pi]$, then
	$$\Xi_j^1=\Big|\frac{t_0}{2}\sin(\xi_1+\eta_1)+ \frac{t_1}{2}\sin(\xi_1+\eta_1+V_1)e^{{\bf i}2V_j}\  + \ \frac{t_2}{2}\sin(\xi_1+\eta_1+W_1)e^{{\bf i}2W_j}\Big|,$$
	we develop this term as
	\begin{equation}
		\begin{aligned}\label{Lemm:Bessel2:2bb7}
			\Xi_j^1\ =\ & \Big|\Big(\frac{t_0}{2}+\frac{t_1}{2}\cos(V_1)e^{{\bf i}2V_j}+\frac{t_2}{2}\cos(W_1)e^{{\bf i}2W_j}\Big)\sin(\xi_1+\eta_1)\\
			& + \Big(\frac{t_1}{2}\sin(V_1)e^{{\bf i}2V_j}+\frac{t_2}{2}\sin(W_1)e^{{\bf i}2W_j}\Big)\cos(\xi_1+\eta_1)\Big|.
		\end{aligned}
	\end{equation}
	Next, we will perform an a priori estimate on $\Xi_j^1$ and $\Xi_j^2,$ to obtain uniform lower bounds independent of $\xi_1+\eta_1$ and $\eta_1$ of those quantities. 
	
	Setting 
	\begin{equation}
		\begin{aligned}\label{Lemm:Bessel2:2bb8}
			\cos(\Upsilon_*(j))
			\ = \ & {\Big|\frac{t_0}{2}+\frac{t_1}{2}\cos(V_1)e^{{\bf i}2V_j}+\frac{t_2}{2}\cos(W_1)e^{{\bf i}2W_j}\Big|}\\
			&\ \times \Big[\Big|\frac{t_0}{2}+\frac{t_1}{2}\cos(V_1)e^{{\bf i}2V_j}+\frac{t_2}{2}\cos(W_1)e^{{\bf i}2W_j}\Big|^2\\
			& \  +\Big|\frac{t_1}{2}\sin(V_1)e^{{\bf i}2V_j}+\frac{t_2}{2}\sin(W_1)e^{{\bf i}2W_j}\Big|^2\Big]^{-\frac12},\\
			\sin(\Upsilon_*(j))\
			= \ & {\Big|\frac{t_1}{2}\sin(V_1)e^{{\bf i}2V_j}+\frac{t_2}{2}\sin(W_1)e^{{\bf i}2W_j}\Big|}\\
			&\ \times \Big[\Big|\frac{t_0}{2}+\frac{t_1}{2}\cos(V_1)e^{{\bf i}2V_j}+\frac{t_2}{2}\cos(W_1)e^{{\bf i}2W_j}\Big|^2\\
			& \ +\Big|\frac{t_1}{2}\sin(V_1)e^{{\bf i}2V_j}+\frac{t_2}{2}\sin(W_1)e^{{\bf i}2W_j}\Big|^2\Big]^{-\frac12},
		\end{aligned}
	\end{equation}
	with $\Upsilon_*\in[0,\pi/2]$, we then find
	\begin{equation}
		\begin{aligned}\label{Lemm:Bessel2:2bb9a}
			\Xi_j^1
			\ =\ 
			& \Big[\Big|\frac{t_0}{2}+\frac{t_1}{2}\cos(V_1)e^{{\bf i}2V_j}+\frac{t_2}{2}\cos(W_1)e^{{\bf i}2W_j}\Big|^2\\
			&+\Big|\frac{t_1}{2}\sin(V_1)e^{{\bf i}2V_j}+\frac{t_2}{2}\sin(W_1)e^{{\bf i}2W_j}\Big|^2\Big]^\frac12\\
			&\times \Big|\cos(\Upsilon_*)\sin(\xi_1+\eta_1)e^{{\bf i}\aleph
				_1^j}+\sin(\Upsilon_*)\cos(\xi_1+\eta_1)e^{{\bf i}\aleph
				_2^j}\Big|.
		\end{aligned}
	\end{equation}
	
	Let us now study the factor containing $\Upsilon_*$ on the right hand side
	\begin{equation}
		\begin{aligned}\label{Lemm:Bessel2:2bb9b}
			& \Big|\cos(\Upsilon_*)\sin(\xi_1+\eta_1)e^{{\bf i}\aleph^j
				_1}+\sin(\Upsilon_*)\cos(\xi_1+\eta_1)e^{{\bf i}\aleph^j
				_2}\Big|^2\\
			=\	&[\cos(\Upsilon_*)\sin(\xi_1+\eta_1)\cos(\aleph^j
			_1)+\sin(\Upsilon_*)\cos(\xi_1+\eta_1)\cos(\aleph^j
			_2)]^2\\
			&+[\cos(\Upsilon_*)\sin(\xi_1+\eta_1)\sin(\aleph^j
			_1)+\sin(\Upsilon_*)\cos(\xi_1+\eta_1)\sin(\aleph^j
			_2)]^2\\
			=\	&\cos^2(\Upsilon_*)\sin^2(\xi_1+\eta_1)+\sin^2(\Upsilon_*)\cos^2(\xi_1+\eta_1)\\
			&+2\cos(\Upsilon_*)\sin(\Upsilon_*)\sin(\xi_1+\eta_1)\cos(\xi_1+\eta_1)\cos(\aleph^j
			_1-\aleph^j
			_2).
		\end{aligned}
	\end{equation}
	Let us consider the case when $\cos(\aleph^j
	_1-\aleph^j
	_2)\ge0$. The above quantity can be bounded from below as
	\begin{equation}
		\begin{aligned}\label{Lemm:Bessel2:2bb9c}
			& \Big|\cos(\Upsilon_*)\sin(\xi_1+\eta_1)e^{{\bf i}\aleph^j
				_1}+\sin(\Upsilon_*)\cos(\xi_1+\eta_1)e^{{\bf i}\aleph^j
				_2}\Big|^2\\
			=\	&[\cos^2(\Upsilon_*)\sin^2(\xi_1+\eta_1)+\sin^2(\Upsilon_*)\cos^2(\xi_1+\eta_1)][1-\cos(\aleph
			_1^j-\aleph
			_2^j)]\\
			& + \sin^2(\xi_1+\eta_1+\Upsilon_*)\cos(\aleph^j
			_1-\aleph^j
			_2)\\
			\ge\ & \min\{\cos^2(\Upsilon_*),\sin^2(\Upsilon_*)\}[1-\cos(\aleph^j
			_1-\aleph^j
			_2)]+ \sin^2(\xi_1+\eta_1+\Upsilon_*)\cos(\aleph^j
			_1-\aleph
			_2^j).
		\end{aligned}
	\end{equation}
	We then deduce
	\begin{equation}
		\begin{aligned}\label{Lemm:Bessel2:2bb9a}
			&\Xi_j^1
			\ \ge \ 
			\Big[\Big|\frac{t_0}{2}+\frac{t_1}{2}\cos(V_1)e^{{\bf i}2V_j}+\frac{t_2}{2}\cos(W_1)e^{{\bf i}2W_j}\Big|^2\\
			&+\Big|\frac{t_1}{2}\sin(V_1)e^{{\bf i}2V_j}+\frac{t_2}{2}\sin(W_1)e^{{\bf i}2W_j}\Big|^2\Big]^\frac12\min\{\cos^2(\Upsilon_*),\sin^2(\Upsilon_*)\}^\frac12[1-\cos(\aleph^j
			_1-\aleph^j
			_2)]^\frac12.
		\end{aligned}
	\end{equation}
	
	When $\cos(\aleph^j
	_1-\aleph^j
	_2)<0$, we can also estimate
	\begin{equation}
		\begin{aligned}\label{Lemm:Bessel2:2bb9c}
			& \Big|\cos(\Upsilon_*)\sin(\xi_1+\eta_1)e^{{\bf i}\aleph^j
				_1}+\sin(\Upsilon_*)\cos(\xi_1+\eta_1)e^{{\bf i}\aleph^j
				_2}\Big|^2\\
			=\	&\cos^2(\Upsilon_*)\sin^2(\xi_1+\eta_1)+\sin^2(\Upsilon_*)\cos^2(\xi_1+\eta_1)\\
			&-2\cos(\Upsilon_*)\sin(\Upsilon_*)\sin(\xi_1+\eta_1)\cos(\xi_1+\eta_1)|\cos(\aleph^j
			_1-\aleph^j
			_2)|\\
			=\	&[\cos^2(\Upsilon_*)\sin^2(\xi_1+\eta_1)+\sin^2(\Upsilon_*)\cos^2(\xi_1+\eta_1)][1-|\cos(\aleph
			_1^j-\aleph
			_2^j)|]\\
			& + \sin^2(\xi_1+\eta_1-\Upsilon_*)|\cos(\aleph^j
			_1-\aleph^j
			_2)|\\
			\ge\ & \min\{\cos^2(\Upsilon_*),\sin^2(\Upsilon_*)\}[1-|\cos(\aleph^j
			_1-\aleph^j
			_2)|]+ \sin^2(\xi_1+\eta_1-\Upsilon_*)|\cos(\aleph^j
			_1-\aleph
			_2^j)|.
		\end{aligned}
	\end{equation}
	As a consequence, we bound
	\begin{equation}
		\begin{aligned}\label{Lemm:Bessel2:2bb9c:1}
			&\Xi_j^1
			\ \gtrsim \ 
			\min\{\cos^2(\Upsilon_*),\sin^2(\Upsilon_*)\}^\frac12 \Big[\Big|\frac{t_0}{2}+\frac{t_1}{2}\cos(V_1)e^{{\bf i}2V_j}+\frac{t_2}{2}\cos(W_1)e^{{\bf i}2W_j}\Big|^2\\
			&+\Big|\frac{t_1}{2}\sin(V_1)e^{{\bf i}2V_j}+\frac{t_2}{2}\sin(W_1)e^{{\bf i}2W_j}\Big|^2\Big]^\frac12 [1-|\cos(\aleph^j
			_1-\aleph^j
			_2)|]^\frac12,
		\end{aligned}
	\end{equation}	
	and similarly
	\begin{equation}
		\begin{aligned}\label{Lemm:Bessel2:2bb9c:2}
			&\Xi_j^2
			\ \gtrsim \ 
			\min\{\cos^2(\Upsilon_*),\sin^2(\Upsilon_*)\}^\frac12\Big[\Big|\frac{t_0}{2}+\frac{t_1}{2}\cos(V_1)e^{{\bf i}2V_j}+\frac{t_2}{2}\cos(W_1)e^{{\bf i}2W_j}\Big|^2\\
			&+\Big|\frac{t_1}{2}\sin(V_1)e^{{\bf i}2V_j}+\frac{t_2}{2}\sin(W_1)e^{{\bf i}2W_j}\Big|^2\Big]^\frac12 [1-|\cos(\aleph^j
			_1-\aleph^j
			_2)|]^\frac12.
		\end{aligned}
	\end{equation}
	Those are the lower bounds needed for $\Xi_j^1$ and $\Xi_j^2$. Note that, when estimating $\tilde{\mathfrak{F}}$, we can improve the above lower bounds as  $\sin^2(\xi_1+\eta_1-\Upsilon_*)$ is also  bounded from below thanks to the cut-off function $\Psi_4$.
	
	\medskip
	{\bf Step 1: Estimating $\mathfrak{A}_j^a$.}  We divide this step into three smaller steps.
	
	\smallskip
	{\bf Step 1.1: Preliminary bounds on the derivatives in $\eta_j$ of $\mathfrak{B}^c_j$. }
	Setting $	\Xi_j^0=\Xi_j^1+	\Xi_j^2$, $\mathfrak{B}^a_j =	\Xi_j^0\mathfrak{B}^c_j$,  
	we now compute the derivatives in $\eta_j$ of $\mathfrak{B}^c_j$, and provide some preliminary bounds on them, before going into the details of estimating $\mathfrak{A}_j^a$. We have
	\begin{equation}
		\begin{aligned}\label{Lemm:Bessel2:2aa1}
			\partial_{\eta_j}\mathfrak{B}^c_j \ = \ & \ 	2\mathrm{Im}\Big[e^{{\bf i}2\xi_j+{\bf i}2\eta_j}\mathfrak{C}_j^3\Big]\ - \ 		2\mathrm{Im}\Big[e^{{\bf i}2\eta_j}\mathfrak{C}_j^4\Big]\end{aligned}
	\end{equation}
	and
	\begin{equation}
		\begin{aligned}\label{Lemm:Bessel2:2aa2}
			\partial_{\eta_j\eta_j}\mathfrak{B}^c_j \ = \ &	 4\mathrm{Re}\Big[e^{{\bf i}2\xi_j+{\bf i}2\eta_j}\mathfrak{C}_j^3\Big]\ - \ 	 4\mathrm{Re}\Big[e^{{\bf i}2\eta_j}\mathfrak{C}_j^4\Big],\end{aligned}
	\end{equation}
	with $\mathfrak{C}_j^1=	\Xi_j^0\mathfrak{C}_j^3$, $\mathfrak{C}_j^2=	\Xi_j^0\mathfrak{C}_j^4$ and $\Xi_j^3$ is defined to be the quantity that satisfies the two identities 
	$$\Xi_j^1=\Xi_j^0\Xi_j^3, \qquad \Xi_j^2=\Xi_j^0(1-\Xi_j^3).$$
	
	In our consideration, the quantity $	\Xi_j^0$ is big,  $	\Xi_j^0>1$, leading to the oscillation of the integral.

	From \eqref{Lemm:Bessel2:2aa2}, we observe that, for all $m\in\mathbb{N}$
	\begin{equation}
		\begin{aligned}\label{Lemm:Bessel2:2aa1:1}
			\Big|\partial_{\eta_j}^{2m+1}\mathfrak{B}^c_j	\Big| \ = \ & \ 		\Big|2^{2m+1}\mathrm{Im}\Big[e^{{\bf i}2\xi_j+{\bf i}2\eta_j}\mathfrak{C}_j^3\Big]\ - \ 		2^{2m+1}\mathrm{Im}\Big[e^{{\bf i}2\eta_j}\mathfrak{C}_j^4\Big]	\Big| = \ 2^{2m}	\Big|\partial_{\eta_j}\mathfrak{B}^c_j	\Big|\end{aligned}
	\end{equation}
	and
	\begin{equation}
		\begin{aligned}\label{Lemm:Bessel2:2aa2:2}
			\Big|	\partial_{\eta_j}^{2m}\mathfrak{B}^c_j 	\Big|\ = \ &		\Big| 2^{2m}\mathrm{Re}\Big[e^{{\bf i}2\xi_j+{\bf i}2\eta_j}\mathfrak{C}_j^3\Big]\ - \ 	 2^{2m}\mathrm{Re}\Big[e^{{\bf i}2\eta_j}\mathfrak{C}_j^4\Big]	\Big|= 2^{2m}	\Big|\mathfrak{B}^c_j	\Big|.\end{aligned}
	\end{equation}
	
	We denote by $\eta_j^*$ the solution of
	$
	\partial_{\eta_j}\mathfrak{B}^c_j \ = \  0,
	$
	which is equivalent to
	\begin{equation}\label{Lemm:Bessel2:2bb1}
		\Xi_j^3 \sin(2(\xi_j+\eta_j^*+\Upsilon_j^1)) \ = \ (1-\Xi_j^3) \sin(2(\eta_j^*+\Upsilon_j^2)).
	\end{equation}
	Noticing that \eqref{Lemm:Bessel2:2bb1} always has a finite number of solutions, independent of the choices of $\xi_j,\xi_1,\eta_1$, except the case when both of the following identities happen $\Xi_j^3=\frac12$ and $2\xi_j+2\Upsilon_j^1-2\Upsilon_j^2 = m\pi$, for $m\in\mathbb{Z}$. We restrict our domain of defining for $\xi_j$ to
	\begin{equation}\label{Lemm:Bessel2:2bb1:1}\mathbb{T}_{\xi_j}=\Big\{\xi_j\in[-\pi,\pi],|2\xi_j+2\Upsilon_j^1-2\Upsilon_j^2 - m\pi|> \epsilon_{\xi_j}\Big\}\end{equation} with $m=0,\pm1,\pm2,\pm3,\pm4,\pm5,\pm6$ and the small constant $\epsilon_{\xi_j}>0$ to be determined later. For $\xi_j\in  \mathbb{T}_{\xi_j}$, equation \eqref{Lemm:Bessel2:2bb1} becomes
	\begin{equation}\label{Lemm:Bessel2:2bb1:1}\begin{aligned}
			& \Xi_j^3 [\sin(2\eta_j^*)\cos(2\xi_j+2\Upsilon_j^1)+\cos(2\eta_j^*)\sin(2\xi_j+2\Upsilon_j^1)] \\
			\ = \ & (1-\Xi_j^3) [\sin(2\eta_j^*)\cos(2\Upsilon_j^2)+\cos(2\eta_j^*)\sin(2\Upsilon_j^2)],\end{aligned}
	\end{equation}
	which is equivalent to
	\begin{equation}\label{Lemm:Bessel2:2bb1:2}\begin{aligned}
			& \sin(2\eta_j^*)[\Xi_j^3 \cos(2\xi_j+2\Upsilon_j^1) -(1-\Xi_j^3) \cos(2\Upsilon_j^2)] \\
			\ = \ & \cos(2\eta_j^*)[(1-\Xi_j^3) \sin(2\Upsilon_j^2)-\Xi_j^3\sin(2\xi_j+2\Upsilon_j^1)].\end{aligned}
	\end{equation}
	Equation \eqref{Lemm:Bessel2:2bb1:2} has a unique solution when $\Xi_j^3 \cos(2\xi_j+2\Upsilon_j^1) -(1-\Xi_j^3) \cos(2\Upsilon_j^2)\ne 0$ and $(1-\Xi_j^3) \sin(2\Upsilon_j^2)-\Xi_j^3\sin(2\xi_j+2\Upsilon_j^1)\ne 0$. It has at most $3$ solutions when one of the two coefficients $\Xi_j^3 \cos(2\xi_j+2\Upsilon_j^1) -(1-\Xi_j^3) \cos(2\Upsilon_j^2)$ and $(1-\Xi_j^3) \sin(2\Upsilon_j^2)-\Xi_j^3\sin(2\xi_j+2\Upsilon_j^1)$ is $0$ while the other one is different from $0$. Therefore, for $\xi_j\in  \mathbb{T}_{\xi_j}$, $\mathfrak{B}^c_j $ has at most $3$ stationary points.

	To make sure that $\Xi_j^0$ is the dominant parameter in our computations, we will put a constraint on the choice of $\epsilon_{\xi_j}$
	\begin{equation}\label{Lemm:Bessel2:2bb6:1}
		1 \gg\frac{1}{\Xi_j^0\epsilon_{\xi_j}^2}.
	\end{equation}  We define by $\mathfrak{I}_{\eta_j}$ the set of all solutions $\eta_j^*$ of \eqref{Lemm:Bessel2:2bb1}.  
	We compute, for $\eta_j^*\in \mathfrak{I}_{\eta_j}$
	
	\begin{equation}\begin{aligned}\label{Lemm:Bessel2:2bb4}
			& |\Xi_j^3\cos(2(\xi_j+\eta_j^*+\Upsilon_j^1))\ - \ (1-\Xi_j^3) \cos(2(\eta_j^*+\Upsilon_j^2))||\sin(2\eta_j^*+2\Upsilon_j^2)|\\
			=\ &	|(1-\Xi_j^3)\sin(2\eta_j^*+2\Upsilon_j^2)\cos(2(\eta_j^*+\Upsilon_j^2)) -\Xi_j^3\sin(2\eta_j^*+2\Upsilon_j^2)\cos(2(\xi_j+\eta_j^*+\Upsilon_j^1))|\\
			=\ &	|\Xi_j^3 \sin(2\xi_j+2\eta_j^*+2\Upsilon_j^1) \cos(2(\eta_j^*+\Upsilon_j^2))
			-\Xi_j^3\sin(2\eta_j^*+2\Upsilon_j^2)\cos(2(\xi_j+\eta_j^*+\Upsilon_j^1))|\\
			=\ &	{\Xi_j^3 |\sin(2\xi_j+2\Upsilon_j^1-2\Upsilon_j^2)|},
		\end{aligned}
	\end{equation}
	which yields
	\begin{equation}\label{Lemm:Bessel2:2bb4:1}
		|\Xi_j^3\cos(2(\xi_j+\eta_j^*+\Upsilon_j^1))-  (1-\Xi_j^3) \cos(2(\eta_j^*+\Upsilon_j^2))|
		\ge 	{\Xi_j^3 |\sin(2\xi_j+2\Upsilon_j^1-2\Upsilon_j^2)|} >0.
	\end{equation}
	A similar argument also gives
	\begin{equation}\begin{aligned}\label{Lemm:Bessel2:2bb4:2}
			& |\Xi_j^3\cos(2(\xi_j+\eta_j^*+\Upsilon_j^1))\ - \ (1-\Xi_j^3)\cos(2(\eta_j^*+\Upsilon_j^2))|\\
			\ge \  &	{(1-\Xi_j^3) |\sin(2\xi_j+2\Upsilon_j^1-2\Upsilon_j^2)|} >0.
		\end{aligned}
	\end{equation}
	As a consequence, we have the bound
	\begin{equation}\label{Lemm:Bessel2:2bb4:3}
		|\Xi_j^3\cos(2(\xi_j+\eta_j^*+\Upsilon_j^1))-  (1-\Xi_j^3) \cos(2(\eta_j^*+\Upsilon_j^2))|
		\ge 	{\frac12 |\sin(2\xi_j+2\Upsilon_j^1-2\Upsilon_j^2)|} >0,
	\end{equation}
	leading to
	\begin{equation}\label{Lemm:Bessel2:2bb5}
		|	\partial_{\eta_j\eta_j}\mathfrak{B}^c_j(\eta_j^*)|
		\ge	2{|\sin(2\xi_j+2\Upsilon_j^1-2\Upsilon_j^2)|}.
	\end{equation} 
	
	Let us also define $\mathfrak{J}_{\eta_j}^*$ to be the set of solutions of 
	$
	\mathfrak{B}^c_j \ = \  0,
	$
	which is equivalent to
	\begin{equation}\label{Lemm:Bessel2:2bb1:1:1}
		\Xi_j^3 \cos(2(\xi_j+\eta_j^{**}+\Upsilon_j^1)) \ = \ (1-\Xi_j^3) \cos(2(\eta_j^{**}+\Upsilon_j^2)).
	\end{equation}
	The same argument as above also shows that, for $\xi_j\in  \mathbb{T}_{\xi_j}$, $\mathfrak{B}^c_j=0$ has at most $3$ solutions. Moreover, for all $\eta_j^{**}\in \mathfrak{J}_{\eta_j}^*$
	\begin{equation}\label{Lemm:Bessel2:2bb5:b:1}
		|	\partial_{\eta_j}\mathfrak{B}^c_j(\eta_j^{**})|
		\ge	 {|\sin(2\xi_j+2\Upsilon_j^1-2\Upsilon_j^2)|}.
	\end{equation} 
	
	Next, we will show another estimate on the phase $\mathfrak{B}^c_j$. Let $c\in(0,1)$ be an arbitrary constant and suppose that 
	\begin{equation}\label{Lemm:Bessel2:2bb5:a:14}
		|\Xi_j^3 \cos(2(\xi_j+\eta_j+\Upsilon_j^1)) \ - \ (1-\Xi_j^3) \cos(2(\eta_j+\Upsilon_j^2))|\le \frac{c}{2}{|\sin(2\xi_j+2\Upsilon_j^1-2\Upsilon_j^2)|},  
	\end{equation}
	for some $\eta_j$ in $[-\pi,\pi]$, we will show that
	\begin{equation}\begin{aligned}\label{Lemm:Bessel2:2bb5:b}
			|	\partial_{\eta_j}\mathfrak{B}^c_j(\eta_j)| \ge &\	{(1-c) |\sin(2\xi_j+2\Upsilon_j^1-2\Upsilon_j^2)|}>0.
		\end{aligned}
	\end{equation}
	We compute, 
	
	\begin{equation}\begin{aligned}\label{Lemm:Bessel2:2bb5:a:15}
			& |\Xi_j^3\sin(2(\xi_j+\eta_j+\Upsilon_j^1))\ - \ (1-\Xi_j^3) \sin(2(\eta_j+\Upsilon_j^2))||\cos(2\eta_j+2\Upsilon_j^2)|\\
			=\ &	|\Xi_j^3\cos(2\eta_j+2\Upsilon_j^2)\sin(2(\xi_j+\eta_j+\Upsilon_j^1)) -(1-\Xi_j^3)\cos(2\eta_j+2\Upsilon_j^2)\sin(2(\eta_j+\Upsilon_j^2))|\\
			\ge \ &	|\Xi_j^3\cos(2\eta_j+2\Upsilon_j^2)\sin(2(\xi_j+\eta_j+\Upsilon_j^1))-\Xi_j^3 \cos(2(\xi_j+\eta_j+\Upsilon_j^1))\sin(2(\eta_j+\Upsilon_j^2))|\\
			& -|(1-\Xi_j^3)\cos(2\eta_j+2\Upsilon_j^2)-\Xi_j^3 \cos(2(\xi_j+\eta_j+\Upsilon_j^1))||\sin(2(\eta_j+\Upsilon_j^2))|\\
			\ge \ &	{\Xi_j^3 |\sin(2\xi_j+2\Upsilon_j^1-2\Upsilon_j^2)|}-\frac{c|\sin(2(\eta_j+\Upsilon_j^2))|}{2}{|\sin(2\xi_j+2\Upsilon_j^1-2\Upsilon_j^2)|},
		\end{aligned}
	\end{equation}
	which yields
	\begin{equation}\begin{aligned}\label{Lemm:Bessel2:2bb5:a:16}
			& |\Xi_j^3\sin(2(\xi_j+\eta_j+\Upsilon_j^1))-  (1-\Xi_j^3) \sin(2(\eta_j+\Upsilon_j^2))|\\
			\ge &\	{\Xi_j^3 |\sin(2\xi_j+2\Upsilon_j^1-2\Upsilon_j^2)|}-\frac{c|\sin(2(\eta_j+\Upsilon_j^2))|}{2}{|\sin(2\xi_j+2\Upsilon_j^1-2\Upsilon_j^2)|}
			.\end{aligned}
	\end{equation}
	A similar argument also gives
	\begin{equation}\begin{aligned}\label{Lemm:Bessel2:2bb5:a:17}
			& |\Xi_j^3\sin(2(\xi_j+\eta_j+\Upsilon_j^1))-  (1-\Xi_j^3) \sin(2(\eta_j+\Upsilon_j^2))|\\
			\ge &\	{(1-\Xi_j^3) |\sin(2\xi_j+2\Upsilon_j^1-2\Upsilon_j^2)|}-\frac{c|\sin(2(\eta_j+\Upsilon_j^2))|}{2}{|\sin(2\xi_j+2\Upsilon_j^1-2\Upsilon_j^2)|}
			.\end{aligned}
	\end{equation}
	Combining \eqref{Lemm:Bessel2:2bb5:a:16} and \eqref{Lemm:Bessel2:2bb5:a:17} yields
	\begin{equation}\begin{aligned}\label{Lemm:Bessel2:2bb5:a:18}
			& |\Xi_j^3\sin(2(\xi_j+\eta_j+\Upsilon_j^1))-  (1-\Xi_j^3) \sin(2(\eta_j+\Upsilon_j^2))|\\
			\ge &\	{\frac12 |\sin(2\xi_j+2\Upsilon_j^1-2\Upsilon_j^2)|}-\frac{c|\sin(2(\eta_j+\Upsilon_j^2))|}{2}{|\sin(2\xi_j+2\Upsilon_j^1-2\Upsilon_j^2)|}\\
			\ge &\	{\frac{1-c}{2} |\sin(2\xi_j+2\Upsilon_j^1-2\Upsilon_j^2)|},
		\end{aligned}
	\end{equation}
	which implies  \eqref{Lemm:Bessel2:2bb5:b}.

	\medskip	
	{\bf Step 1.2: Splitting  $\mathfrak{A}_j^a$. }

	Let $\eta_j^*$ be a point in $\mathfrak{J}_{\eta_j}$, we then write
	\begin{equation}\label{Lemm:Bessel2:2bb5:1}\begin{aligned}
			\mathfrak{B}^c_j(\eta_j)\ = & \ \mathfrak{B}^c_j(\eta_j^*)\ + \ \int_{0}^1\mathrm{d}s \frac{\partial\mathfrak{B}^c_j(s(\eta_j-\eta_j^*)+\eta_j^*) }{\partial s}\\
			\ = & \ \mathfrak{B}^c_j(\eta_j^*)\ + \ (\eta_j-\eta_j^*)\int_{0}^1\mathrm{d}s \partial_{\eta_j}{\mathfrak{B}^c_j(s(\eta_j-\eta_j^*)+\eta_j^*)},
		\end{aligned}
	\end{equation}
	and
	\begin{equation}\label{Lemm:Bessel2:2bb5:2}\begin{aligned}
			\partial_{\eta_j}\mathfrak{B}^c_j(s(\eta_j-\eta_j^*)+\eta_j^*)\ = \ & \partial_{\eta_j}\mathfrak{B}^c_j(\eta_j^*)\ + \ \int_{0}^s\mathrm{d}s' \frac{\partial\mathfrak{B}^c_j(s'(\eta_j-\eta_j^*)+\eta_j^*) }{\partial s'}\\
			\ =\ &  (\eta_j-\eta_j^*)\int_{0}^s\mathrm{d}s' \partial_{\eta_j\eta_j}{\mathfrak{B}^c_j(s'(\eta_j-\eta_j^*)+\eta_j^*) },\end{aligned}
	\end{equation}
	which imply
	\begin{equation}\label{Lemm:Bessel2:2bb5:3}
		\mathfrak{B}^c_j(\eta_j)\ =\ \mathfrak{B}^c_j(\eta_j^*)\ + \ (\eta_j-\eta_j^*)^2\int_{0}^1\int_{0}^s\mathrm{d}s\mathrm{d}s' \partial_{\eta_j\eta_j}{\mathfrak{B}^c_j(s'(\eta_j-\eta_j^*)+\eta_j^*) }.
	\end{equation}
	As $|	\partial_{\eta_j\eta_j}\mathfrak{B}^c_j(\eta_j^*)|
	\ge	2|\sin(2\xi_j+2\Upsilon_j^1-2\Upsilon_j^2)|$ we set $\partial_{\eta_j\eta_j}\mathfrak{B}^c_j(\eta_j^*)=\sigma_{\eta_j^*}|\partial_{\eta_j\eta_j}\mathfrak{B}^c_j(\eta_j^*)|$, where $\sigma_{\eta_j^*}$ is either $1$ or $-1$. There exist constants $c',\delta_{\eta_j^*},\delta_{\eta_j^*}'>0$ such that for  all $\eta_j\in[-\pi,\pi]\cap (\eta_j^*-\delta_{\eta_j^*},\eta_j^*+\delta_{\eta_j^*}')$, we have $\partial_{\eta_j\eta_j}\mathfrak{B}^c_j(s'(\eta_j-\eta_j^*)+\eta_j^*)=\sigma_{\eta_j^*}|\partial_{\eta_j\eta_j}\mathfrak{B}^c_j(s'(\eta_j-\eta_j^*)+\eta_j^*)|$ and
	\begin{equation}\label{Lemm:Bessel2:2bb5:3:1}
		|	\partial_{\eta_j\eta_j}\mathfrak{B}^c_j(\eta_j)|
		\ge 	{c'}{2}|\sin(2\xi_j+2\Upsilon_j^1-2\Upsilon_j^2)|,
	\end{equation}
	for all $s\in[0,1], s'\in[0,s]$. Moreover, $	|	\partial_{\eta_j\eta_j}\mathfrak{B}^c_j(\eta_j)|=2c'|\sin(2\xi_j+2\Upsilon_j^1-2\Upsilon_j^2)|$ when $\eta_j\in\{\eta_j^*-\delta_{\eta_j^*},\eta_j^*+\delta_{\eta_j^*}'\}$.

	We find
	\begin{equation}\label{Lemm:Bessel2:2bb5:4}\begin{aligned}
			\mathfrak{B}^c_j(\eta_j)\ = \ & \mathfrak{B}^c_j(\eta_j^*)\ + \ \sigma_{\eta_j^*}(\eta_j-\eta_j^*)^2\int_{0}^1\int_{0}^s\mathrm{d}s\mathrm{d}s' |\partial_{\eta_j\eta_j}{\mathfrak{B}^c_j(s'(\eta_j-\eta_j^*)+\eta_j^*) }|\\
			\ = \ & \mathfrak{B}^c_j(\eta_j^*)\ + \ \sigma_{\eta_j^*}(\eta_j-\eta_j^*)^2|\mathfrak{G}_{\eta_j^*}(\eta_j-\eta_j^*)|^2,		\end{aligned}
	\end{equation}
	where $\mathfrak{G}_{\eta_j^*}(\eta_j-\eta_j^*)=\sqrt{\sigma_{\eta_j^*}\int_{0}^1\int_{0}^s\mathrm{d}s\mathrm{d}s'\partial_{\eta_j\eta_j}{\mathfrak{B}^c_j(s'(\eta_j-\eta_j^*)+\eta_j^*) }}>0$ is a smooth function with $\eta_j\in\mathbb{T}_{\eta_j^*}:=[-\pi,\pi]\cap (\eta_j^*-\delta_{\eta_j^*},\eta_j^*+\delta_{\eta_j^*}')$. We then define a new variable $y_{\eta_j*}= (\eta_j-\eta_j^*)\mathfrak{G}_{\eta_j^*}(\eta_j-\eta_j^*)$. When $c'$ is close to $1$, according to the inverse function theorem, there exists a neighborhood  $U_{\eta_j^*}$ of the origin $0$ and a smooth function $\psi_{\eta_j^*}: C^\infty_c(U_{\eta_j^*})\to \mathbb{T}_{\eta_j^*}$ such that $\psi_{\eta_j^*}(y_{\eta_j^*})=\eta_j$ and the function $\psi_{\eta_j^*}$ is bijective. The function $ \psi_{\eta_j^*}$ is  the inverse of $(\eta_j-\eta_j^*)\mathfrak{G}_{\eta_j^*}(\eta_j-\eta_j^*)$.
	
	We now consider $\eta_j\in 	\mathbb{T}_{\eta_j}':=[-\pi,\pi]\backslash\cup_{\eta_j^*\in\mathfrak{J}_{\eta_j}}	\big([-\pi,\pi]\cap (\eta_j^*-\delta_{\eta_j^*},\eta_j^*+\delta_{\eta_j^*}')\big)$. As $\partial_{\eta_j}\mathfrak{B}^c_j(\eta_j)\ne 0$ for $\eta_j\in 	\mathbb{T}_{\eta_j}'$, the function $\mathfrak{B}^c_j(\eta_j)$ is monotone on any interval $[\alpha',\beta']\subset \mathbb{T}_{\eta_j}'$. Therefore, $\partial_{\eta_j\eta_j}\mathfrak{B}^c_j(\eta_j)$ is also monotone on any interval $[\alpha',\beta']\subset \mathbb{T}_{\eta_j}'$, by \eqref{Lemm:Bessel2:2aa2:2}. Since  $	|	\partial_{\eta_j\eta_j}\mathfrak{B}^c_j(\eta_j)|={c'}{2}|\sin(2\xi_j+2\Upsilon_j^1-2\Upsilon_j^2)|$ when $\eta_j\in\{\eta_j^*-\delta_{\eta_j^*},\eta_j^*+\delta_{\eta_j^*}'\}$, we deduce that $|	\partial_{\eta_j\eta_j}\mathfrak{B}^c_j(\eta_j)|\le {c'}{2}|\sin(2\xi_j+2\Upsilon_j^1-2\Upsilon_j^2)|$ for any $\eta_j\in 	\mathbb{T}_{\eta_j}'$. Applying \eqref{Lemm:Bessel2:2bb5} and \eqref{Lemm:Bessel2:2bb5:b} for $c=c'$, we find
	\begin{equation}\begin{aligned}\label{Lemm:Bessel2:2bb5:5}
			|	\partial_{\eta_j}\mathfrak{B}^c_j(\eta_j)|> &\	{(1-c') |\sin(2\xi_j+2\Upsilon_j^1-2\Upsilon_j^2)|}.
		\end{aligned}
	\end{equation}
	
	We now split 
	\begin{equation}\begin{aligned}\label{Lemm:Bessel2:2bb5:6}
			\mathfrak{A}_j^a \ = \ &  \sum_{\eta_j^*\in\mathfrak{J}_{\eta_j}}\int_{\mathbb{T}_{\eta_j^*}}\mathrm{d}\eta_j e^{{\bf i}\mathfrak{B}^a_j} \ + \ \int_{\mathbb{T}_{\eta_j}'}\mathrm{d}\eta_j e^{{\bf i}\mathfrak{B}^a_j} \ = \ \sum_{\eta_j^*\in\mathfrak{J}_{\eta_j}}\mathfrak{A}_{\mathbb{T}_{\eta_j}^*} \ + \ \mathfrak{A}_{j,o}^a.
		\end{aligned}
	\end{equation}

	{\bf Step 1.3: Stationary phase estimates of  $\mathfrak{A}_{\mathbb{T}_{\eta_j}^*}$. }

	We employ the change of variables  $\eta_j\to y_{\eta_j^*}$  to rewrite $\mathfrak{A}_{\mathbb{T}_{\eta_j}^*}$ as

	\begin{equation}\begin{aligned}\label{Lemm:Bessel2:2bb5:7}
			\mathfrak{A}_{\mathbb{T}_{\eta_j}^*} \ = \ &     
			\ = \   e^{{\bf i}\mathfrak{B}^a_j(\eta_j^*)}\int_{U_{\eta_j}}\mathrm{d}y_{\eta_j^*} e^{{\bf i}\sigma_{\eta_j^*}\Xi_j^0|y_{\eta_j^*}|^2}\psi_{\eta_j^*}'(y_{\eta_j^*}).
		\end{aligned}
	\end{equation}
	
	By Plancherel's theorem, we can write 
	\begin{equation}
		\label{Lemm:Bessel2:2bb5:8}
		\begin{aligned}
			& \int_{U_{\eta_j}}\mathrm{d}y_{\eta_j^*} e^{{\bf i}\sigma_{\eta_j^*}\Xi_j^0|y_{\eta_j^*}|^2}\psi_{\eta_j^*}'(y_{\eta_j^*}) \ = \ \frac{e^{{\bf i}\pi^2{\sigma_{\eta_j^*}}/4}}{\sqrt{\Xi_j^0}}\int_{\mathbb{R}}\mathrm{d}k e^{-\frac{{\bf i} \pi^2}{\sigma_{\eta_j^*} \Xi_j^0}k^2}\widehat{\psi_{\eta_j^*}'}(k),
		\end{aligned}
	\end{equation}
	where $\widehat{\psi_{\eta_j^*}'}$ is the  Fourier transform of $\psi_{\eta_j^*}'$.
	
	By Taylor's theorem, 
	\begin{equation}
		\label{Lemm:Bessel2:2bb5:9}
		\begin{aligned}
			e^{-\frac{{\bf i} \pi^2}{\sigma_{\eta_j^*} \Xi_j^0}k^2}
			&\ = \ &1+ \sum_{m=1}^\infty\Big(\frac{-{\bf i}\pi^2|k|^2}{\Xi_j^0\sigma_{\eta_j^*}}\Big)^m,
		\end{aligned}
	\end{equation}
	uniformly in $k,|\sigma_{\eta_j^*}|$ and $\Xi_j^0$. We can rewrite \eqref{Lemm:Bessel2:2bb5:8} as
	\begin{equation}
		\label{Lemm:Bessel2:2bb5:10}
		\begin{aligned}
			&\int_{U_{\eta_j}}\mathrm{d}y_{\eta_j^*} e^{{\bf i}\sigma_{\eta_j^*}\Xi_j^0|y_{\eta_j^*}|^2}\psi_{\eta_j^*}'(y_{\eta_j^*})\   =\ \frac{e^{{\bf i}\pi^2{\sigma_{\eta_j^*}}/4}}{\sqrt{\Xi_j^0}}\int_{\mathbb{R}}\mathrm{d}k\Big[1\ + \ \sum_{m=1}^\infty\Big(\frac{-{\bf i}\pi^2|k|^2}{\Xi_j^0\sigma_{\eta_j^*}}\Big)^n\Big]\widehat{\psi_{\eta_j^*}'}(k)\\
			& \ = \ \frac{e^{{\bf i}\pi^2{\sigma_{\eta_j^*}}/4}}{\sqrt{\Xi_j^0}}\int_{\mathbb{R}}\mathrm{d}k\widehat{\psi_{\eta_j^*}'}(k) \ +\ \frac{e^{{\bf i}\pi^2{\sigma_{\eta_j^*}}/4}}{\sqrt{\Xi_j^0}}\int_{\mathbb{R}}\mathrm{d}k \sum_{m=1}^\infty\Big(\frac{-{\bf i}\pi^2k^2}{\Xi_j^0\sigma_{\eta_j^*}}\Big)^m\widehat{\psi_{\eta_j^*}'}(k)\\
			& \ = \ \frac{e^{{\bf i}\pi^2{\sigma_{\eta_j^*}}/4}}{\sqrt{\Xi_j^0}}\int_{\mathbb{R}}\mathrm{d}k\widehat{\psi_{\eta_j^*}'}(k) \ +\ \sum_{m=1}^\infty\frac{e^{{\bf i}\pi^2{\sigma_{\eta_j^*}}/4}}{\sqrt{\Xi_j^0}}\Big(\frac{{\bf i}\pi^2}{\Xi_j^0\sigma_{\eta_j^*}}\Big)^m\int_{\mathbb{R}}\mathrm{d}k \widehat{\psi_{\eta_j^*}^{(2m+1)}}(k)\\
			& \ = \ \frac{e^{{\bf i}\pi^2{\sigma_{\eta_j^*}}/4}}{\sqrt{\Xi_j^0}}{\psi_{\eta_j^*}'}(\eta_j^*) \ +\ \sum_{m=1}^\infty\frac{e^{{\bf i}\pi^2{\sigma_{\eta_j^*}}/4}}{\sqrt{\Xi_j^0}}\Big(\frac{{\bf i}\pi^2}{\Xi_j^0\sigma_{\eta_j^*}}\Big)^m{\psi_{\eta_j^*}^{(2m+1)}}(\eta_j^*).
		\end{aligned}
	\end{equation}
	Our next step will be to compute explicitly the coefficients ${\psi_{\eta_j^*}'}(\eta_j^*)$ and ${\psi_{\eta_j^*}^{(2m+1)}}(\eta_j^*)$.  We first compute
	\begin{equation}
		\label{Lemm:Bessel2:2bb5:11}
		\begin{aligned}
			\psi_{\eta_j^*}'(\eta_j)\ =\ 	& \frac{1}{\partial_{\eta_j}[(\eta_j-\eta_j^*)\mathfrak{G}_{\eta_j^*}(\eta_j-\eta_j^*)]}\ =\ 	 \frac{1}{[(\eta_j-\eta_j^*)\mathfrak{G}'_{\eta_j^*}(\eta_j-\eta_j^*)+\mathfrak{G}_{\eta_j^*}(\eta_j-\eta_j^*)]},
		\end{aligned}
	\end{equation}
	leading to
	\begin{equation}
		\label{Lemm:Bessel2:2bb5:11:1}
		\begin{aligned}
			\psi_{\eta_j^*}'(\eta_j^*)\ =\ 	& 	 \frac{1}{\mathfrak{G}_{\eta_j^*}(0)}\ =\ \frac{\sqrt2}{\sqrt{|\partial_{\eta_j\eta_j}\mathfrak{B}^c_j(\eta_j^*)|}}.
		\end{aligned}
	\end{equation}
	By Fa\`a di Bruno's formula, we find, for all $n\in\mathbb{N}, n\ge 1$
	\begin{equation}
		\label{Lemm:Bessel2:2bb5:12}
		\begin{aligned}
			\psi_{\eta_j^*}^{(n+1)}(\eta_j)\ = \	&  	 \partial_{\eta_j}^{(n)}\Big[\frac{1}{[(\eta_j-\eta_j^*)\mathfrak{G}'_{\eta_j^*}(\eta_j-\eta_j^*)+\mathfrak{G}_{\eta_j^*}(\eta_j-\eta_j^*)]}\Big] \\
			=\	& \sum_{\substack{1m_1+2m_2+\cdots+nm_n=n\\ m_1,\cdots,m_n\in\mathbb{Z},m_1,\cdots,m_n\ge0}}\frac{n!}{m_1!\cdots m_n!}\\
			&\times\frac{(-1)^{m_1+\cdots+m_n}(m_1+\cdots+m_n)!}{[(\eta_j-\eta_j^*)\mathfrak{G}'_{\eta_j^*}(\eta_j-\eta_j^*)+\mathfrak{G}_{\eta_j^*}(\eta_j-\eta_j^*)]^{m_1+\cdots+m_n+1}}\\
			&\times\prod_{i=1}^n\left(\frac{[(\eta_j-\eta_j^*)\mathfrak{G}^{(i+1)}_{\eta_j^*}(\eta_j-\eta_j^*)+(i+1)\mathfrak{G}^{(i)}_{\eta_j^*}(\eta_j-\eta_j^*)]}{i!}\right)^{m_i},
		\end{aligned}
	\end{equation}	
	which implies
	\begin{equation}
		\label{Lemm:Bessel2:2bb5:13}
		\begin{aligned}
			\psi_{\eta_j^*}^{(n+1)}(\eta_j^*)
			=\	& \sum_{\substack{1m_1+2m_2+\cdots+nm_n=n\\ m_1,\cdots,m_n\in\mathbb{Z},m_1,\cdots,m_n\ge0}}\frac{n!}{m_1!\cdots m_n!}\\
			&\times\frac{(-1)^{m_1+\cdots+m_n}(m_1+\cdots+m_n)!}{[\mathfrak{G}_{\eta_j^*}(0)]^{m_1+\cdots+m_n+1}}\prod_{i=1}^n\left(\frac{[(i+1)\mathfrak{G}^{(i)}_{\eta_j^*}(0)]}{i!}\right)^{m_i}.
		\end{aligned}
	\end{equation}

	We next compute, using again Fa\`a di Bruno's formula
	
	\begin{equation}
		\label{Lemm:Bessel2:2bb5:14}\begin{aligned}
			\mathfrak{G}_{\eta_j^*}^{(i)}(\eta_j-\eta_j^*)\ 
			=\int_{0}^1\int_{0}^s\mathrm{d}s\mathrm{d}s'\ &\sum_{\substack{1n_1+2n_2+\cdots+in_i=i\\n_1,\cdots,n_i\in\mathbb{Z},n_1,\cdots,n_i\ge0}}\frac{i!}{n_1!\cdots n_i!}\\
			&\times\Big(\frac12\Big)\cdots\Big(\frac32-i\Big)\Big(\sigma_{\eta_j^*}\Big)^i\big[|\partial_{\eta_j\eta_j}\mathfrak{B}^c_j(s'(\eta_j-\eta_j^*)+\eta_j^*)|\big]^{\frac12-i}\\
			&\times\prod_{l=1}^i\left(\frac{|\partial_{\eta_j}^{l+2}\mathfrak{B}^c_j(s'(\eta_j-\eta_j^*)+\eta_j^*)|}{l!}\right)^{n_l},
		\end{aligned}
	\end{equation}	
	which yields, for $i\ge 1$
	\begin{equation}
		\label{Lemm:Bessel2:2bb5:15}\begin{aligned}
			\mathfrak{G}_{\eta_j^*}^{(i)}(0)\ 
			=\ &\frac12\sum_{\substack{1n_1+2n_2+\cdots+in_i=i\\n_1,\cdots,n_i\in\mathbb{Z},n_1,\cdots,n_i\ge0}}\frac{i!}{n_1!\cdots n_i!}\\
			&\times\Big(\frac12\Big)\cdots\Big(\frac32-i\Big)\Big(\sigma_{\eta_j^*}\Big)^i\big[|\partial_{\eta_j\eta_j}\mathfrak{B}^c_j(\eta_j^*)|\big]^{\frac12-i}\prod_{l=1}^i\left(\frac{|\partial_{\eta_j}^{l+2}\mathfrak{B}^c_j(\eta_j^*)|}{l!}\right)^{n_l}.
		\end{aligned}
	\end{equation}					
	
	Plugging \eqref{Lemm:Bessel2:2bb5:15} into \eqref{Lemm:Bessel2:2bb5:13}, we obtain
	\begin{equation}
		\label{Lemm:Bessel2:2bb5:16}
		\begin{aligned}
			\psi_{\eta_j^*}^{(n+1)}(\eta_j^*)
			=\	& \sum_{\substack{1m_1+2m_2+\cdots+nm_n=n\\ m_1,\cdots,m_n\in\mathbb{Z},m_1,\cdots,m_n\ge0}}\frac{n!\sqrt2}{m_1!\cdots m_n!}\frac{(-1)^{m_1+\cdots+m_n}(m_1+\cdots+m_n)!}{\Big[\sqrt{|\partial_{\eta_j\eta_j}\mathfrak{B}^c_j(\eta_j^*)|}\Big]^{m_1+\cdots+m_n+1}}\\
			&\times\prod_{i=1}^n\left[\sum_{\substack{1n_1+2n_2+\cdots+in_i=i\\n_1,\cdots,n_i\in\mathbb{Z},n_1,\cdots,n_i\ge0}}\frac{i+1}{n_1!\cdots n_i!}\right.\\
			\				&\left.\times\Big(\frac12\Big)\cdots\Big(\frac32-i\Big)\Big(\sigma_{\eta_j^*}\Big)^i\big[|\partial_{\eta_j\eta_j}\mathfrak{B}^c_j(\eta_j^*)|\big]^{\frac12-i}\prod_{l=1}^i\left(\frac{|\partial_{\eta_j}^{l+2}\mathfrak{B}^c_j(\eta_j^*)|}{l!}\right)^{n_l}\right]^{m_i}.
		\end{aligned}
	\end{equation}	
	From \eqref{Lemm:Bessel2:2aa1:1}, we deduce that $\partial_{\eta_j}^{l+2}\mathfrak{B}^c_j(\eta_j^*)=0$ when $l$ is odd. Moreover, when $l$ is even, $|\partial_{\eta_j}^{l+2}\mathfrak{B}^c_j(\eta_j^*)|=2^l|\partial_{\eta_j\eta_j}\mathfrak{B}^c_j(\eta_j^*)|.$ We claim $\prod_{l=1}^i\Big(\frac{|\partial_{\eta_j}^{l+2}\mathfrak{B}^c_j(\eta_j^*)|}{l!}\Big)^{n_l}=0$ when $i$ is odd. This can be shown as follows. As $i$ is odd and $i=1n_1+2n_2+\cdots+in_i$, there is an odd index $l$ such that $n_l\ne 0$, leading to  $\partial_{\eta_j}^{l+2}\mathfrak{B}^c_j(\eta_j^*)=0$  and hence  $\prod_{l=1}^i\Big(\frac{|\partial_{\eta_j}^{l+2}\mathfrak{B}^c_j(\eta_j^*)|}{l!}\Big)^{n_l}=0$. Moreover, to make sure that this product is not zero, $n_l,m_i$ need to be zero when $l,i$ are odd. 
	When $n$ is odd, $n=2m-1$, $m\in\mathbb{N}, m\ge 1$, since $1m_1+2m_2+\cdots+nm_n=n$, there exists an odd index $i$ such that $m_i$ is odd, and thus $m_i\ne 0$, yielding 
	\begin{equation}
		\label{Lemm:Bessel2:2bb5:17:1}
		\psi_{\eta_j^*}^{(2m)}(\eta_j^*)\ =\ 0, \forall  m\in\mathbb{N}, m\ge 1.
	\end{equation}
	Now, 				
	for $n=2m$, we find
	\begin{equation}
		\label{Lemm:Bessel2:2bb5:17}
		\begin{aligned}
			&\psi_{\eta_j^*}^{(2m+1)}(\eta_j^*)
			=\	 \sum_{\substack{2m_2+4m_4+\cdots+2mm_{2m}=2m\\ m_2,\cdots,m_{2m}\in\mathbb{Z},m_2,\cdots,m_{2m}\ge0}}\frac{(2m)!\sqrt2}{m_2!\cdots m_{2m}!}\\
			&\times\frac{(-1)^{m_2+m_4+\cdots+m_{2m}}(m_2+m_4+\cdots+m_{2m})!}{\Big[\sqrt{|\partial_{\eta_j\eta_j}\mathfrak{B}^c_j(\eta_j^*)|}\Big]^{m_2+m_4+\cdots+m_{2m}+1}}\\
			&\times\prod_{i=2, i \text{ is even}}^{2m}\left[\sum_{\substack{2n_2+4n_4+\cdots+in_i=i\\n_2,n_4,\cdots,n_i\in\mathbb{Z},n_2,n_4,\cdots,n_i\ge0}}\frac{(i+1)}{n_1!\cdots n_i!}\right.\\
			\				&\left.\times\Big(\frac12\Big)\cdots\Big(\frac32-i\Big)\Big(\sigma_{\eta_j^*}\Big)^i\big[|\partial_{\eta_j\eta_j}\mathfrak{B}^c_j(\eta_j^*)|\big]^{\frac12-i}\prod_{l=2, l \text{ is even}}^i\left(\frac{|2^l\partial_{\eta_j\eta_j}\mathfrak{B}^c_j(\eta_j^*)|}{l!}\right)^{n_l}\right]^{m_i}.
		\end{aligned}
	\end{equation}	
	By the H\"older sum inequality, we bound, for $m_i>1$,
	\begin{equation}
		\label{Lemm:Bessel2:2bb5:17:1}
		\begin{aligned}
			&\left|\sum_{\substack{2n_2+4n_4+\cdots+in_i=i\\n_2,n_4,\cdots,n_i\in\mathbb{Z},n_2,n_4,\cdots,n_i\ge0}}\frac{(i+1)}{n_1!\cdots n_i!}\right.\\
			\				&\left.\times\Big(\frac12\Big)\cdots\Big(\frac32-i\Big)\Big(\sigma_{\eta_j^*}\Big)^i\big[|\partial_{\eta_j\eta_j}\mathfrak{B}^c_j(\eta_j^*)|\big]^{\frac12-i}\prod_{l=2, l \text{ is even}}^i\left(\frac{|2^l\partial_{\eta_j\eta_j}\mathfrak{B}^c_j(\eta_j^*)|}{l!}\right)^{n_l}\right|\\
			\le\	&\left[\sum_{\substack{2n_2+4n_4+\cdots+in_i=i\\n_2,n_4,\cdots,n_i\in\mathbb{Z},n_2,n_4,\cdots,n_i\ge0}}\left(\frac{(i+1)}{n_1!\cdots n_i!}\Big|\frac12\Big|\cdots\Big|\frac32-i\Big|\right)^{\frac{m_i}{m_i-1}}\right]^{\frac{m_i-1}{m_i}}\\
			\				&\times\left[\sum_{\substack{2n_2+4n_4+\cdots+in_i=i\\n_2,n_4,\cdots,n_i\in\mathbb{Z},n_2,n_4,\cdots,n_i\ge0}}\left(\big[|\partial_{\eta_j\eta_j}\mathfrak{B}^c_j(\eta_j^*)|\big]^{\frac12-i}\prod_{l=2, l \text{ is even}}^i\left(\frac{|2^l\partial_{\eta_j\eta_j}\mathfrak{B}^c_j(\eta_j^*)|}{l!}\right)^{n_l}\right)^{m_i}\right]^{\frac{1}{m_i}}.
		\end{aligned}
	\end{equation}	
	
	Plugging \eqref{Lemm:Bessel2:2bb5:17:1} into \eqref{Lemm:Bessel2:2bb5:17}, we find
	\begin{equation}
		\label{Lemm:Bessel2:2bb5:18}
		\begin{aligned}
			&|\psi_{\eta_j^*}^{(2m+1)}(\eta_j^*)|
			\le\	 \left|\sum_{\substack{2m_2+4m_4+\cdots+2mm_{2m}=2m\\ m_2,\cdots,m_{2m}\in\mathbb{Z},m_2,\cdots,m_{2m}\ge0}}\frac{(2m)!(m_2+m_4+\cdots+m_{2m})!}{m_2!\cdots m_{2m}!}\right.\\
			&\times\left.\prod_{i=2, i \text{ is even}}^{2m}\left[\sum_{\substack{2n_2+4n_4+\cdots+in_i=i\\n_2,n_4,\cdots,n_i\in\mathbb{Z},n_2,n_4,\cdots,n_i\ge0}}\frac{(i+1)}{n_1!\cdots n_i!}\Big|\frac12\Big|\cdots\Big|\frac32-i\Big|\prod_{l=2, l \text{ is even}}^i\left(\frac{2^l}{l!}\right)^{n_l}\right]^{m_i}\right|\\
			&\times|\mathfrak{C}_{\eta_j^*}^1|^m \left|\sum_{\substack{2m_2+4m_4+\cdots+2mm_{2m}=2m\\ m_2,\cdots,m_{2m}\in\mathbb{Z},m_2,\cdots,m_{2m}\ge0}}\prod_{i=2, i \text{ is even}}^{2m}\left(\sum_{\substack{2n_2+4n_4+\cdots+in_i=i\\n_2,n_4,\cdots,n_i\in\mathbb{Z},n_2,n_4,\cdots,n_i\ge0}}\big[|\partial_{\eta_j\eta_j}\mathfrak{B}^c_j(\eta_j^*)|\big]^{(\frac12-i)m_i}\right.\right.\\
			&\left.\left.\  \times{\prod_{l=2, l \text{ is even}}^i}|\partial_{\eta_j\eta_j}\mathfrak{B}^c_j(\eta_j^*)|^{n_lm_i}\Big[\sqrt{|\partial_{\eta_j\eta_j}\mathfrak{B}^c_j(\eta_j^*)|}\Big]^{-m_2-m_4-\cdots-m_{2m}-1}\right)\right|,
		\end{aligned}
	\end{equation}	
	for some constant  $\mathfrak{C}^1_{\eta_j^*}>0$ independent of $m$.
	There exists an explicit constant $\mathfrak{C}^2_{\eta_j^*}>0$ independent of $m$, such that the constant in front of the term containing only $\mathfrak{B}^c_j(\eta_j^*)$  in \eqref{Lemm:Bessel2:2bb5:18} is bounded as
	\begin{equation}
		\label{Lemm:Bessel2:2bb5:19}
		\begin{aligned}
			& \left|\sum_{\substack{2m_2+4m_4+\cdots+2mm_{2m}=2m\\ m_2,\cdots,m_{2m}\in\mathbb{Z},m_2,\cdots,m_{2m}\ge0}}\frac{(2m)!(m_2+m_4+\cdots+m_{2m})!}{m_2!\cdots m_{2m}!}\right|\\
			&\times\prod_{i=2, i \text{ is even}}^{2m}\left[\sum_{\substack{2n_2+4n_4+\cdots+in_i=i\\n_2,n_4,\cdots,n_i\in\mathbb{Z},n_2,n_4,\cdots,n_i\ge0}}\frac{(i+1)}{n_1!\cdots n_i!}\right.\\
			\				&\left.\times\Big|\frac12\Big|\cdots\Big|\frac32-i\Big|\prod_{l=2, l \text{ is even}}^i\left(\frac{2^l}{l!}\right)^{n_l}\right]^{m_i}|\mathfrak{C}_{\eta_j^*}^1|^m\ \le \  |\mathfrak{C}^2_{\eta_j^*}|^{m}.
		\end{aligned}
	\end{equation}	
	We only need to estimate the term containing only $\mathfrak{B}^c_j(\eta_j^*)$  in \eqref{Lemm:Bessel2:2bb5:18}. This quantity can be simplified as
	\begin{equation}
		\label{Lemm:Bessel2:2bb5:20}				\begin{aligned}
			&	\left|\sum_{\substack{2m_2+4m_4+\cdots+2mm_{2m}=2m\\ m_2,\cdots,m_{2m}\in\mathbb{Z},m_2,\cdots,m_{2m}\ge0}}\left[\prod_{i=2, i \text{ is even}}^{2m}\left(\sum_{\substack{2n_2+4n_4+\cdots+in_i=i, i \text{ is even}\\n_2,n_4,\cdots,n_i\in\mathbb{Z},n_2,n_4,\cdots,n_i\ge0}}\big[|\partial_{\eta_j\eta_j}\mathfrak{B}^c_j(\eta_j^*)|\big]^{(\frac12-i)m_i}\right.\right.\right.\\
			&\left.\left.\left.\  \times{\prod_{l=2, l \text{ is even}}^i}|\partial_{\eta_j\eta_j}\mathfrak{B}^c_j(\eta_j^*)|^{n_lm_i}\Big[\sqrt{|\partial_{\eta_j\eta_j}\mathfrak{B}^c_j(\eta_j^*)|}\Big]^{-m_2-m_4-\cdots-m_{2m}-1}\right)\right]\right|\\
			\lesssim \	&|\mathfrak{C}^3_{\eta_j^*}|^{m}\left|\sum_{\substack{2m_2+4m_4+\cdots+2mm_{2m}=2m\\ m_2,\cdots,m_{2m}\in\mathbb{Z},m_2,\cdots,m_{2m}\ge0}}\left(\sum_{\substack{2n_2+4n_4+\cdots+in_i=i\\n_2,n_4,\cdots,n_i\in\mathbb{Z},n_2,n_4,\cdots,n_i\ge0,  i \text{ is even}}}\right.\right.\\
			&\left.\left.\left[\frac{1}{\big|\partial_{\eta_j\eta_j}\mathfrak{B}^c_j(\eta_j^*)\big|}\right]^{\sum_{i=2, i \text{ is even}}^{2m}(-\frac12+i-n_2-\cdots-n_i)m_i+\frac{m_2+\cdots+m_{2m}+1}{2}}\right)\right|,				\end{aligned}
	\end{equation}for some constant  $\mathfrak{C}^3_{\eta_j^*}>0$ independent of $m$.
	We observe that 
	$$\sum_{i=2, i \text{ is even}}^{2m}(-\frac12+i-n_2-\cdots-n_i)m_i+\frac{m_2+\cdots+m_{2m}}{2}=\sum_{i=2, i \text{ is even}}^{2m}(i-n_2-\cdots-n_i)m_i.$$
	
	As $2n_2+4n_4+\cdots+in_i=i$, we deduce $i-n_2-\cdots-n_i\le{i}$. Since $2m_2+4m_4+\cdots+2mm_{2m}=2m$, then $\sum_{i=2, i \text{ is even}}^{2m}(i-n_2-\cdots-n_i)m_i\le \sum_{i=2, i \text{ is even}}^{2m}im_i\le 2m$. Thus, there exists an explicit constant $\mathfrak{C}^4_{\eta_j^*}>0$ such that we could bound \eqref{Lemm:Bessel2:2bb5:20} by 	$(\mathfrak{C}^4_{\eta_j^*})^m\Big[\frac{1}{\big|\partial_{\eta_j\eta_j}\mathfrak{B}^c_j(\eta_j^*)\big|}\Big]^{{2m+1}}$, leading to
	\begin{equation}
		\label{Lemm:Bessel2:2bb5:21}			
		|\psi_{\eta_j^*}^{(2m+1)}(\eta_j^*)|\le 	(\mathfrak{C}^2_{\eta_j^*}\mathfrak{C}^4_{\eta_j^*})^m\left[\frac{1}{\big|\partial_{\eta_j\eta_j}\mathfrak{B}^c_j(\eta_j^*)\big|}\right]^{{2m+1}}.				
	\end{equation}
	We then deduce from \eqref{Lemm:Bessel2:2bb4:3} 
	\begin{equation}
		\label{Lemm:Bessel2:2bb5:22}
		\begin{aligned}
			&\Big|\int_{U_{\eta_j}}\mathrm{d}y_{\eta_j^*} e^{{\bf i}\sigma_{\eta_j^*}\Xi_j^0|y_{\eta_j^*}|^2}\psi_{\eta_j^*}'(y_{\eta_j^*})\Big|\   \le \ \frac{\sqrt2}{\sqrt{\Xi_j^0|\partial_{\eta_j\eta_j}\mathfrak{B}^c_j(\eta_j^*)|}}+\sum_{m=1}^\infty\frac{(\mathfrak{C}^5_{\eta_j^*})^m}{\big|\partial_{\eta_j\eta_j}\mathfrak{B}^c_j(\eta_j^*)\big|^{{2m+1}}|\Xi_j^0|^{\frac{2m+1}{2}}}\\
			&\le \ \frac{1}{\sqrt{\Xi_j^0|\sin(2\xi_j+2\Upsilon_j^1-2\Upsilon_j^2)|}}+\sum_{m=1}^\infty\frac{(\mathfrak{C}^5_{\eta_j^*})^m}{\big[2|\sin(2\xi_j+2\Upsilon_j^1-2\Upsilon_j^2)|\big]^{{2m+1}}|\Xi_j^0|^{\frac{2m+1}{2}}},
		\end{aligned}
	\end{equation}
	in which, the constant $\mathfrak{C}^5_{\eta_j^*}>0$ is universal and explicit.
	By \eqref{Lemm:Bessel2:2bb6:1}, when $\Xi_j^0$ is sufficiently large, we find
	\begin{equation}
		\label{Lemm:Bessel2:2bb5:Final1}				\begin{aligned}
			\Big|\int_{U_{\eta_j}}\mathrm{d}y_{\eta_j^*} e^{{\bf i}\sigma_{\eta_j^*}\Xi_j^0|y_{\eta_j^*}|^2}\psi_{\eta_j^*}'(y_{\eta_j^*})\Big|\   \lesssim & \ \frac{1}{\sqrt{\Xi_j^0|\sin(2\xi_j+2\Upsilon_j^1-2\Upsilon_j^2)|}}\\
			&\ \ 	+\frac{1}{[\Xi_j^0|\sin(2\xi_j+2\Upsilon_j^1-2\Upsilon_j^2)|^2]^\frac32},				\end{aligned}
	\end{equation}
	where the constant on the right hand side is explicit.
	\medskip
	
	{\bf Step 1.4: Stationary phase estimates of  $\mathfrak{A}_{j,o}^a$. }  
	
	We first perform an integration by parts on $\mathfrak{A}_{j,o}^a$
	
	\begin{equation}\begin{aligned}
			\label{Lemm:Bessel2:2bb5:23}
			\mathfrak{A}_{j,o}^a \ = \ & \int_{\mathbb{T}_{\eta_j}'}\mathrm{d}\eta_j e^{{\bf i}\mathfrak{B}^a_j} \ = \ \int_{\mathbb{T}_{\eta_j}'}\mathrm{d}\eta_j \frac{1}{{\bf i}\partial_{\eta_j}\mathfrak{B}^a_j}\partial_{\eta_j}\big(e^{{\bf i}\mathfrak{B}^a_j}\big) \\
			\ = \ & -\frac{1}{\Xi_j^0}\int_{\mathbb{T}_{\eta_j}'}\mathrm{d}\eta_j \partial_{\eta_j}\Big(\frac{1}{{\bf i}\partial_{\eta_j}\mathfrak{B}^c_j}\Big)e^{{\bf i}\Xi_j^0\mathfrak{B}^c_j} \ + \ \frac{1}{\Xi_j^0}\frac{e^{{\bf i}\Xi_j^0\mathfrak{B}^c_j} }{{\bf i}\partial_{\eta_j}\mathfrak{B}^c_j}\Big|_{\partial\mathbb{T}_{\eta_j}'}\\
			\ = \ & -\frac{{\bf i}}{\Xi_j^0}\int_{\mathbb{T}_{\eta_j}'}\mathrm{d}\eta_j \frac{\partial_{\eta_j\eta_j}\mathfrak{B}^c_j}{|\partial_{\eta_j}\mathfrak{B}^c_j|^2}e^{{\bf i}\Xi_j^0\mathfrak{B}^c_j} \ + \ \frac{{\bf i}}{\Xi_j^0}\frac{e^{{\bf i}\Xi_j^0\mathfrak{B}^c_j} }{{\bf i}\partial_{\eta_j}\mathfrak{B}^c_j}\Big|_{\partial\mathbb{T}_{\eta_j}'},
		\end{aligned}
	\end{equation}
	with the notice that we have set $\mathfrak{B}^a_j=\mathfrak{B}^c_j\Xi_j^0.$
	
	Using \eqref{Lemm:Bessel2:2bb5:5}, we bound
	\begin{equation}\begin{aligned}
			\label{Lemm:Bessel2:2bb5:Final2}
			|\mathfrak{A}_{j,o}^a| 
			\ \le \ & \left|\frac{{\bf i}}{\Xi_j^0}\int_{\mathbb{T}_{\eta_j}'}\mathrm{d}\eta_j \frac{\partial_{\eta_j\eta_j}\mathfrak{B}^c_j}{|\partial_{\eta_j}\mathfrak{B}^c_j|^2}e^{{\bf i}\Xi_j^0\mathfrak{B}^c_j} \right| \ + \ \left|\frac{{\bf i}}{\Xi_j^0}\frac{e^{{\bf i}\Xi_j^0\mathfrak{B}^c_j} }{{\bf i}\partial_{\eta_j}\mathfrak{B}^c_j}\Big|_{\partial\mathbb{T}_{\eta_j}'}\right|\\
			\ \lesssim \  & \frac{1}{|\Xi_j^0||\sin(2\xi_j+2\Upsilon_j^1-2\Upsilon_j^2)|^2} \ + \ \frac{1}{|\Xi_j^0||\sin(2\xi_j+2\Upsilon_j^1-2\Upsilon_j^2)|},
		\end{aligned}
	\end{equation}
	where the constant on the right hand side is explicit.
	
	\medskip
	
	{\bf Step 1.5: Final estimate on  $\mathfrak{A}_j^a$. }
	
	Combining \eqref{Lemm:Bessel2:2bb5:6}, \eqref{Lemm:Bessel2:2bb5:Final1} and \eqref{Lemm:Bessel2:2bb5:Final2}, we find
	\begin{equation}\begin{aligned}
			\label{Lemm:Bessel2:2bb5:24}
			|\mathfrak{A}_{j}^a| 
			\ \lesssim \  & \frac{1}{|\Xi_j^0||\sin(2\xi_j+2\Upsilon_j^1-2\Upsilon_j^2)|^2} \ + \ \frac{1}{|\Xi_j^0||\sin(2\xi_j+2\Upsilon_j^1-2\Upsilon_j^2)|}\\
			& \ + \ \frac{1}{\sqrt{\Xi_j^0|\sin(2\xi_j+2\Upsilon_j^1-2\Upsilon_j^2)|}}+\frac{1}{[\Xi_j^0|\sin(2\xi_j+2\Upsilon_j^1-2\Upsilon_j^2)|^2]^\frac32},
		\end{aligned}
	\end{equation}
	which, under the condition \eqref{Lemm:Bessel2:2bb6:1} and the assumption that $\Xi_j^0$ is large, can be estimated as
	\begin{equation}\begin{aligned}
			\label{Lemm:Bessel2:2bb5:25}
			|\mathfrak{A}_{j}^a| 
			\ \lesssim \  & \frac{1}{|\Xi_j^0||\sin(2\xi_j+2\Upsilon_j^1-2\Upsilon_j^2)|^2} \ + \  \frac{1}{\sqrt{\Xi_j^0|\sin(2\xi_j+2\Upsilon_j^1-2\Upsilon_j^2)|}},
		\end{aligned}
	\end{equation}
	where the constant on the right hand side is explicit. This inequality can be combined with \eqref{Lemm:Bessel2:2bb9c:1}-\eqref{Lemm:Bessel2:2bb9c:2}, yielding
	\begin{equation}\begin{aligned}\label{Lemm:Bessel2:2bb5:Final}
			|\mathfrak{A}_{j}^a| \
			\lesssim\ & |\mathfrak{A}_{j}^a|_\infty:= \frac{[1-\cos(\aleph
				_1^j-\aleph
				_2^j)]^{-\frac14}}{\sqrt{|\sin(2\xi_j+2\Upsilon_j^1-2\Upsilon_j^2)|}}\min\Big\{\Big|\frac{t_0}{2}+\frac{t_1}{2}\cos(V_1)e^{{\bf i}2V_j}+\frac{t_2}{2}\cos(W_1)e^{{\bf i}2W_j}\Big|,\\
			&\ \ \ \ \ \ \ \ \ \ \ \ \ \ \ \ \ \ \Big|\frac{t_1}{2}\sin(V_1)e^{{\bf i}2V_j}+\frac{t_2}{2}\sin(W_1)e^{{\bf i}2W_j}\Big|\Big\}^{-\frac12}\\
			&+
			\ \frac{[1-\cos(\aleph
				_1^j-\aleph
				_2^j)]^{-\frac12}}{{|\sin(2\xi_j+2\Upsilon_j^1-2\Upsilon_j^2)|}^2}\min\Big\{\Big|\frac{t_0}{2}+\frac{t_1}{2}\cos(V_1)e^{{\bf i}2V_j}+\frac{t_2}{2}\cos(W_1)e^{{\bf i}2W_j}\Big|,\\
			&\ \ \ \ \ \ \ \ \ \ \ \ \ \ \ \ \ \ \Big|\frac{t_1}{2}\sin(V_1)e^{{\bf i}2V_j}+\frac{t_2}{2}\sin(W_1)e^{{\bf i}2W_j}\Big|^2\Big\}^{-1},
		\end{aligned}
	\end{equation} 
	where the constant on the right hand side is explicit.
	
	\medskip	
	
	{\bf Step 2: Studying the first phase $\mathfrak{B}_1^a$.} 
	
	\smallskip
	{\bf Step 2.1: Estimating $|\partial_{\eta_1}\mathfrak{A}_{j}^a|$}. 		
	We first bound, using \eqref{204new}
	\begin{equation}\begin{aligned}
			& |\partial_{\eta_1}\mathfrak{A}_{j}^a| \ =  \ \Big|\int_{-\pi}^\pi\mathrm{d}\eta_j \partial_{\eta_1}\mathfrak{B}^a_j e^{{\bf i}\mathfrak{B}^a_j}\Big|\ \leq \ \int_{-\pi}^\pi\mathrm{d}\eta_j \Big|\partial_{\eta_1}\mathfrak{B}^a_j\Big|\\
			\le & \ \int_{-\pi}^\pi\mathrm{d}\eta_j \Big|\partial_{\eta_1}\Big[\mathrm{Re}\Big(e^{{\bf i}2\xi_j+{\bf i}2\eta_j}\mathfrak{C}_j^1\Big)\ - \ \mathrm{Re}\Big(e^{{\bf i}2\eta_j}\mathfrak{C}_j^2\Big)\Big]\Big|
			\\
			\le &  \ \int_{-\pi}^\pi\mathrm{d}\eta_j\Big|\Big[\mathrm{Re}\Big\{e^{{\bf i}2\xi_j+{\bf i}2\eta_j}\Big[\Big(\frac{t_0}{2}+\frac{t_1}{2}\cos(V_1)e^{{\bf i}2V_j}+\frac{t_2}{2}\cos(W_1)e^{{\bf i}2W_j}\Big)\cos(\xi_1+\eta_1)\\
			&\ - \Big(\frac{t_1}{2}\sin(V_1)e^{{\bf i}2V_j}+\frac{t_2}{2}\sin(W_1)e^{{\bf i}2W_j}\Big)\sin(\xi_1+\eta_1)\Big]\Big\}-  \mathrm{Re}\Big(e^{{\bf i}2\eta_j}\Big[\Big(\frac{t_0}{2}\\
			&\ +\frac{t_1}{2}\cos(V_1)e^{{\bf i}2V_j}+\frac{t_2}{2}\cos(W_1)e^{{\bf i}2W_j}\Big)\cos(\eta_1)\\
			&\ - \Big(\frac{t_1}{2}\sin(V_1)e^{{\bf i}2V_j}+\frac{t_2}{2}\sin(W_1)e^{{\bf i}2W_j}\Big)\sin(\eta_1)\Big]\Big\}\\
			\ \lesssim  & \ \Big|\frac{t_0}{2}+\frac{t_1}{2}\cos(V_1)e^{{\bf i}2V_j}+\frac{t_2}{2}\cos(W_1)e^{{\bf i}2W_j}\Big|+\Big|\frac{t_1}{2}\sin(V_1)e^{{\bf i}2V_j}+\frac{t_2}{2}\sin(W_1)e^{{\bf i}2W_j}\Big|=:\mathfrak{R}_0^{0,j}.
		\end{aligned}
	\end{equation}
	Therefore, we can set $\partial_{\eta_1}\mathfrak{B}^a_j=\mathfrak{R}_0^{0,j}\mathfrak{D}_{j}^a$ and rewrite  $\partial_{\eta_1}\mathfrak{A}_{j}^a$ as
	$
	\partial_{\eta_1}\mathfrak{A}_{j}^a \ = \ \mathfrak{R}_0^{0,j}\int_{-\pi}^\pi\mathrm{d}\eta_j \mathfrak{D}^a_j e^{{\bf i}\mathfrak{B}^a_j}.
	$
	The same argument of Step 1 can be reused in precisely the same manner to obtain
	$$
	\Big|\int_{-\pi}^\pi\mathrm{d}\eta_j \mathfrak{D}^a_j e^{{\bf i}\mathfrak{B}^a_j}\Big| \
	\lesssim\ |\mathfrak{A}_{j}^a|_\infty,
	$$ yielding 
	\begin{equation}\label{Lemm:Bessel2:2bb5:FinalX}
		|\partial_{\eta_1}\mathfrak{A}_{j}^a| \
		\lesssim\  \mathfrak{R}_0^{0,j}|\mathfrak{A}_{j}^a|_\infty,
	\end{equation} 
	where the constant on the right hand side is explicit.
	
	\medskip
	
	{\bf Step 2.2: The final estimates}.		
	We write, using \eqref{206new}
	\begin{equation}
		\begin{aligned}
			& \mathfrak{B}^a_1(\xi_1,\eta_1) =    {t_0}\Big[\Big(\frac34+\frac{d-1}{2}\Big)\sin(\xi_1+\eta_1)-\frac14\sin(3\xi_1+3\eta_1)\Big]\\
			& +{t_1}\Big[\Big(\frac34+\frac{d-1}{2}\Big)\sin(\xi_1+\eta_1+V_1)-\frac14\sin(3\xi_1+3\eta_1+3V_1)\Big]\\
			&  +{t_2}[\Big(\frac34+\frac{d-1}{2}\Big)\sin(\xi_1+\eta_1+W_1)-\frac14\sin(3\xi_1+3\eta_1+3W_1)\Big]\\
			&-{t_0}\Big[\Big(\frac34+\frac{d-1}{2}\Big)\sin(\eta_1)-\frac14\sin(3\eta_1)\Big]-{t_1}\Big[\Big(\frac34+\frac{d-1}{2}\Big)\sin(\eta_1+V_1)-\frac14\sin(3\eta_1+3V_1)\Big]\\
			& -{t_2}\Big[\Big(\frac34+\frac{d-1}{2}\Big)\sin(\eta_1+W_1)-\frac14\sin(3\eta_1+3W_1)\Big],
		\end{aligned}
	\end{equation}
	which could be expressed as
	\begin{equation}
		\begin{aligned}
			\mathfrak{B}^a_1(\xi_1,\eta_1)\ = & \  \mathrm{Im}\Big[e^{{\bf i}\xi_1+{\bf i}\eta_1}\Big(\frac34+\frac{d-1}{2}\Big)\Big({t_0}+{t_1}e^{{\bf i}V_1}+{t_2}e^{{\bf i}W_1}\Big)\Big]\\
			&\ -  \mathrm{Im}\Big[e^{{\bf i}3\xi_1+{\bf i}3\eta_1}\Big(\frac{t_0}{4}+\frac{t_1}{4}e^{{\bf i}3V_1}+\frac{t_2}{4}e^{{\bf i}3W_1}\Big)\Big]\\
			& \ - \mathrm{Im}\Big[e^{{\bf i}\eta_1}\Big(\frac34+\frac{d-1}{2}\Big)\Big({t_0}+{t_1}e^{{\bf i}V_1}+{t_2}e^{{\bf i}W_1}\Big)\Big]\\
			&\ +  \mathrm{Im}\Big[e^{{\bf i}3\eta_1}\Big(\frac{t_0}{4}+\frac{t_1}{4}e^{{\bf i}3V_1}+\frac{t_2}{4}e^{{\bf i}3W_1}\Big)\Big].
		\end{aligned}
	\end{equation}
	We set
	\begin{equation}
		\begin{aligned}
			\mathfrak{C}^1_1 \ = \mathfrak{C}^1_{1,re} + {\bf i } \mathfrak{C}^1_{1,im} \ =  & \ \Big(\frac34+\frac{d-1}{2}\Big)\Big({t_0}+{t_1}e^{{\bf i}V_1}+{t_2}e^{{\bf i}W_1}\Big),\\ \ \ \ \mathfrak{C}^2_1\ = \mathfrak{C}^2_{1,re} + {\bf i } \mathfrak{C}^2_{1,im}  \ = & \ \frac{t_0}{4}+\frac{t_1}{4}e^{{\bf i}3V_1}+\frac{t_2}{4}e^{{\bf i}3W_1},
		\end{aligned}
	\end{equation}
	and compute the derivative in $\eta_1$ of $\mathfrak{B}_1^a$
	\begin{equation}
		\begin{aligned}
			\partial_{\eta_1}\mathfrak{B}^a_1(\xi_1,\eta_1)\ = & \  \cos(\xi_1+\eta_1)\mathfrak{C}_{1,re}^1-\sin(\xi_1+\eta_1)\mathfrak{C}^1_{1,im}\ - \ 3\cos(3\xi_1+3\eta_1)\mathfrak{C}_{1,re}^2\\
			& +3\sin(3\xi_1+3\eta_1)\mathfrak{C}^2_{1,im}-  \cos(\eta_1)\mathfrak{C}_{1,re}^1+\sin(\eta_1)\mathfrak{C}^1_{1,im}\\
			& \ +3\cos(3\eta_1)\mathfrak{C}_{1,re}^2-3\sin(3\eta_1)\mathfrak{C}^2_{1,im},
		\end{aligned}
	\end{equation}
	leading to
	\begin{equation}
		\begin{aligned}
			\partial_{\eta_1}\mathfrak{B}^a_1(\xi_1,\eta_1)\ = & \  \mathrm{Re}\Big[e^{{\bf i}\xi_1+{\bf i}\eta_1}\mathfrak{C}^1_1\Big] -  \mathrm{Re}\Big[3e^{{\bf i}3\xi_1+{\bf i}3\eta_1}\mathfrak{C}^2_1\Big] \ -\ \mathrm{Re}\Big[e^{{\bf i}\eta_1}\mathfrak{C}^1_1\Big] +  \mathrm{Re}\Big[3e^{{\bf i}3\eta_1}\mathfrak{C}^2_1\Big]
			.
		\end{aligned}
	\end{equation}
	For the sake of simplicity, we set
	\begin{equation}\begin{aligned}
			\mathfrak{C}^1_1\  = & \ \mathfrak{R}^0_1 \ = \ |\mathfrak{R}^0_1|e^{{\bf i}\Upsilon_1^1}, \ \ \ \mathfrak{C}^2_1\  = \ \mathfrak{R}^0_2/3\ = \ e^{{\bf i}3\Upsilon_1^2}|\mathfrak{R}^0_2|/3,\\
			\mathfrak{R}^0_0 \  = & \ |\mathfrak{R}^0_1| + \  |\mathfrak{R}^0_2|, \ \ \ \mathfrak{R}^0_3=\frac{|\mathfrak{R}^0_1|}{\mathfrak{R}^0_0 },\ \ \
			\mathfrak{B}^a_1(\xi_1,\eta_1) \ = \ \mathfrak{R}^0_0\mathfrak{B}^b_1(\xi_1,\eta_1),
		\end{aligned}
	\end{equation}
	and continue to compute the derivative
	\begin{equation}
		\begin{aligned}
			\partial_{\eta_1}\mathfrak{B}^b_1(\xi_1,\eta_1)\ = & \  \mathfrak{R}^0_3\mathrm{Re}\Big[e^{{\bf i}\xi_1+{\bf i}\eta_1+{\bf i}\Upsilon_1^1}\Big] -  (1-\mathfrak{R}^0_3)\mathrm{Re}\Big[e^{{\bf i}3\xi_1+{\bf i}3\eta_1+{\bf i}3\Upsilon_1^2}\Big]\\
			& \ -\mathfrak{R}^0_3\mathrm{Re}\Big[e^{{\bf i}\eta_1+{\bf i}\Upsilon_1^1}\Big] +  (1-\mathfrak{R}^0_3)\mathrm{Re}\Big[e^{{\bf i}3\eta_1+{\bf i}3\Upsilon_1^1}\Big]\\
			\ = & \  \mathfrak{R}^0_3\mathrm{Re}\Big[e^{{\bf i}\eta_1+{\bf i}\Upsilon_1^1}(e^{{\bf i}\xi_1}-1)\Big] -  (1-\mathfrak{R}^0_3)\mathrm{Re}\Big[e^{{\bf i}3\eta_1+{\bf i}3\Upsilon_1^2}(e^{3{\bf i}\xi_1}-1)\Big]\\
			\ = & \  -2\sin(\xi_1')\mathfrak{R}^0_3\mathrm{Im}\Big[e^{{\bf i}\xi_1'+{\bf i}\eta_1+{\bf i}\Upsilon_1^1}\Big] \\
			& \ + 2\sin(\xi_1')[3-4\sin^2(\xi_1')](1-\mathfrak{R}^0_3)\mathrm{Im}\Big[e^{{\bf i}3\xi_1'+{\bf i}3\eta_1+{\bf i}3\Upsilon_1^2}\Big]\\
			\ = & \  -2\sin(\xi_1')\mathfrak{R}^0_3\sin(\eta_1'+\Upsilon_1^1)
			\\
			& \ + 2\sin(\xi_1')[3-4\sin^2(\xi_1')](1-\mathfrak{R}^0_3)\sin(3\eta_1'+3\Upsilon_1^2)
			,
		\end{aligned}
	\end{equation}
	where $\xi_1'=\xi_1/2$, $\eta_1'=\eta_1+\xi_1'$. For the sake of correctness, we restrict the domain of $\xi_1$ to  $$\mathbb{T}_{\xi_1}=\Big\{\xi_1\in[-\pi,\pi] \mbox{ such that }\Big|\sin(\xi_1/2)\Big|>0\Big\}.$$ However, those singular points in $\mathbb{T}\backslash \mathbb{T}_{\xi_1}$ will be eliminated by the integration in $\xi_1$, as we will see later.  
	%
	%
	
	We	restrict the domain of $\eta_1$ to
	$$\mathbb{T}_{\eta_1}=\Big\{\eta_1\in[-\pi,\pi]\mbox{ such that }$$
	$$\Big|\mathfrak{R}^0_3\sin(\eta_1'+\Upsilon_1^1)-[3-4\sin^2(\xi_1')](1-\mathfrak{R}^0_3)\sin(3\eta_1'+3\Upsilon_1^2)\Big|
	>\epsilon_{\eta_1}
	\Big\},$$
	where $\epsilon_{\eta_1}>0$ is a small parameter to be fixed later. 
	
	We will now estimate
	\begin{equation}
		\mathfrak{A}_1^b \ = \ \int_{[-\pi,\pi]\backslash\mathbb{T}_{\eta_1}}\mathrm{d}\eta_1 e^{\mathfrak{R}^0_0{\bf i}\mathfrak{B}^b_1}\prod_{j=2}^d\mathfrak{A}_j^a(\xi_1,\eta_1,\xi_j), \mbox{ and } 
		\mathfrak{A}_1^c \ = \ \int_{\mathbb{T}_{\eta_1}}\mathrm{d}\eta_1 e^{\mathfrak{R}^0_0{\bf i}\mathfrak{B}^b_1}\prod_{j=2}^d\mathfrak{A}_j^a(\xi_1,\eta_1,\xi_j) .
	\end{equation}
	The first quantity can be trivially bounded as
	\begin{equation}\label{Lemm:Bessel2:2bb5:d:1}
		|\mathfrak{A}_1^b| \ \lesssim \ \epsilon_{\eta_1}\prod_{j=2}^d|\mathfrak{A}_j^a|_\infty.
	\end{equation}
	We now develop the second one by simply doing integration by parts
	\begin{equation}\begin{aligned}
			\mathfrak{A}_1^c \ = \ & \int_{\mathbb{T}_{\eta_1}}\mathrm{d}\eta_1 \frac{1}{{\bf i}\mathfrak{R}^0_0\partial_{\eta_1}\mathfrak{B}^b_1}\partial_{\eta_1}\big(e^{\mathfrak{R}^0_0{\bf i}\mathfrak{B}^b_1}\big)\prod_{j=2}^d\mathfrak{A}_j^a(\xi_1,\eta_1,\xi_j) \\
			\ = \  & -\sum_{j'=2}^d\int_{\mathbb{T}_{\eta_1}}\mathrm{d}\eta_1 \frac{e^{\mathfrak{R}^0_0{\bf i}\mathfrak{B}^b_1}}{{\bf i}\mathfrak{R}^0_0\partial_{\eta_1}\mathfrak{B}^b_1}\partial_{\eta_1}\big(\mathfrak{A}_{j'}^a(\xi_1,\eta_1,\xi_j) \big)\prod_{j=2,j\ne j'}^d\mathfrak{A}_j^a(\xi_1,\eta_1,\xi_j)\\
			&- \int_{\mathbb{T}_{\eta_1}}\mathrm{d}\eta_1 \partial_{\eta_1}\Big(\frac{1}{{\bf i}\mathfrak{R}^0_0\partial_{\eta_1}\mathfrak{B}^b_1}\Big)e^{\mathfrak{R}^0_0{\bf i}\mathfrak{B}^b_1}\prod_{j=2}^d\mathfrak{A}_j^a(\xi_1,\eta_1,\xi_j) \\
			&+ \frac{e^{\mathfrak{R}^0_0{\bf i}\mathfrak{B}^b_1}}{{\bf i}\mathfrak{R}^0_0\partial_{\eta_1}\mathfrak{B}^b_1}\prod_{j=2}^d\mathfrak{A}_j^a(\xi_1,\eta_1,\xi_j)\Big|_{\partial \mathbb{T}_{\eta_1}},\end{aligned}
	\end{equation}
	leading to \begin{equation}\begin{aligned}
			\mathfrak{A}_1^c \ = \  & -\sum_{j'=2}^d\int_{\mathbb{T}_{\eta_1}}\mathrm{d}\eta_1 \frac{e^{\mathfrak{R}^0_0{\bf i}\mathfrak{B}^b_1}}{{\bf i}\mathfrak{R}^0_0\partial_{\eta_1}\mathfrak{B}^b_1}\partial_{\eta_1}\big(\mathfrak{A}_{j'}^a(\xi_1,\eta_1,\xi_j) \big)\prod_{j=2,j\ne j'}^d\mathfrak{A}_j^a(\xi_1,\eta_1,\xi_j)\\
			&+ \int_{\mathbb{T}_{\eta_1}}\mathrm{d}\eta_1 \frac{\partial_{\eta_1\eta_1}\mathfrak{B}^b_1}{{\bf i}\mathfrak{R}^0_0|\partial_{\eta_1}\mathfrak{B}^b_1|^2}e^{\mathfrak{R}^0_0{\bf i}\mathfrak{B}^b_1}\prod_{j=2}^d\mathfrak{A}_j^a(\xi_1,\eta_1,\xi_j) \\
			&+ \frac{e^{\mathfrak{R}^0_0{\bf i}\mathfrak{B}^b_1}}{{\bf i}\mathfrak{R}^0_0\partial_{\eta_1}\mathfrak{B}^b_1}\prod_{j=2}^d\mathfrak{A}_j^a(\xi_1,\eta_1,\xi_j)\Big|_{\partial \mathbb{T}_{\eta_1}}, \end{aligned}
	\end{equation}
	which can be bounded as
	\begin{equation}\begin{aligned}
			|\mathfrak{A}_1^c|
			\ \lesssim \  & \sum_{j'=2}^d\prod_{j=2,j\ne j'}^d\Big|\mathfrak{A}_j^a(\xi_1,\eta_1,\xi_j)\Big|_{\infty}\int_{\mathbb{T}_{\eta_1}}\mathrm{d}\eta_1 \frac{1}{\mathfrak{R}^0_0|\partial_{\eta_1}\mathfrak{B}^b_1|}\Big|\partial_{\eta_1}\big(\mathfrak{A}_{j'}^a(\xi_1,\eta_1,\xi_j)\Big| 
			\\
			&+ \prod_{j=2}^d\Big|\mathfrak{A}_j^a(\xi_1,\eta_1,\xi_j)\Big|_{\infty}\int_{\mathbb{T}_{\eta_1}}\mathrm{d}\eta_1 \frac{1}{\mathfrak{R}^0_0|\partial_{\eta_1}\mathfrak{B}^b_1|^2}
			+ \prod_{j=2}^d\Big|\mathfrak{A}_j^a(\xi_1,\eta_1,\xi_j)\Big|_{\infty}\frac{1}{\mathfrak{R}^0_0|\partial_{\eta_1}\mathfrak{B}^b_1|}\Big|_{\partial \mathbb{T}_{\eta_1}}\\
			\ \lesssim \  & \frac{1}{\mathfrak{R}^0_0\epsilon_{\eta_1}}\sum_{j'=2}^d\prod_{j=2,j\ne j'}^d\Big|\mathfrak{A}_j^a(\xi_1,\eta_1,\xi_j)\Big|_{\infty} \int_{\mathbb{T}_{\eta_1}}\mathrm{d}\eta_1\Big|\partial_{\eta_1}\big(\mathfrak{A}_{j'}^a(\xi_1,\eta_1,\xi_j)\Big| 
			\\
			&+ \frac{1}{\mathfrak{R}^0_0\epsilon_{\eta_1}^2}\prod_{j=2}^d\Big|\mathfrak{A}_j^a(\xi_1,\eta_1,\xi_j)\Big|_{\infty}
			+ \frac{1}{\mathfrak{R}^0_0\epsilon_{\eta_1}}\prod_{j=2}^d\Big|\mathfrak{A}_j^a(\xi_1,\eta_1,\xi_j)\Big|_{\infty}\\
			\ \lesssim \  &  \sum_{j=2}^d\Big[\Big|\frac{t_0}{2}+\frac{t_1}{2}\cos(V_1)e^{{\bf i}2V_j}+\frac{t_2}{2}\cos(W_1)e^{{\bf i}2W_j}\Big|+\Big|\frac{t_1}{2}\sin(V_1)e^{{\bf i}2V_j}\\
			&+\frac{t_2}{2}\sin(W_1)e^{{\bf i}2W_j}\Big|\Big]^\frac12\frac{1}{|\mathfrak{R}^0_0\epsilon_{\eta_1}|}\prod_{j=2}^d\Big|\mathfrak{A}_j^a(\xi_1,\eta_1,\xi_j)\Big|_{\infty}\\
			&+ \frac{1}{\mathfrak{R}^0_0\epsilon_{\eta_1}^2}\prod_{j=2}^d\Big|\mathfrak{A}_j^a(\xi_1,\eta_1,\xi_j)\Big|_{\infty}
			+ \frac{1}{\mathfrak{R}^0_0\epsilon_{\eta_1}}\prod_{j=2}^d\Big|\mathfrak{A}_j^a(\xi_1,\eta_1,\xi_j)\Big|_{\infty},
		\end{aligned}
	\end{equation}
	yielding
	\begin{equation}\begin{aligned}\label{Lemm:Bessel2:2bb5:d:2}
			|\mathfrak{A}_1^c|
			\ \lesssim \  &  \prod_{j=2}^d\Big|\mathfrak{A}_j^a(\xi_1,\eta_1,\xi_j)\Big|_{\infty}\left[\frac{1}{\mathfrak{R}^0_0\epsilon_{\eta_1}^2}+\frac{1}{\mathfrak{R}^0_0\epsilon_{\eta_1}}+\sum_{j=2}^d\frac{\mathfrak{R}_0^{0,j}}{\mathfrak{R}^0_0\epsilon_{\eta_1}}\right].
		\end{aligned}
	\end{equation}
	We deduce from \eqref{Lemm:Bessel2:2bb5:d:1} and \eqref{Lemm:Bessel2:2bb5:d:2} that
	\begin{equation}\begin{aligned}\label{Lemm:Bessel2:2bb5:d:3}
			|\mathfrak{A}_1^b+\mathfrak{A}_1^c|
			\ \lesssim \  &  \prod_{j=2}^d\Big|\mathfrak{A}_j^a(\xi_1,\eta_1,\xi_j)\Big|_{\infty}\left[\frac{1}{\mathfrak{R}^0_0\epsilon_{\eta_1}^2}+\frac{1}{\mathfrak{R}^0_0\epsilon_{\eta_1}}+\sum_{j=2}^d\frac{\mathfrak{R}_0^{0,j}}{\mathfrak{R}^0_0\epsilon_{\eta_1}}+\epsilon_{\eta_1}\right].
		\end{aligned}
	\end{equation}
	We balance the quantity $\epsilon_{\eta_1}$ in the above estimate and suppose that $\mathfrak{R}^0_0$ is sufficiently large, we obtain
	\begin{equation}\begin{aligned}\label{Lemm:Bessel2:7:1}
			|\mathfrak{A}_1^b+\mathfrak{A}_1^c| \ = \ &  \sum_{j=2}^d\Big|\int_{[-\pi,\pi]}\mathrm{d}\eta_1 e^{\mathfrak{R}^0_0{\bf i}\mathfrak{B}^b_1}\prod_{j=2}^d\mathfrak{A}_j^a(\xi_1,\eta_1,\xi_j)\Big|\\ 
			\ \lesssim \  &  \Big[\Big|{t_0}+{t_1}\cos(V_1)e^{{\bf i}2V_j}+{t_2}\cos(W_1)e^{{\bf i}2W_j}\Big|+\Big|\frac{t_1}{2}\sin(V_1)e^{{\bf i}2V_j}\\
			&+{t_2}\sin(W_1)e^{{\bf i}2W_j}\Big|\Big]^\frac12\frac{1}{|{\mathfrak{R}^0_0}|^\frac12}\prod_{j=2}^d\Big|\mathfrak{A}_j^a(\xi_1,\eta_1,\xi_j)\Big|_{\infty}.
		\end{aligned}
	\end{equation}
	A straightforward estimate also gives
	\begin{equation}\begin{aligned}\label{Lemm:Bessel2:7:2}
			|\mathfrak{A}_1^b+\mathfrak{A}_1^c| \ \lesssim \ &  \prod_{j=2}^d\Big|\mathfrak{A}_j^a(\xi_1,\eta_1,\xi_j)\Big|_{\infty},
		\end{aligned}
	\end{equation}
	which, in combination with \eqref{Lemm:Bessel2:7:1} leads to
	\begin{equation}\begin{aligned}\label{Lemm:Bessel2:7}
			|\mathfrak{A}_1^b+\mathfrak{A}_1^c|
			\ \lesssim \  &  \min\Big\{1,\sum_{j=2}^d\Big[\Big|\frac{t_0}{2}+\frac{t_1}{2}\cos(V_1)e^{{\bf i}2V_j}+\frac{t_2}{2}\cos(W_1)e^{{\bf i}2W_j}\Big|+\Big|\frac{t_1}{2}\sin(V_1)e^{{\bf i}2V_j}\\
			&+\frac{t_2}{2}\sin(W_1)e^{{\bf i}2W_j}\Big|\Big]^\frac12\frac{1}{|{\mathfrak{R}^0_0}|^\frac12}\Big\}\prod_{j=2}^d\Big|\mathfrak{A}_j^a(\xi_1,\eta_1,\xi_j)\Big|_{\infty}.
		\end{aligned}
	\end{equation}
	In \eqref{Lemm:Bessel2:7}, we will  take the square $ |\mathfrak{A}_j^a|_\infty^2$ and integrate in $\xi_j$, leading to an integral of $\frac{1}{{|\sin(2\xi_j+2\Upsilon_j^1-2\Upsilon_j^2)|^4}}$, which is not integrable. To avoid this, we  restrict the domain of $\xi_j$ to $\mathbb{T}_{\xi_j}$,  defined in \eqref{Lemm:Bessel2:2bb1:1} and $[-\pi,\pi]\backslash \mathbb{T}_{\xi_j}$. On the domain $\mathbb{T}_{\xi_j}$ we have the bound 
	\begin{equation}\begin{aligned}\label{Lemm:Bessel2:2bb11}
			&\Big|\int_{\mathbb{T}_{\xi_j}}\mathrm{d}\xi_j 	|\mathfrak{A}_j^a|_\infty^2\Big| \lesssim\\ 
			&  \Big\{\Big\langle [1-\cos(\aleph^j
			_1-\aleph^j
			_2)]^{\frac{1}{2}}\min\Big\{\Big|\frac{t_0}{2}+\frac{t_1}{2}\cos(V_1)e^{{\bf i}2V_j}+\frac{t_2}{2}\cos(W_1)e^{{\bf i}2W_j}\Big|,\\
			&\ \ \ \ \Big|\frac{t_1}{2}\sin(V_1)e^{{\bf i}2V_j}+\frac{t_2}{2}\sin(W_1)e^{{\bf i}2W_j}\Big|\Big\}\Big\rangle^{-\frac12}\Big\}^2\int_{\mathbb{T}_{\xi_j}}\mathrm{d}\xi_j\frac{1}{|\sin(2\xi_j+2\Upsilon_j^1-2\Upsilon_j^2)|}\\
			&+ \Big\{\Big\langle [1-\cos(\aleph^j
			_1-\aleph^j
			_2)]^{\frac{1}{2}}\min\Big\{\Big|\frac{t_0}{2}+\frac{t_1}{2}\cos(V_1)e^{{\bf i}2V_j}+\frac{t_2}{2}\cos(W_1)e^{{\bf i}2W_j}\Big|,\\
			&\ \ \ \ \Big|\frac{t_1}{2}\sin(V_1)e^{{\bf i}2V_j}+\frac{t_2}{2}\sin(W_1)e^{{\bf i}2W_j}\Big|\Big\}\Big\rangle^{-\frac12}\Big\}^4\int_{\mathbb{T}_{\xi_j}}\mathrm{d}\xi_j\frac{1}{|\sin(2\xi_j+2\Upsilon_j^1-2\Upsilon_j^2)|^4}.
		\end{aligned}
	\end{equation}
	On the other hand, we can see that
	\begin{equation}\begin{aligned}\label{Lemm:Bessel2:2bb11}
			\Big|\int_{[-\pi,\pi]\backslash\mathbb{T}_{\xi_j}}\mathrm{d}\xi_j |\mathfrak{A}_j^a|^2_\infty\Big| \lesssim\ & \epsilon_{\xi_j}.	\end{aligned}
	\end{equation}
	Balancing $\epsilon_{\xi_j}$, under the constraint \eqref{Lemm:Bessel2:2bb6:1}, we  obtain the estimate
	\begin{equation}\begin{aligned}\label{Lemm:Bessel2:2bb12a}
			&\Big|\int_{[-\pi,\pi]}\mathrm{d}\xi_j |\mathfrak{A}_j^a|^2_\infty\Big|\\
			\lesssim\  &  \Big\langle[1-\cos(\aleph^j
			_1-\aleph^j
			_2)]^{\frac{1}{2}}\min\Big\{\Big|\frac{t_0}{2}+\frac{t_1}{2}\cos(V_1)e^{{\bf i}2V_j}+\frac{t_2}{2}\cos(W_1)e^{{\bf i}2W_j}\Big|,\\
			&\Big|\frac{t_1}{2}\sin(V_1)e^{{\bf i}2V_j}+\frac{t_2}{2}\sin(W_1)e^{{\bf i}2W_j}\Big|\Big\}\Big\rangle^{-(\frac12-)}.		\end{aligned}
	\end{equation}
	Combining \eqref{Lemm:Bessel2:2a:b1},  \eqref{Lemm:Bessel2:7}  and \eqref{Lemm:Bessel2:2bb12a} and integrate in $\xi_1$, we finally obtain
	\begin{equation}
		\label{Lemm:Bessel2:2}\begin{aligned}
			& \|\mathfrak{F}(\cdot,t_0,t_1,t_2) \|_{l^4}^4\ = \ \sum_{m\in\mathbb{Z}^d}|\mathfrak{F}(m,t_0,t_1,t_2) |^4\lesssim\min\Big\{1,\sum_{j=2}^d\Big[\Big|{t_0}+{t_1}\cos(V_1)e^{{\bf i}2V_j}\\
			&\ \ \ +{t_2}\cos(W_1)e^{{\bf i}2W_j}\Big|+\Big|{t_1}\sin(V_1)e^{{\bf i}2V_j}+{t_2}\sin(W_1)e^{{\bf i}2W_j}\Big|\Big]\\
			&\times\Big[\Big|{t_0}+{t_1}e^{{\bf i}3V_1}+{t_2}e^{{\bf i}3W_1}\Big|+(2d+1)\Big|{t_0}+{t_1}e^{{\bf i}V_1}+{t_2}e^{{\bf i}W_1}\Big|\Big]^{-1}\Big\}\\
			&\times \prod_{j=2}^{d} \Big\langle[1-|\cos(\aleph^j
			_1-\aleph^j
			_2)|]^{\frac{1}{2}}\min\Big\{\Big|\frac{t_0}{2}+\frac{t_1}{2}\cos(V_1)e^{{\bf i}2V_j}+\frac{t_2}{2}\cos(W_1)e^{{\bf i}2W_j}\Big|,\\
			&\ \ \ \Big|\frac{t_1}{2}\sin(V_1)e^{{\bf i}2V_j}+\frac{t_2}{2}\sin(W_1)e^{{\bf i}2W_j}\Big|\Big\}\Big\rangle^{-(\frac12-)}.\end{aligned}
	\end{equation}
	
	{\bf (ii)  We will now prove estimate \eqref{Lemm:Bessel3:2} for $\tilde{\mathfrak{F}}$}. To this end, we will need finer estimates, in comparison to \eqref{Lemm:Bessel2:2bb9c:1}-\eqref{Lemm:Bessel2:2bb9c:2} by making use of $\check\Psi$. Starting from \eqref{Lemm:Bessel2:2bb9c}, we have 	\begin{equation}
		\begin{aligned}\label{Lemm:Bessel2:2bb9c:1:0}
			& \Big|\cos(\Upsilon_*)\sin(\xi_1+\eta_1)e^{{\bf i}\aleph^j
				_1}+\sin(\Upsilon_*)\cos(\xi_1+\eta_1)e^{{\bf i}\aleph^j
				_2}\Big|^2\\
			\ge\ & \sin^2(\xi_1+\eta_1-\Upsilon_*)|\cos(\aleph^j
			_1-\aleph
			_2^j)|.
		\end{aligned}
	\end{equation}
	Now, as $t_1/t_2=r_*=(1+\epsilon_{r_*})\tilde{r}_l$ (see \eqref{Lemm:Angle:1}), we find $\sin^2(\xi_1+\eta_1-\Upsilon_*)=\sin^2(\xi_1+\eta_1-\Upsilon_*((1+\epsilon_{r_*})\tilde{r}_l))$. By the mean value theorem, we have
	\begin{equation}
		\begin{aligned}\label{Lemm:Bessel2:2bb9c:1:1}
			& \sin^2(\xi_1+\eta_1-\Upsilon_*((1+\epsilon_{r_*})\tilde{r}_l)) - \sin^2(\xi_1+\eta_1-\Upsilon_*(\tilde{r}_l)) \\
			= \ & (r_*- \tilde{r}_l) \frac{\mathrm{d} }{\mathrm{d}r_*}\sin^2(\xi_1+\eta_1-\Upsilon_*((1+\epsilon_{r_*}')\tilde{r}_l)),\end{aligned}
	\end{equation}
	which can be developed as 
	\begin{equation}
		\begin{aligned}\label{Lemm:Bessel2:2bb9c:2:1}
			& \sin^2(\xi_1+\eta_1-\Upsilon_*((1+\epsilon_{r_*})\tilde{r}_l)) - \sin^2(\xi_1+\eta_1-\Upsilon_*(\tilde{r}_l)) \\
			= \ & (r_*- \tilde{r}_l)\frac{\mathrm{d} }{\mathrm{d}r_*}\Big[\sin(\xi_1+\eta_1)\cos(\Upsilon_*((1+\epsilon_{r_*}')\tilde{r}_l))-\cos(\xi_1+\eta_1)\sin(\Upsilon_*((1+\epsilon_{r_*}')\tilde{r}_l))\Big]^2\\
			= \ & \epsilon_{r_*}\tilde{r}_l\Big[\sin^2(\xi_1+\eta_1)\frac{\mathrm{d} }{\mathrm{d}r_*}\cos^2(\Upsilon_*((1+\epsilon_{r_*}')\tilde{r}_l))+\cos^2(\xi_1+\eta_1)\frac{\mathrm{d} }{\mathrm{d}r_*}\sin^2(\Upsilon_*((1+\epsilon_{r_*}')\tilde{r}_l))\\
			& \ - \ \sin(2\xi_1+2\eta_1)\frac{\mathrm{d} }{\mathrm{d}r_*}[\sin(\Upsilon_*((1+\epsilon_{r_*}')\tilde{r}_l))\cos(\Upsilon_*((1+\epsilon_{r_*}')\tilde{r}_l))]\Big].
		\end{aligned}
	\end{equation}
	Applying \eqref{Lemm:Angle:6}-\eqref{Lemm:Angle:7}-\eqref{Lemm:Angle:8} to \eqref{Lemm:Bessel2:2bb9c:2:1}, we find
	\begin{equation}
		\label{Lemm:Bessel2:2bb9c:3}\begin{aligned}
			& \Big|\sin^2(\xi_1+\eta_1-\Upsilon_*((1+\epsilon_{r_*})\tilde{r}_l)) - \sin^2(\xi_1+\eta_1-\Upsilon_*(\tilde{r}_l))\Big|\\
			\lesssim\ & |\epsilon_{r_*}||\tilde{r}_l|\langle\ln\lambda\rangle^{\mathfrak{C}_{\aleph_1^j,\aleph_2^j}^4}\lesssim \epsilon_{r_*}'\langle\ln\lambda\rangle^{-\mathfrak{C}_{\aleph_1^j,\aleph_2^j}^{4'}}.
		\end{aligned}
	\end{equation} for some constants $\mathfrak{C}_{\aleph_1^j,\aleph_2^j}^4,$ $\mathfrak{C}_{\aleph_1^j,\aleph_2^j}^{4'}>0$ and for a suitable choice  $$|\epsilon_{r_*}|=|\epsilon_{r_*}'|\mathcal{O} \Big(\langle\ln\lambda\rangle^{-\mathfrak{C}_{\aleph_1^j,\aleph_2^j}^{4'}-\mathfrak{C}_{\aleph_1^j,\aleph_2^j}^4-\mathfrak{C}_{\aleph_1^j,\aleph_2^j}^{5}}\Big),$$ 
	and thanks to the trivial bound $$|\tilde{r}_l|\lesssim \langle\ln\lambda\rangle^{\mathfrak{C}_{\aleph_1^j,\aleph_2^j}^{5}}, \ \ \ \ \mathfrak{C}_{\aleph_1^j,\aleph_2^j}^{5}>0.$$
	As
	\begin{equation}
		\label{Lemm:Bessel2:2bb9c:4}\begin{aligned}
			& \Big| \sin^2(\xi_1+\eta_1-\Upsilon_*(\tilde{r}_l))\Big| \gtrsim \langle\ln\lambda\rangle^{-\mathfrak{C}_{\aleph_1^j,\aleph_2^j}^{6}},
		\end{aligned}
	\end{equation}
	for some constant $\mathfrak{C}_{\aleph_1^j,\aleph_2^j}^6>0$,
	due to the cut-off function $\tilde\Psi_4$, we deduce
	\begin{equation}
		\label{Lemm:Bessel2:2bb9c:5}
		\sin^2(\xi_1+\eta_1-\Upsilon_*) \gtrsim \langle\ln\lambda\rangle^{-\mathfrak{C}_{\aleph_1^j,\aleph_2^j}^{7}},
	\end{equation}
	for $\epsilon_{r_*}'$ sufficiently small and for some constant $\mathfrak{C}_{\aleph_1^j,\aleph_2^j}^7>0$.

	As a consequence, we bound (see \eqref{Lemm:Bessel2:2bb8})
	\begin{equation}
		\begin{aligned}\label{Lemm:Bessel2:2bb9c:1:1:1}
			&\Xi_j^1
			\ \gtrsim \ 
			\Big[\Big|\frac{t_0}{2}+\frac{t_1}{2}\cos(V_1)e^{{\bf i}2V_j}+\frac{t_2}{2}\cos(W_1)e^{{\bf i}2W_j}\Big|^2\\
			&+\Big|\frac{t_1}{2}\sin(V_1)e^{{\bf i}2V_j}+\frac{t_2}{2}\sin(W_1)e^{{\bf i}2W_j}\Big|^2\Big]^\frac12 |\cos(\aleph^j
			_1-\aleph^j
			_2)|^\frac12  \langle\ln\lambda\rangle^{\mathfrak{C}_{\aleph_1^j,\aleph_2^j}^{8}},
		\end{aligned}
	\end{equation}
	and
	\begin{equation}
		\begin{aligned}\label{Lemm:Bessel2:2bb9c:1:1:1:1}
			&\Xi_j^2
			\ \gtrsim \ 
			\Big[\Big|\frac{t_0}{2}+\frac{t_1}{2}\cos(V_1)e^{{\bf i}2V_j}+\frac{t_2}{2}\cos(W_1)e^{{\bf i}2W_j}\Big|^2\\
			&+\Big|\frac{t_1}{2}\sin(V_1)e^{{\bf i}2V_j}+\frac{t_2}{2}\sin(W_1)e^{{\bf i}2W_j}\Big|^2\Big]^\frac12 |\cos(\aleph^j
			_1-\aleph^j
			_2)|^\frac12   \langle\ln\lambda\rangle^{\mathfrak{C}_{\aleph_1^j,\aleph_2^j}^{8}},
		\end{aligned}
	\end{equation}
	for some constant $\mathfrak{C}_{\aleph_1^j,\aleph_2^j}^8>0$, 
	and we use them to replace \eqref{Lemm:Bessel2:2bb9c:1}-\eqref{Lemm:Bessel2:2bb9c:2}. The same argument as above gives 
	\begin{equation}
		\label{Lemm:Bessel2:2:1}\begin{aligned}
			& \|\tilde{\mathfrak{F}}(\cdot,t_0,t_1,t_2) \|_{l^4}^4\ = \ \sum_{m\in\mathbb{Z}^d}|\tilde{\mathfrak{F}}(m,t_0,t_1,t_2) |^4\\ \lesssim&\ \langle\ln\lambda\rangle^{\mathfrak{C}_{\tilde{\mathfrak{F}},4}}\prod_{j=2}^{d} \Big\langle|\cos(\aleph^j
			_1-\aleph^j
			_2)|^{\frac{1}{2}}\Big\{\Big|\frac{t_0}{2}+\frac{t_1}{2}\cos(V_1)e^{{\bf i}2V_j}+\frac{t_2}{2}\cos(W_1)e^{{\bf i}2W_j}\Big|+\\
			&\ \ \ \Big|\frac{t_1}{2}\sin(V_1)e^{{\bf i}2V_j}+\frac{t_2}{2}\sin(W_1)e^{{\bf i}2W_j}\Big|\Big\}\Big\rangle^{-(\frac12-)},\end{aligned}
	\end{equation}
	for some constant $\mathfrak{C}_{\tilde{\mathfrak{F}}}>0.$
\end{proof}

\begin{lemma}\label{Lemm:Bessel4} Under assumption \eqref{Lemm:Bessel:aa}, there exists a universal constant $\mathfrak{C}_{\mathfrak{F},3}>0$ independent of $t_0,t_1,t_2,\aleph_1,\aleph_2$, such that 
	\begin{equation}
		\label{Lemm:Bessel4:1}\begin{aligned}
			&	\|\mathfrak{F}(\cdot,t_0,t_1,t_2) \|_{l^3} \ \le  \mathfrak{C}_{\mathfrak{F},3}\prod_{j=2}^{d}\Big\langle\min\Big\{\Big|{t_0}+{t_1}\cos(V_1)e^{{\bf i}2V_j}+{t_2}\cos(W_1)e^{{\bf i}2W_j}\Big|,\\
			& \ \Big|{t_1}\sin(V_1)e^{{\bf i}2V_j}+{t_2}\sin(W_1)e^{{\bf i}2W_j}\Big|\Big\}[1-|\cos(\aleph^j
			_1-\aleph^j
			_2)|]^{\frac{1}{2}}\Big\rangle^{-(\frac{1}{12}-)}\\
			& \ \times \min\Big\{1,\sum_{j=2}^d\Big[\Big|{t_0}+{t_1}\cos(V_1)e^{{\bf i}2V_j}  +{t_2}\cos(W_1)e^{{\bf i}2W_j}\Big|\\
			&\ \ \ +\Big|{t_1}\sin(V_1)e^{{\bf i}2V_j}+{t_2}\sin(W_1)e^{{\bf i}2W_j}\Big|\Big]\\
			&\ \ \ \times\Big[\Big|{t_0}+{t_1}e^{{\bf i}3V_1}+{t_2}e^{{\bf i}3W_1}\Big|+(2d+1)\Big|{t_0}+{t_1}e^{{\bf i}V_1}+{t_2}e^{{\bf i}W_1}\Big|\Big]^{-1}\Big\}^\frac16.\end{aligned}
	\end{equation}
	In addition, \eqref{Lemm:Bessel4:1} also holds true for $\tilde{\mathfrak{F}}$.
	Suppose further that \eqref{Lemm:Angle:0}-\eqref{Lemm:Angle:4}-\eqref{Lemm:Angle:5} hold true, and 	\begin{equation}
		\label{Lemm:Bessel4:1:a}r_*=t_1/t_2=(1+\epsilon_{r_*})\tilde{r}_l\end{equation} for $l=1,2,3$ in which $\tilde{r}_l$ are defined in Lemma \ref{Lemm:Angle}. When \begin{equation}
		\label{Lemm:Bessel4:1:b} |\epsilon_{r_*}|=|\epsilon_{r_*}'|\langle
		\ln\lambda\rangle^{-c}\end{equation}  for an explicit constant $c>0$ depending only on the cut-off functions, and $\epsilon_{r_*}'$ is sufficiently small but independent of $\lambda$ and the cut-off functions,  then we have the estimate 
	\begin{equation}
		\label{Lemm:Bessel4:2}\begin{aligned}
			&	\|\tilde{\mathfrak{F}}(\cdot,t_0,t_1,t_2) \|_{l^3} \ \le  \mathfrak{C}_{\mathfrak{F},3}\langle \ln\lambda\rangle^{\mathfrak{C}_{\mathfrak{F},3}'}\prod_{j=2}^{d}\Big\langle\Big\{\Big|{t_0}+{t_1}\cos(V_1)e^{{\bf i}2V_j}+{t_2}\cos(W_1)e^{{\bf i}2W_j}\Big|+\\
			&+ \ \Big|{t_1}\sin(V_1)e^{{\bf i}2V_j}+{t_2}\sin(W_1)e^{{\bf i}2W_j}\Big|\Big\}|\cos(\aleph^j
			_1-\aleph^j
			_2)|^{\frac{1}{2}}\Big\rangle^{-(\frac{1}{12}-)},\end{aligned}
	\end{equation}
	for universal constants $\mathfrak{C}_{\mathfrak{F},3},\mathfrak{C}_{\mathfrak{F},3}'>0.$
\end{lemma}

\begin{proof}
	Observe that $\frac{1}{3} =\frac{\theta}{4} +\frac{1-\theta}{2}$ with $\theta=\frac{2}{3}$. By interpolating between $l_4$ and  $l_2$,  with the  use  of $\theta$,  we  get  \eqref{Lemm:Bessel4:1}.

\end{proof}

\subsection{The role of the cut-off functions on the dispersive estimates} We follow the same notations used in Section \ref{Subsec:DispersiveEstimates}.
Let us consider the dispersion relation
$\omega(k) \ = \ \sin(2\pi k^1)\Big[\sin^2(2\pi k^1) + \cdots + \sin^2(2\pi k^d)\Big]$ considered in \eqref{Nearest}, with its equivalent form    $\omega(\xi) \ = \sin(\xi_1)\Big[\sin^2(\xi_1) + \cdots + \sin^2(\xi_d)\Big]$, defined in \eqref{Nearest1}. We define a similar function as in  \eqref{extendedBesselfunctions1}, where $t_1,t_2$ are set to be zero. However, the  cut-off function $\check\Psi$ is replaces by $\check\Phi$, which is the component of $\check\Psi$ that removes the singular set  \begin{equation*}\label{Setstar}\begin{aligned}\mathfrak{S}'=\Big\{2\pi k=(2\pi k^1,\cdots,2\pi k^d)\in [-\pi,\pi]^d \Big|	&\ k^j=0,\pm\frac14,\pm\frac12, j=1,\cdots,d \Big\}.\end{aligned}\end{equation*}
We write
\begin{equation}
	\label{Improved:extendedBesselfunctions1}
	\mathfrak{F}^{O}(m,t_0) \ = \ \int_{[-\pi,\pi]^d} \mathrm{d}\xi \check{\Phi}e^{{\bf i}m\cdot \xi} e^{{\bf i}t_0\omega(\xi)}\mathfrak{K}(\xi),
\end{equation}
for $m=(m_1,\cdots,m_d)\in\mathbb{Z}^d$. The kernel $\mathfrak{K}(\xi)$ can be either $1$ or $|\sin(\xi_1)|$, with $\xi=(\xi_1,\cdots,\xi_d)\in[-\pi,\pi]^d$. 

We rewrite \eqref{Improved:extendedBesselfunctions1} as follows

\begin{equation}
	\label{Improved:extendedBesselfunctions1a}
	\begin{aligned}
		&\mathfrak{F}^{O}(m,t_0)  =   \int_{-\pi}^\pi\mathrm{d}\xi_1 \mathfrak{K}(\xi)\check{\Phi}(\xi) \exp\Big({{\bf i}m_1\xi_1}\Big)\exp\Big({\bf i}t_0\sin^3(\xi_1)\Big)\\
		&\ \  \times\Big[\prod_{j=2}^d\int_{-\pi}^\pi\mathrm{d}\xi_j \exp\Big({{\bf i}m_j\xi_j}\Big) \exp\Big({\bf i}t_0\sin(\xi_1)\sin^2(\xi_j)\Big)\Big].
	\end{aligned}
\end{equation}
\begin{lemma}\label{Lemm:Improved:Bessel2} There exists a universal constant $\mathfrak{C}_{\mathfrak{F}^{cut},2}$ independent of $t_0$ and $\lambda$, such that 
	\begin{equation}
		\label{Lemm:Improved:Bessel2:1}
		\|\mathfrak{F}^{O}(\cdot,t_0) \|_{l^2} \ = \ \left( \sum_{m\in\mathbb{Z}^d}|\mathfrak{F}^{O}(m,t_0) |^2\right)^\frac12 \ \le \ \mathfrak{C}_{\mathfrak{F}^{cut},2}.
	\end{equation}
\end{lemma}
\begin{proof} 	By the Plancherel theorem, we obtain
	\begin{equation}
		\label{Lemm:Improved:Bessel2:2}
		\begin{aligned}
			\sum_{m\in\mathbb{Z}^d}|\mathfrak{F}^{O}(m,t_0) |^2 
			\ = \ &  \int_{[-\pi,\pi]^d} \mathrm{d}\xi \Big|e^{{\bf i}t_0\omega(\xi)}\mathfrak{K}(\xi) \check{\Phi}(\xi) \Big|^2,
		\end{aligned}
	\end{equation}
	which is a  bounded quantity. The conclusion of the lemma then follows.
\end{proof}
\begin{lemma}\label{Lemm:Improved:Bessel3}
	There exist  universal constants $\mathfrak{C}_{\mathfrak{F}^{O},4},\mathfrak{C}_{\mathfrak{F}^{O},4'}>0$ independent of $t_0$ and $\lambda$, such that
	\begin{equation}
		\label{Lemm:Improved:Bessel3:1}\begin{aligned}
			&\|\mathfrak{F}^{O}(\cdot,t_0) \|_{l^4}  \le |\ln|\lambda||^{\mathfrak{C}_{\mathfrak{F}^{O},4'}}\mathfrak{C}_{\mathfrak{F}^{O},4}\langle t_0\rangle^{-{\big(\frac{
						d-1}{8}-\big)}}.\end{aligned}
	\end{equation}
\end{lemma}
\begin{remark}\label{remark:Improved:Bessel3}
	Differently  from Lemmas \ref{Lemm:Bessel2} and \ref{Lemm:Bessel3}, the constant $|\ln|\lambda||^{\mathfrak{C}_{\mathfrak{F}^{O},4'}}$ depends on $\lambda$ (see also Remark \ref {remark:Bessel2}). However, it produces a factor which is a power of $|\ln|\lambda||$. This factor will be negligible in  estimates of the Feynman diagrams.  
\end{remark}
\begin{proof} We only prove the case when $\mathfrak{K}(\xi)=1$. The case when $\mathfrak{K}(\xi)=|\sin(\xi_1)|$ can be proved by precisely the same argument. The proof is based on the same strategy used in the proof of Lemma \ref{Lemm:Bessel3}.
	We observe that
	\begin{equation}
		\label{Lemm:Improved:Bessel2:1a}
		\begin{aligned}
			|\mathfrak{F}^{O}(m,t_0) |^2
			\ =  \ & \int_{[-\pi,\pi]^d}\mathrm{d}\xi e^{{\bf i}m\cdot\xi} \int_{-\pi}^{\pi}\mathrm{d}\eta_1 \exp\Big({\bf i}t_0\sin^3(\xi_1+\eta_1) -{\bf i}t_0\sin^3(\eta_1)\Big)\\
			& \times\Big[\prod_{j=2}^d\int_{-\pi}^\pi\mathrm{d}\eta_j \check{\Phi}(\eta)\check{\Phi}(\xi+\eta)\exp\Big({\bf i}t_0\sin(\xi_1+\eta_1)\sin^2(\xi_j+\eta_j)\\
			&- {\bf i}t_0\sin(\eta_1)\sin^2(\eta_j)\Big)\Big], 
		\end{aligned}
	\end{equation} 
	and  the Plancherel theorem
	\begin{equation}
		\label{Lemm:Improved:Bessel2:2}\begin{aligned}
			& \|\mathfrak{F}^{O}(\cdot,t_0) \|_{l^4}^4\ = \ \sum_{m\in\mathbb{Z}^d}|\mathfrak{F}^{O}(m,t_0) |^4\\  
			\ = \ & \int_{[-\pi,\pi]^d}\mathrm{d}\xi \Big| \int_{-\pi}^{\pi}\mathrm{d}\eta_1 \exp\Big({\bf i}t_0\sin^3(\xi_1+\eta_1)-{\bf i}t_0\sin^3(\eta_1)\Big)\\
			& \  \times\Big[\prod_{j=2}^d\int_{-\pi}^\pi\mathrm{d}\eta_j \check{\Phi}(\eta)\check{\Phi}(\xi+\eta)\exp\Big({\bf i}t_0\sin(\xi_1+\eta_1)\sin^2(\xi_j+\eta_j)\\
			&\ \ - {\bf i}t_0\sin(\eta_1)\sin^2(\eta_j)\Big]\Big|^2.\end{aligned}
	\end{equation}
	Note that the cut-off function $\check{\Phi}(\xi)$ can be constructed by cut-off functions of each component 
	\begin{equation}
		\check{\Phi}(\xi)=\prod_{j=1}^d \check{F}_1(\xi_j), \ \ \ \  \xi=(\xi_1,\cdot,\xi_d),
	\end{equation}
	where $\check{F}_1(\xi_j)$ are cut-off functions at $\xi_j$.
	We can eliminate some components, for the sake of simplicity,  without affecting the correctness of the lemma. As in the proof of Lemma \ref{Lemm:Bessel3}, we also define for $2\le j\le d$
	\begin{equation}
		\begin{aligned}
			\mathfrak{A}_j(\xi_1,\eta_1,\xi_j) \ = & \ \int_{-\pi}^\pi\mathrm{d}\eta_j \check{F}_1(\xi_j+\eta_j)\check{F}_1(\eta_j) (\check{F}_1(\xi_1+\eta_1))^{\frac{1}{d}}e^{{\bf i}\mathfrak{B}_j(\xi_1,\eta_1,\xi_j,\eta_j)},
		\end{aligned}
	\end{equation}
	in which
	\begin{equation}
		\begin{aligned}
			& \mathfrak{B}_j(\xi_1,\eta_1,\xi_j,\eta_j) =  \ t_0\sin(\xi_1+\eta_1)\sin^2(\xi_j+\eta_j)
			\ - \ t_0\sin(\eta_1)\sin^2(\eta_j),
		\end{aligned}
	\end{equation}
	and
	\begin{equation}
		\begin{aligned}
			& \mathfrak{B}_1(\xi_1,\eta_1) =  \ t_0\sin^3(\xi_1+\eta_1) -t_0\sin^3(\eta_1),
		\end{aligned}
	\end{equation}
	we find
	\begin{equation}
		\label{Lemm:Improved:Bessel2:2a}\begin{aligned}
			\sum_{m\in\mathbb{Z}^d}|\mathfrak{F}^{O}(m,t_0) |^4
			\ = \ & \int_{[-\pi,\pi]^d}\mathrm{d}\xi  \Big|\int_{-\pi}^{\pi}\mathrm{d}\eta_1 \prod_{j=2}^d\mathfrak{A}_j(\xi_1,\eta_1,\xi_j)\\
			&\times (\check{F}_1(\xi_1+\eta_1))^{\frac{1}{d}}\check{F}_1(\eta_1) e^{{\bf i}\mathfrak{B}_1(\xi_1,\eta_1) }\Big|^2.\end{aligned}
	\end{equation}
	The phase $\mathfrak{B}_j$ can be split as the sum of  
	\begin{equation}
		\begin{aligned}
			& \mathfrak{B}^a_j(\xi_1,\eta_1,\xi_j,\eta_j)
			\ =   \ -\ \mathrm{Re}\Big[e^{{\bf i}2\xi_j+{\bf i}2\eta_j}\mathfrak{C}_j^1\Big]\ + \ \mathrm{Re}\Big[e^{{\bf i}2\eta_j}\mathfrak{C}_j^2\Big],
		\end{aligned}
	\end{equation}
	with
	$$\mathfrak{C}_j^1=\frac{t_0}{2}\sin(\xi_1+\eta_1), \mbox{ and } \mathfrak{C}_j^2=\frac{t_0}{2}\sin(\eta_1),$$
	and
	\begin{equation}
		\begin{aligned}
			&\mathfrak{B}^b_j(\xi_1,\eta_1)
			=   \frac{t_0}{2}\sin(\xi_1+\eta_1)    
			-  \frac{t_0}{2}\sin(\eta_1).
		\end{aligned}
	\end{equation}
	The oscillatory integral $\mathfrak{A}_j$ can now be written
	\begin{equation}
		\mathfrak{A}_j \ = \ e^{{\bf i}\mathfrak{B}^b_j }\int_{-\pi}^\pi\mathrm{d}\eta_j \check{F}_1(\xi_j+\eta_j)(\check{F}_1(\xi_1+\eta_1))^{\frac{1}{d}}\check{F}_1(\eta_j) e^{{\bf i}\mathfrak{B}^a_j} \ = \  e^{{\bf i}\mathfrak{B}^b_j }\mathfrak{A}_j^a.
	\end{equation}
	We
	combine the phases ${\mathfrak{B}^b_j }$ and ${\mathfrak{B}_1}$  
	\begin{equation}
		\begin{aligned}
			& \mathfrak{B}^a_1=\mathfrak{B}_1 + \sum_{j=2}^{d}\mathfrak{B}^b_j \  =  \ t_0\Big[\sin^3(\xi_1+\eta_1)+\frac{d-1}{2}\sin(\xi_1+\eta_1)\Big]\\
			& -t_0\Big[\sin^3(\eta_1)+\frac{d-1}{2}\sin(\eta_1)\Big].
		\end{aligned}
	\end{equation}
	As a result, we write
	\begin{equation}
		\label{Lemm:Improved:Bessel2:2a:a}\begin{aligned}
			\sum_{m\in\mathbb{Z}^d}|\mathfrak{F}^{O}(m,t_0) |^4
			\ = \ & \int_{[-\pi,\pi]^d}\mathrm{d}\xi  \Big|\int_{-\pi}^{\pi}\mathrm{d}\eta_1 \prod_{j=2}^d\mathfrak{A}_j^a(\xi_1,\eta_1,\xi_j)\\
			&\times (\check{F}_1(\xi_1+\eta_1))^{\frac{1}{d}}\check{F}_1(\eta_1)e^{{\bf i}\mathfrak{B}_1^a(\xi_1,\eta_1) }\Big|^2.\end{aligned}
	\end{equation}
	{\bf Step 1: Splitting $\mathfrak{A}_j^a$.} 
	Using the stationary phase estimate, we denote by $\eta_j^*$ the solution of $
	\partial_{\eta_j}\mathfrak{B}^a_j  =   0,$
	which is equivalent to
	\begin{equation}\label{Lemm:Improved:Bessel2:2bb1}
		\mathfrak{C}_j^1 \sin(2(\xi_j+\eta_j^*)) \ = \ \mathfrak{C}_j^2 \sin(2(\eta_j^*)).
	\end{equation}
	First, we consider the case $\sin(2(\eta_j^*))=0$. This case can be eliminated due to the cut-off function $\check{F}_1(\eta_j)$. Therefore, $\sin(2(\eta_j^*))\ne 0$ and we compute
	
	\begin{equation}\begin{aligned}\label{Lemm:Improved:Bessel2:2bb2}
			\partial_{\eta_j\eta_j}\mathfrak{B}^a_j(\eta_j^*) \ = \ & 4\mathfrak{C}_j^1 \cos(2(\xi_j+\eta_j^*)) \ - \ 4\mathfrak{C}_j^2 \cos(2(\eta_j^*))\\
			\ = \ &\ 4\frac{\mathfrak{C}_j^1 \cos(2(\xi_j+\eta_j^*))\sin(2(\eta_j^*)) \ - \ \mathfrak{C}_j^2 \cos(2(\eta_j^*))\sin(2(\eta_j^*))}{\sin(2(\eta_j^*))}\\
			\ = \ &\ 4\frac{\mathfrak{C}_j^1 \cos(2(\xi_j+\eta_j^*))\sin(2(\eta_j^*)) \ - \ \mathfrak{C}_j^1 \cos(2(\eta_j^*))\sin(2(\xi_j+\eta_j^*))}{\sin(2(\eta_j^*))}\\
			\ = \ &\ -4\frac{\mathfrak{C}_j^1 \sin(2\xi_j)}{\sin(2(\eta_j^*))} \ = \ -2\frac{{t_0}\sin(\xi_1+\eta_1) \sin(2\xi_j)}{\sin(2(\eta_j^*))}.
		\end{aligned}
	\end{equation}
	We deduce the bound
	\begin{equation}\begin{aligned}\label{Lemm:Improved:Bessel2:2bb2:1:1}
			|\partial_{\eta_j\eta_j}\mathfrak{B}^a_j(\eta_j^*)| \ \ge  \ & 2|{{t_0}\sin(\xi_1+\eta_1) \sin(2\xi_j)}| \ge \mathfrak{C}^2_{cut,j}|\ln(\lambda)|^{-\mathfrak{C}^1_{cut,j}}|\sin(2\xi_j)||t_0|,
		\end{aligned}
	\end{equation}
	for some constants $\mathfrak{C}_{cut,j}^1,\mathfrak{C}^2_{cut,j}>0$, when $\xi_1+\eta_1$  belongs to the support of $\check{F}_1$.

	Similarly as  \eqref{Lemm:Bessel2:2bb1:1}, we also restrict our domain of defining $\xi_j$ to
	\begin{equation}\label{Lemm:Improved:Bessel2:2bb2:1}\mathbb{T}_{\xi_j}=\Big\{\xi_j\in[-\pi,\pi],|2\xi_j  - m\pi|> \epsilon_{\xi_j}>0\Big\},\end{equation} 
	for $m=0,\pm 1, \pm 2$. We then obtain
	\begin{equation}\begin{aligned}\label{Lemm:Improved:Bessel2:2bb2:1}
			|\partial_{\eta_j\eta_j}\mathfrak{B}^a_j(\eta_j^*)| \ \ge  \ & 2|{{t_0}\sin(\xi_1+\eta_1) \sin(2\xi_j)}| \ge \mathfrak{C}^2_{cut,j}|\ln(\lambda)|^{-\mathfrak{C}^1_{cut,j}}|t_0|\epsilon_{\xi_j}.
		\end{aligned}
	\end{equation}
	We also add a constraint on the choice of $\epsilon_{\xi_j}$	\begin{equation}\label{Lemm:Improved:Bessel2:2bb2:1:2}
		1\gg \frac{1}{\mathfrak{C}^2_{cut,j}|\ln(\lambda)|^{-\mathfrak{C}^1_{cut,j}}|t_0|\epsilon_{\xi_j}^2}.
	\end{equation}

	Note that there is a factor of $|\ln|\lambda||^{\mathfrak{C}_{\mathfrak{F}^{O},4'}}$ in front of the constant on the right hand side of \eqref{Lemm:Improved:Bessel3:1}, which is what we want to prove, we only need to consider the case when $\mathfrak{C}^2_{cut,j}|\ln(\lambda)|^{-\mathfrak{C}^1_{cut,j}}|t_0| \gg 1$. Indeed, for smaller values of $t_0$, we always have $\|\mathfrak{F}^{O}(\cdot,t_0) \|_{l^4} \lesssim1$, which automatically yields \eqref{Lemm:Improved:Bessel3:1}.
	
	Denote by $\mathfrak{J}_{\eta_j}$ the set of all stationary points of $\mathfrak{B}_j^a$. By the same argument used in Step 1.1 of the proof of Lemma \ref{Lemm:Bessel3}, we can prove that 
	\begin{equation}\label{Lemm:Improved:Bessel2:2bb5:3:1}
		|\mathfrak{J}_{\eta_j}|\le 3.\end{equation}
	Set $\mathfrak{B}^c_j(\eta_j)=\mathfrak{B}^a_j(\eta_j)/t_0$. Let $\eta_j^*$ be a point in $\mathfrak{J}_{\eta_j}$, we follow equation \eqref{Lemm:Bessel2:2bb5:3} in Step 1.2 of Lemma \ref{Lemm:Bessel3} to write
	\begin{equation}\label{Lemm:Improved:Bessel2:2bb5:3}
		\mathfrak{B}^c_j(\eta_j)\ =\ \mathfrak{B}^c_j(\eta_j^*)\ + \ (\eta_j-\eta_j^*)^2\int_{0}^1\int_{0}^s\mathrm{d}s\mathrm{d}s' \partial_{\eta_j\eta_j}{\mathfrak{B}^c_j(ss'(\eta_j-\eta_j^*)+\eta_j^*) }.
	\end{equation}
	From \eqref{Lemm:Improved:Bessel2:2bb2:1}, $|	\partial_{\eta_j\eta_j}\mathfrak{B}^c_j(\eta_j^*)|
	\ge	2|{\sin(\xi_1+\eta_1) \sin(2\xi_j)}|$, we set $\partial_{\eta_j\eta_j}\mathfrak{B}^c_j(\eta_j^*)=\sigma_{\eta_j^*}|\partial_{\eta_j\eta_j}\mathfrak{B}^c_j(\eta_j^*)|$, where $\sigma_{\eta_j^*}$ is either $1$ or $-1$. Thus, 
	there exist constants $c',\delta_{\eta_j^*},\delta_{\eta_j^*}'>0$ such that for  all $\eta_j\in[-\pi,\pi]\cap (\eta_j^*-\delta_{\eta_j^*},\eta_j^*+\delta_{\eta_j^*}')$, we have $\partial_{\eta_j\eta_j}\mathfrak{B}^c_j(ss'(\eta_j-\eta_j^*)+\eta_j^*)=\sigma_{\eta_j^*}|\partial_{\eta_j\eta_j}\mathfrak{B}^c_j(ss'(\eta_j-\eta_j^*)+\eta_j^*)|$ for all $s\in[0,1], s'\in[0,s]$ and
	\begin{equation}\label{Lemm:Bessel2:Improved:2bb5:3:1}
		|	\partial_{\eta_j\eta_j}\mathfrak{B}^c_j(\eta_j)|
		\ge 	2{c'}|{\sin(\xi_1+\eta_1) \sin(2\xi_j)}|.
	\end{equation}
	In addition, $	|	\partial_{\eta_j\eta_j}\mathfrak{B}^c_j(\eta_j)|={c'}{2}|{\sin(\xi_1+\eta_1) \sin(2\xi_j)}|$ when $\eta_j\in\{\eta_j^*-\delta_{\eta_j^*},\eta_j^*+\delta_{\eta_j^*}'\}$.

	We obtain, by a similar argument used to get \eqref{Lemm:Bessel2:2bb5:4}
	\begin{equation}\label{Lemm:Bessel2:Improved:2bb5:4}\begin{aligned}
			\mathfrak{B}^c_j(\eta_j)\ = \ & \mathfrak{B}^c_j(\eta_j^*)\ + \ \sigma_{\eta_j^*}(\eta_j-\eta_j^*)^2\int_{0}^1\int_{0}^s\mathrm{d}s\mathrm{d}s' |\partial_{\eta_j\eta_j}{\mathfrak{B}^c_j(ss'(\eta_j-\eta_j^*)+\eta_j^*) }|\\
			\ = \ & \mathfrak{B}^c_j(\eta_j^*)\ + \ \sigma_{\eta_j^*}(\eta_j-\eta_j^*)^2|\mathfrak{G}_{\eta_j^*}(\eta_j-\eta_j^*)|^2,		\end{aligned}
	\end{equation}
	in which $\mathfrak{G}_{\eta_j^*}(\eta_j-\eta_j^*)=\sqrt{\sigma_{\eta_j^*}\partial_{\eta_j\eta_j}{\mathfrak{B}^c_j(ss'(\eta_j-\eta_j^*)+\eta_j^*) }}>0$ is a smooth function with $\eta_j\in\mathbb{T}_{\eta_j^*}:=[-\pi,\pi]\cap (\eta_j^*-\delta_{\eta_j^*},\eta_j^*+\delta_{\eta_j^*}')$. We follow Step 1.2 in the proof of  Lemma \ref{Lemm:Bessel3} to define a new variable $y_{\eta_j*}= (\eta_j-\eta_j^*)\mathfrak{G}_{\eta_j^*}(\eta_j-\eta_j^*)$ and infer the existence of  a neighborhood of $U_{\eta_j^*}$ of the origin $0$ and a smooth function $\psi_{\eta_j^*}: C^\infty_c(U_{\eta_j^*})\to \mathbb{T}_{\eta_j^*}$ such that $\psi_{\eta_j^*}(y_{\eta_j^*})=\eta_j$. When $c'$ is closed to $1$, 
	the function $\psi_{\eta_j^*}$ is bijective and is  the inverse of $(\eta_j-\eta_j^*)\mathfrak{G}_{\eta_j^*}(\eta_j-\eta_j^*)$. Moreover, for any $\eta_j\in 	\mathbb{T}_{\eta_j}':=[-\pi,\pi]\backslash\cup_{\eta_j^*\in\mathfrak{J}_{\eta_j}}\Big(	[-\pi,\pi]\cap (\eta_j^*-\delta_{\eta_j^*},\eta_j^*+\delta_{\eta_j^*}')\Big)$, since $\partial_{\eta_j}\mathfrak{B}^c_j(\eta_j)\ne 0$ for $\eta_j\in 	\mathbb{T}_{\eta_j}'$, the function $\mathfrak{B}^c_j(\eta_j)$ is monotone on any interval $[\alpha',\beta']\subset \mathbb{T}_{\eta_j}'$. Therefore, $\partial_{\eta_j\eta_j}\mathfrak{B}^c_j(\eta_j)=-4\mathfrak{B}^c_j(\eta_j)$ is also monotone on any interval $[\alpha',\beta']\subset \mathbb{T}_{\eta_j}'$. Since  $	|	\partial_{\eta_j\eta_j}\mathfrak{B}^c_j(\eta_j)|={c'}{2}|{\sin(\xi_1+\eta_1) \sin(2\xi_j)}|$ when $\eta_j\in\{\eta_j^*-\delta_{\eta_j^*},\eta_j^*+\delta_{\eta_j^*}'\}$, we deduce that $|	\partial_{\eta_j\eta_j}\mathfrak{B}^c_j(\eta_j)|\le {c'}{2}|{\sin(\xi_1+\eta_1) \sin(2\xi_j)}|$ for any $\eta_j\in 	\mathbb{T}_{\eta_j}'$.

	Repeating the same argument used to prove \eqref{Lemm:Bessel2:2bb5:5} in the proof of  Lemma \ref{Lemm:Bessel3}, we obtain
	\begin{equation}\begin{aligned}\label{Lemm:Improved:Bessel2:2bb5:5}
			|	\partial_{\eta_j}\mathfrak{B}^c_j(\eta_j)|> &\	{(1-c')|{\sin(\xi_1+\eta_1) \sin(2\xi_j)}|}\ge  (1-c')\frac12\mathfrak{C}^2_{cut,j}|\ln(\lambda)|^{-\mathfrak{C}^1_{cut,j}}\epsilon_{\xi_j},
		\end{aligned}
	\end{equation}
	for all $\eta_j\in 	\mathbb{T}_{\eta_j}'$. Inequality \eqref{Lemm:Improved:Bessel2:2bb5:5} gives a uniform lower bound for $|	\partial_{\eta_j}\mathfrak{B}^c_j(\eta_j)|$ on the set of non-stationary points $\mathbb{T}_{\eta_j}'$.
	
	We follow \eqref{Lemm:Bessel2:2bb5:6} and split 
	\begin{equation}\begin{aligned}\label{Lemm:Improved:Bessel2:2bb5:6}
			\mathfrak{A}_j^a \ = \ &  \sum_{\eta_j^*\in\mathfrak{J}_{\eta_j}}\int_{\mathbb{T}_{\eta_j^*}}\mathrm{d}\eta_j e^{{\bf i}\mathfrak{B}^a_j} \check{F}_1(\xi_j+\eta_j)(\check{F}_1(\xi_1+\eta_1))^{\frac{1}{d}}\check{F}_1(\eta_j)\\
			& \ + \ \int_{\mathbb{T}_{\eta_j}'}\mathrm{d}\eta_j e^{{\bf i}\mathfrak{B}^a_j}\check{F}_1(\xi_j+\eta_j)(\check{F}_1(\xi_1+\eta_1))^{\frac{1}{d}}\check{F}_1(\eta_j)  \ = \ \sum_{\eta_j^*\in\mathfrak{J}_{\eta_j}}\mathfrak{A}_{\mathbb{T}_{\eta_j^*}} \ + \ \mathfrak{A}_{j,o}^a.
		\end{aligned}
	\end{equation}
	
	{\bf Step 2: Estimating $\mathfrak{A}_{\mathbb{T}_{\eta_j^*}}$.} 
	
	We employ the change of variables $\eta_j\to y_{\eta_j^*}$ and the same computations as in \eqref{Lemm:Bessel2:2bb5:7}-\eqref{Lemm:Bessel2:2bb5:10} to get

	\begin{equation}\begin{aligned}\label{Lemm:Improved:Bessel2:2bb5:7}
			& 	\mathfrak{A}_{\mathbb{T}_{\eta_j^*}} \ = \ \int_{\mathbb{T}_{\eta_j^*}}\mathrm{d}\eta_j e^{{\bf i}\mathfrak{B}^a_j}\check{F}_1(\xi_j+\eta_j)(\check{F}_1(\xi_1+\eta_1))^{\frac{1}{d}}\check{F}_1(\eta_j) \\
			&
			\ = \  e^{{\bf i}\mathfrak{B}^a_j(\eta_j^*)}\int_{U_{\eta_j}}\mathrm{d}y_{\eta_j^*} e^{{\bf i}\sigma_{\eta_j^*}t_0|y_{\eta_j^*}|^2}\psi_{\eta_j^*}'(y_{\eta_j^*})\check{F}_1(\xi_j+\psi_{\eta_j^*}(y_{\eta_j^*}))(\check{F}_1(\xi_1+\eta_1))^{\frac{1}{d}}\check{F}_1(\psi_{\eta_j^*}(y_{\eta_j^*})) \\
			& \ = \ e^{{\bf i}\mathfrak{B}^a_j(\eta_j^*)}\frac{e^{{\bf i}\pi^2{\sigma_{\eta_j^*}}/4}}{\sqrt{|t_0|}}{\psi_{\eta_j^*}'}(\eta_j^*)\check{F}_1(\xi_j+\eta_j^*)(\check{F}_1(\xi_1+\eta_1))^{\frac{1}{d}}\check{F}_1(\eta_j^*)\\
			& \ \ \ +\ e^{{\bf i}\mathfrak{B}^a_j(\eta_j^*)}\sum_{m=1}^\infty\frac{e^{{\bf i}\pi^2{\sigma_{\eta_j^*}}/4}}{\sqrt{|t_0|}}\Big(\frac{{\bf i}\pi^2}{t_0\sigma_{\eta_j^*}}\Big)^m\partial_{\eta_j}^{2m}[{\psi_{\eta_j^*}'}(\eta_j^*)\check{F}_1(\xi_j+\eta_j^*)\check{F}_1(\eta_j^*)](\check{F}_1(\xi_1+\eta_1))^{\frac{1}{d}}\\
			& \ = \ e^{{\bf i}\mathfrak{B}^a_j(\eta_j^*)}\frac{e^{{\bf i}\pi^2{\sigma_{\eta_j^*}}/4}}{\sqrt{|t_0|}}{\psi_{\eta_j^*}'}(\eta_j^*)\check{F}_1(\xi_j+\eta_j^*)(\check{F}_1(\xi_1+\eta_1))^{\frac{1}{d}}\check{F}_1(\eta_j^*)\ +\ e^{{\bf i}\mathfrak{B}^a_j(\eta_j^*)}\sum_{m=1}^\infty\frac{e^{{\bf i}\pi^2{\sigma_{\eta_j^*}}/4}}{\sqrt{|t_0|}}\\
			& \ \ \ \times\Big(\frac{{\bf i}\pi^2}{t_0\sigma_{\eta_j^*}}\Big)^m\left\{\sum_{i=0}^{2m}\psi_{\eta_j^*}^{(i+1)}(\eta_j^*)\partial_{\eta_j}^{2m-i}[\check{F}_1(\xi_j+\eta_j^*)\check{F}_1(\eta_j^*)]\right\}(\check{F}_1(\xi_1+\eta_1))^{\frac{1}{d}},
		\end{aligned}
	\end{equation}
	where, in the last line, we have used the product rule.
	
	Now, by precisely the same arguments used to obtain \eqref{Lemm:Bessel2:2bb5:11:1}, \eqref{Lemm:Bessel2:2bb5:17:1} and \eqref{Lemm:Bessel2:2bb5:21}, we can obtain precisely the same results, namely
	\begin{equation}
		\label{Lemm:Improved:Bessel2:2bb5:11:1}
		\begin{aligned}
			\psi_{\eta_j^*}'(\eta_j^*)\ = \ 	& 	  \frac{1}{\sqrt{|\partial_{\eta_j\eta_j}\mathfrak{B}^c_j(\eta_j^*)|}},
		\end{aligned}
	\end{equation}
	\begin{equation}
		\label{Lemm:Improved:Bessel2:2bb5:17:1}
		\psi_{\eta_j^*}^{(2i)}(\eta_j^*)\ =\ 0,\ \ \ \  \forall  i\in\mathbb{N}, i\ge 1,
	\end{equation}
	and
	\begin{equation}
		\label{Lemm:Improved:Bessel2:2bb5:21}			
		|\psi_{\eta_j^*}^{(2i+1)}(\eta_j^*)|\le 	(\mathfrak{C}^3_{cut,j})^i\left[\frac{1}{\big|\partial_{\eta_j\eta_j}\mathfrak{B}^c_j(\eta_j^*)\big|}\right]^{{2i+1}},	\ \ \ \  \forall  i\in\mathbb{N}, i\ge 1,			
	\end{equation}
	for some constant $\mathfrak{C}^3_{cut,j}>0$.
	
	Moreover, from the definition of $\check{F}_1$, we deduce the existence of a constant  $\mathfrak{C}^4_{cut,j}>0$ such that
	\begin{equation}
		\label{Lemm:Improved:Bessel2:2bb5:22}			
		|\partial_{\eta_j}^{i}[\check{F}_1(\xi_j+\eta_j^*)\check{F}_1(\eta_j^*)]| \ \le \ |\ln(\lambda)|^{-i\mathfrak{C}^4_{cut,j}},	\ \ \ \  \forall  i\in\mathbb{N}, i\ge 1.	
	\end{equation}
	Combining \eqref{Lemm:Improved:Bessel2:2bb5:7}, \eqref{Lemm:Improved:Bessel2:2bb5:17:1}, \eqref{Lemm:Improved:Bessel2:2bb5:21} and \eqref{Lemm:Improved:Bessel2:2bb5:22}	 yields
	\begin{equation}\begin{aligned}\label{Lemm:Improved:Bessel2:2bb5:23}
			& 	|\mathfrak{A}_{\mathbb{T}_{\eta_j^*}}| \ \le \ \frac{1}{\sqrt{|t_0|}}\frac{1}{\sqrt{|\partial_{\eta_j\eta_j}\mathfrak{B}^c_j(\eta_j^*)|}}\\
			&\ +\ \sum_{m=1}^\infty\frac{|\mathfrak{C}^5_{cut,j}|^{2m+1}}{\sqrt{|t_0|}^{2m+1}}\left\{\sum_{i=0, i \text{ is even}}^{2m}(\mathfrak{C}^3_{cut,j})^{{i}}\left[\frac{1}{\big|\partial_{\eta_j\eta_j}\mathfrak{B}^c_j(\eta_j^*)\big|}\right]^{{i+1}}|\ln(\lambda)|^{(2m-i)\mathfrak{C}^4_{cut,j}}\right\},
		\end{aligned}
	\end{equation}
	for an explicit constant $\mathfrak{C}^5_{cut,j}>0$. Inequality \eqref{Lemm:Improved:Bessel2:2bb5:23}, together with \eqref{Lemm:Improved:Bessel2:2bb2:1}, implies
	\begin{equation}\begin{aligned}\label{Lemm:Improved:Bessel2:2bb5:24}
			& 	|\mathfrak{A}_{\mathbb{T}_{\eta_j}^*}| 
			\ \le \ \frac{1}{\sqrt{|t_0|}}\frac{1}{\sqrt{\mathfrak{C}^2_{cut,j}|\ln(\lambda)|^{-\mathfrak{C}^1_{cut,j}}
					|\sin(2{\xi_j})|}}\\
			&\ +\ \sum_{m=1}^\infty\frac{|\mathfrak{C}^5_{cut,j}|^{2m+1}}{\sqrt{|t_0|}^{2m+1}}\left\{\sum_{i=0, i \text{ is even}}^{2m}(\mathfrak{C}^3_{cut,j})^{{i}}\left[\frac{1}{\big|\mathfrak{C}^2_{cut,j}|\ln(\lambda)|^{-\mathfrak{C}^1_{cut,j}}|\sin(2{\xi_j})|\big|}\right]^{{i+1}}|\ln(\lambda)|^{(2m-i)\mathfrak{C}^4_{cut,j}}\right\}\\
			&\ \le \ \frac{\mathcal{C}^o_{cut,j}}{\sqrt{|t_0|\mathfrak{C}^2_{cut,j}|\ln(\lambda)|^{-\mathfrak{C}^1_{cut,j}}|\sin(2{\xi_j})|}}\ +\ \sum_{m=1}^\infty\frac{|\mathcal{C}^o_{cut,j}|^{2m+1}}{\Big[{|t_0|\mathfrak{C}^2_{cut,j}|\ln(\lambda)|^{-\mathfrak{C}^1_{cut,j}}|\sin(2{\xi_j})|}\Big]^{\frac{2m+1}{2}}}\\
			\ \le \ &\frac{1}{\sqrt{|t_0|}}\frac{1}{\sqrt{\mathfrak{C}^2_{cut,j}|\ln(\lambda)|^{-\mathfrak{C}^1_{cut,j}}\epsilon_{\xi_j}}}\\
			&\ +\ \sum_{m=1}^\infty\frac{|\mathfrak{C}^5_{cut,j}|^{2m+1}}{\sqrt{|t_0|}^{2m+1}}\left\{\sum_{i=0, i \text{ is even}}^{2m}(\mathfrak{C}^3_{cut,j})^{{i}}\left[\frac{1}{\big|\mathfrak{C}^2_{cut,j}|\ln(\lambda)|^{-\mathfrak{C}^1_{cut,j}}\epsilon_{\xi_j}\big|}\right]^{{i+1}}|\ln(\lambda)|^{(2m-i)\mathfrak{C}^4_{cut,j}}\right\}\\
			&\ \le \ \frac{\mathcal{C}^o_{cut,j}}{\sqrt{|t_0|\mathfrak{C}^2_{cut,j}|\ln(\lambda)|^{-\mathfrak{C}^1_{cut,j}}\epsilon_{\xi_j}}}\ +\ \sum_{m=1}^\infty\frac{|\mathcal{C}^o_{cut,j}|^{2m+1}}{\Big[{|t_0|\mathfrak{C}^2_{cut,j}|\ln(\lambda)|^{-\mathfrak{C}^1_{cut,j}}\epsilon_{\xi_j}}\Big]^{\frac{2m+1}{2}}},
		\end{aligned}
	\end{equation}
	in which the   constants $\mathcal{C}^o_{cut,j}>0$ is explicit.
	
	\smallskip
	{\bf Step 3: Estimating $\mathfrak{A}_{\mathbb{T}_{\eta_j,o}}$.} 
	Note that \eqref{Lemm:Bessel2:2bb5:23} still holds true in our case. Thus, replacing $\Xi_j^0$ by $t_0$, we deduce a similar estimate with \eqref{Lemm:Bessel2:2bb5:Final2}
	\begin{equation}\begin{aligned}
			\label{Lemm:Improved:Bessel2:2bb5:Final2}
			|\mathfrak{A}_{j,o}^a| 
			\ \le \ & \left|\frac{{\bf i}}{t_0}\int_{\mathbb{T}_{\eta_j}'}\mathrm{d}\eta_j \frac{\partial_{\eta_j\eta_j}\mathfrak{B}^c_j}{|\partial_{\eta_j}\mathfrak{B}^c_j|^2}e^{{\bf i}t_0\mathfrak{B}^c_j} \right| \ + \ \left|\frac{{\bf i}}{t_0}\frac{e^{{\bf i}t_0\mathfrak{B}^c_j} }{{\bf i}\partial_{\eta_j}\mathfrak{B}^c_j}\Big|_{\partial\mathbb{T}_{\eta_j}'}\right|\\
			\ \le \  & \frac{\mathfrak{C}^o_{cut,j}}{|t_0||\mathfrak{C}^2_{cut,j}|\ln(\lambda)|^{-\mathfrak{C}^1_{cut,j}}|\sin(2{\xi_j})||^2} \ + \ \frac{\mathfrak{C}^o_{cut,j}}{|t_0||\mathfrak{C}^2_{cut,j}|\ln(\lambda)|^{-\mathfrak{C}^1_{cut,j}}|\sin(2{\xi_j})||},
		\end{aligned}
	\end{equation}
	in which the   constants on the right hand side are  explicit. 
	
	Combining \eqref{Lemm:Improved:Bessel2:2bb5:24} and 
	\eqref{Lemm:Improved:Bessel2:2bb5:Final2}, we find
	\begin{equation}\label{factordx}\begin{aligned}
			|\mathfrak{A}_j^a | \ \le \ & \frac{\mathfrak{C}^o_{cut,j}}{|t_0||\mathfrak{C}^2_{cut,j}|\ln(\lambda)|^{-\mathfrak{C}^1_{cut,j}}|\sin(2{\xi_j})||^2} \ + \ \frac{\mathfrak{C}^o_{cut,j}}{|t_0|\mathfrak{C}^2_{cut,j}|\ln(\lambda)|^{-\mathfrak{C}^1_{cut,j}}|\sin(2{\xi_j})|}\\
			& + \ 	\frac{3\mathcal{C}^o_{cut,j}}{\sqrt{|t_0|\mathfrak{C}^2_{cut,j}|\ln(\lambda)|^{-\mathfrak{C}^1_{cut,j}}|\sin(2{\xi_j})|}}\ +\ \sum_{m=1}^\infty\frac{3|\mathcal{C}^o_{cut,j}|^{2m+1}}{\Big[{|t_0|\mathfrak{C}^2_{cut,j}|\ln(\lambda)|^{-\mathfrak{C}^1_{cut,j}}|\sin(2{\xi_j})|}\Big]^{\frac{2m+1}{2}}}.			\end{aligned}
	\end{equation}
	where, we have used \eqref{Lemm:Improved:Bessel2:2bb5:3:1}, meaning that there are at most $3$ terms of the type $\mathfrak{A}_{\mathbb{T}_{\eta_j^*}}$.

	Under assumption \eqref{Lemm:Improved:Bessel2:2bb2:1:2}, we find
	
	\begin{equation}\label{factordx2}\begin{aligned}
			|\mathfrak{A}_j^a | \ \lesssim \ & |\mathfrak{A}_j^a |_\infty:= \frac{\mathfrak{C}^o_{cut,j}}{|t_0||\mathfrak{C}^2_{cut,j}|\ln(\lambda)|^{-\mathfrak{C}^1_{cut,j}}|\sin(2{\xi_j})||^2} \ + \ 	\frac{3\mathcal{C}^o_{cut,j}}{\sqrt{|t_0|\mathfrak{C}^2_{cut,j}|\ln(\lambda)|^{-\mathfrak{C}^1_{cut,j}}|\sin(2{\xi_j})|}}.			\end{aligned}
	\end{equation}

	%
	
	{\bf Step 4: The final estimate.} Finally, we estimate 
	\begin{equation}
		\label{Lemm:Improved:Step4:1}\begin{aligned}
			\sum_{m\in\mathbb{Z}^d}|\mathfrak{F}^{O}(m,t_0) |^4
			\ = \ & \int_{[-\pi,\pi]^d}\mathrm{d}\xi  \Big|\int_{-\pi}^{\pi}\mathrm{d}\eta_1 \prod_{j=2}^d\mathfrak{A}_j^a(\xi_1,\eta_1,\xi_j)(\check{F}_1(\xi_1+\eta_1))^{\frac{1}{d}}\check{F}_1(\eta_1)e^{{\bf i}\mathfrak{B}_1^a(\xi_1,\eta_1) }\Big|^2\\
			\ \le \ & \int_{[-\pi,\pi]^d}\mathrm{d}\xi  \Big|\int_{-\pi}^{\pi}\mathrm{d}\eta_1 \prod_{j=2}^d|\mathfrak{A}_j^a(\xi_1,\eta_1,\xi_j)|\Big|^2\\
			\ \lesssim \ & \int_{[-\pi,\pi]^d}\mathrm{d}\xi    \prod_{j=2}^d\left[\sup_{\eta_1\in[-\pi,\pi]}|\mathfrak{A}_j^a|\right]^2\\
			\ \lesssim \ &   \prod_{j=2}^d\int_{[-\pi,\pi]}\mathrm{d}\xi_j  \left[\sup_{\eta_1\in[-\pi,\pi]}|\mathfrak{A}_j^a|\right]^2.\end{aligned}
	\end{equation}
	
	On the domain $\mathbb{T}_{\xi_j}$, we bound, by \eqref{factordx2}
	\begin{equation}\label{factordx2:1}\begin{aligned}
			&	\int_{\mathbb{T}_{\xi_j}}\mathrm{d}\xi_j	\left[\sup_{\eta_1\in[-\pi,\pi]}|\mathfrak{A}_j^a|\right]^2 \ \lesssim\ \\
			\ \lesssim\	& \int_{\mathbb{T}_{\xi_j}}\mathrm{d}\xi_j\left[ \frac{\mathfrak{C}^o_{cut,j}}{|t_0||\mathfrak{C}^2_{cut,j}|\ln(\lambda)|^{-\mathfrak{C}^1_{cut,j}}|\sin(2{\xi_j})||^2} \ + \ 	\frac{3\mathcal{C}^o_{cut,j}}{\sqrt{|t_0|\mathfrak{C}^2_{cut,j}|\ln(\lambda)|^{-\mathfrak{C}^1_{cut,j}}|\sin(2{\xi_j})|}}\right]^2	
			\\
			\ \lesssim\	& \int_{\mathbb{T}_{\xi_j}}\mathrm{d}\xi_j\left[ \frac{1}{|t_0|^2|\mathfrak{C}^2_{cut,j}|\ln(\lambda)|^{-\mathfrak{C}^1_{cut,j}}|\sin(2{\xi_j})||^4} \ + \ 	\frac{1}{{|t_0|\mathfrak{C}^2_{cut,j}|\ln(\lambda)|^{-\mathfrak{C}^1_{cut,j}}|\sin(2{\xi_j})|}}\right]	\\
			\ \lesssim\	& \int_{\mathbb{T}_{\xi_j}}\mathrm{d}\xi_j \frac{1}{|t_0|^2|\mathfrak{C}^2_{cut,j}|\ln(\lambda)|^{-\mathfrak{C}^1_{cut,j}}|^4|\sin(2{\xi_j})|^{1-\epsilon}|\sin(2{\xi_j})|^{3+\epsilon}}\\ 
			& \ + \ \int_{\mathbb{T}_{\xi_j}}\mathrm{d}\xi_j	\frac{1}{|t_0|\mathfrak{C}^2_{cut,j}|\ln(\lambda)|^{-\mathfrak{C}^1_{cut,j}}|\sin(2{\xi_j})|^{1-\epsilon}|\sin(2{\xi_j})|^{\epsilon}}\\	
			\ \lesssim_\epsilon \ &  \frac{1}{|t_0|^2|\mathfrak{C}^2_{cut,j}|\ln(\lambda)|^{-\mathfrak{C}^1_{cut,j}}|^4|\epsilon_{\xi_j}|^{3+\epsilon}} \ + \ 	\frac{1}{{|t_0|\mathfrak{C}^2_{cut,j}|\ln(\lambda)|^{-\mathfrak{C}^1_{cut,j}}|\epsilon_{\xi_j}|^{\epsilon}}},\end{aligned}
	\end{equation}
	for any $0<\epsilon<1.$
	Moreover, it is straightforward to see that
	\begin{equation}
		\label{factordx2:2}\begin{aligned}
			\int_{[-\pi,\pi]\backslash\mathbb{T}_{\xi_j}}\mathrm{d}\xi_j	\left[\sup_{\eta_1\in[-\pi,\pi]}|\mathfrak{A}_j^a|\right]^2 \ \lesssim \ & \epsilon_{\xi_j}.\end{aligned}
	\end{equation}
	Combining \eqref{factordx2:1} and \eqref{factordx2:2}, and balancing $\epsilon_{\xi_j}$ under the constraint \eqref{Lemm:Improved:Bessel2:2bb2:1:2}, we obtain
	\begin{equation}\label{factordx2:3}\begin{aligned}
			\int_{[-\pi,\pi]}\mathrm{d}\xi_j	\left[\sup_{\eta_1\in[-\pi,\pi]}|\mathfrak{A}_j^a|\right]^2 
			\ \lesssim \ &  \frac{1}{|t_0|^{\frac12-}|\mathfrak{C}^2_{cut,j}|\ln(\lambda)|^{-\mathfrak{C}^1_{cut,j}}|^{1-}}.\end{aligned}
	\end{equation}
	Plugging \eqref{factordx2:3} into \eqref{Lemm:Improved:Step4:1}, we find
	
	\begin{equation}
		\label{Lemm:Improved:Step4:2}\begin{aligned}
			\sum_{m\in\mathbb{Z}^d}|\mathfrak{F}^{O}(m,t_0) |^4
			\lesssim \ &   |\ln|\lambda||^{\mathcal{C}_{\mathfrak{F}^{O},4'}}\mathcal{C}_{\mathfrak{F}^{O},4}\langle{t_0}\rangle^{-\big(\frac{(
					d-1)}{2}-\big)}.\end{aligned}
	\end{equation}

	Thus, we obtain
	\begin{equation}
		\label{Lemm:Improved:Bessel2:final}\begin{aligned}
			& \|\mathfrak{F}^{O}(\cdot,t_0) \|_{l^4}^4\ \le  \ |\ln|\lambda||^{\mathcal{C}_{\mathfrak{F}^{O},4'}}\mathcal{C}_{\mathfrak{F}^{O},4}\langle{t_0}\rangle^{-\big(\frac{
					d-1}{2}-\big)},\end{aligned}
	\end{equation}
	yielding
	\begin{equation}
		\label{Lemm:Improved:Bessel2:final:a}\begin{aligned}
			& \|\mathfrak{F}^{O}(\cdot,t_0) \|_{l^4}\ \le  \ |\ln|\lambda||^{\mathcal{C}_{\mathfrak{F}^{O},4'}/4}\mathcal{C}_{\mathfrak{F}^{O},4}^\frac14\langle{t_0}\rangle^{-{\big(\frac{
						d-1}{8}-\big)}},\end{aligned}
	\end{equation}
	which is the conclusion of  Lemma \ref{Lemm:Improved:Bessel3}. 
\end{proof}
\begin{lemma}\label{Lemm:Improved:Bessel4} There exist universal constants $\mathfrak{C}_{\mathfrak{F}^{O},3}, \mathfrak{C}_{\mathfrak{F}^{O},3'}>0$ independent of $t_0$ and $\lambda$, such that
	\begin{equation}
		\label{Lemm:Improved:Bessel4:1}\begin{aligned}
			&	\|\mathfrak{F}^{O}(\cdot,t_0) \|_{l^3} \ \le  |\ln|\lambda||^{\mathfrak{C}_{\mathfrak{F}^{O},3'}}\mathfrak{C}_{\mathfrak{F}^{O},3}\langle{t_0}\rangle^{-{\big(\frac{
						d-1}{12}-\big)}}.\end{aligned}
	\end{equation}
\end{lemma}

\begin{proof}
	By interpolating between $l_4$ and  $l_2$,   we  get  \eqref{Lemm:Improved:Bessel4:1}.

\end{proof}

\subsection{Free momentum estimates}

We follow the notation of Section \ref{Subsec:DispersiveEstimates} and let $\tilde\Psi_1( \xi,V)$ be the component of  $\check\Psi_1( \xi,V,W)$ that concerns only $V$. We set $\xi=2\pi k$, $V=2\pi k^*$ and denote $\tilde\Psi_1( \xi,V)=\mu_1(k,k^*)$.

\begin{lemma}[Degree-one vertex estimate]\label{lemma:degree1vertex}   Let $H(k)$ be any function  in $L^2(\mathbb{T}^d)$. Let  $\diamond_\ell:\mathbb{R}_+\to\mathbb{R}_+$ be a non-negative function such that $\diamond_\ell\in   L^{\mathscr{M}'}(\mathbb{R}_+) $ where $1/\mathscr{M}'+1/\mathscr{M}=1$ and $\mathscr M$ is a sufficiently large constant.  For any $k^*\in\mathbb{T}^d$, $d\ge 2$,  and $\sigma \in\{\pm 1\}$, we have 
	\begin{equation}
		\label{eq:degree1vertex}
		\int_{\mathbb{R}_+}\mathrm{d}\nu \diamond_\ell(\nu)\Big|	\int_{\mathbb{T}^d}\mathrm{d}ke^{\mathbf{i}\nu(\omega_\infty (k)+\sigma \omega_\infty (k^*+k))}\mu_1(k,k^*)H\Big|\ \lesssim \ \langle \ln|\lambda| \rangle^{2+c\eth}\|H\|_{L^2}\|\diamond_\ell\|_{L^{\mathscr{M}'}},
	\end{equation}
	in which the constants on the right hand side  are universal and $\eth$ is associated to the definition of the cut-off function.  
	
\end{lemma}

\begin{proof} We denote $\langle\ln|\lambda| \rangle^{\eth}$  by $d(\mathfrak{S})^6$.
	We first bound
	\begin{equation}
		\label{eq:degree1vertex:E7:a:1}\begin{aligned}
			& 	\int_{\mathbb{R}_+}\mathrm{d}\nu \diamond_\ell(\nu)\Big|	\int_{\mathbb{T}^d}\mathrm{d}ke^{\mathbf{i}\nu(\omega_\infty (k)+\sigma \omega_\infty (k^*+k))}\mu_1(k,k^*)H\Big|\\
			\lesssim\ &  \left[\int _{0}^{\infty }\mathrm{d}\nu \left|\int_{\mathbb{T}^d}\mathrm{d}k \mu_1H e^{{\bf i}\nu(\omega_\infty(k)+\sigma\omega_\infty(k^*+k))}\right|^\mathscr{M}\right]^\frac{1}{\mathscr{M}}\left[\int _{-\infty }^{\infty }\mathrm{d}\nu |\diamond_\ell(\nu)|^{\mathscr{M}'} \right]^\frac{1}{\mathscr{M}'}.\end{aligned}
	\end{equation}

	We now estimate the integral  on the right hand side of \eqref{eq:degree1vertex:E7:a:1} using a $TT^*$ argument. We choose ${G}(\nu)$ to be a test function in $L^p(\mathbb{R})$ with $p=\mathscr{M}'$ and $q=\mathscr{M}$,  and develop
	
	\begin{equation}
		\label{eq:degree1vertex:E7:a:2}\begin{aligned}
			&\  \Big|\int_{\mathbb{R}_+} {\mathrm{d}\nu}	 \int_{\mathbb{T}^d}\mathrm{d}k \mu_1H(k)e^{{\bf i}\nu(\omega_\infty(k)+\sigma\omega_\infty(k^*+k))}{G}(\nu)\Big|.
		\end{aligned}
	\end{equation}
	We  study the $L^2$-norm with respect to $k$
	\begin{equation}
		\label{eq:degree1vertex:E7:a:3}\begin{aligned}
			& \Big\|\int_{\mathbb{R}_+} {\mathrm{d}\nu}	 \mu_1 e^{{\bf i}\nu(\omega_\infty(k)+\sigma\omega_\infty(k^*+k))}{G}(\nu)\Big\|_{L^2}^2,
		\end{aligned}
	\end{equation}
	which we will expand and bound as follows
	\begin{equation}
		\label{eq:degree1vertex:E7:a:4}\begin{aligned}
			& \Big|\int_{\mathbb{R}_+} {\mathrm{d}\nu}\int_{\mathbb{R}_+} {\mathrm{d}\nu'}	\int_{\mathbb{T}^d}\mathrm{d}k\mu_1^2 e^{{\bf i}\nu(\omega_\infty(k)+\sigma\omega_\infty(k^*+k))}{G}(\nu) e^{-{\bf i}\nu'(\omega_\infty(k)+\sigma\omega_\infty(k^*+k))}\overline{{G}(\nu')}\Big|\\
			\lesssim\	& \left[\int_{\mathbb{R}_+} {\mathrm{d}\nu}|{G}(\nu)|^p\right]^\frac2p \left[\int_{\mathbb{R}_+} {\mathrm{d}\nu}|\mathfrak{F}_\omega(\nu)|^\frac{q}{2}\right]^\frac{2}{q},
		\end{aligned}
	\end{equation}
	where
	\begin{equation}
		\label{eq:degree1vertex:E7:a:5}\begin{aligned}
			\mathfrak{F}_\omega(\nu) \ = \ & \int_{\mathbb{T}^d}\mathrm{d}k \mu_1^2(k)\  e^{{\bf i}\nu(\omega_\infty(k)+\sigma\omega_\infty(k^*+k))}.
		\end{aligned}
	\end{equation}
	
	To estimate the oscillatory integral $\int_{\mathbb{T}^d}\mathrm{d}k \mu_1^2 e^{{\bf i}\nu(\omega_\infty(k)+\sigma\omega_\infty(k^*+k))}$, we employ the same method used to prove Lemma \ref{Lemm:Bessel3}. To this end, we set $\xi=2\pi k$, $\xi^*=2\pi k^*$  and write the phase as
	$$e^{{\bf i}\nu(\omega_\infty(k)+\sigma\omega_\infty(k^*+k))}=e^{{\bf i}\nu Z(\xi)},$$
	in which
	\begin{equation}\label{eq:degree1vertex:E7:1}
		\begin{aligned}
			Z(\xi) \ = & \ \sin(\xi_1)\Big[\sin^2(\xi_1)+\cdots+\sin^2(\xi_d)\Big]\\
			&\ + \  \sigma\sin(\xi_1+\xi^*_1)\Big[\sin^2(\xi_1+\xi^*_1)+\cdots+\sin^2(\xi_d+\xi^*_d)\Big],
		\end{aligned}
	\end{equation}
	with $\xi=(\xi_1,\cdots,\xi_d)$, $\xi^*=(\xi_1^*,\cdots,\xi_d^*)$.
	This phase can be written as the sum \begin{equation}\label{sumZ}
		Z=Z_1+\sum_{j=2}^dZ_j,
	\end{equation}
	in which
	\begin{equation}\label{eq:degree1vertex:E7:2}
		\begin{aligned}
			Z_1 \ = & \ \Big(\frac{3}{4}+\frac{d-1}{2}\Big)\sin(\xi_1) \ - \ \frac{1}{4}\sin(3\xi_1)\\
			&\ + \  \sigma\Big(\frac{3}{4}+\frac{d-1}{2}\Big)\sin(\xi_1+\xi^*_1) \ - \ \frac{\sigma}{4}\sin(3\xi_1+3\xi^*_1),
		\end{aligned}
	\end{equation}
	and
	\begin{equation}\label{eq:degree1vertex:E7:3}
		\begin{aligned}
			Z_j \ = & \ -\frac12\sin(\xi_1)\cos(2\xi_j)\
			\ - \  \frac{\sigma}{2}\sin(\xi_1+\xi^*_1)\cos(2\xi_j+2\xi^*_j).
		\end{aligned}
	\end{equation} 
	Therefore 
	\begin{equation}\label{eq:degree1vertex:E7:3:a}
		\begin{aligned}
			& 	\left|\int_{\mathbb{T}^d}\mathrm{d}k \mu_1 e^{{\bf i}\nu(\omega_\infty(k)+\sigma\omega_\infty(k^*+k))}\right| \
			\lesssim \ \left|\int_{[-\pi,\pi]}\mathrm{d}\xi_1 \mu_1e^{{\bf i}\nu Z_1}\prod_{j=2}^d\int_{[-\pi,\pi]}\mathrm{d}\xi_j e^{{\bf i}\nu Z_j}\right|.
		\end{aligned}
	\end{equation}
	From the 	construction of the cut-off function $\mu_1$ in Section \ref{Subsec:DispersiveEstimates}, we can assume that  
	
	\begin{equation}\label{eq:degree1vertex:E7:3:a:1}\mu_1=\mu_1(\xi^*,\xi_1).	\end{equation}
	Setting 
	\begin{equation}\label{yj}
		\begin{aligned}
			\mathcal{Y}_j\	& = \ 	\mu_1^{\frac{1}{d}}(\xi^*,\xi_1)\int_{[-\pi,\pi]}\mathrm{d}\xi_j e^{{\bf i}\nu Z_j},
		\end{aligned}
	\end{equation}
	we divide the proof into several steps. 
	
	\smallskip
	
	{\bf Step 1: Stationary points of $Z_j$.}
	
	\smallskip
	
	We now study the stationary point $\xi_j^o$ of $Z_j$, which means $\partial_{\xi_j}Z_j(\xi_j^o)=0$.
	
	By setting $$\sin(\xi_1)
	\ + \  {\sigma}\sin(\xi_1+\xi^*_1)e^{{\bf i}2\xi^*_j}=\Big|\sin(\xi_1)
	\ + \  {\sigma}\sin(\xi_1+\xi^*_1)e^{{\bf i}2\xi^*_j}\Big|e^{{\bf i}2\upsilon},$$
	with $\upsilon\in[-\pi,\pi]$, we compute
	\begin{equation}\label{eq:degree1vertex:E7:13}
		\begin{aligned}
			0\  = \ \partial_{\xi_j}Z_j \ = & \ \mathrm{Im}\Big[\sin(\xi_1)e^{{\bf i}2\xi_j^o}\Big]\
			\ + \  {\sigma}\mathrm{Im}\Big[\sin(\xi_1+\xi^*_1)e^{{\bf i}2\xi_j^o+{\bf i}2\xi^*_j}\Big]\\
			\ = & \ \mathrm{Im}\Big[e^{{\bf i}2\xi_j^o}\Big(\sin(\xi_1)
			\ + \  {\sigma}\sin(\xi_1+\xi^*_1)e^{{\bf i}2\xi^*_j}\Big)\Big]\\
			\ = & \ \mathrm{Im}\Big[e^{{\bf i}2\xi_j^o+{\bf i}2\upsilon}\Big|\sin(\xi_1)
			\ + \  {\sigma}\sin(\xi_1+\xi^*_1)e^{{\bf i}2\xi^*_j}\Big|\Big],
		\end{aligned}
	\end{equation} 
	and
	\begin{equation}\label{eq:degree1vertex:E7:14}
		\begin{aligned}
			\partial_{\xi_j\xi_j}Z_j \ = & \ 2\mathrm{Re}\Big[e^{{\bf i}2\xi_j^o+{\bf i}2\upsilon}\Big|\sin(\xi_1)
			\ + \  {\sigma}\sin(\xi_1+\xi^*_1)e^{{\bf i}2\xi^*_j}\Big|\Big].
		\end{aligned}
	\end{equation}
	If $\sigma=1$, we set $\tilde\xi_1^*=\xi_1^*$, otherwise for $\sigma=-1$, we set   $\tilde\xi_1^*=\xi_1^*+\pi$. Observing that 
	\begin{equation}\label{eq:degree1vertex:E7:15}
		\begin{aligned}
			& 2\Big|\sin(\xi_1)
			\ + \  {\sigma}\sin(\xi_1+\xi^*_1)e^{{\bf i}2\xi^*_j}\Big|\ = \  2\Big|\sin(\xi_1)
			\ + \ \sin(\xi_1+\tilde\xi_1^*)e^{{\bf i}2\xi^*_j}\Big|\\
			\ge \ & 2\Big||\sin(\xi_1)|-|\sin(\xi_1+\tilde\xi_1^*)|\Big| \
			\ge	\ \Big[|\sin(\xi_1)|+|\sin(\xi_1+\tilde\xi_1^*)|\Big]\Big||\sin(\xi_1)|-|\sin(\xi_1+\tilde\xi_1^*)|\Big|\\
			\ge \ & \Big|\sin^2(\xi_1)
			\ - \  \sin^2(\xi_1+\tilde\xi_1^*)\Big| \ = \ \Big|\sin(\tilde\xi_1^*)\sin(2\xi_1+\tilde\xi_1^*)\Big|,
		\end{aligned}
	\end{equation}
	we deduce
	\begin{equation}\label{eq:degree1vertex:E7:16}
		\begin{aligned}
			\Big|\sin(\xi_1)
			\ + \  {\sigma}\sin(\xi_1+\xi^*_1)e^{{\bf i}2\xi^*_j}\Big|
			\ge \ & \frac12\Big|\sin(\tilde\xi_1^*)\sin(2\xi_1+\tilde\xi_1^*)\Big|.
		\end{aligned}
	\end{equation}
	Moreover, we also have
	\begin{equation}\label{eq:degree1vertex:E7:17}
		\begin{aligned}
			& \Big|\sin(\xi_1)
			\ + \  {\sigma}\sin(\xi_1+\xi^*_1)e^{{\bf i}2\xi^*_j}\Big|^2\\
			= \ & \sin^2(\xi_1)+\sin^2(\xi_1+\tilde\xi_1^*)+2\sin(\xi_1)\sin(\xi_1+\tilde\xi_1^*)\cos(2\xi_j^*)\\
			= \ & \Big(\sin(\xi_1)\cos(2\xi_j^*)+ \sin(\xi_1+\tilde\xi_1^*)\Big)^2\ + \sin^2(\xi_1)\sin^2(2\xi_j^*).
		\end{aligned}
	\end{equation}
	
	%
	%
	%
	%
	Combining \eqref{eq:degree1vertex:E7:16} and \eqref{eq:degree1vertex:E7:17}, we obtain
	\begin{equation}\label{eq:degree1vertex:E7:18}
		\begin{aligned}
			\Big|\sin(\xi_1)
			\ + \  {\sigma}\sin(\xi_1+\xi^*_1)e^{{\bf i}2\xi^*_j}\Big|
			\ge \ & \frac18\Big|\sin(\tilde\xi_1^*)\sin(2\xi_1+\tilde\xi_1^*)\Big|\\
			& \ + \ \frac12\Big|\sin(\xi_1)\sin(2\xi_j^*)\Big|
		\end{aligned}
	\end{equation}
	which is $0$ when either $\xi_1^*=0$, $\pi$ or $\xi^*_j$ belongs to the singular manifold $  \mathfrak{S} $ and are removed by the cut-off function $	\mu_1^{\frac{1}{d}}$. Therefore $\partial_{\xi_j}Z=0$ when $\sin(2\xi_j^o+2\upsilon)=0$ under the constraint that $\xi^*_j,\xi_1^*$ are outside of the singular manifold. In this case, $\cos(2\xi_j^o+2\upsilon)=1$ and the following estimate then holds true
	\begin{equation}
		\label{eq:degree1vertex:E7:Zj}\begin{aligned}
			|\partial_{\xi_j\xi_j}Z_j(\xi_j^o)|\ = \ &\Big|\sin(\xi_1)
			\ + \  {\sigma}\sin(\xi_1+\xi^*_1)e^{{\bf i}2\xi^*_j}\Big|\\
			\ \ge\ & \frac18\Big|\sin(\tilde\xi_1^*)\sin(2\xi_1+\tilde\xi_1^*)\Big|\
			\ + \ \frac18\Big|\sin(\xi_1)\sin(2\xi_j^*)\Big|.\end{aligned}
	\end{equation}
	
	We now write
	\begin{equation}\label{eq:degree1vertex:E7:Zj:1}\begin{aligned}
			Z_j(\xi_j)\ = & \ Z_j(\xi_j^o)\ + \ \int_{0}^1\mathrm{d}s \frac{\partial Z_j(s(\xi_j-\xi_j^o)+\xi_j^o) }{\partial s}\\
			\ = & \ Z_j(\xi_j^o)\ + \ (\xi_j-\xi_j^o)\int_{0}^1\mathrm{d}s \partial_{\xi_j}{Z_j(s(\xi_j-\xi_j^o)+\xi_j^o)},
		\end{aligned}
	\end{equation}
	and
	\begin{equation}\label{eq:degree1vertex:E7:Zj:2}\begin{aligned}
			\partial_{\xi_j}Z_j(s(\xi_j-\xi_j^o)+\xi_j^o)\ = \ & \partial_{\xi_j}Z_j(\xi_j^o)\ + \ \int_{0}^s\mathrm{d}s' \frac{\partial\partial_{\xi_j}Z_j(ss'(\xi_j-\xi_j^o)+\xi_j^o) }{\partial s'}\\
			\ =\ &  (\xi_j-\xi_j^o)\int_{0}^s\mathrm{d}s' \partial_{\xi_j\xi_j}{Z_j(ss'(\xi_j-\xi_j^o)+\xi_j^o) },\end{aligned}
	\end{equation}
	which imply
	\begin{equation}\label{eq:degree1vertex:E7:Zj:3}
		Z_j(\xi_j)\ =\ Z_j(\xi_j^o)\ + \ (\xi_j-\xi_j^o)^2\int_{0}^1\int_{0}^s\mathrm{d}s\mathrm{d}s' \partial_{\xi_j\xi_j}{Z_j(ss'(\xi_j-\xi_j^o)+\xi_j^o) }.
	\end{equation}
	Due to \eqref{eq:degree1vertex:E7:Zj}, we set $\partial_{\xi_j\xi_j}Z_j(\xi_j^o)=\sigma_{\xi_j^o}|\partial_{\xi_j\xi_j}Z_j(\xi_j^o)|$, where $\sigma_{\xi_j^o}$ is either $1$ or $-1$. There exist constants $\delta_{\xi_j^o},\delta_{\xi_j^o}'>0$ such that for  all $\xi_j\in[-\pi,\pi]\cap (\xi_j^o-\delta_{\xi_j^o},\xi_j^o+\delta_{\xi_j^o}')$, we have $\partial_{\xi_j\xi_j}Z_j(ss'(\xi_j-\xi_j^o)+\xi_j^o)=\sigma_{\xi_j^o}|\partial_{\xi_j\xi_j}Z_j(ss'(\xi_j-\xi_j^o)+\xi_j^o)|$ and
	\begin{equation}\label{eq:degree1vertex:E7:Zj:4}\begin{aligned}
			|\partial_{\xi_j\xi_j}Z_j(\xi_j)|
			\	\ge\  &	c\Big|\sin(\xi_1)
			\ + \  {\sigma}\sin(\xi_1+\xi^*_1)e^{{\bf i}2\xi^*_j}\Big|	\\
			\	\ge\  &	c\frac18\Big|\sin(\tilde\xi_1^*)\sin(2\xi_1+\tilde\xi_1^*)\Big|
			\ + \ c\frac18\Big|\sin(\xi_1)\sin(2\xi_j^*)\Big|,\end{aligned}
	\end{equation}
	for all $s\in[0,1], s'\in[0,s]$ and for some constant $0<c<1$. Moreover, $	|	\partial_{\xi_j\xi_j}Z_j(\xi_j)|=c\Big|\sin(\xi_1)
	\ + \  {\sigma}\sin(\xi_1+\xi^*_1)e^{{\bf i}2\xi^*_j}\Big|$ when $\xi_j\in\{\xi_j^o-\delta_{\xi_j^o},\xi_j^o+\delta_{\xi_j^o}'\}$.

	We obtain
	\begin{equation}\label{eq:degree1vertex:E7:Zj:5}\begin{aligned}
			Z_j(\xi_j)\ = \ & Z_j(\xi_j^o)\ + \ \sigma_{\xi_j^o}(\xi_j-\xi_j^o)^2\int_{0}^1\int_{0}^s\mathrm{d}s\mathrm{d}s' |\partial_{\xi_j\xi_j}Z_j(ss'(\xi_j-\xi_j^o)+\xi_j^o)|\\
			\ = \ & Z_j(\xi_j^o)\ + \ \sigma_{\xi_j^o}(\xi_j-\xi_j^o)^2|\mathfrak{G}_{\xi_j^o}(\xi_j-\xi_j^o)|^2,		\end{aligned}
	\end{equation}
	in which $\mathfrak{G}_{\xi_j^o}(\xi_j-\xi_j^o)=\sqrt{\sigma_{\xi_j^o}\partial_{\xi_j\xi_j}{Z_j(ss'(\xi_j-\xi_j^o)+\xi_j^o) }}>0$ is a smooth function with $\xi_j\in\mathbb{T}_{\xi_j^*}:=[-\pi,\pi]\cap (\xi_j^o-\delta_{\xi_j^o},\xi_j^o+\delta_{\xi_j^o}')$. We  define a new variable $y_{\xi_j^o}= (\xi_j-\xi_j^o)\mathfrak{G}_{\xi_j^o}(\xi_j-\xi_j^o)$ and infer the existence of  a neighborhood of $U_{\xi_j^o}$ of the origin $0$ and a smooth function $\psi_{\xi_j^o}: C^\infty_c(U_{\xi_j^o})\to \mathbb{T}_{\xi_j^o}$ such that $\psi_{\xi_j^o}(y_{\xi_j^o})=\xi_j$. The function $\psi_{\xi_j^o}$ is bijective and is  the inverse of $(\xi_j-\xi_j^o)\mathfrak{G}_{\xi_j^o}(\xi_j-\xi_j^o)$. Moreover, for any $\xi_j\in 	\mathbb{T}_{\xi_j}':=[-\pi,\pi]\backslash\Big(\cup_{\xi_j^o\in\mathfrak{J}_{\xi_j}}	[-\pi,\pi]\cap (\xi_j^o-\delta_{\xi_j^o},\xi_j^o+\delta_{\xi_j^o}')\Big)$, since $\partial_{\xi_j}Z_j(\xi_j)\ne 0$ for $\xi_j\in 	\mathbb{T}_{\xi_j}'$, the function $Z_j(\xi_j)$ is monotone on any interval $[\alpha',\beta']\subset \mathbb{T}_{\xi_j}'$. Therefore, $\partial_{\xi_j\xi_j}Z_j(\xi_j)=-4Z_j(\xi_j)$ is also monotone on any interval $[\alpha',\beta']\subset \mathbb{T}_{\xi_j}'$.

	Since  $	|	\partial_{\xi_j\xi_j}Z_j(\xi_j)|=c\Big|\sin(\xi_1)
	\ + \  {\sigma}\sin(\xi_1+\xi^*_1)e^{{\bf i}2\xi^*_j}\Big|$ when $\xi_j\in\{\xi_j^o-\delta_{\xi_j^o},\xi_j^o+\delta_{\xi_j^o}'\}$, we deduce that $|\partial_{\xi_j\xi_j}Z_j(\xi_j)|\le c\Big|\sin(\xi_1)
	\ + \  {\sigma}\sin(\xi_1+\xi^*_1)e^{{\bf i}2\xi^*_j}\Big|$ for any $\xi_j\in 	\mathbb{T}_{\xi_j}'$. Thus, for any $\xi_j\in 	\mathbb{T}_{\xi_j}'$
	
	\begin{equation}\label{eq:degree1vertex:E7:Zj:6}
		|\cos(2\xi_j+2\upsilon)|\Big|\sin(\xi_1)
		\ + \  {\sigma}\sin(\xi_1+\xi^*_1)e^{{\bf i}2\xi^*_j}\Big|\le c\Big|\sin(\xi_1)
		\ + \  {\sigma}\sin(\xi_1+\xi^*_1)e^{{\bf i}2\xi^*_j}\Big|,
	\end{equation}
	which implies
	\begin{equation}\label{eq:degree1vertex:E7:Zj:7}
		\cos^2(2\xi_j+2\upsilon)\Big|\sin(\xi_1)
		\ + \  {\sigma}\sin(\xi_1+\xi^*_1)e^{{\bf i}2\xi^*_j}\Big|^2 \le {c^2}\Big|\sin(\xi_1)
		\ + \  {\sigma}\sin(\xi_1+\xi^*_1)e^{{\bf i}2\xi^*_j}\Big|^2.
	\end{equation}
	Thus,
	\begin{equation}\label{eq:degree1vertex:E7:Zj:8}\begin{aligned}
			&\sin^2(2\xi_j+2\upsilon)\Big|\sin(\xi_1)
			\ + \  {\sigma}\sin(\xi_1+\xi^*_1)e^{{\bf i}2\xi^*_j}\Big|^2\\
			\ge \ & {[{1-c^2}]}\Big|\sin(\xi_1)
			\ + \  {\sigma}\sin(\xi_1+\xi^*_1)e^{{\bf i}2\xi^*_j}\Big|,\end{aligned}
	\end{equation}
	yielding, by \eqref{eq:degree1vertex:E7:18},
	\begin{equation}\label{eq:degree1vertex:E7:Zj:9}\begin{aligned}
			|\partial_{\xi_j\xi_j}Z_j(\xi_j)|
			\ge \ & \frac{\sqrt{1-c^2}}{16}\Big[\Big|\sin(\tilde\xi_1^*)\sin(2\xi_1+\tilde\xi_1^*)\Big|
			+ \Big|\sin(\xi_1)\sin(2\xi_j^*)\Big|\Big],\end{aligned}
	\end{equation}
	for any $\xi_j\in 	\mathbb{T}_{\xi_j}'$. 
	
	Similar with \eqref{Lemm:Bessel2:2bb5:6}, we split 
	\begin{equation}\begin{aligned}\label{eq:degree1vertex:E7:Zj:10}
			\mathcal{Y}_j \ = \ &  \sum_{\xi_j^o\in\mathfrak{J}_{\xi_j}}\int_{\mathbb{T}_{\xi_j^o}}\mathrm{d}\xi_j e^{{\bf i}\nu Z_j} \ + \ \int_{\mathbb{T}_{\xi_j}'}\mathrm{d}\xi_j e^{{\bf i}\nu Z_j} \ = \ \sum_{\xi_j^o\in\mathfrak{J}_{\xi_j}}\mathcal{Y}_{{\xi_j^o}} \ + \ \mathcal{Y}_{j,o},
		\end{aligned}
	\end{equation}
	where $\mathfrak{J}_{\xi_j}$ is the set of stationary points of $Z_j$. Then $|\mathfrak{J}_{\xi_j}|\le 3$. 
	We remark that $\sin(\xi_1)$ appears in all $\partial_{\xi_j\xi_j}Z$. A similar phenomenon was also noticed in  \cite{lukkarinen2007asymptotics}.
	\smallskip
	
	{\bf Step 2: Stationary phase estimates of  $\mathcal{Y}_{{\xi_j^o}}$. }  
	\smallskip
	
	We now follow the same steps as in  \eqref{Lemm:Bessel2:2bb5:7}-\eqref{Lemm:Bessel2:2bb5:8}-\eqref{Lemm:Bessel2:2bb5:9}-\eqref{Lemm:Bessel2:2bb5:10}, to write

	\begin{equation}\begin{aligned}\label{eq:degree1vertex:E7:Zj:11}
			\mathcal{Y}_{{\xi_j^o}} \ = \ &    \mu_1^{\frac{1}{d}}(\xi^*,\xi_1)\int_{U_{\xi_j}}\mathrm{d}y_{\xi_j^o} e^{{\bf i}\nu Z_j(\xi_j^o)}
			\ = \   \mu_1^{\frac{1}{d}}(\xi^*,\xi_1)e^{{\bf i}\nu Z_j(\xi_j^o)}\int_{U_{\xi_j}}\mathrm{d}y_{\xi_j^o} e^{{\bf i}\sigma_{\xi_j^o}\nu|y_{\xi_j^o}|^2}\psi_{\xi_j^o}'(y_{\xi_j^o}).
		\end{aligned}
	\end{equation}
	
	By Plancherel's theorem
	\begin{equation}
		\label{eq:degree1vertex:E7:Zj:12}
		\begin{aligned}
			& \int_{U_{\xi_j}}\mathrm{d}y_{\xi_j^o} e^{{\bf i}\sigma_{\xi_j^o}\nu|y_{\xi_j^o}|^2}\psi_{\xi_j^o}'(y_{\xi_j^o}) \ = \ \frac{e^{{\bf i}\pi^2{\sigma_{\xi_j^o}}/4}}{\sqrt{|\nu|}}\int_{\mathbb{R}}\mathrm{d}k e^{-\frac{{\bf i} \pi^2}{\sigma_{\xi_j^o} \nu}k^2}\widehat{\psi_{\xi_j^o}'}(k),
		\end{aligned}
	\end{equation}
	where $\widehat{\psi_{\eta_j^*}'}$ is the  Fourier transform of $\psi_{\eta_j^*}'$.
	
	By Taylor's theorem, 
	\begin{equation}
		\label{eq:degree1vertex:E7:Zj:13}
		\begin{aligned}
			e^{-\frac{{\bf i} \pi^2}{\sigma_{\xi_j^o} \nu}k^2}
			&\ = \ &1+ \sum_{m=1}^\infty\Big(\frac{-{\bf i}\pi^2|k|^2}{\nu\sigma_{\xi_j^o}}\Big)^m,
		\end{aligned}
	\end{equation}
	uniformly in $k,|\sigma_{\xi_j^o}|$ and $\nu$. We can rewrite \eqref{eq:degree1vertex:E7:Zj:12} as
	\begin{equation}
		\label{eq:degree1vertex:E7:Zj:13}
		\begin{aligned}
			&\int_{U_{\xi_j}}\mathrm{d}y_{\xi_j^o} e^{{\bf i}\sigma_{\xi_j^o}\nu|y_{\xi_j^o}|^2}\psi_{\xi_j^o}'(y_{\xi_j^o})\   =\ \frac{e^{{\bf i}\pi^2{\sigma_{\xi_j^o}}/4}}{\sqrt{|\nu|}}\int_{\mathbb{R}}\mathrm{d}k\Big[1\ + \ \sum_{m=1}^\infty\Big(\frac{-{\bf i}\pi^2|k|^2}{\nu\sigma_{\xi_j^o}}\Big)^n\Big]\widehat{\psi_{\xi_j^o}'}(k)\\
			& \ = \ \frac{e^{{\bf i}\pi^2{\sigma_{\xi_j^o}}/4}}{\sqrt{|\nu|}}\int_{\mathbb{R}}\mathrm{d}k\widehat{\psi_{\xi_j^o}'}(k) \ +\ \frac{e^{{\bf i}\pi^2{\sigma_{\xi_j^o}}/4}}{\sqrt{|\sigma_{\xi_j^o}|\nu}}\int_{\mathbb{R}}\mathrm{d}k \sum_{m=1}^\infty\Big(\frac{-{\bf i}\pi^2k^2}{\nu\sigma_{\xi_j^o}}\Big)^m\widehat{\psi_{\xi_j^o}'}(k)\\
			& \ = \ \frac{e^{{\bf i}\pi^2{\sigma_{\xi_j^o}}/4}}{\sqrt{|\nu|}}\int_{\mathbb{R}}\mathrm{d}k\widehat{\psi_{\xi_j^o}'}(k) \ +\ \sum_{m=1}^\infty\frac{e^{{\bf i}\pi^2{\sigma_{\xi_j^o}}/4}}{\sqrt{|\nu|}}\Big(\frac{{\bf i}\pi^2}{\nu\sigma_{\xi_j^o}}\Big)^m\int_{\mathbb{R}}\mathrm{d}k \widehat{\psi_{\eta_j^*}^{(2m+1)}}(k)\\
			& \ = \ \frac{e^{{\bf i}\pi^2{\sigma_{\eta_j^*}}/4}}{\sqrt{|\nu|}}{\psi_{\xi_j^o}'}(\xi_j^o) \ +\ \sum_{m=1}^\infty\frac{e^{{\bf i}\pi^2{\sigma_{\xi_j^o}}/4}}{\sqrt{|\nu|}}\Big(\frac{{\bf i}\pi^2}{\nu\sigma_{\xi_j^o}}\Big)^m{\psi_{\xi_j^o}^{(2m+1)}}(\xi_j^o).
		\end{aligned}
	\end{equation}
	Next, we will compute explicitly the coefficients ${\psi_{\xi_j^o}'}(\xi_j^o)$ and ${\psi_{\xi_j^o}^{(2m+1)}}(\xi_j^o)$.  We first compute
	\begin{equation}
		\label{eq:degree1vertex:E7:Zj:14}
		\begin{aligned}
			\psi_{\xi_j^o}'(\xi_j)\ =\ 	& \frac{1}{\partial_{\xi_j}[(\xi_j-\xi_j^o)\mathfrak{G}_{\xi_j^o}(\xi_j-\xi_j^o)]}\ =\ 	 \frac{1}{[(\xi_j-\xi_j^o)\mathfrak{G}'_{\xi_j^o}(\xi_j-\xi_j^o)+\mathfrak{G}_{\xi_j^o}(\xi_j-\xi_j^o)]},
		\end{aligned}
	\end{equation}
	yielding
	\begin{equation}
		\label{eq:degree1vertex:E7:Zj:15}
		\begin{aligned}
			\psi_{\xi_j^o}'(\xi_j^o)\ =\ 	& 	 \frac{1}{\mathfrak{G}_{\xi_j^o}(0)}\ =\ \frac{1}{\sqrt{|\partial_{\xi_j\xi_j}Z_j(\xi_j^o)|}}.
		\end{aligned}
	\end{equation}
	By the same argument that leads to \eqref{Lemm:Bessel2:2bb5:13} and \eqref{Lemm:Bessel2:2bb5:15} , we find
	\begin{equation}
		\label{eq:degree1vertex:E7:Zj:16}
		\begin{aligned}
			\psi_{\xi_j^o}^{(n+1)}(\xi_j^o)
			=\	& \sum_{\substack{1m_1+2m_2+\cdots+nm_n=n\\ m_1,\cdots,m_n\in\mathbb{Z},m_1,\cdots,m_n\ge0}}\frac{n!}{m_1!\cdots m_n!}\\
			&\times\frac{(-1)^{m_1+\cdots+m_n}(m_1+\cdots+m_n)!}{[\mathfrak{G}_{\eta_j^*}(0)]^{m_1+\cdots+m_n+1}}\prod_{i=1}^n\left(\frac{[(i+1)\mathfrak{G}^{(i)}_{\xi_j^o}(0)]}{i!}\right)^{m_i}.
		\end{aligned}
	\end{equation}	
	and
	\begin{equation}
		\label{eq:degree1vertex:E7:Zj:17}\begin{aligned}
			\mathfrak{G}_{\xi_j^o}^{(i)}(0)\ 
			=\ &\sum_{\substack{1n_1+2n_2+\cdots+in_i=i\\n_1,\cdots,n_i\in\mathbb{Z},n_1,\cdots,n_i\ge0}}\frac{i!}{n_1!\cdots n_i!}\\
			&\times\Big(\frac12\Big)\cdots\Big(\frac32-i\Big)\Big(\sigma_{\xi_j^o}\Big)^i\big[|\partial_{\xi_j\xi_j}Z_j(\xi_j^o)|\big]^{\frac12-i}\prod_{l=1}^i\left(\frac{|\partial_{\xi_j}^{l+2}Z_j(\xi_j^o)|}{l!}\right)^{n_l}.
		\end{aligned}
	\end{equation}					
	
	Plugging \eqref{eq:degree1vertex:E7:Zj:17} into \eqref{eq:degree1vertex:E7:Zj:16}, we obtain
	\begin{equation}
		\label{eq:degree1vertex:E7:Zj:18}
		\begin{aligned}
			\psi_{\xi_j^o}^{(n+1)}(\xi_j^o)
			=\	& \sum_{\substack{1m_1+2m_2+\cdots+nm_n=n\\ m_1,\cdots,m_n\in\mathbb{Z},m_1,\cdots,m_n\ge0}}\frac{n!}{m_1!\cdots m_n!}\frac{(-1)^{m_1+\cdots+m_n}(m_1+\cdots+m_n)!}{\Big[\sqrt{|\partial_{\xi_j\xi_j}Z_j(\xi_j^o)|}\Big]^{m_1+\cdots+m_n+1}}\\
			&\times\prod_{i=1}^n\left[\sum_{\substack{1n_1+2n_2+\cdots+in_i=i\\n_1,\cdots,n_i\in\mathbb{Z},n_1,\cdots,n_i\ge0}}\frac{(i+1)}{n_1!\cdots n_i!}\right.\\
			\				&\left.\times\Big(\frac12\Big)\cdots\Big(\frac32-i\Big)\Big(\sigma_{\xi_j^o}\Big)^i\big[|\partial_{\xi_j\xi_j}Z_j(\xi_j^o)|\big]^{\frac12-i}\prod_{l=1}^i\left(\frac{|\partial_{\xi_j}^{l+2}Z_j(\xi_j^o)|}{l!}\right)^{n_l}\right]^{m_i}.
		\end{aligned}
	\end{equation}	
	As $\xi_j^o$ is a stationary point of $Z_j$, we deduce that $|\partial_{\xi_j}^{l+2}Z_j(\xi_j^o)|=2^{(l+1)}|\partial_{\xi_j}Z_j(\xi_j^o)|=0$ when $l$ is odd. Now, when $l$ is even, $|\partial_{\xi_j}^{l+2}Z_j(\xi_j^o)|=2^l|\partial_{\xi_j\xi_j}Z_j(\xi_j^o)|.$ We will show that $\prod_{l=1}^i\Big(\frac{|\partial_{\xi_j}^{l+2}Z_j(\xi_j^o)|}{l!}\Big)^{n_l}=0$ when $i$ is odd. As $i$ is odd and $i=1n_1+2n_2+\cdots+in_i$, there is an odd index $l$ such that $n_l\ne 0$, yielding $\partial_{\xi_j}^{l+2}Z_j(\xi_j^o)=0$  and thus  $\prod_{l=1}^i\Big(\frac{|\partial_{\xi_j}^{l+2}Z_j(\xi_j^o)|}{l!}\Big)^{n_l}=0$. In addition, to make sure that this product is not zero, $n_l,m_i$ need to be zero when $l,i$ are odd. 
	When $n$ is odd, $n=2m-1$, $m\in\mathbb{N}, m\ge 1$, since $1m_1+2m_2+\cdots+nm_n=n$, there exists an odd index $i$ such that $m_i$ is odd, and hence $m_i\ne 0$, leading to $
	\psi_{\xi_j^o}^{(2m)}(\xi_j^o)\ =\ 0, \forall  m\in\mathbb{N}, m\ge 1.$
	When $n=2m$, we find
	\begin{equation}
		\label{eq:degree1vertex:E7:Zj:20}
		\begin{aligned}
			&\psi_{\xi_j^o}^{(2m+1)}(\xi_j^o)
			=\	 \sum_{\substack{2m_2+4m_4+\cdots+2mm_{2m}=2m\\ m_2,\cdots,m_{2m}\in\mathbb{Z},m_2,\cdots,m_{2m}\ge0}}\frac{(2m)!}{m_1!\cdots m_{2m}!}\\
			&\times\frac{(-1)^{m_2+m_4+\cdots+m_{2m}}(m_2+m_4+\cdots+m_{2m})!}{\Big[\sqrt{|\partial_{\xi_j\xi_j}Z_j(\xi_j^o)|}\Big]^{m_2+m_4+\cdots+m_{2m}+1}}\\
			&\times\prod_{i=2, i \text{ is even}}^{2m}\left[\sum_{\substack{2n_2+4n_4+\cdots+in_i=i\\n_2,n_4,\cdots,n_i\in\mathbb{Z},n_2,n_4,\cdots,n_i\ge0}}\frac{(i+1)}{n_1!\cdots n_i!}\right.\\
			\				&\left.\times\Big(\frac12\Big)\cdots\Big(\frac32-i\Big)\Big(\sigma_{\xi_j^o}\Big)^i\big[|\partial_{\xi_j\xi_j}Z_j(\xi_j^o)|\big]^{\frac12-i}\prod_{l=2, l \text{ is even}}^i\left(\frac{|2^l\partial_{\xi_j\xi_j}Z_j(\xi_j^o)|}{l!}\right)^{n_l}\right]^{m_i}.
		\end{aligned}
	\end{equation}	
	Similarly as in   \eqref{Lemm:Bessel2:2bb5:18}, we can obtain the bound
	\begin{equation}
		\label{eq:degree1vertex:E7:Zj:21}
		\begin{aligned}
			&|\psi_{\xi_j^o}^{(2m+1)}(\xi_j^o)|
			\le\	 \left|\sum_{\substack{2m_2+4m_4+\cdots+2mm_{2m}=2m\\ m_2,\cdots,m_{2m}\in\mathbb{Z},m_2,\cdots,m_{2m}\ge0}}\frac{(2m)!(m_2+m_4+\cdots+m_{2m}+1)!}{m_1!\cdots m_{2m}!}\right.\\
			&\times\prod_{i=2, i \text{ is even}}^{2m}\left.\left[\sum_{\substack{2n_2+4n_4+\cdots+in_i=i\\n_2,n_4,\cdots,n_i\in\mathbb{Z},n_2,n_4,\cdots,n_i\ge0}}\frac{(i+1)}{n_1!\cdots n_i!}\Big|\frac12\Big|\cdots\Big|\frac32-i\Big|\prod_{l=2, l \text{ is even}}^i\left(\frac{2^l}{l!}\right)^{n_l}\right]^{m_i}\right|\\
			&\times |\mathfrak{C}^1_{\eta_j^*}|^{m}\left|\sum_{\substack{2m_2+4m_4+\cdots+2mm_{2m}=2m\\ m_2,\cdots,m_{2m}\in\mathbb{Z},m_2,\cdots,m_{2m}\ge0}}\Big\{\prod_{i=2, i \text{ is even}}^{2m}\Big[\sum_{\substack{2n_2+4n_4+\cdots+in_i=i\\n_2,n_4,\cdots,n_i\in\mathbb{Z},n_2,n_4,\cdots,n_i\ge0}}\big[|\partial_{\xi_j\xi_j}Z_j(\xi_j^o)|\big]^{(\frac12-i)m_i}\right.\\
			&\left.\  \times{\prod_{l=2, l \text{ is even}}^i}|\partial_{\xi_j\xi_j}Z_j(\xi_j^o)|^{n_lm_i}\Big[\sqrt{|\partial_{\xi_j\xi_j}Z_j(\xi_j^o)|}\Big]^{-m_2-m_4-\cdots-m_{2m}-1}\Big]\Big\}\right|,
		\end{aligned}
	\end{equation}	
	for some constant $\mathfrak{C}^1_{\xi_j^o}>0$ independent of $m$.
	Following the same lines of computations of \eqref{Lemm:Bessel2:2bb5:20}-\eqref{Lemm:Bessel2:2bb5:21}, we deduce the existence of an explicit constant $\mathfrak{C}^2_{\xi_j^o}>0$ independent of $m$, such that the constant in front of the term containing only $Z_j(\xi_j^o)$  in \eqref{eq:degree1vertex:E7:Zj:21} is bounded by $  |\mathfrak{C}^2_{\xi_j^o}|^m$.
	We only need to estimate the term containing only $Z_j(\xi_j^o)$  in \eqref{eq:degree1vertex:E7:Zj:21}. This quantity can be simplified as
	\begin{equation}
		\label{eq:degree1vertex:E7:Zj:22}				\begin{aligned}
			&	\left|\sum_{\substack{2m_2+4m_4+\cdots+2mm_{2m}=2m\\ m_2,\cdots,m_{2m}\in\mathbb{Z},m_2,\cdots,m_{2m}\ge0}}\Big\{\prod_{i=2, i \text{ is even}}^{2m}\Big[\sum_{\substack{2n_2+4n_4+\cdots+in_i=i, i \text{ is even}\\n_2,n_4,\cdots,n_i\in\mathbb{Z},n_2,n_4,\cdots,n_i\ge0}}\big[|\partial_{\xi_j\xi_j}Z_j(\xi_j^o)|\big]^{(\frac12-i)m_i}\right.\\
			&\ \ \ \ \ \ \ \ \left.\  \times{\prod_{l=2, l \text{ is even}}^i}|\partial_{\xi_j\xi_j}Z_j(\xi_j^o)|^{n_lm_i}\Big[\sqrt{|\partial_{\xi_j\xi_j}Z_j(\xi_j^o)|}\Big]^{-m_2-m_4-\cdots-m_{2m}-1}\Big]\Big\}\right|\\
			\le \	&|\mathfrak{C}^1_{1,\xi_j^o}|^m\left|\sum_{\substack{2m_2+4m_4+\cdots+2mm_{2m}=2m\\ m_2,\cdots,m_{2m}\in\mathbb{Z},m_2,\cdots,m_{2m}\ge0}}\left(\sum_{\substack{2n_2+4n_4+\cdots+in_i=i\\n_2,n_4,\cdots,n_i\in\mathbb{Z},n_2,n_4,\cdots,n_i\ge0}}\right.\right.\\
			&\ \ \ \ \ \ \ \ \ \ \ \ \ \ \ \ \left.\left.\left[\frac{1}{\big|\partial_{\xi_j\xi_j}Z_j(\xi_j^o)\big|}\right]^{\sum_{i=2, i \text{ is even}}^{2m}(-\frac12+i-n_2-\cdots-n_i)m_i+\frac{m_2+\cdots+m_{2m}+1}{2}}\right)\right|,				\end{aligned}
	\end{equation}
	for some constant $\mathfrak{C}^1_{1,\xi_j^o}>0$.
	Repeating the argument used to bound  \eqref{Lemm:Bessel2:2bb5:20}-\eqref{Lemm:Bessel2:2bb5:21}, we deduce the existence of a constant  $\mathfrak{C}^2_{1,\xi_j^o}>0$ such that
	\begin{equation}
		\label{eq:degree1vertex:E7:Zj:23}		
		|\psi_{\xi_j^o}^{(2m+1)}(\xi_j^o)|\le 	(\mathfrak{C}^2_{1,\xi_j^o})^m\left[\frac{1}{\big|\partial_{\xi_j\xi_j}Z_j(\xi_j^o)\big|}\right]^{{2m+1}}.				
	\end{equation}
	We then deduce from \eqref{eq:degree1vertex:E7:Zj:13}
	\begin{equation}
		\label{eq:degree1vertex:E7:Zj:24}
		\begin{aligned}
			&\Big|\int_{U_{\xi_j}}\mathrm{d}y_{\xi_j^o} e^{{\bf i}\sigma_{\xi_j^o}\nu|y_{\xi_j^o}|^2}\psi_{\xi_j^o}'(y_{\xi_j^o})\Big|\   \le \ \frac{1}{\sqrt{|\nu||\partial_{\xi_j\xi_j}Z_j(\xi_j^o)|}}+\sum_{m=1}^\infty\frac{(\mathfrak{C}^3_{\xi_j^o})^m}{\big|\partial_{\xi_j\xi_j}Z_j(\xi_j^o)\big|^{{2m+1}}|\nu|^{\frac{2m+1}{2}}},
		\end{aligned}
	\end{equation}
	for some explicit constant $\mathfrak{C}^3_{\xi_j^o}>0$, which, in combination with \eqref{eq:degree1vertex:E7:Zj}, yields
	\begin{equation}
		\label{eq:degree1vertex:E7:Zj:25}
		\begin{aligned}
			\mathcal{Y}_{\xi_j^o}\   \le \ &\frac{\mu_1^{\frac{1}{d}}}{\sqrt{|\nu|\Big|\frac18\Big|\sin(\tilde\xi_1^*)\sin(2\xi_1+\tilde\xi_1^*)\Big|
					\ + \ \frac18\Big|\sin(\xi_1)\sin(2\xi_j^*)\Big|\Big|}}\\
			&\ \ +\sum_{m=1}^\infty\frac{\mu_1^{\frac{1}{d}}(\mathfrak{C}^3_{\xi_j^o})^m}{\Big|\frac18\Big|\sin(\tilde\xi_1^*)\sin(2\xi_1+\tilde\xi_1^*)\Big|
				\ + \ \frac18\Big|\sin(\xi_1)\sin(2\xi_j^*)\Big|\Big|^{{2m+1}}|\nu|^{\frac{2m+1}{2}}}.
		\end{aligned}
	\end{equation}
	
	Now, 
	observing  \begin{equation}
		\label{eq:degree1vertex:E7:Zj:25:a}\Big|\sin(\tilde\xi_1^*)\sin(2\xi_1+\tilde\xi_1^*)\Big|
		+  \Big|\sin(\xi_1)\sin(2\xi_j^*)\Big|\gtrsim|d(\mathfrak{S})|^2,\end{equation}
	we deduce from \eqref{eq:degree1vertex:E7:Zj:25} that
	\begin{equation}
		\label{eq:degree1vertex:E7:Zj:26}
		\begin{aligned}
			\mathcal{Y}_{\xi_j^o}\   \le \ &\frac{1}{\sqrt{|\nu||d(\mathfrak{S})|^2}}\ +\ \sum_{m=1}^\infty\frac{(\mathfrak{C}^3_{\xi_j^o})^m}{|d(\mathfrak{S})|^{{2(2m+1)}}|\nu|^{\frac{2m+1}{2}}}.
		\end{aligned}
	\end{equation}
	
	{\bf Step 3: Estimating $\mathcal{Y}_{j}$ in \eqref{yj}}.
	
	\smallskip
	
	Using \eqref{eq:degree1vertex:E7:Zj:10} we first perform an integration by parts on $\mathcal{Y}_{j,o}$
	
	\begin{equation}\begin{aligned}
			\label{eq:degree1vertex:E7:Zj:27}
			\mathcal{Y}_{j,o}\ = \ & \mu_1^{\frac{1}{d}}(\xi^*,\xi_1)\int_{\mathbb{T}_{\xi_j}'}\mathrm{d}\xi_j e^{{\bf i}\nu Z_j} \ = \ \mu_1^{\frac{1}{d}}(\xi^*,\xi_1)\int_{\mathbb{T}_{\xi_j}'}\mathrm{d}\xi_j \frac{1}{{\bf i}\nu\partial_{\xi_j}Z_j}\partial_{\xi_j}\big(e^{{\bf i}\nu Z_j}\big) \\
			\ = \ & -\mu_1^{\frac{1}{d}}(\xi^*,\xi_1)\int_{\mathbb{T}_{\xi_j}'}\mathrm{d}\xi_j \partial_{\xi_j}\Big(\frac{1}{{\bf i}\nu\partial_{\xi_j}Z_j}\Big)e^{{\bf i}\nu Z_j} \ + \ \mu_1^{\frac{1}{d}}(\xi^*,\xi_1) \frac{e^{{\bf i} \nu Z_j} }{{\bf i}\nu\partial_{\xi_j}Z_j}\Big|_{\partial\mathbb{T}_{\xi_j}'}\\
			\ = \ & \mu_1^{\frac{1}{d}}(\xi^*,\xi_1){{\bf i}}\int_{\mathbb{T}_{\xi_j}'}\mathrm{d}\xi_j \frac{\partial_{\xi_j\xi_j}Z_j}{\nu|\partial_{\xi_j}Z_j|^2}e^{{\bf i}\nu Z_j} \ + \ \mu_1^{\frac{1}{d}}(\xi^*,\xi_1)\frac{e^{{\bf i}\nu Z_j} }{{\bf i}\nu\partial_{\xi_j}Z_j}\Big|_{\partial\mathbb{T}_{\xi_j}'},
		\end{aligned}
	\end{equation}
	where we have used \eqref{eq:degree1vertex:E7:3:a:1}.
	Using \eqref{eq:degree1vertex:E7:Zj:9}, we bound
	\begin{equation}\begin{aligned}
			\label{eq:degree1vertex:E7:Zj:28}
			|\mathcal{Y}_{j,o}| 
			\ \le \ & \mu_1^{\frac{1}{d}}(\xi^*,\xi_1)\left|{{\bf i}}\int_{\mathbb{T}_{\xi_j}'}\mathrm{d}\xi_j \frac{\partial_{\xi_j\xi_j}Z_j}{\nu|\partial_{\xi_j}Z_j|^2}e^{{\bf i}\nu Z_j} \right| \ + \ \mu_1^{\frac{1}{d}}(\xi^*,\xi_1)\left|\frac{e^{{\bf i}\nu Z_j} }{{\bf i}\nu\partial_{\xi_j}Z_j}\Big|_{\partial\mathbb{T}_{\xi_j}'}\right|\\
			\ \lesssim \  & \frac{\mu_1^{\frac{1}{d}}(\xi^*,\xi_1)}{|\nu|\Big|\frac{\sqrt{1-c^2}}{16}\Big|\sin(\tilde\xi_1^*)\sin(2\xi_1+\tilde\xi_1^*)\Big|
				\ + \ \frac{\sqrt{1-c^2}}{16}\Big|\sin(\xi_1)\sin(2\xi_j^*)\Big|\Big|^2} \\
			& \ + \ \frac{\mu_1^{\frac{1}{d}}(\xi^*,\xi_1)}{|\nu|\Big|\frac{\sqrt{1-c^2}}{16}\Big|\sin(\tilde\xi_1^*)\sin(2\xi_1+\tilde\xi_1^*)\Big|
				\ + \ \frac{\sqrt{1-c^2}}{16}\Big|\sin(\xi_1)\sin(2\xi_j^*)\Big|\Big|},
		\end{aligned}
	\end{equation}
	where the constant on the right hand side is explicit. Inequality \eqref{eq:degree1vertex:E7:Zj:28}, together with \eqref{eq:degree1vertex:E7:Zj:25:a}, yields
	\begin{equation}\begin{aligned}
			\label{eq:degree1vertex:E7:Zj:29}
			|\mathcal{Y}_{j,o}| 
			\ \lesssim \  & \frac{1}{|\nu||d(\mathfrak{S})|^4
			}  + \ \frac{1}{|\nu||d(\mathfrak{S})|^2}.
		\end{aligned}
	\end{equation}
	
	Combining \eqref{eq:degree1vertex:E7:Zj:26} and 
	\eqref{eq:degree1vertex:E7:Zj:29}, we find
	\begin{equation}\begin{aligned}
			\label{eq:degree1vertex:E7:Zj:30:1}
			|\mathcal{Y}_{j}| 
			\ \lesssim \  & \frac{1}{|\nu||d(\mathfrak{S})|^4}  + \ \frac{1}{|\nu||d(\mathfrak{S})|^2} \ + \ \frac{1}{\sqrt{|\nu||d(\mathfrak{S})|^2}}\ +\ \sum_{m=1}^\infty\frac{(\mathfrak{C}^3_{\xi_j^o})^m}{|d(\mathfrak{S})|^{{2(2m+1)}}|\nu|^{\frac{2m+1}{2}}},
		\end{aligned}
	\end{equation}
	where the constants on the right hand side are explicit. By the trivial bound $|\mathcal{Y}_j|\le 2\pi$, we deduce from \eqref{eq:degree1vertex:E7:Zj:30:1} that
	\begin{equation}\begin{aligned}
			\label{eq:degree1vertex:E7:Zj:30}
			|\mathcal{Y}_{j}| 
			\ \lesssim \  & \min\left\{1,\frac{1}{|\nu||d(\mathfrak{S})|^4}  + \ \frac{1}{|\nu||d(\mathfrak{S})|^2} \ + \ \frac{1}{\sqrt{|\nu||d(\mathfrak{S})|^2}}\ +\ \sum_{m=1}^\infty\frac{(\mathfrak{C}^3_{\xi_j^o})^m}{|d(\mathfrak{S})|^{{2(2m+1)}}|\nu|^{\frac{2m+1}{2}}}\right\}.
		\end{aligned}
	\end{equation}
	
	{\bf Step 4: Final estimate.}
	
	We deduce from \eqref{eq:degree1vertex:E7:3:a} that

	\begin{equation}	\label{eq:degree1vertex:E7:Zj:32}
		\begin{aligned}
			\tilde{Y}	:=\ & 	\left|\int_{\mathbb{T}^d}\mathrm{d}k \mu_1 e^{{\bf i}\nu(\omega(k)+\sigma\omega(k^*+k))}\right| \ = \ 	\left|\int_{[-\pi,\pi]}\mathrm{d}\xi_1 \mu_1 e^{{\bf i}\nu Z_1}\prod_{j=2}^d\int_{[-\pi,\pi]}\mathrm{d}\xi_j e^{{\bf i}\nu Z_j}\right| \\
			\ \le \	& \left|\int_{[-\pi,\pi]}\mathrm{d}\xi_1\prod_{j=2}^d|\mathcal{Y}_j|\right|\\ \ \lesssim \ & \min\left\{1,\frac{1}{|\nu||d(\mathfrak{S})|^4}  + \ \frac{1}{|\nu||d(\mathfrak{S})|^2} \ + \ \frac{1}{\sqrt{|\nu||d(\mathfrak{S})|^2}}\ +\ \sum_{m=1}^\infty\frac{(\mathfrak{C}^3_{\xi_j^o})^m}{|d(\mathfrak{S})|^{{2(2m+1)}}|\nu|^{\frac{2m+1}{2}}}\right\},
		\end{aligned}
	\end{equation}
	in which, we have used \eqref{eq:degree1vertex:E7:Zj:30} for $\mathcal{Y}_2$. For $\mathcal{Y}_j$ with $3\le j\le d$, we have simply used the bound $|\mathcal{Y}_j|\le 2\pi$.
	We bound
	\begin{equation}
		\label{eq:degree1vertex:E10}
		\begin{aligned}
			&  \left(\int_{0}^{1}\mathrm{d}\nu  \tilde{Y}+2\int_1^{|d(\mathfrak{S})|^{-2}|\ln(\lambda)|^{\epsilon}}\mathrm{d}\nu{\tilde{Y}}^\frac{q}{2}+2\int_{|d(\mathfrak{S})|^{-2}|\ln|\lambda||^{\epsilon}}^\infty\mathrm{d}\nu{\tilde{Y}}^\frac{q}{2}\right)^\frac1q\\
			\lesssim \ 	& \left[{\langle\mathrm{ln}|\lambda|^2\rangle} \int_{-1}^{1}\mathrm{d}\nu  +2\int_1^{|d(\mathfrak{S})|^{-2}|\ln|\lambda||^{\epsilon}}\mathrm{d}\nu +2\int_{|d(\mathfrak{S})|^{-2}|\ln|\lambda||^{\epsilon}}^\infty\mathrm{d}\nu\right.\\
			&\left.\times\min\left\{1,\frac{1}{|\nu||d(\mathfrak{S})|^4}  + \ \frac{1}{|\nu||d(\mathfrak{S})|^2} \ + \ \frac{1}{\sqrt{|\nu||d(\mathfrak{S})|^2}}\ +\ \sum_{m=1}^\infty\frac{(\mathfrak{C}^3_{\xi_j^o})^m}{|d(\mathfrak{S})|^{{2(2m+1)}}|\nu|^{\frac{2m+1}{2}}}\right\}^\frac{q}{2}\right]^\frac1q,
		\end{aligned}
	\end{equation}
	where we have used the bounds $\tilde{Y}\le 1$ for the integral from $-1$ to $1$. The parameter $\epsilon>0$ is small. The first and second integrals can simply be bounded as
	\begin{equation}
		\label{eq:degree1vertex:E10:1}
		\begin{aligned}
			\int_{0}^{1}\mathrm{d}\nu  +2\int_1^{|d(\mathfrak{S})|^{-2}|\ln|\lambda||^{\epsilon}}\mathrm{d}\nu \
			\lesssim	&\ 1 + |d(\mathfrak{S})|^{-2}|\ln|\lambda||^{\epsilon}.
		\end{aligned}
	\end{equation}
	We now estimate the last integral by setting $q=4$. By the change of variable $\nu\to \bar{\nu}=|d(\mathfrak{S})|^{2}|\ln(\lambda)|^{-\epsilon}\nu$, we find
	\begin{equation}
		\label{eq:degree1vertex:E11}
		\begin{aligned}
			&\int_{|d(\mathfrak{S})|^{-2}|\ln|\lambda||^{\epsilon}}^\infty\mathrm{d}\nu\left\{\frac{1}{|\nu||d(\mathfrak{S})|^4}  +  \frac{1}{|\nu||d(\mathfrak{S})|^2}\right.\\
			&\ \ \ \ \ \ \ \ \ \ \ \ \ \ \ \ \ \  \ \ \ \ \ \ \ \ \ \ \ \ \ \ \ \ \ \  \left. +  \frac{1}{\sqrt{|\nu||d(\mathfrak{S})|^2}} + \sum_{m=1}^\infty\frac{(\mathfrak{C}^3_{\xi_j^o})^m}{|d(\mathfrak{S})|^{{2(2m+1)}}|\nu|^{\frac{2m+1}{2}}}\right\}^2\\
			=\	&\int_{1}^\infty|d(\mathfrak{S})|^{-2}|\ln|\lambda||^{\epsilon}\mathrm{d}\bar{\nu}\left\{\frac{|\ln|\lambda||^{-\epsilon}}{\bar{\nu}|d(\mathfrak{S})|^2}  + \frac{|\ln|\lambda||^{-\epsilon}}{\bar{\nu}} +  \frac{|\ln|\lambda||^{-\epsilon/2}}{\sqrt{\bar{\nu}}} \right.\\
			&\ \ \ \ \ \ \ \ \ \ \ \ \ \ \ \ \ \  \ \ \ \ \ \ \ \ \ \ \ \ \ \ \ \ \ \  \left. + \sum_{m=1}^\infty\frac{(\mathfrak{C}^3_{\xi_j^o})^m|\ln|\lambda||^{-(2m+1)\epsilon/2}}{[\bar{\nu}|d(\mathfrak{S})|^2]^{\frac{2m+1}{2}}}\right\}^2
			.
		\end{aligned}
	\end{equation}
	Observing that $\lambda <<|d(\mathfrak{S})|^2,$  and suppose that $\lambda$ is sufficiently small, we deduce 
	\begin{equation}
		\label{eq:degree1vertex:E12}	\begin{aligned}
			&\frac{|\ln|\lambda||^{-\epsilon}}{\bar{\nu}|d(\mathfrak{S})|^2}  + \frac{|\ln|\lambda||^{-\epsilon}}{\bar{\nu}} +  \frac{|\ln|\lambda||^{-\epsilon/2}}{\sqrt{\bar{\nu}}} +\sum_{m=1}^\infty\frac{(\mathfrak{C}^3_{\xi_j^o})^m|\ln|\lambda||^{-(2m+1)\epsilon/2}}{[\bar{\nu}|d(\mathfrak{S})|^2]^{\frac{2m+1}{2}}}  \\
			=\ &\frac{|\ln|\lambda||^{-\epsilon}}{\bar{\nu}|d(\mathfrak{S})|^2}  + \frac{|\ln|\lambda||^{-\epsilon}}{\bar{\nu}} +  \frac{|\ln|\lambda||^{-\epsilon/2}}{\sqrt{\bar{\nu}}} +\mathcal{O}(|\ln|\lambda||^{-3\epsilon/2}[\bar{\nu}|d(\mathfrak{S})|^2]^{-\frac32})\\
			\lesssim \ &	\frac{|\ln|\lambda||^{-\epsilon}}{\bar{\nu}|d(\mathfrak{S})|^2}  +  \frac{|\ln|\lambda||^{-\epsilon/2}}{\sqrt{\bar{\nu}}}.\end{aligned}
	\end{equation}
	Combining \eqref{eq:degree1vertex:E11} and \eqref{eq:degree1vertex:E12} yields
	\begin{equation}
		\label{eq:degree1vertex:E13}
		\begin{aligned}
			&\int_{|d(\mathfrak{S})|^{-2}|\ln(\lambda)|^{\epsilon}}^\infty\mathrm{d}\nu\left\{\frac{1}{|\nu||d(\mathfrak{S})|^4}  +  \frac{1}{|\nu||d(\mathfrak{S})|^2}\right.\\
			&\ \ \ \ \ \ \ \ \ \ \ \ \ \ \ \ \ \  \ \ \ \ \ \ \ \ \ \ \ \ \ \ \ \ \ \  \left. +  \frac{1}{\sqrt{|\nu||d(\mathfrak{S})|^2}} + \sum_{m=1}^\infty\frac{(\mathfrak{C}^3_{\xi_j^o})^m}{|d(\mathfrak{S})|^{{2m+1}}|\nu|^{\frac{2m+1}{2}}}\right\}^2\\
			\lesssim \	&\int_{1}^\infty\mathrm{d}\bar{\nu}|d(\mathfrak{S})|^{-2}|\ln|\lambda||^{\epsilon}\left\{\frac{|\ln|\lambda||^{-\epsilon}}{\bar{\nu}|d(\mathfrak{S})|^2}  +  \frac{|\ln|\lambda||^{-\epsilon/2}}{\sqrt{\bar{\nu}}}\right\}^2
			\\
			\lesssim \	&\int_{1}^\infty\mathrm{d}\bar{\nu}|d(\mathfrak{S})|^{-2}|\ln|\lambda||^{\epsilon}\left\{\frac{|\ln|\lambda||^{-2\epsilon}}{\bar{\nu}^2|d(\mathfrak{S})|^4}  +  \frac{|\ln|\lambda||^{-\epsilon}}{{\bar{\nu}}}\right\}\ \lesssim \ \frac{{\langle\mathrm{ln}|\lambda|\rangle}^2}{|d(\mathfrak{S})|^6}.
		\end{aligned}
	\end{equation}
	We finally obtain the conclusion of the lemma.
\end{proof}
\subsection{The role of the kernel on the dispersive estimates}\label{KernelDispersive}  We follow the same notations used in Section \ref{Subsec:DispersiveEstimates}.
Let us consider the dispersion relation
\eqref{Nearest} with its equivalent form   \eqref{Nearest1}. We define a similar function as   \eqref{Improved:extendedBesselfunctions1}, but the cut-off function $\check{\Phi}$ is replaced by the kernel $|\sin(\xi_1)|^n$, $n\ge 0, n\in\mathbb{N}$.    The function reads
\begin{equation}
	\label{Ker:extendedBesselfunctions1}
	\mathfrak{F}^{Ker}(m,t_0) \ = \ \int_{[-\pi,\pi]^d} \mathrm{d}\xi |\sin(\xi_1)|^ne^{{\bf i}m\cdot \xi} e^{{\bf i}t_0\omega(\xi)},
\end{equation}
for $\xi=(\xi_1,\cdots,\xi_d)\in\mathbb{R}^d$, $m=(m_1,\cdots,m_d)\in\mathbb{Z}^d$, $t_0\in\mathbb{R}$. We also suppose that $d\ge 2$. 
\begin{lemma}\label{Lemm:Ker:Bessel2} There exists a universal constant $\mathfrak{C}_{\mathfrak{F}^{Ker},2}$ independent of $t_0$ and $\lambda$, such that, 
	\begin{equation}
		\label{Lemm:Ker:Bessel2:1}
		\|\mathfrak{F}^{Ker}(\cdot,t_0) \|_{l^2} \ = \ \left( \sum_{m\in\mathbb{Z}^d}|\mathfrak{F}^{Ker}(m,t_0) |^2\right)^\frac12 \ \le \ \mathfrak{C}_{\mathfrak{F}^{Ker},2}.
	\end{equation}
\end{lemma}
\begin{proof} 	The proof of the lemma follows precisely the same argument as the one for  Lemma \ref{Lemm:Improved:Bessel2}.
\end{proof}
\begin{lemma}\label{Lemm:Ker:Bessel3}
	There exists a  universal constant $\mathfrak{C}_{\mathfrak{F}^{Ker},4}>0$ independent of $t_0$ and $\lambda$, such that when $2n+5<d$,
	\begin{equation}
		\label{Lemm:Ker:Bessel3:1}\begin{aligned}
			&\|\mathfrak{F}^{Ker}(\cdot,t_0) \|_{l^4}  \le \mathfrak{C}_{\mathfrak{F}^{Ker},4}\langle t_0\rangle^{-{\frac{2n}{7}-\frac{13}{30}}},\end{aligned}
	\end{equation}
	otherwise
	\begin{equation}
		\label{Lemm:Ker:Bessel3:2}\begin{aligned}
			&\|\mathfrak{F}^{Ker}(\cdot,t_0) \|_{l^4}  \le \mathfrak{C}_{\mathfrak{F}^{Ker},4}\langle t_0\rangle^{-{\frac{2n}{7}}},\end{aligned}
	\end{equation}
\end{lemma}
\begin{proof}
	We observe that
	\begin{equation}
		\label{Lemm:Ker:Bessel2:1a}
		\begin{aligned}
			|\mathfrak{F}^{Ker}(m,t_0) |^2
			\ =  \ & \int_{[-\pi,\pi]^d}\mathrm{d}\xi e^{{\bf i}m\cdot\xi} \int_{-\pi}^{\pi}\mathrm{d}\eta_1 |\sin(\xi_1+\eta_1)|^n|\sin(\eta_1)|^n\\
			&\times\exp\Big({\bf i}t_0\sin^3(\xi_1+\eta_1) -{\bf i}t_0\sin^3(\eta_1)\Big)\\
			& \times\Big[\prod_{j=2}^d\int_{-\pi}^\pi\mathrm{d}\eta_j \exp\Big({\bf i}t_0\sin(\xi_1+\eta_1)\sin^2(\xi_j+\eta_j)- {\bf i}t_0\sin(\eta_1)\sin^2(\eta_j)\Big)\Big], 
		\end{aligned}
	\end{equation} 
	and  the Plancherel theorem
	\begin{equation}
		\label{Lemm:Ker:Bessel2:2}\begin{aligned}
			& \|\mathfrak{F}^{Ker}(\cdot,t_0) \|_{l^4}^4\ = \ \sum_{m\in\mathbb{Z}^d}|\mathfrak{F}^{Ker}(m,t_0) |^4\\  
			\ = \ & \int_{[-\pi,\pi]^d}\mathrm{d}\xi \Big| \int_{-\pi}^{\pi}\mathrm{d}\eta_1 |\sin(\xi_1+\eta_1)|^n|\sin(\eta_1)|^n \exp\Big({\bf i}t_0\sin^3(\xi_1+\eta_1)-{\bf i}t_0\sin^3(\eta_1)\Big)\\
			& \  \times\Big[\prod_{j=2}^d\int_{-\pi}^\pi\mathrm{d}\eta_j \exp\Big({\bf i}t_0\sin(\xi_1+\eta_1)\sin^2(\xi_j+\eta_j)- {\bf i}t_0\sin(\eta_1)\sin^2(\eta_j)\Big)\Big]\Big|^2.\end{aligned}
	\end{equation}
	Similarly  as  the proof of Lemma \ref{Lemm:Improved:Bessel3}, we also define for $2\le j\le d$
	\begin{equation}
		\begin{aligned}
			\mathfrak{A}_j(\xi_1,\eta_1,\xi_j) \ = & \ \int_{-\pi}^\pi\mathrm{d}\eta_j e^{{\bf i}\mathfrak{B}_j(\xi_1,\eta_1,\xi_j,\eta_j)},
		\end{aligned}
	\end{equation}
	in which
	\begin{equation}
		\begin{aligned}
			& \mathfrak{B}_j(\xi_1,\eta_1,\xi_j,\eta_j) =  \ t_0\sin(\xi_1+\eta_1)\sin^2(\xi_j+\eta_j)
			\ - \ t_0\sin(\eta_1)\sin^2(\eta_j),
		\end{aligned}
	\end{equation}
	and
	\begin{equation}
		\begin{aligned}
			& \mathfrak{B}_1(\xi_1,\eta_1) =  \ t_0\sin^3(\xi_1+\eta_1) -t_0\sin^3(\eta_1),
		\end{aligned}
	\end{equation}
	we find
	\begin{equation}
		\label{Lemm:Ker:Bessel2:2a}\begin{aligned}
			\sum_{m\in\mathbb{Z}^d}|\mathfrak{F}^{Ker}(m,t_0) |^4
			\ = \ & \int_{[-\pi,\pi]^d}\mathrm{d}\xi  \Big|\int_{-\pi}^{\pi}\mathrm{d}\eta_1 \prod_{j=2}^d\mathfrak{A}_j(\xi_1,\eta_1,\xi_j)\\
			&\times |\sin(\xi_1+\eta_1)|^n|\sin(\eta_1)|^n e^{{\bf i}\mathfrak{B}_1(\xi_1,\eta_1) }\Big|^2.\end{aligned}
	\end{equation}
	The phase $\mathfrak{B}_j$ is the sum of  
	\begin{equation}
		\begin{aligned}
			& \mathfrak{B}^a_j(\xi_1,\eta_1,\xi_j,\eta_j)
			\ =   \ -\ \mathrm{Re}\Big[e^{{\bf i}2\xi_j+{\bf i}2\eta_j}\mathfrak{C}_j^1\Big]\ + \ \mathrm{Re}\Big[e^{{\bf i}2\eta_j}\mathfrak{C}_j^2\Big],
		\end{aligned}
	\end{equation}
	with
	$$\mathfrak{C}_j^1=\frac{t_0}{2}\sin(\xi_1+\eta_1), \mbox{ and } \mathfrak{C}_j^2=\frac{t_0}{2}\sin(\eta_1),$$
	and
	\begin{equation}
		\begin{aligned}
			&\mathfrak{B}^b_j(\xi_1,\eta_1)
			=   \frac{t_0}{2}\sin(\xi_1+\eta_1)    
			-  \frac{t_0}{2}\sin(\eta_1).
		\end{aligned}
	\end{equation}
	The oscillatory integral $\mathfrak{A}_j$ is now  written
	\begin{equation}
		\mathfrak{A}_j \ = \ e^{{\bf i}\mathfrak{B}^b_j }\int_{-\pi}^\pi\mathrm{d}\eta_j e^{{\bf i}\mathfrak{B}^a_j} \ = \  e^{{\bf i}\mathfrak{B}^b_j }\mathfrak{A}_j^a.
	\end{equation}
	We also
	combine the phases ${\mathfrak{B}^b_j }$ and ${\mathfrak{B}_1}$  
	\begin{equation}
		\begin{aligned}
			& \mathfrak{B}^a_1=\mathfrak{B}_1 + \sum_{j=2}^{d}\mathfrak{B}^b_j \  =  \ t_0\Big[\sin^3(\xi_1+\eta_1)+\frac{d-1}{2}\sin(\xi_1+\eta_1)\Big]\\
			& -t_0\Big[\sin^3(\eta_1)+\frac{d-1}{2}\sin(\eta_1)\Big].
		\end{aligned}
	\end{equation}
	And finally, we also write
	\begin{equation}
		\label{Lemm:Ker:Bessel2:2a:a}\begin{aligned}
			\sum_{m\in\mathbb{Z}^d}|\mathfrak{F}^{Ker}(m,t_0) |^4
			\ = \ & \int_{[-\pi,\pi]^d}\mathrm{d}\xi  \Big|\int_{-\pi}^{\pi}\mathrm{d}\eta_1 \prod_{j=2}^d\mathfrak{A}_j^a(\xi_1,\eta_1,\xi_j)\\
			&\times |\sin(\xi_1+\eta_1)|^n|\sin(\eta_1)|^ne^{{\bf i}\mathfrak{B}_1^a(\xi_1,\eta_1) }\Big|^2.\end{aligned}
	\end{equation}
	We now divide the proof into smaller steps.
	
	\smallskip
	
	{\bf Step 1: Estimating $|\sin(\xi_1+\eta_1)|\mathfrak{A}_j^a(\xi_1,\eta_1,\xi_j)$}. We develop 
	\begin{equation}
		\begin{aligned}\label{Lemm:Ker:Bessel2:3}
			\mathfrak{B}^a_j(\xi_1,\eta_1,\xi_j,\eta_j)\ =  & \ -\frac{t_0}{2}\Big[\sin(\xi_1+\eta_1)\cos(2\xi_j+2\eta_j)\ -\sin(\eta_1)\cos(2\eta_j)\Big]\\
			\ =  & \ -\frac{t_0}{2}\Big[\sin(\xi_1+\eta_1)\cos(2\xi_j)\cos(2\eta_j)-\sin(\xi_1+\eta_1)\sin(2\xi_j)\sin(2\eta_j)\\
			&\ -\sin(\eta_1)\cos(2\eta_j)\Big]\\
			\ = &\ -t'\sin(\xi_1+\eta_1)\sin(2\xi_j)\sin(2\eta_j) + t'[\sin(\xi_1+\eta_1)\cos(2\xi_j)\\
			&\ -\sin(\eta_1)]\cos(2\eta_j),
		\end{aligned}
	\end{equation}
	in which $t'=-t_0/2$.

	Next, we set 
	\begin{equation}\begin{aligned}\label{Lemm:Ker:Bessel2:4}
			\cos(G_{j}) \ = \ & \frac{-\sin(\xi_1+\eta_1)\sin(2\xi_j)}{\sqrt{\sin^2(\xi_1+\eta_1)\sin^2(2\xi_j)+[\sin(\xi_1+\eta_1)\cos(2\xi_j)-\sin(\eta_1)]^2}},\\
			\sin(G_{j}) \ = \ & \frac{\sin(\xi_1+\eta_1)\cos(2\xi_j)-\sin(\eta_1)}{\sqrt{\sin^2(\xi_1+\eta_1)\sin^2(2\xi_j)+[\sin(\xi_1+\eta_1)\cos(2\xi_j)-\sin(\eta_1)]^2}},
		\end{aligned}
	\end{equation}
	and rewrite \eqref{Lemm:Ker:Bessel2:3} as
	\begin{equation}
		\begin{aligned}\label{Lemm:Ker:Bessel2:5}
			\mathfrak{B}^a_j(\xi_1,\eta_1,\xi_j,\eta_j)
			\ =  & t'\sqrt{\sin^2(\xi_1+\eta_1)\sin^2(2\xi_j)+[\sin(\xi_1+\eta_1)\cos(2\xi_j)-\sin(\eta_1)]^2}\sin(2\eta_j+G_j).
		\end{aligned}
	\end{equation}
	We estimate 
	\begin{equation}
		\begin{aligned}\label{Lemm:Ker:Bessel2:6}
			\mathfrak{B}_j': =	 \sqrt{\sin^2(\xi_1+\eta_1)\sin^2(2\xi_j)+[\sin(\xi_1+\eta_1)\cos(2\xi_j)-\sin(\eta_1)]^2} \ \ge \ |\sin(\xi_1+\eta_1)\sin(2\xi_j)|,
		\end{aligned}
	\end{equation}
	which means the right hand side of \eqref{Lemm:Ker:Bessel2:6} vanishes when $\sin(\xi_1+\eta_1)\sin(2\xi_j)=0$. Those singular points will be later integrated out. We now write $|\sin(\xi_1+\eta_1)|\mathfrak{A}_i^a(\xi_1,\eta_1,\xi_j)$ as
	\begin{equation}\label{Lemm:Ker:Bessel2:7}\begin{aligned}
			&	|\sin(\xi_1+\eta_1)|\mathfrak{A}_j^a(\xi_1,\eta_1,\xi_j) \ = \  \int_{-\pi}^\pi\mathrm{d}\eta_j 	|\sin(\xi_1+\eta_1)| e^{{\bf i}t'\mathfrak{B}_j'\sin(2\eta_j+G_j)}.\end{aligned}
	\end{equation}
	Let $\epsilon_{\eta_j}>0$ be a positive constant. We split $[-\pi,\pi]$ as the union of
	$$\mathbb{T}_{\eta_j}\ = \ \Big\{\eta_j\in[-\pi,\pi] \Big| |\cos(2\eta_j+G_j)|\ge  \epsilon_{\eta_j}\Big\}$$
	and
	$$\mathbb{T}'_{\eta_j}\ = \ [-\pi,\pi]\backslash \mathbb{T}_{\eta_j}.$$
	On $\mathbb{T}_{\eta_j}$, we compute
	\begin{equation}\label{Lemm:Ker:Bessel2:7:1}\begin{aligned}
			&	 \int_{\mathbb{T}_{\eta_j}}\mathrm{d}\eta_j 	|\sin(\xi_1+\eta_1)| e^{{\bf i}t'\mathfrak{B}_j'\sin(2\eta_j+G_j)}\\
			&	 \ = \  \int_{\mathbb{T}_{\eta_j}}\mathrm{d}\eta_j 	\frac{|\sin(\xi_1+\eta_1)|}{\partial_{\eta_j}[{\bf i}t'\mathfrak{B}_j'\sin(2\eta_j+G_j)]} \partial_{\eta_j}[{\bf i}t'\mathfrak{B}_j'\sin(2\eta_j+G_j)]e^{{\bf i}t'\mathfrak{B}_j'\sin(2\eta_j+G_j)} \\ 
			& \ = \  \int_{\mathbb{T}_{\eta_j}}\mathrm{d}\eta_j 	\frac{|\sin(\xi_1+\eta_1)|}{2{\bf i}t'\mathfrak{B}_j'\cos(2\eta_j+G_j)} \partial_{\eta_j}\Big[e^{{\bf i}t'\mathfrak{B}_j'\sin(2\eta_j+G_j)} \Big],\end{aligned}
	\end{equation}
	which, by a standard integration by parts argument, can be expressed as
	\begin{equation}\label{Lemm:Ker:Bessel2:8}\begin{aligned}
			&	\int_{\mathbb{T}_{\eta_j}}\mathrm{d}\eta_j 	|\sin(\xi_1+\eta_1)| e^{{\bf i}t'\mathfrak{B}_j'\sin(2\eta_j+G_j)}\\
			\ =  \ & 	\frac{|\sin(\xi_1+\eta_1)|}{2{\bf i}t'\mathfrak{B}_j'\cos(2\eta_j+G_j)} e^{{\bf i}t'\mathfrak{B}_j'\sin(2\eta_j+G_j)} \Big|_{\partial\mathbb{T}_{\eta_j}}   -  \int_{\mathbb{T}_{\eta_j}}\mathrm{d}\eta_j 	\partial_{\eta_j}\Big[\frac{|\sin(\xi_1+\eta_1)|}{2{\bf i}t'\mathfrak{B}_j'\cos(2\eta_j+G_j)} \Big]e^{{\bf i}t'\mathfrak{B}_j'\sin(2\eta_j+G_j)} \\
			\ =  \ & 	\frac{|\sin(\xi_1+\eta_1)|}{2{\bf i}t'\mathfrak{B}_j'\cos(2\eta_j+G_j)} e^{{\bf i}t'\mathfrak{B}_j'\sin(2\eta_j+G_j)} \Big|_{\partial\mathbb{T}_{\eta_j}}   -   \int_{\mathbb{T}_{\eta_j}}\mathrm{d}\eta_j 	\frac{|\sin(\xi_1+\eta_1)|\sin(2\eta_j+G_j)}{{\bf i}t'\mathfrak{B}_j'\cos^2(2\eta_j+G_j)} e^{{\bf i}t'\mathfrak{B}_j'\sin(2\eta_j+G_j)}.\end{aligned}
	\end{equation}
	We can simply bound  \eqref{Lemm:Ker:Bessel2:8} as follows
	\begin{equation}\label{Lemm:Ker:Bessel2:9}\begin{aligned}
			&	\Big|\int_{\mathbb{T}_{\eta_j}}\mathrm{d}\eta_j 	|\sin(\xi_1+\eta_1)| e^{{\bf i}t'\mathfrak{B}_j'\sin(2\eta_j+G_j)}\Big|\\
			\ \le   \ &  	     \int_{\mathbb{T}_{\eta_j}}\mathrm{d}\eta_j 	\frac{|\sin(\xi_1+\eta_1)|}{|t'\mathfrak{B}_j'\cos^2(2\eta_j+G_j)|} \ + \ \frac{|\sin(\xi_1+\eta_1)|}{2|t'\mathfrak{B}_j'\cos(2\eta_j+G_j)|}  \Big|_{\partial\mathbb{T}_{\eta_j}} \\
			\ \lesssim\ &  	    \frac{|\sin(\xi_1+\eta_1)|}{|t'\mathfrak{B}_j'\epsilon_{\eta_j}^2|} \ + \ \frac{|\sin(\xi_1+\eta_1)|}{2|t'\mathfrak{B}_j'\epsilon_{\eta_j}|} .\end{aligned}
	\end{equation}
	Combining \eqref{Lemm:Ker:Bessel2:6} and \eqref{Lemm:Ker:Bessel2:9}, we obtan
	\begin{equation}\label{Lemm:Ker:Bessel2:10}\begin{aligned}
			\Big|\int_{\mathbb{T}_{\eta_j}}\mathrm{d}\eta_j 	|\sin(\xi_1+\eta_1)| e^{{\bf i}t'\mathfrak{B}_j'\sin(2\eta_j+G_j)}\Big|
			\ \lesssim\ &  	   	\frac{1}{|t'\sin(2\xi_j)\epsilon_{\eta_j}^2|} \ + \ \frac{1}{|t'\sin(2\xi_j)\epsilon_{\eta_j}|}    	.\end{aligned}
	\end{equation}
	Moreover, a simple computation also gives
	\begin{equation}\label{Lemm:Ker:Bessel2:11}\begin{aligned}
			\Big|\int_{\mathbb{T}_{\eta_j}'}\mathrm{d}\eta_j 	|\sin(\xi_1+\eta_1)| e^{{\bf i}t'\mathfrak{B}_j'\sin(2\eta_j+G_j)}\Big|
			\ \lesssim\ &  	   	\epsilon_{\eta_j}   	.\end{aligned}
	\end{equation}
	Combining \eqref{Lemm:Ker:Bessel2:10}-\eqref{Lemm:Ker:Bessel2:11} and balancing $\epsilon_{\eta_j} $, we obtain
	\begin{equation}\label{Lemm:Ker:Bessel2:12}\begin{aligned}
			\Big||\sin(\xi_1+\eta_1)|\mathfrak{A}_j^a(\xi_1,\eta_1,\xi_j)\Big|
			\ \lesssim\ &  	   	\frac{1}{|t'\sin(2\xi_j)|^\frac13} \ + \ \frac{1}{|t'\sin(2\xi_j)|^{\frac23}}  \\
			\ \lesssim\ &  	   	\frac{1}{|t_0\sin(2\xi_j)|^\frac13} \ + \ \frac{1}{|t_0\sin(2\xi_j)|^{\frac23}}.\end{aligned}
	\end{equation}

	{\bf Step 2: Estimating $|\sin(\eta_1)|\mathfrak{A}_j^a(\xi_1,\eta_1,\xi_j)$}. We follow precisely all computations of Step 1, except the use of \eqref{Lemm:Ker:Bessel2:6}. We modify \eqref{Lemm:Ker:Bessel2:6} as follows

	\begin{equation}
		\begin{aligned}\label{Lemm:Ker:Bessel2:13}
			\mathfrak{B}_j'\ =	\ & \sqrt{\sin^2(\xi_1+\eta_1)\sin^2(2\xi_j)+[\sin(\xi_1+\eta_1)\cos(2\xi_j)-\sin(\eta_1)]^2} \\
			\ =	\ & \sqrt{\sin^2(\xi_1+\eta_1)\sin^2(2\xi_j)+\sin^2(\xi_1+\eta_1)\cos^2(2\xi_j)+\sin^2(\eta_1)-2\sin(\xi_1+\eta_1)\cos(2\xi_j)\sin(\eta_1)} \\
			\ =	\ & \sqrt{\sin^2(\xi_1+\eta_1)+\sin^2(\eta_1)-2\sin(\xi_1+\eta_1)\cos(2\xi_j)\sin(\eta_1)} 
			\\
			\ =	\ & \sqrt{[\sin(\xi_1+\eta_1)-\cos(2\xi_j)\sin(\eta_1)]^2+\sin^2(2\xi_j)\sin^2(\eta_1)} \\
			\ \ge \ & {|\sin(2\xi_j)\sin(\eta_1)|}.
		\end{aligned}
	\end{equation}
	By using \eqref{Lemm:Ker:Bessel2:13} instead of \eqref{Lemm:Ker:Bessel2:6}, we obtain a similar estimate as \eqref{Lemm:Ker:Bessel2:12}
	\begin{equation}\label{Lemm:Ker:Bessel2:14}\begin{aligned}
			\Big||\sin(\eta_1)|\mathfrak{A}_j^a(\xi_1,\eta_1,\xi_j)\Big|
			\ \lesssim\ &  	   	\frac{1}{|t'\sin(2\xi_j)|^\frac13} \ + \ \frac{1}{|t'\sin(2\xi_j)|^{\frac23}}  \\
			\ \lesssim\ &  	   	\frac{1}{|t_0\sin(2\xi_j)|^\frac13} \ + \ \frac{1}{|t_0\sin(2\xi_j)|^{\frac23}}    	.\end{aligned}
	\end{equation}
	
	{\bf Step 3: Estimating $\mathfrak{A}_j^a(\xi_1,\eta_1,\xi_j)$}. From \eqref{Lemm:Ker:Bessel2:10}, we deduce
	\begin{equation}\label{Lemm:Ker:Bessel2:15}\begin{aligned}
			\Big|\int_{\mathbb{T}_{\eta_j}}\mathrm{d}\eta_j 	 e^{{\bf i}t'\mathfrak{B}_j'\sin(2\eta_j+G_j)}\Big|
			\ \lesssim\ &  	   	\frac{1}{|t'||\sin(\xi_1+\eta_1)||\sin(2\xi_j)\epsilon_{\eta_j}^2|} \ + \ \frac{1}{|t'||\sin(\xi_1+\eta_1)||\sin(2\xi_j)\epsilon_{\eta_j}|}    	.\end{aligned}
	\end{equation}
	Moreover, a simple computation also gives
	\begin{equation}\label{Lemm:Ker:Bessel2:15}\begin{aligned}
			\Big|\int_{\mathbb{T}_{\eta_j}'}\mathrm{d}\eta_j 	 e^{{\bf i}t'\mathfrak{B}_j'\sin(2\eta_j+G_j)}\Big|
			\ \lesssim\ &  	   	\epsilon_{\eta_j}   	.\end{aligned}
	\end{equation}
	Combining \eqref{Lemm:Ker:Bessel2:14}-\eqref{Lemm:Ker:Bessel2:15} and balancing $\epsilon_{\eta_j} $, we obtain
	\begin{equation}\label{Lemm:Ker:Bessel2:16}\begin{aligned}
			\Big|\mathfrak{A}_j^a(\xi_1,\eta_1,\xi_j)\Big|
			\ \lesssim\ &  	   	\frac{1}{|t'\sin(\xi_1+\eta_1)\sin(2\xi_j)|^\frac13} \ + \ \frac{1}{|t'\sin(\xi_1+\eta_1)\sin(2\xi_j)|^{\frac23}}  \\
			\ \lesssim\ &  	   	\frac{1}{|t_0\sin(\xi_1+\eta_1)\sin(2\xi_j)|^\frac13} \ + \ \frac{1}{|t_0\sin(\xi_1+\eta_1)\sin(2\xi_j)|^{\frac23}}    	.\end{aligned}
	\end{equation}
	By using \eqref{Lemm:Ker:Bessel2:13} instead of \eqref{Lemm:Ker:Bessel2:6}, we obtain
	\begin{equation}\label{Lemm:Ker:Bessel2:17}\begin{aligned}
			\Big|\mathfrak{A}_j^a(\xi_1,\eta_1,\xi_j)\Big|
			\ \lesssim\ &  	   	\frac{1}{|t'\sin(\eta_1)\sin(2\xi_j)|^\frac13} \ + \ \frac{1}{|t'\sin(\eta_1)\sin(2\xi_j)|^{\frac23}}  \\
			\ \lesssim\ &  	   	\frac{1}{|t_0\sin(\eta_1)\sin(2\xi_j)|^\frac13} \ + \ \frac{1}{|t_0\sin(\eta_1)\sin(2\xi_j)|^{\frac23}}    	.\end{aligned}
	\end{equation}
	
	{\bf Step 4: The final estimates}. We first bound \eqref{Lemm:Ker:Bessel2:2a:a} as
	
	\begin{equation}
		\label{Lemm:Ker:Bessel2:18}\begin{aligned}
			\sum_{m\in\mathbb{Z}^d}|\mathfrak{F}^{Ker}(m,t_0) |^4
			\ 
			\lesssim \ & \int_{[-\pi,\pi]^d}\mathrm{d}\xi  \Big|\int_{-\pi}^{\pi}\mathrm{d}\eta_1 \prod_{j=2}^{n+1}\Big[|\sin(\xi_1+\eta_1)||\mathfrak{A}_j^a(\xi_1,\eta_1,\xi_j)|\Big]\\
			&\times \prod_{j=n+2}^{2n+1}\Big[|\sin(\eta_1)|\mathfrak{A}_j^a(\xi_1,\eta_1,\xi_j)|\Big]\\
			&\times |\mathfrak{A}_{2n+2}^a(\xi_1,\eta_1,\xi_{2n+2})||\mathfrak{A}_{2n+3}^a(\xi_1,\eta_1,\xi_{2n+3})|\\
			&\times |\mathfrak{A}_{2n+4}^a(\xi_1,\eta_1,\xi_{2n+4})||\mathfrak{A}_{2n+5}^a(\xi_1,\eta_1,\xi_{2n+5})|\Big|^2,\end{aligned}
	\end{equation}
	in which the other terms of the type $|\mathfrak{A}_j^a(\xi_1,\eta_1,\xi_j)|$ are simply bounded by constants. Using H\"older's inequality, we can bound \eqref{Lemm:Ker:Bessel2:18} as
	\begin{equation}
		\label{Lemm:Ker:Bessel2:20}\begin{aligned}
			\sum_{m\in\mathbb{Z}^d}|\mathfrak{F}^{Ker}(m,t_0) |^4
			\ 
			\lesssim \ & \int_{[-\pi,\pi]^d}\mathrm{d}\xi  \int_{-\pi}^{\pi}\mathrm{d}\eta_1 \prod_{j=2}^{n+1}\Big[|\sin(\xi_1+\eta_1)||\mathfrak{A}_j^a(\xi_1,\eta_1,\xi_j)|\Big]^2\\
			&\times \prod_{j=n+2}^{2n+1}\Big[|\sin(\eta_1)|\mathfrak{A}_j^a(\xi_1,\eta_1,\xi_j)|\Big]^2\\
			&\times |\mathfrak{A}_{2n+2}^a(\xi_1,\eta_1,\xi_{2n+2})|^2|\mathfrak{A}_{2n+3}^a(\xi_1,\eta_1,\xi_{2n+3})|^2\\
			&\times |\mathfrak{A}_{2n+4}^a(\xi_1,\eta_1,\xi_{2n+4})|^2|\mathfrak{A}_{2n+5}^a(\xi_1,\eta_1,\xi_{2n+5})|^2\\
			\ 
			\lesssim \ &   \int_{-\pi}^{\pi}\mathrm{d}\xi_1  \int_{-\pi}^{\pi}\mathrm{d}\eta_1 \prod_{j=2}^{n+1}\int_{-\pi}^{\pi}\mathrm{d}\xi_j \Big[|\sin(\xi_1+\eta_1)||\mathfrak{A}_j^a(\xi_1,\eta_1,\xi_j)|\Big]^2\\
			&\times \prod_{j=n+2}^{2n+1}\int_{-\pi}^{\pi}\mathrm{d}\xi_j\Big[|\sin(\eta_1)|\mathfrak{A}_j^a(\xi_1,\eta_1,\xi_j)|\Big]^2\\
			&\times \int_{-\pi}^{\pi}\mathrm{d}\xi_{2n+2}|\mathfrak{A}_{2n+2}^a(\xi_1,\eta_1,\xi_{2n+2})|^2\int_{-\pi}^{\pi}\mathrm{d}\xi_{2n+3}|\mathfrak{A}_{2n+3}^a(\xi_1,\eta_1,\xi_{2n+3})|^2\\
			&\times\int_{-\pi}^{\pi}\mathrm{d}\xi_{2n+4}	|\mathfrak{A}_{2n+4}^a(\xi_1,\eta_1,\xi_{2n+4})|^2\int_{-\pi}^{\pi}\mathrm{d}\xi_{2n+5}|\mathfrak{A}_{2n+5}^a(\xi_1,\eta_1,\xi_{2n+5})|^2.\end{aligned}
	\end{equation}
	
	Let $\epsilon_{\xi_j}>0$ be a positive constant and set
	$$\mathbb{T}_{\epsilon_{\xi_j}}=\Big\{\xi_j\in[-\pi,\pi] ~~ \Big| ~~ |\sin(2\xi_j)|>\epsilon_{\xi_j}\Big\}.$$

	We now bound the quantity that involves $\xi_j$, $j=2,\cdots,n+1$ in the product on the right hand side of \eqref{Lemm:Ker:Bessel2:20} 
	\begin{equation}
		\label{Lemm:Ker:Bessel2:21a}\begin{aligned}
			&	\int_{-\pi}^{\pi}\mathrm{d}\xi_j \Big[|\sin(\xi_1+\eta_1)||\mathfrak{A}_j^a(\xi_1,\eta_1,\xi_j)|\Big]^2\\
			\lesssim\	 & \int_{\mathbb{T}_{\epsilon_{\xi_j}}}\mathrm{d}\xi_j   \left[\frac{1}{|t_0\sin(2\xi_j)|^\frac13} \ + \ \frac{1}{|t_0\sin(2\xi_j)|^{\frac23}}  \right]^2 \ + \ \epsilon_{\xi_j}\\
			\lesssim\	  & \int_{\mathbb{T}_{\epsilon_{\xi_j}}}\mathrm{d}\xi_j  \left[\frac{1}{|t_0\sin(2\xi_j)|^\frac23} \ + \ \frac{1}{|t_0\sin(2\xi_j)|^{\frac43}}  \right] \ + \ \epsilon_{\xi_j}\end{aligned}
	\end{equation}
	Observing that $\frac{1}{|\sin(2\xi_j)|^\frac23}$ is always integrable on $[-\pi,\pi]$, we then bound
	\begin{equation}
		\label{Lemm:Ker:Bessel2:21b}\begin{aligned}
			&	\int_{-\pi}^{\pi}\mathrm{d}\xi_j \Big[|\sin(\xi_1+\eta_1)||\mathfrak{A}_j^a(\xi_1,\eta_1,\xi_j)|\Big]^2\\
			\lesssim\	  & \frac{1}{|t_0|^\frac23} \ + \ \int_{\mathbb{T}_{\epsilon_{\xi_j}}}\mathrm{d}\xi_j   \frac{1}{|t_0\sin(2\xi_j)|^{\frac43}}   \ + \ \epsilon_{\xi_j} \lesssim\	   \frac{1}{|t_0|^\frac23} \ + \ \int_{\mathbb{T}_{\epsilon_{\xi_j}}}\mathrm{d}\xi_j   \frac{1}{|t_0\epsilon_{\xi_j}|^{\frac43}}   \ + \ \epsilon_{\xi_j} \\
			\lesssim\	  & \frac{1}{|t_0|^\frac23} \ + \ \frac{1}{|t_0\epsilon_{\xi_j}|^{\frac43}}   \ + \ \epsilon_{\xi_j}.\end{aligned}
	\end{equation}
	Balancing $\epsilon_{\xi_j}$ by choosing $\epsilon_{\xi_j}=\mathcal O(t_0^{-\frac47})$ gives
	\begin{equation}
		\label{Lemm:Ker:Bessel2:21}\begin{aligned}
			&	\int_{-\pi}^{\pi}\mathrm{d}\xi_j \Big[|\sin(\xi_1+\eta_1)||\mathfrak{A}_j^a(\xi_1,\eta_1,\xi_j)|\Big]^2
			\lesssim\	  \frac{1}{\langle t_0\rangle ^\frac47}.\end{aligned}
	\end{equation}
	Similarly, the quantity that involves $\xi_j$, $j=n+2,\cdots,2n+1$ can also be bounded as
	\begin{equation}
		\label{Lemm:Ker:Bessel2:22}\begin{aligned}
			& \int_{-\pi}^{\pi}\mathrm{d}\xi_j \Big[|\sin(\eta_1)||\mathfrak{A}_j^a(\xi_1,\eta_1,\xi_j)|\Big]^2
			\lesssim	  \ \frac{1}{\langle t_0\rangle ^\frac47}.\end{aligned}
	\end{equation}
	Next, we bound the quantity that involves $\xi_{2n+2}$. A similar argument also gives
	\begin{equation}
		\label{Lemm:Ker:Bessel2:23a}\begin{aligned}
			&\int_{-\pi}^{\pi} \mathrm{d}\xi_{2n+2}  |\mathfrak{A}_j^a(\xi_1,\eta_1,\xi_{2n+2})|^2\\
			\lesssim\	&\int_{\mathbb{T}_{\epsilon_{\xi_{2n+2}}}}  \left[	\frac{1}{|t_0\sin(\xi_1+\eta_1)\sin(2\xi_{2n+2})|^\frac13} \ + \ \frac{1}{|t_0\sin(\xi_1+\eta_1)\sin(2\xi_{2n+2})|^{\frac23}}  \right]^2 \ + \ \epsilon_{\xi_{2n+2}}\\
			\lesssim\	&\int_{\mathbb{T}_{\epsilon_{\xi_{2n+2}}}}  \left[	\frac{1}{|t_0\sin(\xi_1+\eta_1)\sin(2\xi_{2n+2})|^\frac23} \ + \ \frac{1}{|t_0\sin(\xi_1+\eta_1)\sin(2\xi_{2n+2})|^{\frac43}}  \right] \ + \ \epsilon_{\xi_{2n+2}} \\
			\lesssim\	& \left[	\frac{1}{|t_0\sin(\xi_1+\eta_1)|^\frac23} \ + \ \frac{1}{|t_0\epsilon_{\xi_{2n+2}}\sin(\xi_1+\eta_1)|^{\frac43}}  \right] \ + \ \epsilon_{\xi_{2n+2}}.\end{aligned}
	\end{equation}
	Balancing $\epsilon_{\xi_{2n+2}}$ by choosing $\epsilon_{\xi_{2n+2}}=|t_0\sin(\xi_1+\eta_1)|^{-\frac47}$ gives 
	\begin{equation}
		\label{Lemm:Ker:Bessel2:23b}\begin{aligned}
			&\int_{-\pi}^{\pi} \mathrm{d}\xi_{2n+2} 	 |\mathfrak{A}_j^a(\xi_1,\eta_1,\xi_{2n+2})|^2\
			\lesssim\		\frac{1}{|t_0\sin(\xi_1+\eta_1)|^\frac23} \ + \ \frac{1}{|t_0\sin(\xi_1+\eta_1)|^\frac47}   .\end{aligned}
	\end{equation}
	Next, we estimate $$|\mathfrak{A}_j^a(\xi_1,\eta_1,\xi_{2n+3})|^2\le |\mathfrak{A}_j^a(\xi_1,\eta_1,\xi_{2n+3})|^\epsilon |\mathfrak{A}_j^a(\xi_1,\eta_1,\xi_{2n+3})|^{2-\epsilon}\le |\mathfrak{A}_j^a(\xi_1,\eta_1,\xi_{2n+3})|^\epsilon,$$ for some constant $\epsilon>0$ to be fixed later. Integrating this inequality gives 
	\begin{equation}
		\label{Lemm:Ker:Bessel2:23c}\begin{aligned}
			&\int_{-\pi}^{\pi} \mathrm{d}\xi_{2n+3} 	 |\mathfrak{A}_j^a(\xi_1,\eta_1,\xi_{2n+3})|^2\
			\lesssim\	\int_{-\pi}^{\pi} \mathrm{d}\xi_{2n+3} 	 |\mathfrak{A}_j^a(\xi_1,\eta_1,\xi_{2n+3})|^\epsilon\\
			\lesssim\	&\int_{-\pi}^{\pi} \mathrm{d}\xi_{2n+3} 	\left[\frac{1}{|t_0\sin(\xi_1+\eta_1)\sin(2\xi_{2n+3})|^\frac13} \ + \ \frac{1}{|t_0\sin(\xi_1+\eta_1)\sin(2\xi_{2n+3})|^{\frac23}}  \right]^\epsilon  \\
			\lesssim\	&\int_{-\pi}^{\pi} \mathrm{d}\xi_{2n+3} 	\left[\frac{1}{|t_0\sin(\xi_1+\eta_1)\sin(2\xi_{2n+3})|^{\frac{\epsilon }{3}}} \ + \ \frac{1}{|t_0\sin(\xi_1+\eta_1)\sin(2\xi_{2n+3})|^{\frac{2\epsilon }{3}}}  \right] 
			\\
			\lesssim\	&	\frac{1}{|t_0\sin(\xi_1+\eta_1)|^{\frac{\epsilon }{3}}} \ + \ \frac{1}{|t_0\sin(\xi_1+\eta_1)|^{\frac{2\epsilon }{3}}}\end{aligned}
	\end{equation}
	under the constraint that $0<\epsilon<\frac32$. 
	Combining \eqref{Lemm:Ker:Bessel2:23b}-\eqref{Lemm:Ker:Bessel2:23c} yields
	\begin{equation}
		\label{Lemm:Ker:Bessel2:23d}\begin{aligned}
			&\int_{-\pi}^{\pi} \mathrm{d}\xi_{2n+2} 	 |\mathfrak{A}_j^a(\xi_1,\eta_1,\xi_{2n+2})|^2\int_{-\pi}^{\pi} \mathrm{d}\xi_{2n+3} 	 |\mathfrak{A}_j^a(\xi_1,\eta_1,\xi_{2n+3})|^2\\
			\lesssim\	&	\frac{1}{|t_0\sin(\xi_1+\eta_1)|^\frac{2+\epsilon }{3}} \ + \ \frac{1}{|t_0\sin(\xi_1+\eta_1)|^{\frac47+\frac{2\epsilon }{3}}}   .\end{aligned}
	\end{equation}
	Choosing $\epsilon<\frac{9}{14}$, then $\frac{2+\epsilon }{3}<1$ and $\frac47+\frac{2\epsilon }{3}<\frac47+\frac{3}{7}=1$, 
	we estimate
	\begin{equation}
		\label{Lemm:Ker:Bessel2:23e}\begin{aligned}
			&\int_{-\pi}^{\pi} \mathrm{d}\xi_{1} \int_{-\pi}^{\pi} \mathrm{d}\xi_{2n+2} 	 |\mathfrak{A}_j^a(\xi_1,\eta_1,\xi_{2n+2})|^2\int_{-\pi}^{\pi} \mathrm{d}\xi_{2n+3} 	 |\mathfrak{A}_j^a(\xi_1,\eta_1,\xi_{2n+3})|^2\\
			\lesssim\	&\int_{-\pi}^{\pi} \mathrm{d}\xi_{1}\Big[\frac{1}{|t_0\sin(\xi_1+\eta_1)|^\frac{2+\epsilon }{3}} \ + \ \frac{1}{|t_0\sin(\xi_1+\eta_1)|^{\frac47+\frac{2\epsilon }{3}}}  \Big]\\
			\lesssim\	& \frac{1}{|t_0|^\frac{2+\epsilon }{3}} \ + \ \frac{1}{|t_0|^{\frac47+\frac{2\epsilon }{3}}}.\end{aligned}
	\end{equation}
	Choosing $\epsilon=\frac{3}{5}<\frac{9}{14}$, then $\frac{2+\epsilon }{3}=\frac{13}{15}$ and $\frac47+\frac{2\epsilon }{3}=\frac{34}{35},$ we find
	\begin{equation}
		\label{Lemm:Ker:Bessel2:23}\begin{aligned}
			&\int_{-\pi}^{\pi} \mathrm{d}\xi_{1} \int_{-\pi}^{\pi} \mathrm{d}\xi_{2n+2} 	 |\mathfrak{A}_j^a(\xi_1,\eta_1,\xi_{2n+2})|^2\int_{-\pi}^{\pi} \mathrm{d}\xi_{2n+3} 	 |\mathfrak{A}_j^a(\xi_1,\eta_1,\xi_{2n+3})|^2
			\lesssim\ \frac{1}{\langle |t_0|\rangle^\frac{13}{15}}.\end{aligned}
	\end{equation}
	And finally, by the same argument, we  estimate the quantity that involves $\xi_{2n+4}$ and $\xi_{2n+5}$
	\begin{equation}
		\label{Lemm:Ker:Bessel2:24}\begin{aligned}
			& \int_{-\pi}^{\pi} \mathrm{d}\eta_{1} \int_{-\pi}^{\pi} \mathrm{d}\xi_{2n+4} 	 |\mathfrak{A}_j^a(\xi_1,\eta_1,\xi_{2n+4})|^2\int_{-\pi}^{\pi} \mathrm{d}\xi_{2n+5} 	 |\mathfrak{A}_j^a(\xi_1,\eta_1,\xi_{2n+5})|^2
			\lesssim\ \frac{1}{\langle |t_0|\rangle^\frac{13}{15}}.\end{aligned}
	\end{equation}
	Combining \eqref{Lemm:Ker:Bessel2:20}-\eqref{Lemm:Ker:Bessel2:21}-\eqref{Lemm:Ker:Bessel2:22}-\eqref{Lemm:Ker:Bessel2:23}-\eqref{Lemm:Ker:Bessel2:24}, we obtain
	\begin{equation}
		\label{Lemm:Ker:Bessel2:25}\begin{aligned}
			\sum_{m\in\mathbb{Z}^d}|\mathfrak{F}^{Ker}(m,t_0) |^4
			\ \lesssim	  \ \frac{1}{\langle t_0\rangle ^{\frac{8n}{7}+\frac{26}{15}}},\end{aligned}
	\end{equation}
	meaning
	\begin{equation}
		\label{Lemm:Ker:Bessel2:26}\begin{aligned}
			\left(\sum_{m\in\mathbb{Z}^d}|\mathfrak{F}^{Ker}(m,t_0) |^4\right)^\frac14
			\ \lesssim	  \ \frac{1}{\langle t_0\rangle ^{\frac{2n}{7}+\frac{13}{30}}},\end{aligned}
	\end{equation}
	which yields the conclusion \eqref{Lemm:Ker:Bessel3:1} of the lemma. The second estimate \eqref{Lemm:Ker:Bessel3:2} is simply an easy consequence.
\end{proof}
\begin{lemma}\label{Lemm:Ker:Bessel4} There exists a universal constant $\mathfrak{C}_{\mathfrak{F}^{Ker},3}>0$ independent of $t_0$ and $\lambda$, such that
	\begin{equation}
		\label{Lemm:Ker:Bessel4:1}\begin{aligned}
			&	\|\mathfrak{F}^{Ker}(\cdot,t_0) \|_{l^3} \ \le \mathfrak{C}_{\mathfrak{F}^{Ker},3}\langle{t_0}\rangle^{{-{\frac{4n}{21}-\frac{13}{45}}}} \mbox { for } d>2n+5,\\
			&	\|\mathfrak{F}^{Ker}(\cdot,t_0) \|_{l^3} \ \le \mathfrak{C}_{\mathfrak{F}^{Ker},3}\langle{t_0}\rangle^{{-{\frac{4n}{21}}}} \mbox { for } d\le2n+5.	\end{aligned}
	\end{equation}
\end{lemma}

\begin{proof}
	The proof is the same with that of Lemma \ref{Lemm:Improved:Bessel4}. 
\end{proof}

\subsection{A convolution estimate involving the kernel}  We follow the same notations used in Section \ref{Subsec:DispersiveEstimates}.
The goal of this subsection is to estimate the following quantity
\begin{equation}
	\label{AnotherKernel:1}	\begin{aligned} 
		{\mathfrak{F}}^{Kern}(x,t) \ = \ & \int_{\mathbb{T}^d} \mathrm{d}k |\sin(2\pi k_0^1)|^2  |\sin(2\pi k^1)|^2|\sin(2\pi (k_0^1+k^1)|^2 e^{{\bf i}2\pi x\cdot k}e^{{\bf i}t\omega(k) +{\bf i}t\omega(-k_0-k)}, 	\end{aligned}
\end{equation}
Similar as in \eqref{Nearest1}-\eqref{extendedBesselfunctions1a:1}, we  define $\xi=(\xi_1,\cdots,\xi_d)=2\pi k\in [-\pi,\pi]^d$, $\xi^*=(\xi_1^*,\cdots,\xi_d^*)=2\pi k_0\in [-\pi,\pi]^d$, and rewrite \eqref{AnotherKernel:1} as
\begin{equation}
	\label{AnotherKernel:2}	\begin{aligned} 
		{\mathfrak{F}}^{Kern}(x,t) \ = \ & \int_{[-\pi,\pi]^d} \mathrm{d}\xi |\sin(\xi^*_1)|^2  |\sin(\xi_1)|^2|\sin(\xi^*_1+\xi_1)|^2 e^{{\bf i}x\cdot \xi}e^{{\bf i}t\omega(\xi) +{\bf i}t\omega(-\xi^*-\xi)}. 	\end{aligned}
\end{equation}
Next, we will prove an estimate that involves the $l^4$-norm of ${\mathfrak{F}}^{Kern}.$
\begin{lemma}\label{Lemm:AnotherKernel2} There exists a universal constant $\mathfrak{C}_{\tilde{\mathfrak{F}}^{Kern},2}>0$ independent of $t$, such that 
	\begin{equation}
		\label{Lemm:AnotherKernel2:1}
		\|{\mathfrak{F}}^{Kern}(\cdot,t) \|_{l^4}  \ \le \ \mathfrak{C}_{\tilde{\mathfrak{F}}^{Kern},2}  \frac{|\sin(\xi^*_1)|^2}{\langle |t||\sin(\xi_1^*)|\rangle^{\frac{1}{10}-}}.
	\end{equation}
\end{lemma}
\begin{proof}
	We first write
	\begin{equation}
		\label{Lemm:AnotherKernel2:E1}
		\begin{aligned}
			&	|{\mathfrak{F}}^{Kern}(x,t) |^2  
			\ =  \  \int_{[-\pi,\pi]^d}\mathrm{d}\xi e^{{\bf i}x\cdot\xi} \int_{-\pi}^{\pi}\mathrm{d}\eta_1  |\sin(\xi^*_1)|^4  |\sin(\xi_1+\eta_1)|^2|\sin(\xi^*_1+\xi_1+\eta_1)|^2\\
			&\ \  |\sin(\eta_1)|^2|\sin(\xi^*_1+\eta_1)|^2\exp\Big({\bf i}t\sin^3(\xi_1+\eta_1)  - {\bf i}t\sin^3(\xi_1+\eta_1+ \xi^*_1) -{\bf i}t\sin^3(\eta_1)+{\bf i}t\sin^3(\eta_1+ \xi^*_1)\Big)\\
			& \  \times\Big[\prod_{j=2}^d\int_{-\pi}^\pi\mathrm{d}\eta_j \exp\Big({\bf i}t_0\sin(\xi_1+\eta_1)\sin^2(\xi_j+\eta_j) - {\bf i}t_0\sin(\eta_1)\sin^2(\eta_j) \\
			&\ \ -{\bf i}t\sin(\xi_1+\eta_1+\xi^*_1)\sin^2(\xi_j+\eta_j+\xi^*_j) 
			\  + \ {\bf i}t\sin(\eta_1+\xi^*_1)\sin^2(\eta_j+ \xi^*_j)\Big)\Big]. 
		\end{aligned}
	\end{equation} 
	By the Plancherel theorem, we obtain
	\begin{equation}
		\label{Lemm:AnotherKernel2:E2}\begin{aligned}
			& \|\mathfrak{F}^{Kern}(\cdot,t) \|_{l^4}^4\ = \ \sum_{x\in\mathbb{Z}^d}|\mathfrak{F}^{Kern}(x,t) |^4\\  
			\ = \ & \int_{[-\pi,\pi]^d}\mathrm{d}\xi \Big| \int_{-\pi}^{\pi}\mathrm{d}\eta_1  |\sin(\xi^*_1)|^4  |\sin(\xi_1+\eta_1)|^2|\sin(\xi^*_1+\xi_1+\eta_1)|^2|\sin(\eta_1)|^2|\sin(\xi^*_1+\eta_1)|^2\\
			&\ \times \exp\Big({\bf i}t\sin^3(\xi_1+\eta_1)  - {\bf i}t\sin^3(\xi_1+\eta_1+ \xi^*_1) -{\bf i}t\sin^3(\eta_1)+{\bf i}t\sin^3(\eta_1+ \xi^*_1)\Big)\\
			& \  \times\Big[\prod_{j=2}^d\int_{-\pi}^\pi\mathrm{d}\eta_j \exp\Big({\bf i}t\sin(\xi_1+\eta_1)\sin^2(\xi_j+\eta_j) - {\bf i}t\sin(\eta_1)\sin^2(\eta_j) \\
			&\ \ -{\bf i}t\sin(\xi_1+\eta_1+\xi^*_1)\sin^2(\xi_j+\eta_j+\xi^*_j) 
			\  + \ {\bf i}t\sin(\eta_1+\xi^*_1)\sin^2(\eta_j+ \xi^*_j)\Big)\Big]\Big|^2.\end{aligned}
	\end{equation}
	We define for $2\le j\le d$
	\begin{equation}
		\begin{aligned}
			\mathfrak{A}_j(\xi_1,\eta_1,\xi_j) \ = & \ \int_{-\pi}^\pi\mathrm{d}\eta_j e^{{\bf i}\mathfrak{B}_j(\xi_1,\eta_1,\xi_j,\eta_j)},
		\end{aligned}
	\end{equation}
	where
	\begin{equation}
		\begin{aligned}
			\mathfrak{B}_j(\xi_1,\eta_1,\xi_j,\eta_j)\ =  	& \ t\sin(\xi_1+\eta_1)\sin^2(\xi_j+\eta_j)\
			-t\sin(\xi_1+\eta_1+\xi^*_1)\sin^2(\xi_j+\eta_j+\xi^*_j)\\
			&\ - t\sin(\eta_1)\sin^2(\eta_j) \ + \  t\sin(\eta_1+\xi^*_1)\sin^2(\eta_j+\xi^*_j),
		\end{aligned}
	\end{equation}
	and
	\begin{equation}
		\begin{aligned}
			\mathfrak{B}_1(\xi_1,\eta_1) = & \ t\sin^3(\xi_1+\eta_1) -t\sin^3(\xi_1+\eta_1+\xi^*_1) -t\sin^3(\eta_1)+t\sin^3(\eta_1+\xi^*_1),
		\end{aligned}
	\end{equation}
	yielding
	\begin{equation}
		\label{Lemm:AnotherKernel2:E3}\begin{aligned}
			\sum_{x\in\mathbb{Z}^d}|\mathfrak{F}^{Kern}(x,t) |^2
			\ = \ & \int_{[-\pi,\pi]^d}\mathrm{d}\xi \Big|\int_{-\pi}^{\pi}\mathrm{d}\eta_1|\sin(\xi^*_1)|^4  |\sin(\xi_1+\eta_1)|^2|\sin(\xi^*_1+\xi_1+\eta_1)|^2\\
			&\times |\sin(\eta_1)|^2|\sin(\xi^*_1+\eta_1)|^2\prod_{j=2}^d\mathfrak{A}_j(\xi_1,\eta_1,\xi_j) e^{{\bf i}\mathfrak{B}_1(\xi_1,\eta_1) }\Big|^2.\end{aligned}
	\end{equation}

	We  rewrite the phase $\mathfrak{B}_j$ as 
	the sum of
	\begin{equation}
		\begin{aligned}
			\mathfrak{B}^a_j(\xi_1,\eta_1,\xi_j,\eta_j)
			\ =   \ & -\frac{t}{2}\sin(\xi_1+\eta_1)\cos(2\xi_j+2\eta_j)
			\  + \ \frac{t}{2}\sin(\xi_1+\eta_1+\xi^*_1)\cos(2\xi_j+2\eta_j+2\xi^*_j)\\
			& \ + \ \frac{t}{2}\sin(\eta_1)\cos(2\eta_j)\
			- \ \frac{t}{2}\sin(\eta_1+\xi^*_1)\cos(2\eta_j+2\xi^*_j),
		\end{aligned}
	\end{equation}
	and
	\begin{equation}
		\begin{aligned}
			&\mathfrak{B}^b_j(\xi_1,\eta_1)
			=   \frac{t}{2}\sin(\xi_1+\eta_1)  - \frac{t}{2}\sin(\xi_1+\eta_1+\xi^*_1) 
			-  \frac{t}{2}\sin(\eta_1)\
			+\frac{t}{2}\sin(\eta_1+\xi^*_1).
		\end{aligned}
	\end{equation}
	The oscillatory integral $\mathfrak{A}_j$ can now be written $
	\mathfrak{A}_j \ = \ e^{{\bf i}\mathfrak{B}^b_j }\int_{-\pi}^\pi\mathrm{d}\eta_j e^{{\bf i}\mathfrak{B}^a_j}.
	$
	We rewrite  $\mathfrak{B}^a_j$ as follows
	\begin{equation}
		\begin{aligned}
			\mathfrak{B}^a_j(\xi_1,\eta_1,\xi_j,\eta_j)
			\ =   \ & -\ \mathrm{Re}\Big[e^{{\bf i}2\xi_j+{\bf i}2\eta_j}\Big(\frac{t}{2}\sin(\xi_1+\eta_1)- \frac{t}{2}\sin(\xi_1+\eta_1+\xi^*_1)e^{{\bf i}2\xi^*_j}\Big)\Big]\\
			& + \ \mathrm{Re}\Big[e^{{\bf i}2\eta_j}\Big(
			\frac{t}{2}\sin(\eta_1)\
			-\frac{t}{2}\sin(\eta_1+\xi^*_1)e^{{\bf i}2\xi^*_j}\Big)\Big].
		\end{aligned}
	\end{equation}
	Setting $$\mathfrak{C}_j^1=\frac{t}{2}\sin(\xi_1+\eta_1)- \frac{t}{2}\sin(\xi_1+\eta_1+\xi^*_1)e^{{\bf i}2\xi^*_j}, \ \ \ \mathfrak{C}_j^2=\frac{t}{2}\sin(\eta_1)\
	-\frac{t}{2}\sin(\eta_1+\xi^*_1)e^{{\bf i}2\xi^*_j},$$
	we write 
	\begin{equation}
		\begin{aligned}
			& \mathfrak{B}^a_j(\xi_1,\eta_1,\xi_j,\eta_j)
			\ =   \ -\ \mathrm{Re}\Big[e^{{\bf i}2\xi_j+{\bf i}2\eta_j}\mathfrak{C}_j^1\Big]\ + \ \mathrm{Re}\Big[e^{{\bf i}2\eta_j}\mathfrak{C}_j^2\Big], \ \ \  \mathfrak{A}_j^a \ = \ \int_{-\pi}^\pi\mathrm{d}\eta_j e^{{\bf i}\mathfrak{B}^a_j}.
		\end{aligned}
	\end{equation}
	We combine the phases ${\mathfrak{B}^b_j }$ and ${\mathfrak{B}_1}$  
	\begin{equation}
		\begin{aligned}
			& \mathfrak{B}^a_1=\mathfrak{B}_1 + \sum_{j=2}^{d}\mathfrak{B}^b_j \  =  \ t\Big[\sin^3(\xi_1+\eta_1)+\frac{d-1}{2}\sin(\xi_1+\eta_1)\Big]\\
			& -t\Big[\sin^3(\xi_1+\eta_1+\xi^*_1)+\frac{d-1}{2}\sin(\xi_1+\eta_1+\xi^*_1)\Big]\\
			& -t\Big[\sin^3(\eta_1)+\frac{d-1}{2}\sin(\eta_1)\Big]+t\Big[\sin^3(\eta_1+\xi^*_1)+\frac{d-1}{2}\sin(\eta_1+\xi^*_1)\Big],
		\end{aligned}
	\end{equation}
	and obtain 
	\begin{equation}
		\label{Lemm:AnotherKernel2:E4}\begin{aligned}
			\sum_{x\in\mathbb{Z}^d}|\mathfrak{F}^{Kern}(x,t) |^2
			\ = \ & \int_{[-\pi,\pi]^d}\mathrm{d}\xi \Big|\int_{-\pi}^{\pi}\mathrm{d}\eta_1|\sin(\xi^*_1)|^4  |\sin(\xi_1+\eta_1)|^2|\sin(\xi^*_1+\xi_1+\eta_1)|^2\\
			&\times |\sin(\eta_1)|^2|\sin(\xi^*_1+\eta_1)|^2 \prod_{j=2}^d\mathfrak{A}_j^a(\xi_1,\eta_1,\xi_j) e^{{\bf i}\mathfrak{B}^a_1(\xi_1,\eta_1) }\Big|^2\\
			\ 
			\le \ & \int_{[-\pi,\pi]^d}\mathrm{d}\xi \int_{-\pi}^{\pi}\mathrm{d}\eta_1|\sin(\xi^*_1)|^8  \prod_{j=2}^d\Big|\mathfrak{A}_j^a(\xi_1,\eta_1,\xi_j)\Big|^2	.\end{aligned}
	\end{equation}		
	
	We now set
	\begin{equation}
		\mathfrak{C}_j^1=\Xi_j^1 e^{{\bf i}2\Upsilon_j^1}, \ \ \ \ \mathfrak{C}_j^2=\Xi_j^2 e^{{\bf i}2\Upsilon_j^2},
	\end{equation}
	with $\Xi_j^1,\Xi_j^2\in\mathbb{R}_+$ and $\Upsilon_j^1,\Upsilon_j^2\in[-\pi,\pi]$, then
	\begin{equation}
		\begin{aligned}\label{Lemm:AnotherKernel2:E5}
			\Xi_j^1\ =\ & \Big|\Big(\frac{t}{2}-\frac{t}{2}\cos(\xi^*_1)e^{{\bf i}2\xi^*_j}\Big)\sin(\xi_1+\eta_1)\
			+\ \Big(-\frac{t}{2}\sin(\xi^*_1)e^{{\bf i}2\xi^*_j}\Big)\cos(\xi_1+\eta_1)\Big|.
		\end{aligned}
	\end{equation}
	Next, we will perform an a priori estimate on $\Xi_j^1$. Setting 
	\begin{equation}
		\begin{aligned}\label{Lemm:AnotherKernel2:E6}
			\cos(\Upsilon_*(j))
			\ = \ & {\Big|1-\cos(\xi^*_1)e^{{\bf i}2\xi^*_j}\Big|} \Big[\Big|1-\cos(\xi^*_1)e^{{\bf i}2\xi^*_j}\Big|^2\
			+\ \Big|\sin(\xi^*_1)e^{{\bf i}2\xi^*_j}\Big|^2\Big]^{-\frac12}\\
			\ = \ & \frac{\sqrt{1+\cos^2(\xi^*_1)-2\cos(\xi^*_1)\cos(2\xi^*_j)}}{\sqrt{2-2\cos(\xi^*_1)\cos(2\xi^*_j)}},\\
			\sin(\Upsilon_*(j))\
			= \ & \Big|\sin(\xi^*_1)e^{{\bf i}2\xi^*_j}\Big|\Big[\Big|1-\cos(\xi^*_1)e^{{\bf i}2\xi^*_j}\Big|^2\
			+\ \Big|\sin(\xi^*_1)e^{{\bf i}2\xi^*_j}\Big|^2\Big]^{-\frac12}
			\ = \  \frac{|\sin(\xi^*_1)|}{\sqrt{2-2\cos(\xi^*_1)\cos(2\xi^*_j)}},
		\end{aligned}
	\end{equation}
	with $\Upsilon_*\in[0,\pi/2]$, we then find
	\begin{equation}
		\begin{aligned}\label{Lemm:AnotherKernel2:E7}
			\Xi_j^1
			\ =\ 
			&\frac{|t|}{2} \sqrt{2-2\cos(\xi^*_1)\cos(2\xi^*_j)} \Big|\cos(\Upsilon_*)\sin(\xi_1+\eta_1)e^{{\bf i}\aleph
				_1}+\sin(\Upsilon_*)\cos(\xi_1+\eta_1)e^{{\bf i}\aleph
				_2}\Big|,
		\end{aligned}
	\end{equation}
	for some angles  $\aleph
	_1, \aleph
	_2$.
	As a consequence, similar with \eqref{Lemm:Bessel2:2bb9c:1}-\eqref{Lemm:Bessel2:2bb9c:2}, we bound
	\begin{equation}
		\begin{aligned}\label{Lemm:AnotherKernel2:E8}
			\Xi_j^1
			\ \gtrsim \ &
			|t|	\sqrt{2-2\cos(\xi^*_1)\cos(2\xi^*_j)}\Big[[\cos^2(\Upsilon_*)\sin^2(\xi_1+\eta_1)+\sin^2(\Upsilon_*)\cos^2(\xi_1+\eta_1)]^\frac12\\ 
			&\times [1-|\cos(\aleph
			_1^j-\aleph
			_2^j)|]^\frac12 + |\sin(\xi_1+\eta_1-\Upsilon_*)||\cos(\aleph^j
			_1-\aleph
			_2^j)|^\frac12 \Big],
		\end{aligned}
	\end{equation}	
	and, by the same argument, we also bound
	\begin{equation}
		\begin{aligned}\label{Lemm:AnotherKernel2:E9}
			\Xi_j^2
			\ \gtrsim \ &
			|t|\sqrt{2-2\cos(\xi^*_1)\cos(2\xi^*_j)}\Big[[\cos^2(\Upsilon_*)\sin^2(\eta_1)+\sin^2(\Upsilon_*)\cos^2(\eta_1)]^\frac12\\ 
			&\times[1-|\cos(\aleph
			_1^j-\aleph
			_2^j)|]^\frac12  + |\sin(\eta_1-\Upsilon_*)||\cos(\aleph^j
			_1-\aleph
			_2^j)|^\frac12 \Big].
		\end{aligned}
	\end{equation}
	Notice that, by the definition \eqref{Lemm:AnotherKernel2:E7}, we find $\aleph
	_2^j=2\xi_j^*$ and  $\aleph
	_1^j=\arctan\Big(\frac{\cos(\xi^*_1)\sin(2\xi^*_j)}{\cos(\xi^*_1)\cos(2\xi^*_j)-1}\Big)$. Observing that
	\begin{equation}\label{Lemm:AnotherKernel2:E10}\begin{aligned}
			|\sin(\aleph_1^j-\aleph_2^j)|& \le |\tan(\aleph_1^j)-\tan(\aleph_2^j)|= \Big|\frac{\cos(\xi^*_1)\sin(2\xi^*_j)}{\cos(\xi^*_1)\cos(2\xi^*_j)-1}-\tan(2\xi^*_j)\Big|\\
			& =\Big|\frac{\sin(2\xi^*_j)}{\cos(2\xi^*_j)(1-\cos(2\xi^*_j)\cos(\xi_1^*))}\Big|\le \Big|\frac{\sin(2\xi^*_j)}{\cos(2\xi^*_j)(1-\cos(2\xi^*_j))}\Big|,
		\end{aligned}
	\end{equation}
	we then deduce
	\begin{equation}\label{Lemm:AnotherKernel2:E11}\begin{aligned}
			|\cos(\aleph_1^j-\aleph_2^j)|& \ge \Big[1- \Big|\frac{\sin(2\xi^*_j)}{\cos(2\xi^*_j)(1-\cos(2\xi^*_j))}\Big|^2\Big]^\frac12.
		\end{aligned}
	\end{equation}	
	Now, by \eqref{Lemm:Angle:9}, we have
	$
	[1-|\cos(\aleph_1^j-\aleph_2^j)|]\ 
	\gtrsim   \ |\sin(2\xi^*_j)|^2.
	$
	We thus bound
	\begin{equation}
		\begin{aligned}\label{Lemm:AnotherKernel2:E13}
			\Xi_j^1
			\ \gtrsim \ &
			|t||\sin(\xi_1^*)||\sin(\xi_1+\eta_1-\Upsilon_*)||\sin(2\xi^*_j)|\\
			&\  + |t|\sqrt{2-2\cos(\xi^*_1)\cos(2\xi^*_j)}|\sin(\xi_1+\eta_1-\Upsilon_*)|\Big[1- \Big|\frac{\sin(2\xi^*_j)}{\cos(2\xi^*_j)(1-\cos(2\xi^*_j))}\Big|^2\Big]^\frac14\\
			\ \gtrsim \ &
			|t||\sin(\xi_1^*)||\sin(\xi_1+\eta_1-\Upsilon_*)||\sin(2\xi^*_j)|\\
			&\  + |t|\sqrt{2-2\cos(\xi^*_1)}|\sin(\xi_1+\eta_1-\Upsilon_*)|\Big[1- \Big|\frac{\sin(2\xi^*_j)}{\cos(2\xi^*_j)(1-\cos(2\xi^*_j))}\Big|^2\Big]^\frac14\\
			\ \gtrsim \ &
			\Xi_j^{1,a} |\sin(2\xi^*_j)|\ +\ \Xi_j^{1,a} \Big[1- \Big|\frac{\sin(2\xi^*_j)}{\cos(2\xi^*_j)(1-\cos(2\xi^*_j))}\Big|^2\Big]^\frac14,
		\end{aligned}
	\end{equation}	
	with $\Xi_j^{1,a} =|t||\sin(\xi_1^*)||\sin(\xi_1+\eta_1-\Upsilon_*)|$   $\le2|t||\sin(\xi_1^*/2)|$ $|\sin(\xi_1+\eta_1-\Upsilon_*)| $. Similarly, we can also bound
	\begin{equation}
		\begin{aligned}\label{Lemm:AnotherKernel2:E14}
			\Xi_j^2
			\ \gtrsim \ &
			\Xi_j^{2,a} |\sin(2\xi^*_j)|\ +\ \Xi_j^{2,a} \Big[1- \Big|\frac{\sin(2\xi^*_j)}{\cos(2\xi^*_j)(1-\cos(2\xi^*_j))}\Big|^2\Big]^\frac14,
		\end{aligned}
	\end{equation}
	with $\Xi_j^{2,a} =|t||\sin(\xi_1^*)||\sin(\eta_1-\Upsilon_*)|$. The two inequalities \eqref{Lemm:AnotherKernel2:E13} and \eqref{Lemm:AnotherKernel2:E14} imply that, when $ |\sin(2\xi^*_j)|$ is small, then $\Big[1- \Big|\frac{\sin(2\xi^*_j)}{\cos(2\xi^*_j)(1-\cos(2\xi^*_j))}\Big|^2\Big]^\frac14\gtrsim1$, thus, we use $\Xi_j^{1,a}$ and $\Xi_j^{2,a}$ as the lower bound for $\Xi_j^{1}$ and $\Xi_j^{2}$. When $ |\sin(2\xi^*_j)|\gtrsim1$, we can still use $\Xi_j^{1,a}$ and $\Xi_j^{2,a}$ as the lower bound for $\Xi_j^{1}$ and $\Xi_j^{2}$. 
	
	We next set $$\mathbb{T}_{\xi_j} : = \Big\{\xi_j\ \ \Big| \ \ |\sin(2\xi_j+2\Upsilon_j^1-2\Upsilon_j^2)| >\epsilon_{\xi_j}\Big\},\ \ \ j=2,\cdots,d.$$
	Under the constraint that
	$
	|t||\sin(\xi_1^*)|\epsilon_{\xi_j}^2>>1,
	$
	by the same strategy used in Step 1 of the Proof of Lemma \ref{Lemm:Bessel3}, we have the following two bounds for $\xi_j\in \mathbb{T}_{\xi_j}$
	\begin{equation}\begin{aligned}
			\label{Lemm:AnotherKernel2:E16}
			|\mathfrak{A}_{j}^a| 
			\ \lesssim \  & \frac{1}{	|t||\sin(\xi_1^*)|\min\{|\sin(\xi_1+\eta_1-\Upsilon_*(j))|,|\sin(\eta_1-\Upsilon_*(j))|\}|\sin(2\xi_j+2\Upsilon_j^1-2\Upsilon_j^2)|^2}\\
			& \ + \  \frac{1}{\sqrt{|t||\sin(\xi_1^*)|\min\{|\sin(\xi_1+\eta_1-\Upsilon_*(j))|,|\sin(\eta_1-\Upsilon_*(j))|\}|\sin(2\xi_j+2\Upsilon_j^1-2\Upsilon_j^2)|}}.
		\end{aligned}
	\end{equation}
	For $\xi_j\notin \mathbb{T}_{\xi_j}$, we simply have 
	\begin{equation}
		\int_{[-\pi,\pi]\backslash\mathbb{T}_{\xi_j}}\mathrm{d}\xi_j |\mathfrak{A}_{j}^a|^2\le \epsilon_{\xi_j}.  
	\end{equation}
	Thus, by balancing $\epsilon_{\xi_j}$, we find
	\begin{equation}\begin{aligned}
			\label{Lemm:AnotherKernel2:E17}
			\int_{[-\pi,\pi]}\mathrm{d}\xi_j|\mathfrak{A}_{j}^a|^2 
			\ \lesssim \  & \frac{1}{	|t|^\frac25|\sin(\xi_1^*)|^\frac25\min\{|\sin(\xi_1+\eta_1-\Upsilon_*(j))|,|\sin(\eta_1-\Upsilon_*(j))|\}^\frac25}.
		\end{aligned}
	\end{equation}
	By the trivial bound $|\mathfrak{A}_{j}^a|\lesssim 1$, we also have, for $0<\epsilon\le 1$
	\begin{equation}\begin{aligned}
			\label{Lemm:AnotherKernel2:E18}
			\int_{[-\pi,\pi]}\mathrm{d}\xi_j|\mathfrak{A}_{j}^a|^2  
			\ \lesssim \ & \Big[\int_{[-\pi,\pi]}\mathrm{d}\xi_j|\mathfrak{A}_{j}^a|^2 \Big]^\epsilon\\
			\ \lesssim \ & \frac{1}{[	|t|^\frac25|\sin(\xi_1^*)|^\frac25\min\{|\sin(\xi_1+\eta_1-\Upsilon_*(j))|,|\sin(\eta_1-\Upsilon_*(j))|\}^\frac25]^\epsilon},\\
			\int_{[-\pi,\pi]}\mathrm{d}\xi_j|\mathfrak{A}_{j}^a|^2  \ \lesssim \ &	\Big[\int_{[-\pi,\pi]}\mathrm{d}\xi_j|\mathfrak{A}_{j}^a|^2\Big]^\epsilon  \\ \lesssim \  & \frac{1}{	[|t|^\frac25|\sin(\xi_1^*)|^\frac25\min\{|\sin(\xi_1+\eta_1-\Upsilon_*(j))|,|\sin(\eta_1-\Upsilon_*(j))|\}^\frac25]^\epsilon}.
		\end{aligned}
	\end{equation}
	From \eqref{Lemm:AnotherKernel2:E4}, we deduce
	\begin{equation}
		\label{Lemm:AnotherKernel2:E19}\begin{aligned}
			\sum_{x\in\mathbb{Z}^d}|\mathfrak{F}^{Kern}(x,t) |^4
			\le \ & |\sin(\xi^*_1)|^8 \int_{-\pi}^{\pi}\mathrm{d}\eta_1 \int_{-\pi}^{\pi}\mathrm{d}\xi_1  \prod_{j=2}^d\int_{[-\pi,\pi]}\mathrm{d}\xi_j \Big|\mathfrak{A}_j^a(\xi_1,\eta_1,\xi_j)\Big|^2.\end{aligned}
	\end{equation}	
	The same argument used in the final step of the Proof  of Lemma \ref{Lemm:Bessel3} can be reiterated,  using \eqref{Lemm:AnotherKernel2:E17} and \eqref{Lemm:AnotherKernel2:E18}, yielding
	\begin{equation}
		\label{Lemm:AnotherKernel2:E20}\begin{aligned}
			\sum_{x\in\mathbb{Z}^d}|\mathfrak{F}^{Kern}(x,t) |^4
			\le \ & |\sin(\xi^*_1)|^8 \frac{1}{\langle |t||\sin(\xi_1^*)|\rangle^{2/5-}}.\end{aligned}
	\end{equation}	
\end{proof}

\subsection{Estimates of the collision operators}
We follow the same notations used in Section \ref{Subsec:DispersiveEstimates}.


\begin{remark}
	\label{Remark:SingularDispersionRelation} As $\omega$ vanishes on a the ghost manifold,   the oscillatory integrals created by the dispersion $\omega$ do not have as much decay in time, in comparison with those created by the dispersion relation of the  Schr\"odinger equation. However, as we could see from Lemma \ref{Lemm:Ker:Bessel3},  thanks to some refined  estimates concerning the kernels, we have sufficient decays in time to  define  the   resonance appearing in the collision operators. This is the main idea of Lemma \ref{Lemma:Resonance1}.  
\end{remark}
\begin{remark}
	In the  definition of $Q_1^{colli}$ below, $F_1 $ is a function of $k_1$, while in the definition of $Q_2^{colli}$, $F_1 $ is a function of $k_0$, as they will be later applied to different types of ladder operators, discussed after inequality \eqref{FinalProof:E10a:1:5}. 
\end{remark}
\begin{lemma}
	\label{Lemma:Resonance1}
	Let $1>\upsilon	>0$ be a positive constant, suppose that $d\ge 2$, we define the following operators, for all $k_0,k_1,k_2\in \mathbb{T}^d$, $\sigma_0,\sigma_1,\sigma_2\in\{\pm1\}$,  and for any constant $T_o>0$
	\begin{equation}\label{Lemma:Resonance1:a}
		\begin{aligned} 
			& 
			F_1,F_2\in L^4(\mathbb{T}^d) \longrightarrow \mathcal{Q}^{colli}_1[F_1,F_2,\sigma_0,\sigma_1,\sigma_2,s](k_0):=\\\
			& e^{{\bf i}s\sigma_0 \omega(k_0)}\iint_{(\mathbb{T}^d)^2}\!\! \mathrm{d} k_1  \mathrm{d} k_2 \,
			\delta(\sigma_0k_0+\sigma_1k_1+\sigma_2k_2) 
			e^{{\bf i}s\sigma_1 \omega(k_1)}\\
			&\times  {F}_1(k_1) F_2(k_2)e^{{\bf i} s \sigma_2 \omega(k_2)} |\sin(2\pi k_1^1)||\sin(2\pi k_2^1)|,
		\end{aligned}
	\end{equation}
	\begin{equation}\label{Lemma:Resonance1:b}
		\begin{aligned} 
			& 
			F_1\in L^\infty(\mathbb{T}^d),F_2\in L^4(\mathbb{T}^d) \longrightarrow \mathcal{Q}^{colli}_2[F_1,F_2,\sigma_0,\sigma_1,\sigma_2,s](k_0):=\\\
			& e^{{\bf i}s\sigma_0 \omega(k_0)}\iint_{(\mathbb{T}^d)^2}\!\! \mathrm{d} k_1  \mathrm{d} k_2 \,
			\delta(\sigma_0k_0+\sigma_1k_1+\sigma_2k_2) 
			e^{{\bf i}s\sigma_1 \omega(k_1)}|\sin(2\pi k_1^1)|\\
			&\times |\sin(2\pi k_0^1)|  {F}_1(k_0) F_2(k_2)e^{{\bf i} s \sigma_2 \omega(k_2)} |\sin(2\pi k_2^1)|,
		\end{aligned}
	\end{equation}
	\begin{equation}\label{Lemma:Resonance1:b}
		\begin{aligned} 
			& 
			F_1\in L^4(\mathbb{T}^d) \longrightarrow \mathcal{Q}^{colli}_3[F_1,\sigma_0,\sigma_1,\sigma_2,s](k_0):=\\\
			& e^{{\bf i}s\sigma_0 \omega(k_0)}\iint_{(\mathbb{T}^d)^2}\!\! \mathrm{d} k_1  \mathrm{d} k_2 \,
			\delta(\sigma_0k_0+\sigma_1k_1+\sigma_2k_2) 
			e^{{\bf i}s\sigma_1 \omega(k_1)}|\sin(2\pi k_0^1)||\sin(2\pi k_1^1)| \\
			&\times  F_1(k_1)e^{{\bf i} s \sigma_2 \omega(k_2)} |\sin(2\pi k_2^1)|.
		\end{aligned}
	\end{equation}
	Moreover, we also define
	
	\begin{equation}
		\label{Lemma:Resonance1:1}\begin{aligned}
			&	F_1,F_2\in L^4(\mathbb{T}^d), \phi\in C^\infty_c(\mathbb{T}^d\times(-\infty,\infty))\longrightarrow  \ \mathbb{O}^1_{\upsilon, \lambda}[F_1,F_2][\phi] :=\\
			&  \int_{-T_o\lambda^{-2}}^{T_o\lambda^{-2}} \frac{\mathrm{d}s}{\pi}e^{-\upsilon |s|}\iiint_{(\mathbb{T}^d)^3}\!\! \mathrm{d} k_0  \mathrm{d} k_1  \mathrm{d} k_2 \Big[\delta(\sigma_0k_0+\sigma_1k_1+\sigma_2k_2)
			e^{{\bf i}s[\sigma_0 \omega(k_0)+\sigma_1 \omega(k_1)+\sigma_2 \omega(k_2)]}\\
			&\times \ \ |\sin(2\pi k_1^1)||\sin(2\pi k_2^1)|	{F}_1(k_1) F_2(k_2)\phi(k_0,s)\Big],
		\end{aligned}	
	\end{equation}
	and
	\begin{equation}
		\label{Lemma:Resonance1:1}\begin{aligned}
			&	F_2\in L^4(\mathbb{T}^d), F_1\in  L^\infty(\mathbb{T}^d), \phi\in C^\infty_c(\mathbb{T}^d\times(-\infty,\infty))\longrightarrow  \ \mathbb{O}^2_{\upsilon, \lambda}[F_1,F_2][\phi] :=\\
			& \int_{-T_o\lambda^{-2}}^{T_o\lambda^{-2}} \frac{\mathrm{d}s}{\pi}e^{-\upsilon |s|} \iiint_{(\mathbb{T}^d)^3}\!\! \mathrm{d} k_0  \mathrm{d} k_1  \mathrm{d} k_2 \, \Big[\delta(\sigma_0k_0+\sigma_1k_1+\sigma_2k_2)
			e^{{\bf i}s[\sigma_0 \omega(k_0)+\sigma_1 \omega(k_1)+\sigma_2 \omega(k_2)]}\\
			&\times \ \ |\sin(2\pi k_0^1)||\sin(2\pi k_1^1)||\sin(2\pi k_2^1)|	{F}_1(k_0) F_2(k_2)\phi(k_0,s)\Big],
		\end{aligned}	
	\end{equation}
	where $k_0=(k_0^1,\cdots,k_0^d),k_1=(k_1^1,\cdots,k_1^d),k_2=(k_2^1,\cdots,k_2^d)\in\mathbb{T}^d.$

	The following claims then hold true.
	\begin{itemize}
		\item[(a)]  
		There exists a constant $\mathscr{M}>1$ such that we  have
		\begin{equation}
			\label{Lemma:Resonance1:2:bis}\begin{aligned}
				& \lim_{\upsilon\to 0} \int_{-T_o\lambda^{-2}}^{T_o\lambda^{-2}} \frac{\mathrm{d}s}{\pi}e^{-\upsilon |s|}\Big\|\mathcal {Q}^{colli}_1[F_1,F_2,\sigma_0,\sigma_1,\sigma_2,s]\Big\|_{L^4}^\mathscr{M}\\
				\ = \ & \int_{-T_o\lambda^{-2}}^{T_o\lambda^{-2}} \frac{\mathrm{d}s}{\pi}\Big\|\mathcal {Q}^{colli}_1[F_1,F_2,\sigma_0,\sigma_1,\sigma_2,s]\Big\|_{L^4}^\mathscr{M},\end{aligned}	
		\end{equation}	
		\begin{equation}
			\label{Lemma:Resonance1:2:bis1}\begin{aligned}
				&\ \int_{-T_o\lambda^{-2}}^{T_o\lambda^{-2}} {\mathrm{d}s}\Big\|\mathcal {Q}^{colli}_1[F_1,F_2,\sigma_0,\sigma_1,\sigma_2,s]\Big\|_{L^4}^\mathscr{M} \lesssim \||\sin(2\pi k^1)|F_1\|_{L^4}^\mathscr{M}  \||\sin(2\pi k^1)|F_2\|_{L^4}^\mathscr{M}. 
			\end{aligned}
		\end{equation}
		and 	\begin{equation}
			\label{Lemma:Resonance1:2:bis1:a}\begin{aligned}
				\lim_{\upsilon\to 0}\mathbb{O}^1_{\upsilon, \lambda}[F_1,F_2][\phi]
				\ =  &\ \mathbb{O}^1_{0, \lambda}[F_1,F_2][\phi].
			\end{aligned}
		\end{equation}	
			\item[(b)] 
			Moreover, we also have
			
			\begin{equation}
				\label{Lemma:Resonance1:2}\begin{aligned}
					& \lim_{\upsilon\to 0} \int_{-T_o\lambda^{-2}}^{T_o\lambda^{-2}} \frac{\mathrm{d}s}{\pi}e^{-\upsilon |s|}\Big\|\mathcal {Q}^{colli}_2[F_1,F_2,\sigma_0,\sigma_1,\sigma_2,s]\Big\|_{L^4}^\mathscr{M}\\
					\ = \ & \int_{-T_o\lambda^{-2}}^{T_o\lambda^{-2}} \frac{\mathrm{d}s}{\pi}\Big\|\mathcal {Q}^{colli}_2[F_1,F_2,\sigma_0,\sigma_1,\sigma_2,s]\Big\|_{L^4}^\mathscr{M},\end{aligned}	
			\end{equation}	
			\begin{equation}
				\label{Lemma:Resonance1:2:bis2}\begin{aligned}
					&\ \int_{-T_o\lambda^{-2}}^{T_o\lambda^{-2}} {\mathrm{d}s}\Big\|\mathcal {Q}^{colli}_2[F_1,F_2,\sigma_0,\sigma_1,\sigma_2,s]\Big\|_{L^4}^\mathscr{M} \lesssim \| |\sin(2\pi k^1)|F_1(k)\|_{L^\infty}^\mathscr{M}\||\sin(2\pi k^1)| F_2\|_{L^4}^\mathscr{M}, 
				\end{aligned}
			\end{equation}
			with $k=(k^1,\cdots,k^d)\in\mathbb{T}^d,$
			and 	\begin{equation}
				\label{Lemma:Resonance1:2:bis2:a}\begin{aligned}
					\lim_{\upsilon\to 0}\mathbb{O}^2_{\upsilon, \lambda}[F_1,F_2][\phi]
					\ =  &\ \mathbb{O}^2_{0, \lambda}[F_1,F_2][\phi].
				\end{aligned}
			\end{equation}	
			\item[(c)] For $\phi=e^{-\ell|s|}$ and for any constant $c>0$, we set 
			\begin{equation}
				\label{Lemma:Resonance1:5}\begin{aligned}
					& \mathbb{O}_{\upsilon, \lambda}^{1,c}[F_1,F_2][\phi]\	=  \  \iiint_{\big\{k_0,k_1,k_2\in\mathbb{T}^d\big||\sigma_0 \omega(k_0)+\sigma_1 \omega(k_1)+\sigma_2 \omega(k_2)|<c\big\}}\mathrm{d} k_0  \mathrm{d} k_1  \mathrm{d} k_2 \,
					\\
					& \ \ \frac{	1}{\pi[(\upsilon+\ell)^2+(\sigma_0 \omega(k_0)+\sigma_1 \omega(k_1)+\sigma_2 \omega(k_2))^2]} \Big\{(\upsilon+\ell)-e^{-T_o\lambda^{-2}(\upsilon+\ell)}\Big[\cos(\lambda^{-2}T_o(\sigma_0 \omega(k_0)\\
					& \ \ +\sigma_1 \omega(k_1)+\sigma_2 \omega(k_2))(\upsilon+\ell)	- (\sigma_0 \omega(k_0)+\sigma_1 \omega(k_1)+\sigma_2 \omega(k_2))\\
					& \ \times\sin(\lambda^{-2}T_o(\sigma_0 \omega(k_0)+\sigma_1 \omega(k_1)+\sigma_2 \omega(k_2))\Big]\Big\}\Big[\delta(\sigma_0k_0+\sigma_1k_1\\
					& \ \   +\sigma_2k_2)
					|\sin(2\pi k_1^1)||\sin(2\pi k_2^1)|	{F}_1(k_1) F_2(k_2)\Big].\end{aligned}	
			\end{equation}
			Then for all $c>0$ 
			\begin{equation}
				\label{Lemma:Resonance1:6}\begin{aligned}& \lim_{\upsilon+\ell\to 0}\lim_{\lambda\to0}|\mathbb{O}_{\upsilon, \lambda}^{1}[F_1,F_2][\phi]-\mathbb{O}_{\upsilon, \lambda}^{1,c}[F_1,F_2][\phi]  | \ = \ 0.\end{aligned}
			\end{equation}
			\item[(d)] For $\phi=e^{-\ell|s|}$ and for any constant $c>0$,  we set 
			\begin{equation}
				\label{Lemma:Resonance1:5:1}\begin{aligned}
					&	\mathbb{O}_{\upsilon, \lambda}^{2,c}[F_1,F_2][\phi]\	=  \  \iiint_{\big\{k_0,k_1,k_2\in\mathbb{T}^d\big||\sigma_0 \omega(k_0)+\sigma_1 \omega(k_1)+\sigma_2 \omega(k_2)|<c\big\}}\mathrm{d} k_0  \mathrm{d} k_1  \mathrm{d} k_2 \,
					\\
					&\frac{	1}{\pi[(\upsilon+\ell)^2+(\sigma_0 \omega(k_0)+\sigma_1 \omega(k_1)+\sigma_2 \omega_\infty(k_2))^2]}  \Big\{(\upsilon+\ell)-e^{-T_o\lambda^{-2}(\upsilon+\ell)}\Big[\cos(\lambda^{-2}T_o(\sigma_0 \omega(k_0)\\
					& \ \ +\sigma_1 \omega(k_1)+\sigma_2 \omega(k_2))(\upsilon+\ell)	- (\sigma_0 \omega(k_0)+\sigma_1 \omega(k_1)+\sigma_2 \omega(k_2))\\
					& \ \ \times\sin(\lambda^{-2}T_o(\sigma_0 \omega(k_0)+\sigma_1 \omega(k_1)+\sigma_2 \omega(k_2))\Big]\Big\}\\
					& \ \  \times\ \Big[\delta(\sigma_0k_0+\sigma_1k_1+\sigma_2k_2)
					|\sin(2\pi k_0^1)||\sin(2\pi k_1^1)||\sin(2\pi k_2^1)|	{F}_1(k_0) F_2(k_2)\Big].\end{aligned}	
			\end{equation}
			Then for all $c>0$ 
			\begin{equation}
				\label{Lemma:Resonance1:6:1}\begin{aligned}& \lim_{\upsilon+\ell\to 0}\lim_{\lambda\to0}|\mathbb{O}_{\upsilon, \lambda}^{2}[F_1,F_2][\phi]-\mathbb{O}_{\upsilon, \lambda}^{2,c}[F_1,F_2][\phi]  | \ = \ 0.\end{aligned}
			\end{equation}
			\item[(e)] The following bounds hold true for the collision operators $\mathcal{Q}^{colli}_1$  \begin{equation}\label{Lemma:Resonance1:7}
				\begin{aligned} 
					& 
					\|\mathcal {Q}^{colli}_1[F_1,F_2,\sigma_0,\sigma_1,\sigma_2,s]\|_{L^\infty}\
					\ \lesssim\	&	\Big\||\sin(2\pi k^1)| F_1(k)\Big\|_{L^4}\Big\||\sin(2\pi k^1)| F_2(k)\Big\|_{L^4},
				\end{aligned}
			\end{equation}
			and $\mathcal{Q}^{colli}_2$
			
			\begin{equation}\label{Lemma:Resonance1:8}
				\begin{aligned} 
					\|\mathcal {Q}^{colli}_2[F_1,F_2,\sigma_0,\sigma_1,\sigma_2,s]\|_{L^\infty}\
					\lesssim\			&	\Big\||\sin(2\pi k^1)|F_1(k)\Big\|_{L^\infty}\Big\||\sin(2\pi k^1)| F_2(k)\Big\|_{L^4}.
				\end{aligned}
			\end{equation}
			\item[(f)] 	For any $\infty\ge p> 2$
			\begin{equation}
				\label{Lemma:Resonance1:9}\begin{aligned}
					&\ \Big[\int_{-T_o\lambda^{-2}}^{T_o\lambda^{-2}} {\mathrm{d}s}\Big\|\mathcal {Q}^{colli}_3[F_1,\sigma_0,\sigma_1,\sigma_2,s]\Big\|_{L^p}^\mathscr{M}\Big]^\frac{1}{\mathscr{M}} \lesssim \||\sin(2\pi k^1)|F_1\|_{L^4}. 
				\end{aligned}
			\end{equation}
			
		\end{itemize}
	\end{lemma}
	
	\begin{proof}
		We first rewrite $\mathbb{O}^1_{\upsilon, \lambda}$ as
		\begin{equation}
			\label{Lemma:Resonance1:E1}\begin{aligned}
				\mathbb{O}^1_{\upsilon, \lambda}[F_1,F_2][\phi] \ =  &\ \int_{-T_o\lambda^{-2}}^{T_o\lambda^{-2}} \frac{\mathrm{d}s}{\pi}e^{-\upsilon |s|}\int_{\mathbb{T}^d}\mathrm{d}k_0\mathcal {Q}^{colli}_1[F_1,F_2,\sigma_0,\sigma_1,\sigma_2,s](k_0)\phi(k_0),
			\end{aligned}
		\end{equation}	for all $F_1,F_2\in L^4(\mathbb{T}^d)$, $\phi\in  C^\infty_c(\mathbb{T}^d\times(-\infty,\infty))$. While the constraint on  $\mathscr{M}$ will be specified later, we denote by $\mathscr{M}'$ the conjugate of $\mathscr{M}$ i.e. $\frac{1}{\mathscr{M}}+\frac{1}{\mathscr{M}'}=1$. By H\"older's inequality, we bound
		\begin{equation}
			\label{Lemma:Resonance1:E1}\begin{aligned}
				|\mathbb{O}^1_{\upsilon, \lambda}[F_1,F_2][\phi]| \ \le  &\ \Big[\int_{-T_o\lambda^{-2}}^{T_o\lambda^{-2}} \frac{\mathrm{d}s}{\pi}e^{-2\mathscr{M}\upsilon |s|}\Big\|\mathcal {Q}^{colli}_1[F_1,F_2,\sigma_0,\sigma_1,\sigma_2,s](k_0)\Big\|^\mathscr{M}_{L^4}\Big]^\frac{1}{\mathscr{M}}\Big[\int_{-T_o\lambda^{-2}}^{T_o\lambda^{-2}} \frac{\mathrm{d}s}{\pi}\|\phi\|^{\mathscr{M}'}_{L^\frac43}\Big]^\frac{1}{\mathscr{M}'}\\
				\ \le  &\ \Big[\int_{-T_o\lambda^{-2}}^{T_o\lambda^{-2}} \frac{\mathrm{d}s}{\pi}e^{-2\mathscr{M}\upsilon |s|}\Big\|\mathcal {Q}^{colli}_1[F_1,F_2,\sigma_0,\sigma_1,\sigma_2,s](k_0)\Big\|^\mathscr{M}_{L^4}\Big]^\frac{1}{\mathscr{M}}\Big[\int_{-\infty}^{\infty} \frac{\mathrm{d}s}{\pi}\|\phi\|^{\mathscr{M}'}_{L^\frac43}\Big]^\frac{1}{\mathscr{M}'}.
			\end{aligned}
		\end{equation}

		The norm of the operator $\mathbb{O}^1_{\upsilon, \lambda}[F]$ is bounded as
		\begin{equation}
			\label{Lemma:Resonance1:E1:a}\begin{aligned}
				\Big\|\mathbb{O}^1_{\upsilon, \lambda}[F_1,F_2]\Big\|^\mathscr{M} \ \le \  \mathbb{O}^1_{o,\upsilon, \lambda}[F]:=  &\ \int_{-T_o\lambda^{-2}}^{T_o\lambda^{-2}} \frac{\mathrm{d}s}{\pi}e^{-2\mathscr{M}\upsilon |s|}\Big\|\mathcal {Q}^{colli}_1[F_1,F_2,\sigma_0,\sigma_1,\sigma_2,s]\Big\|_{L^4}^\mathscr{M}.
			\end{aligned}
		\end{equation}
		Using the identity 
		$$	\delta(\sigma_0k_0+\sigma_1k_1+\sigma_2k_2) \ = \ \sum_{y\in\mathbb{Z}^d}e^{{\bf i}2\pi y\cdot(\sigma_0k_0+\sigma_1k_1+\sigma_2k_2)},$$ 
		we  develop
		\begin{equation}\label{Lemma:Resonance1:E1:1}
			\begin{aligned} 
				& 
				\mathcal {Q}^{colli}_1[F_1,F_2,\sigma_0,\sigma_1,\sigma_2,s](k_0) \\
				=\	& e^{{\bf i}s\sigma_0 \omega(k_0)}\iint_{(\mathbb{T}^d)^2}\!\! \mathrm{d} k_1  \mathrm{d} k_2 \,
				\delta(\sigma_0k_0+\sigma_1k_1+\sigma_2k_2) 
				e^{{\bf i}s\sigma_1 \omega(k_1)}|\sin(2\pi k_1^1)|\\
				&\times  {F}_1(k_1) F_2(k_2)e^{{\bf i} s \sigma_2 \omega(k_2)} |\sin(2\pi k_2^1)|\\
				= &  \
				e^{{\bf i}s\sigma_0 \omega(k_0)}	\sum_{y\in\mathbb{Z}^d}e^{{\bf i}2\pi k_0\sigma_0\cdot y}\iint_{(\mathbb{T}^d)^2}\!\!  \mathrm{d} k_1  \mathrm{d} k_2 
				e^{{\bf i}2\pi y \cdot\sigma_1 k_1+{\bf i}s\sigma_1\omega(k_1)}{F}_1(k_1)F_2(-\sigma_2\sigma_0k_0-\sigma_2\sigma_1k_1)[1+|k_2|^2]^{d+2}\\
				&\times  |\sin(2\pi k_1^1)| |\sin(2\pi k_2^1)| e^{{\bf i}2\pi y \cdot\sigma_2k_2+{\bf i}s\sigma_2\omega(k_2)}[1+|\sigma_2\sigma_0k_0+\sigma_2\sigma_1k_1|^2]^{-d-2}
				\\
				= &\  e^{{\bf i}s\sigma_0 \omega(k_0)}\sum_{y\in\mathbb{Z}^d}e^{{\bf i}2\pi k_0\sigma_0\cdot y}	\mathfrak{F}^{F_1F_2}(y,s\sigma_1)\mathfrak{F}^{o}(y,s\sigma_2),
			\end{aligned}
		\end{equation}
		in which 
		\begin{equation}
			\label{Lemma:Resonance1:E1:1:1}
			\begin{aligned} 
				& 	\mathfrak{F}^{F_1F_2}(x,t) \ = \ \int_{\mathbb{T}^d} \mathrm{d}k |\sin(2\pi k^1)||\sin(2\pi (-\sigma_2\sigma_0k_0^1-\sigma_2\sigma_1k^1))| e^{{\bf i}2\pi x\cdot k} e^{{\bf i}t\omega(k)}F_1(k)F_2(-\sigma_2\sigma_0k_0-\sigma_2\sigma_1k)\\
				&\times [1+|\sigma_2\sigma_0k_0+\sigma_2\sigma_1k|^2]^{-d-2} e^{{\bf i}t\sigma_1\sigma_2\omega(-\sigma_2\sigma_0k_0-\sigma_2\sigma_1k)},		\end{aligned}
		\end{equation}
		$$	\mathfrak{F}^{o}(x) \ = \ \int_{\mathbb{T}^d} \mathrm{d}k  e^{{\bf i}2\pi x\cdot k} [1+|k|^2]^{d+2}.$$
		

		Using H\"older's inequality for the $y$ variable, we find
		\begin{equation}\label{Lemma:Resonance1:E1:2}
			\begin{aligned} 
				& 
				\Big\|\iint_{(\mathbb{T}^d)^2}\!\! \mathrm{d} k_1  \mathrm{d} k_2 \,
				\delta(\sigma_0k_0+\sigma_1k_1+\sigma_2k_2) 
				e^{{\bf i}s\sigma_0 \omega(k_0)}	e^{{\bf i}s\sigma_1 \omega(k_1)}|\sin(2\pi k_1^1)|\\
				&\times {F}_1(k_1)F_2(k_2)e^{{\bf i} s \sigma_2\omega(k_2)} |\sin(2\pi k_2^1)| \Big\|_{L^4}\\
				\le\ &  
				\Big\|\sum_{y\in\mathbb{Z}^d}e^{{\bf i}2\pi k_0\sigma_0\cdot y}	\mathfrak{F}^{F_1F_2}(y,s\sigma_1)\mathfrak{F}^{o}(y)\Big\|_{L^4}\ 
				\\
				\le\ &  
				\Big\|\sum_{y\in\mathbb{Z}^d}e^{{\bf i}2\pi k_0\sigma_0\cdot y}	\mathfrak{F}^{F_1F_2}(y,s\sigma_1)\mathfrak{F}^{o}(y)\Big\|_{L^\infty}\ 
				\le\ \|	\mathfrak{F}^{F_1F_2}(y,s\sigma_1)\mathfrak{F}^{o}(y)\|_{l^1} \\
				\le\ &  
				\|	\mathfrak{F}^{F_1F_2}(y,s\sigma_1)\|_{l^8}\|\mathfrak{F}^{o}(y )\|_{l^\frac87}.
			\end{aligned}
		\end{equation}
		We have
		\begin{equation}
			\|\mathfrak{F}^{o}(y )\|_{l^\frac87}  \ \lesssim \ \  1,
		\end{equation}
		which, in combination with \eqref{Lemma:Resonance1:E1:2}, yields
		\begin{equation}\label{Lemma:Resonance1:E1:3}
			\begin{aligned} 
				& 
				\Big\|\iint_{(\mathbb{T}^d)^2}\!\! \mathrm{d} k_1  \mathrm{d} k_2 \,
				\delta(\sigma_0k_0+\sigma_1k_1+\sigma_2k_2) 
				e^{{\bf i}s\sigma_0 \omega(k_0)}|\sin(2\pi k_1^1)|\\
				&\times {F}_1(k_1)F_2(k_2)e^{{\bf i} s \sigma_2 \omega(k_2)} |\sin(2\pi k_2^1)| e^{{\bf i}2\pi k_0\cdot x_0}\Big\|_{L^4}\ 
				\lesssim \    
				\|\mathfrak{F}^{F_1F_2}(\cdot,s\sigma_1)\|_{l^8}.
			\end{aligned}
		\end{equation}
		Next, we bound
		\begin{equation}
			\label{Lemma:Resonance1:E1:3:a}\begin{aligned}
				\mathbb{O}^1_{o,\upsilon, \lambda}[F] \ =  &\ \int_{-T_o\lambda^{-2}}^{T_o\lambda^{-2}} \frac{\mathrm{d}s}{\pi}e^{-2\upsilon |s|}\Big\|\mathcal {Q}^{colli}_1[F_1,F_2,\sigma_0,\sigma_1,\sigma_2,s]\Big\|_{L^4}^\mathscr{M}\\
				\ \lesssim	&\ \int_{-T_o\lambda^{-2}}^{T_o\lambda^{-2}} {\mathrm{d}s}\|	\mathfrak{F}^{F_1F_2}(y,s\sigma_1)\|_{l^8}^\mathscr{M}.
			\end{aligned}
		\end{equation}
		We now estimate the integral of the first term on the right hand side of \eqref{Lemma:Resonance1:E1:3:a} using a $TT^*$ argument. We choose $\tilde{G}(y,s)$ to be a test function in $L^{\mathscr{M}'}([-T_o\lambda^{-2},T_o\lambda^{-2}],l^\frac87(\mathbb{Z}^d))$ i. e. $\int_{-T_o\lambda^{-2}}^{T_o\lambda^{-2}}{\mathrm{d}s}$ $\Big[ \sum_{y\in\mathbb{Z}^d}|\tilde{G}(y,s)|^\frac87\Big]^{7\mathscr{M}'/8}<\infty$ and develop
		
		\begin{equation}
			\label{Lemma:Resonance1:E1:3:b}\begin{aligned}
				&\  \Big|\sum_{y\in\mathbb{Z}^d}\int_{-T_o\lambda^{-2}}^{T_o\lambda^{-2}} {\mathrm{d}s}	\mathfrak{F}^{F_1F_2}(y,s\sigma_1)\tilde{G}(y,s)\Big|\\
				= \	& \Big|\sum_{y\in\mathbb{Z}^d}\int_{-T_o\lambda^{-2}}^{T_o\lambda^{-2}} {\mathrm{d}s}	\int_{\mathbb{T}^d} \mathrm{d}k |\sin(2\pi k^1)|  e^{{\bf i}2\pi x\cdot k} |\sin(-\sigma_2\sigma_0k_0^1-\sigma_2\sigma_1k^1)|e^{{\bf i}s\sigma_2\omega(-\sigma_2\sigma_0k_0-\sigma_2\sigma_1k)}\\
				&\times e^{{\bf i}s\sigma_1\omega(k)}F_1(k)F_2(-\sigma_2\sigma_0k_0-\sigma_2\sigma_1k)[1+|\sigma_2\sigma_0k_0+\sigma_2\sigma_1k|^2]^{-d-2}\tilde{G}(y,s)\Big|
			\end{aligned}
		\end{equation}
		We  study the $L^2$-norm in $k$
		\begin{equation}
			\label{Lemma:Resonance1:E1:3:c}\begin{aligned}
				& \Big\|\sum_{y\in\mathbb{Z}^d}\int_{-T_o\lambda^{-2}}^{T_o\lambda^{-2}} {\mathrm{d}s}	|\sin(2\pi k^1)|  e^{{\bf i}2\pi x\cdot k} e^{{\bf i}s\sigma_1\omega(k)}e^{{\bf i}s\sigma_2\omega(-\sigma_2\sigma_0k_0-\sigma_2\sigma_1k)}\tilde{G}(y,s)\Big\|_{L^2}^2,
			\end{aligned}
		\end{equation}
		which we will expand and bound as follows
		\begin{equation}
			\label{Lemma:Resonance1:E1:3:d}\begin{aligned}
				& \Big|\sum_{y,y'\in\mathbb{Z}^d}\int_{-T_o\lambda^{-2}}^{T_o\lambda^{-2}} {\mathrm{d}s}\int_{-T_o\lambda^{-2}}^{T_o\lambda^{-2}} {\mathrm{d}s'}	\int_{\mathbb{T}^d}\mathrm{d}k|\sin(2\pi k^1)|^2 e^{{\bf i}2\pi y\cdot k} e^{{\bf i}s\sigma_1\omega(k)}e^{{\bf i}s\sigma_2\omega(-\sigma_2\sigma_0k_0-\sigma_2\sigma_1k)}\\
				&\ \times e^{-{\bf i}2\pi y'\cdot k} e^{-{\bf i}s'\sigma_1\omega(k)}e^{-{\bf i}s'\sigma_2\omega(-\sigma_2\sigma_0k_0-\sigma_2\sigma_1k)}\tilde{G}(y,s)\overline{\tilde{G}(y',s')}\Big|\\
				\lesssim\	& \Big|\sum_{y,y'\in\mathbb{Z}^d}\int_{-T_o\lambda^{-2}}^{T_o\lambda^{-2}} {\mathrm{d}s}\int_{-T_o\lambda^{-2}}^{T_o\lambda^{-2}} {\mathrm{d}s'}	\mathfrak{F}^{Ker}(y-y',s-s')\tilde{G}(y,s)\overline{\tilde{G}(y',s')}\Big|\\
				\lesssim\	& \Big|\int_{-T_o\lambda^{-2}}^{T_o\lambda^{-2}} {\mathrm{d}s}\int_{-T_o\lambda^{-2}}^{T_o\lambda^{-2}} {\mathrm{d}s'}	\|\mathfrak{F}^{Kern}(\cdot,s-s')\|_{l^4}\Big[\sum_{y}|\tilde{G}(y,s)|^\frac87\Big]^\frac78\Big[\sum_{y'}|\tilde{G}(y',s')|^\frac87\Big]\Big|^\frac78\\
				\lesssim\	& \left[\int_{-T_o\lambda^{-2}}^{T_o\lambda^{-2}} {\mathrm{d}s}\left[\sum_{y\in\mathbb{Z}^d}|\tilde{G}(y,s)|^\frac87\right]^{7\mathscr{M}'/8}\right]^\frac{2}{\mathscr{M}'}\left[\int_{-2T_o\lambda^{-2}}^{2T_o\lambda^{-2}} {\mathrm{d}s}\|\mathfrak{F}^{Kern}(\cdot,s)\|_{l^4}^\frac{\mathscr{M}}{2}\right]^\frac{2}{\mathscr{M}},
			\end{aligned}
		\end{equation}
		in which 
		\begin{equation}
			\label{Lemma:Resonance1:E16}	\begin{aligned} 
				{\mathfrak{F}}^{Kern}(x,s) \ = \ & \int_{\mathbb{T}^d} \mathrm{d}k |\sin(2\pi k_0^1)|^2  |\sin(2\pi k^1)|^2|\sin(2\pi (\sigma_0k_0^1+\sigma_1k^1))|^2 e^{{\bf i}2\pi x\cdot k}\\
				&\ \ \ \ \ \ \times e^{{\bf i}s\omega(\sigma_1 k) +{\bf i}s\sigma_2\omega(-\sigma_0k_0-\sigma_1k)}. 	\end{aligned}
		\end{equation}
		
		Now, applying the estimate proved in Lemma \ref{Lemm:AnotherKernel2} to the right hand side of \eqref{Lemma:Resonance1:E1:3:d}, we find
		\begin{equation}\label{Lemma:Resonance1:E1:4}
			\begin{aligned} 
				& \left[\int_{-T_o\lambda^{-2}}^{T_o\lambda^{-2}} {\mathrm{d}s}\left[\sum_{y\in\mathbb{Z}^d}|\tilde{G}(y,s)|^\frac87\right]^{7\mathscr{M}'/8}\right]^\frac{2}{\mathscr{M}'}\left[\int_{-2T_o\lambda^{-2}}^{2T_o\lambda^{-2}} {\mathrm{d}s}\frac{1}{\langle s\rangle^{\frac{\mathscr{M}}{10}-}}\right]^{\frac{2}{\mathscr{M}}}\\
				\	\lesssim\ &	\left[\int_{-T_o\lambda^{-2}}^{T_o\lambda^{-2}} {\mathrm{d}s}\left[\sum_{y\in\mathbb{Z}^d}|\tilde{G}(y,s)|^\frac87\right]^{7\mathscr{M}'/8}\right]^\frac{2}{\mathscr{M}'},
			\end{aligned}
		\end{equation}
		which yields, for $\mathscr{M}$ large
		\begin{equation}
			\label{Lemma:Resonance1:E1:4:1}\begin{aligned}
				\int_{-T_o\lambda^{-2}}^{T_o\lambda^{-2}} {\mathrm{d}s}\|	\mathfrak{F}^{F_1F_2}(y,s\sigma_1)\|_{l^8}^\mathscr{M} \	\lesssim\	& \||\sin(2\pi k^1)|F_1\|_{L^4}^\mathscr{M}\||\sin(2\pi k^1)|F_2\|_{L^4}^\mathscr{M}.
			\end{aligned}
		\end{equation}
		
		Combining \eqref{Lemma:Resonance1:E1:3:a} and \eqref{Lemma:Resonance1:E1:4:1}, we find 
		\begin{equation}\label{Lemma:Resonance1:E1:5}\begin{aligned}
				\mathbb{O}^1_{o,\upsilon, \lambda}[F_1,F_2] \ \lesssim  &\ \| |\sin(2\pi k^1)|F_2\|_{L^4} \| |\sin(2\pi k^1)|F_1\|_{L^4}.
			\end{aligned}
		\end{equation}

		The dominated convergence theorem  yields the following limits, 
		\begin{equation}
			\label{Lemma:Resonance1:E2}\begin{aligned}
				& \lim_{\upsilon\to 0} \int_{-T_o\lambda^{-2}}^{T_o\lambda^{-2}} \frac{\mathrm{d}s}{\pi}e^{-\upsilon |s|}\Big\|\mathcal {Q}^{colli}_1[F_1,F_2,\sigma_0,\sigma_1,\sigma_2,s]\Big\|_{L^4}^\mathscr{M}\\
				\ = \ & \int_{-T_o\lambda^{-2}}^{T_o\lambda^{-2}} \frac{\mathrm{d}s}{\pi}\Big\|\mathcal {Q}^{colli}_1[F_1,F_2,\sigma_0,\sigma_1,\sigma_2,s]\Big\|_{L^4}^\mathscr{M},\end{aligned}	
		\end{equation}	
		and 
		\begin{equation}
			\label{Lemma:Resonance1:E2:1}\begin{aligned}
				& \lim_{\lambda\to 0}\lim_{\upsilon\to 0} \int_{-T_o\lambda^{-2}}^{T_o\lambda^{-2}} \frac{\mathrm{d}s}{\pi}e^{-\upsilon |s|}\Big\|\mathcal {Q}^{colli}_1[F_1,F_2,\sigma_0,\sigma_1,\sigma_2,s]\Big\|_{L^4}^\mathscr{M}\\
				\ = \ & \int_{-\infty}^{\infty} \frac{\mathrm{d}s}{\pi}\Big\|\mathcal {Q}^{colli}_1[F_1,F_2,\sigma_0,\sigma_1,\sigma_2,s]\Big\|_{L^4}^\mathscr{M}.\end{aligned}	
		\end{equation}	

	As the norm $\|\mathbb{O}^1_{\upsilon, \lambda}\|$ of the operator $\mathbb{O}^1_{\upsilon, \lambda}$   is bounded by a universal constant
	\begin{equation}
		\label{Lemma:Resonance1:E2b}\begin{aligned}
			\lim_{\upsilon\to 0}\|\mathbb{O}_{\upsilon, \lambda}^1[F_1,F_2]\|(k)
			\ \lesssim \ & \| |\sin(2\pi k^1)|F_2\|_{L^4} \| |\sin(2\pi k^1)|F_1\|_{L^4}\end{aligned}	\end{equation}	
	and
	\begin{equation}
		\label{Lemma:Resonance1:E2b}\begin{aligned}
			\lim_{\lambda\to 0}\lim_{\upsilon\to 0}\|\mathbb{O}^1_{\upsilon, \lambda}[F_1,F_2]\|(k)
			\ \lesssim \ & \| |\sin(2\pi k^1)|F_2\|_{L^4} \| |\sin(2\pi k^1)|F_1\|_{L^4}.\end{aligned}	\end{equation}	
	The conclusions \eqref{Lemma:Resonance1:2:bis}, \eqref{Lemma:Resonance1:2:bis1}, \eqref{Lemma:Resonance1:2:bis1:a} of (a)  follow. 
	

	Let us now prove (b). We now rewrite $\mathbb{O}^2_{\upsilon, \lambda}$ as
	\begin{equation}
		\label{Lemma:Resonance1:E1:bis}\begin{aligned}
			\mathbb{O}^2_{\upsilon, \lambda}[F_1,F_2][\phi] \ =  &\ \int_{-T_o\lambda^{-2}}^{T_o\lambda^{-2}} \frac{\mathrm{d}s}{\pi}e^{-\upsilon |s|}\int_{\mathbb{T}^d}\mathrm{d}k_0\mathcal {Q}^{colli}_2[F_1,F_2,\sigma_0,\sigma_1,\sigma_2,s](k_0)\phi(k_0),
		\end{aligned}
	\end{equation}	for all $F_1,F_2\in C^\infty(\mathbb{T}^d)$ and $\phi\in C_c^\infty(\mathbb{T}^d\times(-\infty,\infty))$. Again, by H\"older's inequality, the norm of the operator $\mathbb{O}^1_{\upsilon, \lambda}[F]$ is bounded as
	\begin{equation}
		\label{Lemma:Resonance1:E1:a:bis}\begin{aligned}
			&	\Big\|\mathbb{O}^2_{\upsilon, \lambda}[F_1,F_2]\Big\|^\mathscr{M} \ \le \  \mathbb{O}^2_{o,\upsilon, \lambda}[F]:=  \ \int_{-T_o\lambda^{-2}}^{T_o\lambda^{-2}} \frac{\mathrm{d}s}{\pi}e^{-\upsilon |s|}\Big\|\mathcal {Q}^{colli}_2[F_1,F_2,\sigma_0,\sigma_1,\sigma_2,s]\Big\|_{L^4}^\mathscr{M}\\
			\ \le   &\ \| |\sin(2\pi k^1)|F_1(k)\|_{L^\infty}^\mathscr{M}\int_{-T_o\lambda^{-2}}^{T_o\lambda^{-2}} \frac{\mathrm{d}s}{\pi}e^{-\upsilon |s|}\Big\|\mathcal {Q}^{colli}_1[1,F_2,\sigma_0,\sigma_1,\sigma_2,s]\Big\|_{L^4}^\mathscr{M}.
		\end{aligned}
	\end{equation}
	As \eqref{Lemma:Resonance1:E1:a:bis} involves $\mathcal{Q}_1^{colli}$, the argument used to prove (a) can be repeated. The only difference is that, instead of \eqref{Lemma:Resonance1:E1:3}, we use
	\begin{equation}\label{Lemma:Resonance1:E1:3bis}
		\begin{aligned} 
			& 
			\Big\|\iint_{(\mathbb{T}^d)^2}\!\! \mathrm{d} k_1  \mathrm{d} k_2 \,
			\delta(\sigma_0k_0+\sigma_1k_1+\sigma_2k_2) 
			e^{{\bf i}s\sigma_0 \omega(k_0)}	e^{{\bf i}s\sigma_1 \omega(k_1)}|\sin(2\pi k_1^1)|\\
			&\times F_2(k_2)e^{{\bf i} s \sigma_2 \omega(k_2)} |\sin(2\pi k_2^1)|e^{{\bf i}2\pi k_0\cdot x_0}\Big\|_{L^4}\
			\lesssim \   	 \|\mathfrak{F}^{1\cdot  F_2}(\cdot,s\sigma_2)\|_{l^8},
		\end{aligned}
	\end{equation}
	which leads to the conclusions of (b).
	
	Next, we only prove (c) as the proof of (d) is quite the same. To this end, we compute 
	\begin{equation}
		\label{Lemma:Resonance1:E3}\begin{aligned}
			& |\mathbb{O}_{\upsilon,\lambda}^{1,c}[F]-\mathbb{O}_{\upsilon,\lambda}^1[F]| \	=  \  \Big|\iiint_{\big\{|\sigma_0 \omega(k_0)+\sigma_1 \omega(k_1)+\sigma_2 \omega(k_2)|\ge c\big\}}\mathrm{d} k_0  \mathrm{d} k_1  \mathrm{d} k_2 \,
			\\
			&\ \ \frac{	1}{\pi[(\upsilon+\ell)^2+(\sigma_0 \omega(k_0)+\sigma_1 \omega(k_1)+\sigma_2 \omega(k_2))^2]}\Big\{(\upsilon+\ell)-e^{-T_o\lambda^{-2}(\upsilon+\ell)}\Big[\cos(\lambda^{-2}T_o(\sigma_0 \omega(k_0)\\
			& \ \ +\sigma_1 \omega(k_1)+\sigma_2 \omega(k_2))(\upsilon+\ell)	- (\sigma_0 \omega(k_0)+\sigma_1 \omega(k_1)+\sigma_2 \omega(k_2))\\
			& \ \ \times\sin(\lambda^{-2}T_o(\sigma_0 \omega(k_0)+\sigma_1 \omega(k_1)+\sigma_2 \omega(k_2))\Big]\Big\}\\
			& \ \  \times\ \Big[\delta(\sigma_0k_0+\sigma_1k_1+\sigma_2k_2)|\sin(2\pi k_0^1)||\sin(2\pi k_1^1)||\sin(2\pi k_2^1)|	{F}_1(k_0) F_2(k_2)\Big]\\
			\ \lesssim \ &\Big|\iiint_{\big\{|\sigma_0 \omega(k_0)+\sigma_1 \omega(k_1)+\sigma_2 \omega(k_2)|\ge c\big\}}\mathrm{d} k_0  \mathrm{d} k_1  \mathrm{d} k_2 \,\frac{	1}{\pi[(\upsilon+\ell)^2+c^2]} \Big\{(\upsilon+\ell)\\
			& \ \ -e^{-T_o\lambda^{-2}(\upsilon+\ell)}\Big[\cos(\lambda^{-2}T_o(\sigma_0 \omega(k_0) +\sigma_1 \omega(k_1)+\sigma_2 \omega(k_2))(\upsilon+\ell)	- (\sigma_0 \omega(k_0)+\sigma_1 \omega(k_1)\\
			& \ \ +\sigma_2 \omega(k_2))\sin(\lambda^{-2}T_o(\sigma_0 \omega(k_0)+\sigma_1 \omega(k_1)+\sigma_2 \omega(k_2))\Big]\Big\} \Big[\delta(\sigma_0k_0+\sigma_1k_1+\sigma_2k_2)\\
			&\ \ \times |\sin(2\pi k_1^1)||\sin(2\pi k_2^1)|{F}_1(k_1){F}_2(k_2)\Big]\Big|,\end{aligned}	
	\end{equation}
	which implies
	\begin{equation}
		\label{Lemma:Resonance1:E4}\begin{aligned}
			\lim_{\lambda\to0}	|\mathbb{O}_{\upsilon,\lambda}^{1,c}[F]-\mathbb{O}^1_{\upsilon,\lambda}[F]| \	\le  \   &\Big|\iiint_{\big\{k'\in\mathbb{T}^d\big||\sigma_0 \omega(k_0)+\sigma_1 \omega(k_1)+\sigma_2 \omega(k_2)|\ge c\big\}}\mathrm{d} k_0  \mathrm{d} k_1  \mathrm{d} k_2 \,
			\\
			&\frac{	(\upsilon+\ell)}{\pi[(\upsilon+\ell)^2+c^2]} \Big[\delta(\sigma_0k_0+\sigma_1k_1+\sigma_2k_2)
			e^{{\bf i}s[\sigma_0 \omega(k_0)+\sigma_1 \omega(k_1)+\sigma_2 \omega(k_2)]}\\
			&\times \ \  |\sin(2\pi k_1^1)||\sin(2\pi k_2^1)|{F}_1(k_1){F}_2(k_2)\Big]\Big|.\end{aligned}	
	\end{equation}
	Now, taking the limit  $\upsilon+\ell\to 0$, we finally obtain
	\begin{equation}
		\label{Lemma:Resonance1:E5}\begin{aligned}
			\lim_{\upsilon+\ell\to 0}\lim_{\lambda\to0}	|\mathbb{O}_{\upsilon,\lambda}^c[F]-\mathbb{O}_{\upsilon,\lambda}[F]| =0,\end{aligned}	
	\end{equation}
	leading to \eqref{Lemma:Resonance1:6}.

	We will now prove (e). From \eqref{Lemma:Resonance1:E1:1}, we develop 
	\begin{equation}\label{Lemma:Resonance1:E7}
		\begin{aligned} 
			& 
			\|\mathcal {Q}^{colli}_1[F_1,F_2,\sigma_0,\sigma_1,\sigma_2,s]\|_{L^\infty}\
			= \ \Big\|\sum_{y\in\mathbb{Z}^d}e^{{\bf i}2\pi k_0\sigma_0\cdot y}	\mathfrak{F}^{F_1F_2}(y,s\sigma_1)\mathfrak{F}^{o}(y,s\sigma_2)\Big\|_{L^\infty}\\
			\lesssim\		&\Big\|\mathfrak{F}^{F_1F_2}(y,s\sigma_1)\mathfrak{F}^{F_2}(y,s\sigma_2)\Big\|_{l^1} \   \\
			\ \lesssim\	&	\Big\| |\sin(2\pi k^1)| F_1\Big\|_{L^4}\Big\| |\sin(2\pi k^1)| F_2 \Big\|_{L^4},
		\end{aligned}
	\end{equation}
	with $k=(k^1,\cdots,k^d)$,
	which implies the following bound on $\mathcal{Q}^{colli}_2$
	
	\begin{equation}\label{Lemma:Resonance1:E8}
		\begin{aligned} 
			&	\|\mathcal {Q}^{colli}_2[F_1,F_2,\sigma_0,\sigma_1,\sigma_2,s]\|_{L^\infty}\
			\\
			\lesssim\			&	\Big\||\sin(2\pi k^1)|F_1(k)\Big\|_{L^\infty}\Big\||\sin(2\pi k^1)|F_2(k)\Big\|_{L^4}.
		\end{aligned}
	\end{equation}
	We will now prove the last inequality (f). We  develop
	\begin{equation}\label{Lemma:Resonance1:E9}
		\begin{aligned} 
			& 
			\mathcal {Q}^{colli}_3[F_1,\sigma_0,\sigma_1,\sigma_2,s](k_0)\\
			=\	& e^{{\bf i}s\sigma_0 \omega(k_0)}\iint_{(\mathbb{T}^d)^2}\!\! \mathrm{d} k_1  \mathrm{d} k_2 \,
			\delta(\sigma_0k_0+\sigma_1k_1+\sigma_2k_2) 
			e^{{\bf i}s\sigma_1 \omega(k_1)}|\sin(2\pi k_1^1)|\\
			&\times  {F}_1(k_1) e^{-{\bf i} s \omega(\sigma_0k_0^1+\sigma_1k_1^1)} |\sin(2\pi (-\sigma_0k_0-\sigma_1k_1))|\\
			= &  
			e^{{\bf i}s\sigma_0 \omega(k_0)}	\sum_{y\in\mathbb{Z}^d}e^{{\bf i}2\pi k_0\sigma_0\cdot y}\iint_{(\mathbb{T}^d)^2}\!\!  \mathrm{d} k_1  \mathrm{d} k_2 
			e^{{\bf i}2\pi y \cdot\sigma_1 k_1+{\bf i}s\sigma_1\omega(k_1)}\\
			&\times  |\sin(2\pi k_0^1)|  |\sin(2\pi k_1^1)| |\sin(2\pi (\sigma_0k_0+\sigma_1k_1))| e^{{\bf i}2\pi y \cdot\sigma_2k_2+{\bf i}s\omega(-\sigma_0k_0-\sigma_1k_1)}{F}_1(k_1) 
			\\
			= &  e^{{\bf i}s\sigma_0 \omega(k_0)}\sum_{y\in\mathbb{Z}^d}e^{{\bf i}2\pi k_0\sigma_0\cdot y}	\tilde{\mathfrak{F}}^{F_1}(y,s\sigma_1)\mathfrak{F}^{o}(y),
		\end{aligned}
	\end{equation}
	in which 
	\begin{equation}
		\label{Lemma:Resonance1:E10}	\begin{aligned} 
			\tilde{\mathfrak{F}}^{F_1}(x,s) \ = \ & \int_{\mathbb{T}^d} \mathrm{d}k |\sin(2\pi k_0^1)|  |\sin(2\pi k^1)||\sin(2\pi (\sigma_0k_0^1+\sigma_1k^1))| e^{{\bf i}2\pi x\cdot k}\\
			&\ \ \ \ \ \ \times e^{{\bf i}s\omega(\sigma_1 k) +{\bf i}s\omega(-\sigma_0k_0-\sigma_1k)}{F}_1(k)[1+|\sigma_2\sigma_0k_0+\sigma_2\sigma_1k|^2]^{-d-2} . 	\end{aligned}
	\end{equation}
	Using H\"older's inequality for the $y$ variable, we find 
	\begin{equation}\label{Lemma:Resonance1:E12}
		\begin{aligned} 
			& 
			\|\mathcal {Q}^{colli}_3[F_1,\sigma_0,\sigma_1,\sigma_2,s]\|_{L^p}\
			\le \ 
			\Big\|\sum_{y\in\mathbb{Z}^d}e^{{\bf i}2\pi k_0\sigma_0\cdot y}	\tilde{\mathfrak{F}}^{F_1}(y,s\sigma_1){\mathfrak{F}}^{o}(y)\Big\|_{L^\infty}
			\\
			\le\ &\|	\tilde{\mathfrak{F}}^{F_1}(y,s ){\mathfrak{F}}^{o}(y)\|_{l^1} \ 
			\le\    
			\|\tilde{\mathfrak{F}}^{F_1}(y,s )\|_{l^8}\|\mathfrak{{F}}^{o}(y)\|_{l^\frac87}  \ 
			\lesssim \     \|	\tilde{\mathfrak{F}}^{F_1}(y,s )\|_{l^8}.
		\end{aligned}
	\end{equation}
	Next, we bound, similarly as above
	\begin{equation}
		\label{Lemma:Resonance1:E13}\begin{aligned}
			&\ \int_{-T_o\lambda^{-2}}^{T_o\lambda^{-2}} {\mathrm{d}s}\Big\|\mathcal {Q}^{colli}_3[F_1,\sigma_0,\sigma_1,\sigma_2,s]\Big\|_{L^p}^\mathscr{M} 
			\ \lesssim	 \ \int_{-T_o\lambda^{-2}}^{T_o\lambda^{-2}} {\mathrm{d}s} \|	\tilde{\mathfrak{F}}^{F_1}(y,s )\|_{l^8}^\mathscr{M}.
		\end{aligned}
	\end{equation}
	We also estimate the integral on the right hand side of \eqref{Lemma:Resonance1:E13} using a $TT^*$ argument by choosing $\tilde{G}(y,s)$ to be a test function in $L^{\mathscr{M}'}([-T_o\lambda^{-2},T_o\lambda^{-2}],l^\frac87(\mathbb{Z}^d))$ and develop
	
	\begin{equation}
		\label{Lemma:Resonance1:E14}\begin{aligned}
			&\  \Big|\sum_{y\in\mathbb{Z}^d}\int_{-T_o\lambda^{-2}}^{T_o\lambda^{-2}} {\mathrm{d}s}	\tilde{\mathfrak{F}}^{F_1}(y,s\sigma_1)\tilde{G}(y,s)\Big|\ = \  \Big|\sum_{y\in\mathbb{Z}^d}\int_{-T_o\lambda^{-2}}^{T_o\lambda^{-2}} {\mathrm{d}s}	\int_{\mathbb{T}^d} \mathrm{d}k|\sin(2\pi k_0^1)| \\
			& \ \ \times |\sin(2\pi k^1)||\sin(2\pi (\sigma_0k_0^1+\sigma_1k^1))| e^{{\bf i}2\pi x\cdot k}  e^{{\bf i}s\omega(\sigma_1 k) +{\bf i}s\omega(-\sigma_0k_0-\sigma_1k)}\tilde{G}(y,s){F}_1(k)[1+|\sigma_2\sigma_0k_0+\sigma_2\sigma_1k|^2]^{-d-2}\Big|,
		\end{aligned}
	\end{equation}
	which is expanded and bounded as follows
	\begin{equation}
		\label{Lemma:Resonance1:E15}\begin{aligned}
			& \Big|\sum_{y,y'\in\mathbb{Z}^d}\int_{-T_o\lambda^{-2}}^{T_o\lambda^{-2}} {\mathrm{d}s}\int_{-T_o\lambda^{-2}}^{T_o\lambda^{-2}} {\mathrm{d}s'}	\int_{\mathbb{T}^d}\mathrm{d}k|\sin(2\pi k_0^1)|^2 |\sin(2\pi k^1)|^2|\sin(2\pi (\sigma_0k_0^1+\sigma_1k^1))|^2 e^{{\bf i}2\pi (y-y')\cdot k} \\
			&\ \times e^{{\bf i}(s-s')\omega(\sigma_1 k) +{\bf i}(s-s')\omega(-\sigma_0k_0-\sigma_1k)}\tilde{G}(y,s)\overline{\tilde{G}(y',s')}\Big|\\
			\lesssim\	& \left[\int_{-T_o\lambda^{-2}}^{T_o\lambda^{-2}} {\mathrm{d}s}\left[\sum_{y\in\mathbb{Z}^d}|\tilde{G}(y,s)|^\frac87\right]^{7\mathscr{M}'/8}\right]^\frac{2}{\mathscr{M}'}\left[\int_{-2T_o\lambda^{-2}}^{2T_o\lambda^{-2}} {\mathrm{d}s}\|\mathfrak{F}^{Kern}(y,s)\|_{l^4}^{\frac{\mathscr{M}}{2}}\right]^\frac{2}{\mathscr{M}}.
		\end{aligned}
	\end{equation}

	By \eqref{Lemm:AnotherKernel2:1}, the same argument used to prove (a) can be applied yielding the conclusion of (f).
\end{proof}
\begin{remark}\label{Remark:Resonance} 	
	Given $k\in \mathbb{T}^d$, the set of all 	$k'\in \mathbb{T}^d$ satisfying $\omega(k)+\omega(k')-\omega(k+k')=0$ forms a resonance manifold. In order to analyze wave kinetic equations, understanding resonance manifolds is an important step. The problem of analyzing resonance manifolds, as well as the near-resonance consideration, has been studied in \cite{GambaSmithBinh,germain2020optimal,nguyen2017quantum,ToanBinh,PomeauBinh,rumpf2021wave,Binh1} for various types of dispersion relations. 
\end{remark}
\subsection{An oscillatory integral bound}
We follow the same notations used in Section \ref{Subsec:DispersiveEstimates}.

\begin{lemma} Let $q$ be a sufficiently large positive constant. Let   $\tilde\varkappa_1,\tilde\varkappa_2,\tilde\varkappa_3$ be  fixed vectors in $\mathbb{T}^d$ and be different from each other; $\alpha_1,\alpha_2,\alpha_3,\alpha_4\in\{\pm1\}$ be fixed. We define the cut-off function $\tilde\Psi(k_1+ \tilde\varkappa_1,\tilde{\varkappa}_2-\tilde{\varkappa}_1,\tilde{\varkappa}_3-\tilde{\varkappa}_1)= \check{\Psi}(2\pi k_1+ 2\pi\tilde\varkappa_1,2\pi\tilde{\varkappa}_2-2\pi\tilde{\varkappa}_1,2\pi\tilde{\varkappa}_3-2\pi\tilde{\varkappa}_1)$. We have the following estimate
	\label{Lemma:Crossing}
	\begin{equation}
		\begin{aligned}\label{Lemma:Crossing:1}
			&  \left[\iint_{\mathbb{R}^2}\mathrm{d}r_1\mathrm{d}r_2e^{-q\lambda^2|r_1-r_2|/2}e^{-q\lambda^2|r_2|/2}\right.\\
			&\times\left.\left|\int_{\mathbb{T}^{d}}\mathrm{d}k_1 e^{-{\bf i}r_1 \alpha_3\omega(k_1+ \tilde\varkappa_3)-{\bf i}(r_1 \alpha_4+r_2\alpha_1)\omega(k_1+ \tilde\varkappa_1)-{\bf i}r_2\alpha_2\omega(k_1+ \tilde\varkappa_2)}\tilde{\Psi}(k_1+ \tilde\varkappa_1,\tilde{\varkappa}_2-\tilde{\varkappa}_1,\tilde{\varkappa}_3-\tilde{\varkappa}_1)\right|^\frac{q}{2}\right]^\frac1q\\
			\lesssim\	& \langle\ln\lambda\rangle^{\mathfrak{C}}\lambda^{-2/q+\bar{\epsilon}/q}\tilde{\mathfrak{H}} ,
		\end{aligned}
	\end{equation}
	for some constant  $ \mathfrak{C} >0$, $\bar{\epsilon}>0$, 
	where
	\begin{equation}
		\begin{aligned}\label{Lemma:Crossing:2}
			\tilde{\mathfrak{H}}: = \ & \prod_{j=2}^{d}  \Big[\bar{\mathbb{C}}_j^{-\frac{1}{2q_j}}+ \bar{\mathbb{A}}_j^{-\frac{1}{2q_j}}\Big]\sqrt{\tilde{F}(\tilde{\varkappa}_2-\tilde{\varkappa}_1)} \sqrt{\tilde{F}(\tilde{\varkappa}_3-\tilde{\varkappa}_1)},
		\end{aligned}
	\end{equation}
	with $\tilde{F}(\tilde{\varkappa}_2-\tilde{\varkappa}_1)$, $\tilde{F}(\tilde{\varkappa}_3-\tilde{\varkappa}_1)$ are components of $\sqrt{\tilde\Psi}$ that depend only on $\tilde{\varkappa}_2-\tilde{\varkappa}_1$ and $\tilde{\varkappa}_3-\tilde{\varkappa}_1$ respectively, $q_2,\cdots,q_{d}$ are sufficiently large positive real numbers in $(1,\infty)$ such that
	$\frac{1}{q_2}+\cdots+\frac{1}{q_{d}}=\frac{2}{q}$, and
	\begin{equation}
		\begin{aligned}\label{Lemma:Crossing:3}
			\bar{\mathbb{A}}_j   
			\ =\ & |\alpha_3\alpha_1\cos(\tilde{\varkappa}_3^1-\tilde{\varkappa}_1^1){\sin(2(\tilde{\varkappa}_3^{j}-\tilde{\varkappa}_1^{j}))} - \alpha_2\alpha_4\cos(\tilde{\varkappa}_2^1-\tilde{\varkappa}_1^1)\sin(2(\tilde{\varkappa}_2^j-\tilde{\varkappa}_1^{j}))\\
			&+ \alpha_3\alpha_2\cos(\tilde{\varkappa}_3^1-\tilde{\varkappa}_1^1){\sin(2(\tilde{\varkappa}_3^{j}-\tilde{\varkappa}_1^{j}))} \cos(\tilde{\varkappa}_2^1-\tilde{\varkappa}_1^1)\cos(2(\tilde{\varkappa}_2^{j}-\tilde{\varkappa}_1^{j}))\\
			& -\alpha_3\alpha_2 \cos(\tilde{\varkappa}_3^1-\tilde{\varkappa}_1^1){\cos(2(\tilde{\varkappa}_3^{j}-\tilde{\varkappa}_1^{j}))}\cos(\tilde{\varkappa}_2^1-\tilde{\varkappa}_1^1)\sin(2(\tilde{\varkappa}_2^j-\tilde{\varkappa}_1^{j}))|, 
		\end{aligned}
	\end{equation}
	and
	\begin{equation}
		\begin{aligned}\label{Lemma:Crossing:4}
			\bar{\mathbb{C}}_j  
			\ =\ & |\sin(\tilde{\varkappa}_3^1-\tilde{\varkappa}_1^1)\sin(\tilde{\varkappa}_2^1-\tilde{\varkappa}_1^1)\sin(2(\tilde{\varkappa}_2^{j}-\tilde{\varkappa}_3^{j}))|.
		\end{aligned}
	\end{equation}
\end{lemma}
\begin{proof}
	We develop the left hand side of \eqref{Lemma:Crossing:1}
	\begin{equation}
		\begin{aligned}\label{eq:crossingraph12}
			& \Big|\int_{\mathbb{T}^{d}}\mathrm{d}k_1 e^{-{\bf i}r_1 \alpha_3\omega(k_1+ \tilde\varkappa_3)-{\bf i}(r_1 \alpha_4+r_2\alpha_1)\omega(k_1+ \tilde\varkappa_1)-{\bf i}r_2\alpha_2\omega(k_1+ \tilde\varkappa_2)}\tilde\Psi\Big|\\
			= \ & \left|\int_{\mathbb{T}^{d}}\mathrm{d}k_1 e^{-{\bf i}r_1 \alpha_3\omega(k_1+ \tilde\varkappa_3)-{\bf i}(r_1 \alpha_4+r_2\alpha_1)\omega(k_1+ \tilde\varkappa_1)-{\bf i}r_2\alpha_2\omega(k_1+ \tilde\varkappa_2)}\right.
			\\
			&\times \left.\left[\sum_{m_1,m_2\in\mathbb{Z}^d}g(m_1,m_2)e^{{\bf i}2\pi [m_1\cdot k_1+m_2\cdot(k_1+ \tilde\varkappa_2)]}\right]\tilde\Psi\right|,
		\end{aligned}
	\end{equation}
	where $g$ is the  function $\mathbf{1}_{m_1=0}\mathbf{1}_{m_2=0}$, which is the inverse Fourier transform of the function $1$. 
	Now, distributing $e^{{\bf i}2\pi [m_1\cdot k_1+m_2\cdot(k_1+ \tilde\varkappa_2)]}$ into the terms inside the integral of $k_1$, we find
	\begin{equation}
		\begin{aligned}\label{eq:crossingraph13}
			& \Big|\int_{\mathbb{T}^{d}}\mathrm{d}k_1 e^{-{\bf i}r_1 \alpha_3\omega(k_1+ \tilde\varkappa_3)-{\bf i}(r_1 \alpha_4+r_2\alpha_1)\omega(k_1+ \tilde\varkappa_1)-{\bf i}r_2\alpha_2\omega(k_1+ \tilde\varkappa_2)}\tilde\Psi\Big|\\
			= \ & \Big|
			\int_{\mathbb{T}^{d}}\mathrm{d}k_1 \sum_{m_1,m_2\in\mathbb{Z}^d}g(m_1,m_2)e^{-{\bf i}r_1 \alpha_3\omega(k_1+ \tilde\varkappa_3)-{\bf i}(r_1 \alpha_4+r_2\alpha_1)\omega(k_1+ \tilde\varkappa_1)-{\bf i}r_2\alpha_2\omega(k_1+ \tilde\varkappa_2)+{\bf i}2\pi [m_1\cdot k_1+m_2\cdot(k_1+ \tilde\varkappa_2)]}\tilde\Psi\Big|\\
			= \ & \Big|
			\int_{\mathbb{T}^{2d}}\mathrm{d}k_1\mathrm{d}k_*\delta(k_*-k-\tilde{\varkappa}_2) \sum_{m_1,m_2\in\mathbb{Z}^d}g(m_1,m_2)\sqrt{\tilde\Psi}(k+\tilde\varkappa_2)\sqrt{\tilde\Psi}(k_*)e^{{\bf i}2\pi m_2 \cdot k_*}\\
			&\times  e^{-{\bf i}r_1 \alpha_3\omega(k_1+ \tilde\varkappa_3)-{\bf i}(r_1 \alpha_4+r_2\alpha_1)\omega(k_1+ \tilde\varkappa_1)-{\bf i}r_2\alpha_2\omega(k_1+ \tilde\varkappa_2)+{\bf i}2\pi m_1\cdot k_1}\Big|
			.
		\end{aligned}
	\end{equation}
	By the identity
	$\sum_{y\in\mathbb{Z}^d}e^{{\bf i}2\pi ({k}_{*}-k_1-\tilde{\varkappa}_2)\cdot y}=\delta({k}_{*}-k_1-\tilde{\varkappa}_2),$
	we deduce from \eqref{eq:crossingraph13} that
	\begin{equation}
		\begin{aligned}\label{eq:crossingraph13:1}
			& \Big|\int_{\mathbb{T}^{d}}\mathrm{d}k_1 e^{-{\bf i}r_1 \alpha_3\omega(k_1+ \tilde\varkappa_3)-{\bf i}(r_1 \alpha_4+r_2\alpha_1)\omega(k_1+ \tilde\varkappa_1)-{\bf i}r_2\alpha_2\omega(k_1+ \tilde\varkappa_2)}{\tilde\Psi}\Big|\\
			= \ & \Big|
			\int_{\mathbb{T}^{2d}}\mathrm{d}k_1\mathrm{d}k_*\sum_{y\in\mathbb{Z}^d}e^{{\bf i}2\pi ({k}_{*}-k_1-\tilde{\varkappa}_2)\cdot y} \sum_{m_1,m_2\in\mathbb{Z}^d}g(m_1,m_2)\sqrt{\tilde\Psi}(k_1+ \tilde\varkappa_2)\\
			&\times  e^{-{\bf i}r_1 \alpha_3\omega(k_1+ \tilde\varkappa_3)-{\bf i}(r_1 \alpha_4+r_2\alpha_1)\omega(k_1+ \tilde\varkappa_1)-{\bf i}r_2\alpha_2\omega(k_1+ \tilde\varkappa_2)+{\bf i}2\pi m_1\cdot k_1}\sqrt{\tilde\Psi}(k_*)e^{{\bf i}2\pi m_2\cdot k_*}\Big|\\
			= \ & \Big|
			\int_{\mathbb{T}^{2d}}\mathrm{d}k_1\mathrm{d}k_*\sum_{y\in\mathbb{Z}^d}e^{-{\bf i}2\pi \tilde{\varkappa}_2\cdot y} \sum_{m_1,m_2\in\mathbb{Z}^d}g(m_1,m_2)\sqrt{\tilde\Psi}(k_1+ \tilde\varkappa_2)\\
			&\times  e^{-{\bf i}r_1 \alpha_3\omega(k_1+ \tilde\varkappa_3)-{\bf i}(r_1 \alpha_4+r_2\alpha_1)\omega(k_1+ \tilde\varkappa_1)-{\bf i}r_2\alpha_2\omega(k_1+ \tilde\varkappa_2)+{\bf i}2\pi (m_1-y)\cdot k_1}\sqrt{\tilde\Psi}(k_*)e^{{\bf i}2\pi (m_2+y)\cdot k_*}\Big|\\
			= \ & \left|
			\sum_{m_1,m_2\in\mathbb{Z}^d}g(m_1,m_2)\sum_{y\in\mathbb{Z}^d}e^{-{\bf i}2\pi \tilde{\varkappa}_2\cdot y}\mathfrak{H}_1(m_1-y)\mathfrak{H}_2(m_2+y)\right|
			,
		\end{aligned}
	\end{equation}	 
	
	in which
	\begin{equation}
		\begin{aligned}\label{eq:crossingraph14}
			\mathfrak{H}_1(m)= \  & \int_{\mathbb{T}^{d}}\mathrm{d}k_1e^{-{\bf i}r_1 \alpha_3\omega(k_1+ \tilde\varkappa_3)-{\bf i}(r_1 \alpha_4+r_2\alpha_1)\omega(k_1+ \tilde\varkappa_1)-{\bf i}r_2\alpha_2\omega(k_1+ \tilde\varkappa_2)+{\bf i}2\pi m\cdot k_1}\sqrt{\tilde\Psi}(k_1+ \tilde\varkappa_2),\\
			\mathfrak{H}_2(m)= \  & \int_{\mathbb{T}^{d}}\mathrm{d}{k}_{*}\sqrt{\tilde\Psi}(k_*)e^{{\bf i}2\pi m\cdot k_*}.
		\end{aligned}
	\end{equation}

	By H\"oder's inequality, applied to the right hand side of \eqref{eq:crossingraph13}, it follows  
	\begin{equation}
		\begin{aligned}\label{eq:crossingraph15}&\Big|\int_{\mathbb{T}^{d}}\mathrm{d}k_1 e^{-{\bf i}r_1 \alpha_3\omega(k_1+ \tilde\varkappa_3)-{\bf i}(r_1 \alpha_4+r_2\alpha_1)\omega(k_1+ \tilde\varkappa_1)-{\bf i}r_2\alpha_2\omega(k_1+ \tilde\varkappa_2)}{\tilde\Psi}\Big|\\
			\lesssim \  & \left|\sum_{m_1,m_2\in\mathbb{Z}^d}g(m_1,m_2)\right|\big\|\mathfrak{H}_1\big\|_4\big\|\mathfrak{H}_2\big\|_\frac43\lesssim \   \big\|\mathfrak{H}_1\big\|_4\big\|\mathfrak{H}_2\big\|_\frac43,
		\end{aligned}
	\end{equation}
	where we have used the fact that $g$ is the inverse Fourier transform of the function $1$. We will now bound $\|\mathfrak{H}_2\big\|_\frac43$
	\begin{equation}
		\begin{aligned}\label{eq:crossingraph14:1}
			\|\mathfrak{H}_2\big\|_\frac43\ = \  &\left\{\sum_{m\in\mathbb{Z}^d} \left|\int_{\mathbb{T}^{d}}\mathrm{d}{k}_{*}{\tilde\Psi}e^{{\bf i}2\pi m\cdot k_*}\right|^\frac43\right\}^\frac34\\
			\ = \ & \left\{\sum_{m\in\mathbb{Z}^d\backslash\{0\}} \left|\int_{\mathbb{T}^{d}}\mathrm{d}{k}_{*}\frac{\sqrt{\tilde\Psi}(k_*)}{|{\bf i}2\pi m|^{2d}}\Delta^{d}\Big(e^{{\bf i}2\pi m\cdot k_*}\Big)\right|^\frac43+\left|\int_{\mathbb{T}^{d}}\mathrm{d}{k}_{*}{{\tilde\Psi}}(k_*) \right|^\frac43\right\}^\frac34\\
			\ = \ & \left\{\sum_{m\in\mathbb{Z}^d\backslash\{0\}} 
			\frac{1}{|2\pi m|^{8d/3}}\left|\int_{\mathbb{T}^{d}}\mathrm{d}{k}_{*}\Delta^{d}\Big[{{\sqrt{\tilde\Psi}}(k_*)}\Big]e^{{\bf i}2\pi m\cdot k_*}\right|^\frac43+\left|\int_{\mathbb{T}^{d}}\mathrm{d}{k}_{*}{{\tilde\Psi}}(k_*)\right|^\frac43\right\}^\frac34\\
			\ \lesssim \ & \left\{\sum_{m\in\mathbb{Z}^d\backslash\{0\}} 
			\frac{1}{|2\pi m|^{8d/3}}\right\}^\frac34\left|\int_{\mathbb{T}^{d}}\mathrm{d}{k}_{*}\Big|\Delta^{d}\Big[{{\sqrt{\tilde\Psi}}(k_*)}\Big]\Big|\right| + \left|\int_{\mathbb{T}^{d}}\mathrm{d}{k}_{*}{{\tilde\Psi}}(k_*)\right|\ \le \ \langle\ln|\lambda|\rangle^{\mathfrak{C}_{\mathfrak{H}_2}}.
		\end{aligned}
	\end{equation}
	with $\mathfrak{C}_{\mathfrak{H}_2}>0$.	
	As $\mathfrak{H}_2$ clearly belongs to $l^\frac43$, whose bound introduces an additional factor of $\langle\ln|\lambda|\rangle^{\mathfrak{C}_{\mathfrak{H}_2}}$, we obtain
	\begin{equation}
		\begin{aligned}\label{eq:crossingraph15:1:1}&\Big|\int_{\mathbb{T}^{d}}\mathrm{d}k_1 e^{-{\bf i}r_1 \alpha_3\omega(k_1+ \tilde\varkappa_3)-{\bf i}(r_1 \alpha_4+r_2\alpha_1)\omega(k_1+ \tilde\varkappa_1)-{\bf i}r_2\alpha_2\omega(k_1+ \tilde\varkappa_2)} {\tilde\Psi}\Big|
			\lesssim \   \big\|\mathfrak{H}_1\big\|_4\langle\ln|\lambda|\rangle^{\mathfrak{C}_{\mathfrak{H}_2}}.
		\end{aligned}
	\end{equation}
	
	In the next step, we  will apply inequality \eqref{Lemm:Bessel3:1} to the  specific case of $\mathfrak{H}_1$.

	To estimate $\mathfrak{H}_1$, we first perform the change of variable $k_1+\tilde{\varkappa}_1\to k_1$ and apply \eqref{Lemm:Bessel3:1} for $t_0=r_1 \alpha_4+r_2\alpha_1$, $t_1=r_2\alpha_2,$ $t_2=r_1\alpha_3$, for $V=\tilde\varkappa_2- \tilde\varkappa_1$ and  $W=\tilde\varkappa_3- \tilde\varkappa_1$,  yielding
	\begin{equation}
		\label{eq:crossingraph15:1}\begin{aligned}
			& \|\mathfrak{H}_1\|_{l^4} \ \lesssim \ 	 \prod_{j=2}^{d}\Big\langle\min\Big\{\Big|r_1 \alpha_4+r_2\alpha_1+r_2\alpha_2\cos(\tilde{\varkappa}_2^1-\tilde{\varkappa}_1^1)e^{{\bf i}2(\tilde{\varkappa}_2^j-\tilde{\varkappa}_1^j)}\\
			&+r_1\alpha_3\cos(\tilde{\varkappa}_3^1-\tilde{\varkappa}_1^1)e^{{\bf i}2(\tilde{\varkappa}_3^j-\tilde{\varkappa}_1^j)}\Big|,\Big|{r_2}\alpha_2\sin(\tilde{\varkappa}_2^1-\tilde{\varkappa}_1^1)e^{{\bf i}2(\tilde{\varkappa}_2^j-\tilde{\varkappa}_1^j)}\\
			& +{r_1}\alpha_3\sin(\tilde{\varkappa}_3^1-\tilde{\varkappa}_1^1)e^{{\bf i}2(\tilde{\varkappa}_3^j-\tilde{\varkappa}_1^j)}\Big|\Big\}|1-|\cos(\aleph^j
			_{1}-\aleph^j
			_{2})||^{\frac{1}{2}}\Big\rangle^{-(\frac{1}{8}-)}{\sqrt{\tilde\Psi}(\tilde{\varkappa}_2-\tilde{\varkappa}_1)}{\sqrt{\tilde\Psi}(\tilde{\varkappa}_3-\tilde{\varkappa}_1)},\end{aligned}
	\end{equation}
	with $\tilde\varkappa_i = (\tilde\varkappa_i^1,\cdots,\tilde\varkappa_i^d)$ for $i=1,2,3$. The quantities $\aleph^j
	_{1}=\aleph^j
	_{1}(V,W)$, $\aleph^j
	_{2}=\aleph^j
	_{2}(V,W)$ are defined in \eqref{Lemm:Bessel}. 
	
	Now, when	$ r_*=t_1/t_2=(1+\epsilon_{r_*})\tilde{r}_l$ for $l=1,2,3$, as in \eqref{Lemm:Bessel3:1:a} in which $\tilde{r}_l$ are defined in Lemma \ref{Lemm:Angle} and \eqref{Lemm:Bessel3:1:b} is satisfied, we have the estimate, using \eqref{Lemm:Bessel3:2}
	\begin{equation}
		\label{eq:crossingraph15:2}\begin{aligned}
			&	\|\mathfrak{H}_1\|_{l^4} \ \lesssim \ 	 \langle\ln|\lambda|\rangle^{\mathfrak{C}_{\mathfrak{H}_1}}\prod_{j=2}^{d}\Big\langle \Big\{\Big|r_1 \alpha_4+r_2\alpha_1+r_2\alpha_2\cos(\tilde{\varkappa}_2^1-\tilde{\varkappa}_1^1)e^{{\bf i}2(\tilde{\varkappa}_2^j-\tilde{\varkappa}_1^j)}\\
			&+r_1\alpha_3\cos(\tilde{\varkappa}_3^1-\tilde{\varkappa}_1^1)e^{{\bf i}2(\tilde{\varkappa}_3^j-\tilde{\varkappa}_1^j)}\Big|+\Big|{r_2}\alpha_2\sin(\tilde{\varkappa}_2^1-\tilde{\varkappa}_1^1)e^{{\bf i}2(\tilde{\varkappa}_2^j-\tilde{\varkappa}_1^j)}\\
			& +{r_1}\alpha_3\sin(\tilde{\varkappa}_3^1-\tilde{\varkappa}_1^1)e^{{\bf i}2(\tilde{\varkappa}_3^j-\tilde{\varkappa}_1^j)}\Big|\Big\}|\cos(\aleph^j
			_{1}-\aleph^j
			_{2})|^{\frac{1}{2}}\Big\rangle^{-(\frac{1}{8}-)}{\sqrt{\tilde\Psi}(\tilde{\varkappa}_2-\tilde{\varkappa}_1)}{\sqrt{\tilde\Psi}(\tilde{\varkappa}_3-\tilde{\varkappa}_1)},\end{aligned}
	\end{equation}
	for some constant $\mathfrak{C}_{\mathfrak{H}_1}>0$.
	We denote
	\begin{equation}
		\mathcal{S}_{r_1,r_2}\ := \	\Big\{(r_1,r_2)\in\mathbb{R}^2 ~~\Big|~~   |(r_*-\tilde{r}_l)/\tilde{r}_l|\le \epsilon_{r_*}  \Big\}, \mbox{ and } \mathcal{S}_{r_1,r_2}'= \mathbb{R}^2\backslash \mathcal{S}_{r_1,r_2}.
	\end{equation}
	Let us recall the definition above $r_*=t_1/t_2$, in which  $t_1=r_2\alpha_2,$ $t_2=r_1\alpha_3$. Thus $r_*$ provides a relation between $r_1,r_2$, which means $\mathcal{S}_{r_1,r_2}$ is well-defined.
	
	We then get

	\begin{equation}
		\begin{aligned}\label{eq:crossingraph17}
			&  \left[\iint_{\mathbb{R}^2}\mathrm{d}r_1\mathrm{d}r_2e^{-q\lambda^2|r_1-r_2|/2}e^{-q\lambda^2|r_2|/2}\right.\\
			&\times\left.\left|\int_{\mathbb{T}^{d}}\mathrm{d}k_1 e^{-{\bf i}r_1 \alpha_3\omega(k_1+ \tilde\varkappa_3)-{\bf i}(r_1 \alpha_4+r_2\alpha_1)\omega(k_1+ \tilde\varkappa_1)-{\bf i}r_2\alpha_2\omega(k_1+ \tilde\varkappa_2)}\Psi\right|^\frac{q}{2}\right]^\frac1q\\
			\lesssim &  \langle\ln|\lambda|\rangle^{\mathfrak{C}_{\mathfrak{H}_1}}\Big\{\int_{\mathbb{R}}\int_{\mathbb{R}}\mathrm{d}r_1\mathrm{d}r_2\Big\{e^{-\lambda^2|r_1-r_2|^2}e^{-\lambda^2|r_2|^2}\chi_{\mathcal{S}_{r_1,r_2}} \prod_{j=2}^{d}\Big\langle\min\Big\{\Big|r_1 \alpha_4+r_2\alpha_1\\
			&\ \ +r_2\alpha_2\cos(\tilde{\varkappa}_2^1-\tilde{\varkappa}_1^1)e^{{\bf i}2(\tilde{\varkappa}_2^j-\tilde{\varkappa}_1^j)}+r_1\alpha_3\cos(\tilde{\varkappa}_3^1-\tilde{\varkappa}_1^1)e^{{\bf i}2(\tilde{\varkappa}_3^j-\tilde{\varkappa}_1^j)}\Big|,\Big|{r_2}\alpha_2\sin(\tilde{\varkappa}_2^1-\tilde{\varkappa}_1^1)e^{{\bf i}2(\tilde{\varkappa}_2^j-\tilde{\varkappa}_1^j)}\\
			&\ \ +{r_1}\alpha_3\sin(\tilde{\varkappa}_3^1-\tilde{\varkappa}_1^1)e^{{\bf i}2(\tilde{\varkappa}_3^j-\tilde{\varkappa}_1^j)}\Big|\Big\}|1-|\cos(\aleph^j
			_{1}-\aleph^j
			_{2})||^{\frac{1}{2}}\Big\rangle^{-(\frac{1}{8}-)} {\sqrt{\tilde\Psi}(\tilde{\varkappa}_2-\tilde{\varkappa}_1)}{\sqrt{\tilde\Psi}(\tilde{\varkappa}_3-\tilde{\varkappa}_1)}\Big\}^\frac{q}{2}\Big\}^{\frac{1}{q}}\\
			&\ \ +\langle\ln|\lambda|\rangle^{\mathfrak{C}_{\mathfrak{H}_1}} \Big\{\int_{\mathbb{R}}\int_{\mathbb{R}}\mathrm{d}r_1\mathrm{d}r_2\Big\{e^{-|r_1-r_2|}e^{-|r_2|}\chi_{\mathcal{S}_{r_1,r_2}'}{\sqrt{\tilde\Psi}(\tilde{\varkappa}_2-\tilde{\varkappa}_1)}{\sqrt{\tilde\Psi}(\tilde{\varkappa}_3-\tilde{\varkappa}_1)}\\
			&\ \ \times \prod_{j=2}^{d}\Big\langle\Big\{\Big|r_1 \alpha_4+r_2\alpha_1+r_2\alpha_2\cos(\tilde{\varkappa}_2^1-\tilde{\varkappa}_1^1)e^{{\bf i}2(\tilde{\varkappa}_2^j-\tilde{\varkappa}_1^j)}+r_1\alpha_3\cos(\tilde{\varkappa}_3^1-\tilde{\varkappa}_1^1)e^{{\bf i}2(\tilde{\varkappa}_3^j-\tilde{\varkappa}_1^j)}\Big|\\
			&\ \ +\Big|{r_2}\alpha_2\sin(\tilde{\varkappa}_2^1-\tilde{\varkappa}_1^1)e^{{\bf i}2(\tilde{\varkappa}_2^j-\tilde{\varkappa}_1^j)}  +{r_1}\alpha_3\sin(\tilde{\varkappa}_3^1-\tilde{\varkappa}_1^1)e^{{\bf i}2(\tilde{\varkappa}_3^j-\tilde{\varkappa}_1^j)}\Big|\Big\}|\cos(\aleph^j
			_{1}-\aleph^j
			_{2})|^{\frac{1}{2}}\Big\rangle^{-(\frac{1}{8}-)} \Big\}^\frac{q}{2}\Big\}^{\frac{1}{q}},
		\end{aligned}
	\end{equation}
	where $\chi_{\mathcal{S}_{r_1,r_2}}$ and $\chi_{\mathcal{S}_{r_1,r_2}'}$ are the characteristic functions of $\mathcal{S}_{r_1,r_2}$ and $\mathcal{S}_{r_1,r_2}'$, which are
	$\chi_{\mathcal{S}_{r_1,r_2}} = 1 \mbox{ when } (r_1,r_2)\in \mathcal{S}_{r_1,r_2},$ $\chi_{\mathcal{S}_{r_1,r_2}} = 0 \mbox{ when } (r_1,r_2)\notin \mathcal{S}_{r_1,r_2},$ $\chi_{\mathcal{S}_{r_1,r_2}'} = 1 \mbox{ when } (r_1,r_2)\in \mathcal{S}_{r_1,r_2}',$ $\chi_{\mathcal{S}_{r_1,r_2}'} = 0 \mbox{ when } (r_1,r_2)\notin \mathcal{S}_{r_1,r_2}', $
	yielding
	\begin{equation}
		\begin{aligned}\label{eq:crossingraph17:a1}
			&  \left[\iint_{\mathbb{R}^2}\mathrm{d}r_1\mathrm{d}r_2e^{-q\lambda^2|r_1-r_2|/2}e^{-q\lambda^2|r_2|/2}\right.\\
			&\times\left.\left|\int_{\mathbb{T}^{d}}\mathrm{d}k_1 e^{-{\bf i}r_1 \alpha_3\omega(k_1+ \tilde\varkappa_3)-{\bf i}(r_1 \alpha_4+r_2\alpha_1)\omega(k_1+ \tilde\varkappa_1)-{\bf i}r_2\alpha_2\omega(k_1+ \tilde\varkappa_2)}\tilde{\Psi}\right|^\frac{q}{2}\right]^\frac1q\\
			\lesssim & \langle\ln|\lambda|\rangle^{\mathfrak{C}_{\mathfrak{H}_1}} \Big\{\int_{\mathbb{R}}\int_{\mathbb{R}}\mathrm{d}r_1\mathrm{d}r_2\Big\{e^{-\lambda^2|r_1-r_2|}e^{-\lambda^2|r_2|}\chi_{\mathcal{S}_{r_1,r_2}} \prod_{j=2}^{d}\Big\{\Big\langle\Big|r_1 \alpha_4+r_2\alpha_1+r_2\alpha_2\cos(\tilde{\varkappa}_2^1-\tilde{\varkappa}_1^1)e^{{\bf i}2(\tilde{\varkappa}_2^j-\tilde{\varkappa}_1^j)}\\
			&\ \   +r_1\alpha_3\cos(\tilde{\varkappa}_3^1-\tilde{\varkappa}_1^1)e^{{\bf i}2(\tilde{\varkappa}_3^j-\tilde{\varkappa}_1^j)}\Big||1-|\cos(\aleph^j
			_{1}-\aleph^j
			_{2})||^{\frac{1}{2}}\Big\rangle^{-(\frac{1}{8}-)}\\
			& \ \  +\ \Big\langle\Big|{r_2}\alpha_2\sin(\tilde{\varkappa}_2^1-\tilde{\varkappa}_1^1)e^{{\bf i}2(\tilde{\varkappa}_2^j-\tilde{\varkappa}_1^j)}+{r_1}\alpha_3\sin(\tilde{\varkappa}_3^1-\tilde{\varkappa}_1^1)e^{{\bf i}2(\tilde{\varkappa}_3^j-\tilde{\varkappa}_1^j)}\Big|\\
			&\ \ \times |1-|\cos(\aleph^j
			_{1}-\aleph^j
			_{2})||^{\frac{1}{2}}\Big\rangle^{-\frac{1}{10}}\Big\}{\sqrt{\tilde\Psi}(\tilde{\varkappa}_2-\tilde{\varkappa}_1)}{\sqrt{\tilde\Psi}(\tilde{\varkappa}_3-\tilde{\varkappa}_1)}\Big\}^\frac{q}{2}\Big\}^{\frac{1}{q}}\\
			&\ \ +\langle\ln|\lambda|\rangle^{\mathfrak{C}_{\mathfrak{H}_1}} \Big\{\int_{\mathbb{R}}\int_{\mathbb{R}}\mathrm{d}r_1\mathrm{d}r_2\Big\{e^{-\lambda^2|r_1-r_2|}e^{-\lambda^2|r_2|}\chi_{\mathcal{S}_{r_1,r_2}'} \prod_{j=2}^{d}\Big\langle\Big\{\Big|r_1 \alpha_4+r_2\alpha_1+\\
			&\ \ r_2\alpha_2\cos(\tilde{\varkappa}_2^1-\tilde{\varkappa}_1^1)e^{{\bf i}2(\tilde{\varkappa}_2^j-\tilde{\varkappa}_1^j)} +r_1\alpha_3\cos(\tilde{\varkappa}_3^1-\tilde{\varkappa}_1^1)e^{{\bf i}2(\tilde{\varkappa}_3^j-\tilde{\varkappa}_1^j)}\Big|+\Big|{r_2}\alpha_2\sin(\tilde{\varkappa}_2^1-\tilde{\varkappa}_1^1)e^{{\bf i}2(\tilde{\varkappa}_2^j-\tilde{\varkappa}_1^j)}\\
			&\ \ +{r_1}\alpha_3\sin(\tilde{\varkappa}_3^1-\tilde{\varkappa}_1^1)e^{{\bf i}2(\tilde{\varkappa}_3^j-\tilde{\varkappa}_1^j)}\Big|\Big\}|\cos(\aleph^j
			_{1}-\aleph^j
			_{2})|^{\frac{1}{2}}\Big\rangle^{-(\frac{1}{8}-)} {\sqrt{\tilde\Psi}(\tilde{\varkappa}_2-\tilde{\varkappa}_1)}{\sqrt{\tilde\Psi}(\tilde{\varkappa}_3-\tilde{\varkappa}_1)}\Big\}^\frac{q}{2}\Big\}^{\frac{1}{q}},
		\end{aligned}
	\end{equation}	in which $\aleph^j_{1},\aleph^j_{2}$ are defined  in \eqref{Lemm:Bessel}. Note that in the above estimate, \eqref{Lemm:Bessel3:1} is used for $\mathcal{S}_{r_1,r_2}$ and  \eqref{Lemm:Bessel3:2} is applied for $\mathcal{S}_{r_1,r_2}'$.

	We then estimate \eqref{eq:crossingraph17} as
	\begin{equation}
		\begin{aligned}\label{eq:crossingraph17:aa:1}
			&  \left[\iint_{\mathbb{R}^2}\mathrm{d}r_1\mathrm{d}r_2e^{-q\lambda^2|r_1-r_2|/2}e^{-q\lambda^2|r_2|/2}\right.\\
			&\times\left.\left|\int_{\mathbb{T}^{d}}\mathrm{d}k_1 e^{-{\bf i}r_1 \alpha_3\omega(k_1+ \tilde\varkappa_3)-{\bf i}(r_1 \alpha_4+r_2\alpha_1)\omega(k_1+ \tilde\varkappa_1)-{\bf i}r_2\alpha_2\omega(k_1+ \tilde\varkappa_2)}\tilde{\Psi}\right|^\frac{q}{2}\right]^\frac1q\\
			\lesssim & \  \langle\ln|\lambda|\rangle^{\mathfrak{C}_{\mathfrak{H}_1}}\Big\{\int_{\mathbb{R}}\int_{\mathbb{R}}\mathrm{d}r_1\mathrm{d}r_2\Big\{e^{-\lambda^2|r_1-r_2|^2}e^{-\lambda^2|r_2|^2}\\
			&\ \ \ \times \chi_{\mathcal{S}_{r_1,r_2}}\prod_{j=2}^{d}\Big\{\Big\langle\Big|r_1 \alpha_4+r_2\alpha_1+r_2\alpha_2\cos(\tilde{\varkappa}_2^1-\tilde{\varkappa}_1^1)e^{{\bf i}2(\tilde{\varkappa}_2^j-\tilde{\varkappa}_1^j)}\\
			&\ \   +r_1\alpha_3\cos(\tilde{\varkappa}_3^1-\tilde{\varkappa}_1^1)e^{{\bf i}2(\tilde{\varkappa}_3^j-\tilde{\varkappa}_1^j)}\Big||1-|\cos(\aleph^j
			_{1}-\aleph^j
			_{2})||^{\frac{1}{2}}\Big\rangle^{-(\frac{1}{8}-)}+\ \Big\langle\Big|{r_2}\alpha_2\sin(\tilde{\varkappa}_2^1-\tilde{\varkappa}_1^1)e^{{\bf i}2(\tilde{\varkappa}_2^j-\tilde{\varkappa}_1^j)}\\
			& \ \  +{r_1}\alpha_3\sin(\tilde{\varkappa}_3^1-\tilde{\varkappa}_1^1)e^{{\bf i}2(\tilde{\varkappa}_3^j-\tilde{\varkappa}_1^j)}\Big||1-|\cos(\aleph^j
			_{1}-\aleph^j
			_{2})||^{\frac{1}{2}}\Big\rangle^{-(\frac{1}{8}-)}\Big\}{\sqrt{\tilde\Psi}(\tilde{\varkappa}_2-\tilde{\varkappa}_1)}{\sqrt{\tilde\Psi}(\tilde{\varkappa}_3-\tilde{\varkappa}_1)}\Big\}^\frac{q}{2}\Big\}^{\frac{1}{q}}\\
			&\ \ +\langle\ln|\lambda|\rangle^{\mathfrak{C}_{\mathfrak{H}_1}} \Big\{\int_{\mathbb{R}}\int_{\mathbb{R}}\mathrm{d}r_1\mathrm{d}r_2\Big\{e^{-\lambda^2|r_1-r_2|^2}e^{-\lambda^2|r_2|^2}\\
			&\ \ \times \chi_{\mathcal{S}_{r_1,r_2}'}\prod_{j=2}^{d}\Big\langle\Big\{\Big|r_1 \alpha_4+r_2\alpha_1+r_2\alpha_2\cos(\tilde{\varkappa}_2^1-\tilde{\varkappa}_1^1)e^{{\bf i}2(\tilde{\varkappa}_2^j-\tilde{\varkappa}_1^j)}\\
			&\ \ +r_1\alpha_3\cos(\tilde{\varkappa}_3^1-\tilde{\varkappa}_1^1)e^{{\bf i}2(\tilde{\varkappa}_3^j-\tilde{\varkappa}_1^j)}\Big|+\Big|{r_2}\alpha_2\sin(\tilde{\varkappa}_2^1-\tilde{\varkappa}_1^1)e^{{\bf i}2(\tilde{\varkappa}_2^j-\tilde{\varkappa}_1^j)}\\
			&\ \ +{r_1}\alpha_3\sin(\tilde{\varkappa}_3^1-\tilde{\varkappa}_1^1)e^{{\bf i}2(\tilde{\varkappa}_3^j-\tilde{\varkappa}_1^j)}\Big|\Big\}|\cos(\aleph^j
			_{1}-\aleph^j
			_{2})|^{\frac{1}{2}}\Big\rangle^{-(\frac{1}{8}-)} {\sqrt{\tilde\Psi}(\tilde{\varkappa}_2-\tilde{\varkappa}_1)}{\sqrt{\tilde\Psi}(\tilde{\varkappa}_3-\tilde{\varkappa}_1)}\Big\}^\frac{q}{2}\Big\}^{\frac{1}{q}},
		\end{aligned}
	\end{equation}
	which implies the following bound for \eqref{eq:crossingraph17}
	\begin{equation}
		\begin{aligned}\label{eq:crossingraph17:1}
			&  \left[\iint_{\mathbb{R}^2}\mathrm{d}r_1\mathrm{d}r_2e^{-q\lambda^2|r_1-r_2|/2}e^{-q\lambda^2|r_2|/2}\right.\\
			&\times\left.\left|\int_{\mathbb{T}^{d}}\mathrm{d}k_1 e^{-{\bf i}r_1 \alpha_3\omega(k_1+ \tilde\varkappa_3)-{\bf i}(r_1 \alpha_4+r_2\alpha_1)\omega(k_1+ \tilde\varkappa_1)-{\bf i}r_2\alpha_2\omega(k_1+ \tilde\varkappa_2)}\tilde{\Psi}\right|^\frac{q}{2}\right]^\frac1q\\
			\lesssim\	 &  \langle\ln|\lambda|\rangle^{\mathfrak{C}_{\mathfrak{H}_1}}\mathbb{S}_A+ \langle\ln|\lambda|\rangle^{\mathfrak{C}_{\mathfrak{H}_1}}\mathbb{S}_B,
		\end{aligned}
	\end{equation}
	where
	
	\begin{equation}
		\begin{aligned}\label{eq:crossingraph17:1:A:1}
			\mathbb{S}_A:=	& \sqrt{\sqrt{\tilde\Psi}(\tilde{\varkappa}_2-\tilde{\varkappa}_1)}\sqrt{\sqrt{\tilde\Psi}(\tilde{\varkappa}_3-\tilde{\varkappa}_1)} \Big\{\int_{\mathbb{R}^2}\mathrm{d}r_1\mathrm{d}r_2\Big\{e^{-\lambda^2|r_1-r_2|}e^{-\lambda^2|r_2|}\chi_{\mathcal{S}_{r_1,r_2}}\prod_{j=2}^{d}|1-|\cos(\aleph^j
			_{1}-\aleph^j
			_{2})||^{\frac{1}{2}}\\
			&\ \  \times \Big\{\Big\langle\Big|r_1 \alpha_4+r_2\alpha_1+r_2\alpha_2\cos(\tilde{\varkappa}_2^1-\tilde{\varkappa}_1^1)e^{{\bf i}2(\tilde{\varkappa}_2^j-\tilde{\varkappa}_1^j)} +r_1\alpha_3\cos(\tilde{\varkappa}_3^1-\tilde{\varkappa}_1^1)e^{{\bf i}2(\tilde{\varkappa}_3^j-\tilde{\varkappa}_1^j)}\Big|\Big\rangle^{-(\frac{1}{8}-)}\\
			& \ \  +\ \Big\langle\Big|{r_2}\alpha_2\sin(\tilde{\varkappa}_2^1-\tilde{\varkappa}_1^1)e^{{\bf i}2(\tilde{\varkappa}_2^j-\tilde{\varkappa}_1^j)}+{r_1}\alpha_3\sin(\tilde{\varkappa}_3^1-\tilde{\varkappa}_1^1)e^{{\bf i}2(\tilde{\varkappa}_3^j-\tilde{\varkappa}_1^j)}\Big| \Big\rangle^{-(\frac{1}{8}-)}\Big\}^\frac{q}{2}\Big\}^{\frac{1}{q}},
		\end{aligned}
	\end{equation}
	and
	
	\begin{equation}
		\begin{aligned}\label{eq:crossingraph17:1A:2}
			&	\mathbb{S}_B:=		\  \Big\{\int_{\mathbb{R}^2}\mathrm{d}r_1\mathrm{d}r_2\Big\{e^{-\lambda^2|r_1-r_2|}e^{-\lambda^2|r_2|}\chi_{\mathcal{S}_{r_1,r_2}'} \prod_{j=2}^{d}\Big\langle\Big\{\Big|r_1 \alpha_4+r_2\alpha_1+r_2\alpha_2\cos(\tilde{\varkappa}_2^1-\tilde{\varkappa}_1^1)e^{{\bf i}2(\tilde{\varkappa}_2^j-\tilde{\varkappa}_1^j)}\\
			&\ \ +r_1\alpha_3\cos(\tilde{\varkappa}_3^1-\tilde{\varkappa}_1^1)e^{{\bf i}2(\tilde{\varkappa}_3^j-\tilde{\varkappa}_1^j)}\Big|+\Big|{r_2}\alpha_2\sin(\tilde{\varkappa}_2^1-\tilde{\varkappa}_1^1)e^{{\bf i}2(\tilde{\varkappa}_2^j-\tilde{\varkappa}_1^j)}\\
			&\ \ +{r_1}\alpha_3\sin(\tilde{\varkappa}_3^1-\tilde{\varkappa}_1^1)e^{{\bf i}2(\tilde{\varkappa}_3^j-\tilde{\varkappa}_1^j)}\Big|\Big\}|\cos(\aleph^j
			_{1}-\aleph^j
			_{2})|^{\frac{1}{2}}\Big\rangle^{-(\frac{1}{8}-)} \Big\}^\frac{q}{2}\Big\}^{\frac{1}{q}}\sqrt{\sqrt{\tilde\Psi}(\tilde{\varkappa}_2-\tilde{\varkappa}_1)}\sqrt{\sqrt{\tilde\Psi}(\tilde{\varkappa}_3-\tilde{\varkappa}_1)}.
		\end{aligned}
	\end{equation}
	Next, we will estimate $\mathbb{S}_A$ and $\mathbb{S}_B$ separately.
	\smallskip

	Let $q_2,\cdots,q_{d}$ be positive real numbers in $(1,\infty)$ such that
	$\frac{1}{q_2}+\cdots+\frac{1}{q_{d}}=\frac{2}{q}$. Observing that the total power of $r_1$ in the denominator is $(d-1)(\frac{q}{16}-)$, by increasing the value of $q$, we can choose $q_2,\cdots,q_{d}$ such that $q_j(\frac{1}{8}-)>  1$ for $j\in\{2,\cdots,d\}$. 	
	Applying H\"older's inequality to the integrals of $r_1,r_2$ in \eqref{eq:crossingraph17:1}, we find
	
	\begin{equation}
		\begin{aligned}\label{eq:crossingraph17:1:1}
			&\ \Big\{\int_{\mathbb{R}^2}\mathrm{d}r_1\mathrm{d}r_2\Big\{\chi_{\mathcal{S}_{r_1,r_2}}e^{-\lambda^2|r_1-r_2|}e^{-\lambda^2|r_2|} \prod_{j=2}^{d}\Big\{\Big\langle\Big|r_1 \alpha_4+r_2\alpha_1+r_2\alpha_2\cos(\tilde{\varkappa}_2^1-\tilde{\varkappa}_1^1)e^{{\bf i}2(\tilde{\varkappa}_2^j-\tilde{\varkappa}_1^j)}\\
			&\ \   +r_1\alpha_3\cos(\tilde{\varkappa}_3^1-\tilde{\varkappa}_1^1)e^{{\bf i}2(\tilde{\varkappa}_3^j-\tilde{\varkappa}_1^j)}\Big||1-|\cos(\aleph^j
			_{1}-\aleph^j
			_{2})||^{\frac{1}{2}}\Big\rangle^{-(\frac{1}{8}-)}\\
			&\  + \ \Big\langle\Big|{r_2}\alpha_2\sin(\tilde{\varkappa}_2^1-\tilde{\varkappa}_1^1)e^{{\bf i}2(\tilde{\varkappa}_2^j-\tilde{\varkappa}_1^j)}+{r_1}\alpha_3\sin(\tilde{\varkappa}_3^1-\tilde{\varkappa}_1^1)e^{{\bf i}2(\tilde{\varkappa}_3^j-\tilde{\varkappa}_1^j)}\Big||1-|\cos(\aleph^j
			_{1}-\aleph^j
			_{2})||^{\frac{1}{2}}\Big\rangle^{-(\frac{1}{8}-)}\Big\}\Big\}^\frac{q}{2}\Big\}^{\frac{2}{q}}\\
			\lesssim	&\   \prod_{j=2}^{d}\Big\{\Big\|\chi_{\mathcal{S}_{r_1,r_2}} e^{-\lambda^2|r_1-r_2|-\lambda^2|r_2|}\Big\langle\Big|r_1 \alpha_4+r_2\alpha_1+r_2\alpha_2\cos(\tilde{\varkappa}_2^1-\tilde{\varkappa}_1^1)e^{{\bf i}2(\tilde{\varkappa}_2^j-\tilde{\varkappa}_1^j)}\\
			&\ \   +r_1\alpha_3\cos(\tilde{\varkappa}_3^1-\tilde{\varkappa}_1^1)e^{{\bf i}2(\tilde{\varkappa}_3^j-\tilde{\varkappa}_1^j)}\Big||1-|\cos(\aleph^j
			_{1}-\aleph^j
			_{2})||^{\frac{1}{2}}\Big\rangle^{-(\frac{1}{8}-)}\Big\|_{L^{q_j}}\\
			&\  + \Big\|\chi_{\mathcal{S}_{r_1,r_2}} e^{-\lambda^2|r_1-r_2|-\lambda^2|r_2|}\Big\langle\Big|{r_2}\alpha_2\sin(\tilde{\varkappa}_2^1-\tilde{\varkappa}_1^1)e^{{\bf i}2(\tilde{\varkappa}_2^j-\tilde{\varkappa}_1^j)}+{r_1}\alpha_3\sin(\tilde{\varkappa}_3^1-\tilde{\varkappa}_1^1)e^{{\bf i}2(\tilde{\varkappa}_3^j-\tilde{\varkappa}_1^j)}\Big|\\
			&\ \ \times |1-|\cos(\aleph^j
			_{1}-\aleph^j
			_{2})||^{\frac{1}{2}}\Big\rangle^{-(\frac{1}{8}-)}\Big\|_{L^{q_j}}\Big\},
		\end{aligned}
	\end{equation}
	in which the norms $L^{q_j}$, with $j\in\{2,\cdots,d\}$, are taken with respect to both  $r_1$, $r_2$ variables. Since we chose $q_j(\frac{1}{8}-)>1$ with $j\in\{2,\cdots,d\}$, all of those norms are bounded. We will  estimate \eqref{eq:crossingraph17:1:1} below.

	Using    inequality \eqref{Lemm:Angle:00} of Lemma \ref{Lemm:Angle} to bound the angle term $|1-|\cos(\aleph^j
	_{1}-\aleph^j
	_{2})||^{\frac{1}{2}}$ in  the $j$-th term in the product on the right hand side of \eqref{eq:crossingraph17:1:1}, we estimate
	\begin{equation}\begin{aligned}\label{eq:crossingraph17:1:2:1}
			&\   \Big\|\Big\langle \chi_{\mathcal{S}_{r_1,r_2}} e^{-\lambda^2(|r_1-r_2|+|r_2|)}\Big|r_1 \alpha_4+r_2\alpha_1+r_2\alpha_2\cos(\tilde{\varkappa}_2^1-\tilde{\varkappa}_1^1)e^{{\bf i}2(\tilde{\varkappa}_2^j-\tilde{\varkappa}_1^j)}\\
			&\    +r_1\alpha_3\cos(\tilde{\varkappa}_3^1-\tilde{\varkappa}_1^1)e^{{\bf i}2(\tilde{\varkappa}_3^j-\tilde{\varkappa}_1^j)}\Big||1-|\cos(\aleph^j
			_{1}-\aleph^j
			_{2})||^{\frac{1}{2}}\Big\rangle^{-(\frac{1}{8}-)}\Big\|_{L^{q_j}}\\
			&\ + \ \Big\|\Big\langle \chi_{\mathcal{S}_{r_1,r_2}} e^{-\lambda^2(|r_1-r_2|+|r_2|)}\Big\langle\Big|{r_2}\alpha_2\sin(\tilde{\varkappa}_2^1-\tilde{\varkappa}_1^1)e^{{\bf i}2(\tilde{\varkappa}_2^j-\tilde{\varkappa}_1^j)}+{r_1}\alpha_3\sin(\tilde{\varkappa}_3^1-\tilde{\varkappa}_1^1)e^{{\bf i}2(\tilde{\varkappa}_3^j-\tilde{\varkappa}_1^j)}\Big|\\
			&\  \times |1-|\cos(\aleph^j
			_{1}-\aleph^j
			_{2})||^{\frac{1}{2}}\Big\rangle^{-(\frac{1}{8}-)}\Big\|_{L^{q_j}}\\
			\lesssim	&\ \Big|\int_{\mathbb{R}^2}\mathrm{d}r_1 \mathrm{d}r_2\chi_{\mathcal{S}_{r_1,r_2}}e^{-\lambda^2(|r_1-r_2|+|r_2|)q_j}\Big\langle\Big|r_1 \alpha_4+r_2\alpha_1+r_2\alpha_2\cos(\tilde{\varkappa}_2^1-\tilde{\varkappa}_1^1)e^{{\bf i}2(\tilde{\varkappa}_2^j-\tilde{\varkappa}_1^j)}\\
			&\    +r_1\alpha_3\cos(\tilde{\varkappa}_3^1-\tilde{\varkappa}_1^1)e^{{\bf i}2(\tilde{\varkappa}_3^j-\tilde{\varkappa}_1^j)}\Big|\Big[\frac{r_1^2+r_2^2}{|\mathcal{A}_jr_1^2+\mathcal{B}_j r_1r_2+\mathcal{C}_jr_2^2|}+1\Big]^{-1}\Big\rangle^{-{q_j}(\frac{1}{8}-)}\Big|^{\frac{1}{q_j}}\\
			&\ +  \Big|\int_{\mathbb{R}^2}\mathrm{d}r_1 \mathrm{d}r_2\chi_{\mathcal{S}_{r_1,r_2}}e^{-\lambda^2(|r_1-r_2|+|r_2|)q_j}\Big\langle\Big|{r_2}\alpha_2\sin(\tilde{\varkappa}_2^1-\tilde{\varkappa}_1^1)e^{{\bf i}2(\tilde{\varkappa}_2^j-\tilde{\varkappa}_1^j)}\\
			&\ +{r_1}\alpha_3\sin(\tilde{\varkappa}_3^1-\tilde{\varkappa}_1^1)e^{{\bf i}2(\tilde{\varkappa}_3^j-\tilde{\varkappa}_1^j)}\Big|\Big[\frac{r_1^2+r_2^2}{|\mathcal{A}_jr_1^2+\mathcal{B}_j r_1r_2+\mathcal{C}_jr_2^2|}+1\Big]^{-1}\Big\rangle^{-{q_j}(\frac{1}{8}-)}\Big|^{\frac{1}{q_j}}.
		\end{aligned}
	\end{equation}
	Let us recall that, following  inequality \eqref{Lemm:Angle:00} of Lemma \ref{Lemm:Angle}, 	\begin{equation}\begin{aligned}\mathcal{A}_j \ = & \  -\alpha_3\sin(\tilde{\varkappa}_3^1-\tilde{\varkappa}_1^1)\sin(2(\tilde{\varkappa}_3^j-\tilde{\varkappa}_1^j)),\
			\ \ \ \ 
			\mathcal{C}_j \ =  \   
			-\alpha_2\sin(\tilde{\varkappa}_2^1-\tilde{\varkappa}_1^1)\cos(2(\tilde{\varkappa}_2^j-\tilde{\varkappa}_1^j)),\\
			\alpha_2\alpha_3\mathcal{B}_j  \ = & \    -\sin[\tilde{\varkappa}_3^1-\tilde{\varkappa}_2^1]\sin[2(\tilde{\varkappa}_3^j-\tilde{\varkappa}_2^j)]\\
			&\ -\alpha_3\sin(\tilde{\varkappa}_3^1-\tilde{\varkappa}_1^1)\sin(2(\tilde{\varkappa}_3^j-\tilde{\varkappa}_1^j)) -\alpha_2\sin(\tilde{\varkappa}_2^1-\tilde{\varkappa}_1^1)\cos(2(\tilde{\varkappa}_2^j-\tilde{\varkappa}_1^j)),\end{aligned}
	\end{equation}
	for all $j=2,\cdots,d$.
	We also denote $\bar{\mathcal{B}}_j=\alpha_1+\alpha_2\cos(\tilde{\varkappa}_2^1-\tilde{\varkappa}_1^1)e^{{\bf i}2(\tilde{\varkappa}_2^j-\tilde{\varkappa}_1^j)}$, $\bar{\mathcal{A}}_j=\alpha_3\cos(\tilde{\varkappa}_3^1-\tilde{\varkappa}_1^1)e^{{\bf i}2(\tilde{\varkappa}_3^j-\tilde{\varkappa}_1^j)} + \alpha_4,$ $\bar{\mathcal{C}}_j=\alpha_3\sin(\tilde{\varkappa}_3^1-\tilde{\varkappa}_1^1)e^{{\bf i}2(\tilde{\varkappa}_3^j-\tilde{\varkappa}_1^j)}$ and $\bar{\mathcal{D}}_j=\alpha_2\sin(\tilde{\varkappa}_2^1-\tilde{\varkappa}_1^1)e^{{\bf i}2(\tilde{\varkappa}_2^j-\tilde{\varkappa}_1^j)}$ for all $j=2,\cdots,d$, and rewrite the right hand side of \eqref{eq:crossingraph17:1:2:1} as
	\begin{equation}\begin{aligned}\label{eq:crossingraph17:1:2:1:c}
			&\ \Big|\int_{\mathbb{R}^2}\mathrm{d}r_1 \mathrm{d}r_2e^{-\lambda^2(|r_1-r_2|+|r_2|)q_j}\chi_{\mathcal{S}_{r_1,r_2}}\Big\langle\Big|r_1 \bar{\mathcal{A}}_j+r_2\bar{\mathcal{B}}_j\Big|\Big[\frac{r_1^2+r_2^2}{|\mathcal{A}_jr_1^2+\mathcal{B}_j r_1r_2+\mathcal{C}_jr_2^2|}+1\Big]^{-1}\Big\rangle^{-{q_j}(\frac{1}{8}-)}\Big|^{\frac{1}{q_j}}\\
			&\ \ \ \ + \ \Big|\int_{\mathbb{R}^2}\mathrm{d}r_1 \mathrm{d}r_2e^{-\lambda^2(|r_1-r_2|+|r_2|)q_j}\chi_{\mathcal{S}_{r_1,r_2}}\\
			&\ \ \ \ \times \Big\langle\Big|r_1 \bar{\mathcal{C}}_j+r_2\bar{\mathcal{D}}_j\Big|\Big[\frac{r_1^2+r_2^2}{|\mathcal{A}_jr_1^2+\mathcal{B}_j r_1r_2+\mathcal{C}_jr_2^2|}+1\Big]^{-1}\Big\rangle^{-{q_j}(\frac{1}{8}-)}\Big|^{\frac{1}{q_j}}\  := \ \mathcal{U}_j^{\frac{1}{q_j}} +\mathcal{V}_j^{\frac{1}{q_j}}.
		\end{aligned}
	\end{equation}
	We will now estimate $\mathcal{U}_j$ and $\mathcal{V}_j$. We develop $\mathcal{U}_j$ as follow
	\begin{equation}\begin{aligned}\label{eq:crossingraph17:1:2:1:d}
			\mathcal{U}_j \ =	&\ \int_{\mathbb{R}^2}\mathrm{d}r_1 \mathrm{d}r_2\Big\langle\Big|r_1 \bar{\mathcal{A}}_j+r_2\bar{\mathcal{B}}_j\Big|\frac{|\mathcal{A}_jr_1^2+\mathcal{B}_j r_1r_2+\mathcal{C}_jr_2^2|}{|\mathcal{A}_jr_1^2+\mathcal{B}_j r_1r_2+\mathcal{C}_jr_2^2|+r_1^2+r_2^2}\Big\rangle^{-{q_j}(\frac{1}{8}-)} e^{-\lambda^2(|r_1-r_2|+|r_2|)q_j}\chi_{\mathcal{S}_{r_1,r_2}}.
		\end{aligned}
	\end{equation}

	Now, due to the cut-off functions $\sqrt{\sqrt{\tilde\Psi}(\tilde{\varkappa}_3^j-\tilde{\varkappa}_1^j)}$ and $\sqrt{\sqrt{\tilde\Psi}(\tilde{\varkappa}_2^j-\tilde{\varkappa}_1^j)}$, it follows that 		$ |\mathcal{A}_j|, |\mathcal{C}_j|\ge \mathfrak{C}_{\aleph_1^j,\aleph_2^j}^1|\ln\lambda|^{-\mathfrak{C}_{\aleph_1^j,\aleph_2^j}^2}>0,$ for some constants $ \mathfrak{C}_{\aleph_1^j,\aleph_2^j}^1,\mathfrak{C}_{\aleph_1^j,\aleph_2^j}^2>0$. We observe that
	\begin{equation}\begin{aligned}\label{eq:crossingraph17:1:2:1:e}
			|\mathcal{A}_jr_1^2+\mathcal{B}_j r_1r_2+\mathcal{C}_jr_2^2| \ = \ 	& |\mathcal{A}_j(r_1-\mathcal{R}_1r_2)(r_1-\mathcal{R}_2r_2)|
			\ \ge \ 	 |\mathcal{A}_j||r_1-\mathrm{Re}(\mathcal{R}_1)r_2||r_1-\mathrm{Re}(\mathcal{R}_2)r_2|,
		\end{aligned}
	\end{equation}
	where $\mathcal{R}_1$ and $\mathcal{R}_2$ are the two solutions, which could be complex, of $\mathcal{A}_j\mathcal{R}^2+\mathcal{B}_j \mathcal{R}+\mathcal{C}_j=0$ and $\mathrm{Re}(\mathcal{R}_1),\mathrm{Re}(\mathcal{R}_2)$ are the real parts of  $\mathcal{R}_1,\mathcal{R}_2$. We suppose that $\mathrm{Re}(\mathcal{R}_1),\mathrm{Re}(\mathcal{R}_2)\ne 0$, as those points making $\mathrm{Re}(\mathcal{R}_1),\mathrm{Re}(\mathcal{R}_2)= 0$ correspond to values of the angles that will be later integrated out. 
	Noticing that $\overline{\mathcal{R}}_1=1/\mathcal{R}_1$, $\overline{\mathcal{R}}_2=1/\mathcal{R}_2$ are also the solutions of $\mathcal{A}_j+\mathcal{B}_j \overline{\mathcal{R}}+\mathcal{C}_j\overline{\mathcal{R}}^2=0$.

Again, by setting $r^*=\frac{r_1}{r_2}$, 
we will now estimate the fraction appearing in the expression of $\mathcal{U}_j$. Let us  recall the definition  $r_*=t_1/t_2$, in which  $t_1=r_2\alpha_2,$ $t_2=r_1\alpha_3$. Thus $r_*$ and $r^*$ are different only in the signs. We divide this computation into two smaller cases.

{\it	First, we consider the case when $|r^*|\le 4\max\{|\mathrm{Re}\mathcal{R}_1|,|\mathrm{Re}\mathcal{R}_2|,|\mathrm{Re}\mathcal{R}_1+\mathrm{Re}\mathcal{R}_2|,1\}.$} We bound

\begin{equation}\begin{aligned}\label{eq:crossingraph17:1:2:5}
\mathcal{U}_j \ \le	& \ \int_{\mathbb{R}^2}\mathrm{d}r_1 \mathrm{d}r_2\Big\langle\Big|r_1 \bar{\mathcal{A}}_j+r_2\bar{\mathcal{B}}_j\Big|\frac{|\mathcal{A}_jr_1^2+\mathcal{B}_j r_1r_2+\mathcal{C}_jr_2^2|}{r_1^2+r_2^2+|\mathcal{A}_jr_1^2+\mathcal{B}_j r_1r_2+\mathcal{C}_jr_2^2|}\Big\rangle^{-{q_j}(\frac{1}{8}-)}\\
&\times e^{-\lambda^2(|r_1-r_2|+|r_2|)q_j}\chi_{\mathcal{S}_{r_1,r_2}}\\
\ \lesssim	& \ \int_{\mathbb{R}^2}\mathrm{d}r_1 \mathrm{d}r_2\Big\langle|r_1 \bar{\mathcal{A}}_j+r_2\bar{\mathcal{B}}_j|\frac{|\mathcal{A}_j(r^*)^2+\mathcal{B}_j r^*+\mathcal{C}_j|}{(r^*)^2+1}\Big\rangle^{-{q_j}(\frac{1}{8}-)}e^{-\lambda^2(|r_1-r_2|+|r_2|)q_j}\chi_{\mathcal{S}_{r_1,r_2}}\\
\ \lesssim	& \ \int_{\mathbb{R}^2}\mathrm{d}r_1 \mathrm{d}r_2\left|\frac{|r_1 \bar{\mathcal{A}}_j+r_2\bar{\mathcal{B}}_j|{|\mathcal{A}_j(r^*)^2+\mathcal{B}_j r^*+\mathcal{C}_j|}+{(r^*)^2+1}}{(r^*)^2+1}\right|^{-\bar{\epsilon}_j}\\
&\times \Big\langle\frac{|r_1 \bar{\mathcal{A}}_j+r_2\bar{\mathcal{B}}_j|{|\mathcal{A}_j(r^*)^2+\mathcal{B}_j r^*+\mathcal{C}_j|}}{(r^*)^2+1}\Big\rangle^{-{q_j}(\frac{1}{8}-)+\bar{\epsilon}_j}e^{-\lambda^2(|r_1-r_2|+|r_2|)q_j}\chi_{\mathcal{S}_{r_1,r_2}}
\\
\ \lesssim	& \ \int_{\mathbb{R}^2}\mathrm{d}r_1 \mathrm{d}r_2\left|\frac{(r^*)^2+1}{|r_1 \bar{\mathcal{A}}_j+r_2\bar{\mathcal{B}}_j|{|\mathcal{A}_jr_*^2+\mathcal{B}_j r_*+\mathcal{C}_j|}+{(r^*)^2+1}}\right|^{\bar{\epsilon}_j}\\
&\times \Big\langle|r_1 \bar{\mathcal{A}}_j+r_2\bar{\mathcal{B}}_j|\frac{|\mathcal{A}_j(r^*)^2+\mathcal{B}_j r^*+\mathcal{C}_j|}{(r^*)^2+1}\Big\rangle^{-{q_j}(\frac{1}{8}-)+\bar{\epsilon}_j} e^{-\lambda^2(|r_1-r_2|+|r_2|)q_j}\chi_{\mathcal{S}_{r_1,r_2}},
\end{aligned}
\end{equation}
for some small constant $\bar{\epsilon}_j>0$. 
By writing $|\mathcal{A}_j(r^*)^2+\mathcal{B}_j r_*+\mathcal{C}_j|=|\mathcal{A}_j(r^*-\mathcal{R}_1)(r^*-\mathcal{R}_2)|\ge |\mathcal{A}_j||r^*-\mathrm{Re}\mathcal{R}_1||r^*-\mathrm{Re}\mathcal{R}_2|$, we deduce from \eqref{eq:crossingraph17:1:2:5} that
\begin{equation}\begin{aligned}\label{eq:crossingraph17:1:2:6}
\mathcal{U}_j
\ \lesssim	& \ \int_{\mathbb{R}^2}\mathrm{d}r_1 \mathrm{d}r_2\left|\frac{(r^*)^2+1}{|r_1 \bar{\mathcal{A}}_j+r_2\bar{\mathcal{B}}_j||\mathcal{A}_j||r^*-\mathrm{Re}\mathcal{R}_1||r^*-\mathrm{Re}\mathcal{R}_2|+{(r^*)^2+1}}\right|^{\bar{\epsilon}_j}\chi_{\mathcal{S}_{r_1,r_2}}\\
&\ \times \Big\langle\frac{|r_1 \bar{\mathcal{A}}_j+r_2\bar{\mathcal{B}}_j||\mathcal{A}_j||r^*-\mathrm{Re}\mathcal{R}_1||r^*-\mathrm{Re}\mathcal{R}_2|}{(r^*)^2+1}\Big\rangle^{-{q_j}(\frac{1}{8}-)+\bar{\epsilon}_j} e^{-\lambda^2(|r_1-r_2|+|r_2|)q_j}\\
\ \lesssim	& \ \int_{\mathbb{R}^2}\mathrm{d}r_1 \mathrm{d}r_2\left|\frac{(r^*)^2+1}{|r_1 \bar{\mathcal{A}}_j+r_2\bar{\mathcal{B}}_j||\mathcal{A}_j||r^*-\mathrm{Re}\mathcal{R}_1||r^*-\mathrm{Re}\mathcal{R}_2|}\right|^{\bar{\epsilon}_j}\chi_{\mathcal{S}_{r_1,r_2}}\\
&\ \times \Big\langle\frac{|r_1 \bar{\mathcal{A}}_j+r_2\bar{\mathcal{B}}_j||\mathcal{A}_j||r^*-\mathrm{Re}\mathcal{R}_1||r^*-\mathrm{Re}\mathcal{R}_2|}{(r^*)^2+1}\Big\rangle^{-{q_j}(\frac{1}{8}-)+\bar{\epsilon}_j} e^{-\lambda^2(|r_1-r_2|+|r_2|)q_j}\\
\ \lesssim	& \ \int_{\mathbb{R}^2}\mathrm{d}r_1 \mathrm{d}r_2\left|\frac{\max\{|\mathrm{Re}\mathcal{R}_1|,|\mathrm{Re}\mathcal{R}_2|,|\mathrm{Re}\mathcal{R}_1+\mathrm{Re}\mathcal{R}_2|\}^2+1}{|r_1\bar{\mathcal{A}}_j+r_2\bar{\mathcal{B}}_j||\mathcal{A}_j||r^*-\mathrm{Re}\mathcal{R}_1||r^*-\mathrm{Re}\mathcal{R}_2|}\right|^{\bar{\epsilon}_j}\chi_{\mathcal{S}_{r_1,r_2}}\\
&\ \times \Big\langle\frac{|r_1\bar{\mathcal{A}}_j+r_2\bar{\mathcal{B}}_j||\mathcal{A}_j||r^*-\mathrm{Re}\mathcal{R}_1||r^*-\mathrm{Re}\mathcal{R}_2|}{\max\{|\mathrm{Re}\mathcal{R}_1|,|\mathrm{Re}\mathcal{R}_2|,|\mathrm{Re}\mathcal{R}_1+\mathrm{Re}\mathcal{R}_2|\}^2+1}\Big\rangle^{-{q_j}(\frac{1}{8}-)+\bar{\epsilon}_j} e^{-\lambda^2(|r_1-r_2|+|r_2|)q_j}.
\end{aligned}
\end{equation}
Now, by inspecting only the imaginary part of $ \bar{\mathcal{A}}_jr_1+ \bar{\mathcal{B}}_jr_2$, we could bound
\begin{equation}\begin{aligned}\label{eq:crossingraph17:1:3:2:1}
\  | \bar{\mathcal{A}}_jr_1+ \bar{\mathcal{B}}_jr_2|\ \ge & \ \ \Big|\alpha_3\cos(\tilde{\varkappa}_3^1-\tilde{\varkappa}_1^1){\sin(2(\tilde{\varkappa}_3^{j}-\tilde{\varkappa}_1^{j}))}r_1\\
&\ +\alpha_2\cos(\tilde{\varkappa}_2^1-\tilde{\varkappa}_1^1){\sin(2(\tilde{\varkappa}_2^{j}-\tilde{\varkappa}_1^{j}))}r_2\Big| =: |r_1A_{j}'+r_2B_{j}'|,
\end{aligned}
\end{equation}
in which  $\bar{\mathcal{A}}_j=A_{j}^o+{\bf i}A_{j}'$ and $\bar{\mathcal{B}}_j=B_{j}^o+{\bf i}B_{j}'$.	
We suppose that $|A_{j}'|\ge|B_{j}'|$ and set ${B}_{j}''=B_{j}'/A_{j}'$. We develop
\begin{equation}\begin{aligned}\label{eq:crossingraph17:1:2:7}
\mathcal{U}_j
\ \lesssim	& \   \int_{\mathbb{R}^2}\mathrm{d}r_1 \mathrm{d}r_2\left|\frac{\max\{|\mathrm{Re}\mathcal{R}_1|,|\mathrm{Re}\mathcal{R}_2|,|\mathrm{Re}\mathcal{R}_1+\mathrm{Re}\mathcal{R}_2|\}^2+1}{|r_1 {{A}}_j'+r_2{{B}}'_j||\mathcal{A}_j||r^*-\mathrm{Re}\mathcal{R}_1||r^*-\mathrm{Re}\mathcal{R}_2|}\right|^{\bar{\epsilon}_j}\chi_{\mathcal{S}_{r_1,r_2}}\\
&\ \times \Big\langle\frac{|r_1\bar{\mathcal{A}}_j+r_2\bar{\mathcal{B}}_j||\mathcal{A}_j||r^*-\mathrm{Re}\mathcal{R}_1||r^*-\mathrm{Re}\mathcal{R}_2|}{\max\{|\mathrm{Re}\mathcal{R}_1|,|\mathrm{Re}\mathcal{R}_2|,|\mathrm{Re}\mathcal{R}_1+\mathrm{Re}\mathcal{R}_2|\}^2+1}\Big\rangle^{-{q_j}(\frac{1}{8}-)+\bar{\epsilon}_j} e^{-\lambda^2(|r_1-r_2|+|r_2|)q_j}.
\end{aligned}
\end{equation}

Due to the cut-off function $\chi_{\mathcal{S}_{r_1,r_2}}$, we have $|r_*-\tilde{r}_j|\ge \epsilon_{r_*}|\tilde{r}_j|$ for $j=1,2,3$, leading to $|r^*-\mathrm{Re}\mathcal{R}_1||r^*-\mathrm{Re}\mathcal{R}_2|\ge \epsilon_{r_*}^2|\mathrm{Re}\mathcal{R}_1||\mathrm{Re}\mathcal{R}_2|,$ which implies
\begin{equation}\begin{aligned}\label{eq:crossingraph17:1:2:7:a}
\mathcal{U}_j
\ \lesssim	& \   \int_{\mathbb{R}^2}\mathrm{d}r_1 \mathrm{d}r_2\left|\frac{\max\{|\mathrm{Re}\mathcal{R}_1|,|\mathrm{Re}\mathcal{R}_2|,|\mathrm{Re}\mathcal{R}_1+\mathrm{Re}\mathcal{R}_2|\}^2+1}{|r_1 {{A}}_j'+r_2{{B}}'_j||\mathcal{A}_j|\epsilon_{r_*}^2|\mathrm{Re}\mathcal{R}_1||\mathrm{Re}\mathcal{R}_2|}\right|^{\bar{\epsilon}_j}\\
&\ \times\Big\langle\frac{|r_1\bar{\mathcal{A}}_j+r_2\bar{\mathcal{B}}_j||\mathcal{A}_j|\epsilon_{r_*}^2|\mathrm{Re}\mathcal{R}_1||\mathrm{Re}\mathcal{R}_2|}{\max\{|\mathrm{Re}\mathcal{R}_1|,|\mathrm{Re}\mathcal{R}_2|,|\mathrm{Re}\mathcal{R}_1+\mathrm{Re}\mathcal{R}_2|\}^2+1}\Big\rangle^{-{q_j}(\frac{1}{8}-)+\bar{\epsilon}_j} e^{-\lambda^2(|r_1-r_2|+|r_2|)q_j}.
\end{aligned}
\end{equation}
We now bound
\begin{equation}\begin{aligned}\label{eq:crossingraph17:1:3:4}
e^{-\lambda^2q_{j}(|r_1-r_2|+|r_2|)} \ \le \ e^{-\lambda^2q_{j}(|r_1|+|r_2|)/2} \ \le \ & e^{-\lambda^2q_{j}(|r_1|+|B_{j}'' r_2|)/2} \ \le \ e^{-\lambda^2q_{j}(|r_1+B_{j}'' r_2|)/2},
\end{aligned}
\end{equation}
and write

\begin{equation}\begin{aligned}\label{eq:crossingraph17:1:3:3}
\Big|r_1A_{j}'+r_2B_{j}'\Big|^{-\bar{\epsilon}_j}\ = \ &  |A_{j}'|^{-\bar{\epsilon}_j} \Big|r_1+r_2B_{j}''\Big|^{-\bar{\epsilon}_j},
\end{aligned}
\end{equation}
which implies
\begin{equation}\begin{aligned}\label{eq:crossingraph17:1:3:5}
\mathcal{U}_j \ \lesssim	&\ |A_{j}'\mathcal{A}_j^*|^{-\bar{\epsilon}_j} \int_{\mathbb{R} }\mathrm{d}r_1\int_{\mathbb{R} }\mathrm{d}r_2e^{-\lambda^2q_{j}|r_1+B_{j}'' r_2|/2} \Big\langle\Big|\bar{\mathcal{A}}_jr_1+\bar{\mathcal{B}}_jr_2\Big|\mathcal{A}_j^*\Big\rangle^{-{q_j}(\frac{1}{8}-)+\bar{\epsilon}_j}\Big|r_1+r_2B_{j}''\Big|^{-\bar{\epsilon}_j}\chi_{\mathcal{S}_{r_1,r_2}},
\end{aligned}
\end{equation}
with $$\mathcal{A}_j^*= \frac{ |\mathcal{A}_j|\epsilon_{r_*}^2|\mathrm{Re}\mathcal{R}_1||\mathrm{Re}\mathcal{R}_2|}{\max\{|\mathrm{Re}\mathcal{R}_1|,|\mathrm{Re}\mathcal{R}_2|,|\mathrm{Re}\mathcal{R}_1+\mathrm{Re}\mathcal{R}_2|\}^2+1}.$$
We estimate, due to the cut-off functions $\sqrt{\sqrt{\tilde\Psi}(\tilde{\varkappa}_3^j-\tilde{\varkappa}_1^j)}$ and $\sqrt{\sqrt{\tilde\Psi}(\tilde{\varkappa}_2^j-\tilde{\varkappa}_1^j)}$
\begin{equation}\begin{aligned}\label{eq:crossingraph17:1:3:5:a}
|\mathcal{A}_j^*A_j'|\ = \ &  |\mathcal{A}_jA_j'|\frac{  \epsilon_{r_*}^2|\mathrm{Re}\mathcal{R}_1||\mathrm{Re}\mathcal{R}_2|}{\max\{|\mathrm{Re}\mathcal{R}_1|,|\mathrm{Re}\mathcal{R}_2|,|\mathrm{Re}\mathcal{R}_1+\mathrm{Re}\mathcal{R}_2|\}^2+1}\\
\ \gtrsim \ & |\mathcal{A}_jA_j'|\epsilon_{r_*}^2|\mathrm{Re}\mathcal{R}_1||\mathrm{Re}\mathcal{R}_2||\ln\lambda|^{\mathfrak{C}_{\aleph_1^j,\aleph_2^j}^4} \ \gtrsim \ |\ln\lambda|^{\mathfrak{C}_{\aleph_1^j,\aleph_2^j}^5},
\end{aligned}
\end{equation}
for some constants $\mathfrak{C}_{\aleph_1^j,\aleph_2^j}^4,\mathfrak{C}_{\aleph_1^j,\aleph_2^j}^5>0$,
Now, by the change of variables $r_1+r_2B_{j}''\to r_1$, we find
\begin{equation}\begin{aligned}\label{eq:crossingraph17:1:3:6}
\mathcal{U}_j \ \lesssim	&\    |A_{j}'\mathcal{A}_j^*|^{-\bar{\epsilon}_j} \int_{\mathbb{R} }\mathrm{d}r_1\int_{\mathbb{R} }\mathrm{d}r_2e^{-\lambda^2q_{j}|r_1|/2}\chi_{\mathcal{S}_{r_1,r_2}}\\
&\ \ \ \times\Big\langle\Big|(r_1-r_2B_{j}'')\bar{\mathcal{A}}_j+r_2\bar{\mathcal{B}}_j\Big|\mathcal{A}_j^*\Big\rangle^{-{q_j}(\frac{1}{8}-)+\bar{\epsilon}_j} \Big|r_1\Big|^{-\bar{\epsilon}_j}
\\
\lesssim	&\ |A_{j}'\mathcal{A}_j^*|^{-\bar{\epsilon}_j}|\mathcal{A}_j^*(\bar{\mathcal{B}}_j-B_{j}''\bar{\mathcal{A}}_j)|^{-1}  \lambda^{-2+2\bar{\epsilon}_j}\ \lesssim	\ |A_{j}'\mathcal{A}_j^*|^{-1-\bar{\epsilon}_j}|A_{j}'\bar{\mathcal{B}}_j-A_{j}'B_{j}''\bar{\mathcal{A}}_j|^{-1}  \lambda^{-2+2\bar{\epsilon}_j},
\end{aligned}
\end{equation}
when $\bar{\mathcal{B}}_{j}-B_{j}''\bar{\mathcal{A}}_{j}\ne 0$ and ${q_j}(\frac{1}{8}-)-\bar{\epsilon}_j>1$. We write
\begin{equation}
\begin{aligned}\label{eq:crossingraph17:1:3:7}
&\bar{\mathbb{A}}_j \ = \ |A'_{j}||\bar{\mathcal{B}}_{j}-B_{j}''\bar{\mathcal{A}}_{j}| \ =  \ |A'_{j}\bar{\mathcal{B}}_{j}-B_{j}'\bar{\mathcal{A}}_{j}| \ =\   |A'_{j}B^o_{j}-B_{j}'A^o_{j}|
\\
\ =\ & |\alpha_3\cos(\tilde{\varkappa}_3^1-\tilde{\varkappa}_1^1){\sin(2(\tilde{\varkappa}_3^{j}-\tilde{\varkappa}_1^{j}))}[\alpha_1 +\alpha_2\cos(\tilde{\varkappa}_2^1-\tilde{\varkappa}_1^1)\cos(2(\tilde{\varkappa}_2^{j}-\tilde{\varkappa}_1^{j}))]\\
& - [\alpha_4+\alpha_3\cos(\tilde{\varkappa}_3^1-\tilde{\varkappa}_1^1){\cos(2(\tilde{\varkappa}_3^{j}-\tilde{\varkappa}_1^{j}))}] \alpha_2\cos(\tilde{\varkappa}_2^1-\tilde{\varkappa}_1^1)\sin(2(\tilde{\varkappa}_2^j-\tilde{\varkappa}_1^{j}))|\\
\ =\ & |\alpha_3\alpha_1\cos(\tilde{\varkappa}_3^1-\tilde{\varkappa}_1^1){\sin(2(\tilde{\varkappa}_3^{j}-\tilde{\varkappa}_1^{j}))} - \alpha_2\alpha_4\cos(\tilde{\varkappa}_2^1-\tilde{\varkappa}_1^1)\sin(2(\tilde{\varkappa}_2^j-\tilde{\varkappa}_1^{j}))\\
&+ \alpha_3\alpha_2\cos(\tilde{\varkappa}_3^1-\tilde{\varkappa}_1^1){\sin(2(\tilde{\varkappa}_3^{j}-\tilde{\varkappa}_1^{j}))} \cos(\tilde{\varkappa}_2^1-\tilde{\varkappa}_1^1)\cos(2(\tilde{\varkappa}_2^{j}-\tilde{\varkappa}_1^{j}))\\
& -\alpha_3\alpha_2 \cos(\tilde{\varkappa}_3^1-\tilde{\varkappa}_1^1){\cos(2(\tilde{\varkappa}_3^{j}-\tilde{\varkappa}_1^{j}))}\cos(\tilde{\varkappa}_2^1-\tilde{\varkappa}_1^1)\sin(2(\tilde{\varkappa}_2^j-\tilde{\varkappa}_1^{j}))|, 
\end{aligned}
\end{equation}
yielding 	
\begin{equation}\begin{aligned}\label{eq:crossingraph17:1:3}
\mathcal{U}_j 
\
\lesssim	&\ |A_{j}'\mathcal{A}_j^*|^{-\bar{\epsilon}_j}|\mathcal{A}_j^*(\bar{\mathcal{B}}_j-B_{j}''\bar{\mathcal{A}}_j)|^{-1}  \lambda^{-2+2\bar{\epsilon}_j}\ \lesssim	\ \mathfrak{C}_{\aleph_1^j,\aleph_2^j}^5|\ln\lambda|^{\mathfrak{C}_{\aleph_1^j,\aleph_2^j}^6} |A_{j}'|^{-\bar{\epsilon}_j} \bar{\mathbb{A}}_j^{-1}  \lambda^{-2+2\bar{\epsilon}_j}\\ \lesssim	&\ \mathfrak{C}_{\aleph_1^j,\aleph_2^j}^5|\ln\lambda|^{\mathfrak{C}_{\aleph_1^j,\aleph_2^j}^7}  \bar{\mathbb{A}}_j^{-1}  \lambda^{-2+2\bar{\epsilon}_j},
\end{aligned}
\end{equation} for some constants
$ \mathfrak{C}_{\aleph_1^j,\aleph_2^j}^5,\mathfrak{C}_{\aleph_1^j,\aleph_2^j}^6,\mathfrak{C}_{\aleph_1^j,\aleph_2^j}^7>0$. A similar computation can be done in the case $|A_{j}'|<|B_{j}'|,$ leading to
\begin{equation}\begin{aligned}\label{eq:crossingraph17:1:3:0}
\mathcal{U}_j 
\lesssim	&\ \mathfrak{C}_{\aleph_1^j,\aleph_2^j}^5|\ln\lambda|^{\mathfrak{C}_{\aleph_1^j,\aleph_2^j}^7} \bar{\mathbb{A}}_j^{-1}  \lambda^{-2+2\bar{\epsilon}_j}.
\end{aligned}
\end{equation}

{\it	
Next, we consider the case when $|r^*|>4\max\{|\mathrm{Re}\mathcal{R}_1|,$ $|\mathrm{Re}\mathcal{R}_2|,$ $ |\mathrm{Re}\mathcal{R}_1+\mathrm{Re}\mathcal{R}_2|,1\}.$} Different from the previous case, the bound of $\frac{(r^*)^2+1}{|r_1 \bar{\mathcal{A}}_j+r_2\bar{\mathcal{B}}_j||\mathcal{A}_j(r^*)^2+\mathcal{B}_j r^*+\mathcal{C}_j|+{(r^*)^2+1}}$ can be obtained in a more straightforward manner, without using  $\chi_{\mathcal{S}_{r_1,r_2}}$. To see this, we develop 
\begin{equation}\begin{aligned}\label{eq:crossingraph17:1:2:2}
&  \frac{1+(r^*)^2}{|\mathcal{A}_j(r^*)^2+\mathcal{B}_j r^*+\mathcal{C}_j|} \ \le \ \frac{1+(r^*)^2}{|\mathcal{A}_j||(r^*-\mathrm{Re}\mathcal{R}_1)(r^*-\mathrm{Re}\mathcal{R}_2)|} \\
\ = \		& \frac{1}{|\mathcal{A}_j|}\frac{1}{\frac{|(r^*-\mathrm{Re}\mathcal{R}_1)(r^*-\mathrm{Re}\mathcal{R}_2)|}{1+(r^*)^2}} \ = \ \frac{1}{|\mathcal{A}_j|}\frac{1}{\Big|1+\frac{-(\mathrm{Re}\mathcal{R}_1+\mathrm{Re}\mathcal{R}_2)r^*+\mathrm{Re}\mathcal{R}_1\mathrm{Re}\mathcal{R}_2-1}{1+(r^*)^2}\Big|}. 
\end{aligned}
\end{equation}
As $|r^*|>4\max\{|\mathrm{Re}\mathcal{R}_1|,|\mathrm{Re}\mathcal{R}_2|,|\mathrm{Re}\mathcal{R}_1+\mathrm{Re}\mathcal{R}_2|,1\},$ we deduce $\Big|\frac{(\mathrm{Re}\mathcal{R}_1+\mathrm{Re}\mathcal{R}_2)r^*}{1+(r^*)^2}\Big|<\frac14$ and $\Big|\frac{\mathrm{Re}\mathcal{R}_1\mathrm{Re}\mathcal{R}_2-1}{1+(r^*)^2}\Big|<\frac12$. Thus $\Big|1+\frac{-(\mathrm{Re}\mathcal{R}_1+\mathrm{Re}\mathcal{R}_2)r^*+\mathrm{Re}\mathcal{R}_1\mathrm{Re}\mathcal{R}_2-1}{1+(r^*)^2}\Big|>1-\frac12-\frac14=\frac14.$ Therefore, we can bound
$\frac{1+(r^*)^2}{|\mathcal{A}_j(r^*)^2+\mathcal{B}_j r^*+\mathcal{C}_j|} \ \le \ \frac{4}{|\mathcal{A}_j|},$
which implies
\begin{equation}\begin{aligned}\label{eq:crossingraph17:1:2:4}
&  \frac{(r^*)^2+1}{|r_1 \bar{\mathcal{A}}_j+r_2\bar{\mathcal{B}}_j||\mathcal{A}_j(r^*)^2+\mathcal{B}_j r^*+\mathcal{C}_j|+{(r^*)^2+1}}\\
\le \ & \frac{(r^*)^2+1}{|r_1 \bar{\mathcal{A}}_j+r_2\bar{\mathcal{B}}_j||\mathcal{A}_j(r^*)^2+\mathcal{B}_j r^*+\mathcal{C}_j|}\ 
\le \   \frac{4}{|r_1 \bar{\mathcal{A}}_j+r_2\bar{\mathcal{B}}_j||\mathcal{A}_j| }.
\end{aligned}
\end{equation}
We thus bound
\begin{equation}\begin{aligned}\label{eq:crossingraph17:1:2:6:a}
\mathcal{U}_j
\ \lesssim	& \ \int_{\mathbb{R}^2}\mathrm{d}r_1 \mathrm{d}r_2\left|\frac{(r^*)^2+1}{|r_1 \bar{\mathcal{A}}_j+r_2\bar{\mathcal{B}}_j||\mathcal{A}_j(r^*)^2+\mathcal{B}_j r^*+\mathcal{C}_j|+{(r^*)^2+1}}\right|^{\bar{\epsilon}_j}\chi_{\mathcal{S}_{r_1,r_2}}\\
&\ \times \Big\langle\frac{|r_1 \bar{\mathcal{A}}_j+r_2\bar{\mathcal{B}}_j||\mathcal{A}_j(r^*)^2+\mathcal{B}_j r^*+\mathcal{C}_j|}{(r^*)^2+1}\Big\rangle^{-{q_j}(\frac{1}{8}-)+\bar{\epsilon}_j} e^{-\lambda^2(|r_1-r_2|+|r_2|)q_j}\\
\ \lesssim	& \ \int_{\mathbb{R}^2}\mathrm{d}r_1 \mathrm{d}r_2\left|\frac{4}{|r_1 \bar{\mathcal{A}}_j+r_2\bar{\mathcal{B}}_j||\mathcal{A}_j| }\right|^{\bar{\epsilon}_j}\chi_{\mathcal{S}_{r_1,r_2}} \Big\langle\frac{|r_1 \bar{\mathcal{A}}_j+r_2\bar{\mathcal{B}}_j||\mathcal{A}_j|}{4}\Big\rangle^{-{q_j}(\frac{1}{8}-)+\bar{\epsilon}_j} e^{-\lambda^2(|r_1-r_2|+|r_2|)q_j}.
\end{aligned}
\end{equation}
The same argument as above also gives \eqref{eq:crossingraph17:1:3:0}.
We observe that the estimates on $\mathcal{V}_j$ are  the same as  the estimates on $\mathcal{U}_j$, except that $r_1 \bar{\mathcal{A}}_j+r_2\bar{\mathcal{B}}_j$ is now replaced by $r_1 \bar{\mathcal{C}}_j+r_2\bar{\mathcal{D}}_j$. Carrying out the same computations as above gives
\begin{equation}\begin{aligned}\label{eq:crossingraph17:1:3:aa}
\mathcal{V}_j 
\
\lesssim 	&\ \mathfrak{C}_{\aleph_1^j,\aleph_2^j}^5|\ln\lambda|^{\mathfrak{C}_{\aleph_1^j,\aleph_2^j}^7}  \bar{\mathbb{C}}_j^{-1}  \lambda^{-2+2\bar{\epsilon}_j},
\end{aligned}
\end{equation} 
in which	$\bar{\mathcal{C}}_j=C_{j}^o+{\bf i}C_{j}'$, $\bar{\mathcal{D}}_j=D_{j}^o+{\bf i}D_{j}'$ and
\begin{equation}
\begin{aligned}\label{eq:crossingraph17:1:3:7:a}
&\bar{\mathbb{C}}_j \ = \    |C'_{j}\bar{\mathcal{D}}_{j}-D_{j}'\bar{\mathcal{C}}_{j}| \ =\   |C'_{j}D^o_{j}-D_{j}'C^o_{j}|
\\
\ =\ & |\alpha_3\sin(\tilde{\varkappa}_3^1-\tilde{\varkappa}_1^1){\cos(2(\tilde{\varkappa}_3^{j}-\tilde{\varkappa}_1^{j}))}\alpha_2\sin(\tilde{\varkappa}_2^1-\tilde{\varkappa}_1^1)\sin(2(\tilde{\varkappa}_2^{j}-\tilde{\varkappa}_1^{j}))\\
& - \alpha_3\sin(\tilde{\varkappa}_3^1-\tilde{\varkappa}_1^1){\sin(2(\tilde{\varkappa}_3^{j}-\tilde{\varkappa}_1^{j}))} \alpha_2\sin(\tilde{\varkappa}_2^1-\tilde{\varkappa}_1^1)\cos(2(\tilde{\varkappa}_2^j-\tilde{\varkappa}_1^{j}))|	\\
\ =\ & |\sin(\tilde{\varkappa}_3^1-\tilde{\varkappa}_1^1){\cos(2(\tilde{\varkappa}_3^{j}-\tilde{\varkappa}_1^{j}))}\sin(\tilde{\varkappa}_2^1-\tilde{\varkappa}_1^1)\sin(2(\tilde{\varkappa}_2^{j}-\tilde{\varkappa}_1^{j}))\\
& - \sin(\tilde{\varkappa}_3^1-\tilde{\varkappa}_1^1){\sin(2(\tilde{\varkappa}_3^{j}-\tilde{\varkappa}_1^{j}))} \sin(\tilde{\varkappa}_2^1-\tilde{\varkappa}_1^1)\cos(2(\tilde{\varkappa}_2^j-\tilde{\varkappa}_1^{j}))|\\
\ =\ & |\sin(\tilde{\varkappa}_3^1-\tilde{\varkappa}_1^1)\sin(\tilde{\varkappa}_2^1-\tilde{\varkappa}_1^1)\sin(2(\tilde{\varkappa}_2^{j}-\tilde{\varkappa}_3^{j}))|.
\end{aligned}
\end{equation}
Combining the above two estimates \eqref{eq:crossingraph17:1:3:0} and \eqref{eq:crossingraph17:1:3:aa} on $\mathcal{U}_j$ and $\mathcal{V}_j$, we obtain
\begin{equation}\begin{aligned}\label{eq:crossingraph17:1:2:1:c:1}
&\ \Big|\int_{\mathbb{R}^2}\mathrm{d}r_1 \mathrm{d}r_2e^{-\lambda^2(|r_1-r_2|+|r_2|)q_j}\chi_{\mathcal{S}_{r_1,r_2}}\Big\langle\Big|r_1 \bar{\mathcal{A}}_j+r_2\bar{\mathcal{B}}_j\Big|\Big[\frac{r_1^2+r_2^2}{|\mathcal{A}_jr_1^2+\mathcal{B}_j r_1r_2+\mathcal{C}_jr_2^2|}+1\Big]^{-1}\Big\rangle^{-{q_j}(\frac{1}{8}-)}\Big|^{\frac{1}{q_j}}\\
&\ + \ \Big|\int_{\mathbb{R}^2}\mathrm{d}r_1 \mathrm{d}r_2e^{-\lambda^2(|r_1-r_2|+|r_2|)q_j}\chi_{\mathcal{S}_{r_1,r_2}} \Big\langle\Big|r_1 \bar{\mathcal{C}}_j+r_2\bar{\mathcal{D}}_j\Big|\Big[\frac{r_1^2+r_2^2}{|\mathcal{A}_jr_1^2+\mathcal{B}_j r_1r_2+\mathcal{C}_jr_2^2|}+1\Big]^{-1}\Big\rangle^{-{q_j}(\frac{1}{8}-)}\Big|^{\frac{1}{q_j}}\\
\lesssim\ 	& \mathfrak{C}_{\aleph_1^j,\aleph_2^j}^8|\ln\lambda|^{\mathfrak{C}_{\aleph_1^j,\aleph_2^j}^9}  [\bar{\mathbb{C}}_j^{-1/q_j}+ \bar{\mathbb{A}}_j^{-1/q_j}] \lambda^{-2/q_j+2\bar{\epsilon}_j/q_j},
\end{aligned}
\end{equation}
for some constants $ \mathfrak{C}_{\aleph_1^j,\aleph_2^j}^8,\mathfrak{C}_{\aleph_1^j,\aleph_2^j}^9>0$.

Plugging \eqref{eq:crossingraph17:1:2:1:c:1} into \eqref{eq:crossingraph17:1:1}, we find

\begin{equation}
\begin{aligned}\label{eq:crossingraph17:1:3:1:1}
&\ \Big\{\int_{\mathbb{R}^2}\mathrm{d}r_1\mathrm{d}r_2\chi_{\mathcal{S}_{r_1,r_2}}\Big\{e^{-\lambda^2|r_1-r_2|}e^{-\lambda^2|r_2|} \prod_{j=2}^{d}\Big\langle\Big|r_1 \alpha_4+r_2\alpha_1+r_2\alpha_2\cos(\tilde{\varkappa}_2^1-\tilde{\varkappa}_1^1)e^{{\bf i}2(\tilde{\varkappa}_2^j-\tilde{\varkappa}_1^j)}\\
&\ \   +r_1\alpha_3\cos(\tilde{\varkappa}_3^1-\tilde{\varkappa}_1^1)e^{{\bf i}2(\tilde{\varkappa}_3^j-\tilde{\varkappa}_1^j)}\Big||1-\cos[(\aleph^j
_{1}-\aleph^j
_{2})]|^{\frac{1}{2}}\Big\rangle\Big\rangle^{-(\frac{1}{8}-)}\\
&\  + \ \Big\langle\Big|{r_2}\alpha_2\sin(\tilde{\varkappa}_2^1-\tilde{\varkappa}_1^1)e^{{\bf i}2(\tilde{\varkappa}_2^j-\tilde{\varkappa}_1^j)}+{r_1}\alpha_3\sin(\tilde{\varkappa}_3^1-\tilde{\varkappa}_1^1)e^{{\bf i}2(\tilde{\varkappa}_3^j-\tilde{\varkappa}_1^j)}\Big| |1-\cos[(\aleph^j
_{1}-\aleph^j
_{2})]|^{\frac{1}{2}}\Big\rangle^{-(\frac{1}{8}-)}\Big\}^\frac{q}{2}\Big\}^\frac{1}{q}\\
\lesssim	&\   \mathfrak{C}_{\aleph}^{1}|\ln\lambda|^{\mathfrak{C}_{\aleph}^{2}}\prod_{j=2}^{d} [\bar{\mathbb{C}}_j^{-\frac{1}{2q_j}}+ \bar{\mathbb{A}}_j^{-\frac{1}{2q_j}}] \lambda^{-1/q+\bar{\epsilon}_j/q}\
\lesssim	\   \mathfrak{C}_{\aleph}^{1}\lambda^{-2/q+2\bar{\epsilon}/q}|\ln\lambda|^{\mathfrak{C}_{\aleph}^{2}}\prod_{j=2}^{d} \Big[\bar{\mathbb{C}}_j^{-\frac{1}{2q_j}}+ \bar{\mathbb{A}}_j^{-\frac{1}{2q_j}}\Big],
\end{aligned}
\end{equation}
for some constants $ \mathfrak{C}_{\aleph }^{1},\mathfrak{C}_{\aleph}^{2}>0$ and we choose $\bar{\epsilon}_j=\bar{\epsilon}>0$. We therefore find the final estimate on $\mathbb{S}_A$

\begin{equation}\begin{aligned}\label{eq:crossingraph17:1:3:1}
\mathbb{S}_A\lesssim\	&\sqrt{\tilde{F}_1(\tilde{\varkappa}_2-\tilde{\varkappa}_1)} \sqrt{\tilde{F}_1(\tilde{\varkappa}_3-\tilde{\varkappa}_1)}\mathfrak{C}_{\aleph}^{1}\lambda^{-2/q+2\bar{\epsilon}/q}|\ln\lambda|^{\mathfrak{C}_{\aleph}^{2}}\prod_{j=2}^{d} \Big[\bar{\mathbb{C}}_j^{-\frac{1}{2q_j}}+ \bar{\mathbb{A}}_j^{-\frac{1}{2q_j}}\Big].
\end{aligned}
\end{equation}
Next, we will estimate $\mathbb{S}_B$.
Under the effect of $\chi_{\mathcal{S}_{r_1,r_2}'}$, the quantity $\frac{1}{|1-|\cos(\aleph_1^j-\aleph_2^j)||^\frac12} $ can be bounded from above by a fixed constant by \eqref{Lemm:Angle:E11}. We can simply bound

\begin{equation}
\begin{aligned}\label{eq:crossingraph17:1A:2:1}
\mathbb{S}_B\lesssim\	&  \Big\{\int_{\mathbb{R}^2}\mathrm{d}r_1\mathrm{d}r_2\Big\{e^{-\lambda^2|r_1-r_2|}e^{-\lambda^2|r_2|}\chi_{\mathcal{S}_{r_1,r_2}'} \prod_{j=2}^{d}\Big\langle\Big\{\Big|r_1 \alpha_4+r_2\alpha_1+r_2\alpha_2\cos(\tilde{\varkappa}_2^1-\tilde{\varkappa}_1^1)e^{{\bf i}2(\tilde{\varkappa}_2^j-\tilde{\varkappa}_1^j)}\\
&\ \ +r_1\alpha_3\cos(\tilde{\varkappa}_3^1-\tilde{\varkappa}_1^1)e^{{\bf i}2(\tilde{\varkappa}_3^j-\tilde{\varkappa}_1^j)}\Big|+\Big|{r_2}\alpha_2\sin(\tilde{\varkappa}_2^1-\tilde{\varkappa}_1^1)e^{{\bf i}2(\tilde{\varkappa}_2^j-\tilde{\varkappa}_1^j)}\\
&\ \ +{r_1}\alpha_3\sin(\tilde{\varkappa}_3^1-\tilde{\varkappa}_1^1)e^{{\bf i}2(\tilde{\varkappa}_3^j-\tilde{\varkappa}_1^j)}\Big|\Big\} \Big\rangle^{-(\frac{1}{8}-)}\Big\}^\frac{q}{2}\Big\}^\frac{1}{q} \sqrt{\sqrt{\tilde\Psi}(\tilde{\varkappa}_2-\tilde{\varkappa}_1)}\\
\lesssim\	&  \Big\{\int_{\mathbb{R}^2}\mathrm{d}r_1\mathrm{d}r_2\Big\{e^{-\lambda^2|r_1-r_2|}e^{-\lambda^2|r_2|}\chi_{\mathcal{S}_{r_1,r_2}'}\Big\{\prod_{j=2}^{d}\Big\langle\Big|r_1 \alpha_4+r_2\alpha_1+r_2\alpha_2\cos(\tilde{\varkappa}_2^1-\tilde{\varkappa}_1^1)e^{{\bf i}2(\tilde{\varkappa}_2^j-\tilde{\varkappa}_1^j)}\\
&\ \ +r_1\alpha_3\cos(\tilde{\varkappa}_3^1-\tilde{\varkappa}_1^1)e^{{\bf i}2(\tilde{\varkappa}_3^j-\tilde{\varkappa}_1^j)}\Big|\Big\rangle^{-(\frac{1}{8}-)}+\prod_{j=2}^{d}\Big\langle\Big|{r_2}\alpha_2\sin(\tilde{\varkappa}_2^1-\tilde{\varkappa}_1^1)e^{{\bf i}2(\tilde{\varkappa}_2^j-\tilde{\varkappa}_1^j)}\\
&\ \ +{r_1}\alpha_3\sin(\tilde{\varkappa}_3^1-\tilde{\varkappa}_1^1)e^{{\bf i}2(\tilde{\varkappa}_3^j-\tilde{\varkappa}_1^j)}\Big| \Big\rangle^{-(\frac{1}{8}-)}\Big\}\Big\}\Big\}^\frac{q}{2}\Big\}^\frac{1}{q} \sqrt{\sqrt{\tilde\Psi}(\tilde{\varkappa}_2-\tilde{\varkappa}_1)}.
\end{aligned}
\end{equation}
The computations used to estimate $\mathbb{S}_A$ can be repeated, leading to 
\begin{equation}\begin{aligned}\label{eq:crossingraph17:1:3:2}
\mathbb{S}_B\lesssim\	&\sqrt{\tilde{F}(\tilde{\varkappa}_2-\tilde{\varkappa}_1)} \sqrt{\tilde{F}(\tilde{\varkappa}_3-\tilde{\varkappa}_1)}\mathfrak{C}_{\aleph}^{1}\lambda^{-2/q+2\bar{\epsilon}/q}|\ln\lambda|^{\mathfrak{C}_{\aleph}^{2}}\prod_{j=2}^{d} \Big[\bar{\mathbb{C}}_j^{-\frac{1}{2q_j}}+ \bar{\mathbb{A}}_j^{-\frac{1}{2q_j}}\Big].
\end{aligned}
\end{equation}
\end{proof}

\section{Feynman diagrams}
\label{Duhamel}

\begin{figure}
	\centering
	\includegraphics[width=.99\linewidth]{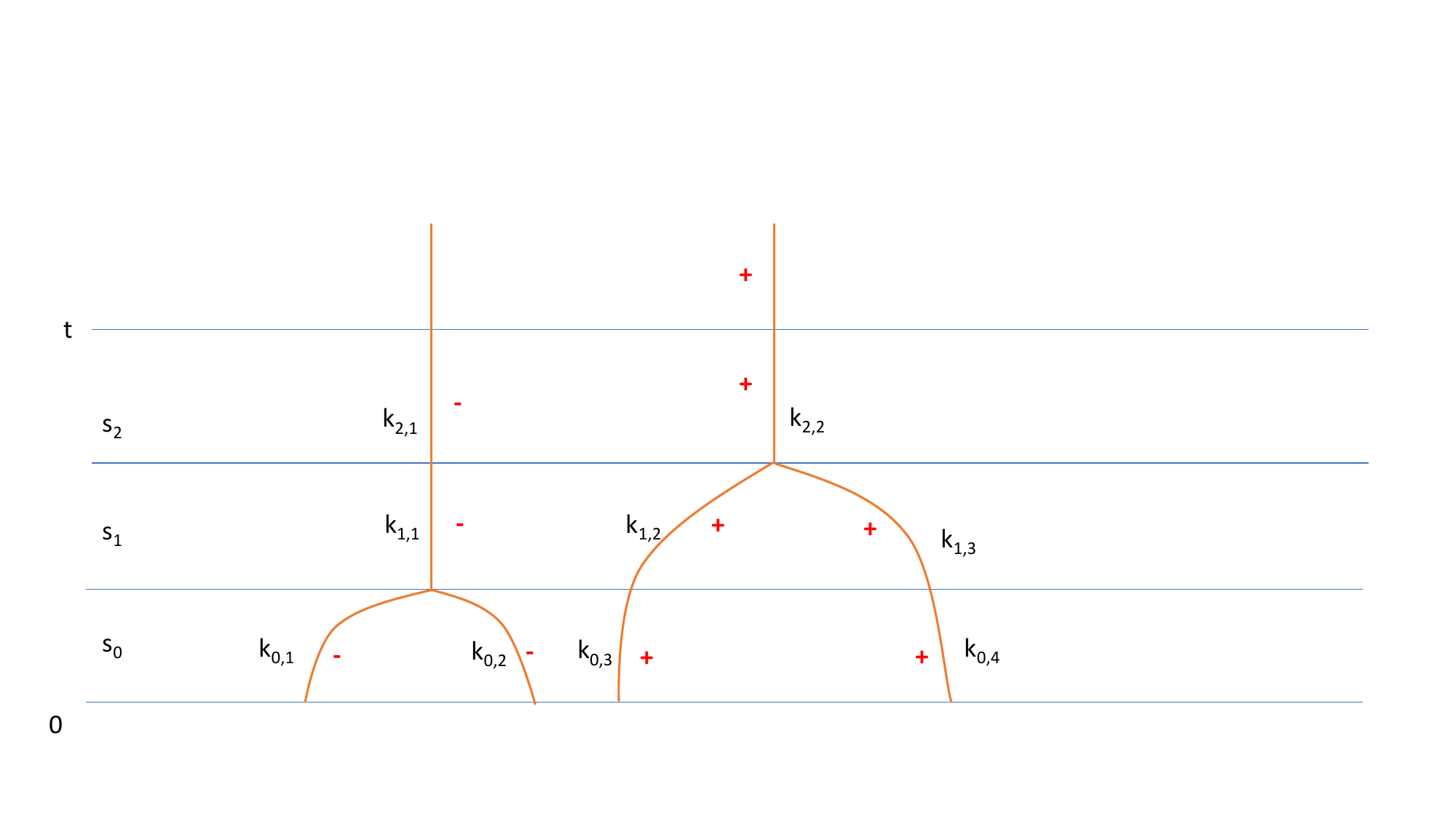}
	\caption{An example of a Feynman diagram. At time slice $s_0$, the edges are $k_{0,1},k_{0,2},k_{0,3},k_{0,4}$, with the signs $-,-,+,+$. At time slice $s_1$, the edges are $k_{1,1},k_{1,2},k_{1,3}$, with the signs $-,+,+$.  At time slice $s_2$, the edges are $k_{2,1},k_{2,2}$, with the signs $-,+$. We have $-k_{1,1}+k_{0,1}+k_{0,2}=0$ and $k_{2,2}-k_{1,2}-k_{1,3}=0$.}
	\label{Fig1}
\end{figure}

In this section, we will discuss the construction of our Feynman diagrams.  Since the time interval $[0,t]$ is divided into $n+1$ time slices, each of length $s_i$, we represent the time slices  from the bottom to the top of the diagram, with the lengths $s_0$, $s_1$, $\dots$, $s_n$, as shown in Figure \ref{Fig1}. As discussed above, the strategy used to get \eqref{StartPointAverageAlphaDuhamel} is repeatedly applied, but only to the term containing $\Phi_{1,i}$, to obtain the full Duhamel expansions. 
Therefore, at each time slice $s_i$, only one Duhamel expansion is allowed, with the term containing $\Phi_{1,i}$.  In this Duhamel expansion, the delta function associated to $\Phi_{1,i}$   means that we combine the two momenta $k_2,k_3$ into the momentum $k_1$.  To represent this delta function on the diagram, at time slice $s_i$, we draw a combination of one couple of the segments of time slice $s_{i-1}$ into one segment of time slice $s_i$. It is straightforward that the number of segments at time slice $s_{i}$ is $2+n-i$, indexing by $k_{i,1},\cdots k_{i,2+n-i}$. We denote those segments of time slice $s_{i-1}$  by $k_{i-1,\rho_i},k_{i-1,\rho_i+1}$ and the one at time slice $s_i$ is denoted by $k_{i,\rho_i}$. The index $\rho_i$ indicates the position where the combination happens. In other words, if  $i$ is the index of time slice $s_i$, the the combination happens at the segment  $\rho_i$.   In this process,  the delta function associated to $\Phi_{1,i}$  implies the identity  $\sigma_{i,\rho_i}k_{i,\rho_i}+\sigma_{i-1,\rho_i}k_{i-1,\rho_i}+\sigma_{i-1,\rho_i+1}k_{i-1,\rho_i+1}=0.$  As we notice from the Duhamel expansions,   there is a ``sign'' parameter $\sigma$  associated to each momentum $k$. The signs for $k_{i-1,\rho_i},k_{i-1,\rho_i+1}$  are denoted by $\sigma_{i-1,\rho_i},\sigma_{i-1,\rho_i+1}$.  Among the $2+(n-i)$ segments of the time slice  $s_{i}$, since  the delta function only enforces the combination at the segment $k_{i,\rho_i}$,  where the 2 segments $\sigma_{i-1,\rho_i},\sigma_{i-1,\rho_i+1}$ merge together, we introduce the following way of indexing the segments in two consecutive time slices $s_{i-1}$ and $s_i$: $k_{i,l}=k_{i-1,l}, \sigma_{i,l}=\sigma_{i-1,l}, $ for  $l\in\{1,\cdots,\rho_i-1\}$ and $k_{i,l}=k_{i-1,l+1}, \sigma_{i,l}=\sigma_{i-1,l+1}, $ for  $l\in\{\rho_i+1,\cdots,2+{n-i}\}$.  

\subsection{Diagrams of $\mathcal{G}^0_n$, $\mathcal{G}^1_n$, $\mathcal{G}^2_n$, $\mathcal{G}^3_n$.}\label{Sec:FirstDiagrams} Below we construct the diagrams corresponding to the expressions   $\mathcal{G}^0_n$, $\mathcal{G}^1_n$, $\mathcal{G}^2_n$, $\mathcal{G}^3_n$, that we introduced  at the end of Section \ref{Sec:Duha}.  We follow the standard definitions in graph theory (cf. \cite{berge2001theory,trees2003translated,tutte2019connectivity,tutte1984graph}).

\begin{figure}
	\centering
	\includegraphics[width=.49\linewidth]{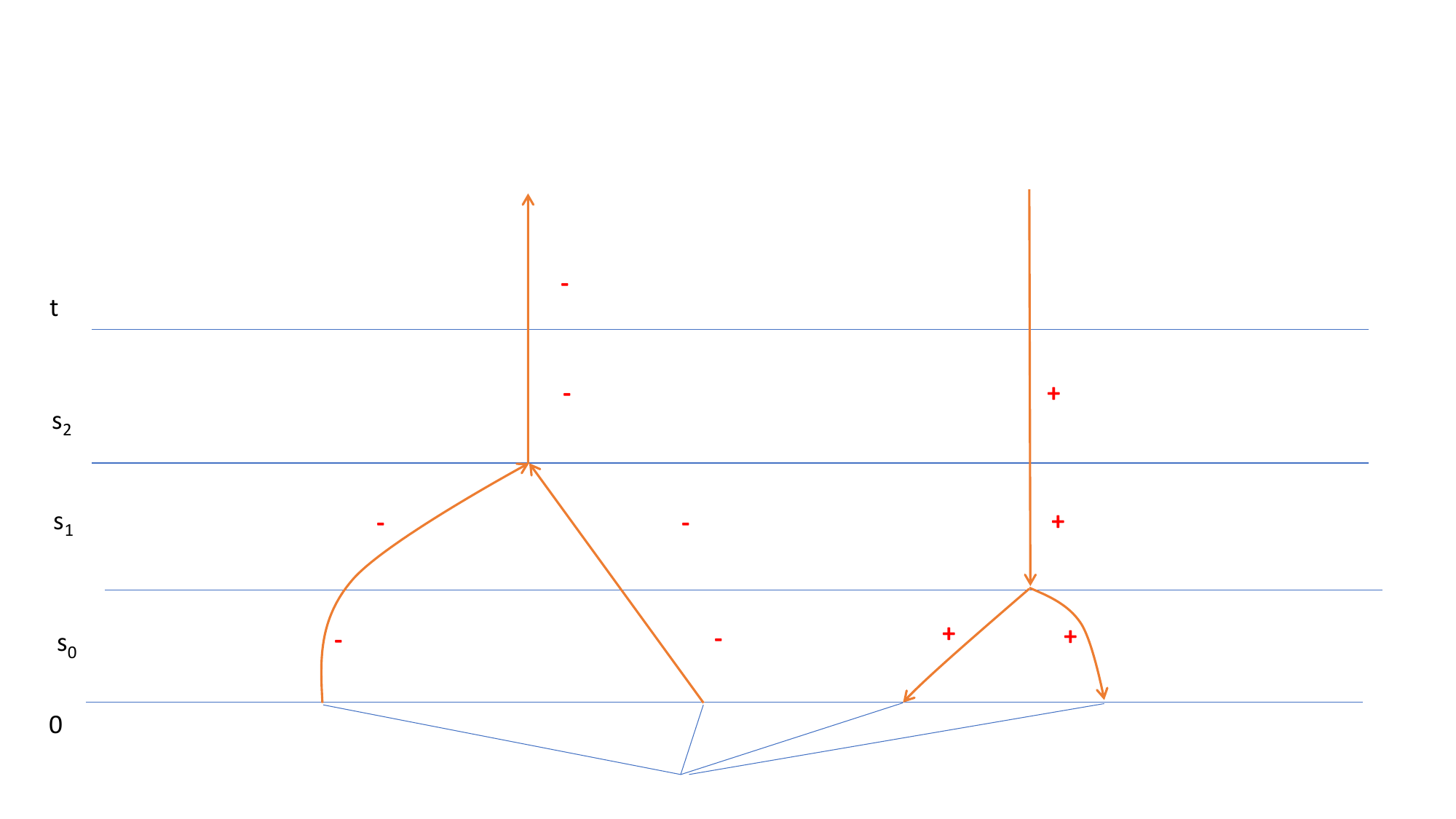}
	\caption{An example of a Feynman diagram with clusters. One new cluster vertex is added at the bottom of the diagram. The arrows correspond to the {\it first assigned orientation.}}
	\label{Fig2}
\end{figure}
\begin{itemize}
	\item[(a)] Even we will be mainly focus on pairing graphs, we also consider the general case that the bottom of the graph has a partition $S$, and for each element $ \{k_{0,j}\sigma_{0,j}\}_{j\in A}$ of the partition $S$, we have the summation $\sum_{j\in A} k_{0,j}\sigma_{0,j}=0$. We call this the cluster partitioning. This is represented on the diagram by using an extra ``cluster vertex'' at the bottom of the graph, for each $A$. This cluster vertex connects itself to all of the vertices $ k_{0,j}$, $j\in A$, appearing in the delta function $\delta\left(\sum_{j\in A} k_{0,j}\sigma_{0,j}\right)$ (see Figure \ref{Fig2}).  

	\item[(b)] In our diagrams, the two identities   $\mathbf{1}(\sigma_{n,2}=1)$ and $\mathbf{1}(\sigma_{n,1}=-1)$ can be repeated by  additionally  assigning $\sigma_{n,1}=-1$ and $\sigma_{n,2}=1$ for the signs of the two topmost segments. This means that the sign of the topmost left segment is always negative and the sign of the topmost right segment is always positive. We assign the value of $k_{n,1}$, the momentum of the topmost segment on the left, to be the same as $k_{n,2}$, the momentum of the topmost segment on the right.  
	\item[(c)]  The remaining signs  can be assigned to the diagrams  using the same rule discussed above, from top to bottom, with the notice that 
	the cluster vertices do not affect those signs. 
	\item[(d)] The total phase of the diagram will not be represented on the diagram, but it could be written in a short form as $\prod_{i=0}^n e^{-\mathrm{i} s_i \vartheta_i}$.
	In Lemma \ref{Lemma:ReallVatheta}, we represent a way to explicitly compute the real part of $\vartheta_i$.

	\item[(e)] As we discuss above, each delta function is associated to one vertex, in which the combination happens. Let us recall that the Duhamel expansion  is continuously applied, but only to the terms containing $\Phi_{1,i}$, to obtain the full Duhamel expansions.  Therefore, the delta functions are associated with the factors $-{\bf i}\lambda \Phi_{1,i}$.  
	\item[(f)] In each diagram, we can see that there are two components, one on the left, associated to the momentum $k_{n,1}$ and one on the right, associated to the momentum $k_{n,2}$. In this diagram, when all cluster vertices and their edges are removed, the two components are disconnected. Since the left component is associated to $k_{n,1}$, whose sign is minus, we call it the ``minus diagram''. The right one is called the ``plus diagram''. The plus and minus diagrams are indeed connected only by the cluster vertices.
	Without loss of generality, we assume that the first splitting always happens on the plus diagram.

\end{itemize} 
\begin{lemma}\label{Lemma:ReallVatheta}
	The real part of $\vartheta_i$ can also be computed as follows \begin{equation}\label{GraphSec:E3}\begin{aligned}
			\mathrm{Re}\vartheta_i\ =\ &\sum_{l=i+1}^n \mathfrak{X}(\sigma_{l,\rho_l},k_{l,\rho_l},\sigma_{l-1,\rho_l},  k_{l-1,\rho_l},\sigma_{l-1,\rho_l+1}, k_{l-1,\rho_l+1})\\
			=\ &\sum_{l=1}^{n-i+2}\sigma_{i,l}\omega(k_{i,l})-2\omega(k_{n,1}).\end{aligned}\end{equation}
\end{lemma}
\begin{proof}
	We will show the claim by an induction argument with respect to $i$. 
	
	For $i=n-1$, we have
	$$\mathrm{Re}\vartheta_{n-1}=\sigma_{n,2}\omega(k_{n,2})+\sigma_{n-1,2}\omega(k_{n-1,2})+\sigma_{n-1,3}\omega(k_{n-1,3}).$$
	On the other hand, 
	$$\sum_{l=1}^{3}\sigma_{n-1,l}\omega(k_{n-1,l})=\sigma_{n-1,1}\omega(k_{n-1,1})+\sigma_{n-1,2}\omega(k_{n-1,2})+\sigma_{n-1,3}\omega(k_{n-1,3}).$$
	Now, since $k_{n,1}=k_{n,2}=k_{n-1,1}$, and  $\sigma_{n,1}=-\sigma_{n,2}=\sigma_{n-1,1}=-1$, we find
	$$\mathrm{Re}\vartheta_{n-1}=\sum_{l=1}^{3}\sigma_{n-1,l}\omega(k_{n-1,l})-2\omega(k_{n,1}).$$
	
	Suppose that the claim is true for $i\ge m$,
	we will show that it is also true for $i=m-1$. Let us compute
	\begin{equation*}\begin{aligned}
			\mathrm{Re}\vartheta_{m-1}\ = \ & \mathrm{Re}\vartheta_{m}+\mathfrak{X}(\sigma_{m,\rho_m},k_{m,\rho_m},\sigma_{m-1,\rho_m},  k_{m-1,\rho_m},\sigma_{m-1,\rho_m+1}, k_{m-1,\rho_m+1})\\
			\ = \ & \sum_{l=1}^{n-m+2}\sigma_{m,l}\omega(k_{m,l})-2\omega(k_{n,1})\\
			& +\mathfrak{X}(\sigma_{m,\rho_m},k_{m,\rho_m},\sigma_{m-1,\rho_m},  k_{m-1,\rho_m},\sigma_{m-1,\rho_m+1}, k_{m-1,\rho_m+1})\\
			\ = \ & \sum_{l=1}^{n-m+2}\sigma_{m,l}\omega(k_{m,l})-2\omega(k_{n,1})\\
			& +\sigma_{m,\rho_m}\omega(k_{m,\rho_m})+\sigma_{m-1,\rho_m}\omega(  k_{m-1,\rho_m})+\sigma_{m-1,\rho_m+1}\omega(k_{m-1,\rho_m+1})\\
			\ = \ & \sum_{l=1}^{n-m+3}\sigma_{m-1,l}\omega(k_{m-1,l})-2\omega(k_{n,1}),\end{aligned}\end{equation*}
	where, in the last equality, we have used the fact that $k_{m,l}=k_{m-1,l}, \sigma_{m,l}=\sigma_{m-1,l}, $ for  $l\in\{1,\cdots,\rho_m-1\}$ and $k_{m,l}=k_{m-1,m+1}, \sigma_{m,l}=\sigma_{m-1,m+1}, $ for  $l\in\{\rho_m+1,\cdots,2+{n-m}\}$.  
	
	Therefore, the identity \eqref{GraphSec:E3} is proved.
\end{proof}
\begin{remark} We anticipate here that the term 
	$-2\omega(k_{n,1})$ in \eqref{GraphSec:E3} will not be integrated, hence the important quantities in $\mathrm{Re}\vartheta_i$ are $\sigma_{i,l}\omega(k_{i,l})$.
\end{remark}
\subsection{How to integrate $\mathcal{G}^0_n$, $\mathcal{G}^1_n$, $\mathcal{G}^2_n$, $\mathcal{G}^3_n$? Construction of integrated graphs.}\label{Sec:MomentumGraph}

In this section we will introduce the way to integrate the diagrams discussed in the previous section, by adding an extra orientation to these diagrams. The orientation allows us to know which edges are integrated first, which edges are integrated next. We call them integrated graphs.
We denote an integrated graph by $\mathfrak{G}=(\mathfrak{V},\mathfrak{E})$, in which $\mathfrak{V}$ and $\mathfrak{E}$ are the sets of vertices and edges.  We first start with the construction of the integrated graph.
\smallskip

{\bf (A) The construction of  integrated graphs}. This graph depends on the following parameters of  $\mathcal{G}^0_n$, $\mathcal{G}^1_n$, $\mathcal{G}^2_n$, $\mathcal{G}^3_n$,  the number of interacting vertices $n$ and the vector  $\rho$ that encodes where the splittings happen. All of these parameters come from the expressions of $\mathcal{G}^0_n$, $\mathcal{G}^1_n$, $\mathcal{G}^2_n$, $\mathcal{G}^3_n$.  Given the parameters $S,n,\rho$, the integrated graph can be reconstructed using the following scheme (see Figure \ref{Fig4}). 
\begin{figure}
	\centering
	\includegraphics[width=.49\linewidth]{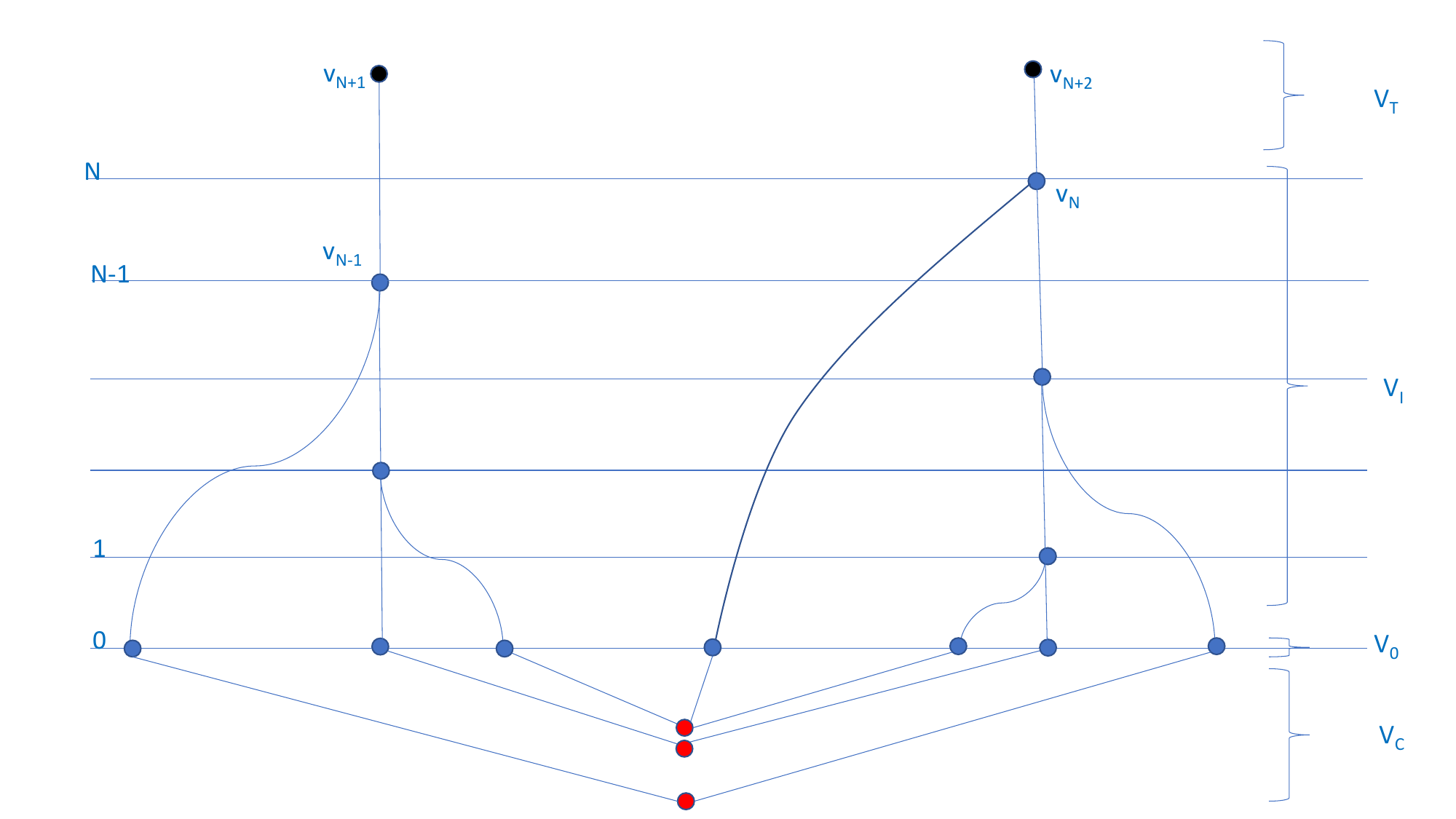}
	\caption{An example of an integrated graph. The sets $\mathfrak{V}_T$, $\mathfrak{V}_I$, $\mathfrak{V}_0$, $\mathfrak{V}_C$ and the vertices $v_{N+1}$, $v_{N+2}$  are marked on the graph.}
	\label{Fig4}
\end{figure}
\begin{itemize}
	\item[(a)] We start with $\mathfrak{G}^0=(\mathfrak{V}^0,\mathfrak{E}^0)$, where $\mathfrak{V}^0$ contains two initial vertices $v_{n+2}$ and $v_{n+1}$. We suppose that $v_{n+1}$ belongs to the minus tree.
	\item[(b)] At the first iteration, we attach one new edge $e_1$ to $v_{n+2}$. This edge belongs to the plus tree and one of its vertices is $v_{n+2}$. We label the other vertex by $v_{n}$. 
	\item[(c)]  In the second iteration, we attach a new edge $e_2$ to $v_{n+1}$. This edge then belongs to the minus tree and one of its vertices is $v_{n+1}$. We label the other vertex of $e_2$ by $v_{n-1}$. 
	\item[(d)] Since we assume that the first combination happens on the plus diagram, the edge $e_1$ will be decomposed into another two edges. In the next step, we attach these two edges to  $v_{n}$, from left to right. 
	\item[(e)] Using $\rho_{n-1}$, we could locate the next  vertex where the combination happens. This vertex could be $v_{n-1}$, or a different vertex denoted by $v_{n-2}$. We keep attaching two new edges to either $v_{n-1}$ or $v_{n-2}$ in the next step.
	\item[(f)] Repeating the above procedure, we can label all of the     vertices associated to the delta functions. There are in total $n$ of them. Those vertices are labeled $v_{n}, v_{n-1}, \cdots v_1$ from the top to the bottom of the diagrams. 
	\item[(g)] The labeled vertices should be in $\mathfrak{V}_H=\{v_1,\cdots,v_{n+2}\}$. The set $\mathfrak{V}_H$ is called ``higher time  vertex set''. The vertices from $v_1$ to $v_n$ are call ``interacting vertices'' and the set of all of them is called ``interacting vertex set'' $$\mathfrak{V}_I=\mathfrak{V}_H\backslash \mathfrak{V}_T=\mathfrak{V}_H\backslash \{v_{n+1},v_{n+2}\},$$
	where  $\mathfrak{V}_T$ denotes the top vertices $\{v_{n+1},v_{n+2}\}$. The set of the other vertices is denoted by $\mathfrak{V}_0$ and called ``zero time vertex set''. For each cluster $A$ in $S$, we call the associated vertex by ``cluster vertex'' and label it by $u_A$. The set $\mathfrak{V}_C=\{u_A\}_{A\in S}$ is called ``cluster vertex set''. Then the vertex set of the final graph $\mathfrak{G}$ is $\mathfrak{V}=\mathfrak{V}_T\cup\mathfrak{V}_H\cup\mathfrak{V}_0\cup\mathfrak{V}_C.$  The edges can be added in a natural way. 
	\item[(h)] Define the function $\mathcal{T}:\mathfrak{V}\to [0,n+2]$, in which $\mathcal{T}(v)=j$ if $v=v_j$ for $j\in\{1,\cdots,n+2\}$,  and $\mathcal{T}(v)=0$ in the other cases. 
	\item[(k)] For a given pair of vertices $(v,v')$, if we can find a set of vertices $\{v_{i_j}\}_{j=1}^m$ such that $v_{i_1}=v$, $v_{i_m}=v'$, and $v_{i_l}$ is connected to $v_{i_{l+1}}$ by an edge, we call $\{v_{i_j}\}_{j=1}^m$ a ``path'' that connects $v$ and $v'$. 
	\item[(l)] For any $v\in\mathfrak{V}$, we denote the set of all edges attached to $v$ by $\mathfrak{E}(v) = \{e\in\mathfrak{E} \  | \ v\in e\}$.  For $v\in\mathfrak{V}_I$, we define two new sets $\mathfrak{E}_+(v)$ and $\mathfrak{E}_-(v)$, that satisfy $\mathfrak{E}(v)=\mathfrak{E}_+(v)\cup \mathfrak{E}_-(v)$. In this definition, $\mathfrak{E}_+(v)=\{e\}$, where $e$ is the first edge attached to $v$ in the construction described  before and $\mathfrak{E}_-(v)=\mathfrak{E}(v)\backslash\{e\}$. 
	
	\item[(m)] 	The cluster $\mathfrak{V}_c^{(j)}$ is defined to be a clustering of the edges intersecting with the time slice $j$. The construction $\mathfrak{V}_c^{(j)}$ is done using an iterative procedure in which the interacting vertices are added to the graph, from the bottom to the top. The additional vertex $v_j$ combines the two edges in $\mathfrak{E}_-(v_j)$ into the new one in $\mathfrak{E}_+(v_j)$. The two edges of $\mathfrak{E}_-(v_j)$ belong to some cluster in the previous iteration $\mathfrak{V}_c^{(j-1)}$. The cluster $\mathfrak{V}_c^{(j)}$ is constructed by joining all clusters of $\mathfrak{V}_c^{(j-1)}$ and replacing the three edges by the one in $\mathfrak{E}_+(v_j)$. The rest remains the same.   
\end{itemize} 
The above process gives us an unoriented graph $(\mathfrak{V},\mathfrak{E}).$

{\bf (B) Embedding the delta functions into the integrated graphs - The first assigned orientation of a diagram.}

For an interacting vertex $v_i$, denote by $k_0,k_1,k_2$ and $\sigma_0,\sigma_1,\sigma_2$ the  momenta associated to $v_i$ and the corresponding signs, we also need the following definition of ``the total phase of $v_i$'', which is an abbreviation of the definition given in \eqref{Def:Omega} 
\begin{equation}
	\label{Def:Xvi}
	\mathfrak{X}_i \ = \ \mathfrak{X}(v_i) \ = \ \sigma_0\omega(k_0)+\sigma_1\omega(k_1)+
	\sigma_2\omega(k_2).
\end{equation}

Moreover, for $v\in \mathfrak{V}_I\cup\mathfrak{V}_C$, 
there is a delta function associated to it
$$
\delta\left(\sum_{e\in\mathfrak{E}_{in}(v)}k_e-\sum_{e\in\mathfrak{E}_{out}(v)}k_e\right),
$$ in which $\mathfrak{E}_{in}(v)$ denotes the set of all the edges $e$ in $\mathfrak{E}(v)$ such that the sign of their momenta $k_e$ is always $+1$, while $\mathfrak{E}_{out}(v)$ denotes the set of all the edges $e$ in $\mathfrak{E}(v)$ such that the sign of their momenta $k_e$ is always $-1$.

{\it The first assigned orientation for a diagram}. For a segment $e$ associated to an interacting vertex $v$ and belonging to $\mathfrak{E}_-(v)$, we assign to this edge the orientation of going inward the vertex $v$ if the sign $\sigma_e$ of the associated  momentum $k_e$ is $+1$, otherwise, we assign to this edge the orientation of going outward the vertex. If  $e$ belongs to $\mathfrak{E}_+(v)$, we assign to this edge the orientation of going inward the vertex $v$ if the sign $\sigma_e$ of the associated  momentum $k_e$ is $-1$, otherwise, we assign to this edge the orientation of going outward the vertex  (see Figure \ref{Fig2}).  We also assign orientations for the edges of the cluster vertices as follows. Suppose that $e$ is a segment associated  to some cluster vertex. One of the vertices of $e$ should belong to $\mathfrak{V}_0$, suppose that it is $v'$. The vertex $v'$ should, again, belong to another segment $e'$, which is connected to another interacting vertex $v''\in\mathfrak{V}_I$. This edge $e'$ has already a orientation, coming from the previous way of assigning the orientations. If the orientation $e'$ is going inward the vertex $v'$, then the orientation of $e$ is going outward the vertex $v'$. If the orientation $e'$ is going outward the vertex $v'$, then the orientation of $e$ is going inward the vertex $v'$ (see Figure \ref{Fig2}).

For any $v\in\mathfrak{V}_I$, let $e$ be an edge in $\mathfrak{E}_-(v)$,  $v'$ be the other vertex of the edge $e$.  Denote by $\sigma_e$ the sign associated to the edge $e$, with respect to the vertex $v$. Then the sign associated to the edge $e$, with respect to the vertex $v'$ is the opposite of $\sigma_e$, which is $-\sigma_e$. We then define a ``sign'' mapping $\sigma_v \ : \ \mathfrak{E}(v) \to \{-1,1\}$ as follow
\begin{equation}\label{GraphSec:E5}
	\begin{aligned}
		\sigma_v (e) \  = & \ \sigma_e,\\
		\sigma_{v'} (e) \  = & \ -\sigma_e.
	\end{aligned}
\end{equation}

\begin{remark}\label{Eplusminus} Though $$\mathfrak{E}_{in}(v)\cup\mathfrak{E}_{out}(v)=\mathfrak{E}_+(v)\cup\mathfrak{E}_-(v)=\mathfrak{E}(v),$$ the pair of sets $\mathfrak{E}_{in}(v)$, $\mathfrak{E}_{out}(v)$ is completely different from the pair of sets $\mathfrak{E}_+(v)$, $\mathfrak{E}_-(v)$. The pair of sets $\mathfrak{E}_{in}(v)$, $\mathfrak{E}_{out}(v)$ is related to the first assigned orientation.
\end{remark}

Our aim is to integrate out all the delta functions. To understand how one can integrate those graphs, we would need the following construction of ``free edges''. 

\smallskip

{\bf (C) The construction of free edges scheme - The second assigned orientation of the diagram.}  Before introducing the scheme, we will need the definition of ``cycles'', following Berge \cite{berge2001theory} and Serre \cite{trees2003translated} (see Figure \ref{Fig20}). 

\begin{definition}[Cycles] Let us consider a graph of $n$ interacting vertices and a set of vertices $\{v_{i}\}_{i\in\mathfrak{I}}$ of $\mathfrak{V}_H\cup\mathfrak{V}_0$. We consider two cases.
	
	{\it Case 1: The set of vertices $\{v_{i}\}_{i\in\mathfrak{I}}$  does not contain both of the two top vertices $v_{n+1}$ and $v_{n+2}$.} If there exists a set of cluster vertices $\{u_A\}_{A\in S'\subset S}$ such that we could go from one vertex $v_{i_1}$ of the set $\{v_{i}\}_{i\in\mathfrak{I}}$  to all of the vertices of $\{v_{i}\}_{i\in\mathfrak{I}} \cup \{u_A\}_{A\in S'\subset S}$ and back to $v_{i_1}$ via the edges  of the graph, we call it a cycle. Suppose that $v_l$ is the top most interacting vertex in the cycle, we say that this is a cycle of the vertex $v_l$.
	
	{\it Case 2: The set of vertices $\{v_{i}\}_{i\in\mathfrak{I}}$   contains both of the two top vertices $v_{n+1}$ and $v_{n+2}$.}
	We ignore $v_{n+1}$ and $v_{n+2}$ and introduce a ``virtual vertex'' $v_*$, that connects both of the two  vertices $v_n$ and $v_{n-1}$.  If there exists a set of cluster vertices $\{u_A\}_{A\in S'\subset S}$ such that we could go from $v_{*}$ to all of the vertices of $\{v_*\}\cup \{v_{i}\}_{i\in\mathfrak{I}}\cup \{u_A\}_{A\in S'\subset S}$ and back to $v_{*}$ via the edges  of the graph, we call it a cycle. Suppose that $v_l$ is the top most  vertex in the cycle, we say that this is a cycle of the vertex $v_l$. Thus, $v_l$ can be either the virtual vertex $v_*$ or an interacting vertex. 
	\begin{figure}
		\centering
		\includegraphics[width=.49\linewidth]{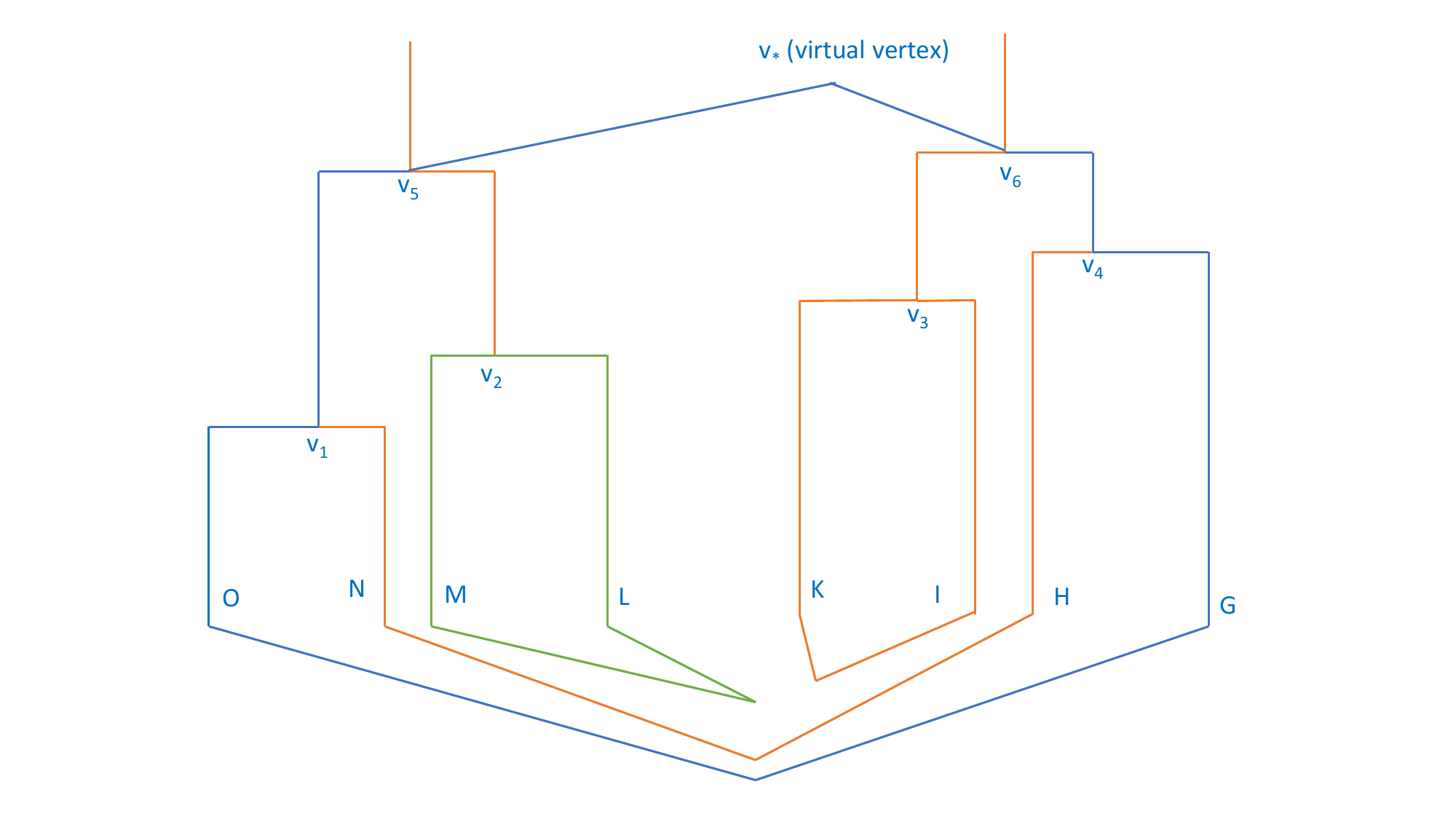}
		\caption{In the picture, $\{v_2,M,L\}$ forms a cycle. Moreover, $v_*$ is the virtual vertex and $\{v_*,v_6,v_4,G,O,v_1,v_5\}$ also forms a cycle.}
		\label{Fig20}
	\end{figure}
\end{definition}
{\it The scheme to obtain free edges.} 
Using the definition of cycles,  we have the following scheme to choose which edges are free and which edges are not. The role of free and non-free edges will be explained in ``(D) Integrating the integrated edges'' below.
\begin{itemize}
	\item[(a)] First, we start with the vertices in $\mathfrak{V}_0$ and the cluster vertices in $\mathfrak{V}_C$ and all of the edges that are connecting them. Those edges are non-free edges or integrated edges. In the following, we will see that we only integrate non-free edges and leave free edges to the end.
	\item[(b)] Now, we continue with the first vertex $v_{1}$. Among the two edges attached to $v_{1}$ that are connected to the vertices in  $\mathfrak{V}_0$, the right edge is set to be a non-free edge, which is then an integrated edge. If the other one forms a cycle with $v_{1}$ and the vertices  in $\mathfrak{V}_C$, we set it free. Otherwise, we set it non-free, and it is again an integrated edge.  
	\item[(c)] We continue the procedure in a recursive manner. At step $i$, we consider the two edges associated to $v_{i}$ that connect $v_{i}$ with the vertices of $\mathfrak{V}_0$ and $\{v_j\}_{j=1}^{i-1}$. If the right one forms a cycle with $\mathfrak{V}_0$ and $\{v_j\}_{j=1}^{i}$, we set it free, otherwise we set it non-free, which then becomes an integrated edge. Then we continue the same procedure with the left one.
	\item[(d)] The procedure is carried on until $i$ reaches $i=n$. 
	\item[(e)] The momenta associated to the free edges are called ``free momenta''. In the sequel, sometimes, we use the edge and its associated momentum for the same roles. We then use the terminologies ``free momenta'' and ``free edges'' for the same purpose. 
	\item[(f)]  	In our construction, it is possible that there is a free edge attached to the virtual vertex $v_*$. If this free edge is on the minus diagram, we call it the ``virtually free edge attached to $v_{n+1}$''. Since $v_{n+1}$ does not play any important role in our construction of free edges, we  also call this free edge the ``virtually free edge attached to $v_{n-1}$''. The associated momentum is called the ``virtually free momentum attached to $v_{n+1}$ (and, equivalently, $v_{n-1}$)''.  If this free edge is on the minus diagram, we call it the ``virtually free edge attached to $v_{n+2}$  (and, equivalently,  $v_{n}$)''.   The associated momentum is called the ``virtually free momentum attached to $v_{n+2}$ (and, equivalently,  $v_{n}$)''. The cycle in this case is the cycle of $v_*$. By our construction, the virtually free edge is always the one on the left, which is  the virtually free edge attached to $v_{n+1}$ (or the virtually free edge attached to $v_{n-1}$, equivalently).
\end{itemize}

{\it The second assigned orientation of the diagram: } The diagram caries a natural orientation, which is defined as follows. For any integrated edge $e=\{v',v\}$, we say that the orientation of the edge is from $v'$ to $v$ if $v$ belongs to the path which does not contain any free edge from $v'$ to the virtual vertex $v_*$. We then define the partial order $v'\succ v$. If $v\succ v'$, we then say $v' 	\prec v$. The graph $\mathfrak{G}$  with the partial order $\succ $ is denoted by $\mathfrak{G}^\succ$. 
We also define the  set 
\begin{equation}\label{Pv}
	\mathfrak{P}(v)=\{v' | v' 	\succ v, v'\ne v_*,v_{n+1},v_{n+2}\}.
\end{equation} 
\begin{remark}\label{SecondAssignedDir} Note that the first assigned orientation is for the graph $\mathfrak{E}$, and is related to the pair of sets $\mathfrak{E}_{in}(v)$, $\mathfrak{E}_{out}(v)$ (see Remark \ref{Eplusminus}), while the second assigned orientation is for the graph $\mathfrak{G}^\succ$ that contains only the integrated edges.
\end{remark}

{\bf (D) Integrating the integrated edges.}
We set $\mathfrak{F}$ to be the set of free edges and $\mathfrak{E}'$ to be the set of integrated edges.  Following the second orientation of the diagram, we could integrate all the delta functions, using the  integrated edges in $\mathfrak{E}'$.  We also need the following definition.

\begin{definition}[Degree of a Vertex]
	For any $v\in\mathfrak{V}_I$, let $\mathfrak{F}(v)$ denote the set of free edges attached to $v$, then $\mathfrak{F}(v)=\mathfrak{E}(v)\cap \mathfrak{F}$. We call the number of free edges in $\mathfrak{F}(v)$ the degree of the interacting vertex and denote it by $\mathrm{deg}v$.
\end{definition}

\begin{figure}
	\centering
	\includegraphics[width=.49\linewidth]{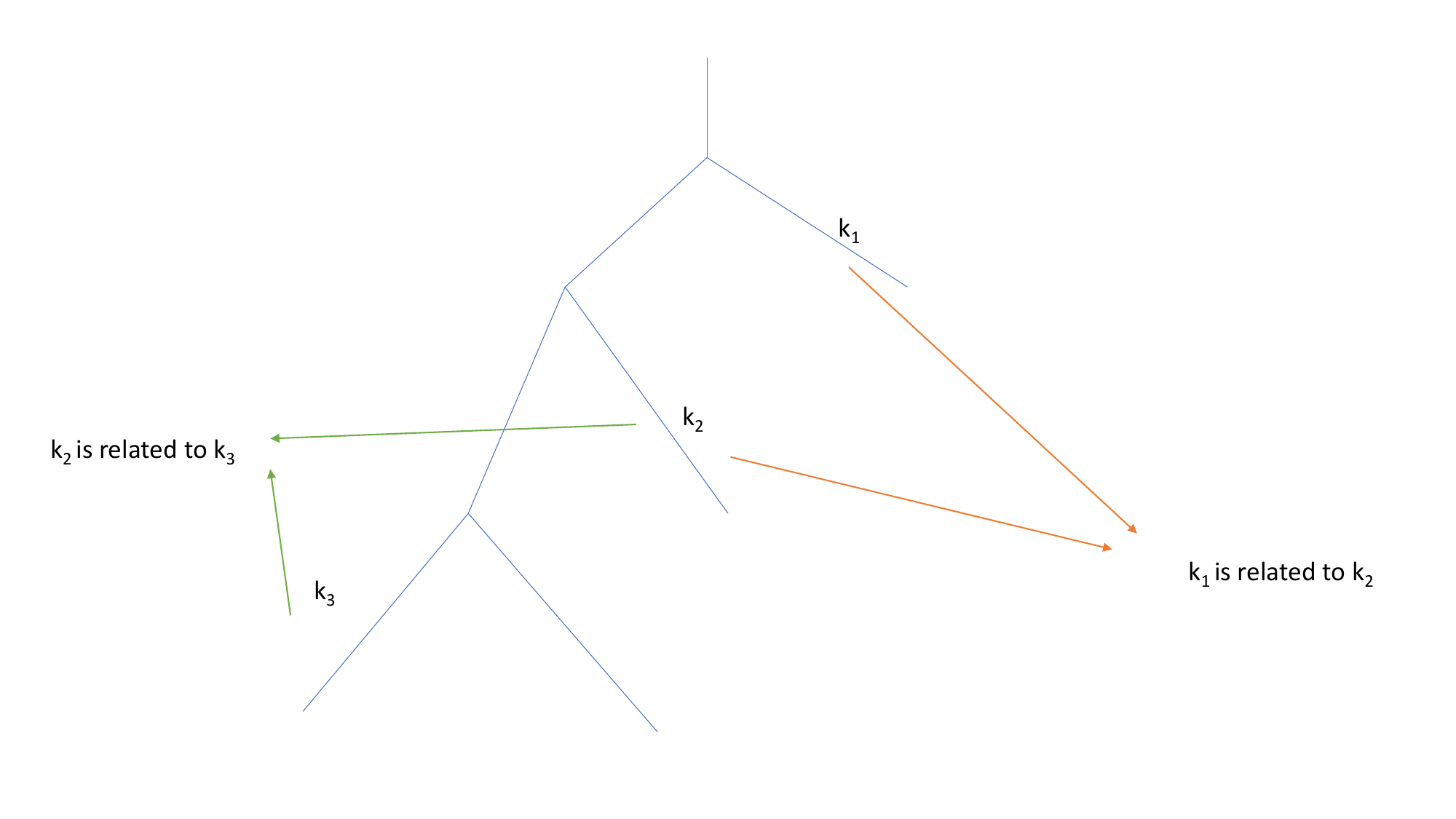}
	\caption{Graphical representation of related momenta. }
	\label{Fig44}
\end{figure}
\subsection{Properties of integrated graphs - Summations of free edges, 1-Separation, 2-Separation}\label{Sec:PropMometum}

Below, we prove properties of the Feynman diagrams constructed above. Many of the properties proved in this subsection are similar by  those of the Feynman diagrams for the nonlinear Schr\"odinger equation studied by Lukkarinen and Spohn \cite{LukkarinenSpohn:WNS:2011}. Several results in this subsection are inspired by  Tutte \cite{tutte2019connectivity,tutte1984graph}. 
\begin{lemma}\label{Lemma:ProductSignFirstDirection}
	If $e=\{v,v'\}\in \mathfrak{E}$ does not intersect $\mathfrak{V}_T$, then $\sigma_v(e)\sigma_{v'}(e)=-1.$
\end{lemma}
\begin{proof}
	Let us consider the orientation of the edge $e$, following the first assigned orientation for the diagram. Without loss of generality, suppose that the first orientation of $e$ is going from $v$ to $v'$, then $\sigma_{v}(e)=1$ and $\sigma_{v'}(e)=-1$, thus $\sigma_v(e)\sigma_{v'}(e)=-1.$ 
\end{proof}

\begin{lemma}\label{Lemma:cycleFreeEdge}
	Let $e=\{v,v'\}$ be a free edge and suppose that $\mathcal{T}(v)>\mathcal{T}(v')$. Let $v''$ be a vertex belonging to the cycle of the vertex $v$, then either $v\succ v''$ or $v'\succ v''$.

	Moreover, let $e'=\{v,v''\}$, $e''=\{v,v'''\}$ be two integrated edges attached to an arbitrary vertex $v$. Then it cannot happen that $v\succ v''$ and $v\succ v'''$. 
\end{lemma}
\begin{proof}
	\begin{figure}
		\centering
		\includegraphics[width=.49\linewidth]{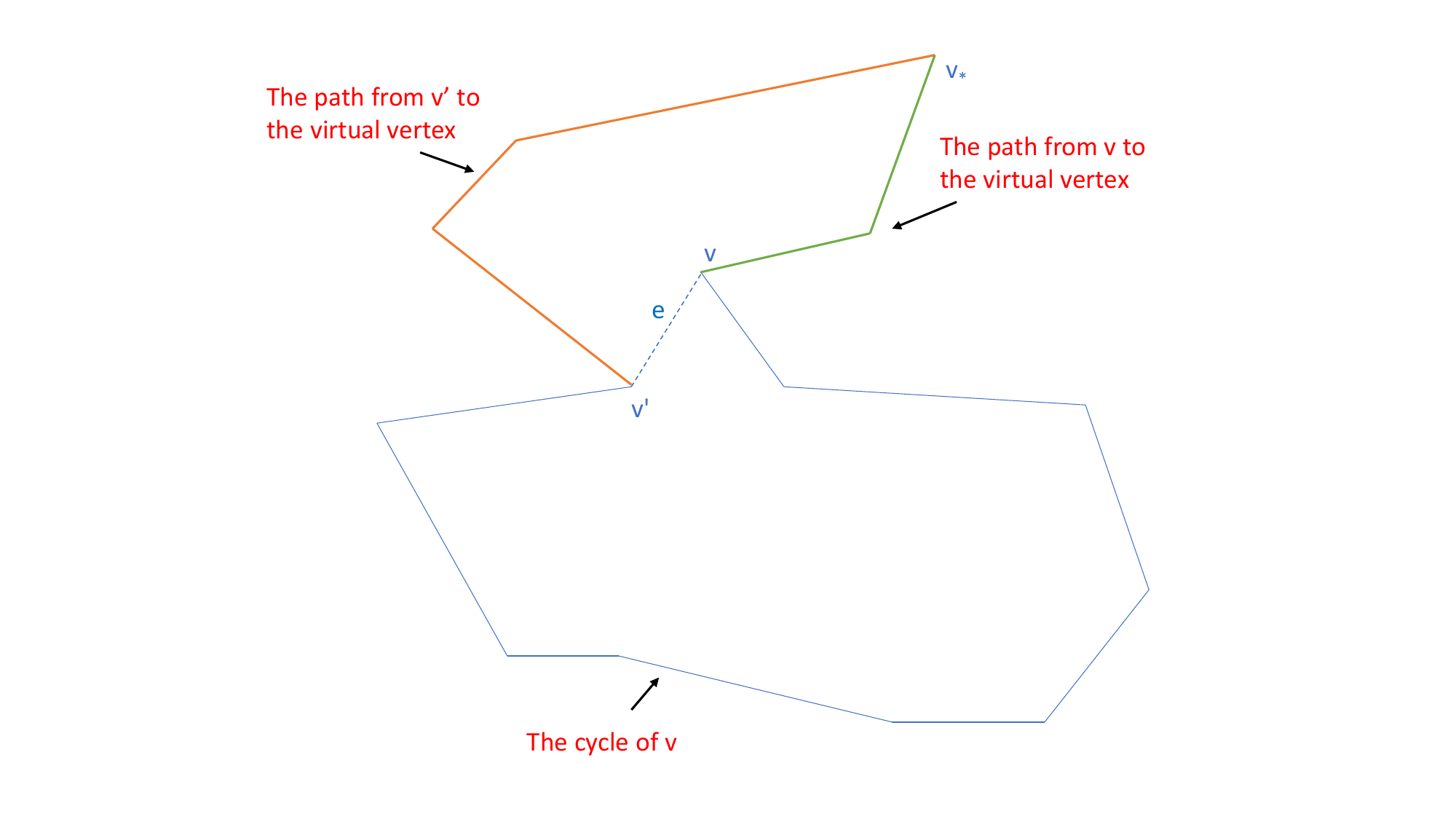}
		\caption{The paths from $v$, $v'$ to the virtual vertex and the cycle of $v$ form another cycle. }
		\label{Fig27}
	\end{figure}
	Let us consider the two paths of $v$ and $v'$ to the virtual vertex $v_*$. Since we can go from $v$ to $v'$ using the integrated edges of the cycle of $v$, the two paths of $v$ and $v'$ have to meet at a vertex $v'''$ in the cycle of $v$; otherwise, together with the cycle of $v$, they will form a full cycle, leading to a contradiction (see Figure \ref{Fig27}) since we already break all cycles using the free edges. As a result, $v\succ v'''$ and $v'\succ v'''$. If $v''\equiv v'''$, then the proof is done. Otherwise, $v''$ belongs to either the path from $v$ to $v'''$ or the path from $v'$ to $v'''$. In the first case, $v\succ v''$ and in the second case, $v'\succ v''$ (see Figure \ref{Fig28}).
	
	Moreover, let $e'=\{v,v''\}$, $e''=\{v,v'''\}$ be  integrated edges attached to an arbitrary vertex $v$ and suppose that $v\succ v''$ and $v\succ v'''$. In this case, it is clear that there is a path from $v$ to $v''$ that ends at the virtual vertex $v_*$. There is also a path from $v$ to $v'''$ and ends at the virtual vertex $v_*$. The two paths form a cycle. This is a contradiction.

\end{proof}

\begin{lemma}\label{lemma:freemomenta}
	A Feynman diagram with $n$ interacting vertices has totally $n+2-|S|$ free momenta. 
\end{lemma}
\begin{proof}
	The number of free edges is indeed the cyclomatic number of the graph (cf. \cite{berge2001theory}). This  is the minimum number of edges that
	must be removed from the graph to break all its cycles,
	making it into a tree or forest. This  concept was first introduced by  Kirchhoff (cf. \cite{hage1996island}). The cyclomatic number  is a function of   the number of edges in the  graph, number of vertices and the number of connected
	components. In our graph, including the virtual vertex and excluding $\{v_{n+1}, v_{n+2}\}$, there are  $3n+4$ edges, $2n+3+|S|$ vertices, and $1$ connected component.    As a result, the number of free edge is computed by the cycle rank formula (cf. \cite{berge2001theory,tutte1984graph}) $$3n+4-(2n+3+|S|)+1=n+2-|S|.$$ 
\end{proof}

\subsubsection{Summations of free edges}
\begin{lemma}\label{Lemma:SumFreeEdges}
	Let $k_e$ be an integrated momentum with $e$ being an edge in $\mathfrak{G}^\succ$, which is the graph $\mathfrak{G}$ associated with the orientation $\succ$ defined above.   Let $v_1$, $v_2$ be the two vertices of $k_e$ and the second orientation of $k_e$ is from $v_1$ to $v_2$, that means $v_1\succ v_2$.   Let  $\mathfrak{P}(v_1)$ be the set defined in \eqref{Pv}.  The following identity then holds true
	\begin{equation}
		\label{Lemma:SumFreeEdges:1}
		k_e \ =  \ \sum_{v\in\mathfrak{P}(v_1)}\sum_{e'\in\mathfrak{F}(v)}\left(-\sigma_{v_1}(e)\sigma_v(e')\right)k_{e'},
	\end{equation}
	in which $k_{e'}$ is the momentum of $e'$.
\end{lemma}
\begin{proof}
	We define $l$ to be the number of elements of the set $\mathfrak{P}(v_1)$ of $v_1$. Our proof will be based on an induction argument with respect to the number $l=|\mathfrak{P}(v_1)|$. Suppose that $l=1$, then $\mathfrak{P}(v_1)=
	\{v_1\}$. We deduce that  $\mathcal{F}(v_1)=\mathfrak{E}(v_1)\backslash\{e\}$. Suppose the contrary, that there is an edge $e'$ in $\mathfrak{E}(v_1)$ such that $e'$ is not free. Let $v_3$ is the other vertex of $e'$, different from $v_1$. If $v_3\succ v_1$, then $v_3\in \mathfrak{P}(v_1)$, that contradicts the fact that $\mathfrak{P}(v_1)=
	\{v_1\}$. Hence, since $e'$ is not free, $v_1\succ v_3$. However, in this case,  we have both $v_1\succ v_3$ and $v_1\succ v_2$, contradicting Lemma \ref{Lemma:cycleFreeEdge}. As a result, there is one  delta function associated to $v_1$, and $\mathcal{F}(v_1)=\mathfrak{E}(v_1)\backslash\{e\}$, leading to
	$$\sum_{e'\in \mathfrak{E}(v_1)} \sigma_{v_1}(e')k_{e'} \ = \ 0.$$
	We then deduce that
	$$k_e \ = \ -\sigma_{v_1}(e)\sum_{e'\in\mathfrak{E}({v_1})\backslash\{e\}}\sigma_{v_1}(e')k_{e'} \ = \ \sum_{e'\in\mathfrak{E}({v_1})\backslash\{e\}}(-\sigma_{v_1}(e)\sigma_{v_1}(e'))k_{e'},$$
	which means \eqref{Lemma:SumFreeEdges:1} holds for $l=1$. Suppose by induction that \eqref{Lemma:SumFreeEdges:1} holds for any $|\mathfrak{P}(v_1)|$ up to $l\ge 1$, we consider an edge $e=\{v_1,v_2\}$ with $|\mathfrak{P}(v_1)|=l+1$ and $v_1\succ v_2$. Now, since $v_1\notin\mathfrak{V}_T$, it follows that 
	$$k_e \ = \ \sum_{e'\in\mathfrak{E}(v_1)\backslash\{e\}}(-\sigma_{v_1}(e)\sigma_{v_1}(e'))k_{e'}.$$
	In this sum, for any edge $e'=\{v_e,v_1\}$,  in $\mathfrak{E}({v_1})\backslash\{e\}$, there are two possibilities. If $e'$ is free, we do not have to worry about it. If $e'$ is not free, we also have two cases, either $v_e\succ v_1$ or $v_1\succ v_e$.  
	\begin{itemize}
		\item If $v_e\succ v_1$, then the set $\mathfrak{P}(v_e)$ contains at most $l$ elements since $v_1$ is on the path from $v_e$ to the virtual vertex $v_*$.  We can use the induction hypothesis applied to  $e'$.
		\item If $v_1\succ v_e$, since we also have $v_1\succ v_2$, this leads to a contradiction with the conclusion of Lemma \ref{Lemma:cycleFreeEdge}.
	\end{itemize}
	Thus,  the identity \eqref{Lemma:SumFreeEdges:1} also holds for $l+1$. This completes the induction proof. 
\end{proof}
\begin{definition}[Dependence of Edges]\label{Def:Dependence}
	In the  formula \eqref{Lemma:SumFreeEdges:1}, the edge $e$ is said to ``depend'' on the edges $e'$ and $k_e$ is said to ``depend'' on $k_{e'}$. If the edge $e$ does not depend on $e'$, we say that $e$ is ``independent'' of $e'$.
\end{definition}

\begin{lemma}\label{Lemma:SumFreeEdges2}
	For any integrated edge $e=(v_1,v_2)\in\mathfrak{E}'$, $v_1\succ v_2$, the following identities hold true
	\begin{equation}
		\begin{aligned}\label{Lemma:SumFreeEdges2:1}
			k_e \ = & \ \sum_{v\in\mathfrak{P}(v_1)}\sum_{e'=(v,v_{e'})\in \mathfrak{F}(v)}\mathbf{1}_{v_{e'}\notin\mathfrak{P}(v_1)}\left(-\sigma_{v_1}(e)\sigma_v(e') \right)k_{e'}\\
			\ = & \ -\sigma_{v_1}(e)\sum_{e'\in\mathfrak{F}}\mathbf{1}(\exists v\in e' \cap \mathfrak{P}(v_1) \mbox{ and } e' \cap \mathfrak{P}(v_1)^c \neq \emptyset)\sigma_v(e')k_{e'},
	\end{aligned}\end{equation}
	in which $k_{e},k_{e'}$ are the momenta of $e$, $e'$. 
	Moreover, if $e'=(v,v')$, such that $e'\neq e$, $v\in\mathfrak{P}(v_1)$ and $v'\notin \mathfrak{P}(v_1)$, then $e'$ is free. 
	
	For any edge $e=(v_1,v_2)\in\mathfrak{E}$, if $e$ is integrated,  we suppose $v_1\succ v_2$ and define  
	\begin{equation}
		\label{Lemma:SumFreeEdges2:0}
		\mathfrak{F}_e \ = \ \{e'\in\mathfrak{F}|\exists v\in e' \cap \mathfrak{P}(v_1) \mbox{ and } e' \cap \mathfrak{P}(v_1)^c \neq \emptyset)\},
	\end{equation}
	then formula \eqref{Lemma:SumFreeEdges2:1} can be expressed under the following form, in which  $\sigma_{e,e'}\in\{\pm 1\}$,
	\begin{equation}
		\label{Lemma:SumFreeEdges2:2}
		k_e \ = \ \sum_{e'\in \mathfrak{F}_e}\sigma_{e,e'}k_{e'}. 
	\end{equation}
	We also denote 
	\begin{equation}
		\label{Lemma:SumFreeEdges2:0:1}
		\mathfrak{F}_{k_e} \ = \ \{k_{e'}\ | \ e'\in\mathfrak{F}, \ \exists v\in e' \cap \mathfrak{P}(v_1) \mbox{ and } e' \cap \mathfrak{P}(v_1)^c \neq \emptyset)\}.
	\end{equation}
	If $e$ is free, we set $\mathfrak{F}_e=\{e\}$ and $\sigma_{e,e}=1$. We have the following version of formula \eqref{Lemma:SumFreeEdges2:2}
	\begin{equation}
		\label{Lemma:SumFreeEdges2:3}
		k_e \ = \ \sum_{e\in \mathfrak{F}_e}\sigma_{e,e}k_{e}. 
	\end{equation}
\end{lemma}
\begin{proof}
	We only prove the first identity, since the second one follows in a straightforward manner.
	According to \eqref{Lemma:SumFreeEdges}
	\begin{equation*}
		k_e \ =  \ \sum_{v\in\mathfrak{P}(v_1)}\sum_{e'\in\mathfrak{F}(v)}\left(-\sigma_{v_1}(e)\sigma_v(e')\right)k_{e'}.
	\end{equation*} 
	Suppose $v,v'\in\mathfrak{P}(v_1)$, such that $e'=(v,v')$ is free. It follows that $-\sigma_{v_1}(e)\sigma_{v}(e')-\sigma_{v_1}(e)\sigma_{v'}(e')=0$. As a result, these pair of edges cancel each other. Summing over edges in $\mathfrak{F}(v)$ and $\mathfrak{F}(v')$ gives the desired identity.
	
	For any edge $e'=(v,v')$, such that $e'\neq e$, $v\in\mathfrak{P}(v_1)$ and $v'\notin \mathfrak{P}(v_1)$, we prove that $e'$ is free using a proof by contradiction. Suppose that $e'$ is integrated, then either $v\succ v'$ or $v' \succ v$. 
	\begin{itemize}
		\item Case 1: $v'\succ v$. By the definition that $v\in \mathfrak{P}(v_1)$,  we have $v\succ v_1$, thus $v'\succ v\succ v_1$ and $v'\in \mathfrak{P}(v_1)$. This contradicts the assumption that $v'\notin \mathfrak{P}(v_1)$.
		\item Case 2: $v\succ v'$. Since $v\in\mathfrak{P}(v_1)$, there exists a vertex $v''\in\mathfrak{P}(v_1)$ such that $v\succ v''$ and the edge $\{v,v''\}$ belongs to the path from $v$ to the virtual vertex $v_*$. At the vertex $v$, we have $v\succ v''$ and $v\succ v'$. This contradicts the conclusion of Lemma \ref{Lemma:cycleFreeEdge}.
	\end{itemize}
	Therefore $e'$ is free. The two formulas \eqref{Lemma:SumFreeEdges2:2} and \eqref{Lemma:SumFreeEdges2:3} then follow in a straightforward manner. 
\end{proof}

\begin{lemma}\label{Lemma:DepedenceOnFreeEdgeInsideAcycle}
	Let $v$ be any interacting vertex. Then $\mathrm{deg}(v)\in\{0,1\}$. 
	If $v\in\mathfrak{V}_I$ and $\mathrm{deg}(v) = 1$, then $\mathfrak{E}_-(v)=\{e,e'\}$, where $e$ is a free edge. Denote by $k_e$ and $k_{e'}$ the momenta associated to $e$ and $e'$. We have $k_{e'} = \pm k_{e} + q$, where  $q$ is independent of $k_e$. Moreover, the  dependence on $k_e$ of the other integrated edges are give below.
	
	\begin{itemize}
		\item[(a)] Consider the integrated edge $e'=(v_1,v_2)\in\mathfrak{E}'$ and denote by $k_{e'}$ its momentum. Then one of the following $3$ possibilities should happen
		$k_{e'} = \pm k_e + \tilde{k}_{e'}$,  or
		$k_{e'}  =   \tilde{k}_{e'}$,
		where $ \tilde{k}_{e'} $ is independent of $k_e$ in all 3 cases. 
		\item[(b)] If $v_1\in\mathfrak{V}_I$, and suppose $e',e''\in\mathfrak{E}(v_1)$, $e'\neq e''$, and denote by $k_{e'},k_{e''}$ their momenta,  then  one of the following  $5$ possibilities should happen
		$k_{e'} + k_{e''}   =   \pm k_e + \tilde{k}_{e',e''}$, or $k_{e'} + k_{e''}  =  \pm 2k_e\ + \tilde{k}_{e',e''},$ or $k_{e'} + k_{e''}  =  \tilde{k}_{e',e''},$
		where $ \tilde{k}_{e',e''}$ is independent of $k_e$ and  $k_{e'}$ in all 5 cases.  
	\end{itemize}
	
\end{lemma}
\begin{proof}
	\begin{figure}
		\centering
		\includegraphics[width=.49\linewidth]{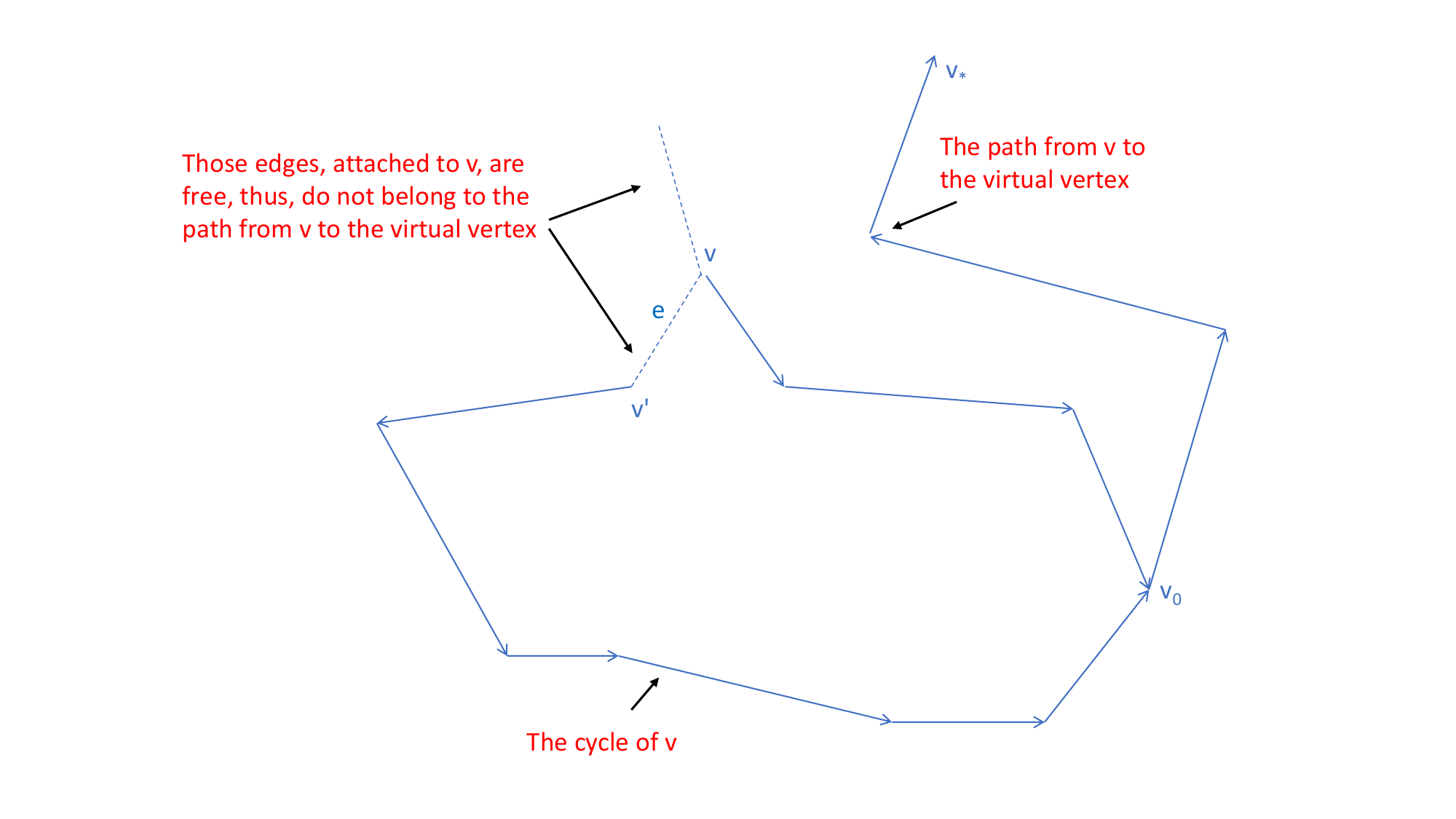}
		\caption{The paths from $v$, $v'$ coincide at $v_0$, a vertex in the cycle of $v$. }
		\label{Fig28}
	\end{figure} The fact that the degree of $v$ has to be smaller than or equal to $1$ is straightforward. We now suppose that $\mathrm{deg}v=1$. Among the integrated edges in $\mathfrak{E}_-(v)$, we denote $e=\{v,v'\}$ to be  free edge.  Since $e$ is a free edge, then according to the proof Lemma \ref{Lemma:cycleFreeEdge}, there is a vertex $v_0$, which belongs to the cycle of $v$, such that the two paths, containing $v_0$, from $v$ and $v'$ to the virtual vertex coincide.   (see Figure \ref{Fig28}). We can suppose that $e'=\{v,v_2\}$ is the first edge in the path from $v$ to $v_0$, then $e'\in \mathfrak{E}_-(v)$ and $k_{e'}=\pm  k_e+q$  by Lemma \ref{Lemma:SumFreeEdges}.  This finishes the proof of the first part of the  lemma. 
	
	We now prove the claim (a). Let us suppose that  $v_1\succ v_2$ and ${e}=\{v,w\}$,  where $w$ is the other vertex, different from $v$. If $v\succ v_1$ or $w\succ v_1$, by \eqref{Lemma:SumFreeEdges:1}, we have $k_{e'}=\pm k_e+ \tilde{k}_{e'}$, where $\tilde{k}_{e'}$ is independent of $k_e$. If there is no oriented path from either $v$ or $w$ to $v_1$ to the virtual vertex $v_*$, then $e\notin \mathfrak{F}_{e'}$ and thus  $k_{e'}=\tilde{k}_{e'}$, where $\tilde{k}_{e'}$ is independent of $k_e$.

	Finally, we prove claim (b). We consider the following cases.
	\begin{itemize}
		\item If both $e'',e'$ are free, we then choose $\tilde{k}_{e'',e'''}$ to be $k_{e'}+k_{e''}$, which is independent of $e$.
		\item If one of them is free, say $e'$. By claim (a), $k_{e''}=\pm k_e+ \tilde{k}_{e''}$, where $\tilde{k}_{e''}$ is independent of $k_e$. Thus, we can choose $\tilde{k}_{e',e''}=k_{e'}+ \tilde{k}_{e''}$.
		\item If both of them are integrated, $e'=\{v,w_1\}$ and $e''=\{v,w_2\}$, with $w_1\neq w_2$. We denote the third edge attach to $v$ by $e'''=\{v,w_3\}$. We then consider two subcases. {\it Subcase 1:  $\sigma_v(e')=\sigma_v(e'')$}.  The momentum $k_{e'''}$ of the third edge $e'''$ is either $k_{e'}+k_{e''}$ or $-k_{e'}-k_{e''}$, due to the delta function associated to the vertex $v$. Applying claim (a) to $k_{e'''}$, we can write $k_{e'''}=\pm k_e+\tilde{k}_{e'''}$, or $\tilde{k}_{e'''}$ with $\tilde{k}_{e'''}$ being independent of $k_e$. Therefore, $k_{e'}+k_{e''}=\mp k_e-\tilde{k}_{e'''}$, or $-\tilde{k}_{e'''}$. We then defined $\tilde{k}_{e,e'}=-\tilde{k}_{e'''}$. {\it Subcase 2:  $\sigma_v(e')=-\sigma_v(e'')$}. Since both $e',e''$ are integrated, by Lemma \ref{Lemma:cycleFreeEdge}, we can assume that $w_1\succ v\succ w_2$. As a result, if both $v_1,v_2$ does not belong to the path from $w_1\succ v\succ w_2$ to the virtual vertex $v_*$, then $k_{e'}+k_{e''}$ does not depend on $k_e$, thus $\tilde{k}_{e',e''}=k_{e'}+k_{e''}$. If either $v_1$ or $v_2$ belongs to this path, without loss of generality, we suppose $v_1\succ w_1\succ v\succ w_2$. Using \eqref{Lemma:SumFreeEdges2:1}, we obtain
		\begin{equation}
			\begin{aligned}\label{Lemma:DepedenceOnFreeEdgeInsideAcycle:E1}
				k_{e'} \ = & \ \sum_{V\in\mathfrak{P}(v)}\sum_{E=(V,v_{E})\in \mathfrak{F}(v)}\mathbf{1}_{v_{E}\notin\mathfrak{P}(v)}\left(-\sigma_{v}(e')\sigma_V(E) \right)k_{E},\\
				k_{e''} \ = & \ \sum_{V\in\mathfrak{P}(w_1)}\sum_{E=(V,v_{E})\in \mathfrak{F}(w_1)}\mathbf{1}_{v_{E}\notin\mathfrak{P}(w_1)}\left(-\sigma_{w_1}(e'')\sigma_V(E) \right)k_{E}
				.
		\end{aligned}\end{equation} 
		In both formulas, $k_e$ appears, but since $\sigma_{v}(e')=\sigma_{w_1}(e'')$, the $k_{e'}+k_{e''}=\pm 2 k_e + \tilde{k}_{e',e''}$, where $\tilde{k}_{e',e''}$ is independent of $k_e$.
	\end{itemize}
	This finishes the proof of claim (b).

\end{proof}

\subsubsection{Singular graphs: The appearance of zero momenta in graphs - 1-Separation}
In this section, we will focus on a very special class of graphs, in which the delta functions enforce some momenta  to be zero. 
\begin{definition}[Singular Graph]\label{Def:SingGraph} A graph  $\mathfrak{G}$ is said to be singular if there  is  an edge $e$ such that its momentum $k_e$ is zero.  The edge $e$ is called the singular edge and the momentum $k_e$ is the singular momentum.
\end{definition}
To understand better singular graphs, we will need the    definition below, inspired by Tutte (cf. \cite{tutte2019connectivity,tutte1984graph}). 
\begin{definition}[1-Separation]\label{Def:1Sep}
	A graph $\mathfrak{G}$ is said to have a 1-separation if $\mathfrak{G}$ is the union   of two components  $\mathfrak{G}_1\cup\mathfrak{G}_2=\mathfrak{G}$ such that the intersection $\mathfrak{G}_1\cap\mathfrak{G}_2$ contains only one edge. 
\end{definition}
\begin{lemma}\label{Lemma:ZeroMomentum}
	For any edge $e\in\mathfrak{E}$, $k_e$ is independent of all free momenta if and only if $\mathfrak{F}_e=\emptyset$, or equivalently the graph is singular with $e$ being the singular edge, which is also equivalent with the fact that the graph  has a 1-separation: it can be decomposed into two components $\mathfrak{G}_1\cup\mathfrak{G}_2=\mathfrak{G}$ such that they are connected only via the edge $e$ (see Figure \ref{Fig31}). As a consequence, consider a cycle of a vertex $v_i$, if one of the vertices of this cycle belongs to $\mathfrak{G}_j$, $j$ can be either $1$ or $2$, then all of the vertices of the cycle belong to $\mathfrak{G}_j$.
	\begin{figure}
		\centering
		\includegraphics[width=.49\linewidth]{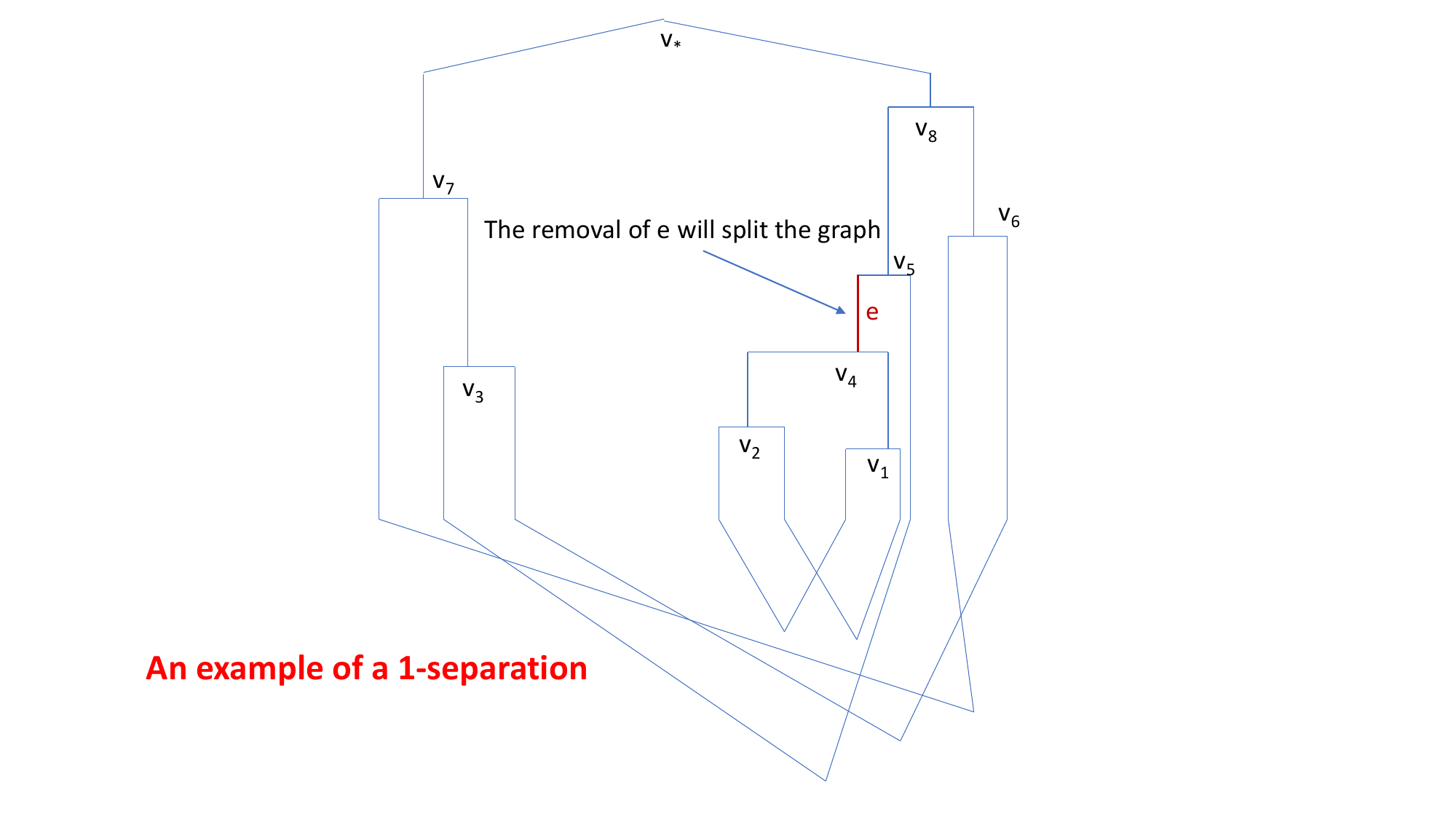}
		\caption{The graph has a 1-separation.}
		\label{Fig31}
	\end{figure}
\end{lemma}
\begin{proof}
	Obviously, $k_e$ is independent of all free momenta if and only if $\mathfrak{F}_e=\emptyset$. It follows from \eqref{Lemma:SumFreeEdges2:1}-\eqref{Lemma:SumFreeEdges2:2} that $k_e$ must be $0$.

	Now we show that if $\mathfrak{F}_e=\emptyset$, then the removal of $e$ will make the number of connected components increase by $1$. Since $\mathfrak{F}_e=\emptyset$, the edge $e=(v_1,v_2)$, $v_1\succ v_2$, is integrated. We suppose the contradiction that the number of connected component does not increase as we remove $e$.  If there is an edge $e'=\{v,v'\}$ such that $v\in \mathfrak{P}(v_1)$ but $v'\notin \mathfrak{P}(v_1)$, then $e'\in\mathfrak{F}$ and thus, $k_e$ depends on $k_{e'}$ and is not $0$. Therefore, if $e'$ is an arbitrary free edge and one of its vertices belongs to $\mathfrak{P}(v_1)$, then the other vertex of $e'$ also belongs to $\mathfrak{P}(v_1)$. Let us now consider any free edge $e'=\{v,v'\},$ whose both vertices belong to $\mathfrak{P}(v_1)$. Since the graph does not split into 2 components when we remove $e$, there is a path, called $\mathcal{P}_{v,v_*}$, that starts from $v$, contains both free and integrated edges, going to the virtual vertex $v_*$, via one of the two topmost vertices $v_{\mathfrak{N}+1},v_{\mathfrak{N}+2}$, and does not contains $e$; otherwise, if we remove $e$, the vertex $v$ will be isolated from the rest of the graph and we have 2 components, contradicting the original assumption. Let us consider the first edge in the path $\mathcal{P}_{v,v_*}$. One of the vertices of this edge is definitely $v$, and we denote the other one by $v''$. There are two cases:
	\begin{itemize}
		\item If the edge $(v,v'')$ is free, then by the previous argument $v''\in\mathfrak{P}(v_1)$. 
		\item If the edge $(v,v'')$ is integrated, then since $v\in\mathfrak{P}(v_1)$, we have $v\succ v_1$. By Lemma \ref{Lemma:cycleFreeEdge}, $v''\succ v$. As a result, $v''\succ v_1$ and thus $v''\in \mathfrak{P}(v_1)$.
	\end{itemize}
	\begin{figure}
		\centering
		\includegraphics[width=.49\linewidth]{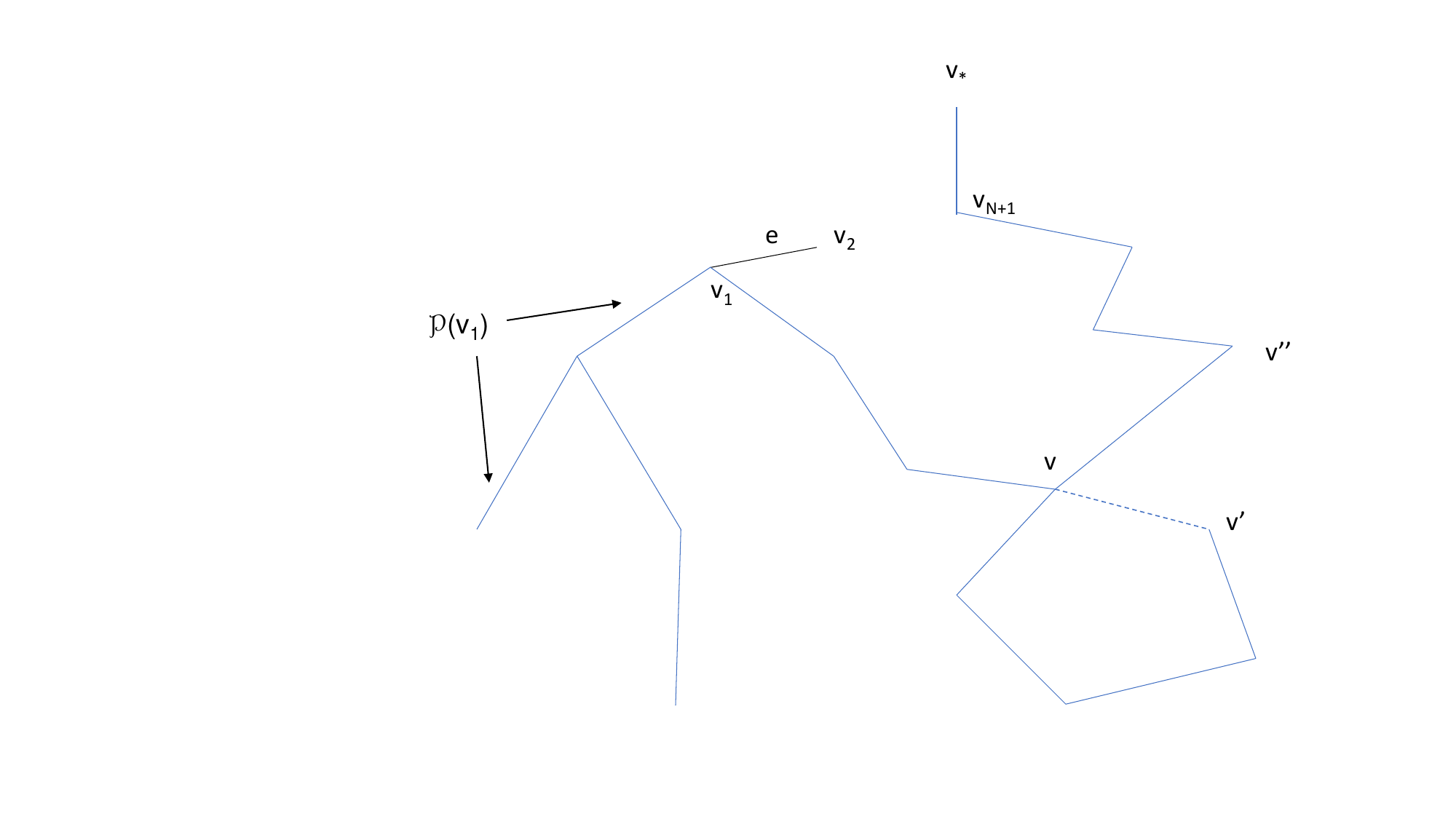}
		\caption{In this picture, on the path going from $v$ to $v_*$,   all of the vertices belong to $\mathfrak{P}(v_1)$. }
		\label{Fig29}
	\end{figure}
	By repeating the above argument, we conclude that all of the vertices of $\mathcal{P}_{v,v_*}$ belong to $\mathfrak{P}(v_1)$ (see Figure \ref{Fig29}). This contradicts the assumption that the path goes to $v_*$ but does not contain $e$. As a result,  the removal of $e$ will make the number of connected components increase by $1$.
	
	\begin{figure}
		\centering
		\includegraphics[width=.49\linewidth]{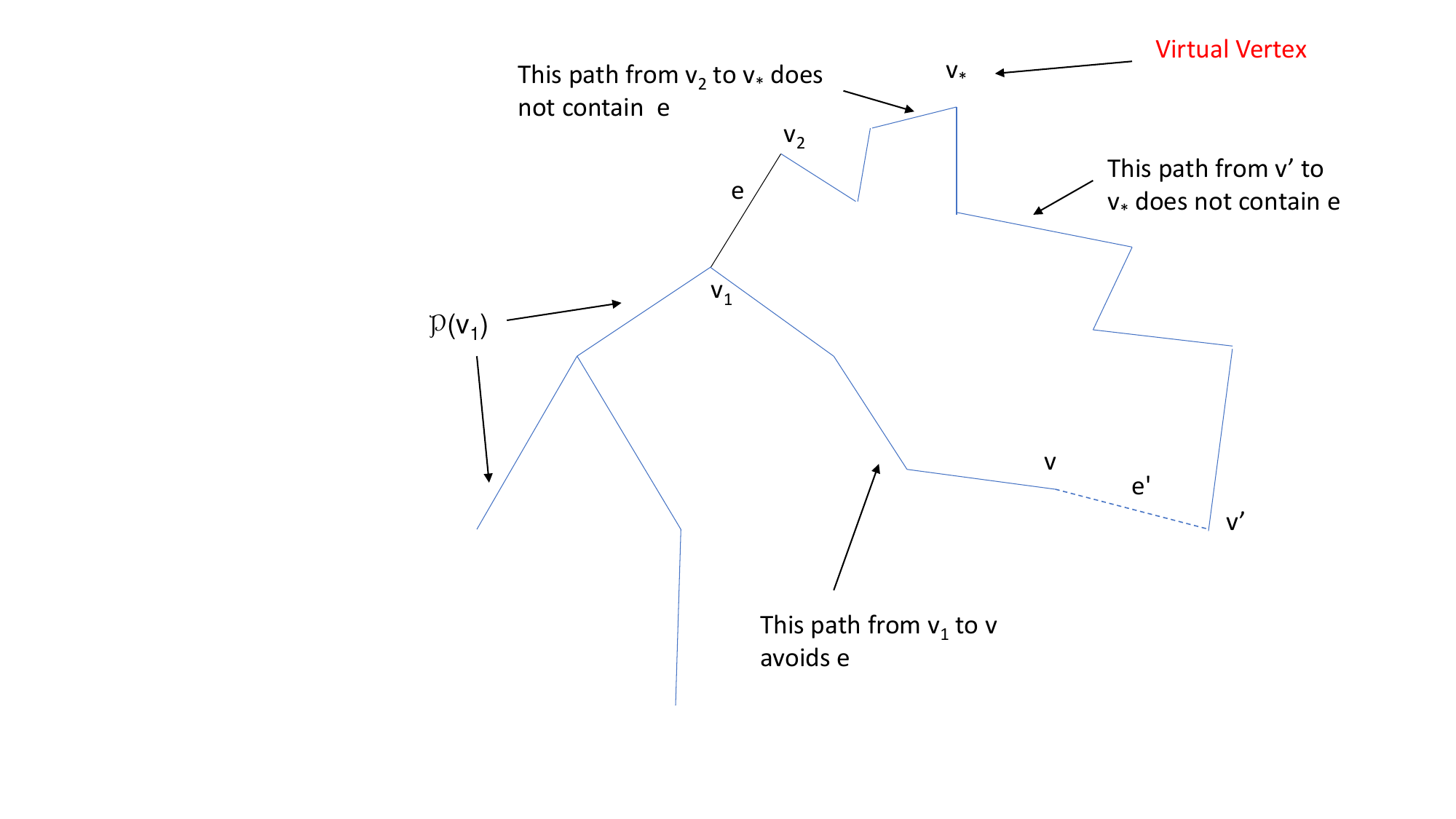}
		\caption{In this picture, the combination of the three paths  $v_2\to v_*$, $v_1\to v$, $v'\to v_*$ and the edge $e'$ form a new path that excludes $e$. }
		\label{Fig30}
	\end{figure}
	Now, suppose that all paths from $v_1$ to $v_2$ contain $e$, that means $v_1$, $v_2$ belong to different components if $e$ is removed, we will have to show that $k_e=0$. Suppose that $k_e$ is not $0$. In this case, $ \mathfrak{F}_e\ne \emptyset$ and we pick $e'=\{v,v'\}\in\mathfrak{F}_e$. By the definition of $\mathfrak{F}_e$, $v\in\mathfrak{P}(v_1)$, $v'\notin \mathfrak{P}(v_1)$. As a result the paths (that can include only integrated edges) from  $v_1$ to $v$ and from the virtual vertex $v_*$ to $v_2$ do not contain $e$. Since $v'\notin \mathfrak{P}(v_1)$, there exists a path  that contains only integrated edge from $v'$ to the virtual vertex $v_*$ and this path does not contain $e$, we can combine the three edges with $e'$ into a new path from $v_1$ to $v_2$ and this path excludes $e$ (see Figure \ref{Fig30}). This is a contradiction, hence, $\mathfrak{F}_e=\emptyset$. 
	
	The final claim of the Lemma, if one of the vertex of the of $v_i$ cycle belongs to $\mathfrak{G}_j$, then all of the vertices of the cycle belong to $\mathfrak{G}_j$, follows from standard results from graph theory (cf. \cite{tutte1984graph}[Theorems III.4 and III.8]).

\end{proof}

\subsubsection{Connectivity of components of graphs -  2-Separation}
We first define the concept of 2-Separation, which is an extension of Definition \ref{Def:1Sep} and  also inspired by Tutte  \cite{tutte2019connectivity,tutte1984graph}.
\begin{definition}[2-Separation]\label{Def:2Sep}
	A graph $\mathfrak{G}$ is said to have a 2-separation if $\mathfrak{G}$ is the union   of two components $\mathfrak{G}_1$ and $\mathfrak{G}_2$: $\mathfrak{G}_1\cup\mathfrak{G}_2=\mathfrak{G}$ such that the intersection $\mathfrak{G}_1\cap\mathfrak{G}_2$ contains exactly two  edges. 
\end{definition}
\begin{lemma}\label{Lemma:RemovalofTwoEdges}
	Suppose that $e,e'\in\mathfrak{E}$, $e\neq e'$, we then have:
	\begin{itemize}
		\item[(i)] Suppose that $\mathfrak{F}_e=\mathfrak{F}_{e'}\ne \emptyset,$ then  if one removes both $e$ and $e'$, the graph $\mathfrak{G}$ is split into  disconnected components (see Figure \ref{Fig32}). As a result, the graph has a 2-separation. 
		\item[(ii)] $\mathfrak{F}_e=\mathfrak{F}_{e'}$ if and only if there is $\sigma\in\{\pm 1\}$ such that $k_e=\sigma k_{e'}$, where $k_e,k_{e'}$ are the momenta associated to $e,e'$.  
	\end{itemize}

\end{lemma}
\begin{proof} {\it (i) The removal of $e,e'$.} We first show that the removal of $e$ and $e'$ splits the graph.
	In the case that both $e,e'$ are  free, since $\mathfrak{F}_e=\mathfrak{F}_{e'}$, we then have $e=e'$. This contradicts the assumption that $e\ne e'$. Thus, at most one of the two edges $e,e'$ is free. We then consider the following two cases.
	
	{\it Case 1: $e'$ is free, $e$ is not free.} If one of the two momenta  $k_e,k_{e'}$ is $0$, then $\mathfrak{F}_e=\mathfrak{F}_{e'}=\emptyset$, leading to a contradiction. As a result, both of them are non zero. We now suppose $e=\{v_1,v_2\}$, with $v_1\succ v_2$, then $\mathfrak{F}_e=\mathfrak{F}_{e'}=\{e'\}$, since $e'$ is free.  Let $e''=\{w,w'\}$ be an edge that has $w\in\mathfrak{P}(v_1)$ and $w'\notin\mathfrak{P}(v_1)$. If $e''\neq e$, we deduce from Lemma \ref{Lemma:SumFreeEdges} that $e''$ is free, then $e''=e'$, and as thus, $e$ is the only edge, besides $e'$, that bridges $\mathfrak{P}(v_1)$ and its complement $\mathfrak{P}(v_1)^c$. As a result, the graph, after the removal of $e$, $e'$, has two disconnected components $\mathfrak{P}(v_1)$ and $\mathfrak{P}(v_1)^c$.
	
	{\it Case 2: Both $e,e'$ are integrated. } We suppose $e'=(v_1',v_2')$, $e=(v_1,v_2)$ with $v_1'\succ v_2'$ and $v_1\succ v_2$. 
	
	{\it Case 2.1. } Let us first consider the  subcase that  $v_1\succ v_1'$, then $\mathfrak{P}(v_1)$ is a subset of $\mathfrak{P}(v_1')$. We proceed with the following two arguments (see Figure \ref{Fig32}).
	\begin{itemize}
		\item Consider a free edge $e''=\{w,w'\}\in\mathfrak{F}_{e}$ and $w\in \mathfrak{P}(v_1)$, $w'\in\mathfrak{P}(v_1)^c$. According to Lemma \ref{Lemma:SumFreeEdges} and the definition of $\mathfrak{F}_{e}$,  such a free edge exits, otherwise, $k_{e}=0$ and $\mathfrak{F}_{e}=\emptyset$. Since $\mathfrak{F}(v_1)\subset\mathfrak{F}(v'_1)$, it follows that $e''\in\mathfrak{F}(v_1')$.  Since $w\in\mathfrak{P}(v_1)\subset \mathfrak{P}(v_1')$, according to Lemma \ref{Lemma:SumFreeEdges} and the fact that $\mathfrak{F}_{e'}=\mathfrak{F}_{e}$, $w'\in \mathfrak{P}(v_1')^c$.  We then deduce that any free edge, that one of its vertices belong to $\mathfrak{P}(v_1)$, has the other vertex in $\mathfrak{P}(v_1')^c$ and at least one such free edge exists, otherwise, $k_e=0$. 
		\item Let us now consider a free edge $e''=\{w,w'\}$ such that $w\in\mathfrak{P}(v'_1)\backslash \mathfrak{P}(v_1)$. We will show that there is only one possibility  that $w'\in\mathfrak{P}(v'_1)\backslash \mathfrak{P}(v_1)$. If $w'\in\mathfrak{P}(v_1)$, since $\mathfrak{P}(v_1)\subset \mathfrak{P}(v_1')$, then $w'\in \mathfrak{P}(v_1')$. Since both $w$ and $w'$ belong to $\mathfrak{P}(v_1')$, the edge $e''$ does not belong to $\mathfrak{F}_{e'}$. However, as $w\notin \mathfrak{P}(v_1)$ but $w'\in\mathfrak{P}(v_1)$, it follows that $e''\in\mathfrak{F}_{e}=\mathfrak{F}_{e'}$, leading to a contradiction. Therefore, $w'\in \mathfrak{P}(v_1)^c$. Now, if $w'\in\mathfrak{P}(v_1')^c$, then since $w\in\mathfrak{P}(v_1')$, the edge $e''$ belongs to $\mathfrak{F}_{e'}=\mathfrak{F}_e$. Therefore, either $w$ or $w'$ belongs to $\mathfrak{P}(v_1)$, leading to the second contradiction. As a consequence, $w'\in \mathfrak{P}(v_1)^c\cap \mathfrak{P}(v_1')=\mathfrak{P}(v'_1)\backslash\mathfrak{P}(v_1).$
	\end{itemize}
	Combining the above two arguments, we can see that $\mathfrak{P}(v'_1)\backslash\mathfrak{P}(v_1)$ and $\mathfrak{P}(v_1)\cup\mathfrak{P}(v_1')^c$ form two disconnected components if we remove $e$ and $e'$. 
	\begin{figure}
		\centering
		\includegraphics[width=.49\linewidth]{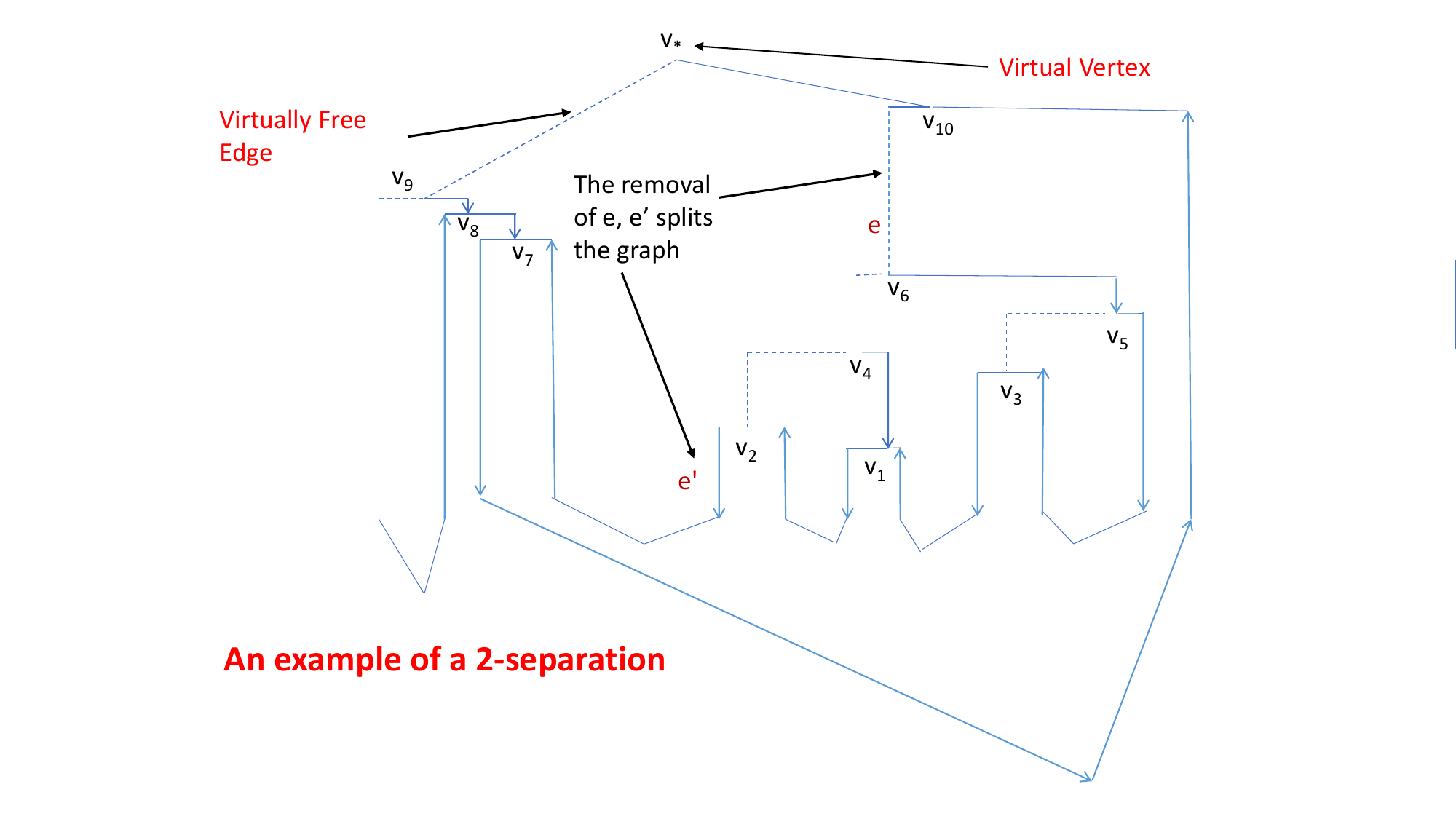}
		\caption{This picture is an example of a graph with a 2-separation. In this picture, the arrows correspond to {\it the second assigned orientation}, that leads to $v_9\succ v_8\succ v_7\succ v_{10}$, $v_6\succ v_5\succ v_3\succ v_1\succ v_{2}\succ v_{7}\succ v_{10}$ and $v_4\succ v_1$.  The dashed edges are free.  The removal of $e$ and $e'$ splits the graph. We also have $k_e= k_{e'}$. Note that Figure \ref{Fig2} corresponds to {\it the first assigned orientation}.}
		\label{Fig32}
	\end{figure}

	{\it Case 2.2. } Now, we consider the subcase that  $\mathfrak{P}(v_1)\cap\mathfrak{P}(v_1')=\emptyset$. Let us consider a free edge $e''=\{w,w'\}$ in $\mathfrak{F}_{e'}=\mathfrak{F}_e$, in which $w\in\mathfrak{P}(v_1)$ and $w'\in\mathfrak{P}(v_1)^c$. If $w'\notin\mathfrak{P}(v_1')$, then $e''\notin\mathfrak{F}_{e'}=\mathfrak{F}_e$ since $w\notin\mathfrak{P}(v_1')$, leading to
	a contradiction. As a consequence, $w'\in\mathfrak{P}(v_1')$. This means that any free edge  $e''=\{w,w'\}$ in $\mathfrak{F}_{e'}=\mathfrak{F}_e$ has one vertex in $\mathfrak{P}(v_1)$ and the other vertex in $\mathfrak{P}(v_1')$ (see Figure \ref{Fig33}). We can see that $e,e'$ connect $\mathfrak{P}(v_1)\cup\mathfrak{P}(v_1')$ and $(\mathfrak{P}(v_1)\cup\mathfrak{P}(v_1'))^c$. We will show that the removal of $e$ and $e'$ will split these two components. Suppose the contrary that there is an edge $e''=\{w,w'\}$ that connects $(\mathfrak{P}(v_1)\cup\mathfrak{P}(v_1'))^c$ with either $\mathfrak{P}(v_1)$  or $\mathfrak{P}(v_1')$. Without loss of generality, we assume that $w\in (\mathfrak{P}(v_1)\cup\mathfrak{P}(v_1'))^c$ and $w'\in \mathfrak{P}(v_1)$. This edge cannot be free, since, otherwise, $w$ has to be in $\mathfrak{P}(v_1')$, as we proved above. As thus, $e''$ is integrated. Since $e''$ is integrated and it is attached to a vertex in $\mathfrak{P}(v_1)$. There is a path from $w$, to $w'\in\mathfrak{P}(v_1)$, then to $v_1$ and finally ends at the virtual vertex $v_*$. Therefore, $w$ also belongs to $\mathfrak{P}(v_1)$, leading to a contradiction. As a consequence, $e$ and $e'$ are the only edges that connect $\mathfrak{P}(v_1)\cup\mathfrak{P}(v_1')$ and $(\mathfrak{P}(v_1)\cup\mathfrak{P}(v_1'))^c$. The removal of $e$ and $e'$ will, therefore, splits those components. 
	\begin{figure}
		\centering
		\includegraphics[width=.49\linewidth]{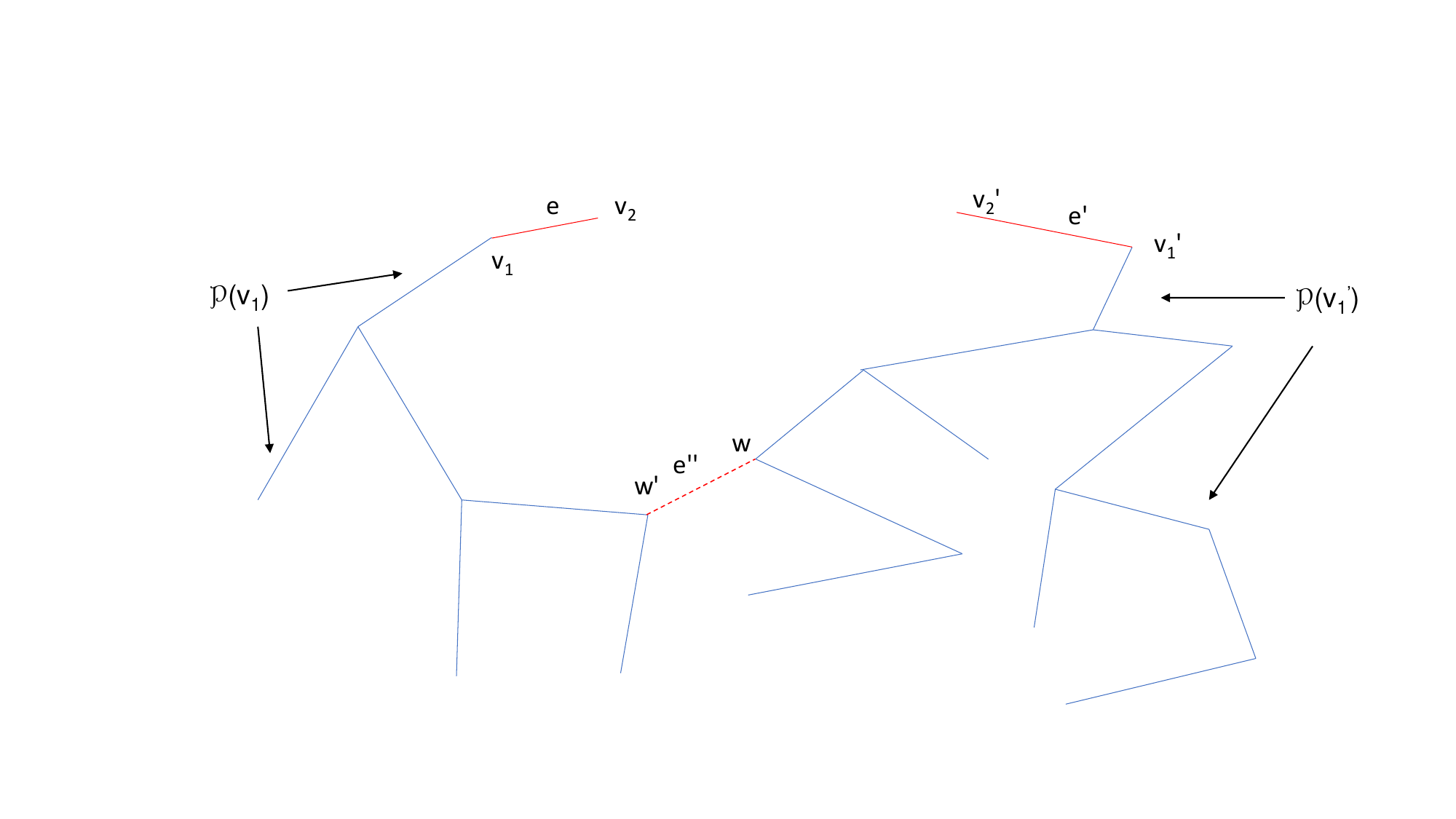}
		\caption{If $e''=\{w,w'\}$ is a free edge and $w\in\mathfrak{P}(v_1')$, then $w'\in\mathfrak{P}(v_1)$. In this case $\sigma_w(e'')=-\sigma_{w'}(e'')$.}
		\label{Fig33}
	\end{figure}

	{\it (ii) The existence of $\sigma$.} If $k_e=\sigma k_{e'}$, it follows that $\mathfrak{F}_e=\mathfrak{F}_{e'}$. 
	We will now show that $\mathfrak{F}_e=\mathfrak{F}_{e'}$ implies $k_e=\sigma k_{e'}$.

	If both of $e, e'$ are free, then  $\mathfrak{F}_e=\mathfrak{F}_{e'}$ implies $e=e'$. If one of them, say $e'$, is free, then $\mathfrak{F}_{e'}=\mathfrak{F}_e=\{e'\}$ and we easily deduce that $e=\pm e'$. 
	
	Now, suppose that both of the two edges are not free and set $e=(v_1,v_2)$, $e'=(v_1',v_2')$ with $v_1\succ v_2$, $v_1'\succ v_2'$.  We write down the formulas for $k_e$ and $k_{e'}$
	\begin{equation}
		\begin{aligned}\label{Lemma:RemovalofTwoEdges:E1}
			k_{e} \ = & \ \sum_{w\in\mathfrak{P}(v_1)}\sum_{e''=(w,w')\in \mathfrak{F}(v_1)}\mathbf{1}_{w'\notin\mathfrak{P}(v)}\left(-\sigma_{w}(e'')\sigma_{v_1}(e) \right)k_{e''},\\
			k_{e'} \ = & \ \sum_{w\in\mathfrak{P}(v_1)}\sum_{e''=(w,w')\in \mathfrak{F}(v_1')}\mathbf{1}_{w'\notin\mathfrak{P}(v_1')}\left(-\sigma_{w}(e'')\sigma_{v_1'}(e') \right)k_{e''}
			.
	\end{aligned}\end{equation} We also consider two cases, that are exactly the same with Cases 2.1 and 2.2 above.

	If  $v_1'\succ v_1$, then $\mathfrak{P}(v_1')\subset \mathfrak{P}(v_1)$, then as argued in Case 2.1, for any free edge $e''=\{w,w'\}\in\mathfrak{F}_e=\mathfrak{F}_{e'}$, we have $w\in\mathfrak{P}(v_1)\cup \mathfrak{P}(v_1')$ and $w'\in\mathfrak{P}(v_1)^c\cup \mathfrak{P}(v_1')^c$. As a result, the constants $\sigma_{w}(e'')$ are the same in the two formulas in \eqref{Lemma:RemovalofTwoEdges:E1}, thus, $k_e=\sigma_{v_1'}(e')\sigma_{v_1}(e)k_{e'}$.
	
	If $\mathfrak{P}(v_1')\cap \mathfrak{P}(v_1)=\emptyset,$ then as argued in Case 2.2, for any free edge $e''=\{w,w'\}\in\mathfrak{F}_e=\mathfrak{F}_{e'}$, one of the vertex, say $w$, belongs to $\mathfrak{P}(v_1)$ and the other vertex, say $w'$, belongs to  $\mathfrak{P}(v_1')$ (see Figure \ref{Fig33}). As a result,  the constants $\sigma_{w}(e'')$ have the opposite orientations  in the two formulas in \eqref{Lemma:RemovalofTwoEdges:E1}, therefore, $k_e=-\sigma_{v_1'}(e')\sigma_{v_1}(e)k_{e'}$.
	
\end{proof}

\subsection{A special class of  integrated graphs: Pairing graphs}
Let us first define the concept of pairing graphs. 
\begin{definition}[Pairing and Non-Pairing Graphs] A graph  $\mathfrak{G}$ is called ``pairing'' if for every $A\in S$, we have $|A|=2$. Otherwise, it is ``non-pairing''.
\end{definition}
\begin{lemma}\label{Lemma:PhaseAtSlice0}
	We  have the following identity for all pairing graphs, for the  phases defined in \eqref{Def:Omega} 
	\begin{equation} \label{Lemma:PhaseAtSlice0:1}
		\vartheta_0\ = \  \sum_{l=1}^n \mathfrak{X}_l \ := \  \sum_{l=1}^n \mathfrak{X}(\sigma_{l,\rho_l},k_{l,\rho_l},\sigma_{l-1,\rho_l},  k_{l-1,\rho_l},\sigma_{l-1,\rho_l+1}, k_{l-1,\rho_l+1}) = \ 0.
	\end{equation}
	Moreover, the phase regulator $\tau_0$ defined in Definition \ref{Def:TauGen} also vanishes
\end{lemma}

\begin{proof}
	Identity \eqref{Lemma:PhaseAtSlice0:1} and the fact that the phase regulator $\tau_0$ vanishes  follow due to the fact that the terms cancel each other pairwise since they are paired.
\end{proof}

\begin{lemma}\label{Lemma:NumberDegreeZeroOne}
	Consider a graph, with $n$ interacting vertices.
	\begin{itemize}
		\item[(i)] If the graph is non-singular and pairing, then the number of degree-zero vertices is $n_0=|S|-1$ and the number of  degree-one vertices is $n_1=n+1-|S|$.
		\item[(ii)] If the graph is non-pairing, and if both of the two edges associated to the top vertices are not virtually free, then  the number of degree-zero vertices is $n_0=|S|-2$ and the number of  degree-one vertices is $n_1=n+2-|S|$.
		\item[(iii)] If the graph is non-pairing, and if there exists an edge among the two edges associated to the top vertices that is virtually free, then the free edge needs to be the one on the left and  the number of degree-zero vertices is $n_0=|S|-1$ and the number of  degree-one vertices is $n_1=n+1-|S|$.
		
	\end{itemize}
\end{lemma}
\begin{proof}
	We first prove $(i)$.
	Let us consider the case when the graph is pairing and non-singular. By Lemma \ref{lemma:freemomenta}, the number of free momenta in this graph is $n+2-|S|$. 
	Suppose that there exists a cluster that connects vertices in $\mathfrak{V}_0$, but those vertices are originated from both the top vertices $v_{n+1}$ and $v_{n+2}$. We will show that one of the two edges attached to $v_{n+1},v_{n+2}$ is free. Suppose that both of them are integrated. Among them, we pick one vertex $v$ originated from $v_{n+1}$ and the other one $v'$ originated from $v_{n+2}$. Since $v$ is originated from $v_{n+1}$, there is a set of vertices $v_{n_1},v_{n_2},\cdots, v_{n_m}$ with $n+1=n_1>n_2>\cdots >n_m$ such that $v_{n_i}$ is connected to $v_{n_{i+1}}$, $i=1,\cdots,m-1$, via an edge and $v_{n_m}$ is connected to $v$ via an edge. Among those edges, let us suppose that $e$ is the first free edge we encounter while going from the bottom to the top. Let us suppose that $e$ connects $v_{n_j}$ and $v_{n_{j+1}}$, if $n_j=n_m$, then we identify $v_{n_{j+1}}$ with $v$.  Let us now consider the cycle of $v_{n_j}$. Since $v$ is inside this cycle, the cycle  goes through another vertex $v''\in\mathfrak{V}_0$, which is connected to a vertex $v_j\in\mathfrak{V}_I$ and originated from $v_{n+1}$. As a result, $v''$ is connected to  $v$ by a cluster and as thus, $v,v',v''$ are in the same cluster (see Figure \ref{Fig6}). The graph is non-pairing, which is a contradiction. Therefore,  one of the two edges attached to $v_{n+1},v_{n+2}$ is a virtually free edge. According to our construction, the left one, which is the one connected to $v_{n+1}$, is the virtually free edge. This means that there are $n+1-|S|$ free edges attached to the interacting vertices of $\mathfrak{V}_I$, and there are $n+1-|S|$ degree-one vertices. The number of degree-zero vertices is then $n-(n+1-|S|)=|S|-1$.
	\begin{figure}
		\centering
		\includegraphics[width=.49\linewidth]{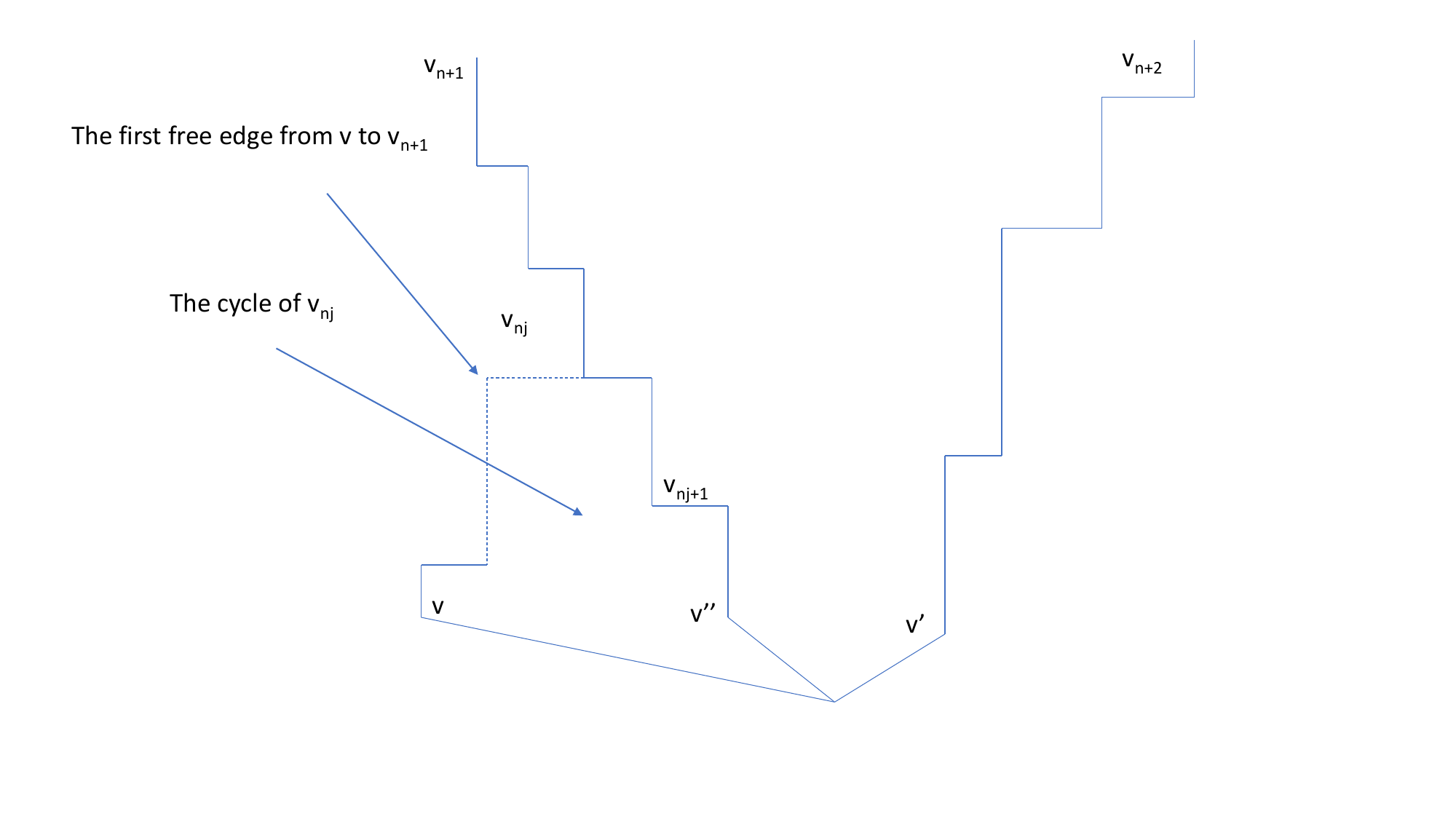}
		\caption{In this picture, the cycle of $v_{n_j}$ belongs to the branch originated from $v_{n+1}$ of the graph. This cycle has two vertices $v, v''$ in $\mathfrak{V}_0$. Then $v,v',v''$ belong to the same cluster and the graph is non-pairing. Note that if $v_{n_j}$ coincides with $v_{n+1}$, then since $v_{n+1}$ is not an interacting vertex, and there is no cycle associated to $v_{n+1}$.}
		\label{Fig6}
	\end{figure}
	
	Now, let us consider the case when there is no cluster that connects vertices in $\mathfrak{V}_0$, that are originated from both the top vertices $v_{n+1}$ and $v_{n+2}$. As a result, the vertices in $\mathfrak{V}_0$ originated from $v_{n+1}$ and $v_{n+2}$ are isolated from each other. As a result, we can remove the edge associated to $v_{n+1}$ and split the graph into two disconnected component. This simply means that this edge is singular and the graph is also singular, leading to a contradiction.
	
	We next prove $(ii)$. By the hypothesis, both of the two edges associated to the two vertices in $\mathfrak{V}_T$ are not virtually free, we then deduce for each of the  $n+2-|S|$ free edges, at least one of its vertices is an interacting vertex. Each free edge is then associated to a degree-one vertex. Thus, there are in total $n_1=n+2-|S|$ degree-one vertices and $n_0=n-n_1=|S|-2$ degree-zero vertices.
	
	Finally, $(iii)$ can be shown by observing that among the two  edges, only one of them can be the virtually free edge, since, by our construction, the addition of the virtually free edge associated to $v_{n+1}$ is to avoid the cycle containing the virtual vertex $v_*$. Once this free edge is added, the non-existence of this cycle is guaranteed, and thus the virtually free edge attached to $v_{n+2}$ is not necessary. Since one of the free edges  is (virtually) attached to $v_{n+1}$, there are $n+1-|S|$ attached to the interacting vertices in $\mathfrak{V}_I$. Each free edge is then associated to a degree-one vertex. Thus, there are in total $n_1=n+1-|S|$ degree-one vertices and $n_0=n-n_1=|S|-1$ degree-zero vertices.
	
	\begin{figure}
		\centering
		\includegraphics[width=.49\linewidth]{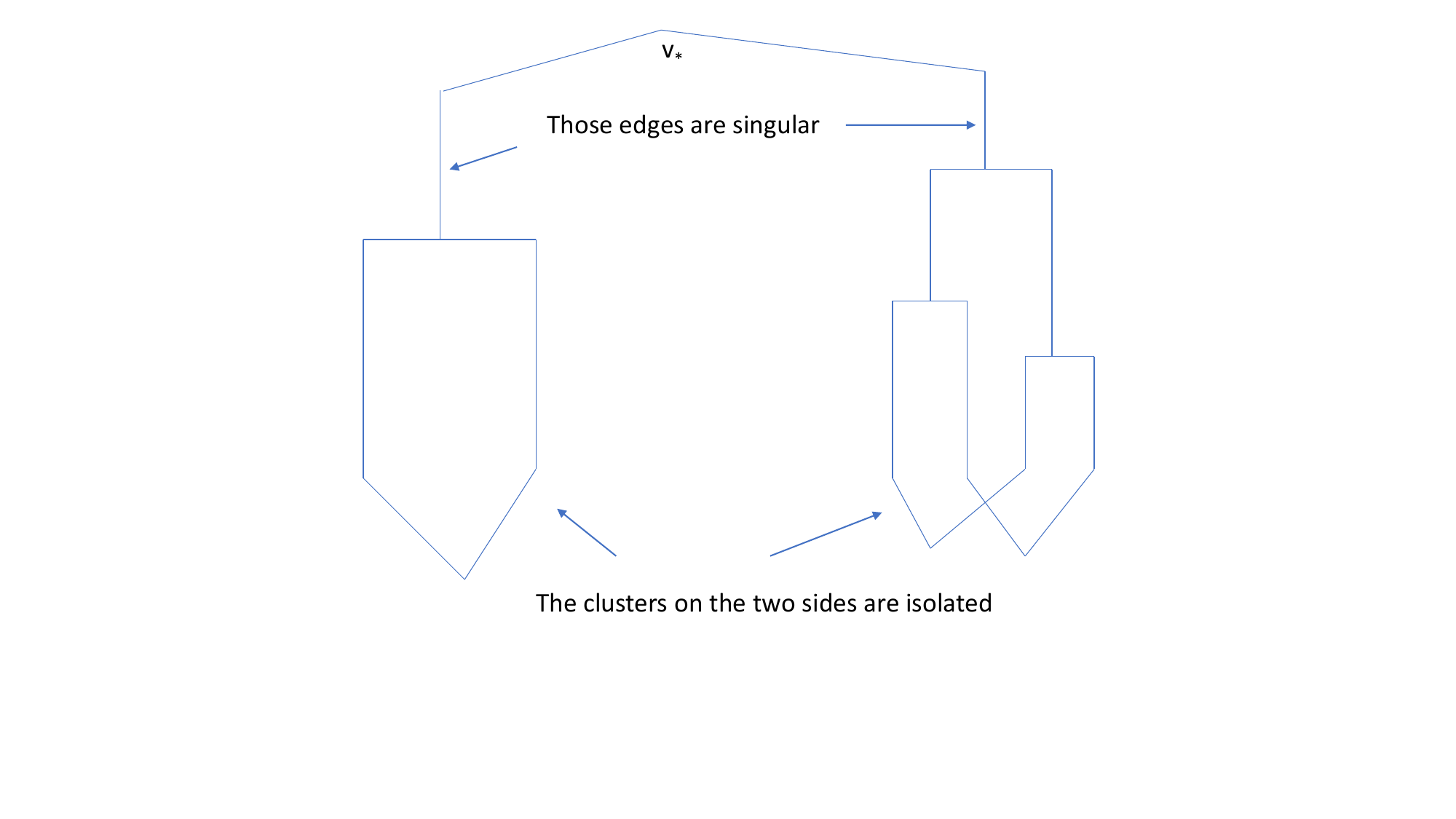}
		\caption{In this picture, the clusters on the left and on the right are isolated. The graph has a 1-separation.}
		\label{Fig34}
	\end{figure}
\end{proof}

We will show later that pairing graphs are the main contribution to the Duhamel multi-layer expansions. However, among all pairing graphs, some are less important than the others. As a result, we need a classification of pairing graphs.
\subsection{Classification of non-singular pairing graphs}
The goal of this section is to give a detailed study on pairing graphs. Among the pairing graphs, we exclude those that are singular. Our classification is based on the concept of the interacting size of an $\mathrm{iC}_{m}$ cycle , which is an extended version of the concept of the size of a $C_{l}$ cycle (defined in Definition \ref{Def:EstradaSizecycle}),  discussed in \cite{estrada2012structure},  due to the distinction between vertices in $\mathfrak{V}_I$ and $\mathfrak{V}_0\cup\mathfrak{V}_C$ in our graphs.  Moreover, differently from classical configurations  in graph theory, the vertices in $\mathfrak{V}_I$  are ordered from the top to the bottom, we then further classify the important  $\mathrm{iC}_{2}$ cycles into $\mathrm{iC}_{2}^d$ and $\mathrm{iC}_{2}^r$ using this ordering. 
\begin{remark} Note that, a pairing graph can also be singular. For instance,  Figure \ref{Fig31} is an example of a pairing, singular graph. 
\end{remark}

To make a classification of non-singular pairing graphs, we will first classify cycles, associated to degree-one vertices, inside a non-singular pairing graph. We then have the classical definition of the size of a cycle, following  Estrada  (cf. \cite{estrada2012structure}).
\begin{definition}[Size of a Cycle]\label{Def:EstradaSizecycle}
	Let $v_{i}$ be a degree-one vertex. Define $\mathcal{I}$ to be the set of all of the vertices in the cycle of $v_i$, including $v_i$. Then the number of elements $l=|\mathcal{I}|$ of $\mathcal{I}$ is defined to be the size of the cycle of $v_i$. Those cycles are  called $C_{l}$ cycles.
\end{definition}
The above definition counts all of the vertices in $\mathfrak{V}_C,\mathfrak{V}_0$ and  $\mathfrak{V}_I$ of a cycle. Since in our proof, we are mostly interested in the interacting vertices, defined in Section \ref{Sec:MomentumGraph},  it is important that we count the number of interacting vertices in a cycle. As thus, the following definition is introduced, to guarantee the classification of $C_l$ cycles based on counting the number of  vertices in $\mathfrak{V}_I$. 
\begin{definition}[Interacting Size of a Cycle]\label{Def:Sizecycle}
	Let $v_{i}$ be a degree-one vertex. Define $\{v_{l}\}_{l\in\mathfrak{I}}$ to be the set of all of the interacting vertices in the cycle of $v_i$, including $v_i$. Then the number of elements $m=|\mathfrak{I}|$ of $\mathfrak{I}$ is defined to be the interacting size of the cycle of $v_i$. The cycles are then called $\mathrm{iC}_{m}$ cycles (see Figure \ref{Fig7} for an illustration).
\end{definition}

We have the following lemmas. 
\begin{lemma}\label{Lemma:Size1cycle}
	Let $v_{i}$ be a degree-one vertex and  $k\in\mathfrak{E}_+(v_{i})$.  Suppose that the interacting size of the cycle of $v_i$ is $1$, then $k=0$ and the graph containing the cycle of $v_{i}$ has a 1-separation and is singular. Therefore, any cycle in a non-singular graph has an interacting size bigger than $1$. 
\end{lemma}
\begin{proof}
	Suppose that $k_1,k_2\in\mathfrak{E}_-(v_{i})$ . Since the cycle containing $v_{i}$ has only one interacting vertex, the two edges $k_1,k_2$ are directly connected. This means  if we remove the edge associated to $k$, the cycle is split from the graph. By Lemma \ref{Lemma:ZeroMomentum}, $k=0$ and the graph has a 1-separation and is singular.
\end{proof}

The following lemma illustrates the structure of vertices inside a cycle.

\begin{lemma}\label{Lemma:VerticeOrderInAcycle}
	Let $v_i$ be a degree-one vertex and define $\{v'_{l}\}_{l\in\mathfrak{I}}$ to be the set of all of the interacting vertices in the cycle of $v_i$, including $v_i$. Then for all $l\in\mathfrak{I}$, we always have $\mathcal{T}(v'_l)\le i$.
\end{lemma}
\begin{proof}
	
	%
	Suppose the contrary, that there exists a vertex $v_j$, $j>i$, such that this vertex has one of its associated edges depending on the free edge  $e=\{v_i,w\}$  associated to $v_i$, and $ i=\mathcal{T}(v_i)>\mathcal{T}(w)$. It follows that $e\cap\mathfrak{V}_T=\emptyset$. We now consider the  paths, with the second assigned orientation, going from $v_i$, $w$ to the virtual vertex $v_*$. These paths must coincide starting from a vertex $w_0$. We have $v_i\succ w_0$ and $w\succ w_0$. Suppose that $e'=(w_1,w_2)$, with $w_1\succ w_2$, satisfies $w\succ w_0\succ w_1$ and $v_i\succ w_0\succ w_1$. Then by Lemma \ref{Lemma:SumFreeEdges}, $k_{e'}=\sigma_{w_1}(e')\sigma_{w}(e)k_e+\dots$ as $w\succ w_0\succ w_1$. Similarly, since $w'\succ w_0\succ w_1$,  then $k_{e'}=\sigma_{w_1}(e')\sigma_{w'}(e)k_e+\dots$. We notice that $\sigma_{w}(e)=-\sigma_{w'}(e)$. Therefore,  $\sigma_{w}k_e$ and $-\sigma_{w'}k_e$ cancel with each other. As a result, $k_{e'}$ is independent of $k_e$.

	From the above arguments, we deduce that the vertex $v_j$ belongs to either the path  $w\succ w_0$ or the path  $v_i\succ w_0$. Let us now revisit the  construction of free edges scheme. In this scheme, the construction of free edges is done from the bottom to the top of the tree. The edge $k_e$ is set free due to the fact that it forms a cycle with a set of edges associated to vertices $\{v_{l'}\}_{l'\in \mathfrak{V}},$ with $l'\le i$ for $l'\in \mathfrak{V}$. As a result, the  paths  $w \succ v_j$, $v_j\succ w_0$, $v_i\succ w_0$ and the vertices $\{v_{l'}\}_{l'\in \mathfrak{V}}$  form a cycle. This means  they coincide, otherwise, we have a contradiction. As a result, $j$ belongs to the set $\mathfrak{V}$, and $j\le i$, contradicting the original assumption that $j>i$. 
\end{proof}

Below, we present a classification of cycles for a non-singular graph.

\begin{figure}
	\centering
	\includegraphics[width=.49\linewidth]{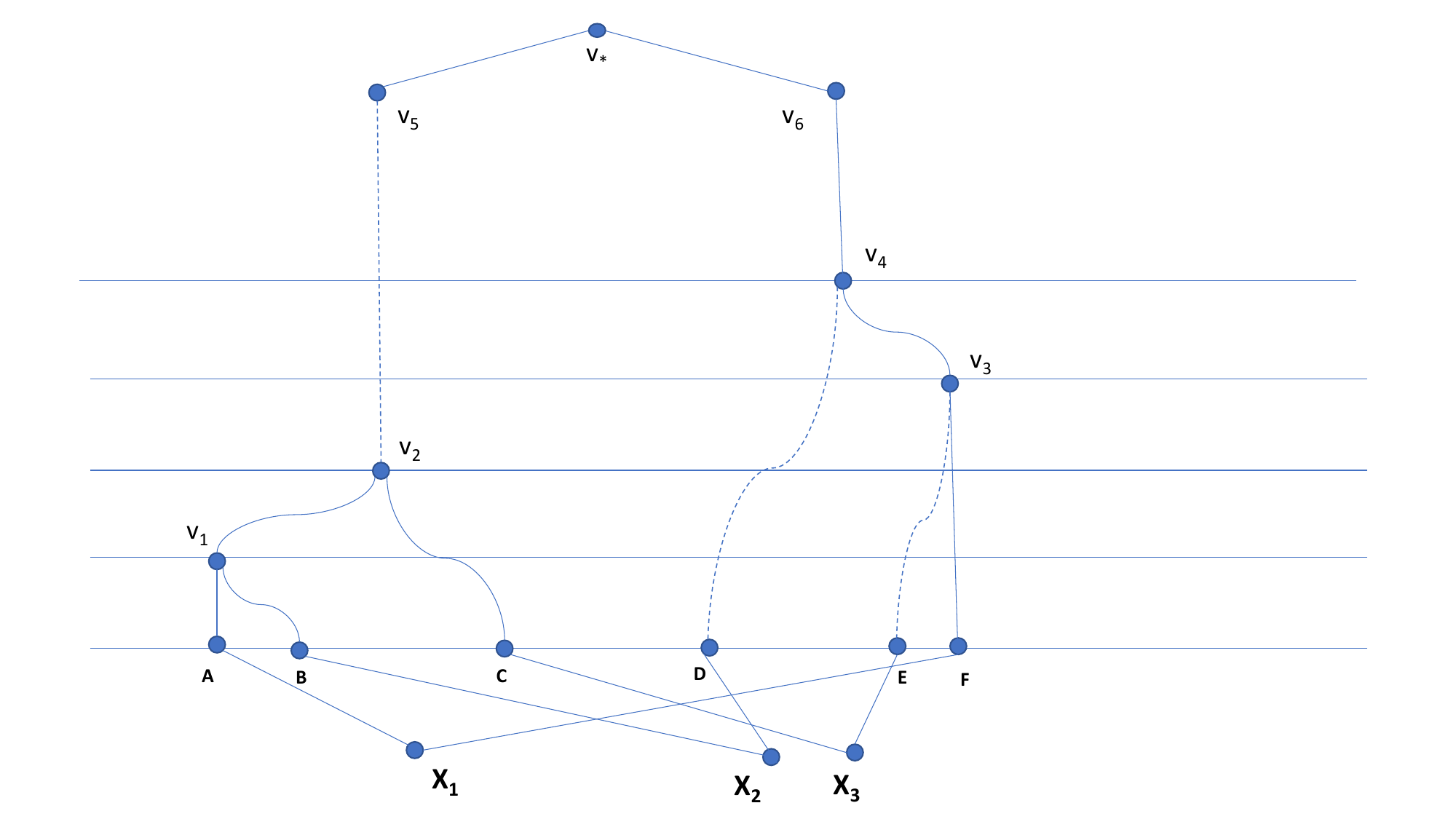}
	\caption{In this picture, the dashed edges are free. The cycles of $v_4$ and $v_3$ have three vertices $\{v_4,v_3,v_1\}$ and $\{v_3,v_2,v_1\}$. Therefore, these cycles have the interacting size of $3$ and they are $\mathrm{iC}_3$ cycles. The cycle of $v_4$ has 9 vertices $\{v_4,v_3,F,X_1,A,v_1,B,X_2,D\}$ and is a $C_9$ cycle. The cycle of $v_3$ has 9 vertices $\{v_3,E,X_3,C,v_2,v_1,A,X_1,F\}$ and is a $C_9$ cycle.  They are long collisions, as described in Definition \ref{Def:Classificationcycles}. }
	\label{Fig8}
\end{figure}
\begin{definition}[Classification of Cycles in Non-singular Graphs]\label{Def:Classificationcycles}
	In a non-singular graph $\mathfrak{G}$, let $v_i$ be a degree-one vertex and $\{v_j\}_{j\in\mathfrak{I}}$ be the sets of all the interacting vertices in its cycle. The cycle of $v_i$ is then an $\mathrm{iC}_{|\mathfrak{I}|}$ cycle.
	\begin{itemize}
		\item  If the size $|\mathfrak{I}|$ of this $\mathrm{iC}_{|\mathfrak{I}|}$ cycle is bigger than or equal to  $3$, the cycle is called a ``long collision'' (see Figure \ref{Fig8}). 
		\item  If the size $|\mathfrak{I}|$ of this $\mathrm{iC}_{|\mathfrak{I}|}$ cycle is exactly  $2$, the cycle is called a ``short collision'' (see Figure \ref{Fig7}).
	\end{itemize}
\end{definition}
\begin{figure}
	\centering
	\includegraphics[width=.49\linewidth]{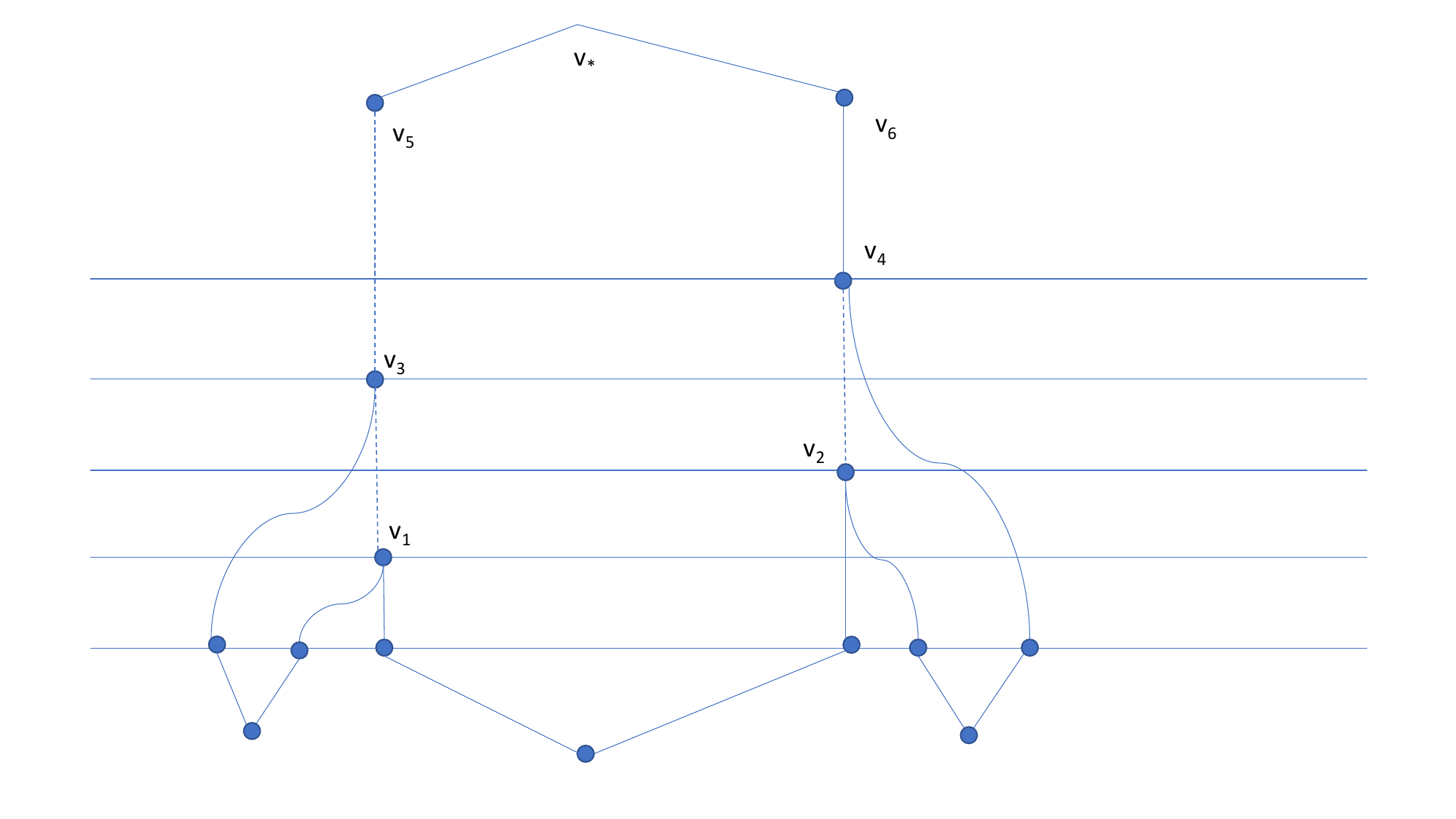}
	\caption{In this picture, the dashed edges are free. The cycles of $v_4$ and $v_3$ have only two interacting vertices $\{v_4,v_2\}$ and $\{v_3,v_1\}$. Therefore, these cycles are $\mathrm{iC}_2$ cycles, which means have the interacting size of $2$, and they are short. They are, in addition, short delayed recollisions  $\mathrm{iC}_2^d$, as described in Definition \ref{Def:ClassificationShortCollision}.}
	\label{Fig7}
\end{figure}
In our graphs, the ordering of the interacting vertices from the top to the bottom plays a very important role. Though $\mathrm{iC}_2$ cycles have only two interacting vertices $v_i,v_j$,  we need to further classify $\mathrm{iC}_2$ cycles based on this ordering. We then have the following classification for  $\mathrm{iC}_2$ short collisions. 
\begin{definition}[Classification of Short Collisions in Non-singular Graphs]\label{Def:ClassificationShortCollision}
	In a non-singular graph $\mathfrak{G}$, let $v_i$ be a degree-one vertex, whose cycle is a short collision and contains two interacting vertices $v_i,v_j$ with $j<i$.
	\begin{itemize}
		\item If $j=i-1$, the cycle is called a ``recollision'' (see Figures \ref{Fig10}-\ref{Fig11}). There are two types of recollision. If a recollision uses only one cluster vertex, we call it a ``single-cluster recollision''. A single-cluster recollision is an $\mathrm{iC}_2$ and a $C_5$ cycle. If a recollision uses two cluster vertices, we call it a ``double-cluster recollision''. A double-cluster recollision is an $\mathrm{iC}_2$ and a $C_8$ cycle. A recollision is called an ``$\mathrm{iC}^r_2$ cycle.''
		\item If $j\ne i-1$, the cycle is called a ``short delayed recollision'' (see Figure \ref{Fig7}). By Lemma \ref{Lemma:VerticeOrderInAcycle}, it follows that $j<i-1$. A short delayed recollision  is called an ``$\mathrm{iC}^d_2$ cycle.''
		\item For both recollisions and short delayed recollisions, each cycle has the same number of interacting vertices (see Figures \ref{Fig7},\ref{Fig10},\ref{Fig11}). 
	\end{itemize}
\end{definition}
\begin{figure}
	\centering
	\includegraphics[height=.50\linewidth]{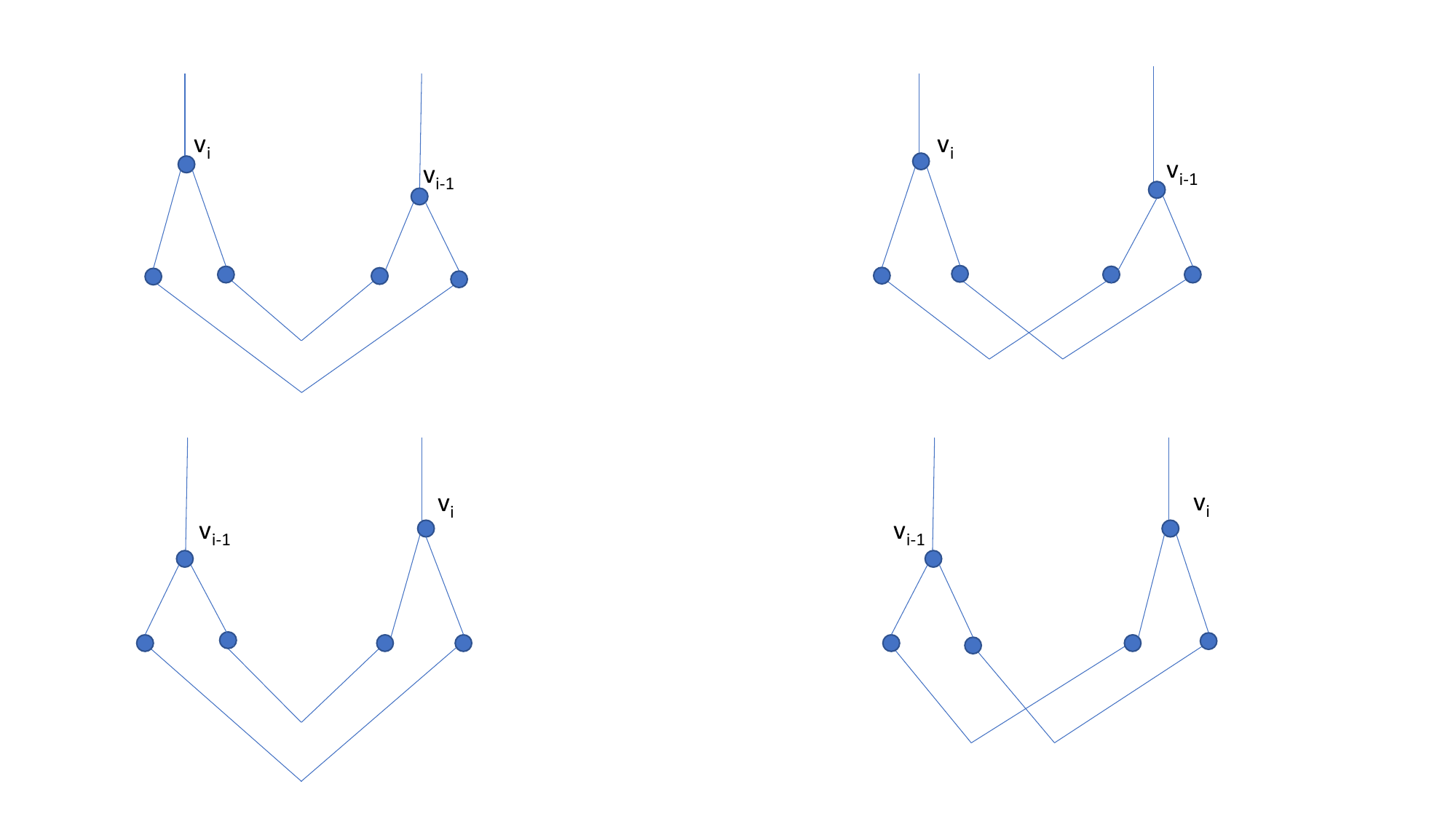}
	\caption{The 4 double-cluster recollisions. Each of them uses two cluster vertices.  These recollisions are not delayed as the one in Figure \ref{Fig7}. They are $\mathrm{iC}_2^r$ and $C_8$ cycles, the skeletons of  $\mathrm{iCL}_2$ ladder graphs.}
	\label{Fig10}
\end{figure}
\begin{figure}
	\centering
	\includegraphics[height=.50\linewidth]{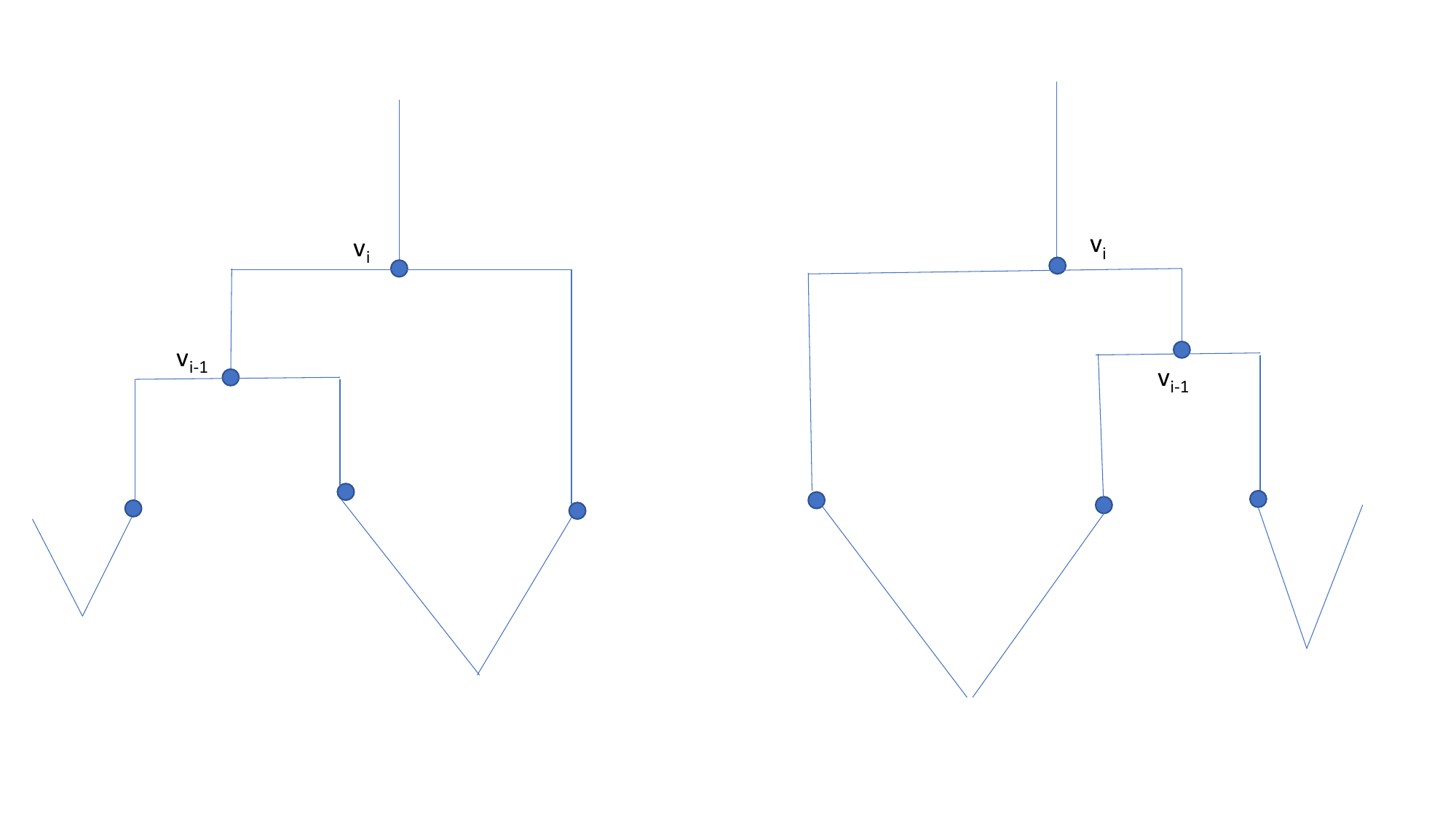}
	\caption{The 2 single-cluster recollisions. Each of them uses only one cluster vertex. These recollisions are not delayed as the one in Figure \ref{Fig7}. They are $\mathrm{iC}_2^r$ and $C_5$ cycles, the skeletons of  $\mathrm{iCL}_2$ ladder graphs.}
	\label{Fig11}
\end{figure}

The following Lemma shows the dependence of the phase of a vertex on an arbitrary free edge.
\begin{lemma}\label{Lemma:XiOnFreeEdge}
	Let $v_{i}$ be an interacting vertex in a non-singular graph and $e\in\mathfrak{E}(v_i)$. Let $e'$ be an arbitrary free edge in $\mathfrak{F}_e$. The followings hold true.
	\begin{itemize}
		\item[(i)] $\nabla_{k_{e'}}\mathfrak{X}(v_{i})\neq 0$.
		\item[(ii)] $\mathfrak{X}(v_{i})$ cannot be independent of all free momenta.
	\end{itemize}
\end{lemma}
\begin{proof} We first prove (i). Denote the two vertices of $e'$ by $w$ and $w'$. 
	Let us suppose the contrary that $\mathfrak{X}(v_i)$ is a constant in $k_{e'}$, hence independent of  $k_{e'}$. Denote the three edges associate to $v_i$ by $e_0,e_1,e_2$ with the momenta $k_0,k_1,k_2$, signs $\sigma_{k_0},\sigma_{k_1},\sigma_{k_2}$, and the three vertices associate to them by $v'_0,v'_1,v'_2$. If the free edge $e'$ depends on all of the 3 momenta $k_0,k_1,k_2$, we consider two cases. In the first case, we suppose that $e'$ is attached to $v_i$, then $k_{e'}$ 
	is one of the 3 momenta $k_0,k_1,k_2$, say $k_0$. Since the other two edges of $v_i$ also depends on $e'$, we have $v_i\succ v_1'$ and $v_i\succ v_2'$, contradicting the conclusion of  
	Lemma \ref{Lemma:cycleFreeEdge}. In the second case, we suppose that $e'$ is not attached to $v_i$. Since $k_0,k_1,k_2$ all depend on $k_{e'}$, for $k_j$, $j=\{0,1,2\}$, we then have one of the following possibilities: 
	(a)  $w$ or $w'$ $\succ v_i$;
	(b) $w$ or $w'$ $\succ v_j'$. The first possibility (a) cannot happens, since suppose without loss of generality that $w\succ v_i$, then we must have $v_i\succ v_j'$ for all  $j\in\{0,1,2\}$ so that $k_0,k_1,k_2$ all depend on $k_{e'}$, leading to a contradiction with the conclusion of  
	Lemma \ref{Lemma:cycleFreeEdge} again. As a result, the second possibility (b) must happen, this means $w$ or $w'$ $\succ v_j'$ for all $j\in\{0,1,2\}$. By the Pigeonhole Principle, we can suppose without loss of generality that there exist $j\ne l\in \{0,1,2\}$ such that $w\succ v_j'$ and  $w\succ v_l'$. This means that $w, v_j',v_l'$ and $v_i$ form a cycle. This is again another contradiction. Therefore,  the free edge $e'$ depends on exactly 2 of the 3 momenta $k_0,k_1,k_2$. 
	
	We now suppose that $k_{e'}$ depends on $k_0$ and $k_1$. By Lemma \ref{Lemma:DepedenceOnFreeEdgeInsideAcycle}, we can then write $k_0=\sigma_0 k_{e'}+k_0'$ and $k_1=\sigma_1 k_{e'}+k_1'$, in which $k_0',k_1'$ are the sum of some other free momenta. Thus, $$\nabla_{k_e'}\Big(\sigma_{k_0}\omega(\sigma_0 k_{e'}+k_0')+\sigma_{k_1}\omega(\sigma_1 k_{e'}+k_1')\Big)=0,$$
	which means 
	$$\sigma_{k_0}\omega(\sigma_0 k_{e'}+k_0')+\sigma_{k_1}\omega(\sigma_1 k_{e'}+k_1')=\mathrm{constant}.$$
	This happens only when $\sigma_{k_0}\sigma_0+\sigma_{k_1}\sigma_1=0$, $\sigma_{k_0}k_0'+\sigma_{k_1}k_1'=0$, and therefore $\sigma_{k_0}k_0+\sigma_{k_1}k_1=0$. By the delta function associated to $v_i$, we deduce that $k_2=0$. Thus, the momentum $k_2$ is singular, that contradicts the hypothesis that the graph is non-singular. Therefore $\nabla_{k_{e'}}\mathfrak{X}(v_{i})\neq 0$. This finishes the proof of $(i)$.
	
	Suppose that $\mathfrak{X}(v_{i})$ is independent of all free momenta. This, by $(i)$, shows that for any edge $e\in\mathfrak{E}(v_i)$, the set $\mathfrak{F}_e$ is empty. This means $k_e=0$ and the edge $e$ is singular, contradicting the hypothesis that the graph is non-singular.
	
\end{proof}
We will need to further classify long collisions. To this end, we first prove the following Lemma.

\begin{lemma}\label{Lemma:VerticesLongCollisions}
	Let $v_{l_1}$ be a degree-one vertex in a non-singular graph. Denote by $k_1,k_2$ its edges in $\mathfrak{E}_-(v_{l_1})$ and $k_0$ the edge in $\mathfrak{E}_+(v_{l_1})$. Denote their signs by $\sigma_{k_0},\sigma_{k_1},\sigma_{k_2}$. Suppose that $k_1$ is the momentum of the free edge. Assume that there exists a vertex $v_{l_0}$ in the cycle of $v_{l_1}$ such that $\mathfrak{X}(v_{l_0})+\mathfrak{X}({v_{l_1}})$ is not a function of $k_1$. Denote the three edges of $v_{l_0}$ by $k_0',k_1',k_2'$, in which $k_0'\in\mathfrak{E}_+(v_{l_0})$ and $k_1',k_2'\in\mathfrak{E}_-(v_{l_0})$.  The signs of those edges are denoted by $\sigma_{k_0'},\sigma_{k_1'},\sigma_{k_2'}$.  Then, the followings hold true.
	\begin{itemize}
		\item[(i)] One of the following identities is satisfied
		\begin{equation}
			\label{Lemma:VerticesLongCollisions:1}
			\sigma_{k_2}k_2+\sigma_{k_2'}k_2'=\sigma_{k_1}k_1+\sigma_{k_1'}k_1'=0,
		\end{equation} 
		\begin{equation}
			\label{Lemma:VerticesLongCollisions:2}
			\sigma_{k_2}k_2+\sigma_{k_1'}k_1'=\sigma_{k_1}k_1+\sigma_{k_2'}k_2'=0,
		\end{equation} 
		\begin{equation}
			\label{Lemma:VerticesLongCollisions:3}
			\sigma_{k_1}k_1+\sigma_{k_0'}k_0'=\sigma_{k_2}k_2+\sigma_{k_2'}k_2'=0,
		\end{equation}
		and
		\begin{equation}
			\label{Lemma:VerticesLongCollisions:4}
			\sigma_{k_1}k_1+\sigma_{k_0'}k_0'=\sigma_{k_2}k_2+\sigma_{k_1'}k_1'=0.
		\end{equation}
		\item[(ii)] There does not exist another vertex $v_j$ in the cycle of $v_{l_1}$ such that $\mathfrak{X}(v_{j})+\mathfrak{X}({v_{l_1}})$  is not a function of $k_1$.  Therefore, $v_{l_0}$ is the only vertex with this property.
		\item[(iii)] Suppose that the cycle is long. If \eqref{Lemma:VerticesLongCollisions:1} or \eqref{Lemma:VerticesLongCollisions:2} happens, we can remove the other vertices in the cycle and join $k_1,k_2$ with $k_1',k_2'$ or with $k_2',k_1'$ to form a short delayed recollision or a recollision.   If \eqref{Lemma:VerticesLongCollisions:3} or \eqref{Lemma:VerticesLongCollisions:4} happens, we can remove the other vertices in the cycle and join $k_1,k_2$ with $k_0',k_1'$ or with $k_0',k_2'$ to form a short delayed recollision or a recollision. In both cases, the graph has 2-separations.
	\end{itemize}

\end{lemma}
\begin{proof}  We divide the proof into several parts.
	
	\smallskip
	
	{\bf Part 1 - The proof of (i): One of the  identities \eqref{Lemma:VerticesLongCollisions:1}-\eqref{Lemma:VerticesLongCollisions:4} is satisfied.}
	
	\smallskip
	\begin{figure}
		\centering
		\includegraphics[width=.49\linewidth]{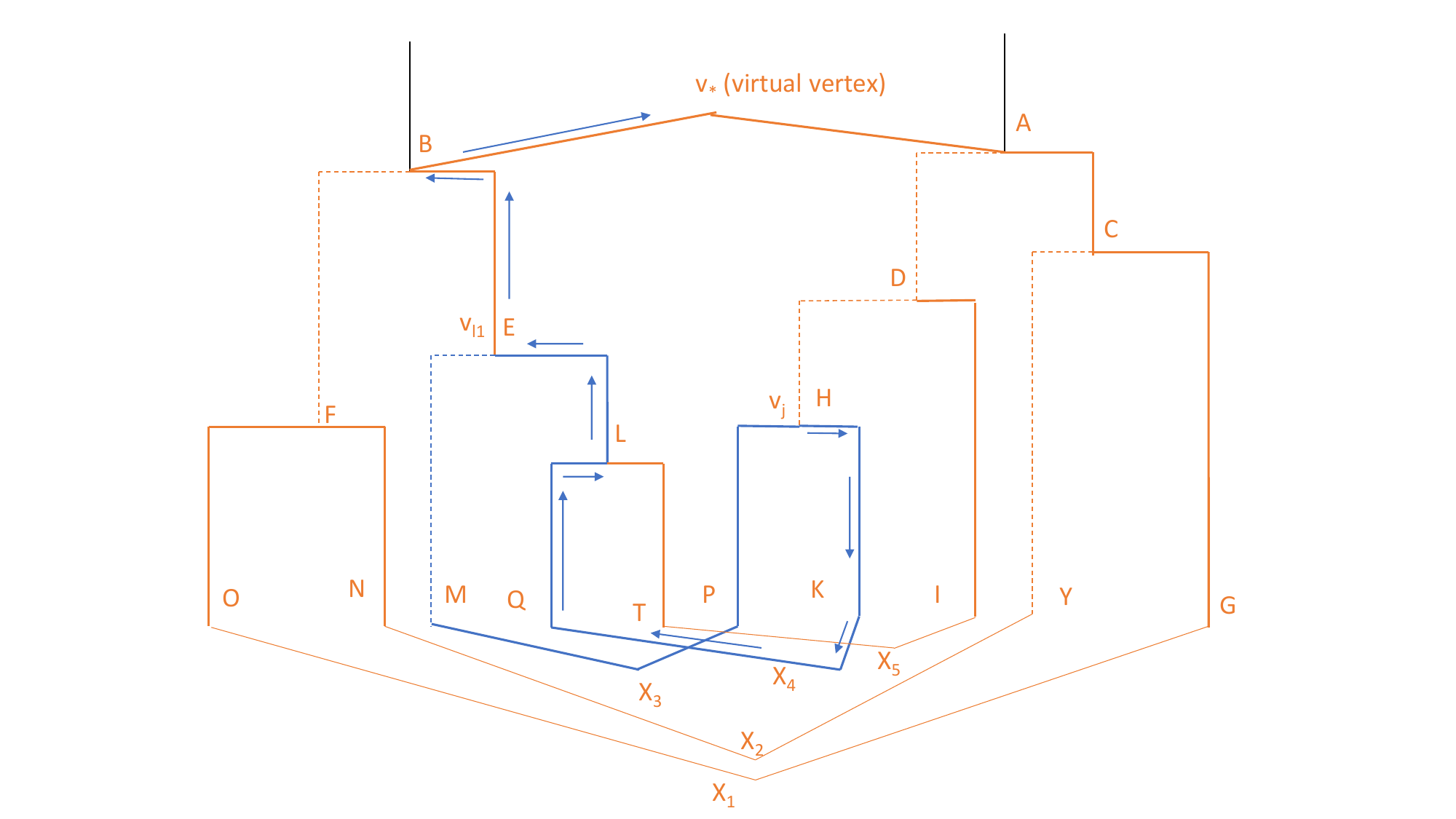}
		\caption{ In this example, the cycle of $v_{l_1}$ is $ELQX_4KHPX_3M$. The free edges are the dashed lines. The path from $v_j$ to the virtual vertex is $HKX_4QLEBv_*$, which is denoted by the arrows. The two vertices $E$ and $M$ of the free edge can be joined into this path by $EBv_*$ and $MX_3PHKX_4QLEBv_*$.}
		\label{Fig23}
	\end{figure}
	
	The cycle of $v_{l_1}$ goes through the vertex $v_{l_0}$ by either the two edges in $\mathfrak{E}_-(v_{l_0})$ or by one edge in $\mathfrak{E}_-(v_{l_0})$ and one edge in $\mathfrak{E}_+(v_{l_0})$. These two edges are then denoted by $h_1,h_2$ and the left-over edge is denoted by $h_0$. As $k_1$ is the free momentum of the cycle, it follows that $\{h_0\}=\{k_0',k_1',k_2'\}\backslash\{h_1,h_2\}$ and $h_1=\sigma_1 k_1+u_1$, $h_2=\sigma_2 k_1+u_2$, in which $u_1,u_2$ are independent of $k_1$ and $\sigma_1,\sigma_2\in\{\pm1\}$. We denote the associated signs by $\sigma_{h_0},\sigma_{h_1},\sigma_{h_2}$. Since $\sigma_{h_0}h_0+\sigma_{h_1}h_1+\sigma_{h_2}h_2=0$, we deduce $h_0=-\sigma_{h_0}(\sigma_{h_1}\sigma_1+\sigma_{h_2}\sigma_2)k_1-\sigma_{h_0}(\sigma_{h_1}u_1+\sigma_{h_2}u_2)$. Therefore, $ -\sigma_{h_0}(\sigma_{h_1}\sigma_1+\sigma_{h_2}\sigma_2)k_1=\pm k_1$ or $0$. Since $\sigma_{h_0}(\sigma_{h_1}\sigma_1+\sigma_{h_2}\sigma_2)$ is an even number, it follows that $\sigma_{h_0}(\sigma_{h_1}\sigma_1+\sigma_{h_2}\sigma_2)=0$. Therefore $h_0$ is independent of $k_1$. 
	\begin{figure}
		\centering
		\includegraphics[width=.49\linewidth]{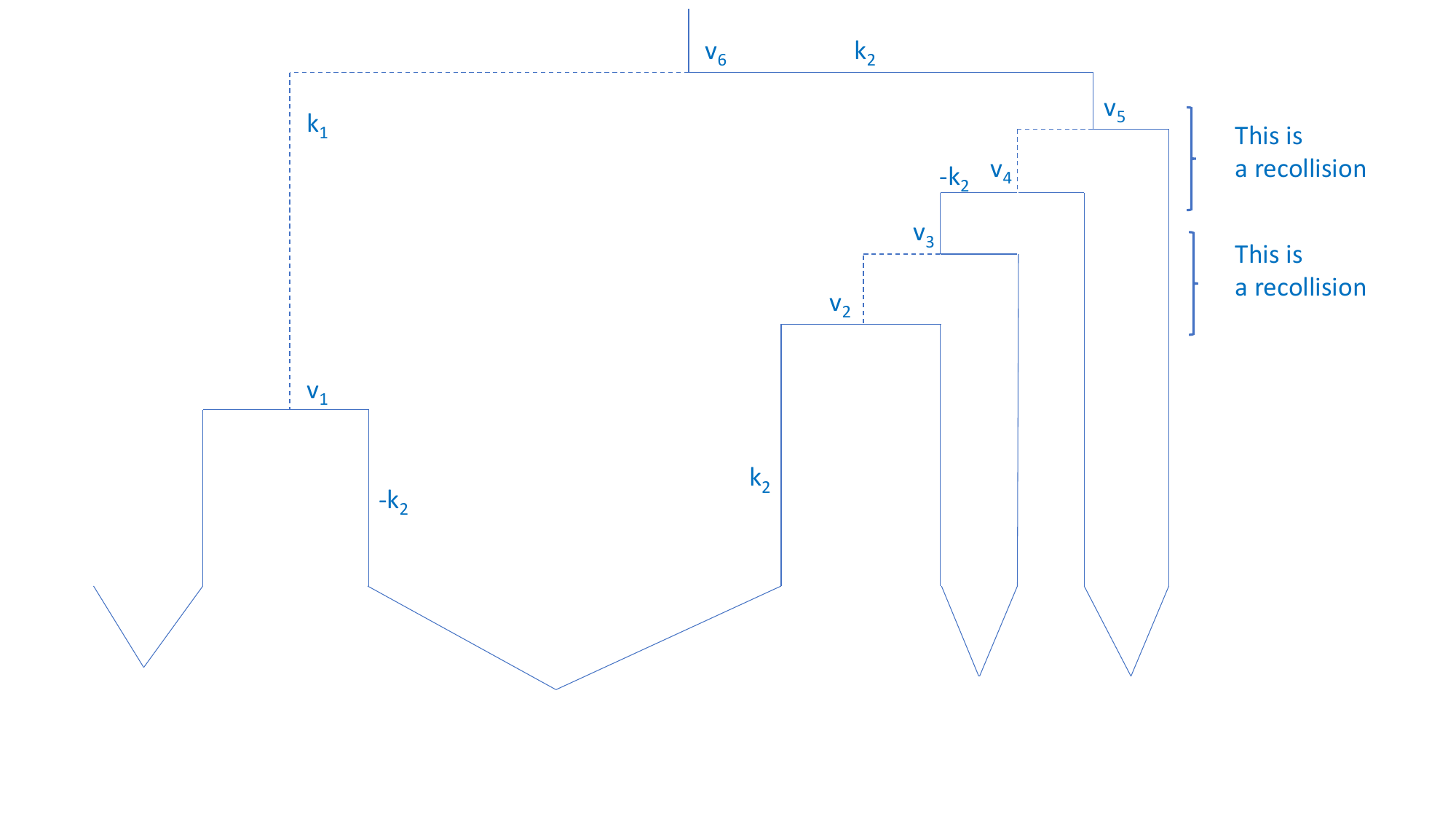}
		\caption{ In this example, after the two recollisions of $v_4,v_5$ and $v_2, v_3$, the momenta is propagated from $v_6$ to $v_2$. As a result, the edge associated with $v_1$ is $-k_2$, that cancels out the $k_2$ edge of $v_6$. Moreover, the momenta $k_1$ appear in both $v_1$ and $v_6$, but with opposite signs, thus, they cancel out each other. Finally, $\mathfrak{X}(v_1)+\mathfrak{X}(v_6)$ is independent of $k_1$. As a result, the cycle of $v_6$ is an  $IC^d_6$ long delayed recollision.}
		\label{Fig25}
	\end{figure}

	We then have
	\begin{equation}\label{Propo:LeadingGraphs:EA1}
		\begin{aligned}
			\mathfrak{X}(v_{l_1})+\mathfrak{X}(v_{l_0}) \ = \  &\sigma_{k_0}\omega(k_0)+\sigma_{k_1}\omega(k_1)+\sigma_{k_2}\omega(k_2)\\
			& \ +\ \sigma_{h_0}\omega(h_0)+\sigma_{h_1}\omega(h_1)+\sigma_{h_2}\omega(h_2)\\
			= \ & \sigma_{k_0}\omega(k_0)+\sigma_{k_1}\omega(k_1)+\sigma_{k_2}\omega(\sigma_{k_2}\sigma_{k_0}k_0-\sigma_{k_2}\sigma_{k_1}k_1)\\
			& \ + \ \sigma_{h_0}\omega(h_0)+\sigma_{h_1}\omega(\sigma_1 k_1+u_1)+\sigma_{h_2}\omega(\sigma_2 k_1+u_2).
	\end{aligned}\end{equation}
	Since the edge $k_0$ belongs to a vertex $v_j$, $j>l_1$, it is independent of $k_1$. We have already proved that $h_0$ is also independent of $k_1$. We then write
	\begin{equation}\label{Propo:LeadingGraphs:EA2}
		\begin{aligned}
			\mathfrak{X}_{l_1}+\mathfrak{X}_{l_0} 
			= \ & \sigma_{k_1}\omega(k_1)+\sigma_{k_2}\omega(\sigma_{k_2}\sigma_{k_0}k_0-\sigma_{k_2}\sigma_{k_1}k_1)\\
			& \ + \ \sigma_{h_1}\omega(\sigma_1 k_1+u_1)+\sigma_{h_2}\omega(\sigma_2 k_1+u_2) \ + \ c_{indepfreeedge},
	\end{aligned}\end{equation}
	in which the constant $c_{indepfreeedge}$ is independent of $k_1$. Therefore, the only quantity that depends on $k_1$ is $\sigma_{k_1}\omega(k_1)+\sigma_{k_2}\omega(\sigma_{k_2}\sigma_{k_0}k_0-\sigma_{k_2}\sigma_{k_1}k_1)+\sigma_{h_1}\omega(\sigma_1 k_1+u_1)+\sigma_{h_2}\omega(\sigma_2 k_1+u_2)$. By our assumption $\mathfrak{X}_{l_1}+\mathfrak{X}_{l_0}$ is not a function of $k_1$, we find $\sigma_{k_1}\omega(k_1)+\sigma_{k_2}\omega(\sigma_{k_2}\sigma_{k_0}k_0-\sigma_{k_2}\sigma_{k_1}k_1)+\sigma_{h_1}\omega(\sigma_1 k_1+u_1)+\sigma_{h_2}\omega(\sigma_2 k_1+u_2)=0$. Since the graph is non-singular, this happens when $k_1=\pm h_1$, $k_2=\pm h_2$ (or $k_2=\pm h_1$ and $k_1=\pm h_2$, by symmetry) and $\sigma_{k_1}\omega(k_1)+\sigma_{h_1}\omega(h_1)=\sigma_{k_2}\omega(k_2)+\sigma_{h_2}\omega(h_2)=0$ (or $\sigma_{k_1}\omega(k_1)+\sigma_{h_2}\omega(h_2)=\sigma_{k_2}\omega(k_2)+\sigma_{h_1}\omega(h_1)=0$, by symmetry).
	\smallskip
	
	{\bf Part 2 - The proof of (ii): There does not exist another vertex $v_j$ in the cycle of $v_{l_1}$ such that $\mathfrak{X}(v_{j})+\mathfrak{X}({v_{l_1}})$ is not a function of $k_1$.}
	\smallskip
	
	From the previous part, we know that $k_1=\pm h_1$, $k_2=\pm h_2$ (or $k_2=\pm h_1$ and $k_1=\pm h_2$, by symmetry) and $\sigma_{k_1}\omega(k_1)+\sigma_{h_1}\omega(h_1)=\sigma_{k_2}\omega(k_2)+\sigma_{h_2}\omega(h_2)=0$ (or $\sigma_{k_1}\omega(k_1)+\sigma_{h_2}\omega(h_2)=\sigma_{k_2}\omega(k_2)+\sigma_{h_1}\omega(h_1)=0$, by symmetry).
	We suppose without loss of generality that the former case happens. In Figure \ref{Fig25}, we give an example that illustrates this case. Now, we will show that $v_{l_0}$ is the only vertex in the cycle of $v_{l_1}$ that creates a delayed recollision. Suppose the contrary that we have another vertex $v_{j}$,  in the cycle of $v_{l_1}$ such that $\mathfrak{X}_{j}+\mathfrak{X}_{l_1}$ is not a function of $k_1$.  The same argument used for $v_{l_0}$ can be repeated, yielding that there are two momenta of $v_j$, denoted by $h'_1,h'_2$, with signs $\sigma_{h_1'},\sigma_{h_2'}$ such that $k_1=\pm h_1'$, $k_2=\pm h_2'$  and $\sigma_{k_1}\omega(k_1)+\sigma_{h_1'}\omega(h_1')=\sigma_{k_2}\omega(k_2)+\sigma_{h_2'}\omega(h_2')=0$.
	\begin{figure}
		\centering
		\includegraphics[width=.49\linewidth]{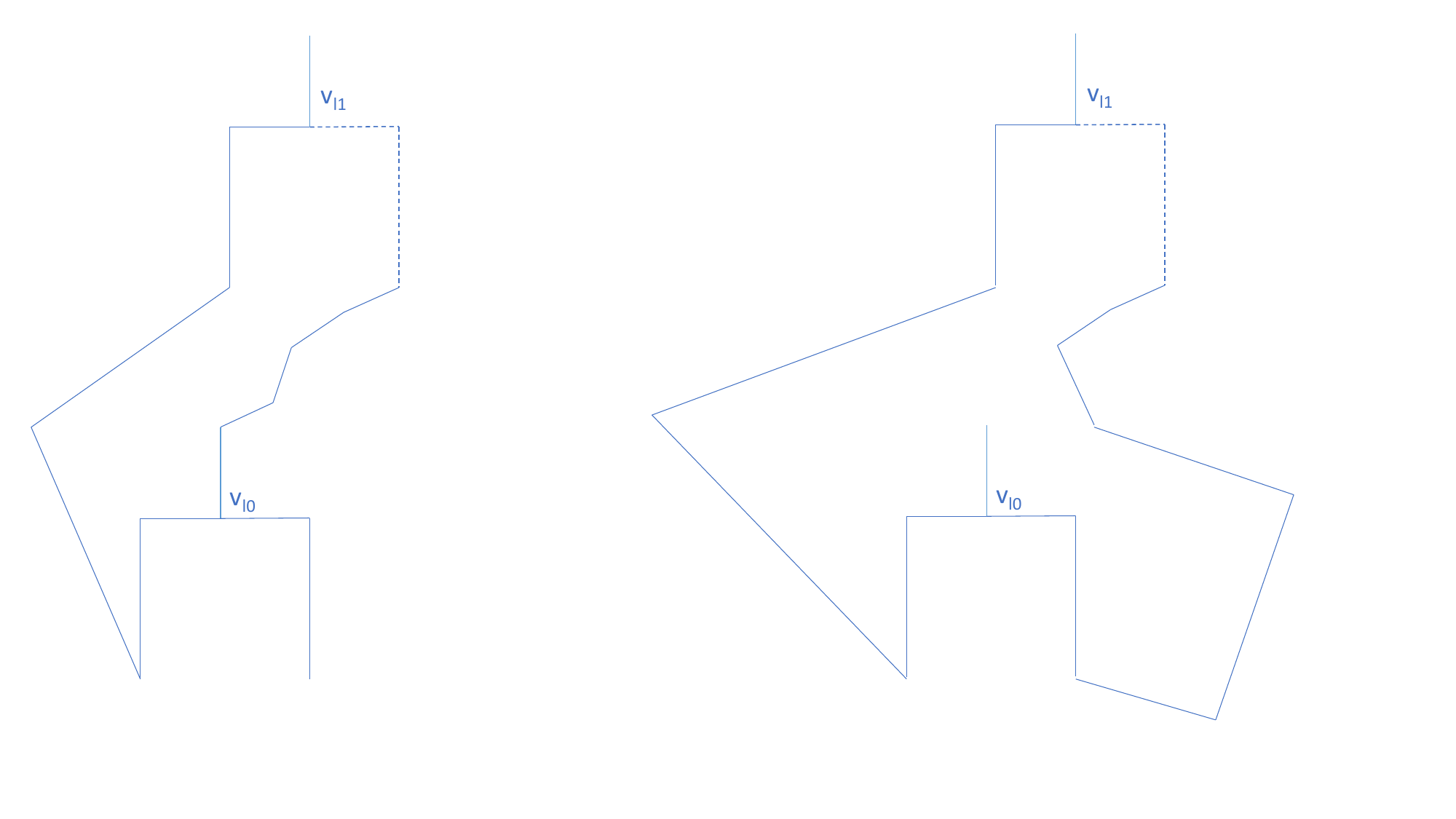}
		\caption{ In this example, the left picture corresponds to the case in which $v_{l_1}$ goes through the vertex $v_{l_0}$ by one edge in $\mathfrak{E}_-(v_{l_0})$ and one edge in $\mathfrak{E}_+(v_{l_0})$. The right picture corresponds to the case when $v_{l_1}$ goes through the vertex $v_{l_0}$ by  the two edges in $\mathfrak{E}_-(v_{l_0})$. The dashed line is the free edge of $v_{l_1}$.}
		\label{Fig24}
	\end{figure}
	We now observe that $\mathfrak{F}_{h_1'}=\mathfrak{F}_{h_1}=\mathfrak{F}_{k_1}=\{k_1\}$. We denote the edge associated to $k_1$ by $e=(v_{l_1},v')$, the edge associated to $h_1$ by $e_*=(v_{h_1},v_{h_1}')$ and the edge associated to $h_1'$ by $e_{**}=(v_{h_1'},v_{h_1'}')$, where the orientation to the virtual vertex $v_*$ are from $v_{h_1}$ to $v_{h_1}'$ and from $v_{h_1'}$ to $v_{h_1'}'.$ Then one of the two vertices $v_{h_1},v_{h_1}'$ is $v_{l_0}$ and one of the two vertices $v_{h_1'},v_{h_1'}'$ is $v_{j}$.
	By Lemma 
	\ref{Lemma:SumFreeEdges}, we have
	\begin{equation}
		\begin{aligned}\label{Propo:LeadingGraphs:EA3}
			h_1 
			\ = & \ -\sigma_{v_{h_1}}(e_*)\sum_{e'\in\mathfrak{F}}\mathbf{1}(\exists v\in e' \cap \mathfrak{P}(v_{h_1}) \mbox{ and } e' \cap \mathfrak{P}(v_{h_1})^c \neq \emptyset)\sigma_v(e')k_{e'},
	\end{aligned}\end{equation}
	and
	\begin{equation}
		\begin{aligned}\label{Propo:LeadingGraphs:EA4}
			h_1 '
			\ = & \ -\sigma_{v_{h_1'}}(e_{**})\sum_{e'\in\mathfrak{F}}\mathbf{1}(\exists v\in e' \cap \mathfrak{P}(v_{h_1'}) \mbox{ and } e' \cap \mathfrak{P}(v_{h_1'})^c \neq \emptyset)\sigma_v(e')k_{e'}.
	\end{aligned}\end{equation}
	\begin{figure}
		\centering
		\includegraphics[width=.49\linewidth]{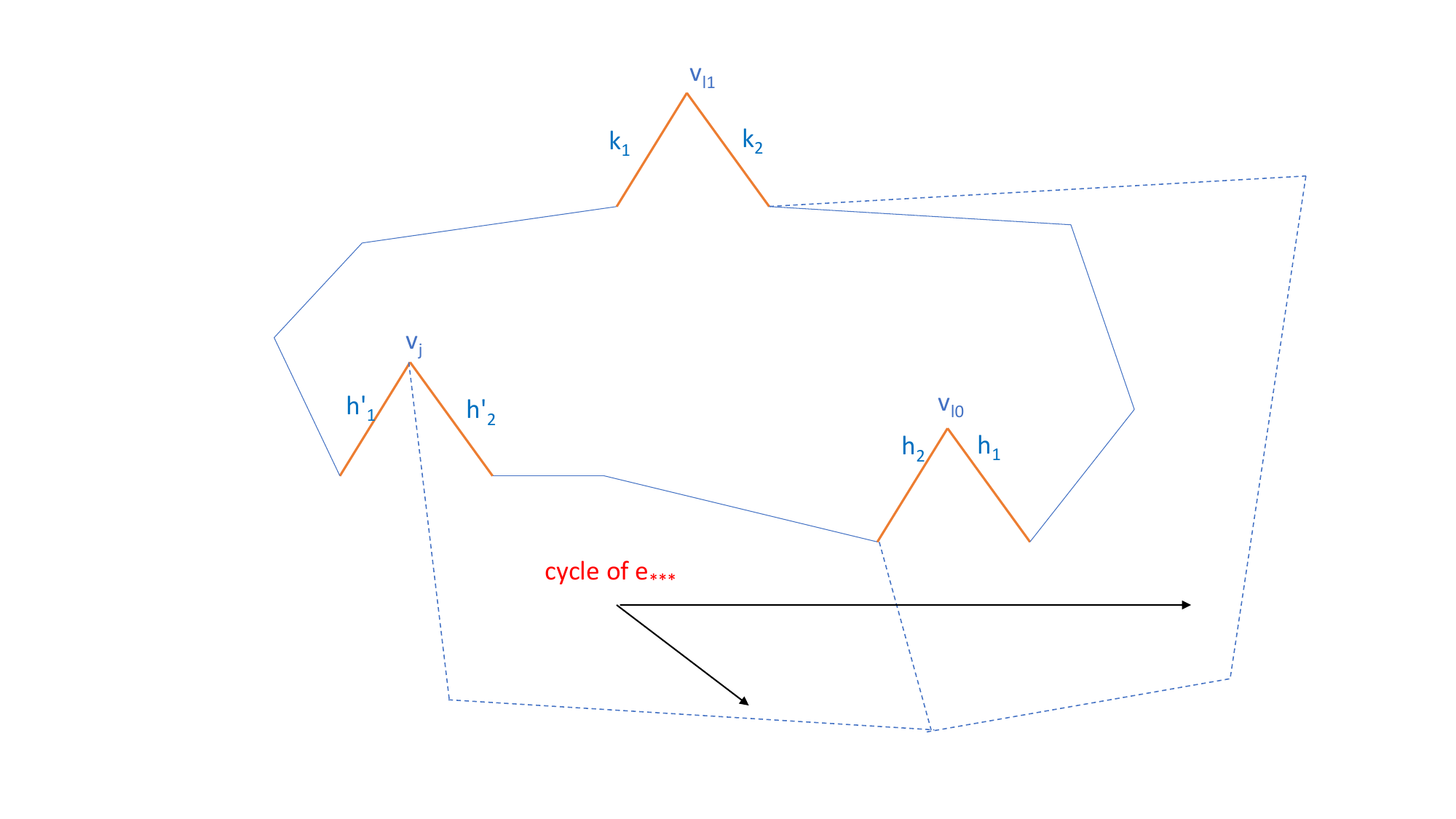}
		\caption{ In this picture, the dashed curves represent the cycle of $e_{***}$. The removal of the pair $(k_1,h_1)$ does not split the graph into two disconnected components due to the cycle of $e_{***}$.}
		\label{Fig26}
	\end{figure}
	
	It is then clear that one vertex of the free edge $e$ belongs to  $\mathfrak{P}(v_{h_1})$ and the other vertex belongs to $ \mathfrak{P}(v_{h_1})^c$, otherwise, if both vertices of $e$ belongs to $\mathfrak{P}(v_{h_1})$, the appearance of $e$ will be canceled out in the sum \eqref{Propo:LeadingGraphs:EA3}. Let $e'=(w,w')$ be another edge, in which $w\in \mathfrak{P}(v_{h_1})$ and $w'\in \mathfrak{P}(v_{h_1})^c$, then $e'$ must be a free edge. Suppose the contrary that $e'$ is not free, since $w'\notin \mathfrak{P}(v_{h_1})$, there exists a unique path from $w'$ to the virtual vertex $v_*$ and the path excludes $v_{h_1}$. This path, and the path from $v_{h_1}$ to the virtual vertex $v_*$ will form a cycle. This contradicts our construction of the free edges, since in our scheme, we break all of the cycles by removing one  edge in each cycle and call them the free edges.   Since $e'$ is free, it coincides with $e$. As a result the set of edges that connect $\mathfrak{P}(v_{h_1})$ and $\mathfrak{P}(v_{h_1})^c$, and different from $e_*$ contains only one element $e$. Similarly,   the set of edges that connect $\mathfrak{P}(v_{h_1'})$ and $\mathfrak{P}(v_{h_1'})^c$, and different from $e_{**}$ contains also one element $e$. Since both $e_*$ and $e_{**}$ belong to the cycle of $e$, we suppose, without loss of generality that there is a path  $e_*\succ e_{**}$. As  the set of edges  connecting $\mathfrak{P}(v_{h_1'})$ and $\mathfrak{P}(v_{h_1'})^c$, and different from $e_{**}$, contains only one element $e$, we deduce that all the paths to the virtual vertex $v_*$ whose components contain edges in $\mathfrak{P}(v_{h_1'})$ and $\mathfrak{P}(v_{h_1'})^c$ include $e_{**}$ as one of their edges. Similarly, all the paths to the virtual vertex $v_*$ whose components contain edges in $\mathfrak{P}(v_{h_1})$ and $\mathfrak{P}(v_{h_1})^c$ include $e_{*}$ as one of their edges. As a result, removing either the pair $\{e_*,e\}$ or the pair $\{e_{**},e\}$ will split the graph into two disconnected component, by Lemma \ref{Lemma:RemovalofTwoEdges} and the graph has a 2-separation. Now, let us consider the three edges associated to $k_2$, $h_2'$ and $h_2''$. Since the three edges $k_2$, $h_2'$ and $h_2''$ depend not only on $k_1$, they also depend on at least another free edge. We name this free edge by $e_{***}$, and denote the associated free momenta by $k_{***}$. Therefore,  for each of the momenta $k_2$, $h_2'$ and  $h_2''$, one of its vertices should belong to the cycle created by $e_{***}$.   This cycle guarantees that the removal of either the pair $\{e_*,e\}$ or the pair $\{e_{**},e\}$ does not split the full graph into two separate components (see Picture \ref{Fig26} for an illustration).  
	
	\smallskip
	
	{\bf Part 3 - The proof of (iii): Forming a recollision or a short delayed recollision.} 
	
	\smallskip
	
	Suppose that  \eqref{Lemma:VerticesLongCollisions:1}  happens, then $\mathfrak{F}_{k_1}=\mathfrak{F}_{k_1'}$, as a result, we can remove the two edges associated to $k_1,k_1'$ and the graph splits into two parts by Lemma \ref{Lemma:RemovalofTwoEdges}. Thus, we can cut the part that does not contains $v_{l_1},v_{l_0}$ out of the graph. Since $\sigma_{k_1}k_1+\sigma_{k_1'}k_1'=0$, we can joint $v_{l_1}$ to $v_{l_0}$, by gluing $k_1$ to $k_1'$ into one new cluster. By the same argument, we can also glue $k_2$ to $k_2'$ into one new cluster (see Figure \ref{Fig15}). In this case, we have either a short delayed recollision or a recollision. The case when  \eqref{Lemma:VerticesLongCollisions:2}  happens can be done by exactly the same argument.
	\begin{figure}
		\centering
		\includegraphics[width=.49\linewidth]{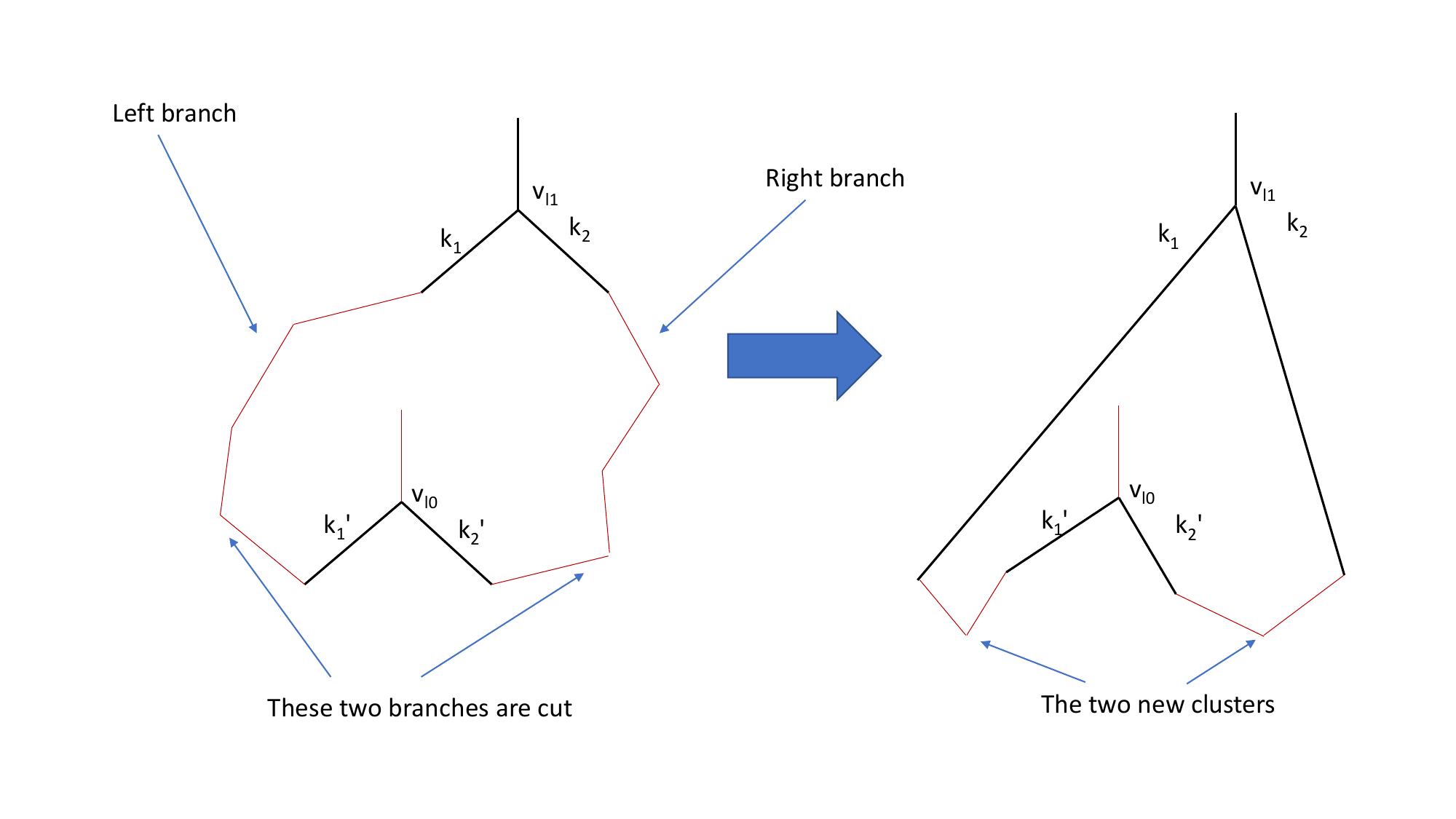}
		\caption{ In this picture, the two parts of the long collision is replaced by the two new clusters, reducing the long collision into a short collision, which is either a short delayed recollision or a recollision.}
		\label{Fig15}
	\end{figure}
	
	If \eqref{Lemma:VerticesLongCollisions:3} happens, then $\mathfrak{F}_{k_1}=\mathfrak{F}_{k_0'}$, as a result, we can remove the two edges associated to $k_1,k_0'$ and the graph splits into two parts by Lemma \ref{Lemma:RemovalofTwoEdges}. Therefore, we can cut the part that does not contains $v_{l_1},v_{l_0}$ out of the graph. Since $\sigma_{k_1}k_1+\sigma_{k_0'}k_0'=0$, we can joint $v_{l_1}$ to $v_{l_0}$, by gluing $k_1$ to $k_1'$ into one new edge $e$. In this case, $\sigma_{v_1}e=\sigma_{k_1}$ and $\sigma_{v_0}e=\sigma_{k_0'}$. By the same argument used for the case when \eqref{Lemma:VerticesLongCollisions:1}  happens, we can glue $k_2$ and $k_1'$ by a new cluster (see Figure \ref{Fig16}). We then have either a short delayed recollision or a recollision. The case when  \eqref{Lemma:VerticesLongCollisions:4}  happens is then similar.
	\begin{figure}
		\centering
		\includegraphics[width=.49\linewidth]{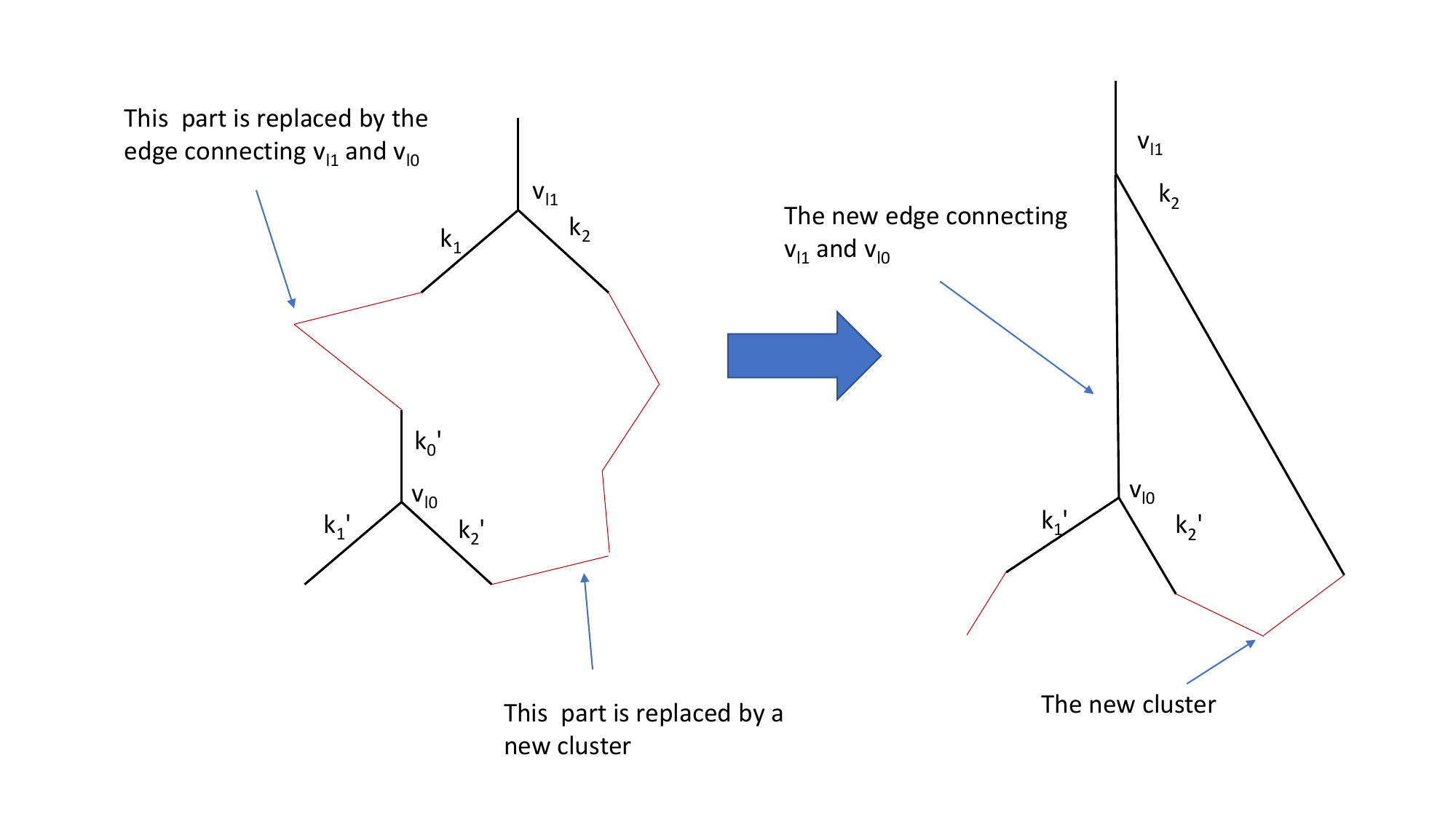}
		\caption{ In this picture, the two parts of the long collision are replaced by the new cluster and the new edge, reducing the long collision into a short collision, which is either a short delayed recollision or a recollision.}
		\label{Fig16}
	\end{figure}
\end{proof}
By Lemma \ref{Lemma:VerticesLongCollisions}, we then have a classification for long collisions. 
\begin{definition}[Classification of Long Collisions in Non-singular Graphs]\label{Def:ClassificationLongCollisions}
	Let $v_{l_1}$ be a degree-one vertex in a non-singular graph and suppose that its $\mathrm{iC}_l$ cycle is long $l>2$. Denote by $k_1$ its free momentum. 
	
	\begin{itemize}
		\item[(i)] If there exists a vertex $v_{l_0}$ in the cycle of $v_{l_1}$ such that $\mathfrak{X}(v_{l_0})+\mathfrak{X}({v_{l_1}})$ is not a function of $k_1$, and $l_1-l_0>1$, the cycle is called a ``long delayed recollision''  and denoted by $\mathrm{iC}_l^d$.
		\item[(ii)] If there does not exist a vertex $v_{l_0}$ in the cycle of $v_{l_1}$ such that $\mathfrak{X}(v_{l_0})+\mathfrak{X}({v_{l_1}})$ is not a function of $k_1$, the cycle is called an ``irreducible long collision''and denoted by $\mathrm{iC}_l^i$.
	\end{itemize}
	
\end{definition}
We have the following observation for a recollision.
\begin{lemma}\label{Lemma:RecollisionXiXi-1}
	Let $v_{i}$ be a degree-one vertex in a non-singular graph and suppose that it is either a  recollision or a delayed recollision. Denote by $k_1$ its free edges. Let $v_{j}$ be the another   vertex in the cycle. Then $\mathfrak{X}(v_{i})+\mathfrak{X}({v_{j}})$ is not a function of $k_1$.
	
\end{lemma}
\begin{proof}
	The proof follows straightforwardly from the computation of $\mathfrak{X}(v_{i})$ and $\mathfrak{X}({v_{j}})$.
\end{proof}
Combining Lemma \ref{Lemma:VerticesLongCollisions} and Lemma \ref{Lemma:RecollisionXiXi-1}, we have the following unified definition for delayed  recollisions.
\begin{definition}[Delayed  Recollisions]\label{Def:DelayAndRecollisions}
	Let $v_{l_1}$ be a degree-one vertex in a non-singular graph. Denote by $k_1$ its free edges. 
	If there exists a vertex $v_{l_0}$ in the $\mathrm{iC}_l$ cycle of $v_{l_1}$ such that $\mathfrak{X}(v_{l_0})+\mathfrak{X}({v_{l_1}})$ is not a function of $k_1$, and if the $\mathrm{iC}_l$ cycle is not a recollision, it is called a ``delayed recollision''. In this case, there are two possibilities.
	\begin{itemize}
		\item[(i)] If the cycle is long $l>2$, it is a ``long delayed recollision''
		\item[(ii)] If the cycle is short $l=2$, it is a ``short delayed recollision''.
	\end{itemize}
\end{definition}

\subsection{Different types of graphs}
Below, we give a definition for different types of graphs.
\begin{definition}[Different types of graphs]
	For a non-singular, pairing graph, we consider all of the degree-one vertices from the bottom to the top of the graph.
	
	If, from the bottom to the top, the first degree-one vertex,  which does not  correspond to an $\mathrm{iC}_2^r$ recollision or to a cycle formed by iteratively applying the recollisions (in Figure \ref{Fig37}, the four vertices $v_3,v_4,v_5,v_6$ are formed by iteratively applying two recollisions), is associated to an $\mathrm{iC}_m^i$  long irreducible  collision $(m> 2)$, we call the graph ``long irreducible''. 
	
	If, from the bottom to the top, the first degree-one vertex, which does not  correspond to a recollision or to a cycle formed by iteratively applying the recollisions, is associated to an $\mathrm{iC}_m^d$ delayed recollision $(m\ge 2)$, we call the graph ``delayed recollisional''. 
\end{definition}

Next, we introduce the concept of long irreducible degree-zero vertex and a modified cut-off function associated to it. 
\begin{definition}[Long-irreducible degree-zero vertex]
	Consider a long irreducible graph. 	We denote the first degree-one vertex,  which does not  correspond to an $\mathrm{iC}_2^r$ recollision, from the bottom to the top, by $v_i$. We denote the two edges in $\mathfrak{E}_-(v_i)$ by $k_1,k_2$ and the other ends of these two edges are $v_i'$ and $v_i''$. We follow the paths starting from $v_i \to v_i'$ and $v_i\to v_i''$ along the cycle of $v_i$. The first degree-zero vertices of each path, that do not belong to any recollision are denoted by $v_i^*$ and $v_i^{**}$. We call the higher one, between the two vertices $v_i^*$ and $v_i^{**}$ in the graph, the long-irreducible degree-zero vertex of the graph (see Figure \ref{Fig43} for an illustration).

\end{definition}

\begin{figure}
	\centering
	\includegraphics[width=.49\linewidth]{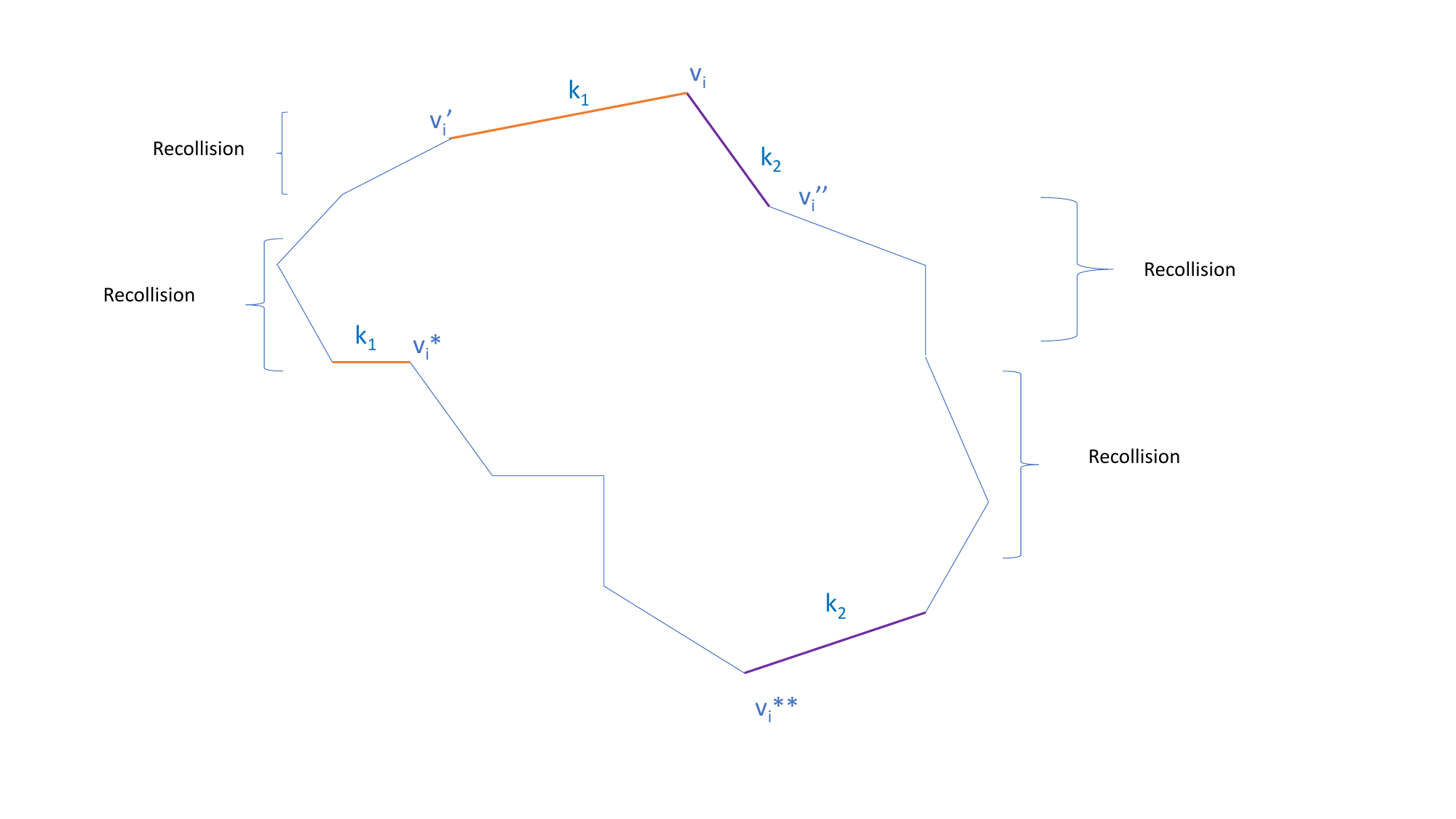}
	\caption{ In this picture, following the path from $v_i$ to $v_i'$, we find $v_i^*$, the first degree-zero vertex that does not belong to any recollision. Following the path from $v_i$ to $v_i''$, we find $v_i^{**}$, the first degree-zero vertex that does not belong to any recollision. As $v_i^*$ is placed higher than $v_{i}^{**}$, $v_i^*$  is the long-irreducible degree-zero vertex of the graph.}
	\label{Fig43}
\end{figure}
We then have the lemma. 
\begin{lemma}
	\label{Lemma:LIDegreeZeroVertex}
	Consider a long irreducible   graph. 	We denote the first degree-one vertex,  which does not  correspond to an $\mathrm{iC}_2^r$ recollision, from the bottom to the top, by $v_i$. We denote the two edges in $\mathfrak{E}_-(v_i)$ by $k_1,k_2$ and the other ends of these two edges are $v_i'$ and $v_i''$. Denote by $v_j$ the  long-irreducible degree-zero vertex of the graph. Suppose that $v_j$ belongs to the path from $v_i$ to $v_i'$. Then the edge that is attached to $v_j$, and belongs to this path, has momentum $k_1$. Moreover, the cut-off function associated to $v_j$ is of type $\Phi_{1,j}$. 
\end{lemma}
\begin{proof}
	As from $v_i'$ to $v_j$, we only have recollisions,  the momentum $k_1$ is propagated to the edge that is attached to $v_j$ and belongs to this path. Moreover, since $v_j$ is placed higher than the other degree-zero vertex emerging from the path from  $v_i$ to $v_i''$, the cut-off function is of type $\Phi_{1,j}$.
\end{proof}

\begin{lemma}
	\label{Lemma:Triplet}
Consider a graph with $n$ time slices, in which the number of vertices at the bottom $\{k_{0,1},\cdots,k_{0,n+2}\}$ is even. Suppose that at each time slice, the sets of the associated momenta $\{k_{i,1},\cdots,k_{i,n+2-i}\}$ are admissible (see Definition \ref{def:distinct}). Then the graph is pairing.
\end{lemma}
\begin{proof}
	Suppose that $i_0$ is the first time slice from the top to the bottom at which the set of associated momenta $\{k_{i_0,1},\cdots,k_{i_0,n+2-i_0}\}$ has  a triple $k_{i_0,j'},k_{i_0,j''},k_{i_0,j'''}$ such that $k_{i_0,j'}+k_{i_0,j''}+k_{i_0,j'''}=0$ appear, while the set $\{k_{i_0,1},\cdots,k_{i_0,n+2-i_0}\}\backslash \{k_{i_0,j'},k_{i_0,j''},k_{i_0,j'''}\}$ contains only pairings. As the splitting happens at $k_{i_0,\rho_{i_0}},k_{i_0-1,\rho_{i_0}},k_{i_0-1,\rho_{i_0}+1},$ there are two cases: either $\rho_{i_0}$ coincides with one of the indices $j',j'',j'''$ or $\rho_{i_0}$ does not coincide  with one of the indices $j',j'',j'''$. Let us consider the first case when $\rho_{i_0}$ coincides with one of the indices $j',j'',j'''$, without loss of generality, we suppose $j'=\rho_{i_0}$. As $k_{i_0,j'}+k_{i_0,j''}+k_{i_0,j'''}=0$, we deduce  $k_{i_0-1,\rho_{i_0}}+k_{i_0-1,\rho_{i_0}+1}+k_{i_0,j''}+k_{i_0,j'''}=0$. As the set $\{k_{i_0-1,1},\cdots,k_{i_0-1,n+3-i_0}\}$ is admissible, we deduce that $\{k_{i_0-1,\rho_{i_0}},k_{i_0-1,\rho_{i_0}+1},k_{i_0,j''},k_{i_0,j'''}\}$ contains  only pairings. We now consider the second case that $\rho_{i_0}$ does not coincide  with one of the indices $j',j'',j'''$, then $k_{i_0,\rho_{i_0}}$  must belong to a paring, we suppose without loss generality that this pairing is $k_{i_0,\rho_{i_0}},k_{i_0,\rho_{i_0}+1}$, yielding $k_{i_0-1,\rho_{i_0}}+k_{i_0-1,\rho_{i_0}+1}+k_{i_0-1,\rho_{i_0}+2} =0$. We suppose that $k_{i_0,j'},k_{i_0,j''},k_{i_0,j'''}$ propagates to  $k_{i_0-1,j'_0},k_{i_0-1,j''_0},k_{i_0-1,j'''_0}$ in the lower time slice $i_0-1$. As $\{k_{i_0-1,j'_0},k_{i_0-1,j''_0},k_{i_0-1,j'''_0}\}$ is different from $\{k_{i_0-1,\rho_{i_0}},k_{i_0-1,\rho_{i_0}+1},k_{i_0-1,\rho_{i_0}+2}\}$, and as the set of momenta in the $i_0-1$ time slice is also admissible, $\{k_{i_0-1,j'_0},k_{i_0-1,j''_0},k_{i_0-1,j'''_0}\}$ and $\{k_{i_0-1,\rho_{i_0}},k_{i_0-1,\rho_{i_0}+1},k_{i_0-1,\rho_{i_0}+2}\}$ create another $3$ pairings. Applying this argument iteratively, we conclude that if the set $\{k_{i,1},\cdots,k_{i,n+2-i}\}$ contains one triplet $k_{i,j'}+k_{i,j''}+k_{i,j'''}=0$ (while the left-over set contains only pairings), then the set of the next time slice $\{k_{i-1,1},\cdots,k_{i-1,n+3-i}\}$ contains only pairings. And next, the set $\{k_{i-2,1},\cdots,k_{i-2,n+4-i}\}$ contains precisely one triplet $k_{i-2,l'}+k_{i-2,l''}+k_{i-2,l'''}=0$ and the left-over set contains only pairing. Finally, the set  $\{k_{0,1},\cdots,k_{0,n+2}\}$ contains only pairings.
\end{proof}
\subsection{The most important graphs: $\mathrm{iCL}_2$  ladder graphs (leading diagrams)}\label{LeadingDiagram}
In the next sections of the paper we will show that graphs with long collisions and delayed recollisions are negligible, as thus, the most important  contribution comes from graphs with only $\mathrm{iC}_2^r$ recollisions. The following lemma shows the form of a graph, whose cycles are only $\mathrm{iC}_2^r$ recollisions. 
\begin{definition}[$\mathrm{iCL}_2$ Recollision Ladder Graphs/Leading Diagrams]\label{Def:Ladder}
	A graph that can be obtained by iteratively adding $M$  recollisions ($\mathrm{iC}_2^r$ cycles) is an ``$\mathrm{iCL}_2$ recollision ladder graph'' (see Figure \ref{Fig35}). We call the $\mathrm{iC}_2^r$ cycles the ``skeletons'' of the $\mathrm{iCL}_2$ ladder graph. 
	\begin{figure}
		\centering
		\includegraphics[width=.49\linewidth]{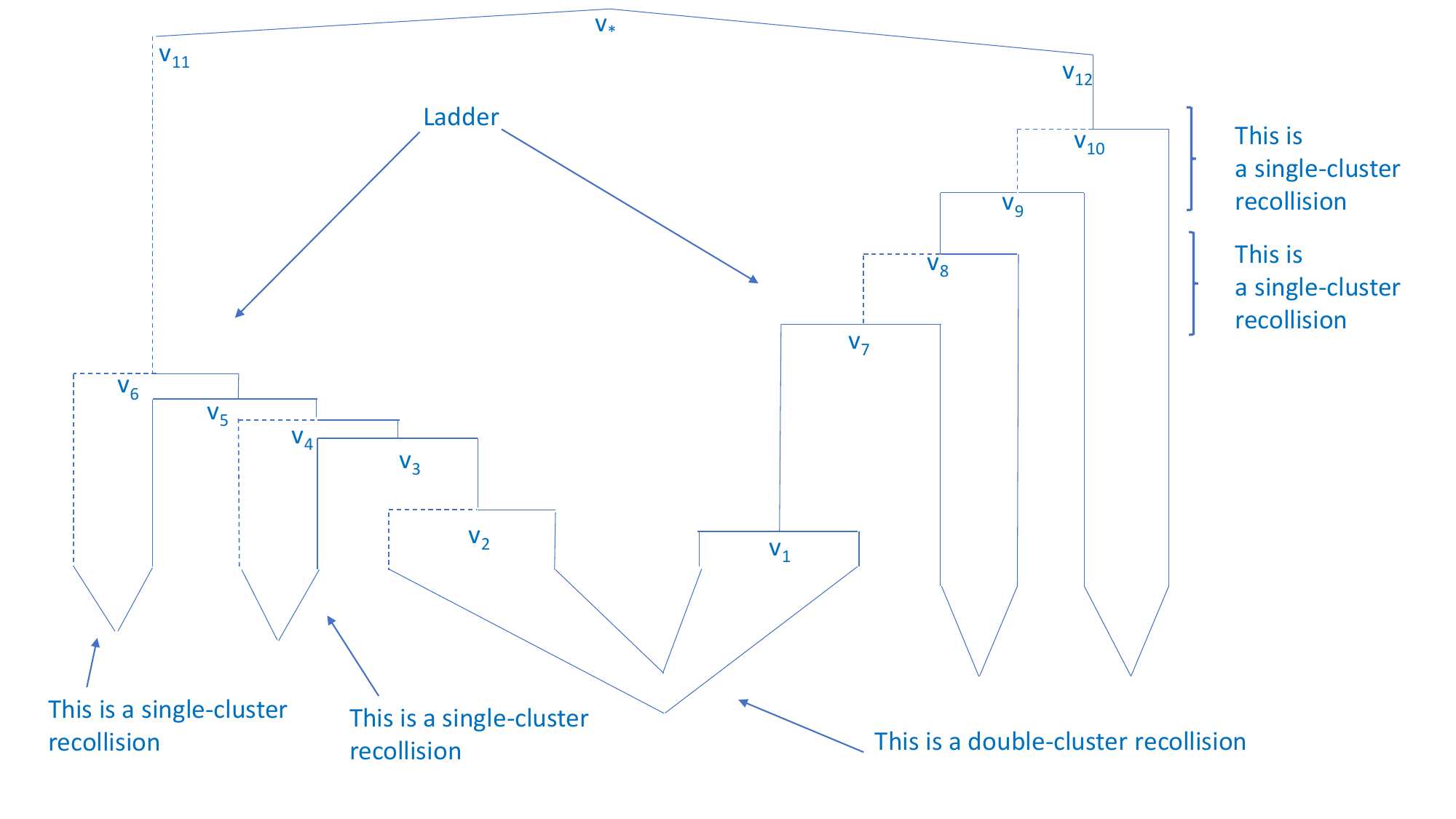}
		\caption{This is an example of a ladder. All of the cycles are recollisions. There are $4$ single-cluster recollisions and $1$ double-cluster recollision.}
		\label{Fig35}
	\end{figure}	
\end{definition}
\begin{remark}\label{Remark:CLn}
	We remark that prism graphs  are $\mathrm{CL}_m$  circular ladder graphs that have one of the prisms as its skeleton. Note that a $\mathrm{CL}_m$ circular ladder graph can be obtained by the Cartesian product of a cycle and an edge.  Our $\mathrm{iCL}_2$ ladder graphs could be understood, in  some sense, as an extension of those graphs using $\mathrm{iC}_2^r$ cycles as the skeletons. However, $\mathrm{iCL}_2$ ladder graphs are obtained by iteratively adding the $\mathrm{iC}_2^r$ cycles, which could be understood as the ``convolutions'' of the skeletons.  
\end{remark}
\begin{lemma}\label{Lemma:Recollisions}
	Given a non-singular pairing graph, suppose that all of its cycles are $\mathrm{iC}_2^r$ recollisions. Then the graph is a $\mathrm{iCL}_2$ ladder. Moreover, the following identity also holds true
	\begin{equation}
		\label{Lemma:Recollisions:1}\begin{aligned}
			&\mathfrak{X}(\sigma_{i,\rho_i},  k_{i,\rho_i},\sigma_{i-1,\rho_i}, k_{i-1,\rho_i},\sigma_{i-1,\rho_i+1}, k_{i-1,\rho_i+1}) \\
			& + \mathfrak{X}(\sigma_{i+1,\rho_{i+1}},  k_{i+1,\rho_{i+1}},\sigma_{i,\rho_{i+1}}, k_{i,\rho_{i+1}},\sigma_{i,\rho_{i+1}}, k_{i,\rho_{i+1}+1}) \ = \ 0,
		\end{aligned}
	\end{equation}
	if $i$ is odd.
\end{lemma}

\begin{proof}

	Suppose the graph has $n$ interacting vertices. Since it is pairing, $n=2m$ is an even number and there are in total $|S|=m+1$ cluster vertices. Using Lemma \ref{Lemma:NumberDegreeZeroOne}, the number of degree-one vertices  is computed by $n+1-|S|=2m+1-m-1=m$. And the number of the degree-zero vertices is $|S|-1=m$. 
	
	Since the numbers of degree-one vertices and degree-zero vertices are equal, the number of degree-zero vertices that do not belong to any cycle is zero. This can be seen as follows. We denote by $l$ the number of the degree-zero vertices that do not belong to any cycle. Since there are in total $m$ degree-one vertices, there are in total $m$ cycles. By the hypothesis, all of the cycles are  $\mathrm{iC}_2^r$ recollisions. Thus, each cycle has only two vertices, one vertex is of degree $1$ and the other one is of degree $0$. Therefore, there are $m$ degree-one vertices and $m$ degree-zero vertices that belong to the cycles of the graph. The number of interacting vertices of the graph is then $2m+l=2m$. Therefore, $l=0$. This shows that each degree-zero vertex belong to a cycle of the graph. 
	
	As a consequence, the graph contains $m$ cycles, each of which has  one degree-one vertex $v_{i}$ at the top and and the next one, $v_{i-1}$ is a degree-zero vertex. 
	
	Now, we can prove by induction in $m$ that the graph can be obtained by iteratively adding the  $\mathrm{iC}_2^r$ recollisions. 
	
	If $m=1$, we only have one cycle in the graph, meaning that the cycle is an  $\mathrm{iC}_2^r$ recollision. Suppose that the claim is true for $m=n$, we will prove that the claim is true for $m=n+1$. Let us consider the  $\mathrm{iC}_l$ cycle, whose interacting size is the smallest among all the cycles of the graph, then the cycle is associated to the degree-one vertex $v_i$. The next vertex in the cycle is $v_{i-1}$ and this is a degree-zero vertex. Suppose that inside the cycle, there is another vertex $v_j$. If this is a degree-zero vertex, then $v_j$ belongs to the cycle of the degree-one vertex $v_{j+1}$, according to the previous argument. Moreover, $v_{j+1}$ and its cycle should also belong to the cycle of $v_i$. Thus the cycle of $v_{j+1}$ has a smaller size than the cycle of $v_i$. This is a contradiction. If $v_j$ is a degree-one vertex, then the cycle of $v_j$ also belongs to the cycle of $v_i$, meaning that the size of the cycle of $v_j$ is smaller than that of $v_i$. This also leads to another contradiction. Therefore, the cycle of $v_i$ has only two vertices $v_i$ and $v_{i-1}$, and it is a  $\mathrm{iC}_2^r$ recollision. We can thus remove this cycle (see Figure \ref{Fig37}), then the graph becomes a graph with $n$ cycles of the type  $\mathrm{iC}_2^r$ and the induction assumption can then be applied. Therefore, the graph can be obtained by iteratively adding $n+1$ recollisions. 
	
	By induction, the graph can be obtained  by iteratively adding $m$ recollisions.
	\begin{figure}
		\centering
		\includegraphics[width=.49\linewidth]{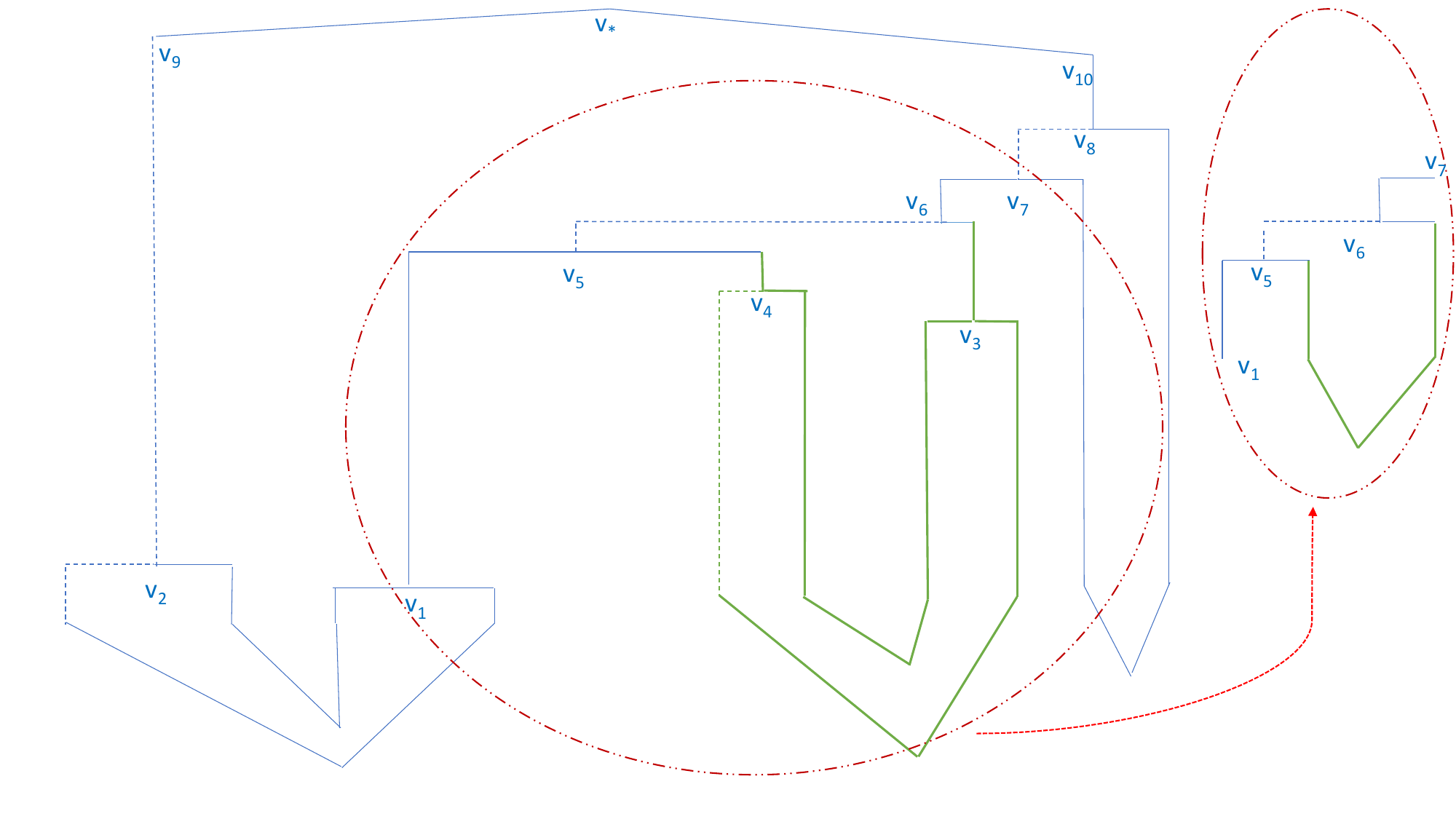}
		\caption{In this graph, the removal of the $\mathrm{iC}_2$ cycle of $\{v_3,v_4\}$, which is a two-cluster recollision, makes the $\mathrm{iC}_4$  cycle of $\{v_5,v_6\}$ become an $\mathrm{iC}_2$  one-cluster recollision. The four vertices $v_3,v_4,v_5,v_6$ are formed by iteratively applied two recollisions. }
		\label{Fig37}
	\end{figure}

	The identity \eqref{Lemma:Recollisions:1} follows straightforwardly from \eqref{Lemma:PhaseAtSlice0:1} and the fact that the sum of the phases of a recollision is $0$.
\end{proof}

\section{General graph estimates}
\label{Sec:GeneralGraph}

\subsection{Estimates of ${Q}^3$} In this subsection, we will provide estimates on $Q^3$, which is defined in Proposition \ref{Proposition:Ampl2} .

\begin{proposition}[The first   estimate on ${Q}^3$]\label{lemma:BasicGEstimate1}
	Let $	\tilde{Q}^3$ be a term in the sum of $Q^{3,nonpair}$. There are constants $\mathfrak{C}_{3,1}$ such that for $1\le n$ and $t=\tau\lambda^{-2}>0$ 
	\begin{equation}
		\begin{aligned} \label{eq:BasicGEstimate3:1}
			\tilde{Q}^3:= & \	\lambda^{n}\mathbf{1}({\sigma_{n,1}=-1})\mathbf{1}({\sigma_{n,2}=1})\Big[  \\&\times \sum_{\substack{\bar\sigma\in \{\pm1\}^{\mathcal{I}_n},\\ \sigma_{i,\rho_i}+\sigma_{i-1,\rho_i}+\sigma_{i-1,\rho_i+1}\ne \pm3,
					\\ \sigma_{i-1,\rho_i}\sigma_{i-1,\rho_i+1}= 1}}\int_0^t\mathrm{d}s_0 e^{{\bf i}s_0\vartheta_0}\int_{(\Lambda^*)^{\mathcal{I}_n}}  \mathrm{d}\bar{k}\Delta_{n,\rho}(\bar{k},\bar\sigma)\left\langle\prod_{i=1}^{n+2}\alpha(k_{0,i},\sigma_{0,i})\right\rangle_{s_0} h^d  \\
			&\times \Phi_{0,1}(\sigma_{0,\rho_1},  k_{0,\rho_1},\sigma_{0,\rho_2},  k_{0,\rho_2})\sigma_{i,\rho_i}\mathcal{M}( k_{1,\rho_1}, k_{0,\rho_1}, k_{0,\rho_1+1})\\
			&\times \prod_{i=2}^n\Big[\sigma_{i,\rho_i}\mathcal{M}( k_{i,\rho_i}, k_{i-1,\rho_i}, k_{i-1,\rho_i+1})\Phi_{1,i}(\sigma_{i-1,\rho_i},  k_{i-1,\rho_i},\sigma_{i-1,\rho_i+1},  k_{i-1,\rho_i+1})\Big] \\
			&\times \int_{(\mathbb{R}_+)^{\{1,\cdots,n\}}}\mathrm{d}\bar{s} \delta\left(t-\sum_{i=0}^ns_i\right)\prod_{i=1}^{n}e^{-s_i \varsigma_{n-i} }\prod_{i=0}^{n}e^{-s_i \tau_{i}}\prod_{i=0}^{n-1} e^{-{\bf i}t_i(s)\mathfrak{X}_i}(1-\mho)\Big],
	\end{aligned}\end{equation}
and
	
	\begin{equation}
		\begin{aligned} \label{eq:BasicGEstimate3}
			\lim_{\lambda\to0} \limsup_{D\to\infty}\Big|\tilde{Q}^3\Big|
			\ 
			= & 0.
	\end{aligned}\end{equation}
 The time integration is defined to be $$\int_{(\mathbb{R}_+)^{\{1,\cdots,n\}}}\mathrm{d}\bar{s} \delta\left(t-\sum_{i=0}^ns_i\right)= \int_{(\mathbb{R}_+)^{\{1,\cdots,n\}}}\mathrm{d}{s}_1\cdots \mathrm{d}{s}_n \delta\left(t-\sum_{i=0}^ns_i\right).$$

	If the graph is pairing and singular then  the right hand side of \eqref{eq:BasicGEstimate3}  becomes $0$. 
\end{proposition}

\begin{proof}
	We divide the proof into several steps.
	
	{\bf Step 1: A Priori Estimates.} 
	As the set $\{(k_{0,1}\sigma_{0,1},\cdots,k_{0,n+2}\sigma_{0,n+2})\}$ is not admissible, as a result, the vector $v_o=\left(\mathbf{1}_{k_{0,j}\sigma_{0,j}}\right)_{j=1}^{n-i+2}$, which can be seen as a vector on  $\mathbb{Z}^{|\Lambda_*^+|}$ after appropriate arrangements,  can be decompose as $v_o=v_o'+ v_o''$ where $v_o'\in \mathcal{V}_\mathscr{E}^1$ and $v_o''\in \mathcal{V}_\mathscr{E}^2$ (see \eqref{VectorSpace2}). We compute
	\begin{equation}
		\label{eq:BasicGEstimate3:A1}
		v_o\mathscr{E}v_o^T\  = \    \sum_{i=1}^{\mathscr{M}_\mathscr{E}} \mho_iv_oe_ie_i^Tv_o^T \  = \ \sum_{i=1}^{\mathscr{M}_\mathscr{E}} \mho_i(v_o'+ v_o'')e_ie_i^T(v_o'+ v_o'')^T  \  = \ \sum_{i=1}^{\mathscr{M}_\mathscr{E}} \mho_i v_o''e_ie_i^T{ v_o''}^T \gtrsim \mathscr{C}_\mho.
	\end{equation}	
	Recalling  
		$
		\tau_0   =   \mathscr{T}_{\left(\mathbf{1}_{k_{0,j}\sigma_{0,j}}\right)_{j=1}^{n-i+2}} $, we then have $	\tau_0\ge \mathscr{C}_\mho$ (see Definition \ref{Def:PhaseRegulator}.)


		{\bf Step 2: Time Estimates.}
	Next, we estimate the time integral. By the definition of $\Delta_{n,\rho}$, there
	are in total $n$ interacting vertices associated to this quantity.  There are also $n$ time slices due to the delta function $\delta\left(t-\sum_{i=0}^ns_i\right)$. We then split the set $\{0,1,\cdots,n\}$ into two sets. One of which is the set of time slice with indices $1\le l<n$ such that $v_{l+1}$ is a degree one vertex.  We denote this set by $I$.  The other set then contains all of the other time slices . We denote this set by $I''$. 
	
	Introducing  the new time slice $s_{n+1}=\bar{s}_{n+1}\lambda^{-2}$, we express the time integration as
	
	\begin{equation}
		\begin{aligned}
			\label{eq:Aestimate0:1}
			&	\int_{{(\mathbb{R}_+)^{{\{0,\cdots,n\}}}}}\mathrm{d}\bar{s}\delta\left(t-\sum_{i=0}^{n}s_i\right) \ = \ \int_{(\mathbb{R}_+)^{I'}}\mathrm{d}\bar{s}'\delta\left(t-\sum_{i\in I'}s_i\right)\int_{(\mathbb{R}_+)^{I}}\mathrm{d}\bar{s}''\delta\left(s_{n+1}-\sum_{i\in I}s_i\right),
		\end{aligned}
	\end{equation}
	where $I'=\{n+1\}\cup I''$,  $\bar{s}'$ denotes the vector whose components are $s_i$ with $i\in I'$,  $\bar{s}''$ denotes the vector whose components are $s_i$ with $i\in I$.

	This means we can write $\tilde{Q}^3(\tau)$ as $\int_0^\tau\mathrm{d}\bar{s}_{n+1} $ $F_1(\tau-\bar{s}_{n+1})\lambda^{-2} F_2(s_{n+1})$, where $F_1(\tau-\bar{s}_{n+1})$ involves  $\delta\left(t-\sum_{i\in I'}s_i\right)$ and all of terms containing the time slices $s_i$ with $i\in I'$ and $F_2(s_{n+1})$ involves $\delta\left(s_{n+1}-\sum_{i\in I}s_i\right),$ and all of terms containing the time slices with $i\in I$. The quantity $F_1$ then has the form

	\begin{equation}
		\begin{aligned}\label{eq:Aestimate0:1:A}
			F_1(\tau-\bar{s}_{n+1}) \ = \ 	(\mathbf{i} \lambda)^n\int_{(\mathbb{R}_+)^{{I}''}}\mathrm{d}\tilde{s}'\delta\left(t-\sum_{i\in {I}''}s_i-s_{n+1}\right)\prod_{i\in I''}e^{-{\bf i}s_i\vartheta_i},
	\end{aligned}\end{equation}
where $\tilde{s}'$ denotes the vectors whose components are included in $(\mathbb{R}_+)^{{I}''}$. This means that $F_1$ involves all time slices of all the degree zero vertices and $F_2$ involves all the other time slices.
	
We now estimate $F_1$. 	  By using the identity,
	\begin{equation}
		\label{eq:Aestimate1a}
		\int_{\mathbb{R}_+^m}\mathrm{d}\bar{s}\delta\left(t-\sum_{i=1}^m s_i\right) \ = \ \frac{t^{m-1}}{(m-1)!},
	\end{equation}
	we then obtain
 
	\begin{equation}
		\begin{aligned}\label{eq:Aestimate2a}
	& 	 \|F_1(\tau-\bar{s}_{n+1})\|_\infty  \
			\
			 \le\ \lambda^n \frac{t^{\tilde{n}}}{\tilde{n}!}, 
	\end{aligned}\end{equation}
	in which  $\tilde{n}=|{I}''|-1.$ Therefore
		\begin{equation}
		\begin{aligned}\label{eq:Aestimate0:1:A:3}
			&\int_0^\tau\mathrm{d}\bar{s}_{n+1}    	 |F_1(\tau-\bar{s}_{n+1})|\lambda^{-2} |F_2(s_{n+1})|\ 
		\le \	  \lambda^{n-2\tilde{n}}\|F_1 \|_{L^\infty}\int_0^\tau\mathrm{d}\bar{s}_{n+1}    	 |F_2( {s}_{n+1})|\lambda^{-2}.
	\end{aligned}\end{equation}

We now study the number of possibilities to assign signs to $\bar{\sigma}$. First, at the bottom, there are $n+2$ momenta $k_{0,1},\cdots, k_{0,n+2}$, and we have at most $2^{n+2}$ ways to choose $\sigma_{0,1},\cdots, \sigma_{0,n+2}$. For each choice of the signs for the zero time slice $\sigma_{0,1},\cdots, \sigma_{0,n+2}$, there are at most two choices for the signs of the momentum which is split in the time slice 1.  Continuing this counting, for each of the time slice $i$, $i=1,\cdots, n$, there are at most 2 choices. In total, we have at most $2^{n+2}2^n=4^{n+1}$ choices of $\bar\sigma$.

As all of the time slices that are associated to the degree-zero vertices have already been removed by \eqref{Propo:ExpectationEqui:2:1}, we will now develop a way to integrate out the momenta and the time slices associated to the degree-one vertices.
The process of integrating the momenta can be carried out as follows. We integrate all of the free momenta from the bottom to the top, starting from time slice $0$. 	This matches  the direction  we used to construct the free momenta. Whenever we meet a degree one vertex, we use Lemma \ref{lemma:degree1vertex}  to integrate the associated time slice as well. Each time, we get a bound with an extra factor of $\langle \ln|\lambda| \rangle^{2+c\eth}$, for some $c>0$. We use the following rough estimate
	\begin{equation}\label{Propo:ExpectationEqui:2:1}
	\begin{aligned}
		 \sup_{k_1,\cdots,k_n\in \Lambda^*}	[\tilde\omega(k_1)\cdots\tilde\omega(k_{n+2})]^\frac12	\left \langle |a_{k_{1}}|^2\cdots|a_{k_{n+2}}|^2 \right\rangle_s \
		\lesssim  &\ \	h^{-d(n+2)}\mathcal{C}_o^{n+2},
	\end{aligned}
\end{equation}
instead of \eqref{Propo:ExpectationEqui:2:1}, that can be proved by the same proof. This means we have the point-wise bound at the bottom of the graph

\begin{equation}
	\label{MomentsFeynman:1}\begin{aligned}
		&	\left|\left\langle\prod_{i=1}^{n+2}a(k_{0,i},\sigma_{0,i})(s_0)\right\rangle\right|\prod_{i=0}^{n}e^{-s_i \tau_{i}}\
		\lesssim	& (\mathcal C_oh^{-d})^{n+2}e^{-s_0\mathscr{C}_\mho}. 
	\end{aligned}
\end{equation}

Therefore, we can bound $\int_0^\tau\mathrm{d}\bar{s}_{n+1}    	 \lambda^{-2} |F_2(s_{n+1})|$ by
	\begin{equation}
		\begin{aligned}\label{eq:Aestimate1}  \int_0^\tau\mathrm{d}\bar{s}_{n+1}    	 |F_2( {s}_{n+1})|\lambda^{-2}\ 
				\le	& \left\|\int_{(\mathbb{R}_+)^{I}}\mathrm{d}\bar{s}''\delta\left(s_{n+1}-\sum_{i\in I}s_i\right)\prod_{i\in I}e^{-{\bf i}s_i \vartheta_i}\mathfrak{W}^*  \right\|_{L^\infty}\\
		 	\le	 &   	 \left\|\int_{(\mathbb{R}_+)^{I}}\mathrm{d}\bar{s}''\delta\left(s_{n+1}-\sum_{i\in I}s_i\right)\prod_{i\in I}^{n}e^{-{\bf i}s_i \mathrm{Re}\vartheta_i}\prod_{i\in I}e^{s_i \mathrm{Im}\vartheta_i}\mathfrak{W}^*e^{-\sum_{i\in I}s_i \varsigma_{n-i}} \right\|_{L^\infty}\\
			\le  & \left\|\int_{(\mathbb{R}_+)^{I}}\mathrm{d}\bar{s}''\delta\left(s_{n+1}-\sum_{i\in I}s_i\right)\mathfrak{W}^*e^{-\sum_{i\in I}s_i \varsigma_{n-i}} \right\|_{L^\infty}\\
			\le &   ( h^{-d})^{n+1}e^{-s_0\mathscr{C}_\mho}\mathfrak{C}_{3,1}^{n},
		\end{aligned}
	\end{equation}
where $\mathfrak{W}^*$ contains all of the momenta integrals of $F_2$,
for some constant $\mathfrak{C}_{3,1}>0$, which depends on powers of $\langle \ln|\lambda| \rangle^{2+c\eth}$, $\lambda^{-1}$ as well as the constant $4^{n+1}$ discussed above.
We therefore have the following  estimate
\begin{equation}
	\begin{aligned}\label{eq:Aestimate0:1:A:4}
		 \int_0^\tau\mathrm{d}\bar{s}_{n+1}      |F_1(\tau-\bar{s}_{n+1})|\lambda^{-2} |F_2(s_{n+1})|\ 
		\le \	& \lambda^{n-2\tilde{n}} ( h^{-d})^{n+1}e^{-s_0\mathscr{C}_\mho}\mathfrak{C}_{3,1}^{n}.
\end{aligned}\end{equation}
 
which vanishes in the limit  $D\to\infty$. 

		If the graph is pairing and singular then, one of the kernel $\mathcal{M}$ vanishes and   the right hand side of \eqref{eq:BasicGEstimate3} then becomes $0$.

\end{proof}

\begin{proposition}[The main estimate on   ${Q}^3$]\label{lemma:Q1FinalEstimate}   There is a constant $\mathfrak{C}_{Q_{3,2}}>0$ such that for $t=\tau\lambda^{-2}>0$, we have

	\begin{equation}
		\begin{aligned} \label{lemma:Q1FinalEstimate:1}
			&	\lim_{\lambda\to0}\limsup_{D\to\infty}	\Big|Q^{3,nonpair}\Big|\ = \ 0
		,
	\end{aligned}\end{equation}
	where the same notations an in  Proposition  \ref{lemma:BasicGEstimate1} have been used.
\end{proposition}
\begin{proof}
    
We observe that in each graph, the vector $\{\rho_i\}_{i=1}^n$ encodes the place where the splitting happens, as defined in \eqref{Def:Rho}.
	We now compute  the number of choices for $\{\rho_i\}_{i=1}^n$. For each $\rho_i\in\{1,\cdots,n-i+2\}$, we have at most $(n-i+2)!$ choices. As a result, the number of choices of $\{\rho_i\}_{i=1}^n$ is  at most 
	\begin{equation}\label{lemma:Q1FinalEstimate:E0}
		\prod_{i=0}^{n-1}(n-i+2)!\le c_{\rho}n^n,
	\end{equation}
	in which $c_\rho$ is a universal constant.  
	Using the result of the previous propositions, we obtain the conclusion.
\end{proof}

\subsection{Estimates of ${Q}^4$}
\begin{proposition}[The first estimate on   ${Q}^4$]\label{lemma:Q4FirstEstimate}  Let $	\tilde{Q}^4$ be a term in the sum of $Q^{4,nonpair}$. There is a constant $\mathfrak{C}_{Q_{4,1}}>0$ such that for $1\le n\le \mathfrak{N}$ and $t=\tau\lambda^{-2}>0$, we have
	\begin{equation}
		\begin{aligned} \label{lemma:Q4FirstEstimate:1}
			\tilde{Q}^4(\tau) :=\	& 
			\lambda^{n}\mathbf{1}({\sigma_{n,1}=-1})\mathbf{1}({\sigma_{n,2}=1})\Big[  \\&  \sum_{\substack{\bar\sigma\in \{\pm1\}^{\mathcal{I}_n},\\ \sigma_{i,\rho_i}+\sigma_{i-1,\rho_i}+\sigma_{i-1,\rho_i+1}\ne \pm3,
					\\ \sigma_{i-1,\rho_i}\sigma_{i-1,\rho_i+1}= 1}}\Big[\int_0^t\mathrm{d}s_0\int_{(\Lambda^*)^{\mathcal{I}_n}}  \mathrm{d}\bar{k} \Delta_{n,\rho}(\bar{k},\bar\sigma) e^{{\bf i}s_0\vartheta_0}\left\langle\prod_{i=1}^{n+2}\alpha(k_{0,i},\sigma_{0,i})\right\rangle_{s_0}  e^{{\bf i}s_0\vartheta_0}\\
			&\times (1-\mho)  \prod_{i=1}^n\Big[\sigma_{i,\rho_i}\mathcal{M}( k_{i,\rho_i}, k_{i-1,\rho_i}, k_{i-1,\rho_i+1}) \\
			& \times\Phi_{1,i}(\sigma_{i-1,\rho_i}, k_{i-1,\rho_i},\sigma_{i-1,\rho_i+1}, k_{i-1,\rho_i+1})\Big] \int_{(\mathbb{R}_+)^{\{1,\cdots,n\}}}\mathrm{d}\bar{s} \delta\left(t-\sum_{i=0}^ns_i\right)\\
			&\times \prod_{i=1}^{n}e^{-s_i[\varsigma_{n-i}+\tau_i]} \prod_{i=1}^{n} e^{-{\bf i}t_i(s)\mathfrak{X}_i}\prod_{i=0}^{n}e^{-s_i \tau_{i}}h^d\Big]\Big],	\end{aligned}\end{equation} and
	\begin{equation}
		\begin{aligned} \label{lemma:Q4FirstEstimate:2}
			&	\lim_{\lambda\to0}\limsup_{D\to\infty} \Big|{\tilde{Q}^4(\tau)}\Big|
			\ 
			=\   0,
	\end{aligned}\end{equation}
	where the same notations as in Proposition \ref{lemma:BasicGEstimate1} have been used.  
\end{proposition}
\begin{proof}
	The proof of the Proposition is the same as  that of Proposition \ref{lemma:BasicGEstimate1}. 
\end{proof}

\begin{proposition}[The main estimate on   ${Q}^4$]\label{lemma:Q4SecondEstimate}
	There is a constant $\mathfrak{C}_{Q_{4,2}}>0$ such that for $\tau=t\lambda^{-2}>0$, we have 
	
	\begin{equation}
		\begin{aligned} \label{lemma:Q4SecondEstimate:1}
			&	\lim_{\lambda\to0}\limsup_{D\to\infty} 
			\Big|Q^{4,nonpair}\Big|
			\ 				=\     0,
	\end{aligned}\end{equation}
	where the same notations as in Proposition \ref{lemma:BasicGEstimate1} have been used.  
\end{proposition}
\begin{proof} The proposition is an application of Proposition \ref{lemma:Q4FirstEstimate}. 
\end{proof}
\subsection{Estimates of ${Q}^2$}

\begin{proposition}[The first estimate on   ${Q}^2$]\label{lemma:BasicQ2Estimate1} Let $	\tilde{Q}^2$ be a term in the sum of $Q^{4,nonpair}$. There is a constant $\mathfrak{C}_{Q_{2,1}}>0$ such that for $[\mathfrak{N}/4]\le n\le \mathfrak{N}$ and $t=\tau\lambda^{-2}>0$, we have
	\begin{equation}
		\begin{aligned} \label{eq:BasicQ2Estimate1:1}
			\tilde{Q}^2(\tau):=\	&\varsigma_{n}	
			\lambda^{n}\mathbf{1}({\sigma_{n,1}=-1})\mathbf{1}({\sigma_{n,2}=1})\Big[  \\
			& \sum_{\substack{\bar\sigma\in \{\pm1\}^{\mathcal{I}_n},\\ \sigma_{i,\rho_i}+\sigma_{i-1,\rho_i}+\sigma_{i-1,\rho_i+1}\ne \pm3,
					\\ \sigma_{i-1,\rho_i}\sigma_{i-1,\rho_i+1}= 1}}\Big[\int_0^{s_0}e^{-(s_0-s_0')\varsigma_n}\mathrm{d}s_0'\int_{(\Lambda^*)^{\mathcal{I}_n}}  \mathrm{d}\bar{k} \Delta_{n,\rho}(\bar{k},\bar\sigma) e^{{\bf i}s_0'\vartheta_0}\left\langle\prod_{i=1}^{n+2}\alpha(k_{0,i},\sigma_{0,i})\right\rangle_{s_0'}\\
			&\times (1-\mho)\prod_{i=1}^n\Big[\sigma_{i,\rho_i}\mathcal{M}( k_{i,\rho_i}, k_{i-1,\rho_i}, k_{i-1,\rho_i+1})\Phi_{1,i}(\sigma_{i-1,\rho_i}, k_{i-1,\rho_i},\sigma_{i-1,\rho_i+1}, k_{i-1,\rho_i+1})\Big] \\
			&\times\int_{(\mathbb{R}_+)^{\{0,1,\cdots,n\}}}\mathrm{d}\bar{s} \delta\left(t-\sum_{i=0}^ns_i\right)\prod_{i=0}^{n}e^{-s_i[\varsigma_{n-i}+\tau_i]} h^d \prod_{i=0}^{n} e^{-{\bf i}t_i(s)\mathfrak{X}_i}\Big]\Big]\end{aligned}\end{equation}
and
	
	\begin{equation}
		\begin{aligned} \label{eq:BasicQ2Estimate1}
			&\lim_{\lambda\to 0}	\limsup_{D\to\infty} \Big|{\tilde{Q}^2(\tau)}\Big|		\
			=\ 0,
	\end{aligned}\end{equation}
	where   the same notations as in  Proposition \ref{lemma:BasicGEstimate1} have been used. 
\end{proposition}
\begin{proof}
	The quantity  $\tilde{Q}^2$comes from the partial time integration. The additional integration $\varsigma_n\int_0^{s_0}\mathrm{d}s_0'$ $e^{-(s_0-s_0')\varsigma_n}$ is needed. This has a bound
	$$\varsigma_n\int_0^{s_0}\mathrm{d}s_0'e^{-(s_0-s_0')\varsigma_n}=\int_0^{\varsigma_ns_0}\mathrm{d}s_0''e^{-(s_0\varsigma_n-s_0'')}=e^{-(s_0\varsigma_n-s_0'')}\Big|_0^{s_0\varsigma_n}\le 1.$$ 
 As a result, the same strategy of  Proposition \ref{lemma:BasicGEstimate1} can be reused, except that the integration of $s_0'$ can be bounded simply by $1$. 
\end{proof}
\begin{proposition}[The main estimate on   ${Q}^2$]\label{lemma:Q2FirstEstimate}
	There is a constant $\mathfrak{C}_{Q_{2,2}}>0$ such that for $t=\tau\lambda^{-2}>0$, we have 
	\begin{equation}
		\begin{aligned} \label{lemma:Q2FirstEstimate:2}
			&\lim_{\lambda\to 0}\limsup_{D\to\infty} 
			\Big|{{Q}^{2,nonpair}(\tau)}\Big|
			\ 
			=\   0.
	\end{aligned}\end{equation}

\end{proposition}

\begin{proof} 
	The proof follows the same argument used in the proof  of Proposition \ref{lemma:Q1FinalEstimate}.

\end{proof}

\subsection{First reduction of graphs}

\begin{proposition}\label{Lemma:Qmain} We define
	\begin{equation}
		\label{Def:Q0}
		Q^0=Q^1+Q^2+Q^3+Q^4.
	\end{equation}
	Suppose that $t>0$ and $t=\mathcal{O}(\lambda^{-2})$ 
	\begin{equation}\label{Lemma:Qmain:1}
		\begin{aligned}
			\lim_{\lambda\to0}\limsup_{D\to\infty} 
			&\Big | {Q}^0-  {Q}^{2,pair}-{Q}^{3,pair}-{Q}^{4,pair} \Big|
			\
			= \   0,\end{aligned}
	\end{equation}
	where, we have used the same notations as in Proposition   \ref{lemma:BasicGEstimate1}.
\end{proposition}
\begin{proof}
	The limit \eqref{Lemma:Qmain:1} is a direct consequence of Propositions \ref{lemma:Q1FinalEstimate},   \ref{lemma:Q4SecondEstimate}, \ref{lemma:Q2FirstEstimate}. 
\end{proof}

\section{Non-singular, pairing graph estimates}

First, we introduce a sufficiently large constant $\mathscr{M}>0$ and its conjugate $\mathscr{M}'$ such that $\frac{1}{\mathscr{M}}+\frac{1}{\mathscr{M}'}=1$. Due to the cut-off functions $\Phi_{1,i}$, we will only have to  consider graphs that have an even number of moments associated to the initial time slice $\prod_{j\in \{0,\cdots,n+2\}}  a_{s_0}(k_j,\sigma_j)$. For those graphs,  we only need to consider the case that for all time slices from $0$ to $n-1$,   all of the sets $\{k_{i,1},\cdots,k_{i,n-i+2}\}$, for $i=0,\cdots,n-1$, are admissible, otherwise the corresponding factor $\tau_i$ is strictly positive and the same argument used for \eqref{eq:Aestimate0:1:A:4} can be repeated and the contribution of the graph vanishes in the early limit of $D\to\infty$.  As $n$ is even, by Lemma \ref{Lemma:Triplet}, those graphs are pairing. In the next subsections, we will show that for long irreducible and delayed recollisional graphs, the estimate obtained in Proposition \ref{lemma:BasicGEstimate1} can be improved by  additional positive powers of $\lambda$. Those graphs are then negligible  when $\lambda$ is sufficiently small and the only graphs that mainly contribute to the total expansions are $\mathrm{iC}_2$ ladders. Therefore, we only need to consider components of the sums ${Q}_1$,  ${Q}_2$ ,  ${Q}_3$ and ${Q}_4$ in Proposition  \ref{Lemma:Qmain} that do not produce long irreducible and delayed recollisional graphs, leading to a reduction of ${Q}_1$,  ${Q}_2$ ,  ${Q}_3$ and ${Q}_4$ to sums of fewer terms: $\mathfrak{Q}_1$, $\mathfrak{Q}_2$ and $\mathfrak{Q}_3$ introduced in Proposition \ref{Propo:ReduceLadderGraph} below.

\subsection{Non-singular, pairing graphs with long irreducible  collisions}
The Proposition below provides some estimates on some components of the quantities ${Q}_1$,  ${Q}_2$ ,  ${Q}_3$ and ${Q}_4$ of Proposition \ref{Lemma:Qmain}, that have long irreducible  graphs. 
\begin{proposition}[Long irreducible  graphs]\label{Propo:CrossingGraphs} Suppose that the corresponding graph is long irreducible. There are constants $\mathfrak{C}_{Q_{LIC}},\mathfrak{C}'_{Q_{LIC}},\mathfrak{C}''_{Q_{LIC}} >0$ such that for $1\le n\le \mathfrak{N}$ and $t=\tau\lambda^{-2}>0$, we set
	\begin{equation}
		\begin{aligned} \label{Propo:CrossingGraphs:1:0}
			{Q}_{LIC}(\tau):=	& \left[
			\lambda^{n}\left[\sum_{\bar\sigma\in \{\pm1\}^{\mathcal{I}_n}}\int_{(\Lambda^*)^{\mathcal{I}_n}}  \mathrm{d}\bar{k}\Delta_{n,\rho}(\bar{k},\bar\sigma) \right. \right.h^d\\
			&\times \mathbf{1}({\sigma_{n,1}=-1})\mathbf{1}({\sigma_{n,2}'=1})\left\langle\prod_{i=1}^{n+2}\alpha(k_{0,i},\sigma_{0,i})\right\rangle_{0} \\
			&\times  \int_{(\Lambda^*)^{\mathcal{I}_n'}}\mathrm{d}\bar{k}\Delta_{n,\rho}(\bar k,\bar\sigma)\tilde{\Phi}_0(\sigma_{0,\rho_1},  k_{0,\rho_1},\sigma_{0,\rho_2},  k_{0,\rho_2})\\
			&\times \sigma_{1,\rho_1}\mathcal{M}( k_{1,\rho_1}, k_{0,\rho_1}, k_{0,\rho_1+1})\prod_{i=2}^n\Big[ \sigma_{i,\rho_i}\mathcal{M}( k_{i,\rho_i}, k_{i-1,\rho_i}, k_{i-1,\rho_i+1})\\
			&\times \Phi_{1,i}(\sigma_{i-1,\rho_i},  k_{i-1,\rho_i},\sigma_{i-1,\rho_i+1},  k_{i-1,\rho_i+1})\Big]\mho \\
			&\left.\left.\times\int_{(\mathbb{R}_+)^{\{0,\cdots,n\}}}\mathrm{d}\bar{s} \delta\left(t-\sum_{i=0}^ns_i\right)\prod_{i=1}^{n}e^{-s_i[\varsigma_{n-i}+\tau_i]}\prod_{i=1}^{n} e^{-{\bf i}t_i(s)\mathfrak{X}_i}\prod_{j\in I_{d1}}\diamond_\ell(s_{j-1})\right]\right],\end{aligned}\end{equation}
	then for any constants $1>T_*>0$ and $T_*>\tau>0$,
	
	\begin{equation}
		\begin{aligned} \label{Propo:CrossingGraphs:1}
			\limsup_{D\to\infty} 
			\Big\|{Q}_{LIC}(\tau)\Big\|_{L^4}
			\
			\le\ &   e^{T_*}\frac{T_*^{1+n}}{({n}_0(n))!}\lambda^{\mathfrak{C}'_{Q_{LIC}}}\mathfrak{C}_{Q_{LIC}}^{n}\langle \ln n\rangle \rangle^{ \mathfrak{C}''_{Q_{LIC}}}.
	\end{aligned}\end{equation}
	where the constants are universal and we have used the same notations as in Proposition \ref{lemma:BasicGEstimate1} and Proposition \ref{Lemma:Qmain}.  The cut-off function $\tilde{\Phi}_0$ can be either $\Phi_{0,1}$ or $\Phi_{1,1}$. If we replace $\left\langle\prod_{i=1}^{n+2}\alpha(k_{0,i},\sigma_{0,i})\right\rangle_{0}$ by $\left\langle\prod_{i=1}^{n+2}\alpha(k_{0,i},\sigma_{0,i})\right\rangle_{s_0} $ or    $\varsigma_n\int_0^{s_0}e^{-(s_0-s_0')\varsigma_n}\mathrm{d}s_0'$ $\left\langle\prod_{i=1}^{n+2}\alpha(k_{0,i},\sigma_{0,i})\right\rangle_{s_0'} ,$ as introduced in the quantities ${Q}_1,{Q}_2,{Q}_3,{Q}_4$  of Proposition \ref{Lemma:Qmain}, the same estimates hold true. 
\end{proposition}

\begin{proof}   By Lemma \ref{Lemma:Triplet}, we deduce that one of the quantities $\tau_i$ is non zero. Therefore, there exists  and a sequence of non-negative function $\diamond_\ell:\mathbb{R}_+\to\mathbb{R}_+$ such that $\diamond_\ell\in L^{\mathscr{M}'}(\mathbb{R}_+)\cap L^{2}(\mathbb{R}_+) $ and $\lim_{\ell\to 0}\|\diamond_\ell-\diamond_0\|_{L^{\mathscr{M}'}(\mathbb{R}_+)}=0$ with $\diamond_0(s)=1$ for all $s\in\mathbb{R}_+$ such that  we can bound $|{Q}_{LIC}|\le |{Q}_{LIC}^\ell|$, with
		\begin{equation}
		\begin{aligned} \label{Propo:CrossingGraphs:1:0:a0}
			{Q}_{LIC}^\ell(\tau):=	& \left[
			\lambda^{n}\left[\sum_{\bar\sigma\in \{\pm1\}^{\mathcal{I}_n}}\int_{(\Lambda^*)^{\mathcal{I}_n}}  \mathrm{d}\bar{k}\Delta_{n,\rho}(\bar{k},\bar\sigma) \right. \right.h^d\\
			&\times \mathbf{1}({\sigma_{n,1}=-1})\mathbf{1}({\sigma_{n,2}'=1})\left\langle\prod_{i=1}^{n+2}\alpha(k_{0,i},\sigma_{0,i})\right\rangle_{0} \\
			&\times  \int_{(\Lambda^*)^{\mathcal{I}_n'}}\mathrm{d}\bar{k}\Delta_{n,\rho}(\bar k,\bar\sigma)\tilde{\Phi}_0(\sigma_{0,\rho_1},  k_{0,\rho_1},\sigma_{0,\rho_2},  k_{0,\rho_2})\\
			&\times \sigma_{1,\rho_1}\mathcal{M}( k_{1,\rho_1}, k_{0,\rho_1}, k_{0,\rho_1+1})\prod_{i=2}^n\Big[ \sigma_{i,\rho_i}\mathcal{M}( k_{i,\rho_i}, k_{i-1,\rho_i}, k_{i-1,\rho_i+1})\\
			&\times \Phi_{1,i}(\sigma_{i-1,\rho_i},  k_{i-1,\rho_i},\sigma_{i-1,\rho_i+1},  k_{i-1,\rho_i+1})\Big]\mho \\
			&\left.\left.\times\int_{(\mathbb{R}_+)^{\{0,\cdots,n\}}}\mathrm{d}\bar{s} \delta\left(t-\sum_{i=0}^ns_i\right)\prod_{i=1}^{n}e^{-s_i[\varsigma_{n-i}+\tau_i]}\prod_{i=1}^{n} e^{-{\bf i}t_i(s)\mathfrak{X}_i}\prod_{j\in I_{d1}}\diamond_\ell(s_{j-1})\right]\right],\end{aligned}\end{equation}
			where $I_{d1}$ is the set of all indices $j$ such that the  vertex $v_j$ is a degree one vertex. 
	Therefore, instead of proving the required estimates for ${Q}_{LIC}(\tau)$, we will prove the required estimates for ${Q}_{LIC}^\ell(\tau)$.
	 We only prove the estimate for $\left\langle\prod_{i=1}^{n+2}\alpha(k_{0,i},\sigma_{0,i})\right\rangle_{0}$, the other ones can be proved by precisely the same argument.
	We observe that $\tilde{\Phi}_0$ is associated to the first interacting vertex $v_1$, which cannot be a degree-one vertex, otherwise, the edge in $\mathfrak{E}_+(v_1)$  is singular.   As $v_1$ is a degree-zero vertex, we can simply  bound $\tilde{\Phi}_0$  by a constant in our estimates. As the graph is pairing, the power of $h$ in the estimates of Proposition \ref{Propo:ExpectationEqui} are naturally incorporated into the delta functions of the pairing graphs. 
	
	According to our definition and Lemma \ref{Lemma:VerticesLongCollisions}, since the graph contains at least one long irreducible cycle, we denote  the set of all degree-one vertices at the top of those long irreducible cycles by $\{v_{l}\}_{l\in I}$. Suppose that $l_1=\min\{l\in I\}$, then  in the  long irreducible $\mathrm{iC}_m^i$ cycle of $v_{l_1}$, there  exists at least one vertex  $v_{i}$ such that $\mathfrak{X}_{i}+\mathfrak{X}_{l_1}$ is a function of the free edge  of $v_{l_1}$. We denote by $l_0$ the largest among such indices. We now divide the proof into several steps. Below we denote ${\Phi}_{1,l}$ as the cut-off function associated to the vertex $v_l$. 
	
	\smallskip{}{}
	{\bf Step 1: Analyzing the cycle of $v_{l_1}$ and the vertices between $v_{l_0}$ and $v_{l_1}$.}

	Suppose that the momenta of the edges in $\mathfrak{E}(v_{l_1})$ are $k_0,k_1,k_2$, in which $k_1,k_2\in\mathfrak{E}_-(v_{l_1})$  and $k_1$ is the free momenta. They are equipped with the signs $\sigma_{k_0},\sigma_{k_1},\sigma_{k_2}$. 
	In view of $l_0<l_1,$  we have
	\begin{equation}
		\label{Propo:CrossingGraphs:E0}
		\mathrm{Re}\vartheta(l_0-1) \ = \ \mathfrak{X}_{l_1} \ + \ \mathfrak{X}_{l_0} \ + \ \sum_{i>l_1}\mathfrak{X}_i,
	\end{equation}
	for $l_1=l_0+1$.

	In the case that $l_1>l_0+1$, we obtain	
	\begin{equation}
		\label{Propo:CrossingGraphs:E1}
		\mathrm{Re}\vartheta(l_0-1) \ = \ \mathfrak{X}_{l_1} \ + \ \mathfrak{X}_{l_0} \ + \  \sum_{i=l_0+1}^{l_1-1}\mathfrak{X}_i \ + \ \sum_{i>l_1}\mathfrak{X}_i.
	\end{equation}

	From our construction, the quantity $\sum_{i>l_1}\mathfrak{X}_i$ is not a function of  $k_1$ as well as any  free edge associated to $v_j$ with $j\le l_1$ by Lemma \ref{Lemma:VerticeOrderInAcycle}. 
	We denote the three edges associated to $v_{l_0}$ by $k_0'\in\mathfrak{E}_+(v_{l_0})$ and $k_1',k_2'\in\mathfrak{E}_-(v_{l_0})$. They are equipped with the signs $\sigma_{k_0'},\sigma_{k_1'},\sigma_{k_2'}$. 
	
	Let us consider a vertex $i$ such that  $l_0+1\le i<l_1$.
	If $\mathrm{deg}v_i=0$, then $\mathfrak{X}_i+\mathfrak{X}_{l_1}$ is independent of $k_1$, otherwise, when  $\mathfrak{X}_i+\mathfrak{X}_{l_1}$ depends on $k_1$, $v_i$ must belong to the cycle of $v_{l_1}$ and this contradicts the assumption that $l_0$ is the largest index in the cycle such that  $\mathfrak{X}_{l_0}+\mathfrak{X}_{l_1}$ depends on the free edge of $v_{l_1}$. If $\mathrm{deg}v_i=1$ and $l_0+1\le i<l_1$, by our assumption, $v_i$ has to  correspond to a  recollision $\mathrm{iC}_2^r$ or to a cycle formed by iteratively applying the recollisions (in Figure \ref{Fig37}, the four vertices $v_3,v_4,v_5,v_6$ are formed by iteratively applied two recollisions). We consider the three cases that happen when $\mathrm{deg}v_i=1$.
	
	\smallskip
	{Case 1: The recollision $\mathrm{iC}_2^r$ of $v_i$ corresponds to a double-cluster recollision and $i>l_0+1$.} 
	In this case, suppose that $k_0''$ is associated to the edge in $\mathfrak{E}_+(v_i)$ and $k_0'''$ is associated to the edge in $\mathfrak{E}_+(v_{i-1})$.  Since the two edges of $\mathfrak{E}_-(v_i)$ are paired with the two edges of $\mathfrak{E}_-(v_{i-1})$, we deduce that $\sigma_{k_0''}k_0''=-\sigma_{k_0'''} k_0'''$ and 
	$\mathfrak{X}_i+\mathfrak{X}_{i-1}  =  \sigma_{k_0}\omega(k_0'')  +  \sigma_{k_0'''}\omega(k_0'') =  0.$

	\smallskip
	{Case 2: The recollision $\mathrm{iC}_2^r$ of $v_i$ corresponds to a single-cluster recollision and $i>l_0+1$.} 
	Denote the three edges associated to $v_i$ by $k_0'',k_1'',k_3''$, in which $k_3''\in\mathfrak{E}_+(v_i)$. Denote the three edges associated to $v_{i-1}$ by $k_0'',k_2'',k_4''$. 
	Due to the pairing property of the single-cluster recollision, $\sigma_{k_2''}k_2''=-\sigma_{k_3''} k_3''.$ As a result,
	$\mathfrak{X}_i+\mathfrak{X}_{i-1}  =  \sigma_{k_3''}\omega(k_3'')  +  \sigma_{k_2''}\omega(k_2'') = 0.$

	\smallskip	
	{Case 3: $v_i$ corresponds to a  recollision and $i=l_0+1$}. Suppose that in the cycle of $v_{l_1}$, there is another vertex $v_{l_2}$, with $l_2>i$ and $\mathrm{deg}v_{l_2}=0$ while $\{v_{l_2},v_{l_2+1}\}$ do not form a recollision within the cycle of $v_{l_1}$. Then since $l_0$ is the largest among the indices of the cycle that make $\mathfrak{X}_i+\mathfrak{X}_{l_1}$ a function of $k_1$, we deduce that $\mathfrak{X}_{l_2}+\mathfrak{X}_{l_1}$ is independent of $k_1$. Therefore, by Lemma \ref{Lemma:VerticesLongCollisions},  $v_{l_2}$ is the only vertex with this property in the cycle and the cycle is then not long irreducible, contradicting our original assumption. As a result, there is no vertex $v_{l_2}$ in the cycle of $v_{l_1}$, with $l_2>i$ and  $\mathrm{deg}v_{l_2}=0$  while $\{v_{l_2},v_{l_2+1}\}$ do not form a recollision within the cycle of $v_{l_1}$. Since $\{v_{i},v_{i-1}\}=\{v_{l_0+1},v_{l_0}\}$ form a recollision, we have $\mathfrak{X}(v_{l_0+1})=-\mathfrak{X}(v_{l_0})$. Since $v_{l_0}$ has the largest index among all of the vertices $v_j$ inside the cycle of $v_{l_1}$ that has the quantity $\mathfrak{X}(v_{j})+\mathfrak{X}(v_{l_1})$ to be a function of $k_1$, we deduce that $\mathfrak{X}(v_{l_0+1})+\mathfrak{X}(v_{l_1})$ is then independent of $k_1$. Therefore, the graph is  a delayed recollision and not an irreducible long collision.   
	In conclusion, Case 3 does not happen in general. 
	Thus, $v_{l_0}$ is a long irreducible vertex. As between $v_{l_0}$, $v_{l_1}$, there are only recollisions, the difference $l_1-l_0$ is an odd number.
	

	In the above process,  if we define $$\mathfrak{I}=\Big\{i\in \{l_0+1,\cdots,l_1-1\} \Big| \mathrm{deg}v_i=1\Big\},$$
	and $$\mathfrak{I}'=\Big\{i\in \{l_0+1,\cdots,l_1-1\}, i+1\notin \mathfrak{I} \Big| \mathrm{deg}v_i=0\Big\},$$
	then 
	$$\vartheta_*=\sum_{i\in \mathfrak{I}}[\mathfrak{X}_i+\mathfrak{X}_{i-1}]+\sum_{i\in \mathfrak{I}'}\mathfrak{X}_i=\sum_{i\in \mathfrak{I}'}\mathfrak{X}_i,$$
	and
	$$\vartheta_{**}=\sum_{i>l_1}\mathfrak{X}_i,$$
	are both independent of $k_1$. Note that in the case of \eqref{Propo:CrossingGraphs:E0}, the quantity $\vartheta_*$ is simply $0$.
	As a result, we can write
	\begin{equation}
		\label{Propo:CrossingGraphs:E2}
		\mathrm{Re}\vartheta(l_0-1) \ = \  \mathfrak{X}_{l_1} + \mathfrak{X}_{l_0}+\vartheta_* + \vartheta_{**},\ \  \mathrm{Re}\vartheta(l_1-1) \ = \  \mathfrak{X}_{l_1} + \vartheta_{**}.
	\end{equation}

	Let us now consider the  $v_{l_0}$ vertex and its three momenta $k_0',k_1',k_2'$. We deduce that two of these momenta,  $k_{j_2}',k_{j_3}'$, ($j_2,j_3\in\{0,1,2\}$), depend on $k_1$. The other edge, denoted by $k_{j_1}'$, ($j_1\in \{0,1,2\}\backslash\{j_2,j_3\}$), is independent of $k_1$. 
	We  suppose that $k'_{j_2}=\sigma'  k_1+ \varkappa_1$, in which $\sigma'$ can be either $1$ or $-1$ and $\varkappa_1$ is independent of $k_1$, and $k'_{j_3}=\sigma''  k_1+ \varkappa_2$, in which $\sigma''$ can be either $1$ or $-1$ and $\varkappa_2$ is independent of $k_1$. It follows from the delta function at the vertex $v_{l_0}$ that $
	\sigma_{k_{j_3}'}k_{j_1}'+\sigma_{k_{j_2}'}k_{j_2}'+\sigma_{k_{j_3}'}k_{j_3}'=0,
	$
	which implies 
	$\sigma_{k_{j_1}'}k_{j_1}'+\sigma_{k_{j_2}'}(\sigma'  k_1+ \varkappa_1)+\sigma_{k_{j_3}'}(\sigma''  k_1+ \varkappa_2)=0.$
	Thus $\sigma_{k_{j_2}'}\sigma'+\sigma_{k_{j_3}'}\sigma''=0,$ and $\sigma_{k_{j_1}}k_{j_1}'+\sigma_{k_{j_2}'} \varkappa_1+\sigma_{k_{j_3}'}\varkappa_2=0.$
	Defining $\tilde\varkappa_1'=-\sigma_{k_{j_2}'} \varkappa_1$, $\tilde\varkappa_2'=\sigma_{k_{j_3}'} \varkappa_2$, we get $\tilde\varkappa_1'-\tilde\varkappa_2'=\sigma_{k_{j_1}'}k_{j_1}'.$
	Plugging the above expressions into $\mathfrak{X}_{l_0}$, we obtain \begin{equation}\label{eq:crossingraph9}
		\begin{aligned}
			\mathfrak{X}_{l_0}\ =\ & \sigma_{k_{j_1}'}\omega(k_{j_1}')+\sigma_{k_{j_2}'}\omega(\sigma'  k_1+ \varkappa_1)+\sigma_{k_{j_3}'}\omega(\sigma''  k_1+ \varkappa_2)\\
			\ =\ & \sigma_{k_{j_1}'}\omega(k_{j_1}')+\omega(\sigma_{k_{j_2}'}\sigma'  k_1- \tilde\varkappa_1')+\omega(\sigma'' \sigma_{k_{j_3}'} k_1+ \tilde\varkappa_2')
			\\
			\ =\ & \sigma_{k_{j_1}'}\omega(k_{j_1}')+\omega(-\sigma'' \sigma_{k_{j_3}'} k_1- \tilde\varkappa_1')+\omega(\sigma'' \sigma_{k_{j_3}'} k_1+ \tilde\varkappa_2')\\
			\ =\ & \alpha_1\omega(k_1+ \tilde\varkappa_1)+\alpha_2\omega(k_1+ \tilde\varkappa_2) + \omega',
		\end{aligned}
	\end{equation}
	where $\alpha_1,\alpha_2\in\{\pm1\}$,  $\tilde\varkappa_1-\tilde\varkappa_2= \pm (\tilde\varkappa_1'-\tilde\varkappa_2')$, $\omega'$ is independent of $k_1$.
	
	We also write the phase of $v_{l_1}$ as follows
	\begin{equation}\label{eq:crossingraph9:1}
		\begin{aligned}
			\mathfrak{X}_{l_1}\ =\ & \sigma_{k_{0}}\omega(k_{0})+\sigma_{k_{1}}\omega(k_1)+\sigma_{k_{2}}\omega(k_2)\
			\ =\  \alpha_3\omega(k_1+ \tilde\varkappa_3)+\alpha_4\omega(k_1+ \tilde\varkappa_4) + \omega'',
		\end{aligned}
	\end{equation}
	where $\alpha_3,\alpha_4\in\{\pm1\}$,  $\tilde\varkappa_3,\tilde\varkappa_4$, $\omega''$ are independent of $k_1$. One of the two momenta $\tilde\varkappa_3,\tilde\varkappa_4$ should be $0$.

	We recall that for any $i\in\{l_0+1,\cdots,l_1-1\}$, if $v_i$ belongs to the cycle of $v_{l_1}$, then either $v_i$ forms a recollision with $v_{i+1}$ or $v_{i-1}$ within the cycle of $v_{l_1}$. Thus $k_1+ \tilde\varkappa_1$ or $k_1+ \tilde\varkappa_2$ should  coincide with either the momentum $k_1+ \tilde\varkappa_3$ or the momentum $k_1+ \tilde\varkappa_4$ of the vertex $v_{l_1}$.  We refer to Figure \ref{Fig42} for an illustration of this situation. Without loss of generality, we suppose  $\tilde\varkappa_4=\tilde\varkappa_1$ and write
	\begin{equation}\label{eq:crossingraph9:2}
		\begin{aligned}
			\mathfrak{X}_{l_1}
			\ =\ & \alpha_3\omega(k_1+ \tilde\varkappa_3)+\alpha_4\omega(k_1+ \tilde\varkappa_1) + \omega''.
		\end{aligned}
	\end{equation}
	\begin{figure}
		\centering
		\includegraphics[width=.49\linewidth]{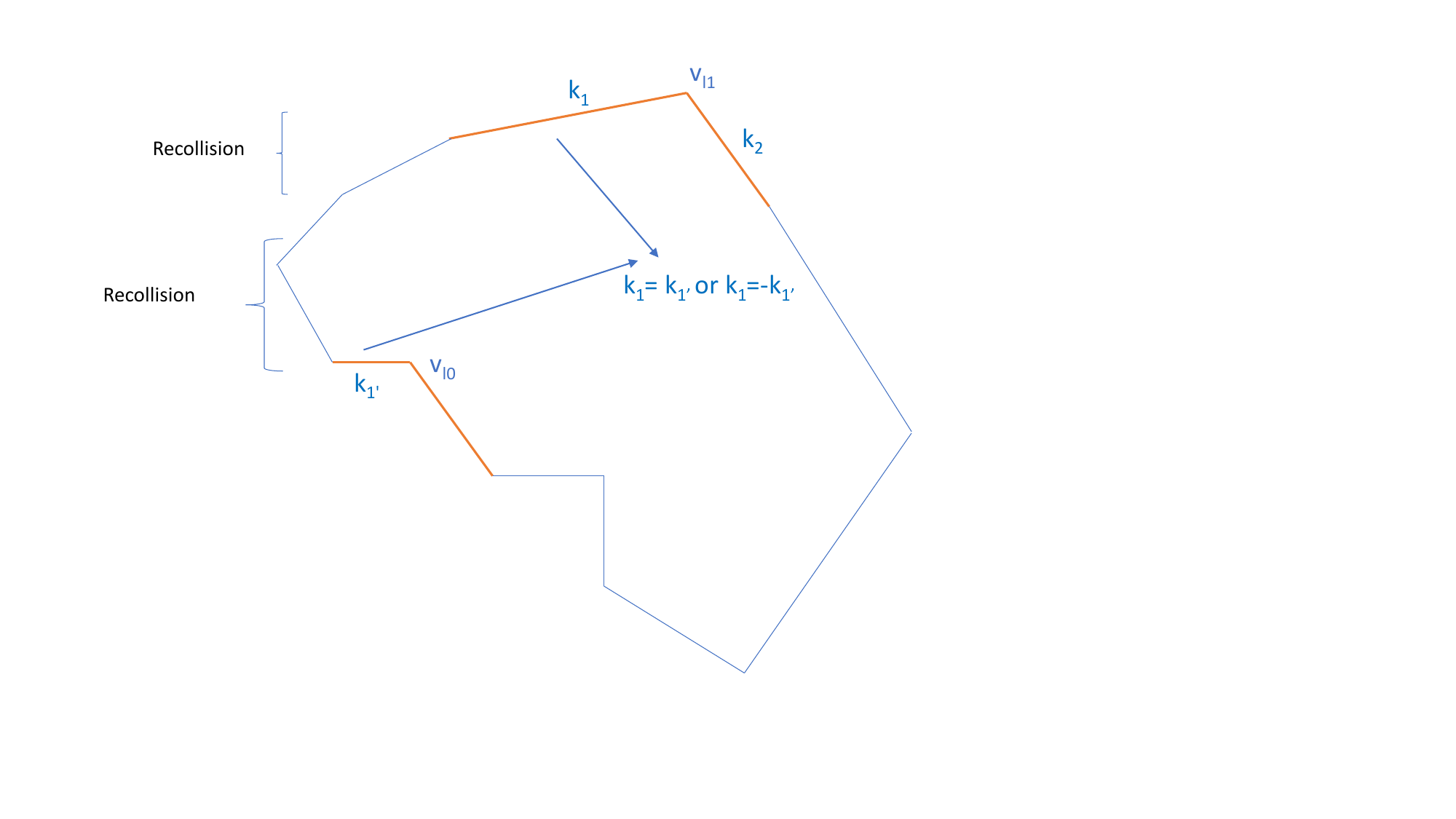}
		\caption{Illustration: $k_1+ \tilde\varkappa_1$ or $k_1+ \tilde\varkappa_2$ should  coincide with either the momentum $k_1+ \tilde\varkappa_3$ or the momentum $k_1+ \tilde\varkappa_4$ of the vertex $v_{l_1}$.}
		\label{Fig42}
	\end{figure}

	After rewriting the phases $\vartheta(l_0-1)$, $\vartheta(l_1-1)$ as in \eqref{Propo:CrossingGraphs:E2}, 	and $\mathfrak{X}_{l_1}$, $\mathfrak{X}_{l_0}$ as in \eqref{eq:crossingraph9:2}-\eqref{eq:crossingraph9}, we will move to Step 2 of the proof.
	
	\smallskip
	{\bf Step 2: Using the strategy of Proposition \ref{lemma:BasicGEstimate1} to estimate the graph.}	
	
 Now, we follow the strategy of Proposition \ref{lemma:BasicGEstimate1} by splitting the set $\{0,1,\cdots,n\}$ into smaller sets.  One of which is the set of time slice indices $1\le l\le n$ such that $v_{l+1}$ is of degree zero excluding $l_{0}-1$.   We denote this set by $I''$.   The other set then contains all of  time slice indices $0\le l<n$ such that $\mathrm{deg}(v_{l+1})=1$ and $n$. We denote this set by $I'$. We  obtain, following the same lines of computations as \eqref{eq:Aestimate1}
\begin{equation}
	\begin{aligned}\label{eq:crossingraph4}
			 & \	\Big|\lambda^{n}\mathbf{1}({\sigma_{n,1}=-1})\mathbf{1}({\sigma_{n,2}=1})\Big[  \\&\times \sum_{\substack{\bar\sigma\in \{\pm1\}^{\mathcal{I}_n},\\ \sigma_{i,\rho_i}+\sigma_{i-1,\rho_i}+\sigma_{i-1,\rho_i+1}\ne \pm3,
					\\ \sigma_{i-1,\rho_i}\sigma_{i-1,\rho_i+1}= 1}}\int_0^t\mathrm{d}s_0\int_{(\Lambda^*)^{\mathcal{I}_n}}  \mathrm{d}\bar{k}\Delta_{n,\rho}(\bar{k},\bar\sigma)\left\langle\prod_{i=1}^{n+2}\alpha(k_{0,i},\sigma_{0,i})\right\rangle_{0} h^d  \\
			&\times \Phi_{0,1}(\sigma_{0,\rho_1},  k_{0,\rho_1},\sigma_{0,\rho_2},  k_{0,\rho_2})\sigma_{i,\rho_i}\mathcal{M}( k_{1,\rho_1}, k_{0,\rho_1}, k_{0,\rho_1+1})\mho\\
			&\times \prod_{i=2}^n\Big[\sigma_{i,\rho_i}\mathcal{M}( k_{i,\rho_i}, k_{i-1,\rho_i}, k_{i-1,\rho_i+1})\Phi_{1,i}(\sigma_{i-1,\rho_i},  k_{i-1,\rho_i},\sigma_{i-1,\rho_i+1},  k_{i-1,\rho_i+1})\Big] \\
			&\times \int_{(\mathbb{R}_+)^{\{1,\cdots,n\}}}\mathrm{d}\bar{s} \delta\left(t-\sum_{i=0}^ns_i\right)\prod_{i=1}^{n}e^{-s_i \varsigma_{n-i} }\prod_{i=0}^{n}e^{-s_i \tau_{i}}\prod_{i=0}^{n-1} e^{-{\bf i}t_i(s)\mathfrak{X}_i}\Big]\Big|\\
		\lesssim	&\ \int_0^\tau\mathrm{d}\bar{s}_{n+1}      |F_1(\tau-\bar{s}_{n+1})|\lambda^{-2} |F_2(s_{n+1})|,
\end{aligned}\end{equation}
which means, in comparison to \eqref{eq:Aestimate1}, $F_1$ does not contain the time slice $l_{0}-1$.
To estimate $F_1$, we use an inequality similar to \eqref{eq:Aestimate2a} 
\begin{equation}
	\begin{aligned}
		& \int_{(\mathbb{R}_+)^{I''\cup\{n+1\}}}\mathrm{d}\bar{s}'\delta\left(t-\sum_{i\in I''\cup\{n+1\}}s_i\right)\prod_{i\in I''}e^{-\varsigma_{n-i}s_i}\\
		\quad \le\ &   \int_{(\mathbb{R}_+)^{I''\cup\{n+1\}}}\mathrm{d}\bar{s}'\delta\left(t-\sum_{i\in I''\cup\{n+1\}}s_i\right)\ \le\   \frac{t^{\frac{n}{2}-1}}{(\frac{n}{2}-1)!} \ = \ \lambda^{-2(\frac{n}{2}-1)} \frac{\tau^{\frac{n}{2}-1}}{(\frac{n}{2}-1)!}, 
\end{aligned}\end{equation} 
and the notation $	\bar{s}'$ stands for the vector whose components are $s_i$ with $i\in I''\cup\{n+1\}$.   The factor $\frac{n}{2}-1$ appears due to the fact that we have moved one degree-zero vertex $v_{l_0}$ from the set $I''$ of ``degree-zero time slices'' to $I'$, which is the set of ``degree-one time slices''.  The quantity $F_2$ contains only the integration on time slices that are associated to degree one vertices. For all of the degree one vertices, except $v_{l_1}$, we simply apply the estimate of  \ref{lemma:degree1vertex}. The terms that are associated to $v_{l_0},v_{l_1}$  are grouped and estimated separately (see Figure \ref{Fig14}).\begin{figure}
		\centering
		\includegraphics[width=.49\linewidth]{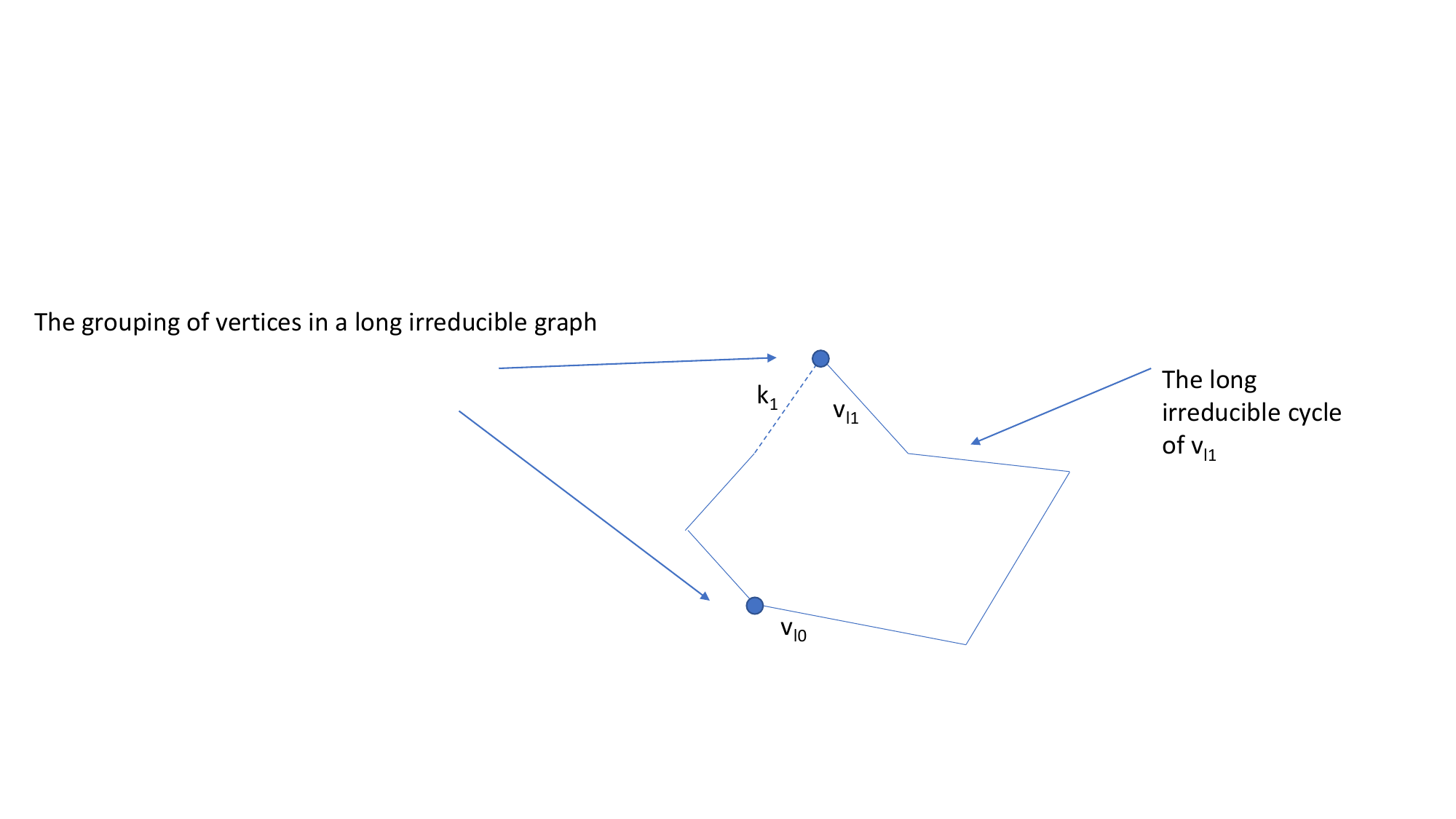}
		\caption{ The grouping of the 2 vertices $v_{l_0},v_{l_1}$ in a long irreducible collision gives an extra power of $\lambda$.}
		\label{Fig14}
	\end{figure}

	\smallskip
	{\bf Step 3: Grouping $v_{l_0}$ and $v_{l_1}$.}

	We now develop the quantities that are associated to $v_{l_0}$ and $v_{l_1}$, which are those that contain $\vartheta_i$ with $i=l_0-1, l_1-1$ respectively. We set $r_1=s_{l_0-1}$ and $r_2=s_{l_1-1}$ and  will need to estimate a quantity of the type
	\begin{equation}
		\begin{aligned}\label{eq:crossingraph6}
			&\int_{\Lambda^{*}}\mathrm{d}k_1
			 {{\Phi}_{1,l_0}(v_{l_0}) }  H(k_1) {\Phi}_{1,l_1}(v_{l_1})e^{-{\bf i}r_1(\mathfrak{X}_{l_1}+\vartheta_{**}) -{\bf i}r_2(\mathfrak{X}_{l_0}+\mathfrak{X}_{l_1}+\vartheta_{*}+\vartheta_{**})}e^{-\varsigma_{n-l_{0}+1}r_1}e^{-\varsigma_{n-l_{1}+1}r_2}\diamond_\ell(r_2),
		\end{aligned}
	\end{equation}
	where $H$ contain  all the quantities depending on $k_1$. Due to the structure of our graph, $H$ is obtained by iteratively applying the recollisions. This is precisely the situation described in Figures \ref{Fig41}, \ref{Fig47} and the collision operator
	\eqref{FinalProof:E10a:1:5} defined later. We therefore refer to the discussion of  Figures \ref{Fig41}, \ref{Fig47} and  
	\eqref{FinalProof:E10a:1:5} for the precise structure of $H$. In general, we suppose $H\in L^4(\mathbb{T}^d)$ (see also \eqref{OperatorRecollisionbis:1} for a similar situation). Note that $k_1$ is the free edge attached to $v_{l_1}$ and $H$ is bounded due to Lemma \ref{Lemma:Cutoffn} and the cut-off function defined in Definition \ref{def:hbar}
	. 
	Using the Lebesgue dominated convergence theorem, we can replace $\Lambda^{*}$ by $\mathbb{T}^d$ in the limit $D\to\infty$ and we still denote the cut-off functions on $\mathbb{T}^d$ by ${\Phi}_{1,l_0}(v_{l_0})$. 	Note that in the above expression, we keep the cut-off function with respect to $v_{l_0}$ while bounding the other one of $v_{l_1}$ simply by $1$. According to our assumption, for any vertex $v_i\ne v_{l_1},v_{l_0},$ and $v_i$ belonging to the cycle of $v_{l_1}$:  (a) if $\mathrm{deg}v_i=1$, then  $v_i,v_{i-1}$ correspond to a recollision, and (b) if $\mathrm{deg}v_i=0$,  then $v_{i+1},v_{i}$ do not correspond to a recollision and hence $\mathfrak{X}(v_{l_1})+\mathfrak{X}(v_{i})$ is a function of $k_1$. Suppose that for any vertex  $v_i\ne v_{l_1},v_{l_0},$ and $v_i$  belonging to the cycle of $v_{l_1}$, only possibility (a) happens,  then  the cycle of $v_{l_1}$ is a delayed recollision and $\mathfrak{X}(v_{l_1})+\mathfrak{X}(v_{l_0})=0$, leading to a contradiction. Therefore, there exists at least another vertex $v_{l'}$ satisfying (b), that means $l'<l_0$ and $\mathfrak{X}(v_{l_1})+\mathfrak{X}(v_{l'})$ is a function of $k_1$.  Hence, $v_{l_0}$ cannot be associated to the first time slice. This means $\Phi_{1,l_0}=\Phi_{1}$.  When $\Phi_{1,l_0}$ is  the function $1$, we then need  a strategy to replace it by a different cut-off function as we will describe below. In this case, the number of momenta in the time slice associated to  $\Phi_{1,l_0}$  is an odd number. By the same argument used in Lemma  \ref{Lemma:Triplet}, we could see that the  splitting of the  time slice $l_0-1$ has to happen at $v_{l_0}$ and the two momenta split at $v_{l_0}$ belong to a triplet, whose sum is zero. Moreover, only the first case discussed in the proof of Lemma  \ref{Lemma:Triplet} happens: one of the momenta in this triplet will be split again  giving two pairings, which  belong to either a recollision or a delayed recollision. We exclude the possibility that this is a delayed recollision as this contradicts the definition of  long irreducible collisional graphs. One end of the recollision is associated to ${\Phi}_{1,l_0}(v_{l_0})$, while the other end is associated to a different cut-off function ${\Phi}_{1,l_0-1}(v_{l_0-1})$ of the next vertex $v_{l_0-1}$ which must be different from the function $1$. As ${\Phi}_{1,l_0}(v_{l_0})$ and ${\Phi}_{1,l_0-1}(v_{l_0-1})$ are associated to the same   recollision, we then replace ${\Phi}_{1,l_0}(v_{l_0})$ by ${\Phi}_{1,l_0-1}(v_{l_0-1}).$ After this replacement, we still denote ${\Phi}_{1,l_0-1}(v_{l_0-1})$ by ${\Phi}_{1,l_0}(v_{l_0})$ as this makes no difference. It is then clear that in this case, $l_0-1$ has to be strictly smaller than $l'$ and $v_{l_0-1}$ is not in the first time slice.

	 We bound, using H\"older's inequality

	\begin{equation}
		\begin{aligned}\label{eq:crossingraph8:1}
			&\int_0^{{\lambda^{-2}}}\int_0^{{\lambda^{-2}}}\mathrm{d}r_{1}\mathrm{d}r_2\left|\int_{\mathbb{T}^d}\mathrm{d}k_1
			{{\Phi}_{1,l_0}(v_{l_0}) }  H(k_1) {\Phi}_{1,l_1}(v_{l_1})\right.\\
			&\left. \times e^{-{\bf i}r_1(\mathfrak{X}_{l_1}+\vartheta_{**}) -{\bf i}r_2(\mathfrak{X}_{l_0}+\mathfrak{X}_{l_1}+\vartheta_{*}+\vartheta_{**})}e^{-\varsigma_{n-l_{0}+1}r_1}e^{-\varsigma_{n-l_{1}+1}r_2}\diamond_\ell(r_2)\right|\\
			\lesssim \ &\iint_{\mathbb{R}^2}\mathrm{d}r_{1}\mathrm{d}r_2e^{-\lambda^2[|r_1|+|r_2|]}\left|\int_{\mathbb{T}^d}\mathrm{d}k_1
			{{\Phi}_{1,l_0}(v_{l_0}) }  H(k_1) {\Phi}_{1,l_1}(v_{l_1})\right.\\
			&\left. \times e^{-{\bf i}r_1(\mathfrak{X}_{l_1}+\vartheta_{**}) -{\bf i}r_2(\mathfrak{X}_{l_0}+\mathfrak{X}_{l_1}+\vartheta_{*}+\vartheta_{**})}e^{-\varsigma_{n-l_{0}+1}r_1}e^{-\varsigma_{n-l_{1}+1}r_2}\diamond_\ell(r_2)\right|\\
			\lesssim \ &\left[\iint_{\mathbb{R}^2}\mathrm{d}r_{1}\mathrm{d}r_2e^{-\lambda^2p|r_1|}|\diamond_\ell(r_2)|^p\right]^\frac1p\left[\iint_{\mathbb{R}^2}\mathrm{d}r_{1}\mathrm{d}r_2e^{-\lambda^2q[|r_1|+|r_2|]}\left|\int_{\mathbb{T}^d}\mathrm{d}k_1
			{{\Phi}_{1,l_0}(v_{l_0}) }  H(k_1) {\Phi}_{1,l_1}(v_{l_1})\right.\right.\\
			&\left.\left. \times e^{-{\bf i}r_1(\mathfrak{X}_{l_1}+\vartheta_{**}) -{\bf i}r_2(\mathfrak{X}_{l_0}+\mathfrak{X}_{l_1}+\vartheta_{*}+\vartheta_{**})}e^{-\varsigma_{n-l_{0}+1}r_1}e^{-\varsigma_{n-l_{1}+1}r_2}\right|^q\right]^\frac1q\\
			\lesssim \ &\lambda^{-\frac2p}\left[\iint_{\mathbb{R}^2}\mathrm{d}r_{1}\mathrm{d}r_2e^{-\lambda^2q[|r_1|+|r_2|]}\left|\int_{\mathbb{T}^d}\mathrm{d}k_1
			{{\Phi}_{1,l_0}(v_{l_0}) }  H(k_1) {\Phi}_{1,l_1}(v_{l_1})\right.\right.\\
			&\left.\left. \times e^{-{\bf i}r_1(\mathfrak{X}_{l_1}+\vartheta_{**}) -{\bf i}r_2(\mathfrak{X}_{l_0}+\mathfrak{X}_{l_1}+\vartheta_{*}+\vartheta_{**})}e^{-\varsigma_{n-l_{0}+1}r_1}e^{-\varsigma_{n-l_{1}+1}r_2}\right|^q\right]^\frac1q,\end{aligned}
		\end{equation}
	in which, we recall that $q=\mathscr{M}>1$ is a sufficiently large constant   and $\frac1p+\frac1q=1$.
	By a $TT^*$ argument, we bound
			\begin{equation}
				\begin{aligned}
					&\iint_{\mathbb{R}^2}\mathrm{d}r_{1}\mathrm{d}r_2\left|\int_{\mathbb{T}^d}\mathrm{d}k_1
					{{\Phi}_{1,l_0}(v_{l_0}) }  H(k_1) {\Phi}_{1,l_1}(v_{l_1})\right.\\
					&\left. \times e^{-{\bf i}r_1(\mathfrak{X}_{l_1}+\vartheta_{**}) -{\bf i}r_2(\mathfrak{X}_{l_0}+\mathfrak{X}_{l_1}+\vartheta_{*}+\vartheta_{**})}e^{-\varsigma_{n-l_{0}+1}r_1}e^{-\varsigma_{n-l_{1}+1}r_2}\diamond_\ell(r_2)\right|\\
			&	\lesssim\  \|H\|_{L^4}\lambda^{-\frac2p}\left[\iint_{\mathbb{R}^2}\mathrm{d}r_1\mathrm{d}r_2e^{-q\lambda^2|r_1|/2}e^{-q\lambda^2|r_2|/2}\left|\int_{\mathbb{T}^{d}}\mathrm{d}k_1 |{\Phi}_{1,l_0}(v_{l_0})|^2e^{-{\bf i}r_1(\mathfrak{X}_{l_1}+\vartheta_{**}) -{\bf i}r_2(\mathfrak{X}_{l_0}+\mathfrak{X}_{l_1}+\vartheta_{*}+\vartheta_{**})}\right|^\frac{q}{2}\right]^\frac1q,	\end{aligned}
	\end{equation}
	 where $r_1\vartheta_*+(r_1+r_2)\vartheta_{**}$ is independent of $k_1$. By the change of variable $r_1+r_2\to r_1$, we find
	\begin{equation}
		\begin{aligned}\label{eq:crossingraph8}
			&\iint_{\mathbb{R}^2}\mathrm{d}r_1\mathrm{d}r_2\left|\int_{\mathbb{T}^d}\mathrm{d}k_1
			{{\Phi}_{1,l_0}(v_{l_0}) }  H(k_1) {\Phi}_{1,l_1}(v_{l_1})e^{-{\bf i}r_1(\mathfrak{X}_{l_1}+\vartheta_{**}) -{\bf i}r_2(\mathfrak{X}_{l_0}+\mathfrak{X}_{l_1}+\vartheta_{*}+\vartheta_{**})}\right|\\
			&\times e^{-\varsigma_{n-l_{0}+1}r_1}e^{-\varsigma_{n-l_{1}+1}r_2}\diamond_\ell(r_2)\\
			 \lesssim & \|H\|_{L^4}\lambda^{-\frac2p} \left[\iint_{\mathbb{R}^2}\mathrm{d}r_1\mathrm{d}r_2e^{-q\lambda^2|r_1-r_2|/2}e^{-q\lambda^2|r_2|/2}\left|\int_{\mathbb{T}^{d}}\mathrm{d}k_1| {\Phi}_{1,l_0}(v_{l_0})|^2e^{-{\bf i}r_1\mathfrak{X}_{l_1} -{\bf i}r_2\mathfrak{X}_{l_0}-{\bf i}(r_2\vartheta_*+r_1\vartheta_{**})}\right|^\frac{q}{2}\right]^\frac1q,
		\end{aligned}
	\end{equation}
	
	and
	\begin{equation}
		\begin{aligned}\label{eq:crossingraph10}
			&\left|\int_{\mathbb{T}^{d}}\mathrm{d}k_1|{\Phi}_{1,l_0}(v_{l_0})|^2e^{-{\bf i}r_1\mathfrak{X}_{l_1}-{\bf i}r_2\mathfrak{X}_{l_0}-{\bf i}\vartheta^*}\right|\
			= \   \left|\int_{\mathbb{T}^{d}}\mathrm{d}k_1 |{\Phi}_{1,l_0}(v_{l_0})|^2\right.\\
			&\times e^{-{\bf i}r_1 (\sigma_{k_1}\omega({k_1})+\sigma_{k_2}\omega(k_2)+\sigma_{k_0}\omega(k_0))} e^{-{\bf i}r_2(\alpha_1\omega(k_1+ \tilde\varkappa_1)+\alpha_2\omega(k_1+ \tilde\varkappa_2) + \omega')}\left. e^{-{\bf i}\vartheta^*}\right|,
		\end{aligned}
	\end{equation}
	where $\vartheta^*=r_2\vartheta_*+r_1\vartheta_{**}$.

	By bounding all the terms that are independent of $k_1$ simply by $1$, we obtain 
	\begin{equation}
		\begin{aligned}\label{eq:crossingraph11}
			&\left|\int_{\mathbb{T}^{d}}\mathrm{d}k_1 |{\Phi}_{1,l_0}(v_{l_0})|^2e^{-{\bf i}r_1\mathfrak{X}_{l_1}-{\bf i}r_2\mathfrak{X}_{l_0}-{\bf i}\vartheta^*}\right|\\
			= \  & \left|\int_{\mathbb{T}^{d}}\mathrm{d}k_1 |{\Phi}_{1,l_0}(v_{l_0})|^2e^{-{\bf i}r_1 (\alpha_3\omega(k_1+ \tilde\varkappa_3)+\alpha_4\omega(k_1+ \tilde\varkappa_1) + \omega'')} e^{-{\bf i}r_2(\alpha_1\omega(k_1+ \tilde\varkappa_1)+\alpha_2\omega(k_1+ \tilde\varkappa_2) + \omega')}e^{-{\bf i}\vartheta^*}\right|\\
			\lesssim \  & \left|\int_{\mathbb{T}^{d}}\mathrm{d}k_1 |{\Phi}_{1,l_0}(v_{l_0})|^2e^{-{\bf i}r_1 (\alpha_3\omega(k_1+ \tilde\varkappa_3)+\alpha_4\omega(k_1+ \tilde\varkappa_1))}e^{-{\bf i}r_2(\alpha_1\omega(k_1+ \tilde\varkappa_1)+\alpha_2\omega(k_1+ \tilde\varkappa_2))}\right|\\
			\lesssim \  & \left|\int_{\mathbb{T}^{d}}\mathrm{d}k_1 |{\Phi}_{1,l_0}(v_{l_0})|^2e^{-{\bf i}r_1 \alpha_3\omega(k_1+ \tilde\varkappa_3)}e^{-{\bf i}(r_1 \alpha_4+r_2\alpha_1)\omega(k_1+ \tilde\varkappa_1)}\right.\left. e^{-{\bf i}r_2\alpha_2\omega(k_1+ \tilde\varkappa_2)}\right|.
		\end{aligned}
	\end{equation}

Using Lemma \ref{Lemma:Crossing}, we finally obtain
	\begin{equation}
		\begin{aligned}\label{eq:crossingraph17:d}
			&\iint_{\mathbb{R}^2}\mathrm{d}r_1\mathrm{d}r_2\left|\int_{\mathbb{T}^d}\mathrm{d}k_1
		{{\Phi}_{1,l_0}(v_{l_0}) }  H(k_1) {\Phi}_{1,l_1}(v_{l_1})e^{-{\bf i}r_1(\mathfrak{X}_{l_1}+\vartheta_{**}) -{\bf i}r_2(\mathfrak{X}_{l_0}+\mathfrak{X}_{l_1}+\vartheta_{*}+\vartheta_{**})}\right.\\
		&\times\left.e^{-\varsigma_{n-l_{0}+1}r_1}e^{-\varsigma_{n-l_{1}+1}r_2}\diamond_\ell(r_2)\right|\\
			&	\lesssim  \|H\|_{L^4}\langle\ln\lambda\rangle^{\mathfrak{C}_{\aleph}^{3}}\lambda^{-2/q+ \bar{\epsilon} }\lambda^{-2/p}\tilde{\mathfrak{H}}	\lesssim  \|H\|_{L^2}\langle\ln\lambda\rangle^{\mathfrak{C}_{\aleph}^{3}}\lambda^{-2+ \bar{\epsilon} }\tilde{\mathfrak{H}},
		\end{aligned}
	\end{equation}
	for some constant  $ \mathfrak{C}_{\aleph}^{3} >0$, $ \bar{\epsilon}>0$.  We recall that
	\begin{equation}
		\begin{aligned}\label{Lemma:Crossing:2}
			\tilde{\mathfrak{H}}: = \ & \prod_{j=2}^{d}  \Big[\bar{\mathbb{C}}_j^{-\frac{1}{2q_j}}+ \bar{\mathbb{A}}_j^{-\frac{1}{2q_j}}\Big]\sqrt{\tilde{F}(\tilde{\varkappa}_2-\tilde{\varkappa}_1)} \sqrt{\tilde{F}(\tilde{\varkappa}_3-\tilde{\varkappa}_1)},
		\end{aligned}
	\end{equation}
	where $\tilde{F}_1(\tilde{\varkappa}_2-\tilde{\varkappa}_1)$, $\tilde{F}_1(\tilde{\varkappa}_3-\tilde{\varkappa}_1)$ are components of $\sqrt{\check\Psi}$ that depend only on $\tilde{\varkappa}_2-\tilde{\varkappa}_1$ and $\tilde{\varkappa}_3-\tilde{\varkappa}_1$ respectively, $q_2,\cdots,q_{d}$ are positive real numbers in $(1,\infty)$ such that
	$\frac{1}{q_2}+\cdots+\frac{1}{q_{d}}=\frac{2}{q}$, and
	\begin{equation*}
		\begin{aligned}\label{Lemma:Crossing:3}
			\bar{\mathbb{A}}_j   
			\ =\ & |\alpha_3\alpha_1\cos(\tilde{\varkappa}_3^1-\tilde{\varkappa}_1^1){\sin(2(\tilde{\varkappa}_3^{j}-\tilde{\varkappa}_1^{j}))} - \alpha_2\alpha_4\cos(\tilde{\varkappa}_2^1-\tilde{\varkappa}_1^1)\sin(2(\tilde{\varkappa}_2^j-\tilde{\varkappa}_1^{j}))\\
			&+ \alpha_3\alpha_2\cos(\tilde{\varkappa}_3^1-\tilde{\varkappa}_1^1){\sin(2(\tilde{\varkappa}_3^{j}-\tilde{\varkappa}_1^{j}))} \cos(\tilde{\varkappa}_2^1-\tilde{\varkappa}_1^1)\cos(2(\tilde{\varkappa}_2^{j}-\tilde{\varkappa}_1^{j}))\\
			& -\alpha_3\alpha_2 \cos(\tilde{\varkappa}_3^1-\tilde{\varkappa}_1^1){\cos(2(\tilde{\varkappa}_3^{j}-\tilde{\varkappa}_1^{j}))}\cos(\tilde{\varkappa}_2^1-\tilde{\varkappa}_1^1)\sin(2(\tilde{\varkappa}_2^j-\tilde{\varkappa}_1^{j}))|, 
		\end{aligned}
	\end{equation*}
	and
	\begin{equation*}
		\begin{aligned}\label{Lemma:Crossing:4}
			\bar{\mathbb{C}}_j  
			\ =\ & |\sin(\tilde{\varkappa}_3^1-\tilde{\varkappa}_1^1){\cos(2(\tilde{\varkappa}_3^{j}-\tilde{\varkappa}_1^{j}))}\sin(\tilde{\varkappa}_2^1-\tilde{\varkappa}_1^1)\sin(2(\tilde{\varkappa}_2^{j}-\tilde{\varkappa}_1^{j}))\\
			& - \sin(\tilde{\varkappa}_3^1-\tilde{\varkappa}_1^1){\sin(2(\tilde{\varkappa}_3^{j}-\tilde{\varkappa}_1^{j}))} \sin(\tilde{\varkappa}_2^1-\tilde{\varkappa}_1^1)\cos(2(\tilde{\varkappa}_2^j-\tilde{\varkappa}_1^{j}))|\\
			\ =\ & |\sin(\tilde{\varkappa}_3^1-\tilde{\varkappa}_1^1)\sin(\tilde{\varkappa}_2^1-\tilde{\varkappa}_1^1)\sin(2(\tilde{\varkappa}_2^{j}-\tilde{\varkappa}_3^{j}))|.
		\end{aligned}
	\end{equation*}
	In this estimate, we obtain a factor of $\lambda^{-2+\bar{\epsilon}/q}$. We find the total factor $$\lambda^{-2(\frac{n}{2}-1)}\lambda^{-2+\bar{\epsilon}/q}=\lambda^{-n}\lambda^{2\bar{\epsilon}/q}.$$  The factor $\lambda^{2\bar{\epsilon}/q}$ caries a positive power of $\lambda$, we will see later that this extra factor guarantees the convergence to $0$ of terms associated to irreducible long graphs.

	\smallskip
	
	{\bf Step 4: Integrating $\tilde{\mathfrak{H}}$, using other free momenta.}	
	
	Next, we will need to integrate  $\tilde\varkappa_1-\tilde\varkappa_2$  and $\tilde\varkappa_1-\tilde\varkappa_3$. 
	It follows from Lemma \ref{Lemma:ZeroMomentum} that $\tilde\varkappa_1-\tilde\varkappa_2=\pm\sigma_{k_{j_1}'}k_{j_1}'$ cannot be inpendent of all free momenta, otherwise $k_{j_1}'$ becomes a singular momentum, contradicting our original assumption. We then suppose that $\tilde\varkappa_1-\tilde\varkappa_2$ depends on the free momentum $k_{e_1}$. Thus, $k_{e_1}$ is either the virtually free momentum of  the virtually free edge or associated to a free edge of a degree-one vertex $v_{l_{e_1}}$.  If $k_{e_1}$ is the virtually free momentum, it is associated to the virtual vertex $v_*$, which is above $v_{l_1}$. If $k_{e_1}$ is associated to a free edge of a degree-one vertex $v_{l_{e_1}}$, we will show that $l_{e_1}>l_1$. Suppose the contrary, then by the definition of $v_{l_1}$, we deduce that $v_{l_{e_1}}$ is the top vertex of a recollision, leading to a contradiction. Among the free momenta that $k_{j_1}'$ depends on, we choose $k_{l_{e_1}}$ to be the free momentum of a degree-one vertex  if possible. Only in the case that the virtually free momentum is the only free momentum that $k_{j_1}'$ depends on, then we choose $k_{l_{e_1}}$ to be the virtually free momentum.
	
	Now, since the momentum of the top edge $k_0$ of $v_{l_1}$ is independent of the free edge $k_1$, it should depend on some other free momentum $k_e$, otherwise, the graph is singular.  We denote the 
	vertex, that $k_e$ is attached to by $v_{l_e}$. The same argument as above shows that  $v_{l_e}$ is at a higher position in the diagram in comparison with  $v_{l_1}$. Among the free momenta that $k_0$ depends on, we choose $k_{l_{e}}$ to be the free momentum of a degree-one vertex $v_{l_{e}}$ if possible. Only in the case that the virtually free momentum is the only free momentum that $k_0$ depends on, then we choose $k_{l_{e}}$ to be the virtually free momentum.

	In the case that both $k_e$ and $k_{e_1}$ are not the virtually free momenta, if  $v_{l_e}$ and $v_{l_{e_1}}$ are distinct and as they are above the cycle of $v_{l_1}$, we can integrate them independently with respect to this cycle.  It is then safe to embed $\tilde{\mathfrak{H}}$ into the integrals of the free momenta $k_e$  and $k_{e_1}$.  
 	We also observe that the singularities created by 	
	$|\bar{\mathbb{A}}_j|^{\frac{1}{q_j}}$ and $|\bar{\mathbb{C}}_j|^{\frac{1}{q_j}}$
	are absorbed into those integrals, since they contribute  powers smaller than $1$ in the denominator of the integrals of $k_{e_1}$ and $k_e$, under the influence of the cut-off functions $\sqrt{\tilde{F}_1(\tilde{\varkappa}_2-\tilde{\varkappa}_1)}$, $ \sqrt{\tilde{F}_1(\tilde{\varkappa}_3-\tilde{\varkappa}_1)}$.
	
	In the case that $v_{l_e}$ and $v_{l_{e_1}}$ coincide and $k_e\equiv k_{e_1}$ is not the virtually free momentum, we only need to integrate once, with respect to the free momentum $k_e$, to remove the singularities. 
	
	In the case that either $k_e$ or $k_{e_1}$ is the virtually free edge,  then either $k_e$ or $k_{e_1}$ is $k_{n,1}$ (see Figure \ref{Fig40}). We then  need to perform the integration
	$
	\int_{\mathbb{T}^d}\mathrm{d}k_{n,1}$ in the $L^2$-norm to eliminate the singularities.   The integration of $\tilde\varkappa_1-\tilde\varkappa_3$ can be done in the standard way. 
	
	Therefore, we gain a factor of $\lambda^{\mathfrak{C}'_{Q_{Long}}}=\lambda^{2\bar{\epsilon}/q}$.  
	Thus the total $\lambda$ power   guarantees the convergence of the whole graph to $0$.

\end{proof}


\subsection{Non-singular, pairing graphs with delayed recollisions}

\begin{proposition}[Non-singular, pairing graphs with delayed recollisions] \label{Propo:LeadingGraphs} Suppose that the corresponding graph is delayed recollisional and $d\ge 2$. There are constants $\mathfrak{C}_{Q_{Delayed}}$, $\mathfrak{C}'_{Q_{Delayed}}$, $\mathfrak{C}''_{Q_{Delayed}}>0$ such that for $1\le n\le \mathfrak{N}$ and $t=\tau\lambda^{-2}>0$, we set
	\begin{equation}
		\begin{aligned} \label{Propo:LeadingGraphs:1:0}
			{Q}_{Delayed}^\ell(\tau):=\	& 
			\lambda^{n}\left[ \sum_{\bar\sigma\in \{\pm1\}^{\mathcal{I}_n}}\int_{(\Lambda^*)^{\mathcal{I}_n}}  \mathrm{d}\bar{k}\Delta_{n,\rho}(\bar{k},\bar\sigma)h^d \right. \\
			&\times \prod_{A=\{j,l\}\in S}\left[\delta(k_{0,l}+k_{0,j})\mathbf{1}({\sigma_{0,l}=-\sigma_{0,j}})\right]\mathbf{1}({\sigma_{n,1}=-1})\mathbf{1}({\sigma_{n,2}'=1})\left\langle\prod_{i=1}^{n+2}\alpha(k_{0,i},\sigma_{0,i})\right\rangle_{0}  \\
			&\times  \int_{(\Lambda^*)^{\mathcal{I}_n'}}\mathrm{d}\bar{k}\Delta_{n,\rho}(\bar k,\bar\sigma)\tilde{\Phi}_0(\sigma_{0,\rho_1},  k_{0,\rho_1},\sigma_{1,\rho_2},  k_{0,\rho_2})\\
			&\times \sigma_{1,\rho_1}\mathcal{M}( k_{1,\rho_1}, k_{0,\rho_1}, k_{0,\rho_1+1})\prod_{i=2}^n\Big[ \sigma_{i,\rho_i}\mathcal{M}( k_{i,\rho_i}, k_{i-1,\rho_i}, k_{i-1,\rho_i+1})\\
			&\times \Phi_{1,i}(\sigma_{i-1,\rho_i},  k_{i-1,\rho_i},\sigma_{i-1,\rho_i+1},  k_{i-1,\rho_i+1})\Big] \mho\\
			&\left.\times\int_{(\mathbb{R}_+)^{\{0,\cdots,n\}}}\mathrm{d}\bar{s} \delta\left(t-\sum_{i=0}^ns_i\right)\prod_{i=1}^{n}e^{-s_i[\tau_{n-i}+\tau_i]}\prod_{i=1}^{n} e^{-{\bf i}t_i(s)\mathfrak{X}_i}\right]\prod_{j\in I_{d1}}\diamond_\ell(s_{j-1}),
	\end{aligned}\end{equation}
then for any constants $1>T_*>0$ and $T_*>\tau>0$,

\begin{equation}
	\begin{aligned} \label{Propo:DelayedGraph:1}
		\limsup_{D\to\infty} 
		\Big\|{Q}_{Delayed}^\ell(\tau)\Big\|_{L^4}
		\
		\le\ &   e^{T_*}\frac{T_*^{{n}}}{({n}/2)!}\lambda^{\mathfrak{C}'_{Q_{Delayed}}}\mathfrak{C}_{Q_{Delayed}}^n}\langle \ln\lambda\rangle^{\mathfrak{C}''_{Q_{Delayed}},
\end{aligned}\end{equation}

and
\begin{equation}
	\begin{aligned} \label{Propo:DelayedGraph:2}
		&\lim_{\lambda\to0}\lim_{\ell\to0} \limsup_{D\to\infty} 
		\Big\|{Q}_{Delayed}^\ell(\tau)-{Q}_{Delayed}^0(\tau)\Big\|_{L^4}=0,
\end{aligned}\end{equation}
	where,   and we have used the same notations as in Propositions \ref{lemma:BasicGEstimate1},  \ref{Lemma:Qmain}, \ref{Propo:CrossingGraphs}.   The cut-off function $\tilde{\Phi}_0$ can be either $\Phi_{0,1}$ or $\Phi_{1,1}$. If we replace $\left\langle\prod_{i=1}^{n+2}\alpha(k_{0,i},\sigma_{0,i})\right\rangle_{0} $ by $\left\langle\prod_{i=1}^{n+2}\alpha(k_{0,i},\sigma_{0,i})\right\rangle_{s_0} $ or    $\varsigma_n\int_0^{s_0}e^{-(s_0-s_0')\varsigma_n}\mathrm{d}s_0'$ $\left\langle\prod_{i=1}^{n+2}\alpha(k_{0,i},\sigma_{0,i})\right\rangle_{s_0'} ,$ as introduced in the quantities ${Q}_1,{Q}_2,{Q}_3,{Q}_4$  of Proposition \ref{Lemma:Qmain}, the same estimates hold true. 
\end{proposition}

\begin{proof}
By Lemma \ref{Lemma:Triplet}, we deduce that one of the quantities $\tau_i$ is non zero. Therefore, there exists  and a sequence of non-negative function $\diamond_\ell:\mathbb{R}_+\to\mathbb{R}_+$ such that $\diamond_\ell\in L^{\mathscr{M}'}(\mathbb{R}_+)\cap L^{2}(\mathbb{R}_+) $ and $\lim_{\ell\to 0}\|\diamond_\ell-\diamond_0\|_{L^{\mathscr{M}'}(\mathbb{R}_+)}=0$ with $\diamond_0(s)=1$ for all $s\in\mathbb{R}_+$ such that  we can bound $|{Q}_{Delayed}|\le |{Q}_{Delayed}^\ell|$, with
	\begin{equation}
	\begin{aligned} \label{Propo:LeadingGraphs:1:0}
		{Q}_{Delayed}^\ell(\tau):=\	& 
		\lambda^{n}\left[ \sum_{\bar\sigma\in \{\pm1\}^{\mathcal{I}_n}}\int_{(\Lambda^*)^{\mathcal{I}_n}}  \mathrm{d}\bar{k}\Delta_{n,\rho}(\bar{k},\bar\sigma)h^d \right. \\
		&\times \prod_{A=\{j,l\}\in S}\left[\delta(k_{0,l}+k_{0,j})\mathbf{1}({\sigma_{0,l}=-\sigma_{0,j}})\right]\mathbf{1}({\sigma_{n,1}=-1})\mathbf{1}({\sigma_{n,2}'=1})\left\langle\prod_{i=1}^{n+2}\alpha(k_{0,i},\sigma_{0,i})\right\rangle_{0}  \\
		&\times  \int_{(\Lambda^*)^{\mathcal{I}_n'}}\mathrm{d}\bar{k}\Delta_{n,\rho}(\bar k,\bar\sigma)\tilde{\Phi}_0(\sigma_{0,\rho_1},  k_{0,\rho_1},\sigma_{1,\rho_2},  k_{0,\rho_2})\\
		&\times \sigma_{1,\rho_1}\mathcal{M}( k_{1,\rho_1}, k_{0,\rho_1}, k_{0,\rho_1+1})\prod_{i=2}^n\Big[ \sigma_{i,\rho_i}\mathcal{M}( k_{i,\rho_i}, k_{i-1,\rho_i}, k_{i-1,\rho_i+1})\\
		&\times \Phi_{1,i}(\sigma_{i-1,\rho_i},  k_{i-1,\rho_i},\sigma_{i-1,\rho_i+1},  k_{i-1,\rho_i+1})\Big] \mho\\
		&\left.\times\int_{(\mathbb{R}_+)^{\{0,\cdots,n\}}}\mathrm{d}\bar{s} \delta\left(t-\sum_{i=0}^ns_i\right)\prod_{i=1}^{n}e^{-s_i[\tau_{n-i}+\tau_i]}\prod_{i=1}^{n} e^{-{\bf i}t_i(s)\mathfrak{X}_i}\right]\prod_{j\in I_{d1}}\diamond_\ell(s_{j-1}).
\end{aligned}\end{equation}
	Therefore, instead of proving the required estimates for ${Q}_{Delayed}(\tau)$, we will prove the required estimates for ${Q}_{Delayed}^\ell(\tau)$.
	We now prove the first estimate  \eqref{Propo:DelayedGraph:1} for ${Q}_{Delayed}^\ell(\tau)$ and $ \left\langle\prod_{j\in \{0,\cdots,n+2\}}  a(k_j,\sigma_j)\right\rangle_{0}$, as the other ones can be treated by precisely the same argument.
	We can suppose that $\Lambda^*$, after taking the limit $D\to\infty$, can be replaced by $\mathbb{T}^d$.  
	\begin{figure}
		\centering
		\includegraphics[width=.49\linewidth]{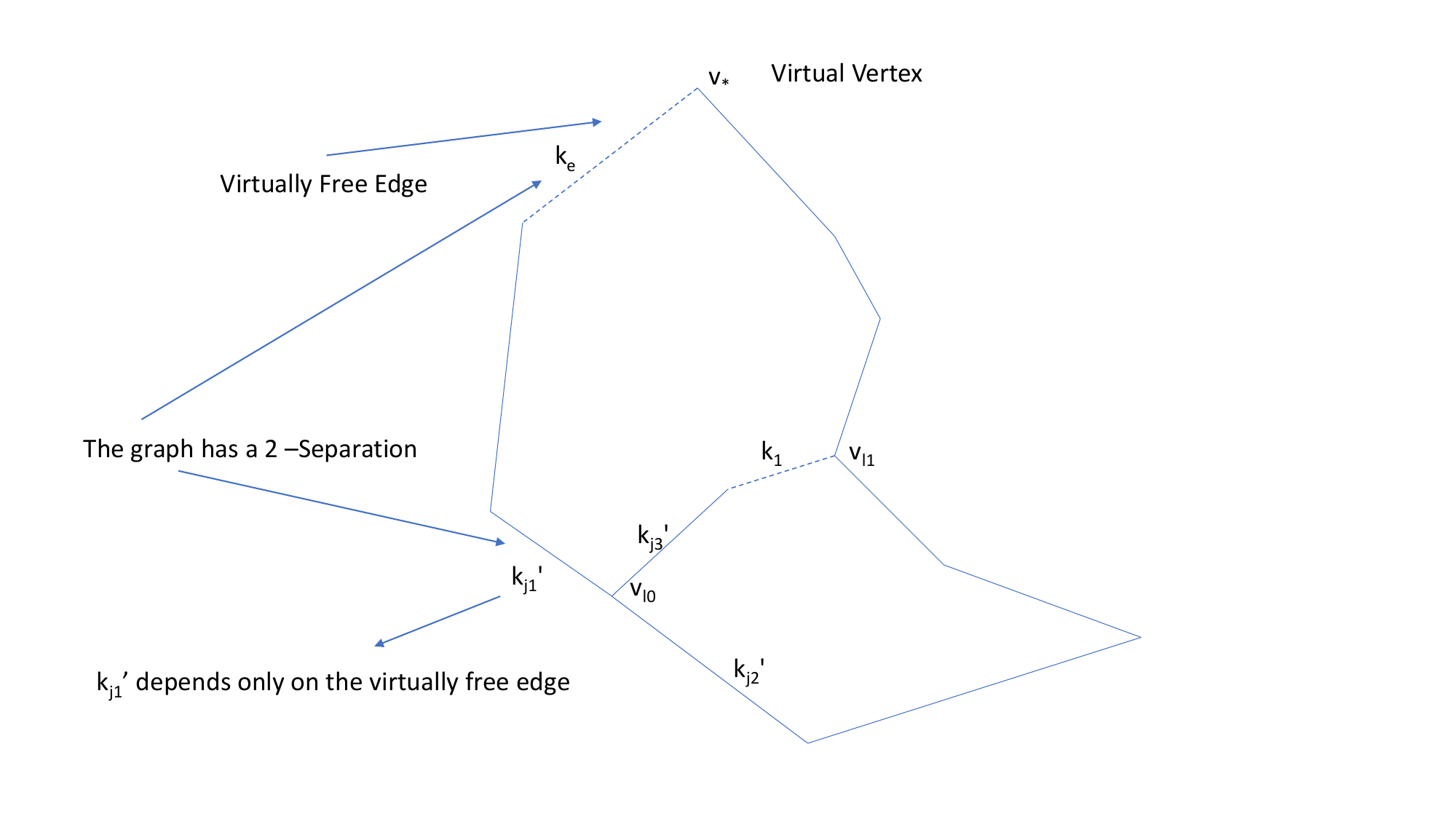}
		\caption{ An example in which $k_e$ is virtually free.}
		\label{Fig40}
	\end{figure}
	In this delayed recollisional graph, let us denote by $l_1$ the index of the first degree-one vertex $v_{l_1}$, such that it does not corresponds to a double-cluster, a single-cluster recollision or a cycle formed by iteratively applying   the recollisions. 
	Since this is the first among such vertices, according to the definition, the vertex corresponds to a delayed recollision. Denote by $k_0$ the edge in $\mathfrak{E}_+(v_{l_1})$ and $k_1,k_2$ the two edges in $\mathfrak{E}_-(v_{l_1})$, in which $k_1$ is the free edge. As usual, we denote the signs of those edges by $\sigma_{k_0},\sigma_{k_1},\sigma_{k_2}$. By the second conclusion  (ii) of Lemma \ref{Lemma:VerticesLongCollisions}  there exists a unique vertex $v_{l_2}$, $l_1-l_2>1,$ within the cycle of $v_{l_1}$ such that $\mathfrak{X}_{l_2}+\mathfrak{X}_{l_1}$ is not a function of $k_1$. In addition, by Lemma \ref{Lemma:VerticeOrderInAcycle}, all quantities $\mathfrak{X}_i$ for $i>l_1$ are not functions of $k_1$.

	Let us now  denote by $\mathfrak{I}_{v_{l_1}}$ the collection of all vertices belonging to the cycle of $v_{l_1}$. Let $v_j$, $j\ne l_1$ be any vertex in  $\mathfrak{I}_{v_{l_1}}$. If $\mathrm{deg}v_j=1$, then $v_j$ corresponds to a recollision and $\mathrm{deg}v_{j-1}=0$. If $\mathrm{deg}v_j=0$, then either $v_j\equiv v_{l_2}$ or $v_j$ has to belong to a recollision within the cycle of $v_{l_1}$, $\mathrm{deg}v_{j+1}=1$ and $v_{j+1}$ also belongs to $\mathfrak{I}_{v_{l_1}}$ . By Lemma \ref{Lemma:VerticesLongCollisions}, $v_{l_2}$ is the only degree-zero vertex that does not belong to any recollisions in the cycle of $v_{l_1}$; while the other parts of the cycle contain only recollisions  or no recollision at all - the cycle is then just a short delayed  recollision (see Figure \ref{Fig8}).

	Let us now study the real part of the  total phase of the first $l_1$ time slices
	\begin{equation}
		\label{DelayedRecollision1}
		\sum_{i=0}^{l_1-1}s_i\mathrm{Re}\vartheta_i \ = \ \tilde\xi\ + \ \sum_{i=0}^{l_1-1}s_i\sum_{j=i+1}^{l_1}\mathfrak{X}_j.
	\end{equation}
	in which the quantity $\tilde\xi$ is independent of all of the free edges of the first $l_1$ time slices.
	
	Let us denote $$\mathfrak{D}_i=\{0\le j<l_1 \ | \ \mathrm{deg}v_j=i\}$$
	for $i=0,1$. 
	We also denote 
	$$\mathfrak{D}_0' \ = \ \{i\in \mathfrak{D}_0\backslash\{l_2\} \ | \ i+1 \in	 \mathfrak{D}_0\}. $$
	We then develop
	\begin{equation}
		\begin{aligned}
			& \sum_{i=0}^{l_1-1}s_i\sum_{j=i+1}^{l_1}\mathfrak{X}_j \ =  \ \sum_{j\in \mathfrak{D}_1} \left[\mathfrak{X}_j\sum_{i=0}^{j-1}s_i \ +\mathfrak{X}_{j-1}\sum_{i=0}^{j-2}s_i \right]\ + \ \mathfrak{X}_{l_1}\sum_{i=l_2}^{l_1-1}s_i\\
			&\ \  + \ (\mathfrak{X}_{l_2}+\mathfrak{X}_{l_1})\sum_{i=0}^{l_2-1}s_i   \ +\ \sum_{j\in\mathfrak{D}_0'}\mathfrak{X}_j\sum_{i=0}^{j-1}s_i. 
		\end{aligned}
	\end{equation}
	We now define
	\begin{equation}\begin{aligned}
			\zeta_0 \ = \ &\ \tilde\xi\ + \ (\mathfrak{X}_{l_2}+\mathfrak{X}_{l_1})\sum_{i=0}^{l_2-1}s_i  \ + \ \sum_{j\in\mathfrak{D}_0'}\mathfrak{X}_j\sum_{i=0}^{j-1}s_i. 
		\end{aligned} 
	\end{equation}
	For all $j\in\mathfrak{D}_0'$ , we have $\mathrm{deg}(v_j)=\mathrm{deg}(v_{j+1})=0$. Hence, by conclusion (iii) of Lemma \ref{Lemma:VerticesLongCollisions}, we deduce that $\mathfrak{X}_j$ is independent of  the cycle of $v_{l_1}$. Note that $\mathfrak{X}_{l_1}$ appears in  $\zeta_0$, together with  $\mathfrak{X}_{l_2}$, and since $\mathfrak{X}_{l_2}+\mathfrak{X}_{l_1}$ is independent of the free edge of  $v_{l_1}$, the quantity $(\mathfrak{X}_{l_2}+\mathfrak{X}_{l_1})\sum_{i=0}^{l_2-1}s_i$ is also independent of the free edge of $v_{l_1}$. As a result, $\zeta_0$
	is independent of the free edge of $v_{l_1}$. 
	
	We then split $\mathfrak{D}_1=\mathfrak{D}_1'\cup\mathfrak{D}_1''$, in which $\mathfrak{D}_1'$ denotes the set of the indices of the degree-one vertices belonging to the cycle of $v_{l_1}$ excluding $l_1$, and  $\mathfrak{D}_1''$ denotes the set of the indices of the degree-one vertices outside of the cycle of $v_{l_1}$. The total phase of the first $l_1$ time slices becomes
	
	\begin{equation}
		\label{DelayedRecollision:A1}\begin{aligned}
			& \sum_{i=0}^{l_1-1}s_i\mathrm{Re}\vartheta_i \ = \  \zeta_0 \ + \ \sum_{j\in \mathfrak{D}_1'} \left[\mathfrak{X}_j\sum_{i=0}^{j-1}s_i \ +\mathfrak{X}_{j-1}\sum_{i=0}^{j-2}s_i \right]\ + \ \mathfrak{X}_{l_1}\sum_{i=l_2}^{l_1-1}s_i\\
			& \ + \ \sum_{j\in \mathfrak{D}_1''} \left[\mathfrak{X}_j\sum_{i=0}^{j-1}s_i\ + \ \mathfrak{X}_{j-1}\sum_{i=0}^{j-2}s_i \right]\\
			\ = \ & \zeta_1 \ + \ \sum_{j\in \mathfrak{D}_1'} \left[\mathfrak{X}_j\sum_{i=0}^{j-1}s_i \ +\mathfrak{X}_{j-1}\sum_{i=0}^{j-2}s_i \right]\ + \ \mathfrak{X}_{l_1}\sum_{i=l_2}^{l_1-1}s_i,
		\end{aligned}
	\end{equation}
	where the new phase $\zeta_1$ is totally independent of the free edge of $v_{l_1}$. 
	This identity means that we can estimate the phases using exactly the same method used in the proof of Proposition \ref{Propo:CrossingGraphs}, except for all of the vertices belonging to the cycle of the degree-one vertex $v_{l_1}$, at which we need a more delicate estimate  than Lemma \ref{lemma:degree1vertex}. This estimate is performed on the total phase of the cycle. 
	
	
	We set
	\begin{equation}
		\label{DelayedRecollision:A1}\begin{aligned}
			\vartheta_{l_1}^{cyc}& \ :=  \ \sum_{j\in \mathfrak{D}_1'} \left[\mathfrak{X}_j\sum_{i=0}^{j-1}s_i \ +\mathfrak{X}_{j-1}\sum_{i=0}^{j-2}s_i \right]\ + \ \mathfrak{X}_{l_1}\sum_{i=l_2}^{l_1-1}s_i\\
			& \ =  \ \sum_{j\in \mathfrak{D}_1'} \mathfrak{X}_js_{j-1} \ + \ \mathfrak{X}_{l_1}\sum_{i=l_2}^{l_1-1}s_i.
		\end{aligned}
	\end{equation}
	Next, we will reduce $ \vartheta_{l_1}^{cyc}$ to a new quantity $	\vartheta_{l_1}^{cy}$
	\begin{equation}
		\label{DelayedRecollision:A3:a}\begin{aligned}
			\vartheta_{l_1}^{cy}& \ :=   \mathfrak{X}_{l_1}\sum_{i=l_2, i+1\notin \mathfrak{D}_1'}^{l_1-1}s_i \ = \ \mathfrak{X}_{l_1}\sum_{i\in\mathfrak{I}^o}s_i,
		\end{aligned}
	\end{equation}
	where $\mathfrak{I}^o$ is  the set of indecies $i=l_2,l_1-1, i+1\notin \mathfrak{D}_1'.$ Since the graph is delayed recollisional, $l_1-1\ne l_2$. As the cycle of $v_{l_1}$ is a delayed recollision, we have $|\mathfrak{I}^o|>1$. This reduction  is described as follows. 
	From \eqref{DelayedRecollision:A1}, since only the phase $\vartheta_{l_1}^{cyc}$ depends on the free momentum of the degree-one vertex $v_{l_1}$, while the other phase $\zeta_1 $ is independent of this momentum, we could perform the standard estimate of Proposition \ref{Propo:CrossingGraphs} for the degree-one and degree-zero vertices outside the cycle of $v_{l_1}$ and leave the estimate of  the part involving $\vartheta_{l_1}^{cyc}$ of this cycle to the final step. However, we can also perform the standard strategy for all of the degree-one vertices included in $\mathfrak{D}_1'$ as well. Thus, comparing to the  strategy of   Proposition \ref{lemma:BasicGEstimate1}, we remove all the degree zero vertices of the set $\mathfrak{I}^o$  from  the total set of degree zero vertices and perform a special treatment for the degree-one vertex $v_{l_1}$. This is the key difference between the treatment of delayed recollisional graphs and the standard graph estimates of Propositions \ref{Propo:CrossingGraphs} and \ref{lemma:BasicGEstimate1}.  To this end,  the set of all time slices $i$ includes the indices $i$ such that $v_{i+1}$ are degree-one vertices satisfying  $i+1\ne l_1$, we follow the standard treatment of degree one vertices of Propositions \ref{Propo:CrossingGraphs} and \ref{lemma:BasicGEstimate1}.  As we  remove all the degree zero vertices of the set $\mathfrak{I}^o$  from  the total set of degree zero vertices, the quantity $F_1$ defined in \eqref{eq:Aestimate0:1:A} is now modified to  include only the time slices associated to zero degree vertices  that are not in $\mathfrak{I}^o$.   The bound  \eqref{eq:Aestimate2a} for this new version of  $F_1$ is now modified by  $\lambda^n \frac{t^{n/2-(|\mathfrak{I}^o|-1)}}{[n/2-(|\mathfrak{I}^o|-1)]!}$. This process also moves all  the degree zero vertices of the set $\mathfrak{I}^o$ into the quantity $F_2$ in \eqref{eq:Aestimate0:1:A:3}.

To estimate the new quantity $F_2$, we  use the standard method of Proposition \ref{Propo:CrossingGraphs}, from the bottom of the graph to the top, until we reach $v_{l_1}$: We integrate all of the free momenta from the bottom to the top, starting from time slice $0$. 	When we meet a degree one vertex, we use Lemma \cite{lemma:degree1vertex} to  get a bound with an extra factor of $\langle \ln|\lambda| \rangle^{2+c\eth}$ for some $c>0$. This process will completely change the form of $\vartheta_{l_1}^{cyc}$, due to the following reason. From the bottom to the top, whenever we meet a degree one-vertex $v_j$, with $j\in \mathfrak{D}_1'$, we will remove $e^{{\bf i}  \mathfrak{X}_js_{j-1}}$ from the total phases.  In other  words, $ \mathfrak{X}_{j}s_{j-1}$ is now removed from $\vartheta_{l_1}^{cyc}$. We refer to Figure \ref{Fig45} for an illustration. After integrating all the degree-one vertices below $v_{l_1}$, we reduce $\vartheta_{l_1}^{cyc}$ to the  new quantity $	\vartheta_{l_1}^{cy}$, that depends on the free edge of the cycle of $v_{l_1}$.

	We will now integrated the free edge of $v_{l_1}$, using the strategy described below instead of using again  Lemma \cite{lemma:degree1vertex}. This strategy will lead to the gain of an extra factor of $\lambda^{c}$ for some $c>0$, guaranteeing the convergence of the whole graph to zero. After this, we can use the standard method of Proposition \ref{lemma:BasicGEstimate1} to integrate all degree-one vertices above $v_{l_1}$, by Lemma \cite{lemma:degree1vertex} and  get a bound with an extra factor of $\langle \ln|\lambda| \rangle^{2+c\eth}$ each time.

	\begin{figure}
		\centering
		\includegraphics[width=.49\linewidth]{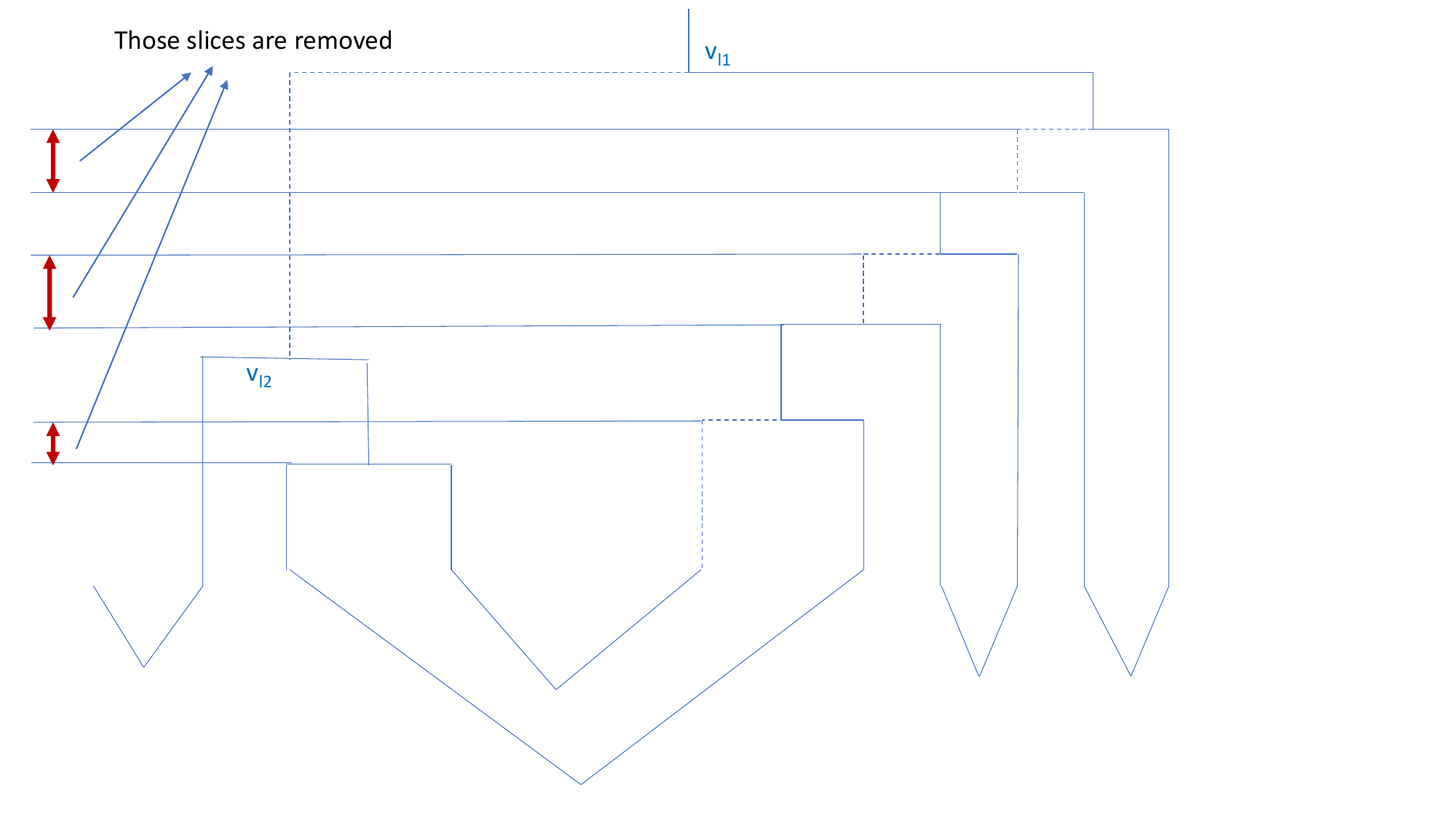}
		\caption{ In the cycle of $v_{l_1}$, which is a long delayed recollision, we remove some time slices to obtain $\vartheta_{l_1}^{cy}$.}
		\label{Fig45}
	\end{figure}
	
	To integrate the free edge of $v_{l_1}$, we isolate only the time slices from $l_2$ to $l_1-1$ that appear in the sum $\vartheta_{l_1}^{cy}$ and set
	\begin{equation}\label{DelayedRecollision:A3:1}
		\begin{aligned}
			\mathfrak{C}_1^{del} :=	&\ \int_{(\mathbb{R}_+)^{\mathfrak{I}^o}}\! \mathrm{d} \bar{s} \,  \delta\left(t-\sum_{i\in \mathfrak{I}^o} s_i\right) \int_{(\mathbb{T}^{d})^{\mathcal{I}_n^o}}\mathrm{d}\bar{k}\Delta^o_{n,\rho}(\bar k,\bar\sigma)\\
			&\times\prod_{i\in \mathfrak{I}^o}\Big[ \sigma_{i,\rho_i}\mathcal{M}( k_{i,\rho_i}, k_{i-1,\rho_i}, k_{i-1,\rho_i+1})\diamond_\ell(s_{l_1-1})\\
			&\times {\Phi}_{1,i}(\sigma_{i-1,\rho_i},  k_{i-1,\rho_i},\sigma_{i-1,\rho_i+1},  k_{i-1,\rho_i+1})\Big] 
			e^{ -{\bf i} \vartheta_{l_1}^{cy}}\mathfrak{F}^{{I}^o}, 
	\end{aligned}\end{equation}
	which is a portion of $	{Q}_{Delayed}\ell(\tau)$ containing the time slices from $0$ to $l_1-1$ that involves the index set $\mathfrak{I}^o$ of the cycle of $v_{l_1}$. The quantity ${(\mathbb{T}^{d})^{\mathcal{I}_n^o}}$ is the part of ${(\mathbb{T}^{d})^{\mathcal{I}_n'}}$ that  contains all the domain of all momenta $k_{i,j}$ of all  the time slices from $0$ to $l_1-1$ that involve the index set $\mathfrak{I}^o$ of the cycle of $v_{l_1}$. The quantity $\Delta^o_{n,\rho}(\bar k,\bar\sigma)$ is the part of $\Delta_{n,\rho}(\bar k,\bar\sigma)$ that  involves all momenta connecting to the index set $\mathfrak{I}^o$ of the cycle of $v_{l_1}$. The quantity $\mathfrak{F}^{{I}^o}$ is the component of $F_2$ that depends on the free edge of $v_{l_1}$ and that appears after we have done the integration of all the degree-one vertex before reaching $v_{l_1}$.  As the other part can be treated as before, we will therefore devote the rest of the proof to estimate $\mathfrak{C}_1^{del}$. We also set
	\begin{equation}\label{DelayedRecollision:A3:2}
		\begin{aligned}
			\mathfrak{C}_2^{del} :=	&\  \int_{(\mathbb{T}^{d})^{\mathcal{I}_n^o}}\mathrm{d}\bar{k}\Delta^o_{n,\rho}(\bar k,\bar\sigma)\sigma_{1,\rho_1} \mathcal{M}( k_{1,\rho_1}, k_{0,\rho_1}, k_{0,\rho_1+1}) \\
			&\times\prod_{i\in \mathfrak{I}^o}\Big[ \sigma_{i,\rho_i}\mathcal{M}( k_{i,\rho_i}, k_{i-1,\rho_i}, k_{i-1,\rho_i+1})\\
			&\times {\Phi}_{1,i}(\sigma_{i-1,\rho_i},  k_{i-1,\rho_i},\sigma_{i-1,\rho_i+1},  k_{i-1,\rho_i+1})\Big] 
			e^{ -{\bf i} \vartheta_{l_1}^{cy}}\mathfrak{F}^{{I}^o}. 
	\end{aligned}\end{equation}

	The quantity $\mathfrak{C}_2^{del}$ contains a delayed recollision in its form. 	To proceed further, we need to understand the structure of the recollisions and the delayed recollision formed by the cycle of $v_{l_1}$ and the  recollisions below this cycle.
	Let us consider a recollision consisting of two vertices $v_j,v_{j-1}$, in which $\mathrm{deg}v_j=1$, $\mathrm{deg}v_{j-1}=0$. Since the total phase of the two time slices associated to $v_j,v_{j-1}$ is $$\mathfrak{X}_j\sum_{i=0}^{j-1}s_i  +\mathfrak{X}_{j-1}\sum_{i=0}^{j-2}s_i \ =\ \mathfrak{X}_js_{j-1},$$ we only need to consider the time slice $s_{j-1}$ and the associated phase $\mathfrak{X}_j$. The recollision contains the following product of delta functions $$\delta(\sigma_{j,\rho_j}k_{j,\rho_j}+\sigma_{j-1,\rho_j} k_{j-1,\rho_j}+\sigma_{j-1,\rho_j+1}k_{j-1,\rho_j+1})$$$$\times\delta(\sigma_{j-1,\rho_{j-1}}k_{j-1,\rho_{j-1}}+\sigma_{j-2,\rho_{j-1} k_{j-2,\rho_{j-1}}}+\sigma_{j-2,\rho_{j-1}+1}k_{j-2,\rho_{j-1}+1}).$$ This product of delta functions finally reduces to only one delta function $$\delta(\sigma_{j,\rho_j}k_{j,\rho_j}+\sigma_{j-1,\rho_j} k_{j-1,\rho_j}+\sigma_{j-1,\rho_j+1}k_{j-1,\rho_j+1})$$ due to the pairings of the momenta. The product $$\mathcal{M}(k_{j,\rho_j}, k_{j-1,\rho_j},k_{j-1,\rho_j+1})\mathcal{M}(k_{j-1,\rho_{j-1}}, k_{j-2,\rho_{j-1}},k_{j-2,\rho_{j-1}+1})$$ is also combined into  $$\mbox{ either }-\sigma_{j,\rho_j}\sigma_{j-1,\rho_j}\mathrm{sign}k_{j,\rho_j}^1\mathrm{sign}k_{j-1,\rho_j}^1|\mathcal{W}(k_{j,\rho_j}, k_{j-1,\rho_j},k_{j-1,\rho_j+1})|^2,$$ 
	$$ -\sigma_{j,\rho_j}\sigma_{j-1,\rho_j+1}\mathrm{sign}k_{j,\rho_j}^1\mathrm{sign}k_{j-1,\rho_j+1}^1|\mathcal{W}(k_{j,\rho_j}, k_{j-1,\rho_j},k_{j-1,\rho_j+1})|^2,$$
	$$\mbox{ or }-\sigma_{j,\rho_j}^2(\mathrm{sign}k_{j,\rho_j}^1)^2|\mathcal{W}(k_{j,\rho_j}, k_{j-1,\rho_j},k_{j-1,\rho_j+1})|^2,$$ due to the pairings of the momenta.

	As a result, if a recollision is of double-cluster type, adding it to a pairing means that we change a factor $\tilde{f}$ to 
	
	\begin{equation}
		\label{DelayedRecollision:A5b}\begin{aligned}
			&
			-\iint_{(\mathbb{T}^d)^2}\!\! \mathrm{d} k_{j-1,\rho_j+1}  \mathrm{d} k_{j-1,\rho_j}\,
			\delta(\sigma_{j,\rho_j}k_{j,\rho_j}+\sigma_{j-1,\rho_j} k_{j-1,\rho_j}+\sigma_{j-1,\rho_j+1}k_{j-1,\rho_j+1}) \\ & 
			\times|\mathcal{W}(k_{j,\rho_j}, k_{j-1,\rho_j},k_{j-1,\rho_j+1})|^2\Phi_{1,j}(v_j) \\
			&\times e^{{\bf i} s_{j-1} \sigma_{j,\rho_j} \omega(k_{j,\rho_j})}e^{{\bf i} s_{j-1} \sigma_{j-1,\rho_j} \omega(k_{j-1,\rho_j})} \\
			&\times \tilde{f}(k_{j-1,\rho_j})e^{{\bf i} s_{j-1} \sigma_{j-1,\rho_j+1} \omega(k_{j-1,\rho_j+1})} \tilde{f}(k_{j-1,\rho_j+1}),
		\end{aligned}
	\end{equation}
in which we denote $\tilde{f}(k_{j-1,\rho_j+1})\delta_{k_{j-1,\rho_j+1}-k_{j-1,\rho_j+1}'}=\langle a_{k_{j-1,\rho_j+1}}a_{k_{j-1,\rho_j+1}'}^* \rangle_0$, $\tilde{f}(k_{j,\rho_j})\delta_{k_{j,\rho_j}-k_{j,\rho_j}'}=\langle a_{k_{j,\rho_j}}a_{k_{j,\rho_j}'}^* \rangle_0$ and $\tilde{f}(k_{j-1,\rho_j})\delta_{k_{j-1,\rho_j}-k_{j-1,\rho_j}'}=\langle a_{k_{j-1,\rho_j}}a_{k_{j-1,\rho_j}'}^* \rangle_0$.
	And if a recollision is of single-cluster type, adding it to an edge means that we change a factor $\tilde{f}$ to   
	\begin{equation}
		\label{DelayedRecollision:A5c}\begin{aligned}
			&
			\iint_{(\mathbb{T}^d)^2}\!\! \mathrm{d} k_{j-1,\rho_j+1}   \mathrm{d} k_{j-1,\rho_j}\,
			\delta(\sigma_{j,\rho_j}k_{j,\rho_j}+\sigma_{j-1,\rho_j} k_{j-1,\rho_j}+\sigma_{j-1,\rho_j+1}k_{j-1,\rho_j+1}) \\ & 
			\times(-\sigma_{j,\rho_j}\sigma_{j-1,\rho_j}\mathrm{sign}k_{j,\rho_j}^1\mathrm{sign}k_{j-1,\rho_j}^1)|\mathcal{W}(k_{j,\rho_j}, k_{j-1,\rho_j},k_{j-1,\rho_j+1})|^2 \Phi_{1,j}(v_j)\\
			&\times e^{{\bf i} s_{j-1} \sigma_{j,\rho_j} \omega(k_{j,\rho_j})}e^{{\bf i} s_{j-1} \sigma_{j-1,\rho_j} \omega(k_{j-1,\rho_j})} \\
			&\times \tilde{f}(k_{j-1,\rho_j+1})e^{{\bf i} s_{j-1} \sigma_{j-1,\rho_j+1} \omega(k_{j-1,\rho_j+1})} \tilde{f}(k_{j,\rho_j}),
		\end{aligned}
	\end{equation}
	and
	\begin{equation}
		\label{DelayedRecollision:A5d}\begin{aligned}
			&
			\iint_{(\mathbb{T}^d)^2}\!\! \mathrm{d} k_{j-1,\rho_j+1}   \mathrm{d} k_{j-1,\rho_j}\,
			\delta(\sigma_{j,\rho_j}k_{j,\rho_j}+\sigma_{j-1,\rho_j} k_{j-1,\rho_j}+\sigma_{j-1,\rho_j+1}k_{j-1,\rho_j+1}) \\ & 
			\times(-\sigma_{j,\rho_j}\sigma_{j-1,\rho_j+1}\mathrm{sign}k_{j,\rho_j}^1\mathrm{sign}k_{j-1,\rho_j+1}^1)|\mathcal{W}(k_{j,\rho_j}, k_{j-1,\rho_j},k_{j-1,\rho_j+1})|^2\Phi_{1,j}(v_j) \\
			&\times e^{{\bf i} s_{j-1} \sigma_{j,\rho_j} \omega(k_{j,\rho_j})}e^{{\bf i} s_{j-1} \sigma_{j-1,\rho_j} \omega(k_{j-1,\rho_j})} \\
			&\times \tilde{f}(k_{j-1,\rho_j})e^{{\bf i} s_{j-1} \sigma_{j-1,\rho_j+1} \omega(k_{j-1,\rho_j+1})} \tilde{f}(k_{j,\rho_j}).
		\end{aligned}
	\end{equation}
	
	Note that by changing $\sigma_{j,\rho_j}k_{j,\rho_j}$ $\to k_{j,\rho_j}$, $\sigma_{j-1,\rho_j} k_{j-1,\rho_j}$ $\to -k_{j-1,\rho_j}$ and $\sigma_{j-1,\rho_j+1}k_{j-1,\rho_j+1}$ $\to -k_{j-1,\rho_j+1}$, the collision operators \eqref{DelayedRecollision:A5b}, \eqref{DelayedRecollision:A5c}, \eqref{DelayedRecollision:A5d}  become
	\begin{equation}
		\label{DelayedRecollision:A5b:1}\begin{aligned}
			&
			-\iint_{(\mathbb{T}^d)^2}\!\! \mathrm{d} k_{j-1,\rho_j+1}  \mathrm{d} k_{j-1,\rho_j}\,
			\delta(k_{j,\rho_j}- k_{j-1,\rho_j}-k_{j-1,\rho_j+1})\Phi_{1,j}(v_j) \\ & 
			\times|\mathcal{W}(k_{j,\rho_j}, k_{j-1,\rho_j},k_{j-1,\rho_j+1})|^2 e^{{\bf i} s_{j-1}  \omega(k_{j,\rho_j})}e^{-{\bf i} s_{j-1}  \omega(k_{j-1,\rho_j})} \\
			&\times \tilde{f}(k_{j-1,\rho_j})e^{-{\bf i} s_{j-1} \omega(k_{j-1,\rho_j+1})} \tilde{f}(k_{j-1,\rho_j+1})\diamond_\ell(s_{l_1-1}),
		\end{aligned}
	\end{equation}
	\begin{equation}
		\label{DelayedRecollision:A5c:1}\begin{aligned}
			&
			\iint_{(\mathbb{T}^d)^2}\!\! \mathrm{d} k_{j-1,\rho_j+1}   \mathrm{d} k_{j-1,\rho_j}\,
			\delta(k_{j,\rho_j}- k_{j-1,\rho_j}-k_{j-1,\rho_j+1}) \Phi_{1,j}(v_j)\\
			& 
			\times(\mathrm{sign}k_{j,\rho_j}^1\mathrm{sign}k_{j-1,\rho_j}^1)|\mathcal{W}(k_{j,\rho_j}, k_{j-1,\rho_j},k_{j-1,\rho_j+1})|^2e^{{\bf i} s_{j-1} \omega(k_{j,\rho_j})} \\
			&\times e^{-{\bf i} s_{j-1}\omega(k_{j-1,\rho_j})} \tilde{f}(k_{j-1,\rho_j+1})e^{-{\bf i} s_{j-1}  \omega(k_{j-1,\rho_j+1})}\tilde{f}(k_{j,\rho_j}),
		\end{aligned}
	\end{equation}
	and
	\begin{equation}
		\label{DelayedRecollision:A5d:1}\begin{aligned}
			&
			\iint_{(\mathbb{T}^d)^2}\!\! \mathrm{d} k_{j-1,\rho_j+1}   \mathrm{d} k_{j-1,\rho_j}\,
			\delta(k_{j,\rho_j}- k_{j-1,\rho_j}-k_{j-1,\rho_j+1}) \Phi_{1,j}(v_j)\\ & 
			\times(\mathrm{sign}k_{j,\rho_j}^1\mathrm{sign}k_{j-1,\rho_j+1}^1)|\mathcal{W}(k_{j,\rho_j}, k_{j-1,\rho_j},k_{j-1,\rho_j+1})|^2 e^{{\bf i} s_{j-1} \omega(k_{j,\rho_j})}\\
			&\times e^{-{\bf i} s_{j-1} \omega(k_{j-1,\rho_j})}\tilde{f}(k_{j-1,\rho_j})e^{-{\bf i} s_{j-1}  \omega(k_{j-1,\rho_j+1})} \tilde{f}(k_{j,\rho_j}).
		\end{aligned}
	\end{equation}
	
	We can generalize \eqref{DelayedRecollision:A5b}-\eqref{DelayedRecollision:A5c}-\eqref{DelayedRecollision:A5d} to the forms
	\begin{equation}\begin{aligned}\label{OperatorRecollision}
			& \mathcal{C}^1_{recol}[g_1,g_2,s,\sigma_0,\sigma_1,\sigma_2](v_j)(k_0)
			=
			-\iint_{(\mathbb{T}^d)^2}\!\! \mathrm{d} k_1  \mathrm{d} k_2\,
			\delta(\sigma_0k_0+\sigma_1k_1+\sigma_2k_2) \\ & \quad
			\times|\mathcal{W}(k_{0}, k_{1},k_{2})|^2\Phi_{1,j}(v_j) e^{-{\bf i} s \sigma_0 \omega(k_0)}
			\\
			&\quad \times e^{-{\bf i} s \sigma_1 \omega(k_1)} g_1(k_1)e^{-{\bf i} s \sigma_2 \omega(k_2)} g_2(k_2),\quad  k_0\in \mathbb{T}^d\,,
		\end{aligned}
	\end{equation}
	and
	\begin{equation}\begin{aligned}\label{OperatorRecollision2}
			& \mathcal{C}^2_{recol}[g_0,g_1,s,\sigma_0,\sigma_1,\sigma_2](v_j)(k_0)
			=
			\iint_{(\mathbb{T}^d)^2}\!\! \mathrm{d} k_1  \mathrm{d} k_2\,
			\delta(\sigma_0k_0+\sigma_1k_1+\sigma_2k_2) \\ & \quad
			\times(-\sigma_0\sigma_2\mathrm{sign}k_{0}^1\mathrm{sign}k_{2}^1)|\mathcal{W}(k_{0}, k_{1},k_{2})|^2\Phi_{1,j}(v_j) e^{-{\bf i} s \sigma_0 \omega(k_0)}\\
			&\quad \times e^{-{\bf i} s \sigma_1 \omega(k_1)} g_1(k_1)e^{-{\bf i} s \sigma_2 \omega(k_2)} g_0(k_0),\quad  k_0\in \mathbb{T}^d\,,
		\end{aligned}
	\end{equation}
	in which  $s\in \mathbb{R}$, $\sigma_0,\sigma_1,\sigma_2\in \{\pm 1\}$ and $k_j =(k_j^1,\cdots,k_j^d)$ for $j=0,1,2$.

	We rewrite the two operators using the abbreviations 
	\begin{equation}
		\label{CollisionVj}
		\mathcal{C}_{recol}[v_j,s_{j-1}]=\mathcal{C}_{recol}[v_j],
	\end{equation}
	which means ``the collision operator for the recollision happpening at the vertex $v_j$ and time slice $s_{j-1}$'' (see Figure \ref{Fig36}).
	Note that the cycle of $v_{l_1}$ contains  only recollisions associated to the degree-one vertices $v_j$ and the operator $\mathcal{C}_{recol}$ can be either $\mathcal{C}_{recol}^1$ or $\mathcal{C}_{recol}^2$.

	\begin{figure}
		\centering
		\includegraphics[width=.49\linewidth]{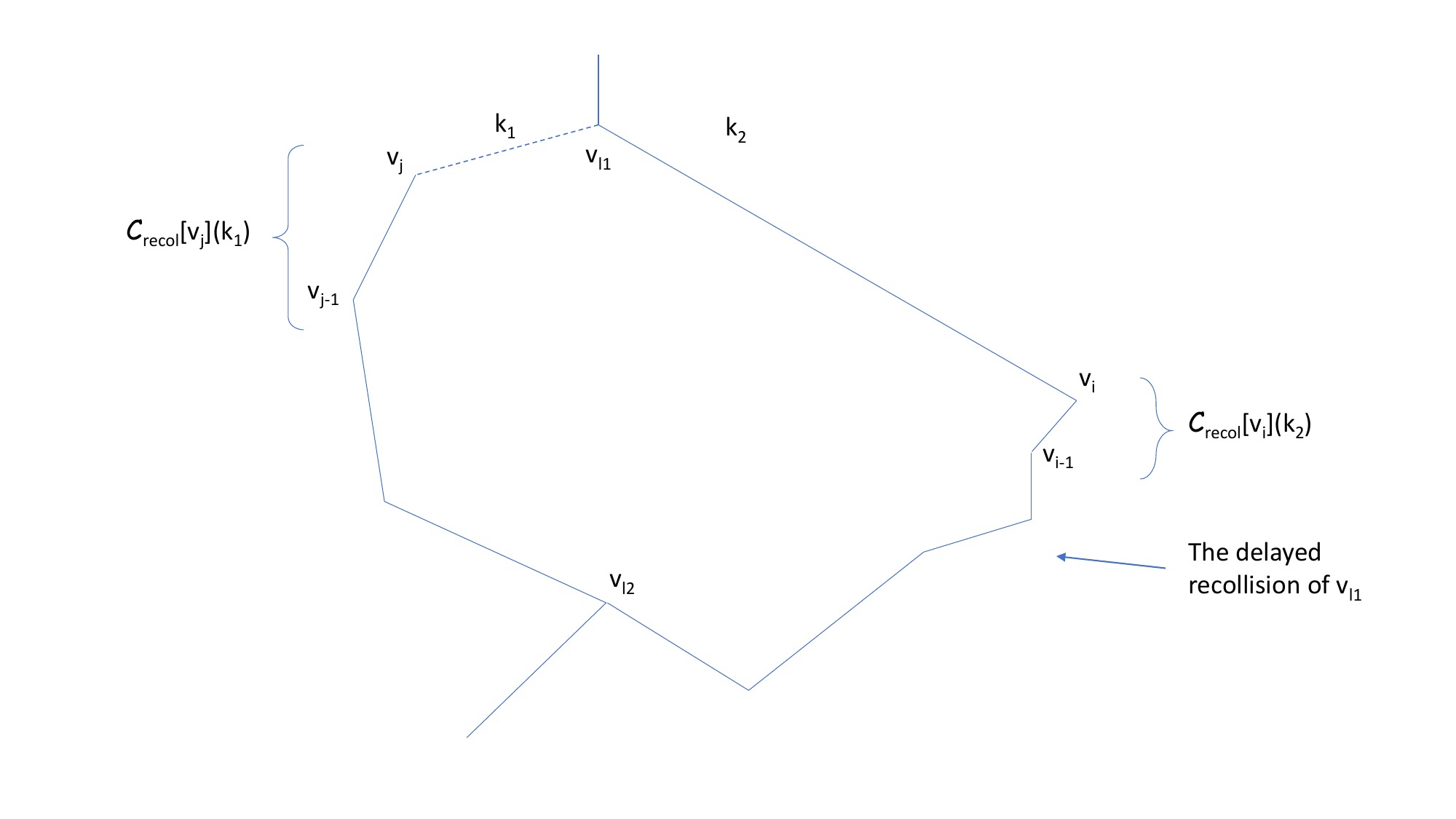}
		\caption{ This is the delayed recollision of $v_{l_1}$. In the picture, the recollisions by $v_i,v_{i-1}$ and $v_j,v_{j-1}$ are grouped into $\mathcal{C}_{recol}[v_i]$ and $\mathcal{C}_{recol}[v_j]$.}
		\label{Fig36}
	\end{figure}

	Now, let us study the delayed recollision created by the cycle of $v_{l_1}$.
	We have already shown that this cycle consists of a series of recollisions and one degree-zero vertex $v_{l_2}$. By the third conclusion  (iii) of Lemma \ref{Lemma:VerticesLongCollisions},  if we remove all of the recollisions in the manner explained in Lemma \ref{Lemma:VerticesLongCollisions}, the cycle then becomes a recollision. As a result, the delayed recollision can also be written in the form of  \eqref{OperatorRecollision} or \eqref{OperatorRecollision2}, in which the time slice $s_{j-1}$ is replaced by the sum $\sum_{i\in\mathfrak{I}^o}s_i$. We then consider two cases.
	
	{\bf Case 1: The delayed recollision can also be written in the form of  \eqref{OperatorRecollision}.} We write
	\begin{equation}
		\begin{aligned}\label{OperatorDelayedRecollision}
			&\mathcal{C}_{recol}\left[G_1,G_2,\sum_{i=l_2}^{l_1-1}s_i,\sigma_{l_1,\rho_{l_1}},\sigma_{l_1-1,\rho_{l_1}},\sigma_{l_1-1,\rho_{l_1}+1}\right](k_{l_1,\rho_{l_1}})
			\nonumber \\
			=
			&\  \iint_{(\mathbb{T}^d)^2}\!\! \mathrm{d} k_{l_1-1,\rho_{l_1}}  \mathrm{d} k_{l_1-1,\rho_{l_1}+1}\Phi_{1,l_2}(v_{l_2})\,\\
			&\times
			\delta(\sigma_{l_1,\rho_{l_1}}k_{l_1,\rho_{l_1}}+\sigma_{l_1-1,\rho_{l_1}}k_{l_1-1,\rho_{l_1}}+\sigma_{l_1-1,\rho_{l_1}}k_{l_1-1,\rho_{l_1}+1}) \\ 
			& 
			\times|\mathcal{M}(k_{l_1,\rho_{l_1}},k_{l_1-1,\rho_{l_1}},k_{l_1-1,\rho_{l_1}+1})|^2 e^{{\bf i}\Big(\sum_{i\in\mathfrak{I}^o}s_i \Big)\sigma_{l_1,\rho_{l_1}} \omega(k_{l_1,\rho_{l_1}})}\\
			&\times e^{{\bf i} \Big(\sum_{i\in\mathfrak{I}^o}s_i \Big) \sigma_{l_1-1,\rho_{l_1}}k_{l_1-1,\rho_{l_1}} \omega(k_{l_1-1,\rho_{l_1}})} G_1(k_{l_1-1,\rho_{l_1}})\\
			&\times e^{{\bf i}\Big(\sum_{i\in\mathfrak{I}^o}s_i \Big)s_i \sigma_{l_1-1,\rho_{l_1}+1} \omega(k_{l_1-1,\rho_{l_1}+1})} G_2(k_{l_1-1,\rho_{l_1}+1}),
		\end{aligned}
	\end{equation}
	in which, the expressions of $G_1,G_2$ will be explained below. We suppose that $\Phi_{1,l_2}$ is not the function $1$, since if it is the function $1$, we can replace it by $\Phi_{1,l_1}$, which is clearly not the function $1$. 
	
	In principle, $G_1,G_2$ should be the products of the recollisions $\mathcal{C}_{recol}[v_j]$. If the vertex $v_j$ is to the top vertex of a double-cluster recollision, then $\mathcal{C}_{recol}[v_j]$ takes the form $\mathcal{C}_{recol}^1[v_j]$, with $G_1,G_2$ being $\tilde{f}$. If the vertex $v_j$ is to the top vertex of a single-cluster recollision, then $\mathcal{C}_{recol}[v_j]$ takes the form $\mathcal{C}_{recol}^2[v_j]$, in which one of the functions $G_1,G_2$ is $\tilde{f}$ and the other one is $1$ (which belongs to a sequence of single-cluster recollisions).  The situation is precisely the same with the case of the leading graphs, described in Figures \ref{Fig41}, \ref{Fig47} and the collision operator
	\eqref{FinalProof:E10a:1:5}. We refer to the discussions of \ref{Fig41}, \ref{Fig47} and the collision operator
	\eqref{FinalProof:E10a:1:5} for a more precise description of the structures of $G_1,G_2$. It is clear that those products are bounded in the $L^4$-norm.  We can thus suppose that 
	\begin{equation}
		\begin{aligned}\label{OperatorRecollisionbis:1}
			\|{G}_1\|_{L^4},\   	  \|{G}_2\|_{L^4}\le \mathfrak{C}_G,
		\end{aligned}
	\end{equation}
	for some universal constant $\mathfrak{C}_G>0$ depending on the bound of all the quantities $\tilde{f}(k_{0,l})$.  Even when at the bottom of the graph, we have   $\left\langle\prod_{i=1}^{n+2}\alpha(k_{0,i},\sigma_{0,i})\right\rangle_{s} $ with $s\ne 0$  then due to Lemma \ref{Lemma:Cutoffn} and the cut-off function defined in Definition \ref{def:hbar}, we still have an $L^4$ bound.

	Next, to simplify the notations, we denote $k_{l_1,\rho_{l_1}}$, $k_{l_1-1,\rho_{l_1}}$, $k_{l_1-1,\rho_{l_1}+1}$ and $\sigma_{l_1,\rho_{l_1}},$ $\sigma_{l_1-1,\rho_{l_1}},$ $\sigma_{l_1-1,\rho_{l_1}+1}$  by $k_0,k_1,k_2$, $\sigma_0,\sigma_1,\sigma_2$. We also  set $S=\sum_{i\in\mathfrak{I}^o}s_i$ and rewrite 
	\begin{equation}\begin{aligned}\label{OperatorRecollisionbis}
			& \mathcal{C}_{recol}\left[G_1,G_2,\sum_{i\in\mathfrak{I}^o}s_i,\sigma_{l_1,\rho_{l_1}},\sigma_{l_1-1,\rho_{l_1}},\sigma_{l_1-1,\rho_{l_1}+1}\right](k_{l_1,\rho_{l_1}})\  \\
			= \ &   
			\iint_{(\mathbb{T}^d)^2}\!\! \mathrm{d} k_1  \mathrm{d} k_2\,
			\delta(\sigma_0k_0+\sigma_1k_1+\sigma_2k_2)\Phi_{1,l_2}(v_{l_2}) \\ & \quad
			\times|\mathcal{M}(k_{0}, k_{1},k_{2})|^2 e^{{\bf i} S \sigma_0 \omega(k_0)}e^{{\bf i} S \sigma_1 \omega(k_1)} G_1(k_1)e^{{\bf i} S \sigma_2 \omega(k_2)} G_2(k_2).
		\end{aligned}
	\end{equation}
	Note that as both $k_1$ and $k_2$ are associated to $v_{l_1}$, the function $\Phi_{1,l_2}(v_{l_2})$ has both $k_1$ and $k_2$ as its variables. With an abuse of notation, we can denote $\Phi_{1,o}^a(k_1)=\Phi_{1,o}^a(k_2)=\sqrt{\Phi_{1,l_2}(v_{l_2})}.$  
	
	Next, we set
	\begin{equation}\begin{aligned}\label{OperatorRecollisionbis:1}
			& \mathcal{C}_{recol}^o\left[G_1,G_2,\sum_{i=l_2}^{l_1-1}s_i,\sigma_{l_1,\rho_{l_1}},\sigma_{l_1-1,\rho_{l_1}},\sigma_{l_1-1,\rho_{l_1}+1}\right](k_{l_1,\rho_{l_1}})\ 
			\\
			= \ &   
			\iint_{(\mathbb{T}^d)^2}\!\! \mathrm{d} k_1  \mathrm{d} k_2\,
			\delta(\sigma_0k_0+\sigma_1k_1+\sigma_2k_2)\Phi_{1,o}^a(k_2) \Phi_{1,o}^a(k_1)\\ & \quad
			\times|\bar\omega(k_1)||\bar\omega(k_2)|e^{{\bf i} S \sigma_1 \omega(k_1)} G_1(k_1)e^{{\bf i} S \sigma_2 \omega(k_2)} G_2(k_2),
		\end{aligned}
	\end{equation}
	where we have used \eqref{Kernel} and \eqref{Kernel2}, and  the quantity $ |\bar\omega(k_0)| $ and the phase $\omega(k_0)$ involving $k_0$ has been removed.  
	We now estimate the  recollision operators \eqref{OperatorRecollisionbis} and \eqref{OperatorRecollisionbis:1}. To this end, we define
	$$
	\mathfrak{H}_1(x,t) :=  \int_{\mathbb{T}^d}\! \mathrm{d}  k \, |\bar\omega(k)||\bar\omega(-\sigma_0 \sigma_2 k_0-\sigma_1\sigma_2 k)|e^{{\bf i} 2 \pi x \cdot k +{\bf i} t\sigma_1 \omega(k)+ {\bf i} t\sigma_2 \omega(-\sigma_0 \sigma_2 k_0-\sigma_1\sigma_2 k)}\Phi_{1,o}^a(k)$$
	$$\times G_1(k)G_2(-\sigma_0 \sigma_2 k_0-\sigma_1\sigma_2 k),
	$$
	and
	$$
	\mathfrak{H}_2(x,t) :=  \int_{\mathbb{T}^d}\! \mathrm{d}  k \, e^{{\bf i} 2 \pi x \cdot k}{\Phi_{1,o}^a(k)}.
	$$

	By the identity 
	$$	\delta(\sigma_0k_0+\sigma_1k_1+\sigma_2k_2) \ = \ \sum_{y\in\mathbb{Z}^d}e^{{\bf i}2\pi y\cdot(\sigma_0k_0+\sigma_1k_1+\sigma_2k_2)},$$
	we write 
	\begin{equation}
		\begin{aligned}
			&	{\mathcal{C}_{recol}^o}\Big[G_1,G_2,\sum_{i=l_2}^{l_1-1}s_i,\sigma_{l_1,\rho_{l_1}},\sigma_{l_1-1,\rho_{l_1}},\sigma_{l_1-1,\rho_{l_1}+1}\Big](k_0)\\
			= &  
			\sum_{y\in\mathbb{Z}^d}e^{{\bf i}2\pi k_0\sigma_0\cdot y}		\iint_{(\mathbb{T}^d)^2}\!\! \mathrm{d} k_1  \mathrm{d} k_2 \Phi_{1,o}^a(k_1)\Phi_{1,o}^a(k_2) |\bar\omega(k_1)||\bar\omega(k_2)|  
			\\
			&\times {G}_2(-\sigma_0\sigma_2k_0-\sigma_1\sigma_2k_2) e^{{\bf i}2\pi y \cdot\sigma_1k_1+{\bf i}S\sigma_1\omega(k_1)}{G}_1(k_1) e^{{\bf i}2\pi y \cdot\sigma_2k_2+{\bf i}S\sigma_2\omega(k_2)}
			\\
			= & \sum_{y\in\mathbb{Z}^d}e^{{\bf i}2\pi k_0\sigma_0\cdot y}	\mathfrak{H}_1(\sigma_1y,S\sigma_1)
			\mathfrak{H}_{2}(\sigma_2y,S\sigma_2).
		\end{aligned}
	\end{equation} 	
	Thus, we can bound
	\begin{equation}
		\begin{aligned}\label{DelayedRecollision:A6:a1}
			&	\Big\|{\mathcal{C}_{recol}^o}\Big[G_1,G_2,\sum_{i=l_2}^{l_1-1}s_i,\sigma_{l_1,\rho_{l_1}},\sigma_{l_1-1,\rho_{l_1}},\sigma_{l_1-1,\rho_{l_1}+1}\Big](k_0)\Big\|_{L^4}\\
			\lesssim\ &  
			\Big\| \sum_{y\in\mathbb{Z}^d}e^{{\bf i}2\pi k_0\sigma_0\cdot y}\mathfrak{H}_1(\sigma_1y,S\sigma_1)
			\mathfrak{H}_{2}(\sigma_2y,S\sigma_2)\Big\|_{L^4}.
		\end{aligned}
	\end{equation} 	
	
	Next
	we bound the  terms on  the right hand side of \eqref{DelayedRecollision:A6:a1} as
	\begin{equation}
		\begin{aligned}\label{DelayedRecollision:A6:a2}
			&	\Big\|{\mathcal{C}^o_{recol}}\Big[G_1,G_2,\sum_{i\in\mathfrak{I}^o}s_i,\sigma_{l_1,\rho_{l_1}},\sigma_{l_1-1,\rho_{l_1}},\sigma_{l_1-1,\rho_{l_1}+1}\Big](k_0)\Big\|_{L^4}\
			\lesssim\ \Big\| {\mathfrak{H}}_{1}(\sigma_1\cdot,S\sigma_1) {\mathfrak{H}}_{2}(\sigma_2\cdot,S\sigma_2)\Big\|_{l^1},
		\end{aligned}
	\end{equation} 	
	which yields, by H\"older's inequality
	\begin{equation}
		\begin{aligned}\label{DelayedRecollision:A6:a3}
			&	\Big\|{\mathcal{C}^o_{recol}}\Big[G_1,G_2,\sum_{i\in\mathfrak{I}^o}s_i,\sigma_{l_1,\rho_{l_1}},\sigma_{l_1-1,\rho_{l_1}},\sigma_{l_1-1,\rho_{l_1}+1}\Big](k_0)\Big\|_{L^4}\
			\lesssim\ \Big\| {\mathfrak{H}}_{1}(\sigma_1\cdot,S\sigma_1)\Big\|_{l^8}\Big\| {\mathfrak{H}}_{2}(\sigma_2\cdot,S\sigma_2)\Big\|_{l^\frac87}.
		\end{aligned}
	\end{equation} 	

	Note that the different norms we use are defined in \eqref{Def:Norm1}-\eqref{Def:Norm2}-\eqref{Def:Norm3}-\eqref{Def:Norm4}. Thus we bound
	\begin{equation}\label{DelayedRecollision:A6:b1}
		\begin{aligned}
			\Big\|\mathfrak{C}_1^{del}\Big\|_{L^4}
			\lesssim\ &  \lambda^{-\frac2q(|\mathfrak{I}^o|-1)}\left[\int_{(\mathbb{R}_+)^{\mathfrak{I}^o}}\! \mathrm{d} \bar{s} \Big\|{\mathcal{C}_{recol}}\Big[G_1,G_2,\sum_{i\in\mathfrak{I}^o}s_i,\sigma_{l_1,\rho_{l_1}},\sigma_{l_1-1,\rho_{l_1}},\sigma_{l_1-1,\rho_{l_1}+1}\Big](k_0)\Big\|_{L^4}^p {\bf \mathbf{1}}\left(\sum_{i\in \mathfrak{I}^o} s_i\le t\right)\right]^{\frac1p}\\
			\lesssim\ &   \lambda^{-\frac2q(|\mathfrak{I}^o|-1)} \left[\int_{(\mathbb{R}_+)^{\mathfrak{I}^o}}\! \mathrm{d} \bar{s} \Big\| {\mathfrak{H}}_{1}(\sigma_1\cdot,S\sigma_1)\Big\|_{l^8}^p{\bf \mathbf{1}}\left(\sum_{i\in \mathfrak{I}^o} s_i\le t\right)\right]^{\frac1p}\Big\| {\mathfrak{H}}_{2}(\sigma_2\cdot,S\sigma_2)\Big\|_{l^\frac87},
	\end{aligned}\end{equation}
	in which $\frac1p+\frac1q=1$ and $q=\mathscr{M}'$. This inequality is attainable due to the fact that all the $|\mathfrak{I}^o|-1$ lower time slices in $\mathfrak{I}^o$ are all  degree $0$.
	
	{\bf Case 2: The delayed recollision can also be written in the form of  \eqref{OperatorRecollision}.} We write
	\begin{equation}
		\begin{aligned}\label{OperatorDelayedRecollision:bis}
			&\mathcal{C}_{recol}\Big[G_1,G_2,\sum_{i=l_2}^{l_1-1}s_i,\sigma_{l_1,\rho_{l_1}},\sigma_{l_1-1,\rho_{l_1}},\sigma_{l_1-1,\rho_{l_1}+1}\Big](k_{l_1,\rho_{l_1}})
			\nonumber \\
			=
			&\  \iint_{(\mathbb{T}^d)^2}\!\! \mathrm{d} k_{l_1-1,\rho_{l_1}}  \mathrm{d} k_{l_1-1,\rho_{l_1}+1}\Phi_{1,l_2}(v_{l_2})
			\delta(\sigma_{l_1,\rho_{l_1}}k_{l_1,\rho_{l_1}}+\sigma_{l_1-1,\rho_{l_1}}k_{l_1-1,\rho_{l_1}}+\sigma_{l_1-1,\rho_{l_1}}k_{l_1-1,\rho_{l_1}+1}) \\ 
			& 
			\times|\mathcal{M}(k_{l_1,\rho_{l_1}},k_{l_1-1,\rho_{l_1}},k_{l_1-1,\rho_{l_1}+1})|^2 e^{{\bf i}\Big(\sum_{i\in\mathfrak{I}^o}s_i \Big)\sigma_{l_1,\rho_{l_1}} \omega(k_{l_1,\rho_{l_1}})} e^{{\bf i} \Big(\sum_{i\in\mathfrak{I}^o}s_i \Big) \sigma_{l_1-1,\rho_{l_1}} \omega(k_{l_1-1,\rho_{l_1}})} \\
			&\times e^{{\bf i}\Big(\sum_{i\in\mathfrak{I}^o}s_i \Big) \sigma_{l_1-1,\rho_{l_1}+1} \omega(k_{l_1-1,\rho_{l_1}+1})}G_1(k_{l_1,\rho_{l_1}}) G_2(k_{l_1-1,\rho_{l_1}+1}),
		\end{aligned}
	\end{equation}
	or
	\begin{equation}
		\begin{aligned}\label{OperatorDelayedRecollision:bis:1}
			&\mathcal{C}_{recol}\Big[G_1,G_2,\sum_{i=l_2}^{l_1-1}s_i,\sigma_{l_1,\rho_{l_1}},\sigma_{l_1-1,\rho_{l_1}},\sigma_{l_1-1,\rho_{l_1}+1}\Big](k_{l_1,\rho_{l_1}})
			\nonumber \\
			=
			&\  \iint_{(\mathbb{T}^d)^2}\!\! \mathrm{d} k_{l_1-1,\rho_{l_1}}  \mathrm{d} k_{l_1-1,\rho_{l_1}+1}\Phi_{1,l_2}(v_{l_2})
			\delta(\sigma_{l_1,\rho_{l_1}}k_{l_1,\rho_{l_1}}+\sigma_{l_1-1,\rho_{l_1}}k_{l_1-1,\rho_{l_1}}+\sigma_{l_1-1,\rho_{l_1}}k_{l_1-1,\rho_{l_1}+1}) \\ 
			& 
			\times|\mathcal{M}(k_{l_1,\rho_{l_1}},k_{l_1-1,\rho_{l_1}},k_{l_1-1,\rho_{l_1}+1})|^2 e^{{\bf i}\Big(\sum_{i\in\mathfrak{I}^o}s_i \Big)\sigma_{l_1,\rho_{l_1}} \omega(k_{l_1,\rho_{l_1}})}G_1(k_{l_1-1,\rho_{l_1}})\\
			&\times e^{{\bf i} \Big(\sum_{i\in\mathfrak{I}^o}s_i \Big) \sigma_{l_1-1,\rho_{l_1}}k_{l_1-1,\rho_{l_1}} \omega(k_{l_1-1,\rho_{l_1}})} e^{{\bf i}\Big(\sum_{i\in\mathfrak{I}^o}s_i \Big)s_i \sigma_{l_1-1,\rho_{l_1}+1} \omega(k_{l_1,\rho_{l_1}})} G_2(k_{l_1-1,\rho_{l_1}+1}).
		\end{aligned}
	\end{equation}
	The estimate \eqref{DelayedRecollision:A6:b1} can be proved by the same arguments, as  thus, \eqref{DelayedRecollision:A6:b1} is  still satisfied in the second case considered here.

	To proceed further, we will use a $TT^*$-type argument, starting from \eqref{DelayedRecollision:A6:b1}. We choose $\tilde G(x,\bar{s})$ to be a test function in $L^{q}([0,t]^{\mathfrak{J}^o},l^\frac87(\mathbb{Z}^d))$, i.e. $\Big[\int_{[0,t]^{\mathfrak{J}^o}}\mathrm{d}\bar{s}$ $\Big[\sum_{x\in\mathbb{Z}^d} |\tilde G(x,\bar{s})|^\frac87\Big]^{7q/8}\Big]^{\frac{1}{q}}<\infty$, 
	and estimate
	\begin{equation}\label{DelayedRecollision:A6:b2}
		\begin{aligned}
			& \left|\sum_{x\in\mathbb{Z}^d} \int_{(\mathbb{R}_+)^{\mathfrak{I}^o}}\! \mathrm{d} \bar{s}  {\mathfrak{H}}_{1}(\sigma_1x,S\sigma_1)\tilde G(x,\bar{s})
			{\bf \mathbf{1}}\left(\sum_{i\in \mathfrak{I}^o} s_i\le t\right) \right|\\
			\ = \ &\Big|\sum_{x\in\mathbb{Z}^d} \int_{(\mathbb{R}_+)^{\mathfrak{I}^o}}\! \mathrm{d} \bar{s}   \int_{\mathbb{T}^d}\! \mathrm{d}  k \,  |\bar\omega(k)||\bar\omega(-\sigma_0 \sigma_2 k_0-\sigma_1\sigma_2 k)|e^{{\bf i} 2 \pi x \cdot k +{\bf i} S\sigma_1 \omega(k)+ {\bf i} S\sigma_2 \omega(-\sigma_0 \sigma_2 k_0-\sigma_1\sigma_2 k)}
			\Phi_{1,o}^a(k)\\
			&\times G_1(k)G_2(-\sigma_0 \sigma_2 k_0-\sigma_1\sigma_2 k)\tilde G(x,\bar{s})
			{\bf \mathbf{1}}\Big(\sum_{i\in \mathfrak{I}^o} s_i\le t\Big) \Big|.
	\end{aligned}\end{equation}
	Next, we will  study the $L^2$ norm in $k$
	\begin{equation}\label{DelayedRecollision:A6:b3}
		\begin{aligned}
			& \left\|\sum_{x\in\mathbb{Z}^d} \int_{(\mathbb{R}_+)^{\mathfrak{I}^o}}\! \mathrm{d} \bar{s}  e^{{\bf i} 2 \pi x \cdot k +{\bf i} S\sigma_1 \omega(k)}|\bar\omega(k)|e^{{\bf i} 2 \pi x \cdot k +{\bf i} S\sigma_1 \omega(k)+ {\bf i} S\sigma_2 \omega(-\sigma_0 \sigma_2 k_0-\sigma_1\sigma_2 k)}\tilde G(x,\bar{s}) |\bar\omega(-\sigma_0 \sigma_2 k_0-\sigma_1\sigma_2 k)|
		\right.\\
		&\times\left.\Phi_{1,o}^a(k)	{\bf \mathbf{1}}\left(\sum_{i\in \mathfrak{I}^o} s_i\le t\right) \right\|_{L^2},
	\end{aligned}\end{equation}
	which can be expanded as 
	\begin{equation}\label{DelayedRecollision:A6:b4}
		\begin{aligned}
			& \left|\int_{\mathbb{T}^d}\mathrm{d}k\sum_{x'\in\mathbb{Z}^d} \int_{(\mathbb{R}_+)^{\mathfrak{I}^o}}\! \mathrm{d} \bar{s}'\sum_{x\in\mathbb{Z}^d} \int_{(\mathbb{R}_+)^{\mathfrak{I}^o}}\! \mathrm{d} \bar{s}  e^{{\bf i} 2 \pi x \cdot k +{\bf i} S\sigma_1 \omega(k)+ {\bf i} S\sigma_2 \omega(-\sigma_0 \sigma_2 k_0-\sigma_1\sigma_2 k)} |\bar\omega(k)|^2|\bar\omega(-\sigma_0 \sigma_2 k_0-\sigma_1\sigma_2 k)|^2 \right.\\
			&\left.\times |\Phi_{1,o}^a(k)|^2 \tilde G(x,\bar{s})e^{-{\bf i} 2 \pi x' \cdot k -{\bf i} S'\sigma_1 \omega(k)- {\bf i} S'\sigma_2 \omega(-\sigma_0 \sigma_2 k_0-\sigma_1\sigma_2 k)} \overline{\tilde G(x',\bar{s}')}
			{\bf \mathbf{1}}\left(\sum_{i\in \mathfrak{I}^o} s_i\le t\right)	{\bf \mathbf{1}}\left(\sum_{i\in \mathfrak{I}^o} s_i'\le t\right)\right|,
	\end{aligned}\end{equation}
	in which $S'=\sum_{i\in \mathfrak{I}^o} s_i'$.
	We  now introduce 
	\begin{equation}
		\begin{aligned}
			\mathfrak{H}_{1}'(x,t)\ :=\  &\int_{\mathbb{T}^d}\! \mathrm{d}  ke^{{\bf i} 2 \pi x \cdot k +{\bf i} S\sigma_1 \omega(k)+ {\bf i} S\sigma_2 \omega(-\sigma_0 \sigma_2 k_0-\sigma_1\sigma_2 k)} |\bar\omega(k)|^2|\bar\omega(-\sigma_0 \sigma_2 k_0-\sigma_1\sigma_2 k)|^2|\Phi_{1,o}^a(k)|^2, \end{aligned}
	\end{equation} 
	    so that we can continue the estimate of the above quantity  
	\begin{equation}\label{DelayedRecollision:A6:b5}
		\begin{aligned}
			& \left|\sum_{x'\in\mathbb{Z}^d} \int_{(\mathbb{R}_+)^{\mathfrak{I}^o}}\! \mathrm{d} \bar{s}'\sum_{x\in\mathbb{Z}^d} \int_{(\mathbb{R}_+)^{\mathfrak{I}^o}}\! \mathrm{d} \bar{s}  | \tilde G(x,\bar{s})||{\tilde G(x',\bar{s}')}|\right.\\
			&\left.\times|\mathfrak{H}'_1(x-x',S\sigma_1-S'\sigma_1)|
			{\bf \mathbf{1}}\left(\sum_{i\in \mathfrak{I}^o} s_i\le t\right) 	{\bf \mathbf{1}}\left(\sum_{i\in \mathfrak{I}^o} s_i'\le t\right) \right|\\
			\lesssim\	& \left[\int_{[0,t]^{\mathfrak{J}^o}}\mathrm{d}\bar{s} \left[\sum_{x\in\mathbb{Z}^d} |\tilde G(x,\bar{s})|^\frac87\right]^{7q/8}\right]^\frac2q \left[\int_{(\mathbb{R})^{\mathfrak{I}^o}}\! \mathrm{d} \bar{s}	{\bf \mathbf{1}}\left( -t\le \sum_{i\in \mathfrak{I}^o} s_i\le t\right) \|\mathfrak{H}'_1(\cdot,S\sigma_1)\|_4^\frac{p}{2} \right]^\frac{2}{p},
	\end{aligned}\end{equation}
	with $\frac{1}{p}+\frac{1}{q}=1$. 
		Since by Lemma \ref{Lemm:AnotherKernel2}, $$\Big\| \mathfrak{H}_1'(\cdot,{S\sigma_1})\Big\|_4\lesssim \frac{\langle \ln |\lambda|\rangle^c}{\langle S\rangle^{-\frac{1}{10}-}},$$
		for some $c>0$ that varies from lines to lines,
	we deduce that 
	\begin{equation}\label{DelayedRecollision:A6:b5}
		\begin{aligned}
			& \Big|\sum_{x'\in\mathbb{Z}^d} \int_{(\mathbb{R}_+)^{\mathfrak{I}^o}}\! \mathrm{d} \bar{s}'\sum_{x\in\mathbb{Z}^d} \int_{(\mathbb{R}_+)^{\mathfrak{I}^o}}\! \mathrm{d} \bar{s}  | \tilde G(x,\bar{s})||{\tilde G(x',\bar{s}')}|\\
			&\times|\mathfrak{H}'_1(x-x',S\sigma_1-S'\sigma_1)|
			{\bf \mathbf{1}}\left(\sum_{i\in \mathfrak{I}^o} s_i\le t\right) 	{\bf \mathbf{1}}\left(\sum_{i\in \mathfrak{I}^o} s_i'\le t\right) \Big|\\
			\lesssim\	& \langle \ln |\lambda|\rangle^c\left[\int_{[0,t]^{\mathfrak{J}^o}}\mathrm{d}\bar{s} \left[\sum_{x\in\mathbb{Z}^d} |\tilde G(x,\bar{s})|^\frac87\right]^{7q/8}\right]^\frac{2}{q} \left[\int_{(\mathbb{R})^{\mathfrak{I}^o}}\! \mathrm{d} \bar{s}	{\bf \mathbf{1}}\left(-t\le \sum_{i\in \mathfrak{I}^o} s_i\le t\right) {\langle S\rangle^{-\frac{p}{20}-}}\right]^\frac2p
	\end{aligned}\end{equation}

	{On the other hand} 
	\begin{equation}
		\begin{aligned}\label{eq:crossingraph14:1}
			\|\mathfrak{H}_{2}\big\|_\frac43\ = \  &\left\{\sum_{m\in\mathbb{Z}^d} \left|\int_{\mathbb{T}^{d}}\mathrm{d}{k}_{*}{\Phi_{1,o}^b(k_*)}e^{{\bf i}2\pi m\cdot k_*}\right|^\frac87\right\}^\frac78\\
			\ = \ & \left\{\sum_{m\in\mathbb{Z}^d\backslash\{0\}} \left|\int_{\mathbb{T}^{d}}\mathrm{d}{k}_{*}\frac{\Phi_{1,o}^b(k_*)}{|{\bf i}2\pi m|^{2d}}\Delta^{d}\Big(e^{{\bf i}2\pi m\cdot k_*}\Big)\right|^\frac87+\left|\int_{\mathbb{T}^{d}}\mathrm{d}{k}_{*}{\Phi_{1,o}^b }(k_*)\right|^\frac87\right\}^\frac78\\
			\ = \ & \left\{\sum_{m\in\mathbb{Z}^d\backslash\{0\}} 
			\frac{1}{|2\pi m|^{16d/7}}\left|\int_{\mathbb{T}^{d}}\mathrm{d}{k}_{*}\Delta^{d}\Big[{{\Phi_{1,o}^b(k_*)}}\Big]e^{{\bf i}2\pi m\cdot k_*}\right|^\frac43+\left|\int_{\mathbb{T}^{d}}\mathrm{d}{k}_{*}{\Phi_{1,o}^b(k_*)} \right|^\frac87\right\}^\frac78\\
			\ \le \ & \left\{\sum_{m\in\mathbb{Z}^d\backslash\{0\}} 
			\frac{1}{|2\pi m|^{16d/7}}\right\}^\frac78\left|\int_{\mathbb{T}^{d}}\mathrm{d}{k}_{*}\Big|\Delta^{d}\Big[\Phi_{1,o}^b(k_*)\Big]\Big|\right|  +\left|\int_{\mathbb{T}^{d}}\mathrm{d}{k}_{*}{\Phi_{1,o}^b(k_*)} \right| \ \le \ \langle\ln|\lambda|\rangle^{\mathfrak{C}_{\mathfrak{H}_2}}.
		\end{aligned}
	\end{equation}
	with $\mathfrak{C}_{\mathfrak{H}_2}>0$, where the last inequality follows from the property of the cut-off function.

	Combining \eqref{DelayedRecollision:A6:b1},\eqref{DelayedRecollision:A6:b5} and \eqref{eq:crossingraph14:1} yields

	%

	\begin{equation}\label{DelayedRecollision:A6:b6}
		\begin{aligned}
			\Big\|\mathfrak{C}_1^{del}\Big\|_{L^2}
			\lesssim\ &   	&\langle\ln|\lambda|\rangle^{\mathfrak{C}_{\mathfrak{H}_{1,2}}}\lambda^{-\frac2q(|\mathfrak{I}^o|-1)} \left[\int_{(\mathbb{R})^{\mathfrak{I}^o}}\! \mathrm{d} \bar{s}	{\mathbf{1}}\left( \sum_{i\in \mathfrak{I}^o} s_i\le t\right) \langle S\rangle^{-\big( \frac{p}{4}-\big)}\right]^\frac{1}{p},
	\end{aligned}\end{equation}
	for some ${\mathfrak{C}_{\mathfrak{H}_{1,2}}}>0.$
	
	By using the same strategy proposed in the proof of Proposition \ref{lemma:BasicGEstimate1}, we have another estimate 
	\begin{equation}\label{DelayedRecollision:A4}\begin{aligned}
			&		\limsup_{D\to\infty} 
			\Big\|{Q}_{Delayed}^\ell\Big\|_{L^2} \\ \le\  & 4^n e^{T_*}\frac{T_*^{\bar{n}}}{(\bar{n})!}\lambda^{n-2\bar{n}}\mathfrak{C}_{Q_{1,1}}^{n}\langle \ln\lambda\rangle^{\mathfrak{C}''_{Q_{Delayed}}}
			\left[\int_{(\mathbb{R})^{\mathfrak{I}^o}}\! \mathrm{d} \bar{s}	{\bf \mathbf{1}}\left( \sum_{i\in \mathfrak{I}^o} s_i\le t\right) \langle S\rangle^{-\big(\frac{p}{4}-\big)}\right]^\frac{2}{p},\end{aligned}
	\end{equation}
	in which $\bar{n}+(|\mathfrak{I}^o|-1)=n/2$. We now estimate, 
	\begin{equation}\begin{aligned}\label{DelayedRecollision:A6a}
			& \int_{(\mathbb{R})^{\mathfrak{I}^o}}\! \mathrm{d} \bar{s}	{\bf \mathbf{1}}\left(\sum_{i\in \mathfrak{I}^o} s_i\le t\right) {\left\langle \sum_{i\in \mathfrak{I}^o} s_i\right\rangle^{-\big(\frac{p}{20}-\big)}}\
			\lesssim \  \int_0^t\!  \mathrm{d} s \, s^{|\mathfrak{I}^o|-1} \langle{s}\rangle^{-\big(\frac{p}{4}-\big)}\\
			&\lesssim \    \int_0^t\!  \mathrm{d} s \, s^{|\mathfrak{I}^o|-1-\frac{9}{16}} \langle{s}\rangle^{-\big(\frac{17}{16}-\big)} \ \lesssim \ 
			t^{|\mathfrak{I}^o|-1-\frac{9}{16}},
	\end{aligned}\end{equation}
	for $p \ge 200$. 
	
	Combining \eqref{DelayedRecollision:A4} and \eqref{DelayedRecollision:A6a} yields 
	\begin{equation}\label{DelayedRecollision:A6}\begin{aligned}
					\limsup_{D\to\infty} 
			\Big\|{Q}_{Delayed}^\ell\Big\|_{L^2}\ \le\  &  e^{T_*}\frac{T_*^{\bar{n}}}{(\bar{n})!}\lambda^{\frac{9}{8p}}\mathfrak{C}_{Q_{1,1}}^{n}  \langle\ln\lambda\rangle^{\mathfrak{C}''_{Q_{Delayed}}},\end{aligned}
	\end{equation}
	which gives a gain of $\lambda^{\frac{9}{8p}}$ in the total bound. The above gain  of $\lambda^{\frac{9}{8p}}$  is attainable   as the graph is delayed recollisional and  $l_1-1\ne l_2$ and this ensures the convergence of the graph to $0$.
\end{proof}
\subsection{Reducing to $\mathrm{iCL}_2$ ladder graphs}
Using the results in the previous subsections, we then have the following improvement of Proposition  \ref{Lemma:Qmain}, in which we restrict the summation of $ Q_1$, $ Q_2$, $Q_3$ and $ Q_4$ only to terms that correspond to $\mathrm{iCL}_2$ ladder  graphs.

\begin{proposition}\label{Propo:ReduceLadderGraph}
	Suppose that $t>0$ and $t=\mathcal{O}(\lambda^{-2})$. 	We define

	\begin{equation}\label{Propo:ReduceLadderGraph:2}
		\mathfrak{Q}_1\ = \ \sum_{n=0}^{\mathfrak{N}-1}\sum_{S\in\mathcal{P}_{pair}(I_{n+2})}\mathfrak{G}_{1,n}'(S,t,k,\sigma,\Gamma),
	\end{equation}
	\begin{equation}\label{Propo:ReduceLadderGraph:3}
		\mathfrak Q_{2} 
		\ = \ \sum_{n=0}^{\mathfrak{N}-1}\sum_{S\in\mathcal{P}_{pair}(\{1,\cdots,n+2\})}\mathfrak{G}_{2,n}'(S,t,k,\sigma,\Gamma),
	\end{equation}
	\begin{equation}\label{Propo:ReduceLadderGraph:3}
		\mathfrak Q_{3} 
		\ = \ \sum_{n=1}^{\mathfrak{N}}\sum_{S\in\mathcal{P}_{pair}(\{1,\cdots,n+2\})}\int_{0}^t\mathrm{d}s_0\mathfrak{G}_{3,n}'(S,s_0,t,k,\sigma,\Gamma),
	\end{equation}
	\begin{equation}\label{Propo:ReduceLadderGraph:4}
	\mathfrak Q_{4} 
	\ = \  \sum_{S\in\mathcal{P}_{pair}^1(\{1,\cdots,\mathfrak{N}+2\})}\mathfrak{G}_{4,\mathfrak{N}}'(S,t,k,\sigma,\Gamma),
\end{equation}
in which $\mathcal{P}^1_{pair}(\{1,\cdots,n+2\})$ denotes the subset of $\mathcal{P}^0_{pair}(\{1,\cdots,n+2\})$, where we only count the terms that correspond to $\mathrm{iCL}_2$ ladder graphs. We denote  terms whose graphs are $\mathrm{iCL}_2$ ladder graphs by $\mathfrak{G}_{1,n}'$, $\mathfrak{G}_{2,n}'$, $\mathfrak{G}_{3,n}'$, $\mathfrak{G}_{4,n}'$. In other words,  $\mathfrak{Q}_1$, $\mathfrak{Q}_2$, $\mathfrak{Q}_3$, $\mathfrak{Q}_4$ are sum of terms coming from $Q^1,Q^2,Q^3,Q^4$ that only correspond to ladder graphs. Define the modifications $ Q_*^\ell$ by adding the additional quantities $\prod_{j\in I_{d1}}e^{-\ell s_{j-1}}$ into all of its graphs as in Proposition \ref{lemma:BasicGEstimate1}.

	in which $\mathcal{P}_{pair}(\{1,\cdots,n+2\})$ denotes the set of all partition $S$ of $\{1,\cdots,n+2\}$, such that for all elements $A\in S$,  then $|A|= 2$ and the graph is not singular. Moreover, in $\mathfrak{G}_{1,n}'(S,t,k,\sigma,\Gamma)$, $\mathfrak{G}_{2,n}'(S,t,k,\sigma,\Gamma)$, $\mathfrak{G}_{3,n}'(S,s_0,t,k,\sigma,\Gamma)$, $\mathfrak{G}_{4,n}'$ we only count the terms that correspond to $\mathrm{iCL}_2$ ladder graphs of ${Q}^{2,pair}$, ${Q}^{3,pair}$, ${Q}^{4,pair}$,  defined in Proposition \ref{Lemma:Qmain}. Therefore, $\mathfrak Q_1$, $\mathfrak Q_2$, $\mathfrak Q_4$  and $\mathfrak Q_3$ are the restriction of ${Q}^{1,pair}$, ${Q}^{2,pair}$, ${Q}^{3,pair}$, ${Q}^{4,pair}$, to $\mathrm{iCL}_2$ ladder graphs. We denote the difference between ${Q}^{1,pair}$, ${Q}^{2,pair}$, ${Q}^{3,pair}$, ${Q}^{4,pair}$ and $\mathfrak Q_1$, $\mathfrak Q_2$, $\mathfrak Q_3$ and $\mathfrak Q_4$ as $\mathfrak Q_1^{c}$, $\mathfrak Q_2^{c}$, $\mathfrak Q_3^{c}$ and $\mathfrak Q_4^{c}$. We have

\begin{equation}\label{Propo:ReduceLadderGraph:1}
	\begin{aligned}
		&\lim_{\lambda\to 0}\lim_{D\to\infty}\Big\|\mathfrak Q_1^{c}+\mathfrak Q_2^{c}+\mathfrak Q_3^{c}+\mathfrak Q_4^{c}\Big\|_{L^2}
		= \  0,\end{aligned}
\end{equation}

\end{proposition}
\begin{proof}The proof follows directly from Proposition \ref{Propo:CrossingGraphs} and Proposition \ref{Propo:LeadingGraphs}.
\end{proof}

\section{The proof of the main theorem}

In Proposition \ref{Propo:ReduceLadderGraph}, we have reduced the whole Duhamel expansions to only three ``ladder'' terms $\mathfrak Q_1$, $\mathfrak Q_2$ and $\mathfrak Q_3$.   Below, we will show that $\mathfrak Q_1$ has the main contribution. The second term $\mathfrak Q_2$ is indeed negligible due to the fact that the partial time integration parameters  $\varsigma_n$ is only non-zero for $n$ going from $\mathfrak{N}/4$ to $\mathfrak{N}$. The last term $\mathfrak Q_3$   vanishes as the truncation function $\Phi_{0,n}$, defined in \eqref{CutoffPhii}, vanishes in the limit of $\lambda\to 0$. The final term $\mathfrak Q_4$ also vanish since it contains only the set of diagrams with $n=\mathfrak N$.

Let us start by  analyzing  $\mathfrak{Q}_1$ in details. We only need to study
\begin{equation}
	\begin{aligned} \label{FinalProof:E1}
		& \ \sum_{n=0}^{[\mathfrak{N}/4]-1}\sum_{S\in\mathcal{P}_{pair}(\{1,\cdots,n+2\})}\mathfrak{G}_{1,n}'(S,t,k,\sigma,\Gamma).
\end{aligned}\end{equation}

Let us consider the analytic expression of an $\mathrm{iCL}_2$ ladder graph,  with $n=2q$ vertices. According to Lemma \ref{Lemma:Recollisions}, this graph is formed by iteratively appling  the $\mathrm{iC}_2^r$ recollisions. Since $n=2q\le [\mathfrak{N}/4]-1$, it follows that $\varsigma_{2 q-i}=0$ for all $i\in\{0,\cdots,2q\}$, due to \eqref{Def:Para3}. Thus, we do not need to worry about the soft-partial time integration parameters  $\varsigma_{2 q-i}=0$ for all $i\in\{0,\cdots,2q\}$.
Now, we notice that   the difference $|\Phi_{1,i}-1|$   vanishes as $\lambda$ tends to $0$ when $i$ is even. We write one term in $\mathfrak{Q}_1$ as
\begin{equation}
	\begin{aligned}\label{FinalProof:E1:a}
		&(-1)^q \lambda^{2 q}
		\int_{(\Lambda^*)^{\mathcal{I}_{2q}}}\mathrm{d}\bar{k} \Delta_{2q,\rho}(\bar k,\bar\sigma) \prod_{i=1}^{2q}\Big[\sigma_{i,\rho_i}\mathcal{M}( k_{i,\rho_i}, k_{i-1,\rho_i}, k_{i-1,\rho_i+1})\Big] \Phi_{1,i}\\
		&\times 
		\prod_{\{i,j\}\in S} \tilde{f}(k_{0,i},0)
		\int_{(\mathbb{R}_+)^{I_{\{0,\cdots,2q\}}}}\! \mathrm{d} \bar{s} \,
		\delta\left(t-\sum_{i=0}^{2 q}
		s_i\right) 
		\prod_{i=1}^{2 q} e^{-s_i \tau_i}
		\prod_{i=1}^{q} e^{-{\bf i} s_{2 i-1} \mathfrak{X}_{2 i}}
		\, \\
		&=\ (-1)^q \lambda^{2 q}
		\int_{(\Lambda^*)^{\mathcal{I}_{2q}}}\mathrm{d}\bar{k} \Delta_{2q,\rho}(\bar k,\bar\sigma) \prod_{i=1}^{2q}\Big[\sigma_{i,\rho_i}\mathcal{M}( k_{i,\rho_i}, k_{i-1,\rho_i}, k_{i-1,\rho_i+1})\Big] \\
		&\times 
		\prod_{\{i,j\}\in S} \tilde{f}(k_{0,i},0)
		\int_{(\mathbb{R}_+)^{I_{\{0,\cdots,2q\}}}}\! \mathrm{d} \bar{s} \,
		\delta\left(t-\sum_{i=0}^{2 q}
		s_i\right) 
		\prod_{i=1}^{q} e^{-{\bf i} s_{2 i-1} \mathfrak{X}_{2 i}}
		\, 
	\end{aligned}
\end{equation}
in which all $\tau_i$ now vanish.  The quantity $\prod_{\{i,j\}\in S} {f}(k_{0,i},0)$ indicates the pairings at the first time slice: The two momenta $k_{0,i}$ and $k_{0,j}$ are paired, forming $\tilde{f}(k_{0,i},0)$.  The sum of all of all of the   quantities \eqref{FinalProof:E1:a} give the   $f_{leading}(\tau)$ in the main theorem. We set $$f_{nonleading}(\tau)=f(\tau)-f_{leading}(\tau),$$
which vanish in the limit of $\lambda\to0$.  We now perform the resonance broadening idea for \eqref{FinalProof:E1:a}, with a broadening size $\ell\ge 0$
\begin{equation}
	\begin{aligned}\label{FinalProof:E1}
		&(-1)^q \lambda^{2 q}
		\int_{(\Lambda^*)^{\mathcal{I}_{2q}}}\mathrm{d}\bar{k} \Delta_{2q,\rho}(\bar k,\bar\sigma) \prod_{i=1}^{2q}\Big[\sigma_{i,\rho_i}\mathcal{M}( k_{i,\rho_i}, k_{i-1,\rho_i}, k_{i-1,\rho_i+1})\Big] \\
		&\times 
		\prod_{\{i,j\}\in S} \tilde{f}(k_{0,i},0)
		\int_{(\mathbb{R}_+)^{I_{\{0,\cdots,2q\}}}}\! \mathrm{d} \bar{s} \,
		\delta\left(t-\sum_{i=0}^{2 q}
		s_i\right) 
		 \prod_{j\in I_{d1}}e^{-s_{j-1}\ell}
		\prod_{i=1}^{q} e^{-{\bf i} s_{2 i-1} \mathfrak{X}_{2 i}}
		\, \\
		&=\ (-1)^q \lambda^{2 q}
		\int_{(\Lambda^*)^{\mathcal{I}_{2q}}}\mathrm{d}\bar{k} \Delta_{2q,\rho}(\bar k,\bar\sigma) \prod_{i=1}^{2q}\Big[\sigma_{i,\rho_i}\mathcal{M}( k_{i,\rho_i}, k_{i-1,\rho_i}, k_{i-1,\rho_i+1})\Big] \\
		&\times 
		\prod_{\{i,j\}\in S} \tilde{f}(k_{0,i},0)
		\int_{(\mathbb{R}_+)^{I_{\{0,\cdots,2q\}}}}\! \mathrm{d} \bar{s} \,
		\delta\left(t-\sum_{i=0}^{2 q}
		s_i\right) 
		\prod_{i=1}^{q} e^{-{\bf i} s_{2 i-1} \mathfrak{X}_{2 i}-s_{2 i-1}\ell}
		\, 
	\end{aligned}
\end{equation}
where the set $I_{d1}$ is defined in Proposition \ref{Propo:CrossingGraphs}.
Let us now study the time integration, which reads 
\begin{equation}
	\begin{aligned}\label{FinalProof:E5}
		&  \int_{(\mathbb{R}_+)^{I_{\{0,\cdots,2q\}}}}  \mathrm{d} \bar{s}
		e^{-{\bf i} s_{2 i-1} \mathfrak{X}_{2 i}-s_{2 i-1}\ell}
		\,
		\delta\left(t-\sum_{j=0}^{q} s_{2j} - \sum_{i=1}^{q} s_{2i-1}\right)
		\\ & 
		=  \int_{(\mathbb{R}_+)^{\{1,\cdots,q\}}}\! \prod_{i=1}^{q}\mathrm{d} s_{2 i-1}
		e^{-{\bf i} s_{2 i-1} \mathfrak{X}_{2 i}-s_{2 i-1}\ell}
		\mathbf{1}\left(\sum_{i=1}^{q} s_{2i-1} \le t\right)
		\frac{1}{q!} \left(t-\sum_{i=1}^{q} s_{2i-1}\right)^q\, ,
	\end{aligned}
\end{equation}
in which we have integrated  all of the variables $s_{2i}$ using \eqref{eq:Aestimate1a}. 
Combining \eqref{FinalProof:E1} and \eqref{FinalProof:E5} yields
\begin{equation}
	\begin{aligned}\label{FinalProof:E6}
		&(-1)^q \lambda^{2 q} 
		\int_{(\Lambda^*)^{\mathcal{I}_{2q}}}\mathrm{d}\bar{k} \Delta_{2q,\rho}(\bar k,\bar\sigma) \prod_{i=1}^{2q}\Big[\sigma_{i,\rho_i}\mathcal{M}( k_{i,\rho_i}, k_{i-1,\rho_i}, k_{i-1,\rho_i+1})\Big]\prod_{\{i,j\}\in S} \tilde{f}(k_{0,i},0) \\
		&\times 
		\int_{(\mathbb{R}_+)^{\{1,\cdots,q\}}}\! \prod_{i=1}^{q}\mathrm{d} s_{2 i-1}
		\prod_{i=1}^{q}  e^{-{\bf i} s_{2 i-1} \mathfrak{X}_{2 i}-s_{2 i-1}\ell}
		\mathbf{1}\left(\sum_{i=1}^{q} s_{2i-1} \le t \right)
		\frac{1}{q!} \left(t-\sum_{i=1}^{q} s_{2i-1}\right)^q.
	\end{aligned}
\end{equation}
Next, we need to sum over all possible ladders and take the limits. In the Feynman diagram associated to \eqref{FinalProof:E6} there are in total $q+1$ pairing clusters.  We denote this sum by $\mathfrak{C}^\lambda_q(t)$.  To write down the explicit form of this sum, we define an input ``$q+1$-correlation function'', which represents all of the $q+1$ pairings at the bottom of the graph \begin{equation}
	\mathfrak{L}_{q+1}(k_{1}',\cdots,k_{q+1}'):=\prod_{\{i,j\}\in S} {f}(k_{0,i},0), \end{equation}
in which $k_{1}',\cdots,k_{q+1}'$ are the momenta $\{k_{0,i}\}_{\{i,j\}\in S}$  with a different way of indexing.  We recall that at the bottom of the graph there are in total $q+1$ pairing clusters $\{i,j\}$, those pairings are denoted by $S$. The quantity $\tilde{f}(k_{0,i},0)$ appears due to the fact that $k_{0,i}
$ and $k_{0,j}$ are paired.
We  define
\begin{equation}
	\begin{aligned}\label{FinalProof:E9}
		&\mathcal{C}^\lambda_{i,q+1,\ell}(s,k_1',\cdots,k'_i,\cdots,k'_{q})=   \int_{\Lambda^*}\int_{\Lambda^*}
		\mathrm{d}k'\mathrm{d}k_{q+1}'\\
		&\times|\mathcal{M}(k_{i}',k',k_{q+1}')|^2  \Big[e^{{\bf i}s(\omega(k')+\omega(k_{q+1}')-\omega(k_{i}')-s\ell} +e^{-{\bf i}s(\omega(k')+\omega(k_{q+1}')-\omega(k_{i}'))-s\ell}\Big]\\
		&\times\delta(k'+k_{q+1}'-k_{i}')\Big( \mathfrak{L}_{q+1}(k_1',\cdots,k',\cdots,k'_{q+1})-\mathfrak{L}_{q+1}(k_1',\cdots,k'_i,\cdots,k'_{q+1})\mathrm{sign}(k_{i}')\mathrm{sign}(k_{q+1}')\\
		& -\mathfrak{L}_{q	+1}(k_1',\cdots,k'_i,\cdots,k')\mathrm{sign}(k_{i}')\mathrm{sign}(k')\Big),
\end{aligned}\end{equation}
in which $k'$ takes the positions of $k'_i$ in $\mathfrak{L}_{q+1}(k_1',\cdots,k',$ $\cdots,k'_{q+1})$ and of $k'_{q+1}$ in $\mathfrak{L}_{q+1}(k_1',\cdots,k'_i,\cdots,k')$, and set 

\begin{equation}
	\begin{aligned}\label{FinalProof:E10:a}
		&\mathcal{C}^\lambda_{q+1,\ell}(s,k_1',\cdots,k'_i,\cdots,k'_{q})=   \sum_{i=1}^{q}\mathcal{C}^\lambda_{i,q+1,\ell}(s,k_{i}').
\end{aligned}\end{equation}
The operator \eqref{FinalProof:E10:a} is created by summing operators of the types  \eqref{DelayedRecollision:A5c:1} and \eqref{DelayedRecollision:A5d:1}, in which   $\Phi_{1,j}$ are already  replaced by $1$ as discussed above. Under the influence of \eqref{FinalProof:E10:a} (and \eqref{DelayedRecollision:A5c:1}, \eqref{DelayedRecollision:A5d:1}), we can see that the number of momenta has been reduced from $q+1$ (in $\mathfrak{L}_{q+1}(k_1',\cdots,k',\cdots,k'_{q+1})$) to $q$ (in $ \mathcal{C}^\lambda_{q+1}(s,k_1',\cdots,k'_i,\cdots,k'_{q})$). Thus, we need to iteratively continue this procedure, until we reach the original momentum. To perform this iteration, we start with $m+1=q$, and define the ``$m+1$-correlation function''
\begin{equation}
	\mathfrak{L}_{m+1}(k_1',\cdots,k'_i,\cdots,k'_{m+1}) \ = \ \mathcal{C}^\lambda_{m+2,\ell}(s,k_1',\cdots,k'_i,\cdots,k'_{m+1}),
\end{equation}
as well as the operators acting on the new correlation functions 
\begin{equation}
	\begin{aligned}\label{FinalProof:E9:bis}
		&\mathcal{C}^\lambda_{i,m+1,\ell}(s,k_1',\cdots,k'_i,\cdots,k'_{m})=   \int_{\Lambda^*}\int_{\Lambda^*}
		\mathrm{d}k'\mathrm{d}k_{m+1}'\\
		&\times|\mathcal{M}(k_{i}',k',k_{q+1}')|^2  \Big[e^{{\bf i}s(\omega(k')+\omega(k_{m+1}')-\omega(k_{i}')-s\ell} +e^{-{\bf i}s(\omega(k')+\omega(k_{m+1}')-\omega(k_{i}'))-s\ell}\Big]\\
		&\times\delta(k'+k_{m+1}'-k_{i}')\Big( \mathfrak{L}_{m+1}(k_1',\cdots,k',\cdots,k'_{m+1})-\mathfrak{L}_{m+1}(k_1',\cdots,k'_i,\cdots,k'_{m+1})\mathrm{sign}(k_{i}')\mathrm{sign}(k_{m+1}')\\
		& -\mathfrak{L}_{m	+1}(k_1',\cdots,k'_i,\cdots,k')\mathrm{sign}(k_{i}')\mathrm{sign}(k')\Big),
\end{aligned}\end{equation}
 $k'$ takes the positions of $k'_i$ in $\mathfrak{L}_{m+1}$ $(k_1',\cdots,k',\cdots,k'_{m+1})$ and of $k'_{m+1}$ in $\mathfrak{L}_{m+1}(k_1',$ $\cdots,k'_i,\cdots,k')$,  and  

\begin{equation}
	\begin{aligned}\label{FinalProof:E10:a:bis}
		&\mathcal{C}^\lambda_{m+1,\ell}(s,k_1',\cdots,k'_i,\cdots,k'_{m})=   \sum_{i=1}^{m}\mathcal{C}^\lambda_{i,m+1,\ell}(s,k_{i}').
\end{aligned}\end{equation}
This  procedure will be iterated all the way to  the original momentum denoted, without loss of generality, by $k_1'$. This means the final operator would be $\mathcal{C}^\lambda_{2}$. Hence, the explicit form of $\mathfrak{C}^\lambda_q(t)$, {sum over all possible ladders},  can be written, at the end of the iteration, as
\begin{equation}
	\begin{aligned}\label{FinalProof:E10}
		\mathfrak{C}^\lambda_{q,\ell}(t)\ = \ &(-1)^q \lambda^{2 q}
		\int_{(\mathbb{R}_+)^{\{1,\cdots,q\}}}\!\mathrm{d} s_{1} \cdots \mathrm{d} s_{2 q-1} 
		[\mathcal{C}^\lambda_{2,\ell}(s_1)\cdots\mathcal{C}^\lambda_{q+1,\ell}(s_{2q-1})]\\
		&\times
		\mathbf{1}\left(\sum_{i=1}^{q} s_{2i-1} \le t \right)
		\frac{1}{q!} \left(t-\sum_{i=1}^{q} s_{2i-1}\right)^q.
	\end{aligned}
\end{equation}

By the change of variables $s_{2i-1}\to \lambda^{-2} s_{2i-1}$ and $t\to \tau\lambda^{-2}$, we write
\begin{equation}
	\begin{aligned}\label{FinalProof:E10a}
		\mathfrak{C}^\lambda_{q,\ell}(\tau)\ = \ &
		(-1)^q \lambda^{-2q}\int_{(\mathbb{R}_+)^{\{1,\cdots,q\}}}\!\mathrm{d} s_{1} \cdots \mathrm{d} s_{2 q-1} 
		[\mathcal{C}^\lambda_{2,\ell}(\lambda^{-2 }s_1)\cdots\mathcal{C}^\lambda_{q+1,\ell}(\lambda^{-2 }s_{2q-1})]\\
		&\times
		\mathbf{1}\left(\sum_{i=1}^{q} s_{2i-1} \le \tau \right)
		\frac{1}{q!} \left(\tau-\sum_{i=1}^{q} s_{2i-1}\right)^q\,.
	\end{aligned}
\end{equation}

As in Proposition \ref{lemma:BasicGEstimate1}, we can pass to the limit $D\to\infty$ and obtain $\mathfrak{C}^\lambda_{q,\infty,\ell}$, introduced in \eqref{FinalProof:E10a:1}.

We  need to  take the limit $\lambda\to 0$ by iteratively applying Lemma \ref{Lemma:Resonance1}. 
We thus set   
\begin{equation}
	\begin{aligned}
		\mathfrak{C}^\lambda_{q,\infty,a,\ell}(\tau) \ := \ &\lambda^{-2q}
		\int_{(\mathbb{R}_+)^{\{1,\cdots,q\}}}\!\mathrm{d} s_{1} \cdots \mathrm{d} s_{2 q-1} 
		\mathcal{C}^\lambda_{2,\infty,\ell}(\lambda^{-2 }s_1)\cdots\mathcal{C}^\lambda_{q+1,\infty,\ell}(\lambda^{-2 }s_{2q-1})\\
		&\ \  \ \times
		\mathbf{1}\left(\sum_{i=1}^{q} s_{2i-1} \le \tau \right) \left(\tau-\sum_{i=1}^{q} s_{2i-1}\right)^q.
	\end{aligned}
\end{equation}

We now bound using H\"older's inequality

\begin{equation}
	\begin{aligned}\label{FinalProof:E10a:1:1c:g}
		&	\|\mathfrak{C}^\lambda_{q,\infty,a,\ell}\|_{L^4(\mathbb{T}^d)}
		\\
		\lesssim\ &\mathfrak{C}_{final,1}^{q}T_*^{q+1}
		\sup_{\tau'\in[0,T_*]}\left[
		\int_{(\mathbb{R}_+)^{\{1,\cdots,q\}}}\!\mathrm{d} s_{1} \cdots \mathrm{d} s_{2 q-1} 
		\lambda^{-2q}\|\mathcal{C}^\lambda_{2,\infty,0}(\lambda^{-2 }s_1)\cdots\mathcal{C}^\lambda_{q+1,\infty,0}(\lambda^{-2 }s_{2q-1})\|_{L^4(\mathbb{T}^d)}^\mathscr{M}(\tau')\right.\\
		&\ \ \ \ \ \ \left. \times
		\mathbf{1}\left(\sum_{i=1}^{q} s_{2i-1} \le \tau'  \right)\right]^\frac{1}{\mathscr{M}},
	\end{aligned}
\end{equation}
for some constant $\mathscr{M}>0$.
We define 
\begin{equation}
	\begin{aligned}\label{FinalProof:E10a:1:3}
		\mathfrak{C}^\lambda_{q,\infty,b}(\tau)\ := \ & 
		\mathcal{C}^\lambda_{2,\infty,0}(\lambda^{-2 }s_1)\cdots\mathcal{C}^\lambda_{q+1,\infty,0}(\lambda^{-2 }s_{2q-1}),
	\end{aligned}
\end{equation}
and
\begin{equation}
	\begin{aligned}\label{FinalProof:E10a:1:4}
		\mathfrak{C}^\lambda_{q,\infty,c}(\tau)\ := \ & \lambda^{-2q}
		\int_{(\mathbb{R}_+)^{\{1,\cdots,q\}}}\!\mathrm{d} s_{1} \cdots \mathrm{d} s_{2 q-1} 
		\Big\|\mathfrak{C}^\lambda_{q,\infty,b}\Big\|_{L^4(\mathbb{T}^d)}^\mathscr{M}
		\mathbf{1}\left(\sum_{i=1}^{q} s_{2i-1} \le \tau \right).
	\end{aligned}
\end{equation}
As discussed above, the operator $\mathfrak{C}^\lambda_{q,\infty,b}(\tau)$   is an iterative application of the $\mathrm{iC}_2^r$ recollisions. Thus,  $\mathfrak{C}^\lambda_{q,\infty,b}(\tau)$  can be iterated by applying ladder operators of the form
\begin{equation}\label{FinalProof:E10a:1:5}
	\begin{aligned} 
		& 
		{Q}^{ladder}_{s,\sigma_0,\sigma_1,\sigma_2}[\mathbf F_0,\mathbf F_1,\mathbf F_2](k_0):=\iint_{(\mathbb{T}^d)^2}\!\! \mathrm{d} k_1  \mathrm{d} k_2 \,
		\delta(\sigma_0k_0+\sigma_1k_1+\sigma_2k_2) e^{{\bf i}s\sigma_0 \omega(k_0)+{\bf i}s\sigma_1 \omega(k_1)+{\bf i} s \sigma_2 \omega(k_2)}\\
		&\times  |\sin(2\pi k_0^1)|\mathbf{F}_0(k_0)\mathbf {F}_1(k_1)\mathbf F_2(k_2) |\sin(2\pi k_1^1)||\sin(2\pi k_2^1)|,
	\end{aligned}
\end{equation}
where $\mathbf F_0,\mathbf F_1,\mathbf F_2$ can be $1, f_0$  or ${Q}^{ladder}_{s',\sigma_0',\sigma_1',\sigma_2'}$, where the parameters $s',\sigma_0',\sigma_1',\sigma_2'$ come from the previous iterations. {These operators }can be classified into two main types. 

\smallskip

{\bf Type 1:} If the time $s$ corresponds to a single-cluster recollision, the operator has the form
\begin{equation}\label{FinalProof:single-cluster}
	\begin{aligned} 
		& 
		{Q}^{ladder}_{s,\sigma_0,\sigma_1,\sigma_2}[\mathbf F_0,1,\mathbf F_2](k_0), \qquad \mbox{ or } \qquad  {Q}^{ladder}_{s,\sigma_0,\sigma_1,\sigma_2}[\mathbf F_0,\mathbf F_2,1](k_0)
	\end{aligned}
\end{equation}
in which  the last function is the product of several operators of the type ${Q}^{ladder}_{s',\sigma_0',\sigma_1',\sigma_2'}$ and ${f}_0$ or $1$, which  come from the previous time slices. The last function, {which we assume without loss of generality} to be  $\mathbf F_2(k)$ has two explicit representations. 
\begin{itemize}
	\item[(i)] The first representation of $\mathbf F_2(k)$   takes the form
	
	\begin{equation}\label{FinalProof:product1}
		\begin{aligned} 
			\mathbf P_1(k) \ =	\ & 
			\mathbf	P_1^0(k)\prod_{l=1}^{m} {Q}^{ladder}_{s^l,\sigma_0^l,\sigma_1^l,\sigma_2^l}\Big[1,1,\mathbf P_1^l\Big](k),
		\end{aligned}
	\end{equation}
	where $\mathbf	P_1^0(k)$ is ${f}_0$, with $m\in\mathbb{N}$, $m\ge 0$ and   the functions $\mathbf	P_1^l(k)$  are either  ${f}_0$ or ladder operators obtained from the previous iterations, also of the types \eqref{FinalProof:product1}-\eqref{FinalProof:product2}, $l=1,\cdots, m$. Figure \ref{Fig41} gives an illustration of this situation. 
	\item[(ii)] The second representation of $\mathbf F_2(k)$   takes the form
	\begin{equation}\label{FinalProof:product2}
		\begin{aligned} 
			\mathbf P_2(k) \ =	\ & 
			{Q}^{ladder}_{s^1,\sigma_0^1,\sigma_1^l,\sigma_2^l}\Big[1,\mathbf P_1^0,\mathbf P_2^0\Big](k)\prod_{l=1}^{m} {Q}^{ladder}_{s^l,\sigma_0^l,\sigma_1^l,\sigma_2^l}\Big[1,1,\mathbf P_2^l\Big](k),
		\end{aligned}
	\end{equation}
\end{itemize}
in which ${Q}^{ladder}_{s^1,\sigma_0^1,\sigma_1^l,\sigma_2^l}[1,\mathbf P_1^0,\mathbf P_2^0]$ corresponds to a double-cluster recollision, which will be described right below and  $m\in\mathbb{N}$, $m\ge 0$. The functions  $\mathbf	P_1^0(k)$, $\mathbf	P_2^l(k)$ are ${f}_0$ or ladder operators obtained from the previous iterations, also of the types \eqref{FinalProof:product1}-\eqref{FinalProof:product2}, $l=0,\cdots, m$.  Figure \ref{Fig47} gives an illustration of this situation. 
\begin{figure}
	\centering
	\includegraphics[width=.49\linewidth]{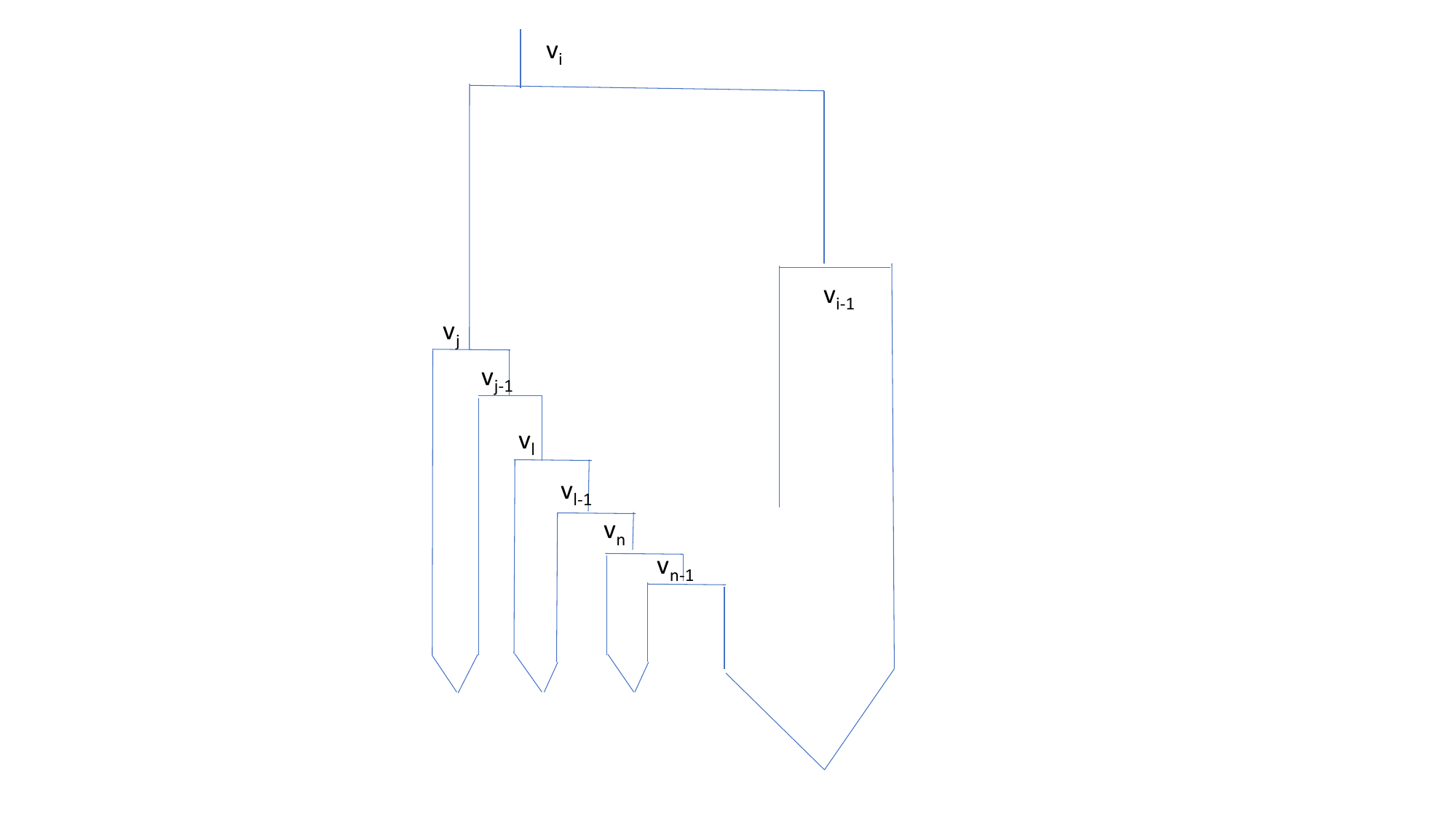}
	\caption{This picture gives an example of the ladder operator ${Q}^{ladder}_{s,\sigma_0,\sigma_1,\sigma_2}[1,1,\mathbf F_2,s]$. In this case $\mathbf F_1$ has to be $1$ as it corresponds to the single-cluster recollision associated to the vertex $v_i$. The function $\mathbf F_2$ a product of the ladder operators formed by the recollisions of the vertices $v_j,v_l,v_n$ and ${f}_0$. }
	\label{Fig41}
\end{figure}

\begin{figure}
	\centering
	\includegraphics[width=.49\linewidth]{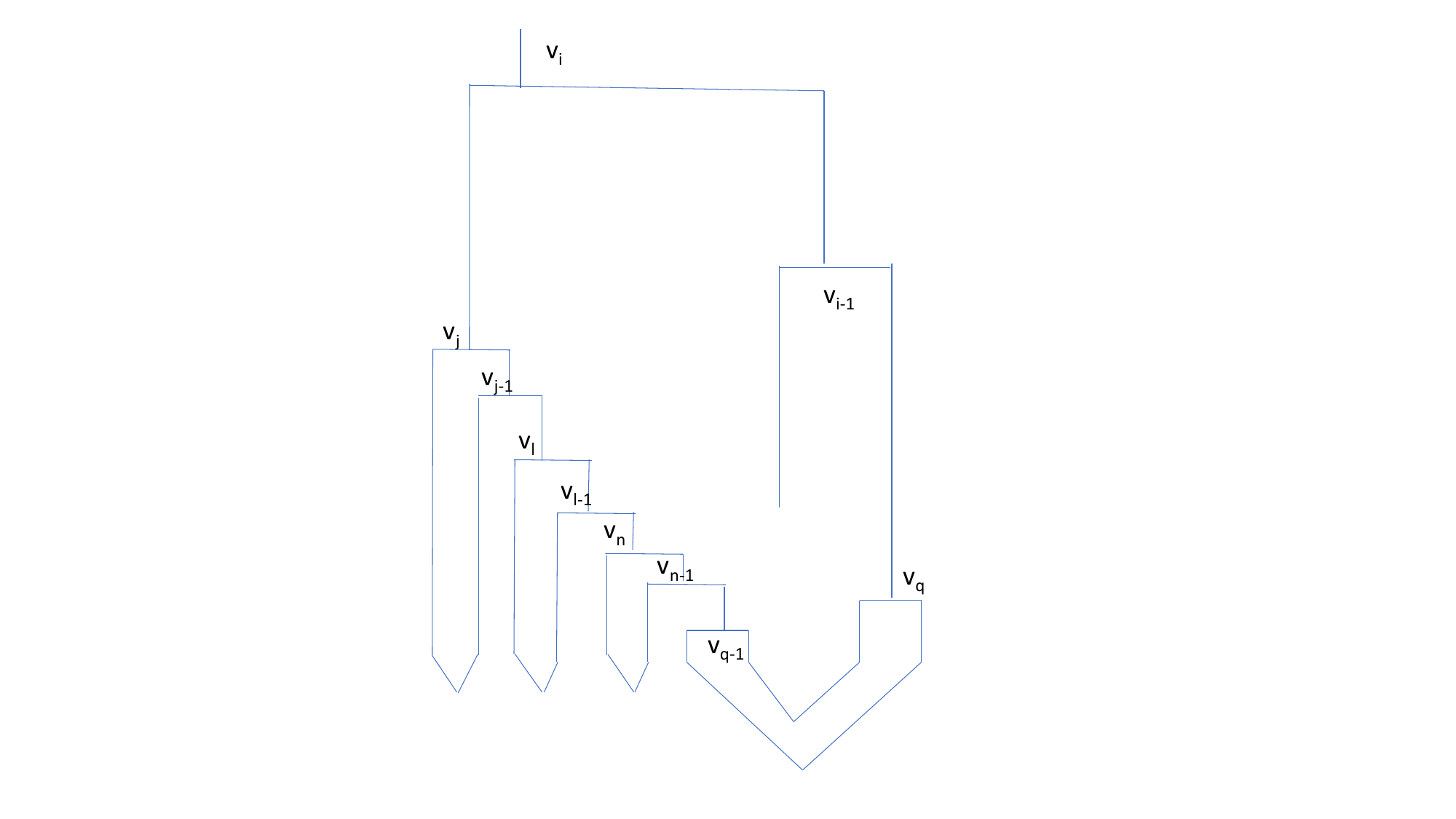}
	\caption{This picture gives an example of the ladder operator ${Q}^{ladder}_{s,\sigma_0,\sigma_1,\sigma_2}[1,1,\mathbf F_2,s]$. In this case $\mathbf F_1$ has to be $1$ as it corresponds to the single-cluster recollision associated to the vertex $v_i$. The function $\mathbf F_2$ a product of the ladder operators formed by the recollisions of the vertices $v_j,v_l,v_n$ and $v_q$. }
	\label{Fig47}
\end{figure}

\smallskip
{\bf Type 2:} If the time $s$ corresponds to a double-cluster recollision, the operator has the form
\begin{equation}\label{FinalProof:double-cluster}
	\begin{aligned} 
		& 
		{Q}^{ladder}_{s,\sigma_0,\sigma_1,\sigma_2}[1,\mathbf F_1,\mathbf F_2](k_0),
	\end{aligned}
\end{equation}
in which the function $\mathbf F_0$ has to be $1$ and the two functions $\mathbf F_1,\mathbf F_2$ also have the forms \eqref{FinalProof:product1} and \eqref{FinalProof:product2}. Figure \ref{Fig46} gives an illustration of this situation.

\begin{figure}
	\centering
	\includegraphics[width=.49\linewidth]{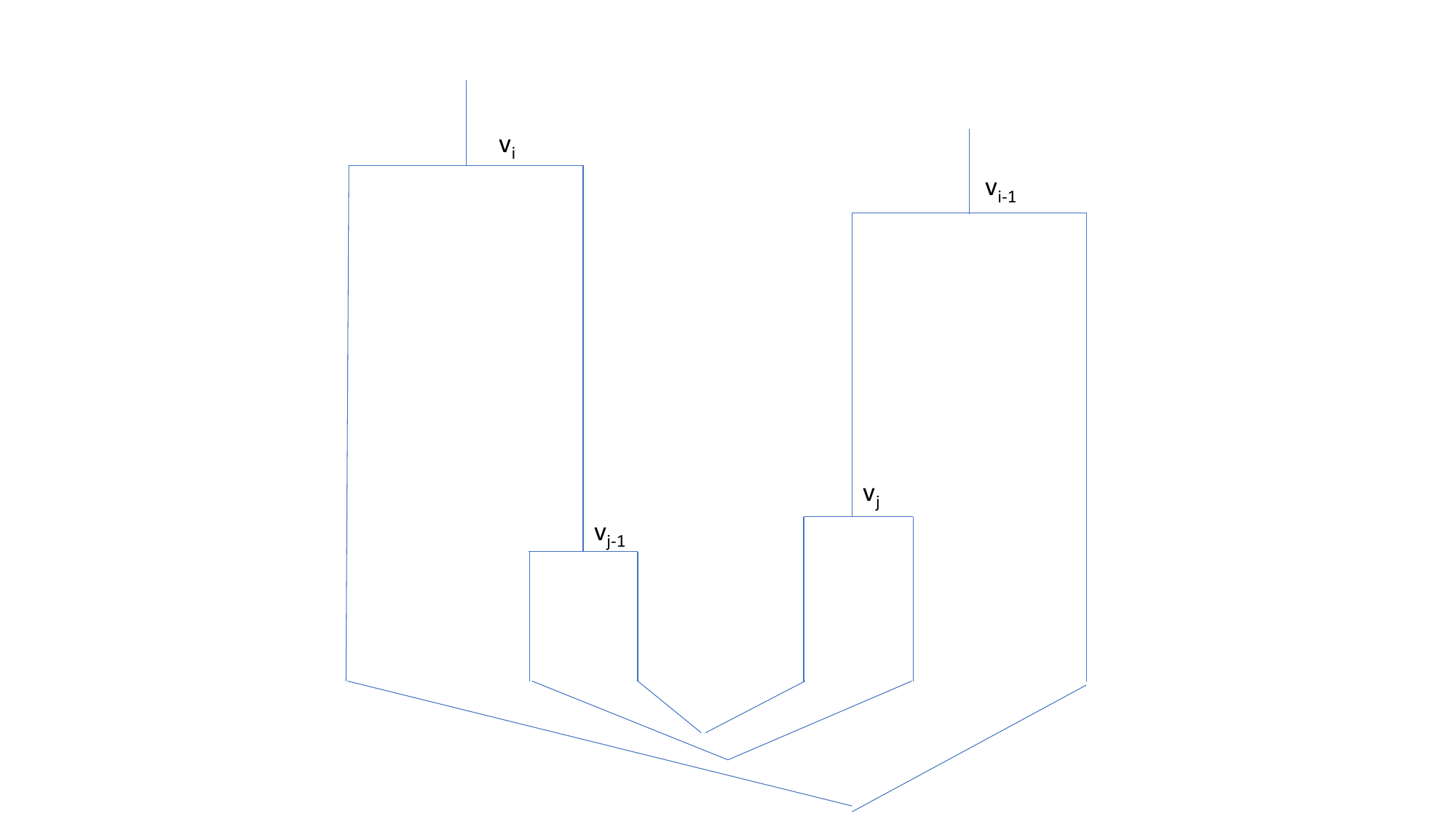}
	\caption{This picture gives an example of the ladder operator ${Q}^{ladder}_{s,\sigma_0,\sigma_1,\sigma_2}[1,\mathbf F_1,\mathbf F_2,s]$. In this case $\mathbf F_0$ has to be $1$ as it corresponds to the double-cluster recollision associated to the vertex $v_i$. The function $\mathbf F_1$ is indeed ${f}_0$  as this corresponds to  a  pairing. The function $\mathbf F_2$ is ${Q}^{ladder}_{s',\sigma_0,\sigma_1,\sigma_2}[1,\mathbf F_1',\mathbf F_2',s]$ as it is associated to the double-cluster recollision of $v_j$. The two functions $\mathbf F_1',\mathbf F_2'$ are ${f}_0$ as this corresponds to  pairings.}
	\label{Fig46}
\end{figure}
The above structure is applied iteratively from the top to the bottom of the graph, yielding the  form of $\mathfrak{C}^\lambda_{q,\infty,c}$.  To obtain the bound for $\mathfrak{C}^\lambda_{q,\infty,c}$, we will interatively apply Lemma \ref{Lemma:Resonance1}  to all of the ladder operator  from the top to the bottom of the graph. The procedure is described as follow.

{\bf Strategy  (I).} If we encounter a ladder operator of the type \eqref{FinalProof:double-cluster},   which is not included in a product of the type \eqref{FinalProof:product2} with $m\ge1$, the general strategy is to apply \eqref{Lemma:Resonance1:2:bis1} for $F_1 =\mathbf F_1$ and $F_2 =\mathbf F_2$.

{\bf Strategy (II).}  If we encounter a ladder operator of the type \eqref{FinalProof:single-cluster}, which is not included in a product of the type \eqref{FinalProof:product1} with $m\ge1$, and if  $\mathbf F_2$ is the function that is different from $1$, the general strategy is to apply \eqref{Lemma:Resonance1:2:bis2} for $F_2=\mathbf F_2$.
\smallskip

{\bf Strategy	(III).} After applying (I) and (II), we will encounter products of the types \eqref{FinalProof:product1} and \eqref{FinalProof:product2}, which we need to control in the  $L^4$-norm.  The strategies to control those products will be discussed below.

\smallskip
{\it Strategy	(III.1).} 	 First, we provide a treatment for the product \eqref{FinalProof:product1}. We bound

\begin{equation}\label{FinalProof:E10a:1:6:3}
	\begin{aligned} 
		\|	\mathbf P_1(k) \|_{L^4(\mathbb{T}^d)} \ \le	\ & 
		\| {f}_0(k) \|_{L^4(\mathbb{T}^d)}\prod_{l=1}^{m}\Big\| {Q}^{ladder}_{s^l,\sigma_0^l,\sigma_1^l,\sigma_2^l}\Big[1,1,\mathbf P_2^l\Big]\Big\|_{L^{\infty}(\mathbb{T}^d)}.
	\end{aligned}
\end{equation}
Next, we apply \eqref{Lemma:Resonance1:9} to bound $\Big\| {Q}^{ladder}_{s^l,\sigma_0^l,\sigma_1^l,\sigma_2^l}\Big[1,1,\mathbf P_2^l\Big](k))\Big\|_{L^{\infty}(\mathbb{T}^d)}$ by $\|\mathbf P_2^l\|_{L^4(\mathbb{T}^d)}$, so that in the next iteration, (I), (II) or (III) can be reused.

\smallskip
{\it Strategy	(III.2).}  Now, we provide a treatment for the product \eqref{FinalProof:product2}. We bound

\begin{equation}\label{FinalProof:E10a:1:6:4}
	\begin{aligned} 
		\|	\mathbf P_2(k)\|_{L^4(\mathbb{T}^d)} \ \le	\ & 
		\Big\| {Q}^{ladder}_{s^1,\sigma_0^1,\sigma_1^l,\sigma_2^l}\Big[1,\mathbf P_1^0,\mathbf P_2^0\Big](k) \Big\|_{L^4(\mathbb{T}^d)}\prod_{l=1}^{m}\Big\| {Q}^{ladder}_{s^l,\sigma_0^l,\sigma_1^l,\sigma_2^l}\Big[1,1,\mathbf P_2^l\Big](k))\Big\|_{L^{\infty}(\mathbb{T}^d)}.
	\end{aligned}
\end{equation}
Next, we apply \eqref{Lemma:Resonance1:9}  to bound $\Big\| {Q}^{ladder}_{s^l,\sigma_0^l,\sigma_1^l,\sigma_2^l}\Big[1,1,\mathbf P_2^l\Big](k)\Big\|_{L^{\infty}(\mathbb{T}^d)}$ by $\|\mathbf P_2^l\|_{L^4(\mathbb{T}^d)}$, so that in the next iteration,  (I), (II) or (III) can be reused. We also bound $	\Big\| {Q}^{ladder}_{s^1,\sigma_0^1,\sigma_1^l,\sigma_2^l}\Big[1,\mathbf P_1^0,\mathbf P_2^0\Big](k) \Big\|_{L^4(\mathbb{T}^d)}$ using strategy (II).

At the end of the procedure, after taking into account all of the possible graph combination of the ladder operators, we obtain the  bound
\begin{equation}
	\begin{aligned}\label{FinalProof:E10b:1}
	\frac{\tau^q}{q!}	\|{\mathfrak{C}^\lambda_{q,\infty,a,\ell}}\|_{L^4(\mathbb{T}^d)} \ \le \ & |\mathfrak{C}_{final,2}|^qT_*^{q},
	\end{aligned}
\end{equation}
for some constants $\mathfrak{C}_{final,2}>0$ independent of $\lambda$.  We set

\begin{equation}
	\begin{aligned}\label{FinalProof:E10a:2}
		\mathfrak{C}^\ell_{q}(k_1')\ = \ &
		\mathcal{C}_{2}\cdots\mathcal{C}_{q+1},
	\end{aligned}
\end{equation}
where
\begin{equation}
	\begin{aligned}\label{FinalProof:E10:a:1}
		&\mathcal{C}_{m+1}(k_1',\cdots,k'_i,\cdots,k'_{m})=   \sum_{i=1}^{m}\mathcal{C}^\lambda_{i,m+1},
\end{aligned}\end{equation}
and
\begin{equation}
	\begin{aligned}\label{FinalProof:E9:1}
		&\mathcal{C}^\lambda_{i,m+1}(s,k_1',\cdots,k'_i,\cdots,k'_{m})=   \iint_{(\mathbb{T}^d)^2}
		\mathrm{d}k'\mathrm{d}k_{m+1}'|\mathcal{M}(k_{i}',k',k_{m+1}')|^2\\
		&\times\frac{1}{\pi}\delta_{\ell}\Big(\omega(k')+\omega(k_{m+1}')-\omega(k_{i}')\Big)\delta(k'+k_{m+1}'-k_{i}')\Big( \mathfrak{L}^{m+1}(k_1',\cdots,k',\cdots,k'_{m+1})\\
		&-\mathfrak{L}^{m+1}(k_1',\cdots,k'_i,\cdots,k'_{m+1})\mathrm{sign}(k_{i}')\mathrm{sign}(k_{m+1}') -\mathfrak{L}^{m+1}(k_1',\cdots,k'_i,\cdots,k')\mathrm{sign}(k_{i}')\mathrm{sign}(k')\Big).
\end{aligned}\end{equation}
In the above expression, $k'$ takes the position of $k'_i$ in $\mathfrak{L}^{m+1}(k_1',\cdots,k',\cdots,k'_{m+1})$ and of $k'_{m+1}$ in $\mathfrak{L}^{m+1}(k_1',\cdots,$ $ k'_i,\cdots,k')$. Moreover,
\begin{equation}
	\mathfrak{L}^{m+1}(k_1',\cdots,k'_i,\cdots,k'_{m}) \ = \ \mathcal{C}_{m+2}(s,k_1',\cdots,k'_i,\cdots,k'_{m}),
\end{equation}
and
\begin{equation}
	\mathfrak{L}^{q+1}(k_{1}',\cdots,k_{q+1}'):=\prod_{\{i,j\}\in S} {f}(k_{0,i},0).\end{equation}

It is straightforward from \eqref{FinalProof:E10b:1} that
$$\frac{\tau^q}{q!}\|\mathcal{C}_{2}\cdots\mathcal{C}_{q+1}\|_{L^4(\mathbb{T}^d)}\ \le \mathcal{C}_o^q,$$
for some universal constant $\mathcal{C}_o>0$, we then find
\begin{equation}\label{FinalProof:E11}
	\sum_{q=0}^\infty \left\| 	\frac{\tau^q}{q!}\mathfrak{C}_{q}^\ell\right\|_{L^4(\mathbb{T}^d)}\lesssim \sum_{q=0}^\infty  (T_* \mathcal{C}_o)^q,
\end{equation}
which converges for $T_*<\frac{1}{\mathcal{C}_o}$.

Now, let us estimate $\mathfrak Q_{2} $, which we write as \begin{equation}\label{FinalProof:E16}
	\mathfrak Q_{2} 
	\ = \ \sum_{n=[\mathfrak{N}/4]}^{\mathfrak{N}-1}\sum_{S\in\mathcal{P}_{pair}(\{1,\cdots,n+2\})}\mathfrak{G}_{2,n}'(S,t,k,\sigma,\Gamma).
\end{equation}
The same estimates used for $\mathfrak{Q}_1$ can be redone for $\mathfrak{Q}_2$, whose resonance broadening sum  is bounded by
\begin{equation}\label{FinalProof:E17}
	\ \lim_{\lambda\to0}\sum_{2q=[\mathfrak{N}/4]}^{\mathfrak{N}}\varsigma' (\mathcal{C}_oT_*)^q,
\end{equation}
where $\varsigma'$ is defined in \eqref{Def:Para3}, 
which tends to $0$ as $\mathfrak{N}$ tends to infinity and $T_*$ sufficiently small.
The last quantities $\mathfrak Q_3$ also goes to $0$ as $\Phi_{0,i}$ goes to $0$ in the limit of $\lambda\to 0$.

As a result, we obtain the limit
\begin{equation}\label{FinalProof:E18}
	 \sum_{q=0}^\infty \frac{\tau^q}{q!}\mathfrak{C}_{q}^\ell,
\end{equation}
and the convergence \eqref{TheoremMain:1}-\eqref{TheoremMain:2} by a standard Chebyshev's inequality. This series is  a solution of
\eqref{Eq:WT1}.

\medskip

\bibliographystyle{plain}

\bibliography{WaveTurbulence}

\end{document}